\documentclass[12pt]{amsart}
\usepackage{graphicx} 
\usepackage{booktabs}
\usepackage{ragged2e}
\usepackage{float}
\usepackage{soul}
\usepackage[citecolor=Navy]{hyperref}
\usepackage[ruled,vlined]{algorithm2e}
\usepackage{amsfonts,amssymb,amsmath,amsthm}
\usepackage{algorithmic}
\usepackage{amssymb,amscd,amsxtra,calc}
\usepackage{cmmib57}
\usepackage{graphicx}
\usepackage{mathrsfs,enumerate,xcolor}

\usepackage[all]{xy}
\usepackage{multirow} 
\usepackage{psfrag,xmpmulti,amscd,color,pstricks, import}
\usepackage{tikz-cd} 

\newtheorem{thm}{Theorem}[section]

\theoremstyle{definition}

\newtheorem*{rem*}{Remark}

\title[BNQN, Newton's flow, Voronoi's diagram, Random equation]{Backtracking New Q-Newton's method, Newton's flow,  Voronoi's diagram and Stochastic root finding}
\date{\today}

\author[Forn\ae ss]{John Erik Forn\ae ss}
\address{Department of Mathematics, NTNU, Norway }
\email{fornaess@gmail.com}

\author[Hu]{Mi Hu}
\address{Department of Mathematics, University of Oslo, Norway}
\email{humihqu@gmail.com}

\author[Truong]{Tuyen Trung Truong}
\address{Department of Mathematics, University of Oslo, Norway}
\email{tuyentt@math.uio.no}

\author[Watanabe]{Takayuki Watanabe}
\address{Chubu University Academy of Emerging Sciences/Center for Mathematical Science and Artificial Intelligence, Chubu University, Japan}
\email{takawatanabe@isc.chubu.ac.jp}

\begin{document}

\maketitle

{\centering\footnotesize To Professor Steven Krantz on his 70th birthday.\par}

\begin{abstract}

A new variant of Newton's method - named Backtracking New Q-Newton's method (BNQN) - which has strong theoretical guarantee, is easy to implement, and has good experimental performance, was recently introduced by the third author. 

Experiments performed previously showed some remarkable properties of the basins of attractions for finding roots of polynomials and meromorphic functions, with BNQN. In general, they look more smooth than that of Newton's method. 

In this paper, we continue to experimentally explore in depth this remarkable phenomenon, and connect BNQN to Newton's flow and Voronoi's diagram. This link poses a couple of challenging puzzles to be explained. Experiments also indicate that BNQN is more robust against random perturbations than Newton's method and Random Relaxed Newton's method.

\end{abstract}

\section{Introduction}

This paper mainly concerns the problem of finding roots of a polynomial or meromorphic function in 1 complex variable $f(z)=0$. 
 
Newton's method is a well known method to solve equations and optimization problems. In the case of a function $f(z)$ concerns here, the method amounts to the update rule: 
\begin{eqnarray*}
z_{n+1}=z_n-\frac{f(z_n)}{f'(z_n)}.
\end{eqnarray*}

One can also apply the higher dimensional version of Newton's method for systems of equations, to the system $G(x,y)=(Re(f(x+iy)),Im(f(x+iy)))=0$ of real and imaginary parts of $f$, where $x,y\in \mathbf{R}$, i.e. to consider the iterative method
\begin{eqnarray*}
z_{n+1}=z_n-JG(z_n)^{-1}.G(z_n), 
\end{eqnarray*}
where $JG$ is the Jacobian matrix of $G$. However, in the case where $f$ is a meromorphic function, this iterative method is the same as that of the usual Newton's method, thanks to Cauchy-Riemann's equations. 

A new variant of Newton's method - named Backtracking New Q-Newton's method (BNQN) - which has strong theoretical guarantee, is easy to implement, and has good experimental performance, was recently introduced by the third author \cite{RefT}. It is a variant of Newton's method for optimization. For details about the method, also comparisons between it and some well known variants of Newton's method (including Newton's method itself), the readers can see \cite{RefT}\cite{RefTT}. In Section 2 we will recall essential information on BNQN, enough for the purpose of this paper.   

Experiments performed previously in \cite{RefT}, indicate a remarkable phenomenon: The basins of attraction, when BNQN is applied to find roots of meromorphic functions, seem to be more smooth than that of Newton's method. In particular, these basins seem not to be fractal, as is often the case for Newton's method \cite{RefN}. Moreover, when finding roots of a polynomial of degree 2 in 1 complex variable, the picture one obtains for basins of attraction of BNQN looks exactly the same as that in the classical Schr\"oder's theorem for Newton's method \cite{RefSE2}\cite{RefAlex}. 

Which leads to a natural and interesting question: Can we explain this phenomenon rigorously, or is it just a lucky occurrence of numerical calculations?

After presentations of BNQN at ICIAM 2023 (Waseda University) and a satellite conference, a couple of candidates, which have a lot of connections in mathematics and physics, for explaining the pictures obtained in \cite{RefT} have emerged: 

- These seem relevant to Newton's flow. 

- These seem relevant to Voronoi's diagrams. 

More details about Newton's flow and Voronoi's diagrams will be recalled in Section 2. Note that Voronoi's diagram for a pair of distinct points in the plane is the same as what produced in Schr\"oder's theorem for Newton's method. Also, Newton's flow for a polynomial of degree 2 produces exactly the same picture. 

In a previous paper \cite{RefFHTW}, we proved that indeed BNQN applied to finding roots of a polynomial $f(z)$ of degree 2 will produce the same picture as in Schr\"oder's theorem for Newton's method. Even so, BNQN is indeed more smooth than Newton's method in this case: On the boundary line, the dynamics of Newton's method is chaotic, while the dynamics of BNQN has a global attraction. Newton's method for optimization, on the other hand, produces a picture with many black blobs of initial points $z_0$ for which the constructed sequence $\{z_n\}$ does not converge to the roots of $f(z)$ but instead to the critical point of $f$, and the basins of attraction even seem to be fractal!   

In this paper, we will explore in depth the possible links mentioned above between BNQN, Newton's flow and Voronoi's diagrams, when finding roots of meromorphic functions. We verify with many different geometric configurations of the roots, that indeed the link is very apparent and well above a random coincidence. We work with three different versions of Newton's flow, with two versions of Newton's method, as well as Random Relaxed Newton's method. For the readers's convenience, we will recall in enough detail on Newton's method and Random Relaxed Newton's method in the next Section. It seems that the picture one obtains from BNQN is a better version of what one obtains with Newton's method for Optimization, and while BNQN is not a flow method it somehow attains many good features of flow methods. From these experiments, another interesting new phenomenon seems to occur when one applies Newton's flow to $f/f'$, even if one only wants to find roots of $f$ and even if $f$ has only roots of multiplicity $1$. 

An interesting possible application of the study in this paper is as follows. In the literature, there is no unique construction of Voronoi's diagrams when two or more points meet together (limits of Voronoi's diagrams) \cite{RefL}. In such a case, we expect that a curved version of a ''canonical" Voronoi's diagram can be obtained by using Newton's flow or BNQN. For example, working with the polynomial $P(z)=(z-a)^2(z-b)(z-c)$ can provide one some insights into how Voronoi's diagram for three points $a$ (with multiplicity 2), $b$ (with multiplicity $1$) and $c$ (with multiplicity 1) should approximately look like. In the experiments, we have several such examples. 

Given that when using computers one cannot avoid dealing with random errors/perturbations, and that solving equations with random perturbations is interesting itself, we also present some experiments concerning this issue. The experiments indicate that BNQN is more robust compared to Newton's method and Random Relaxed Newton's method.

{\bf Plan of the paper:} In Section 2 we recall relevant information concerning the methods, as well as Newton's flow and Voronoi's diagrams, as well as stochastic root finding. In Section 3, we introduce the setting of our experiments and then present the experimental results. In the final section, we draw some conclusions and some challenging puzzles and directions for future's research. 

{\bf Acknowledgments:} Mi Hu and Tuyen Trung Truong are partially supported by Young Research Talents grant 300814 from Research Council of Norway. Takayuki Watanabe is partially supported by JSPS Grant-in-Aid for Early-Career Scientists Grant Number JP 23K13000. Some ideas of the work were discussed/carried out in the inspiring environments of the conferences ICIAM 2023 (Waseda University) and ``Birational geometry and algebraic dynamics" 2023 (University of Tokyo), for which we would like to thank the organisers and participants, in particularly Mark Comerford for suggesting the connection to Voronoi's diagrams.

\section{Preliminaries} In this section, we briefly review some previous results and algorithms, with the emphasis on properties which will be used later.  

\subsection{Newton's method}

Besides the main version in the introduction, here we recall  Newton's method in the optimization setting (Newton's method for optimization) also, to easily compare later to BNQN. Assume that one wants to find (local) minima of an objective function $F:\mathbf{R}^m\rightarrow \mathbf{R}$. Let $\nabla F$ be the gradient of $F$, and $\nabla ^2F$ be the Hessian of $F$. Then Newton's method for the function $F$ is the following iterative algorithm: Choose $z_0\in \mathbf{R}^m$ an initial point, and define: 
$$z_{n+1}=z_n-(\nabla ^2F(z_n))^{-1}.\nabla F(z_n).$$

In this paper, Newton's method for optimization will be applied to the function $F(x,y)=|f(x+iy)|^2/2$, for a meromorphic function $f$ in 1 complex variable $z=x+iy$. 

Both versions of Newton's method are invariant under a linear change of coordinates. However, their behaviour can be very different, see the experimental results for detail. 

Concerning finding roots of polynomials, we note that Newton's method gives rise to an algebraic dynamical system on the Riemann sphere. As such, \cite{RefMc} shows that Newton's method does not guarantee to always find roots of polynomials of degree at least $4$, see the cited paper for the precise statement. For Newton's method for optimization, it has the tendency to converge to the nearest (non-degenerate) critical point of $F$, and hence also has no guarantee for finding roots. See the experimental results for detail.

\subsection{Random relaxed Newton's method} A direct variant of Newton's method is Relaxed Newton's method. It works as follows. We choose a nonzero  complex number $\alpha$, and consider the iterative procedure: 
\begin{eqnarray*}
z_{n+1}=z_n-\alpha \frac{f(z_n)}{f'(z_n)}.
\end{eqnarray*}
Again, this gives rise to an algebraic dynamical system on the Riemann sphere, and hence by \cite{RefMc} does not have guarantee to finding roots of polynomials. However, see the discussion in the subsection on Newton's flow about how the basins of attraction for Relaxed Newton's method behave when $\alpha $ goes to $0$. A surprising fact is that if one adds {\bf randomness} into the design, then the method can find roots of polynomials. The iterative procedure, named Random Relaxed Newton's method, in this case is the following: 
\begin{eqnarray*}
z_{n+1}=z_n-\alpha _n \frac{f(z_n)}{f'(z_n)},
\end{eqnarray*}
where $\alpha _n$ is randomly chosen in an appropriate manner. 

We recall here the relevant result in \cite{RefS} on convergence of Random Relaxed Newton's method when applied to finding roots of polynomials. 

\begin{thm} Let $0.5<\rho <1$ be a constant. Let $P(z)$ be a polynomial in 1 complex variable $z$.

Let $\alpha _n$ be randomly chosen from the uniform distribution on the set $\{\alpha \in \mathbf{C}:~|\alpha _n-1|\leq \rho \}$. 

Then the Random Relaxed Newton's method, applied to finding roots of the polynomial $P(z)$, will converge to a root of $P(z)$ for all - except a finite number of exceptions - initial points $z_0$.  

\label{TheoremSumi}\end{thm}

\subsection{Backtracking New Q-Newton's method}
In this section, we recall enough details on both theoretical and practical aspects of BNQN.

\subsubsection{Heuristics and some delicate issues} A tendency of Newton's method for optimization, applied to find minima of a function $F:\mathbf{R}^m\rightarrow \mathbf{R}$, is that if the initial point $z_0$ is close  to a (non-degenerate) critical point $z^*$ of $F(z)$, then the sequence constructed by Newton's method will converge to $z^*$. Hence, if $z^*$ is a saddle point or local maximum, then Newton's method is not desirable for finding (local) minima. In the special case where $F(z)$ is a quadratic function whose Hessian is invertible, then for every initial point $z_0$, the point $z_1$ constructed by Newton's method will be the unique critical point $z^*=0$.

To overcome this undesirable behaviour of Newton's method, a new variant called New Q-Newton's method (NQN) was recently proposed in \cite{RefTT}. This consists of the following two main ideas: 

- Add a perturbation $\delta ||\nabla F(z)||^{\tau}Id$ to the Hessian $\nabla ^2 F(z)$, where $\tau >0$ is a constant, and $\delta$ is appropriately chosen within a previously chosen set $\{\delta _0,\ldots ,\delta _m\}\subset \mathbf{R}$. Thus, we consider a matrix $A(z)=\nabla ^2F(z)+\delta ||\nabla F(z)||^{\tau}Id$, instead of $\nabla ^2F(z)$ as in the vanilla Newton's method. This has two benefits. First,  it avoids the (minor) difficulty one encounters in Newton's method if $\nabla^2F(z)$ is a singular matrix. Second, more importantly, it turns out that if the $\delta _0,\ldots ,\delta _m$ are randomly chosen, then NQN can avoid saddle points. Since the perturbation $\delta ||\nabla F(z)||^{\tau}Id$ is negligble near non-degenerate local minima, the rate of convergence of NQN is the same as that of Newton's method there.  

- If one mimics Newton's method, in defining $z_{n+1}=z_n-A(z_n)^{-1}.\nabla F(z)$, then the constructed sequence still has the same tendency of convergence to the nearest critical point of $F(z)$. This is remedied in NQN by the following idea: We let $B(z)$ be the matrix with the same eigenvectors as $A(z)$, but whose eigenvalues are all absolute values of the corresponding eigenvalues of $A(z)$. The update rule of NQN is $z_{n+1}=z_n-B(z_n)^{-1}.\nabla F(z_n)$. 

The precise definition of the algorithm NQN is given in the next subsection. If $F(z)$ is a $C^3$ function, then NQN applied to $F(z)$ can avoid saddle points, and has the same rate of convergence near non-degenerate local minima. 

However, NQN has no global convergence guarantee. In \cite{RefT}, a variation of NQN, that is BNQN, was introduced and shown to keep the same good theoretical guarantees as NQN, with the additional bonus that global convergence can be proven for BNQN for very general classes of functions $F(z)$: functions which have at most countably many critical points, or functions which satisfy a Lojasiewicz gradient inequality type. The main idea is to incorporate Armijo's Backtracking line search \cite{RefAr} into NQN. For the readers' convenience, we recall here the idea of Armijo's Backtracking line search. Let $F:\mathbf{R}^m\rightarrow \mathbf{R}$ be a $C^1$ function. Let $z,w\in \mathbf{R}$ such that $<\nabla F(z),w>$ is strictly positive. Then, there exists a positive real number $\gamma$ for which $F(x-\gamma w)-F(x)\leq -<\nabla F(x),\gamma w >/3$. If we choose $\gamma $ by  a backtracking manner (that is, start $\gamma$ from a positive number, and then reduce it exponentially until the above mentioned inequality is satisfied), then the procedure is called Armijo's Backtracking line search.

We notice that there are some delicate issues when doing this incorporation between NQN and Armijo's Backtracking line search, which rely crucially on the fact that we are using a perturbation of the Hessian matrix here. First, in BNQN, the choice of $\delta$ from among $\{\delta _0,\ldots ,\delta _m\}$ is more complicated than NQN. Second, the analog Backtracking line search for Gradient descent is not yet known to be able to avoid saddle points, even though there is a slight variant which can avoid saddle points, and experiments support that the method should be able to avoid saddle point. Third, a priori the learning rate one finds by Armijo's Backtracking line search can be smaller than 1, and hence the rate of convergence can a priori slower than being quadratic.  

The precise definition of BNQN is given in the next subsection. 

\subsubsection{Algorithms: NQN and BNQN} Here we present the basic versions of NQN and BNQN. Many more variations can be found in \cite{RefT}. 

Let $A:\mathbb{R}^m\rightarrow \mathbb{R}^m$ be an invertible {\bf symmetric} square matrix. In particular, it is diagonalisable.  Let $V^{+}$ be the vector space generated by eigenvectors of positive eigenvalues of $A$, and $V^{-}$ the vector space generated by eigenvectors of negative eigenvalues of $A$. Then $pr_{A,+}$ is the orthogonal projection from $\mathbb{R}^m$ to $V^+$, and  $pr_{A,-}$ is the orthogonal projection from $\mathbb{R}^m$ to $V^-$. As usual, $Id$ means the $m\times m$ identity matrix.  

First, we introduce NQN \cite{RefTT}. 

\medskip
{\color{blue}
 \begin{algorithm}[H]
\SetAlgoLined
\KwResult{Find a minimum of $F:\mathbb{R}^m\rightarrow \mathbb{R}$}
Given: $\{\delta_0,\delta_1,\ldots, \delta_{m}\}\subset \mathbb{R}$\  and $\alpha >0$;\\
Initialization: $z_0\in \mathbb{R}^m$\;
 \For{$k=0,1,2\ldots$}{ 
    $j=0$\\
    \If{$\|\nabla f(z_k)\|\neq 0$}{
   \While{$\det(\nabla^2f(z_k)+\delta_j \|\nabla f(z_k)\|^{1+\alpha}Id)=0$}{$j=j+1$}}

$A_k:=\nabla^2f(z_k)+\delta_j \|\nabla f(z_k)\|^{1+\alpha}Id$\\
$v_k:=A_k^{-1}\nabla f(z_k)=pr_{A_k,+}(v_k)+pr_{A_k,-}(v_k)$\\
$w_k:=pr_{A_k,+}(v_k)-pr_{A_k,-}(v_k)$\\
$x_{k+1}:=x_k-w_k$
   }
  \caption{New Q-Newton's method} \label{table:alg}
\end{algorithm}
}
\medskip

BNQN includes a more sophisticated choice of $\delta$ in NQN, together with a combination of Armijo's Backtracking line search. For a symmetric, square real matrix $A$, we define: 
  
  $sp(A)=$ the maximum among $|\lambda |$'s, where $\lambda  $ runs in the set of eigenvalues of $A$, this is usually called the spectral radius in the Linear Algebra literature;
  
  and 
  
  $minsp(A)=$ the minimum among $|\lambda |$'s, where $\lambda  $ runs in the set of eigenvalues of $A$, this number is non-zero precisely when $A$ is invertible.
  
 One can easily check the following more familiar formulas: $sp(A)=\max _{\|e\|=1}\|Ae\|$ and $minsp(A)=\min _{\|e\|=1}\|Ae\|$, using for example the fact that $A$ is diagonalisable.  

We recall that a function $F$ has compact sublevels if for all $C\in \mathbf{R}$ the set $\{z:~F(z)\leq C\}$ is compact. 

\medskip
{\color{blue}
 \begin{algorithm}[H]
\SetAlgoLined
\KwResult{Find a minimum of $F:\mathbb{R}^m\rightarrow \mathbb{R}$}
Given: $\{\delta_0,\delta_1,\ldots, \delta_{m}\} \subset \mathbb{R}$\, $0<\tau $ and $0<\gamma _0\leq 1$;\\
Initialization: $z_0\in \mathbb{R}^m$\;
$\kappa:=\frac{1}{2}\min _{i\not=j}|\delta _i-\delta _j|$;\\
 \For{$k=0,1,2\ldots$}{ 
    $j=0$\\
  \If{$\|\nabla F(z_k)\|\neq 0$}{
   \While{$minsp(\nabla^2F(z_k)+\delta_j \|\nabla F(z_k)\|^{\tau}Id)<\kappa  \|\nabla F(z_k)\|^{\tau}$}{$j=j+1$}}
  
 $A_k:=\nabla^2f(z_k)+\delta_j \|\nabla f(z_k)\|^{\tau}Id$\\
$v_k:=A_k^{-1}\nabla f(z_k)=pr_{A_k,+}(v_k)+pr_{A_k,-}(v_k)$\\
$w_k:=pr_{A_k,+}(v_k)-pr_{A_k,-}(v_k)$\\
$\widehat{w_k}:=w_k/\max \{1,\|w_k\|\}$\\
(If $F$ has compact sublevels, then one can choose $ \widehat{w_k}=w_k$).\\
$\gamma :=\gamma _0$\\
 \If{$\|\nabla f(z_k)\|\neq 0$}{
   \While{$f(z_k-\gamma \widehat{w_k})-f(z_k)>-\gamma \langle\widehat{w_k},\nabla f(z_k)\rangle/3$}{$\gamma =\gamma /3$}}

$z_{k+1}:=z_k-\gamma \widehat{w_k}$
   }
  \caption{Backtracking New Q-Newton's method} \label{table:alg0}
\end{algorithm}
}
\medskip

Algorithm \ref{table:alg0} has two different versions depending whether the objective function $F$ has compact sublevels or not. In \cite{RefFHTW}, we introduced a new variant, named BNQN New Variant, with a new parameter, which includes both these versions as special cases. Indeed, $\theta =0$ in BNQN New Variant is the version of BNQN for functions with compact sublevels, and $\theta =1$ in BNQN New Variant is the version of BNQN for general functions. 

\medskip
{\color{blue}
 \begin{algorithm}[H]
\SetAlgoLined
\KwResult{Find a minimum of $F:\mathbb{R}^m\rightarrow \mathbb{R}$}
Given: $\{\delta_0,\delta_1,\ldots, \delta_{m}\} \subset \mathbb{R}$\, $\theta \geq 0$, $0<\tau $ and $0<\gamma _0\leq 1$;\\
Initialization: $z_0\in \mathbb{R}^m$\;
$\kappa:=\frac{1}{2}\min _{i\not=j}|\delta _i-\delta _j|$;\\
 \For{$k=0,1,2\ldots$}{ 
    $j=0$\\
  \If{$\|\nabla F(z_k)\|\neq 0$}{
   \While{$minsp(\nabla^2F(z_k)+\delta_j \|\nabla F(z_k)\|^{\tau}Id)<\kappa  \|\nabla F(z_k)\|^{\tau}$}{$j=j+1$}}
  
 $A_k:=\nabla^2f(z_k)+\delta_j \|\nabla f(z_k)\|^{\tau}Id$\\
$v_k:=A_k^{-1}\nabla f(z_k)=pr_{A_k,+}(v_k)+pr_{A_k,-}(v_k)$\\
$w_k:=pr_{A_k,+}(v_k)-pr_{A_k,-}(v_k)$\\
$\widehat{w_k}:=w_k/\max \{1,\theta \|w_k\|\}$\\
$\gamma :=\gamma _0$\\
 \If{$\|\nabla f(z_k)\|\neq 0$}{
   \While{$f(z_k-\gamma \widehat{w_k})-f(z_k)>-\gamma \langle\widehat{w_k},\nabla f(z_k)\rangle/3$}{$\gamma =\gamma /3$}}

$z_{k+1}:=z_k-\gamma \widehat{w_k}$
   }
  \caption{Backtracking New Q-Newton's method New Variant} \label{table:alg2}
\end{algorithm}
}
\medskip

\subsubsection{A main theoretical result for finding roots of meromorphic functions by BNQN}

As mentioned previously, BNQN has strong theoretical guarantees for some big classes of functions in any dimension. However, to keep the presentation concise, here we present only one main result relevant to the question pursued in this paper, that of finding roots of meromorphic functions in 1 complex variable. For more general results, the readers can consult \cite{RefT}, \cite{RefTT}. 

\begin{thm} Let $g(z):\mathbf{C}\rightarrow \mathbf{P}^1$ be a non-constant meromorphic function. 
Define a function $F:\mathbf{R}^2\rightarrow [0,+\infty]$ by the formula $F(x,y)=|g(x+iy)|^2/2$. 

Given an initial point $z_0\in \mathbf{C}$, which is not a pole of $g$, we let $\{z_n\}$ be the sequence constructed by BNQN New Variant applied to the function $F$ with initial point $z_0$. If the objective function has compact sublevels (i.e. for all $C\in \mathbf{R}$ the set $\{(x,y)\in \mathbf{R}^2:~F(x,y)\leq C\}$ is compact), we choose $\theta \geq 0$, while in general we choose $\theta >0$. 

1) Any critical point of $F$ is a root of $g(z)g'(z)=0$. 

2) If $z^*$ is a cluster point of $\{z_n\}$ (that is, if it is the limit of a subsequence of $\{z_n\}$), then $z^*$ is a critical point of $F$. Moreover, in this case the whole sequence $\{z_n\}$ converges to $z^*$. 

3) If $F$ has compact sublevels, then $\{z_n\}$ converges. 

4) Assume that the parameters $\delta _0,\delta _1,\delta _2$ in BNQN New Variant are randomly chosen. Assume also that $g(z)$ is generic, in the sense that $\{z\in \mathbf{C}:~g(z)g"(z)=g'(z)=0\}=\emptyset$. There exists an exceptional set $\mathcal{E}\subset \mathbf{C}$ of zero Lebesgue measure so that if $z_0\in \mathbf{C}\backslash \mathcal{E}$, then $\{z_n\}$ must satisfy one of the following two options: 

Option 1: $\{z_n\}$ converges to a root $z^*$ of $g(z)$, and if $\gamma _0=1$ in the algorithm then the rate of convergence is quadratic. 

Option 2: $\lim _{n\rightarrow\infty}|z_n| =+\infty$. 

\label{TheoremMeromorphic}\end{thm}

For BNQN, Theorem \ref{TheoremMeromorphic} is proven in \cite{RefTT}. The stated version for BNQN New Variant is from \cite{RefFHTW}. An example for which part 3 of Theorem \ref{TheoremMeromorphic} can apply is when $g$ is a polynomial or $g=P/Q$ where $P,Q$ are polynomials and $P$ has bigger degree than $Q$ (a special case is when $Q=P'$, in which case the zeros of $g$ are exactly that of $P$, with the advantage that they all have multiplicity $1$).  

\subsubsection{Implementation details} An implementation in Python of BNQN accompanies the paper \cite{RefT}. The implementation is flexible in that one does not need precise values of the gradient and Hessian of the function $F$, but approximate values are good enough. Experiments also show that the performance of BNQN is very stable, with respect to its parameters and to the values of the objective function and its gradient and Hessian matrix. Some experiments later in this paper will illustrate this.   

\subsection{Newton's flow} A general strategy when studying an iterative method is that one also studies its flow counterpart. Usually, the flow counterpart is a smooth version of the (discrete) iterative method.  
In this subsection, we briefly review several variants of Newton's flow. 
Newton's flow means solutions of certain ordinary differential equations which we will define later. 

The authors observed that the basins of attractions for BNQN method look smoother than that for Newton’s method or Backtracking line search for Gradient Descent. 
This observation inspired us to compare the basin structure for BNQN with that for Newton's flow. 
In the following,  we consider three types of continuous-time dynamical systems.

\subsubsection{Newton's flow for $f$}
The first Newton-type flow is the ordinary differential equation of complex-valued function $z(t) = x(t) + i y(t)$ of one real variable $t$: 
\begin{align}\label{eq:ODE}
    \frac{\mathrm{d}z}{\mathrm{d}t} = - \frac{f(z)}{f'(z)}
\end{align}
with initial value $z(0) = z_0$. 
Here, $f(z)$ is a holomorphic function.  
Note that the right-hand side of (\ref{eq:ODE}) has singularity at $\{f' = 0\} \setminus \{f=0\}$. 

If a global solution $z(t)$ for the ODE above exists, 
then it satisfies 
\begin{align*}
\frac{\mathrm{d}}{\mathrm{d}t}f(z(t))
=f'(z(t))\frac{\mathrm{d}z}{\mathrm{d}t} 
= - f(z). 
\end{align*}
Thus, $f(z(t)) = f(z_0)e^{-t}$ 
and $f(z(t)) \to 0$ as $t \to \infty$. 
By this calculation, 
we expect that the solution of the ODE (\ref{eq:ODE}) approximates a zero of the target function $f.$

Neuberger {\it et al.} showed that 
we can interpret the ODE (\ref{eq:ODE}) as the Sobolev gradient of $F(z) =\|f(z)\|^2/2$ relative to the natural Riemannian metric. 
See Appendix A of \cite{RefNeub}. 

From our point of view, 
it is worth noting that the Euler discretization of the ODE (\ref{eq:ODE}) gives the Relaxed Newton's method. 
Let $\alpha >0$ be a constant time step. 
The Euler discretization gives us 
\begin{align*}\label{eq:Euler}
    \frac{z(t+\alpha) - z(t)}{\alpha} = - \frac{f(z(t))}{f'(z(t))}. 
\end{align*}
with initial condition $z(0)=z_0$. 
Thus, we have the recursive formula 
$z(t+\alpha)= z(t) - \alpha f(z(t))/f'(z(t)),$ 
which is the relaxed Newton's method. 
In particular, if the time step $\alpha$ is $1,$
then the Euler discretization is precisely the same as the usual Newton's method. 

It is shown in \cite{RefHK}\cite{RefMe} that if $f$ is a rational function, then the union of basins of attraction for Newton's flow of all the roots of $f$ has full measure. It is also shown that in case $f$ is a polynomial, there is a small value of $\alpha$ for which Relaxed Newton's method has convergence guarantee to find all roots of $f$. Note that the choice of $\alpha$ heavily depends on the polynomial $f$. See \cite{RefBer} for a good review on some results on this direction, as well as dynamics of meromorphic functions in general. 

Compared to the original Newton's method, 
the local convergence rate of Newton's flow is slower.
Namely, 
if $z=a$ is a simple zero of $f$, 
then we have $|z(t) - a| = O(e^{-t})$ as $t \to \infty.$ 
Thus, Newton's flow converges to a simple root with order $1$.  
On the other hand, the original Newton's method converges to a simple root with order $2$ as is well-known. 

From now on, this method will be called Newton's flow. 

\subsubsection{Newton's flow for $f/f'$}
The second Newton-type flow is the ordinary differential equation of complex-valued function $z(t)$ of one real variable $t$: 
\begin{equation}\label{eq:ODEwFirstD}
    \frac{\mathrm{d}z}{\mathrm{d}t} = - \frac{g(z)}{g'(z)}
\end{equation}
with initial value $z(0) = z_0$. 
Here, $g(z) =f(z)/f'(z)$ is a meromorphic function and $f(z)$ is a holomorphic function. The roots of $g$ are exactly the roots of $f$, with the advantage that their multiplicities are all $1$. The ODE (\ref{eq:ODEwFirstD}) is obtained by substituting $g=f/f'$ for $f$ in the ODE (\ref{eq:ODE}). 

By a similar calculation, 
we can show that $g(z(t)) = g(z_0)e^{-t}\to 0$ as $t \to \infty$ 
if a global solution $z(t)$ exists. 
Since the zeros of $g$ is same as the zeros of $f,$ 
we expect the flow $z(t)$ will converge to a root of $f.$ 

The conceptual difference between two ODEs (\ref{eq:ODE}) and (\ref{eq:ODEwFirstD}) is the singular sets. 
The singular set for (\ref{eq:ODE}) is $\{f'=0\}\setminus \{f=0\},$ 
while the singular set for (\ref{eq:ODEwFirstD}) is contained in $\{(f')^2-ff''=0\}.$  

From now on, this version will be named Newton's flow vFraction. 

\subsubsection{Newton's flow for the optimization version} The last Newton's flow is applicable to Newton's method v2, the optimization version of Newton's method. This Newton-type flow is the ordinary differential equations related to the  function $z(t) = x(t) + i y(t)$ of one real variable $t$: 
\begin{align}\label{eq:ODEwHessian}
 \frac{\mathrm{d}}{\mathrm{d}t}
\left( \begin{array}{c}
	 x\\
	 y\\
\end{array} \right)
= - \left( \nabla^2 F(z) \right)^{-1} \nabla F(z)
\end{align}
with initial value $z(0) = z_0$. 
Here, $F(z) =\|f(z)\|^2/2$ and $f(z)$ is a holomorphic function.  
Also, by identifying $z=x+iy$, 
we denote by 
\begin{equation*}
\begin{aligned}
\nabla^2 F =
\left( \begin{array}{cc}
	 F_{xx} & F_{xy}\\
	 F_{xy} & F_{yy}\\
\end{array} \right)
\end{aligned}
\text{ and }
\nabla F =
\left( \begin{array}{c}
	 F_x\\
	 F_y\\
\end{array} \right)
\end{equation*}
where $F_x = \frac{\partial F}{\partial x}, \ F_y = \frac{\partial F}{\partial y}$ and so on. 
The right-hand side of (\ref{eq:ODEwHessian}) is defined outside the singular set $\{\det \nabla^2 F = 0\}$.

If this ODE has a global solution $z(t)$, 
then we have 
$$\frac{\mathrm{d}}{\mathrm{d}t} \nabla F(z(t)) 
= \nabla^2 F(z) \frac{\mathrm{d}z}{\mathrm{d}t} 
= - \nabla F(z(t)). $$
Therefore, the global solution satisfies $\nabla F(z(t)) = e^{-t} \nabla F(z_0)$, 
and hence we have $\nabla F(z(t)) \to 0$ as $t \to \infty$. 
This shows that we can expect to obtain an approximation of points where $\nabla F$ vanishes by solving the ODE.

We are now interested in the local behavior near simple zeros of $f$. 
Suppose $f(a)=0.$ 
We can assume that $a=0$ by changing the coordinate. 
Since $f$ is holomorphic, 
we can write $f(z) = cz + g(z)$ 
where $c \neq 0$ is a complex number and $g(z)$ is a holomorphic function defined near $z=0$ which satisfies $g(0) = g'(0)=0$ and $g''(0) \neq 0.$ 
Let us denote $f(z) = u(z) + i v(z)$ using real functions. 
Then $u(0)=0$, $v(0)=0$, and  $u_x(0) + i v_x(0)= c$ since $f'(z) = u_x(z)+i v_x(z).$ 
Moreover, $(u_{xx}(0), v_{xx}(0))\neq (0, 0)$ since $f''(z) = u_{xx}(z)+i v_{xx}(z).$
Since $F(z) = \|f(z)\|^2/2 = (u^2+v^2)/2$, 
we have 
\begin{equation*}
\begin{aligned}
\nabla F =
\left( \begin{array}{c}
	 u u_x + v v_x\\
	 u u_y + v v_y\\
\end{array} \right), 
\text{ and hence }
\nabla^2 F(0) =
\left( \begin{array}{cc}
	 |c|^2 & 0 \\
	 0 & |c|^2\\
\end{array} \right). 
\end{aligned}
\end{equation*}
Therefore, the ODE (\ref{eq:ODEwHessian}) near $z=0$ is reduced to be  
\begin{equation*}
\begin{aligned}
\frac{\mathrm{d}}{\mathrm{d}t}
\left( \begin{array}{c}
	 x\\
	 y\\
\end{array} \right)
=
-
\left( \begin{array}{c}
	 x\\
	 y\\
\end{array} \right)
+o(\sqrt{x^2+y^2}). 
\end{aligned}
\end{equation*}
Thus, 
the flow that starts near a simple root has a time-global solution and will converge to the simple root. 

From now on, this Newton's flow is named Newton's flow vOptimization. 

\subsubsection{Implementation details}
In this subsubsection, 
we present our implementation of three Newton-type flows.
There are various known methods for solving differential equations numerically, but in this paper, we use the Explicit Runge-Kutta method of order 5(4). 
This method is the default method used in \texttt{scipy.integrate.solve\_ivp} function. 

The Explicit Runge-Kutta method of order 5(4) is discussed in \cite{RefDP}. 
This embedded method automatically adjusts the time step size to enable fast and accurate numerical computation, but for simplicity, we review the classical Runge-Kutta method here. 

We want to numerically solve the autonomous ODE 
$$ \dot{y}=g(y),\quad y(0)=y_{0}$$
with a given initial condition. Here, $y$ and $g$ are $\mathbf{R}^n$-valued, or they are $\mathbf{C}$-valued functions.  
Fix a small time step $h > 0$, and define $Y(0) =y_0$. 
If we already define the value $Y(t)$, 
then define 
$$Y(t+h)=Y(t) + \frac{h}{6}(k_1 + 2k_2 + 2k_3 + k_4)$$
where 
\begin{align*}
    k_1&=g(Y(t)), \\
    k_2&=g(Y(t)+\frac{h}{2}k_1), \\
    k_3&=g(Y(t)+\frac{h}{2}k_2), \\
    k_4&=g(Y(t)+ hk_3). 
\end{align*}
It is known that the truncation error is proportional to $h^5$ for every time step (if the solution $y(t)$ is sufficiently smooth).  
Thus, if we calculate until $t = T,$ 
then the error $|y(T) -Y(T)|$ is proportional to $h^5\cdot T/h = T h^4.$ 
See Section 5.4 and following of \cite{RefBFR} for more details. 

\begin{figure}
    \centering
    \includegraphics[width=3cm]{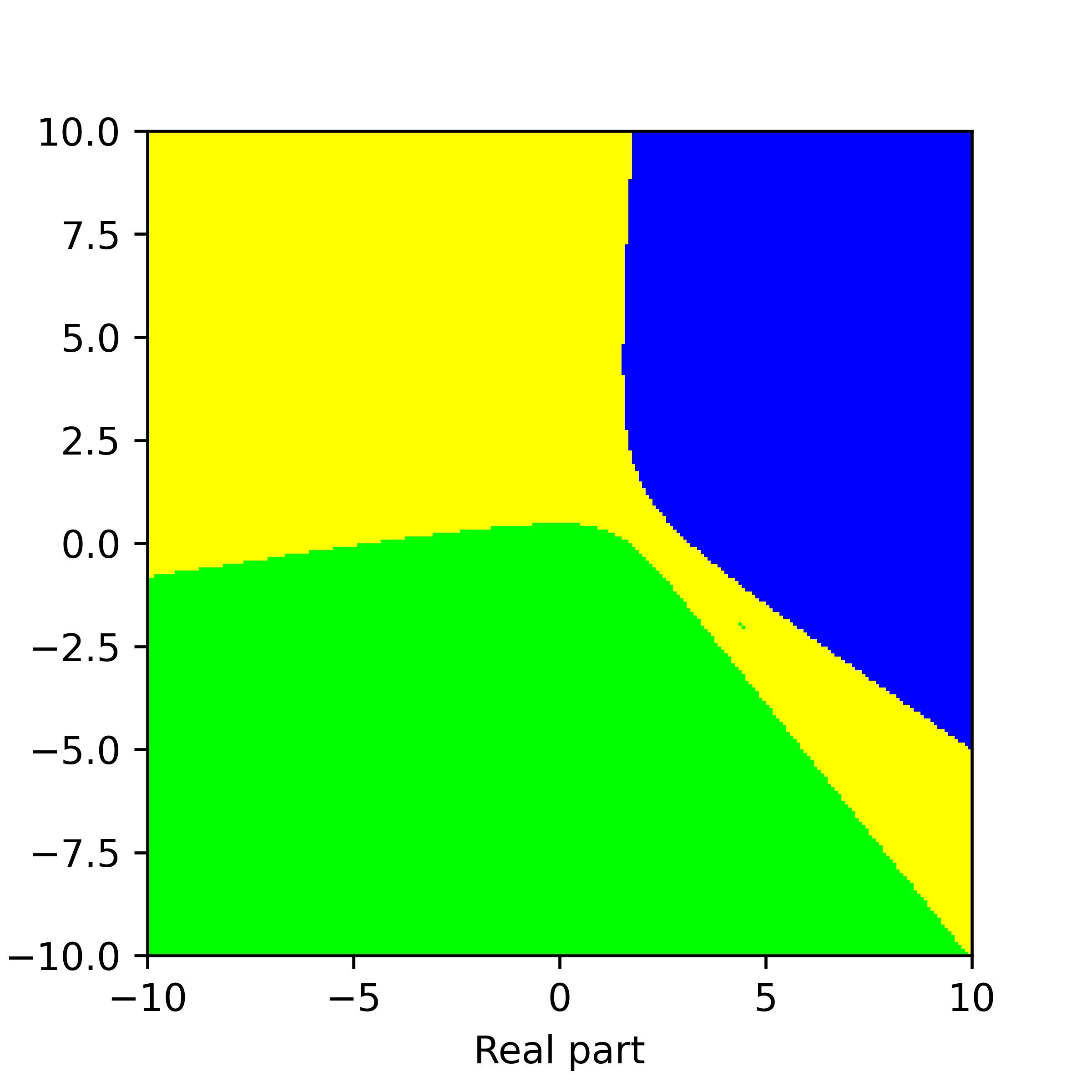}
    \includegraphics[width=3cm]{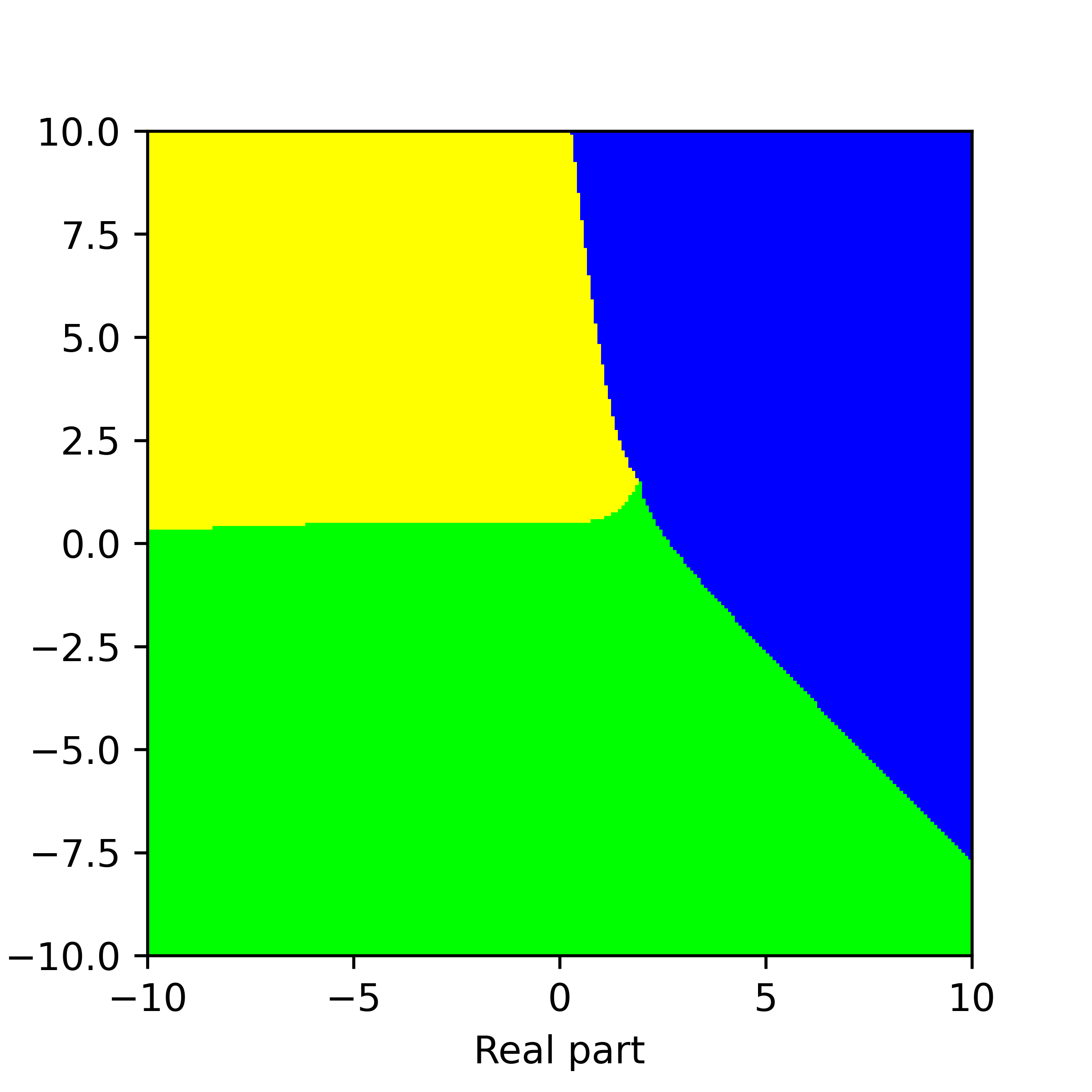}
    \includegraphics[width=3cm]{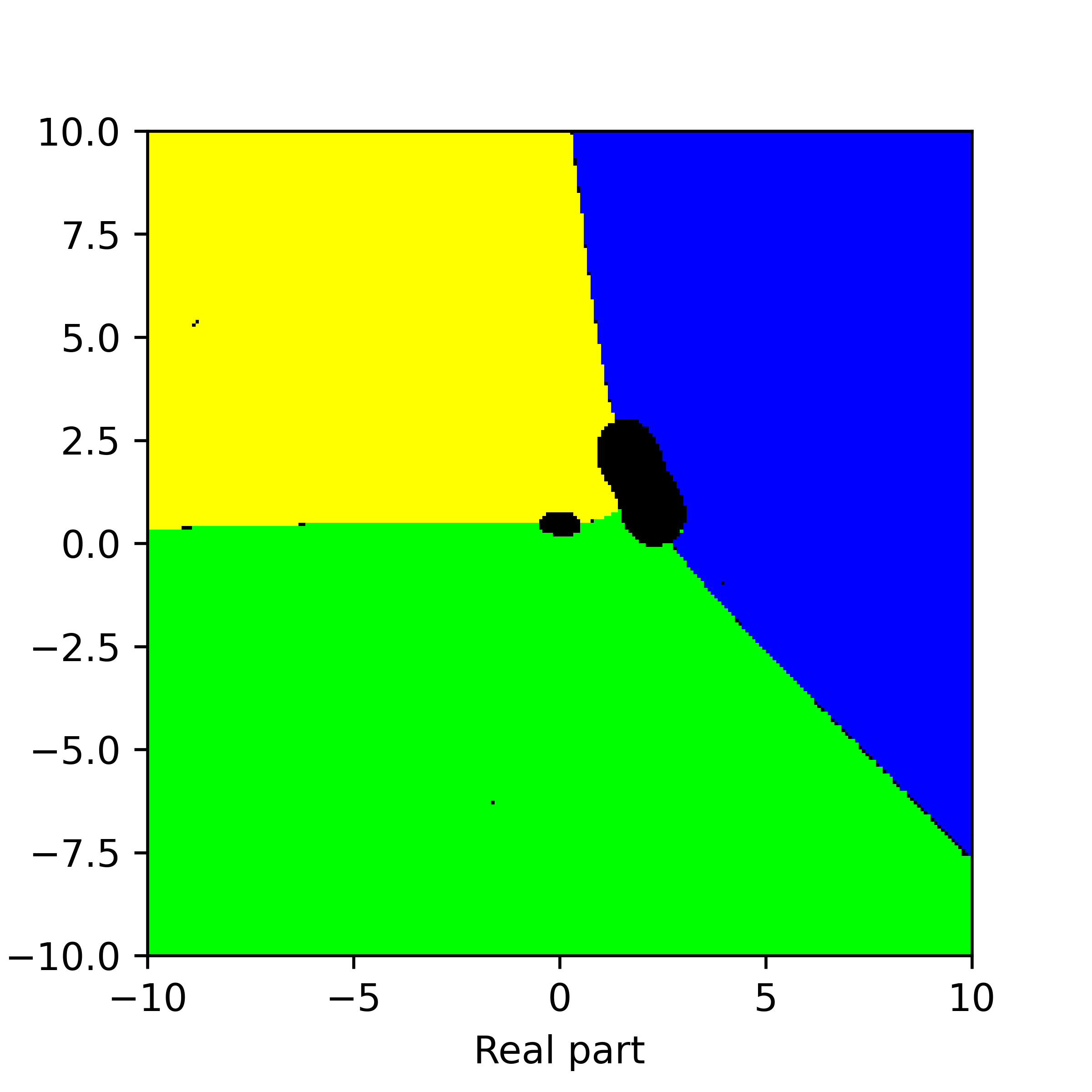}
    \caption{Basins of attraction for the ODEs (\ref{eq:ODE}), (\ref{eq:ODEwFirstD}), and (\ref{eq:ODEwHessian}) from left to right, applied to a polynomial of degree $3$ (the function $f_1$ listed in Section 3). }
    \label{fig:NewtonFlowF1}
\end{figure}

We are interested in the global structure of basins for Newton's flow. 
Therefore, we performed the following experiment for every three ODEs above. 
For several initial values, 
we numerically integrated the trajectories using \texttt{scipy.integrate.solve\_ivp}. 
The time step was $h=0.01$ and the endpoint of integration was $T=100$.  
The same color was applied to the initial values whose trajectory reached within a certain distance from the same root of $f$. 
If the trajectories are away from all the roots, then we colored their initial values black. 
An illustration is shown in Fig. \ref{fig:NewtonFlowF1}.

\begin{figure}
    \centering
    \includegraphics[width=5.5cm]{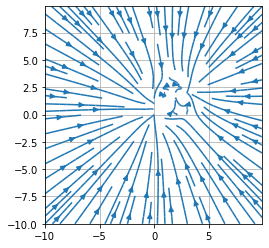}
    \includegraphics[width=5cm]{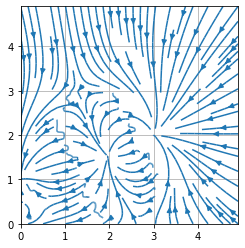}
    \caption{Streamline plots for the ODE (\ref{eq:ODEwHessian}), from Figure \ref{fig:NewtonFlowF1}. The figure on the right is an enlarged image of the figure on the left, whose range is $0 \leq \mathrm{Re} z,\ \mathrm{Im} z \leq 5$.  We can see a stable equilibrium point near  $z = 1.966755+1.516588i$ which does not correspond to the roots of $f_1$, but a critical point of $f_1$. }
    \label{fig:streamlineplotF1}
\end{figure}

\subsection{Voronoi's diagrams} Voronoi's diagrams are a classical object \cite{RefV1}\cite{RefV2}, with an ever increasing domain of applications: geophysics, surface metrology, hydrology, archaeology, dialectometry, political sciences, biology, ecology, ethology, computational chemistry, astrophysics, computational fluid dynamics, computational physics, medical diagnosis, epidemiology, polymer physics, materials science, aviation, architecture, urban planing, mining, robotics, wireless network, computer graphics, machine learning, and a lot more. For example, the skin pattern of a giraffe has Voronoi-type diagrams.  A good place to start with Voronoi's diagrams and applications is \cite{RefVW} and its references. 

There are many versions of Voronoi's diagrams. In this paper, we consider the most basis version, that of a finite set of points $z_1,\ldots ,z_m$ in $\mathbf{R}^2$, with the Euclidean distance $d_E(.,.)$. If $z\in \mathbf{R}^2$ is such that there is a unique $i\in \{1,\ldots ,m\}$ for which $d_E(z,z_i)<\min _{j\not= i}d_E(z,z_j)$, then we say that $z$ belongs to the Voronoi's cell of the point $z_i$. The boundary consists of line intervals where the maximum of $\min _{j}d_E(z,z_j)$ is attained at two or more points $z_i$'s.   

\subsection{Stochastic root finding} Because of various reasons, it is necessary and interesting to consider root finding problem for not only the determnistic case $f(z)=0$, but also the stochastic version $g(z,\xi )$, where $\xi $ is a random variable. 

We assume that the expectation of $g$, with respect to the variable $\xi$, is $f$, that is $E(g(z,.))=f(z)$. 

In stochastic root finding, we aim to test the robustness of a root finding method $IM$ in a stochastic environment. When given a function $h$ in the variable $z$, the root finding method $IM$ will construct a sequence defined by the formula $z_{n+1}=IM(h(z),z_n)$.  More precisely, we will follow the following procedure: 

- Generate a sequence $\{\xi _n\}$ from the distribution for the random variable $\xi$. 

- Choose an initial point $z_0$ and construct a sequence: 

$z_{n+1}=IM(g(z,\xi _n),z_n)$. 

We hope that $\{z_n\}$ will converge to a root of $f(z)$. To ensure the stability, we need the variation of $g$ is small: $Var(g(z,.))$ is small.   

\section{Experimental results} Here we discuss the settings for our experiments, and then present the experiments.

\subsection{Settings} We present the settings of our experiments.

1) The functions: We consider many functions the geometric configurations of whose roots vary (for example, some functions have roots which are vertices of convex polygons, while some functions have some roots inside the convex hull of other roots, some functions have three roots on the same line, and some functions have simple roots while other have non-simple roots). The list of the functions is as follows: 

\begin{eqnarray*}
f_1(z)&=&z(z-i)(z-3-2i),\\
f_2(z)&=&z(z-i)(z-3i),\\
f_3(z)&=&z(z-i)^2,\\
f_4(z)&=&z(z-i)(z-3-2i)(z-1-4i),\\
f_5(z)&=&z(z-i)(z-3-2i)(z-2-4i),\\
f_6(z)&=&z(z-i)^3,\\
f_7(z)&=&z(z-i)(z-1-i)(z-3-2i),\\
f_8(z)&=&z(z-i)(z-2i)(z-3-2i),\\
f_9(z)&=&z^2(z-i)^2,\\
f_{10}(z)&=&z^2(z-i)(z-1-i),\\
f_{11}(z)&=&z^2(z-i)(z-2i),\\
f_{12}(z)&=&z(z-i)(z-5i)(z-3-2i),\\
f_{13}(z)&=&z^2(z-i)(z-5i),\\
f_{14}(z)&=&z(z-2i)(z-3+3i))(z-3-6i)(z-5-2i),\\
f_{15}(z)&=&z(z-2i)(z-3-6i)(z-5-2i)(z-7+i),\\
f_{16}(z)&=&z(z-3-6i)(z-5-2i)(z-7+i)(z-4-3.4i),\\
f_{17}(z)&=&z(z-2i)(z-5-2i)(z-3+3i)(z-2-i),\\
f_{18}(z)&=&z(z-3-6i)(z-5-2i)(z-7+i)(z-2-i),
\end{eqnarray*}

\begin{eqnarray*}
f_{19}(z)&=&z^2(z-5-2i)(z-3+3i)(z-7+i),\\
f_{20}(z)&=&z^2(z-2-i)(z-5-2i)(z-3+3i),\\
f_{21}(z)&=&z^2(z-2-i)(z-5-2i)(z-3-6i),\\
f_{22}(z)&=&z(z-2-i)^2(z-3+3i)(z-3-6i),\\
f_{23}(z)&=&z^2+cos(z)+2sin(z)-1-0.5i,\\
f_{24}(z)&=&f_7(z)e^z,\\
f_{25}(z)&=&f_{17}(z)e^z.
\end{eqnarray*}

The function $f_{23}(z)$ has infinitely many roots. In the region $-10\leq x,y\leq 10$, it has 8 (approximate) roots: $0.01453348 +0.24577632i$, $-1.79690338 -0.16311646i$, $2.65293461 -2.52795741i$, $2.70778504+ 2.4386467i$, $-7.27782023 -4.1230358i$, $-7.26685729  +4.13462414i$, $9.62682067 -4.62305718i$,  and $9.63392763+ 4.61683271i$. 

In some experiments, we also consider functions of the form $g=f/f'$, where $f$ is one of the above functions. 

In stochastic root finding experiments, we consider functions of the form $g(z,\xi )=f(z)+\epsilon \xi (z^3+2z-5)$. Here $f(z)$ is one of the above functions, $\xi$ is a random variable with normal distribution $(E(\xi )=0, ~Var (\xi)=1)$, $\epsilon >0$ is a small number.  We can check that $E(g(z,.))=f(z)$. However, in experiments involving Newton's method vOptimization and BNQN we need to work with the functions $F(z)=|f(z)|^2/2$ and $G(z,\xi )=|g(z,\xi )|^2/2$. One can check that $E(G(z,.))$ is not $F(z)$, but rather $E(G(z,.))=F(z)+\epsilon ^2 |z^3+2z-5|^2$. Therefore, we should choose $\epsilon$ small so that $E(G(z,.))$  is close enough to $F(z)$.  

2) The experiments: 

We will compute basins of attraction when using different methods (iterative or flow) to certain complex-valued functions in 1 complex variable (listed in the previously). We also draw pictures of (reduced) Voronoi's diagrams, where the multiplicity of a point is disregarded.  

3) Colours and interpretations of the results: If the factorisation of a function $f$ is written as (note that the order $z_1$, $z_2$, $\ldots$ is important) $f(z)=(z-z_1)^{n_1}(z-z_2)^{n_2}\ldots $, then we will colour the basins of attraction as follows: that for the first root has the green colour, that for the second root has the yellow colour, that for the third root has the blue colour, that for the fourth root has the red colour, that for the fifth root has the pink colour, that for the sixth root has the cyan colour, that for the seventh root has the orange colour, that for the eighth root has the purple colour. Other points (e.g. points whose corresponding sequence goes to infinity or does not converge to one of the roots) will have the black colour.   

For example, for the function $f_{10}(z)=z^2(z-i)(z-1-i)$, by our convenience, the first root is $0$, the second root is $i$ and the third root is $1+i$. Therefore, the basin of attraction for $0$ will have the green colour, the basin of attraction for $i$ will have the yellow colour, and the basin of attraction for $1+i$ has the blue colour.

In Voronoi's diagrams, we also follow the same convenience. 

4) Implementation of the codes: 

Implementation of BNQN and Newton's flows are as discussed in Section 2. Implementation of Newton's method and Random Relaxed Newton's method is straightforward from their descriptions. 

For BNQN, most of the case we will run BNQN New Variant with $\theta =0$, and this will be reported in the experimental results as BNQN. In some cases, we also run with BNQN New Variant with $\theta =1$, and report it as BNQN v2. 

To avoid overflow errors, we use the mpmath library \cite{Refmp}. 

5) Experimental procedures: 

We will choose a grid of $240\times 240$ points in $[-10,10]\times [-10,10]$, whose center is randomly chosen in a small neighbourhood of $(0,0)$. Each of the points $z_0$ in the grid will be the initial point to run the concerned iterative method. We will run a maximum of $10 000$ iterations, and can stop earlier if either the function value is smaller than a threshold $\epsilon$ or the gradient of the objective cost function is smaller than the threshold (for BNQN and Newton's method vOptimization). If $z_n$ is the last constructed point, we declare that $z_n$ belongs to the basin of attraction for a root $z^*$ if the Euclidean distance $d(z_n,z^*)$ is smaller than a threshold. We choose the threshold to be $10^{-6}$ or smaller. 

This manner of determining the basins of attraction is reasonable. First, we know that if $z_n$ is close enough to a root $z^*$, then the considered iterative methods will converge to the root $z^*$. Second, in the stochastic root finding setting, as explained in Section 3, $E(G(z,.))$ is not the same as $F(z)$, and hence we do not expect that $z_n$ converges precisely to a root of $f(z)$ (or more generally, a minimum of $F(z)$). However, if $E(G(z,.))-F(z)$ is small, then we expect that the minizers of $E(G(z,.))$ are within a small neighbourhood of that of $F(z)$. Since we expect that $z_n$ will converge to minimizers of $E(G(z,.))$, it follows that $z_n$ will converge to within a small neighbourhood of a minimizer of $F(z)$. Explicitly, in $g(z,\xi )=f(z)+\epsilon \xi (z^3+2z-5)$, we choose $\epsilon$ to be $10^{-4}$. In this case, we reduce the threshold to declare that the point $z_n$ belongs to the basin of attraction of a root $z^*$ to be $10 
\epsilon$, which is $10^{-3}$. 

\subsection{Results}

\subsubsection{Polynomials}

\begin{figure}
    \centering
    \includegraphics[width=5cm]{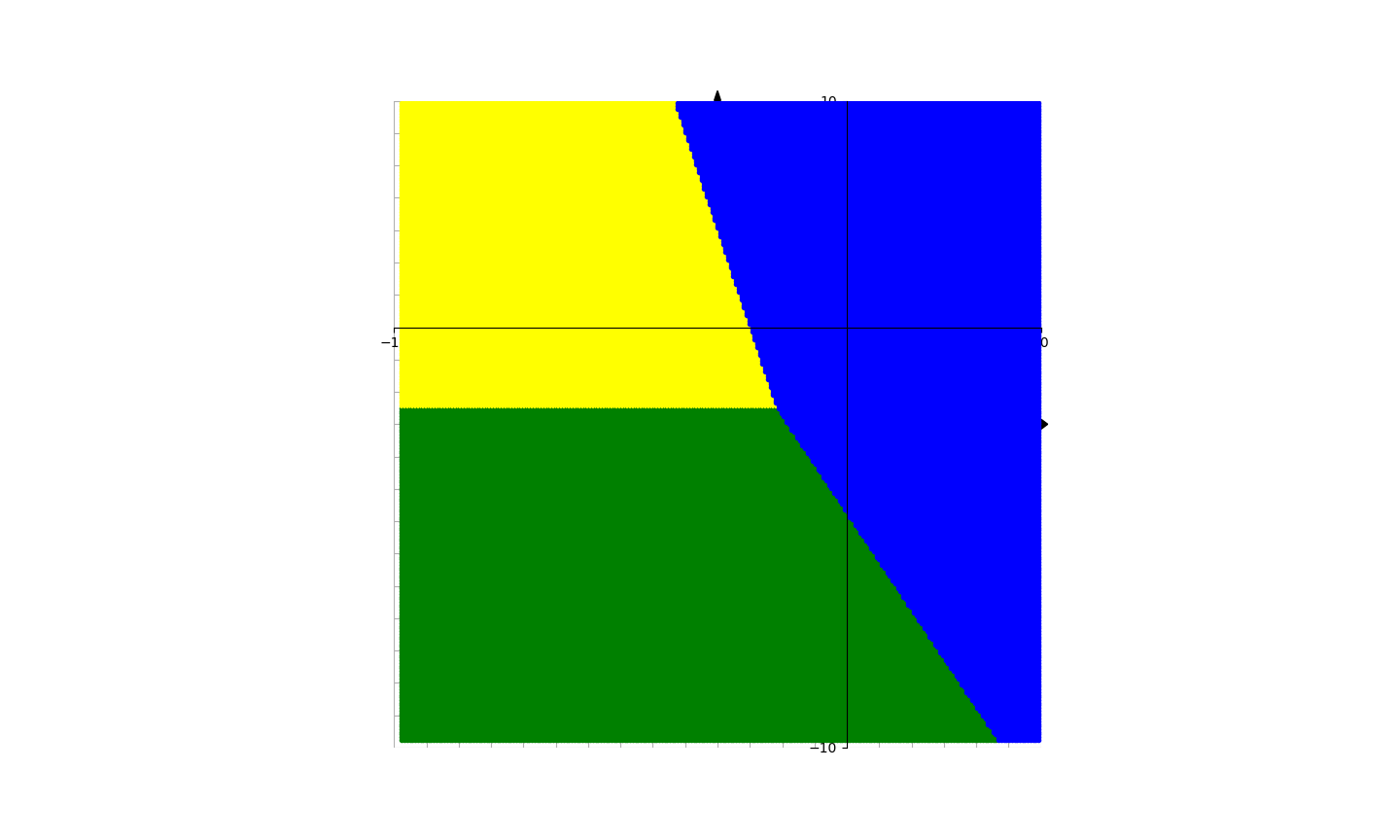}
    \includegraphics[width=3cm]{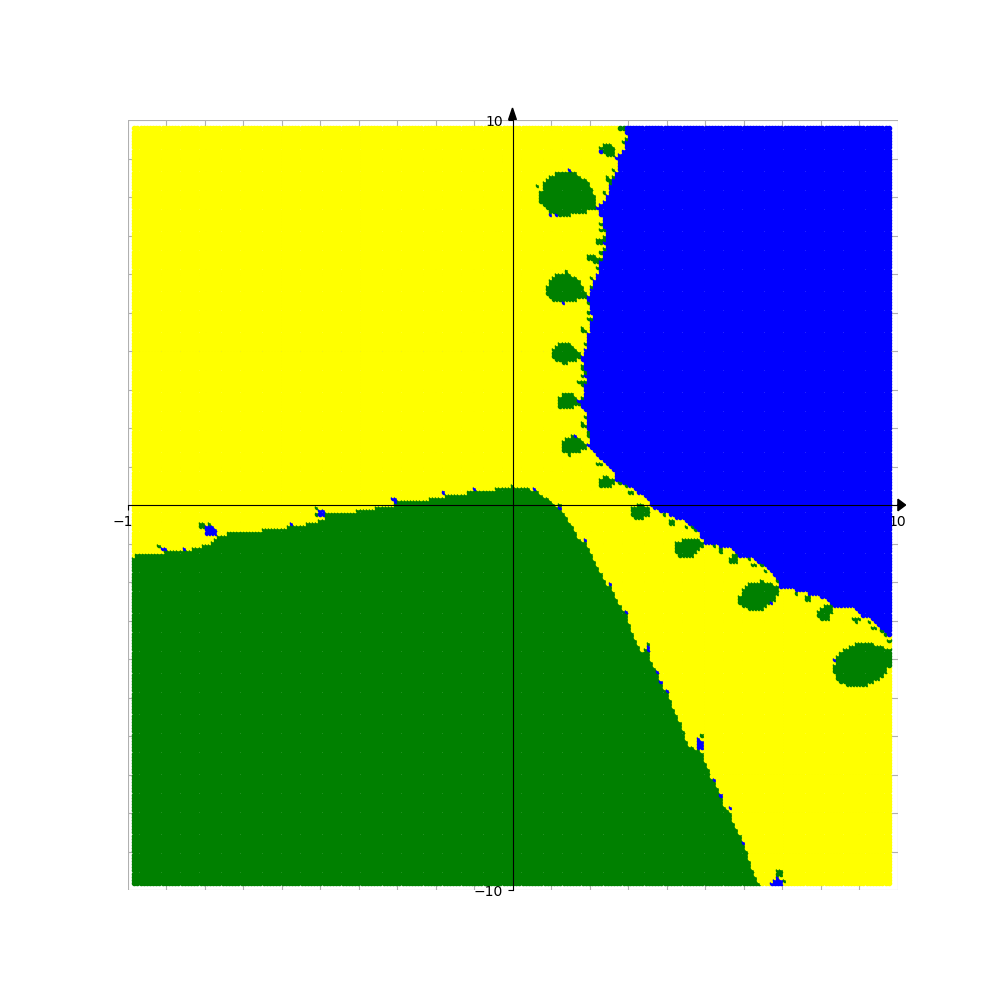}
    \includegraphics[width=3cm]{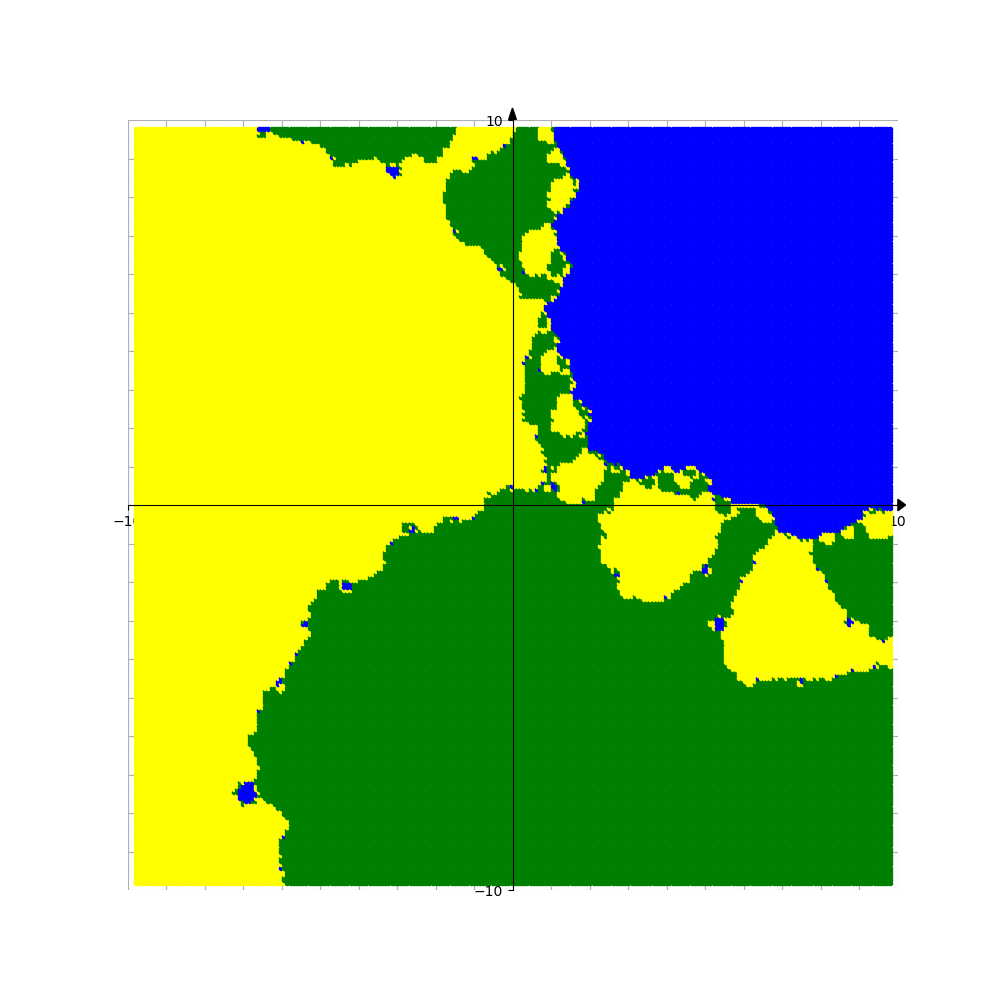}

    \bigskip
    \includegraphics[width=5.5cm]{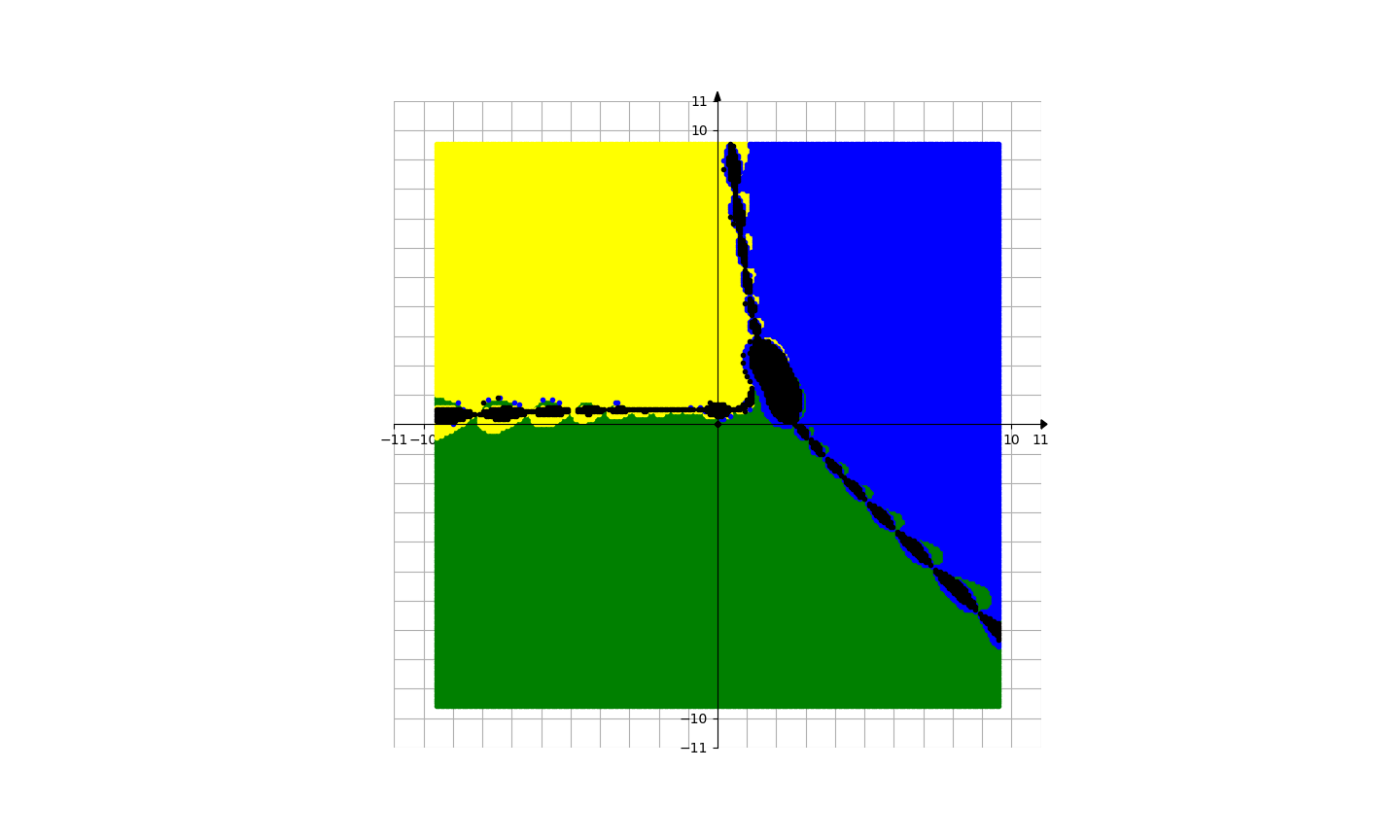}
    \includegraphics[width=3cm]{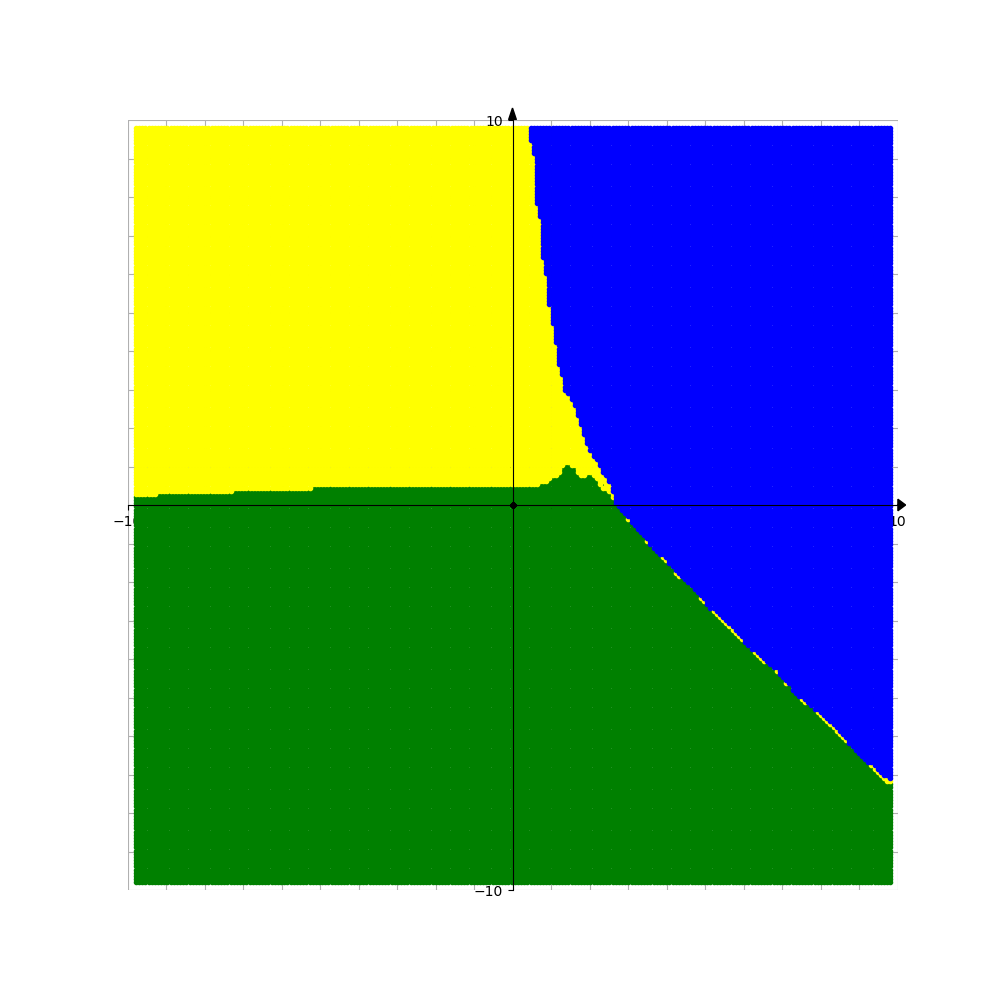}

    \bigskip
    \includegraphics[width=3cm]{NewtonFlowF1.png}
    \includegraphics[width=3cm]{NewtonFlowF1withFirstDerivative.png}
    \includegraphics[width=3cm]{NewtonFlowF1AndDerivative.png}
    
    \caption{Basins of attraction for finding roots of the function $f_1$ by different methods. Pictures are referenced to from top to bottom, from left to right. Row 1: left picture is Voronoi's diagram, central picture is for Newton's method, right picture is for Random Relaxed Newton's method. Row 2: left picture is for Newton's method vOptimization, right picture is for BNQN. Row 3: left picture is for Newton's flow, central picture is for Newton's flow vFraction, right picture is for Newton's flow vOptimization. The black points in some of these pictures are those in the basin of attraction of critical points of $f_1$.}
    \label{fig:f1}
\end{figure}

\begin{figure}
    \centering
    \includegraphics[width=5cm]{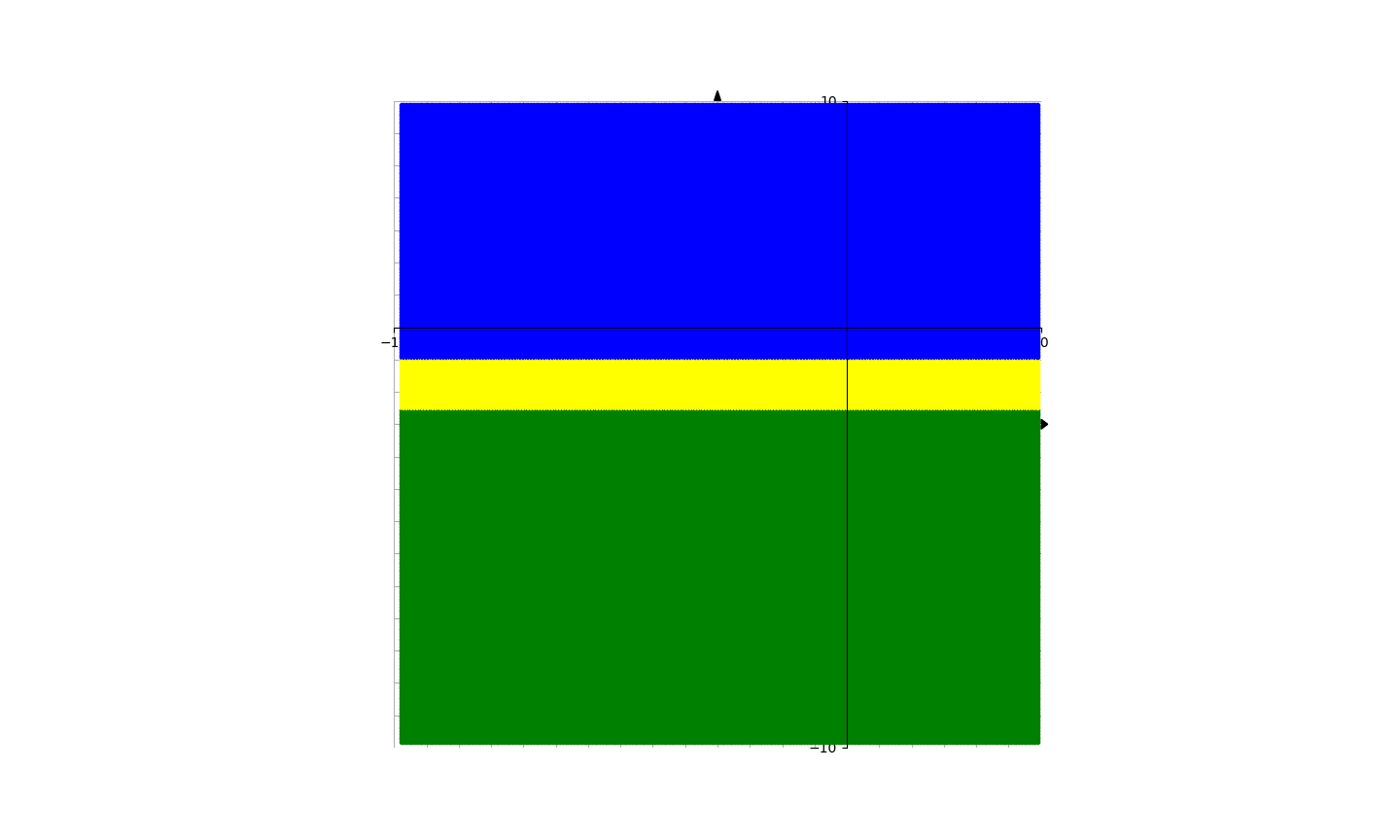}
    \includegraphics[width=3cm]{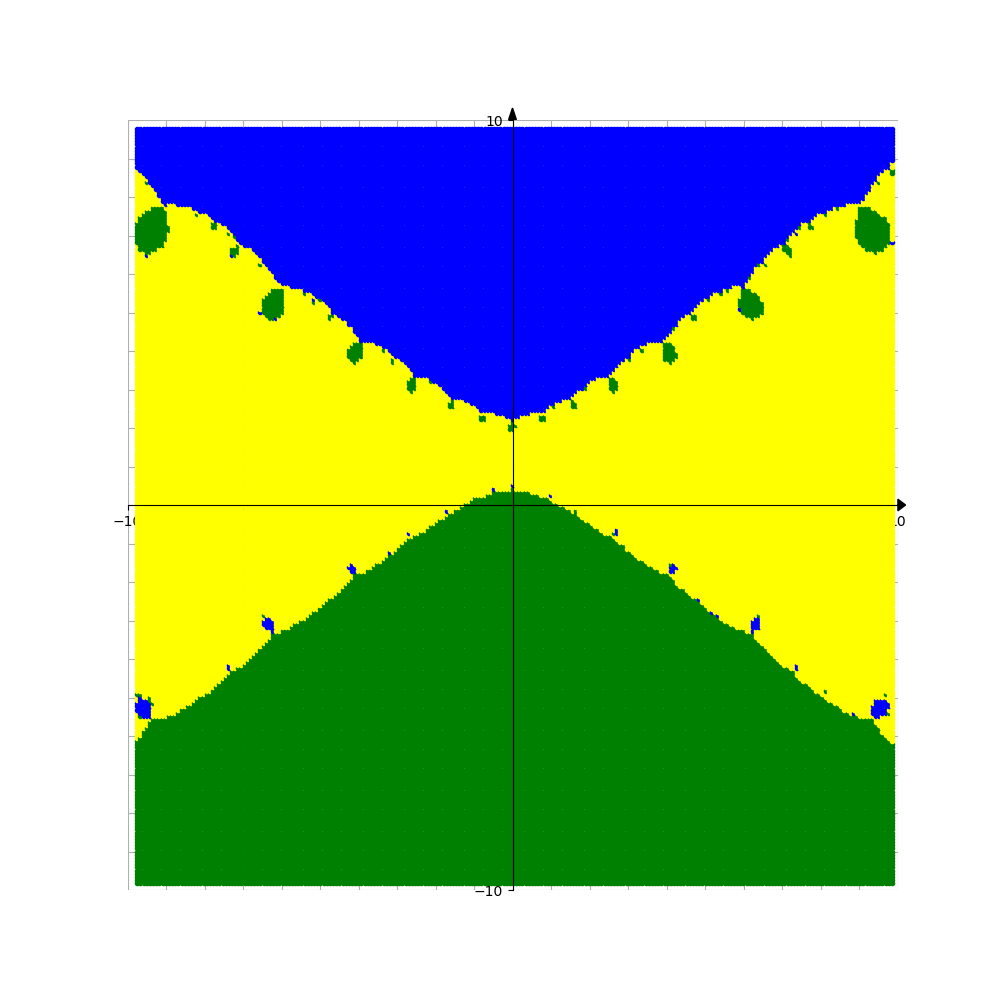}
    \includegraphics[width=3cm]{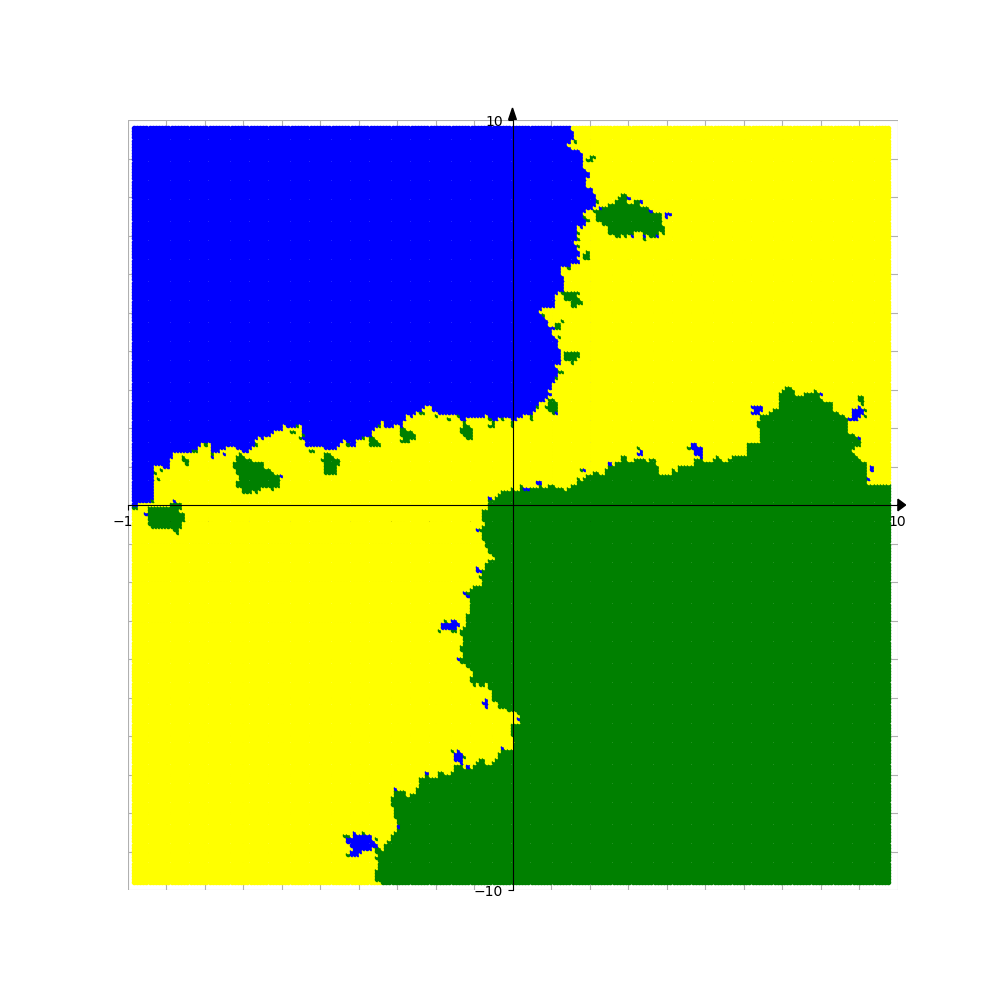}

    \bigskip
    \includegraphics[width=5.5cm]{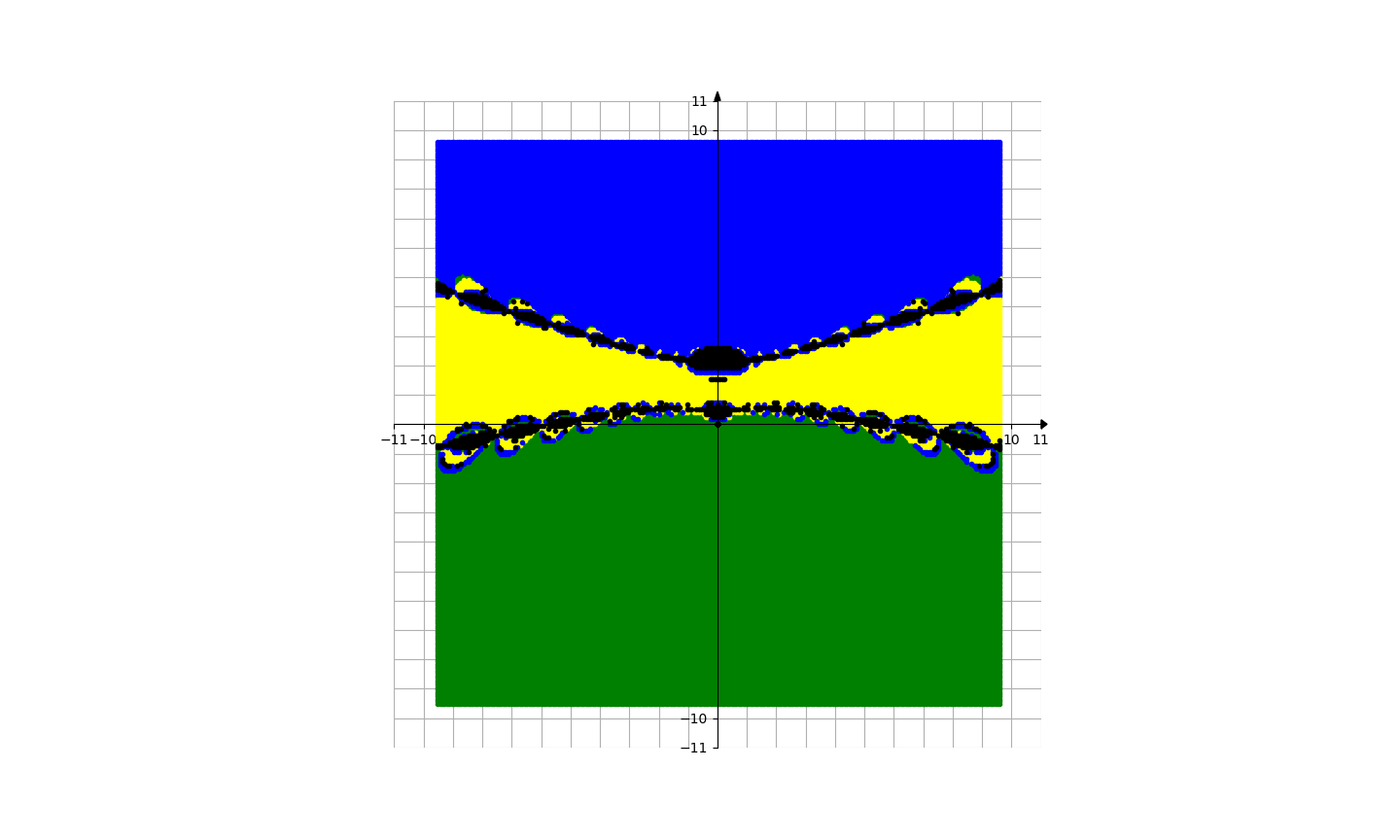}
    \includegraphics[width=3cm]{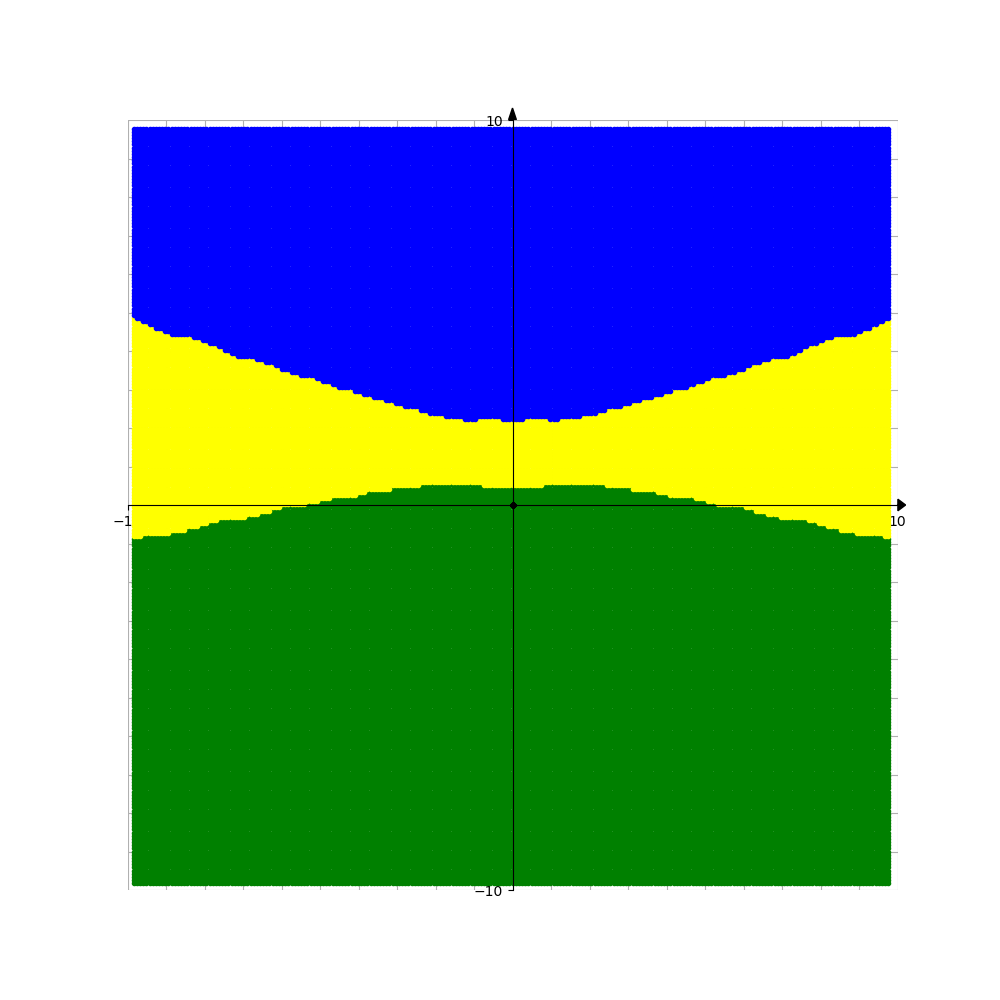}

    \bigskip
    \includegraphics[width=3cm]{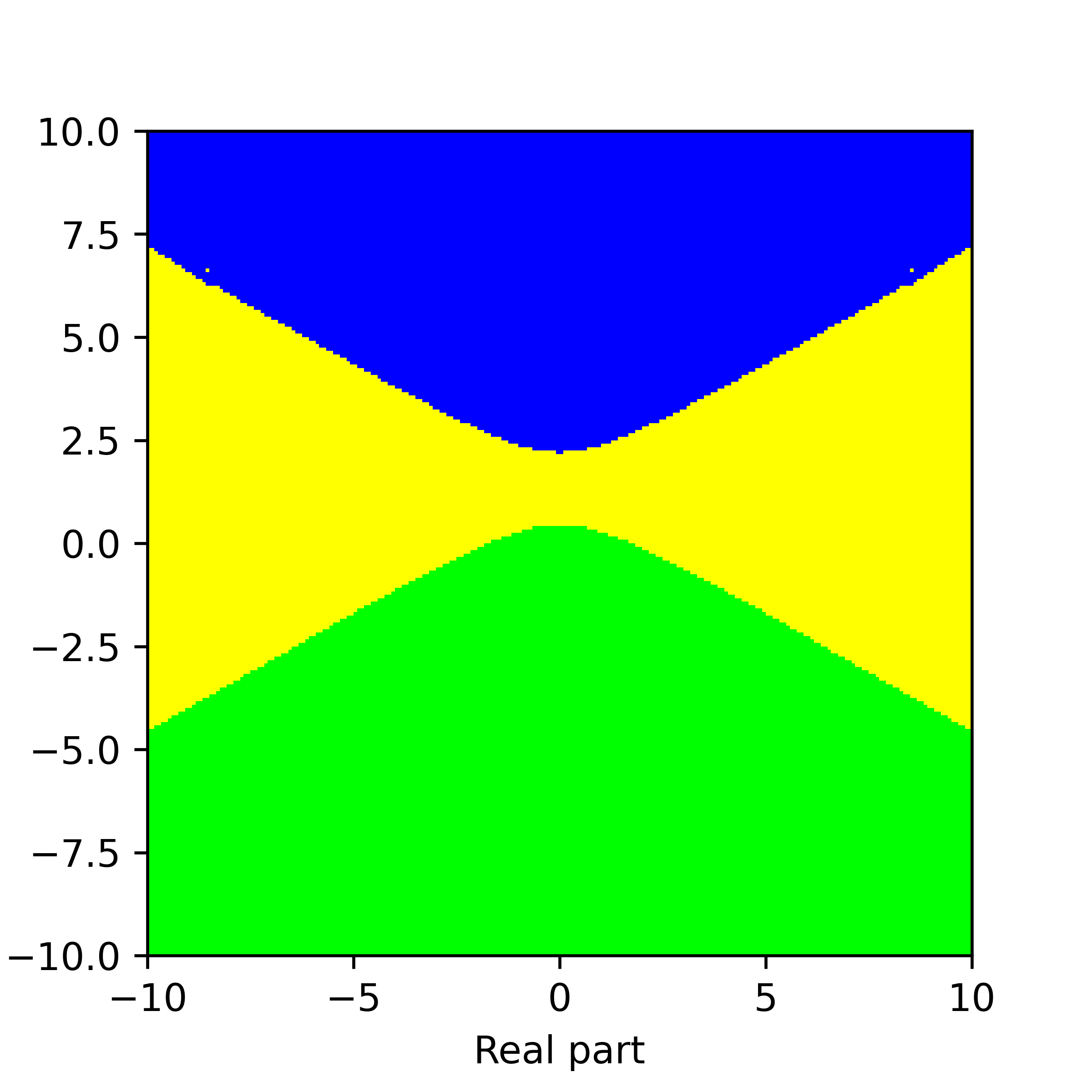}
    \includegraphics[width=3cm]{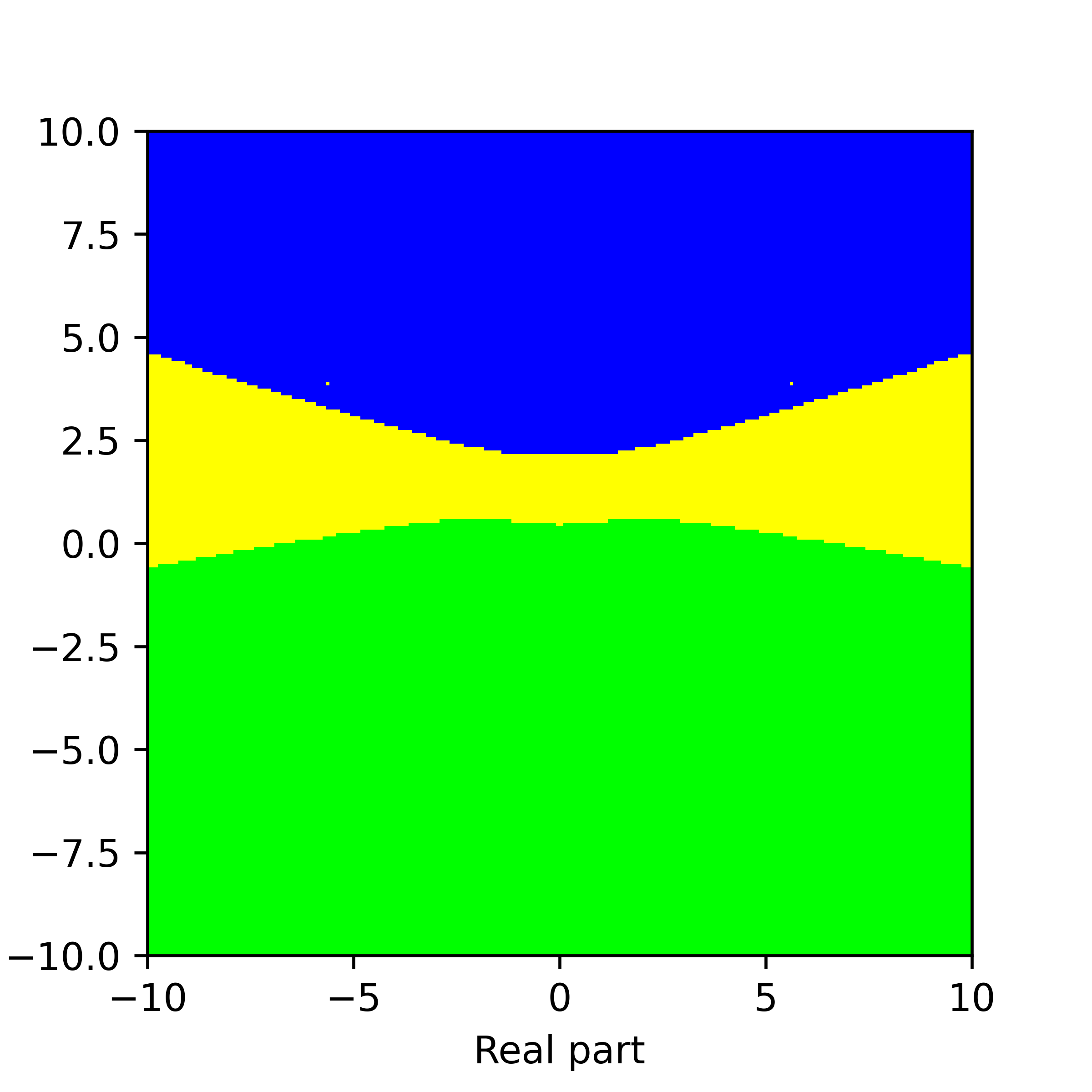}
    \includegraphics[width=3cm]{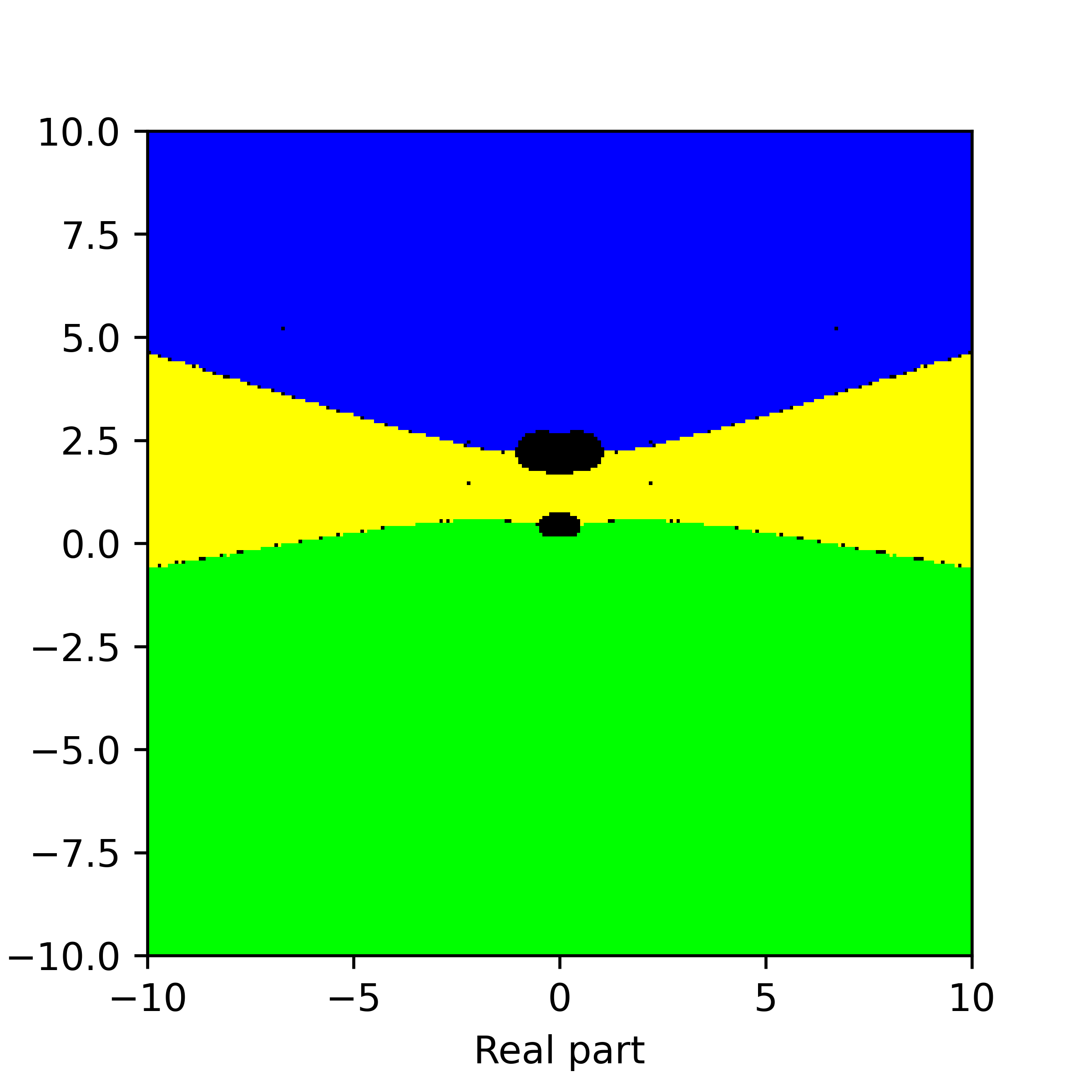}
    
    \caption{Basins of attraction for finding roots of the function $f_2$ by different methods. Pictures are referenced to from top to bottom, from left to right. Row 1: left picture is Voronoi's diagram, central picture is for Newton's method, right picture is for Random Relaxed Newton's method. Row 2: left picture is for Newton's method vOptimization, right picture is for BNQN. Row 3: left picture is for Newton's flow, central picture is for Newton's flow vFraction, right picture is for Newton's flow vOptimization. The black points in some of these pictures are those in the basin of attraction of critical points of $f_2$.}
    \label{fig:f2}
\end{figure}

\begin{figure}
    \centering
    \includegraphics[width=5cm]{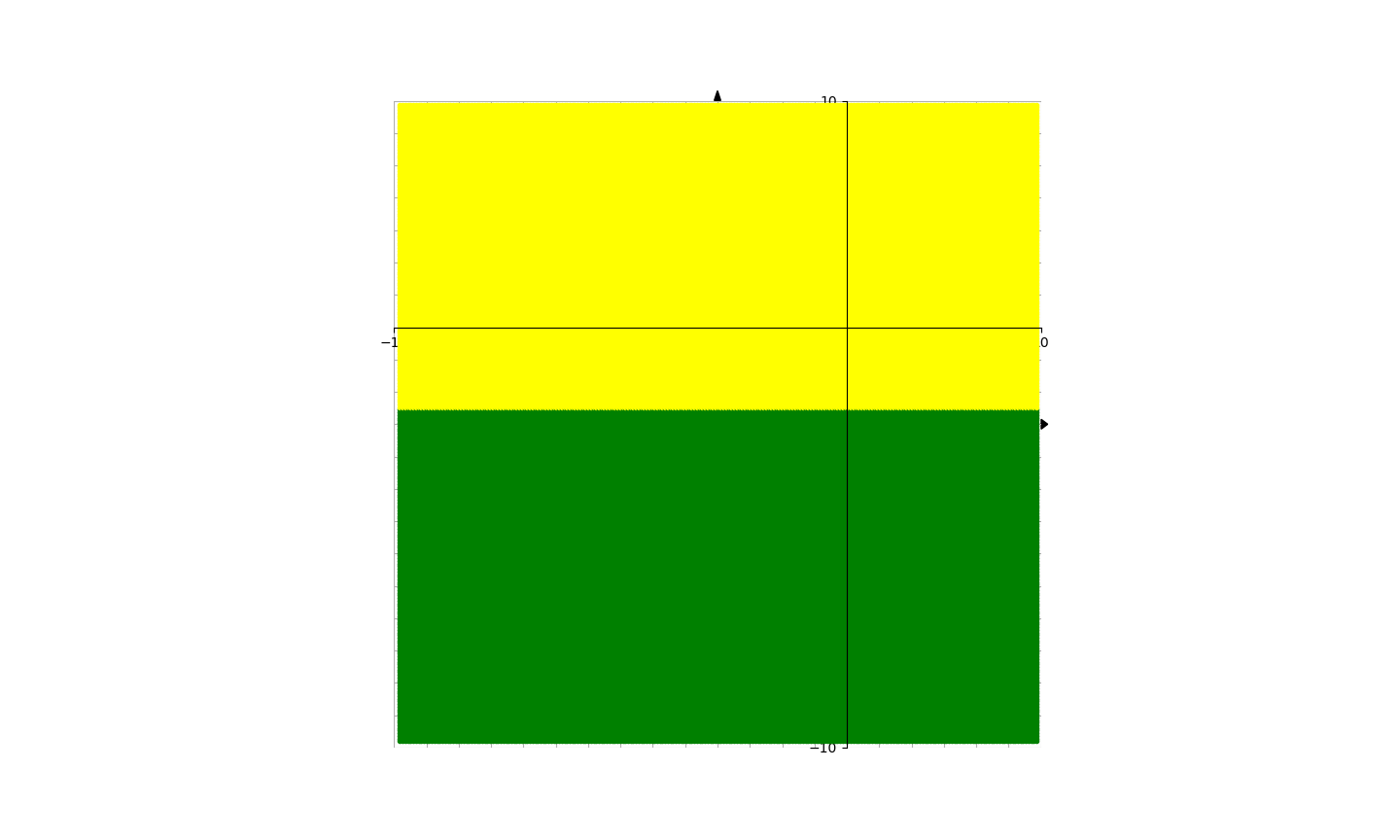}
    \includegraphics[width=3cm]{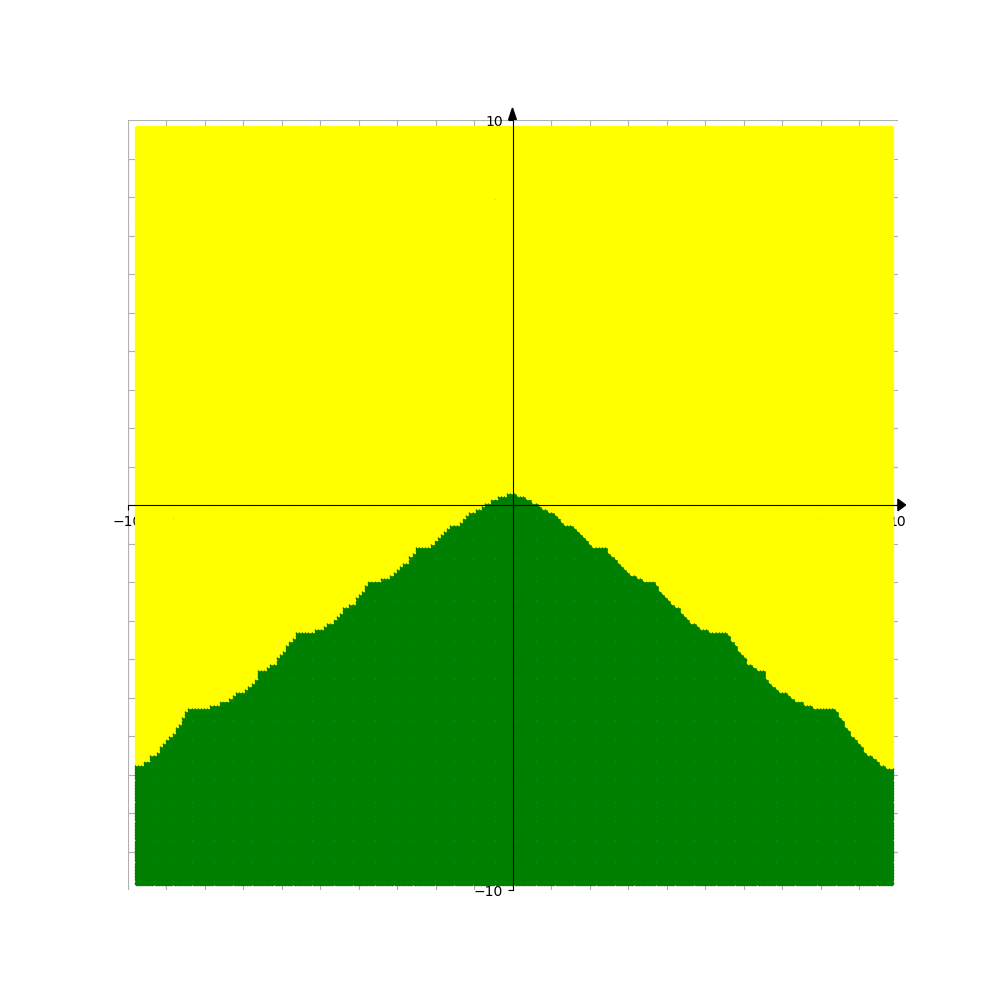}
    \includegraphics[width=3cm]{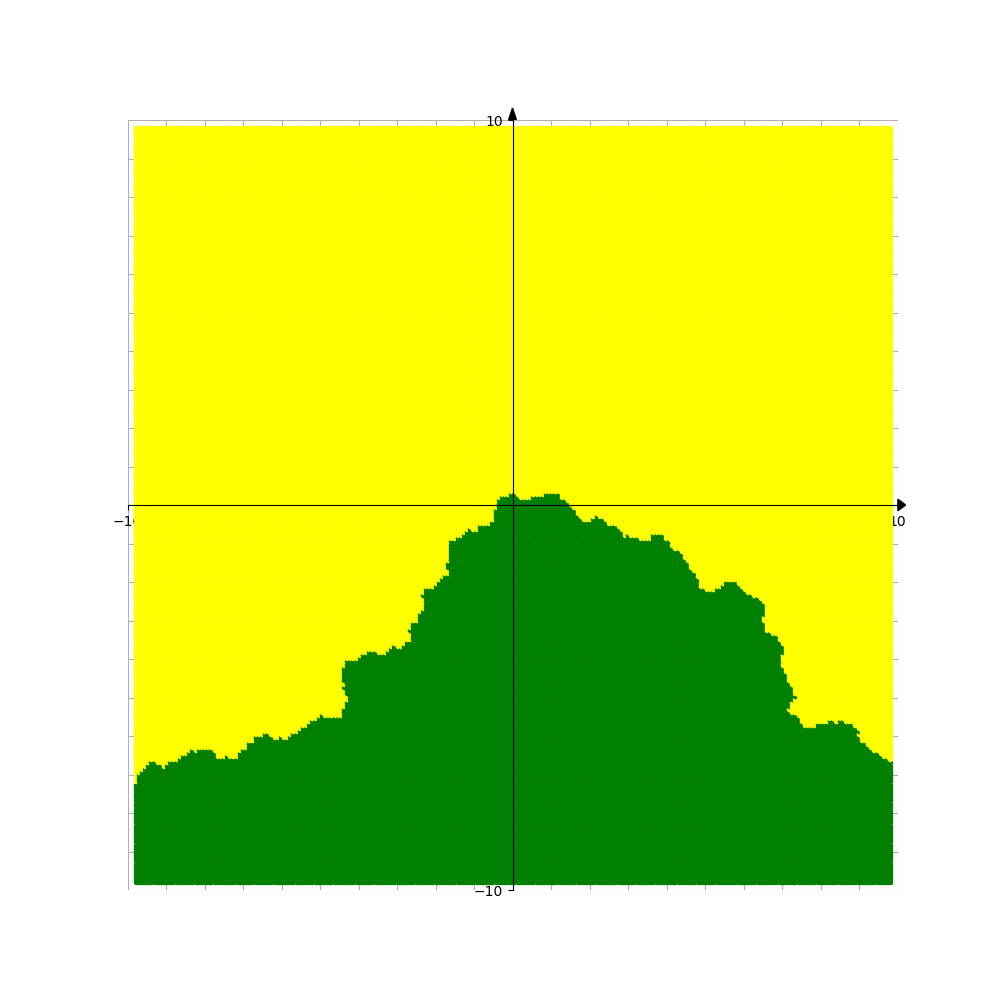}

    \bigskip
    \includegraphics[width=5.5cm]{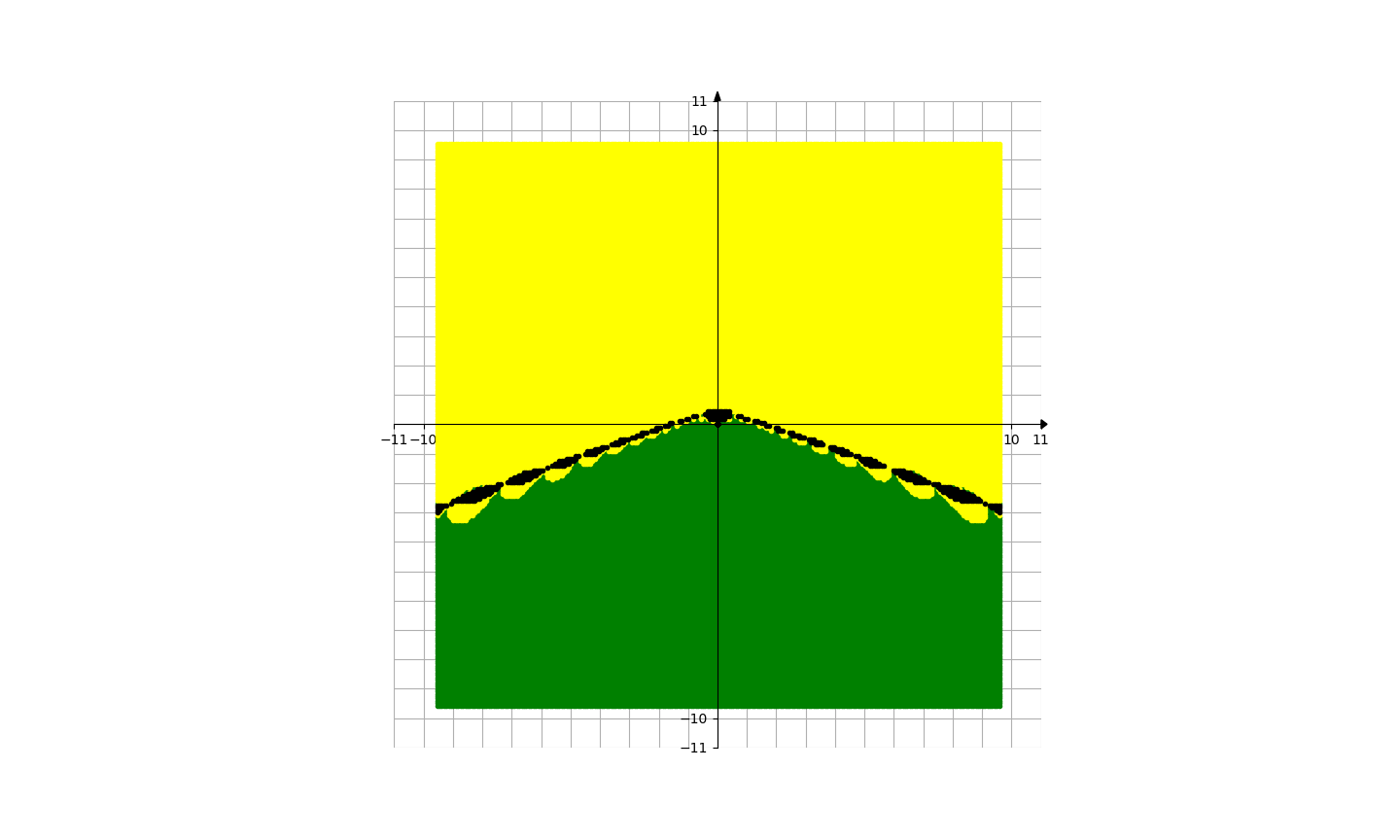}
    \includegraphics[width=3cm]{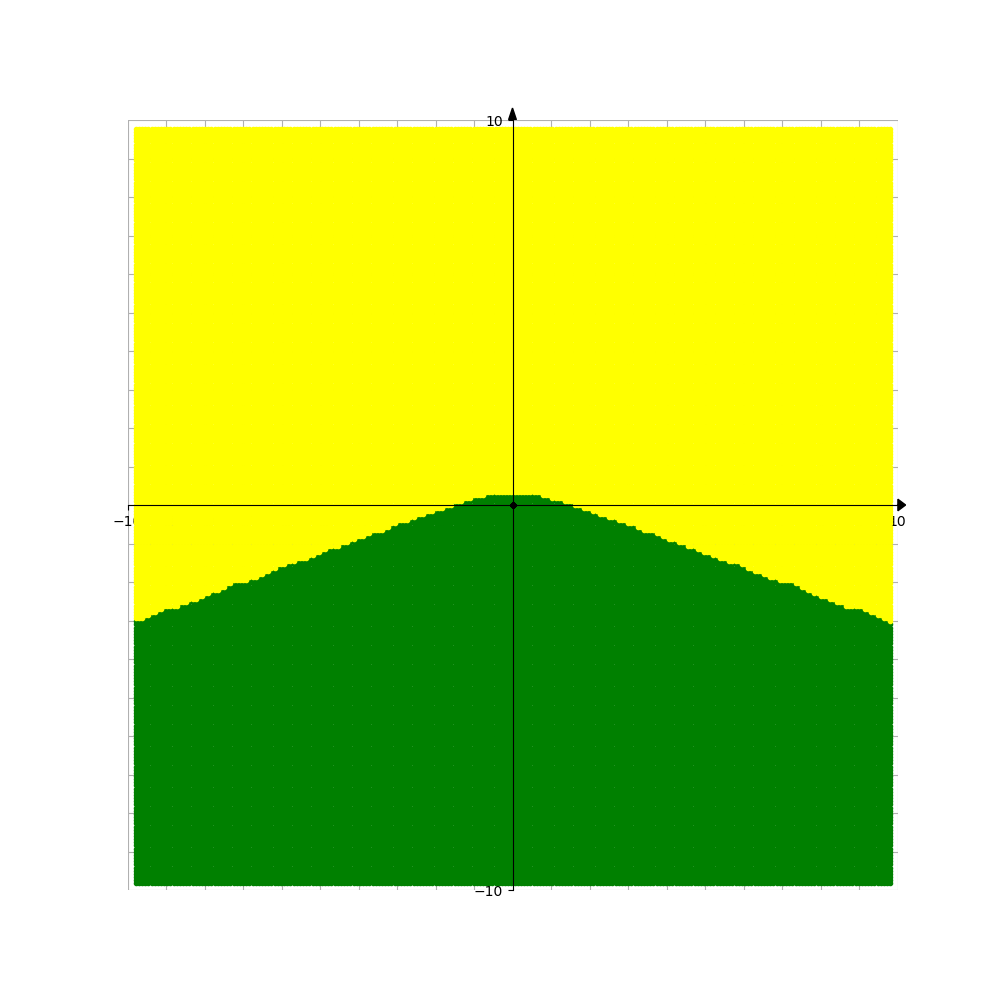}

    \bigskip
    \includegraphics[width=3cm]{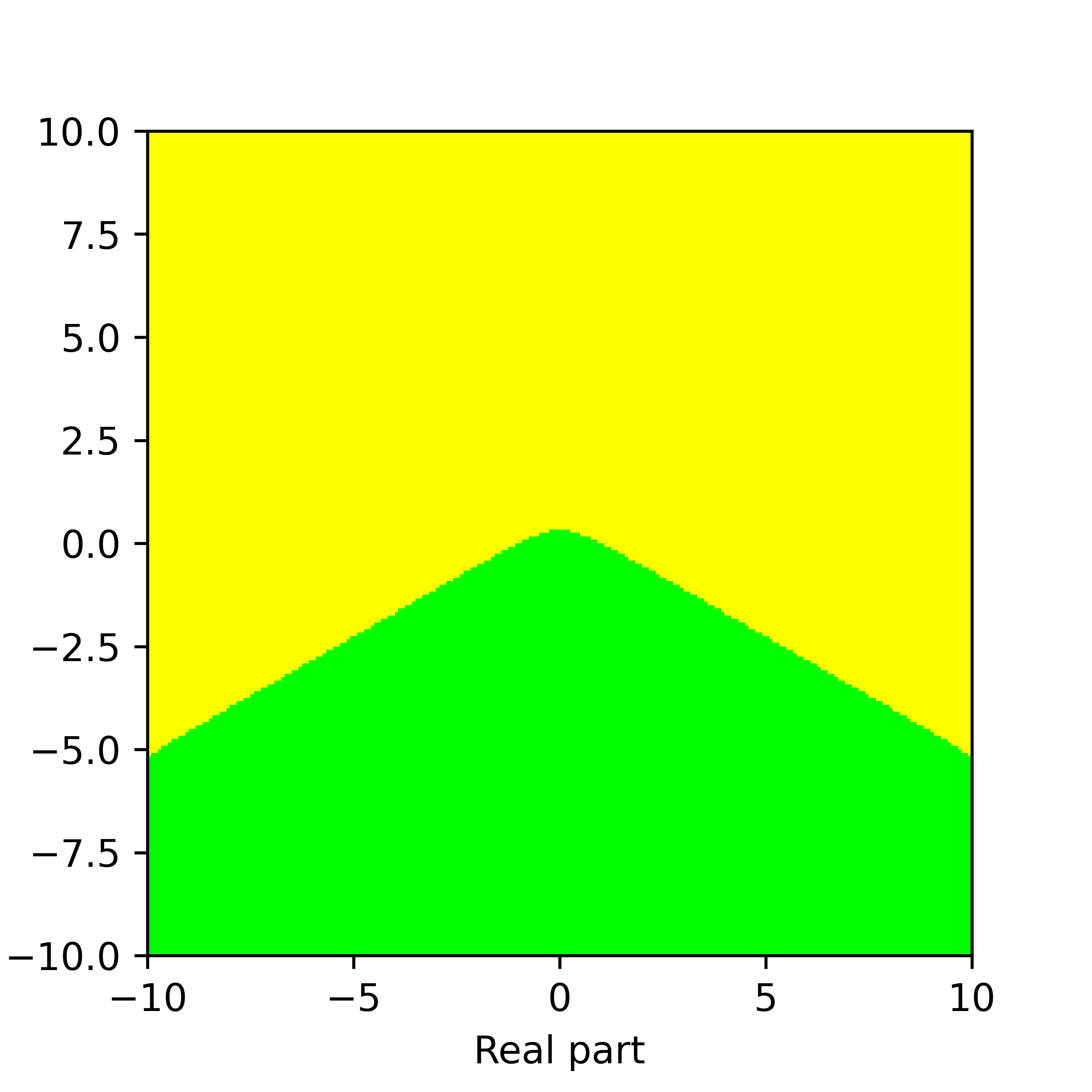}
    \includegraphics[width=3cm]{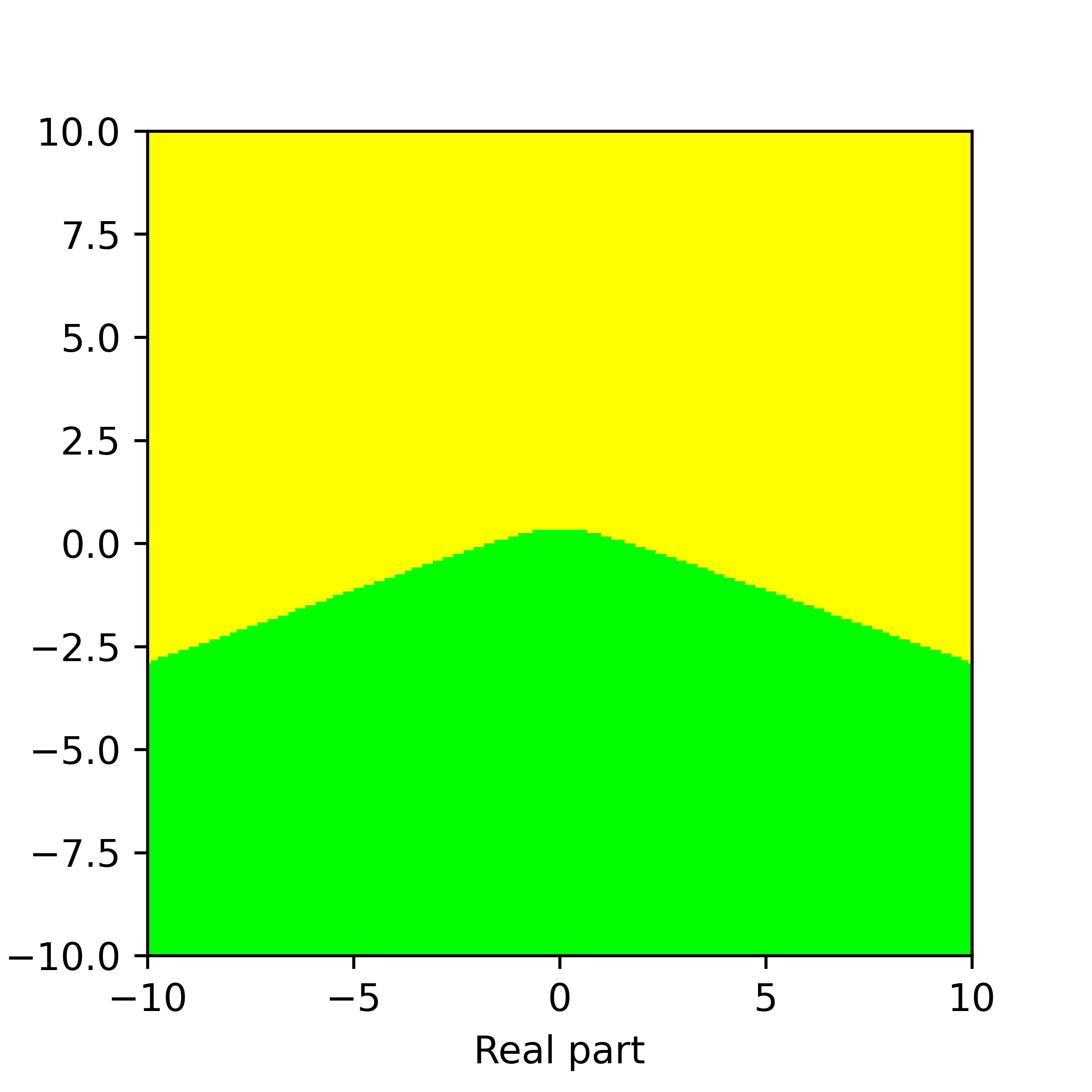}
    \includegraphics[width=3cm]{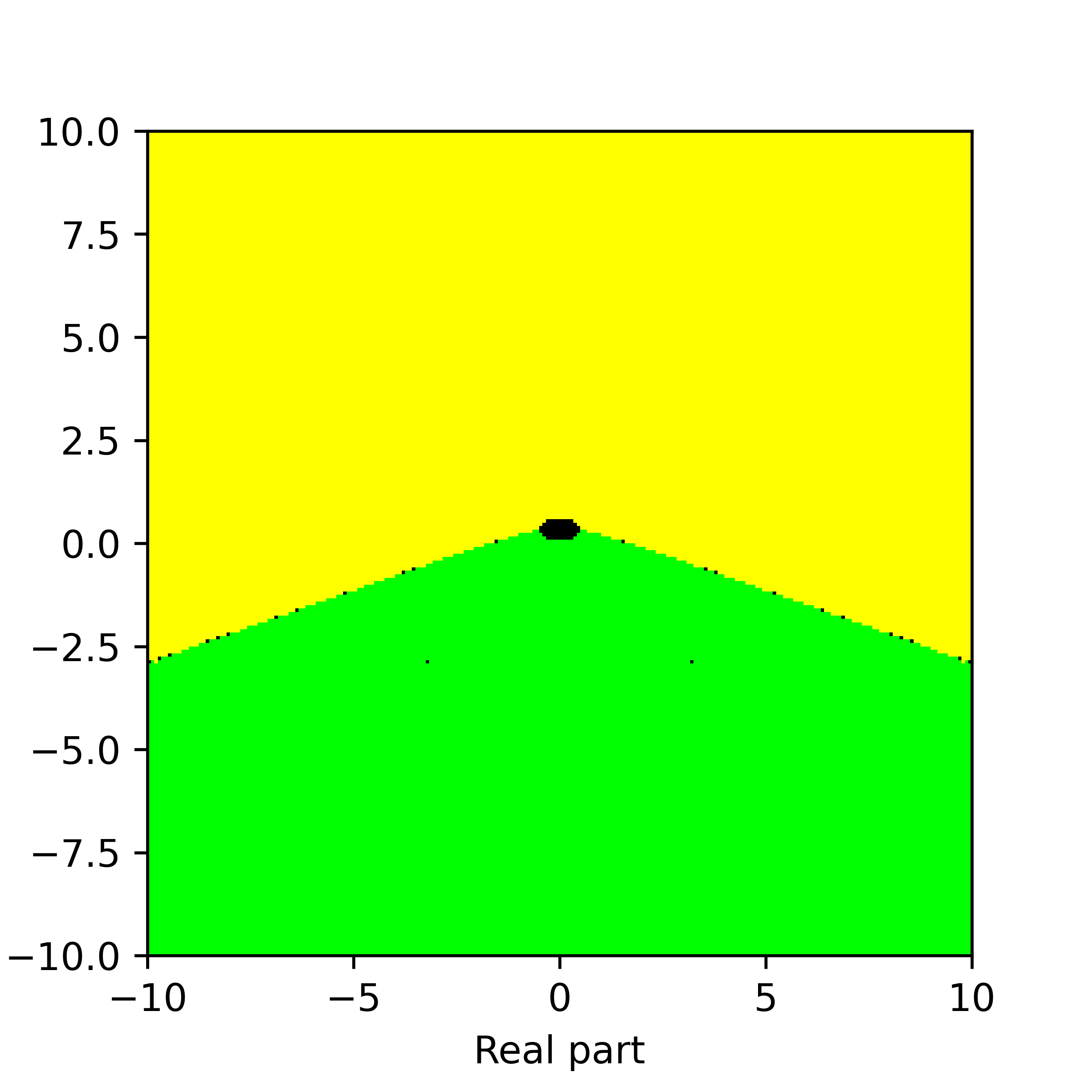}
    
    \caption{Basins of attraction for finding roots of the function $f_3$ by different methods. Pictures are referenced to from top to bottom, from left to right. Row 1: left picture is Voronoi's diagram, central picture is for Newton's method, right picture is for Random Relaxed Newton's method. Row 2: left picture is for Newton's method vOptimization, right picture is for BNQN. Row 3: left picture is for Newton's flow, central picture is for Newton's flow vFraction, right picture is for Newton's flow vOptimization. The black points in some of these pictures are those in the basin of attraction of critical points of $f_3$.}
    \label{fig:f3}
\end{figure}

\begin{figure}
    \centering
    \includegraphics[width=5cm]{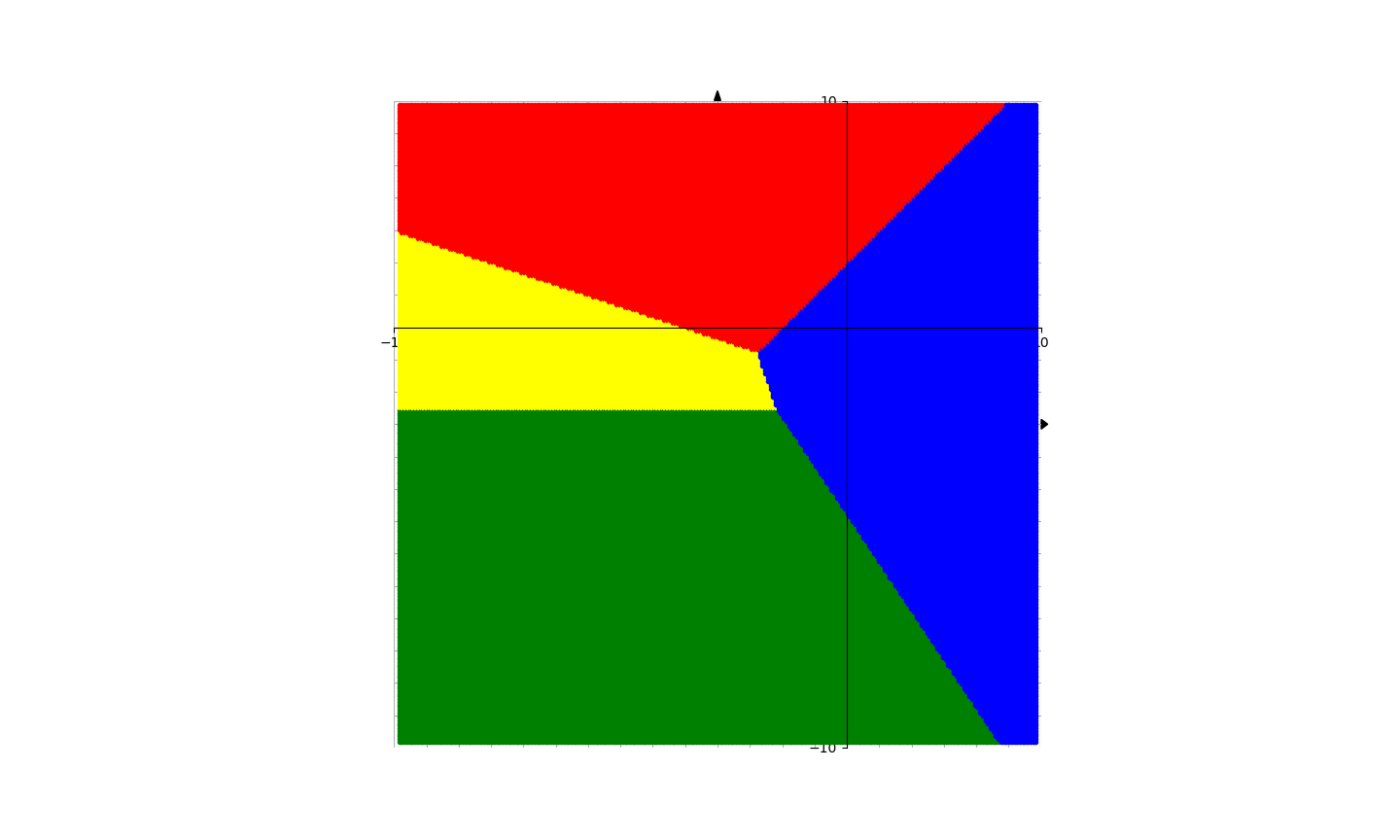}
    \includegraphics[width=3cm]{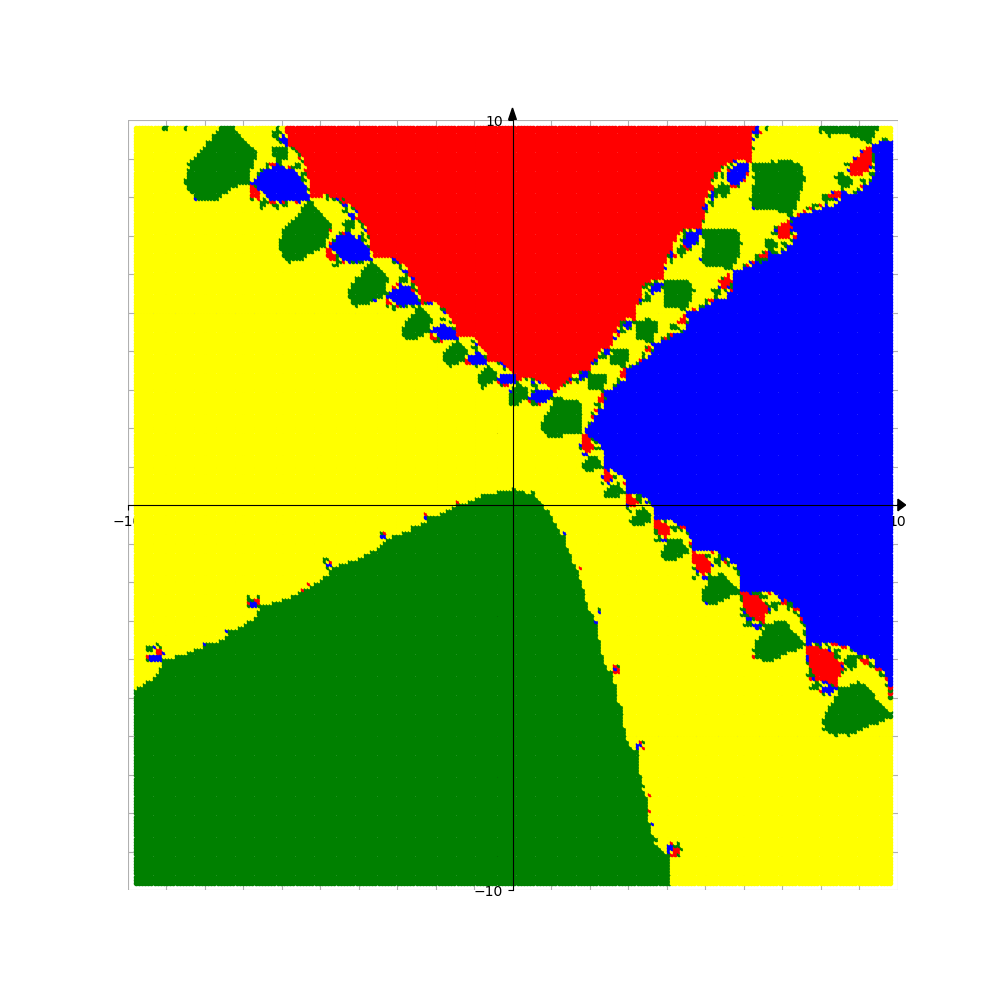}
    \includegraphics[width=3cm]{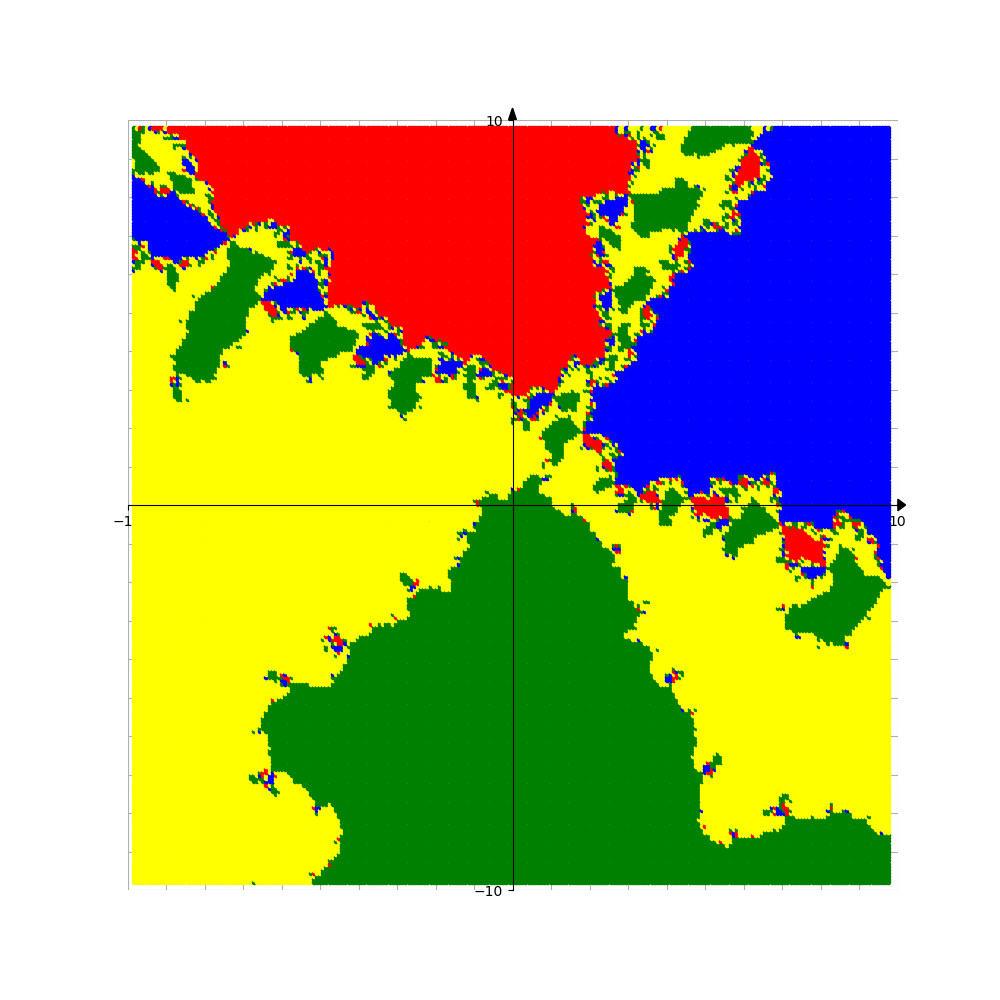}

    \bigskip
    \includegraphics[width=5.5cm]{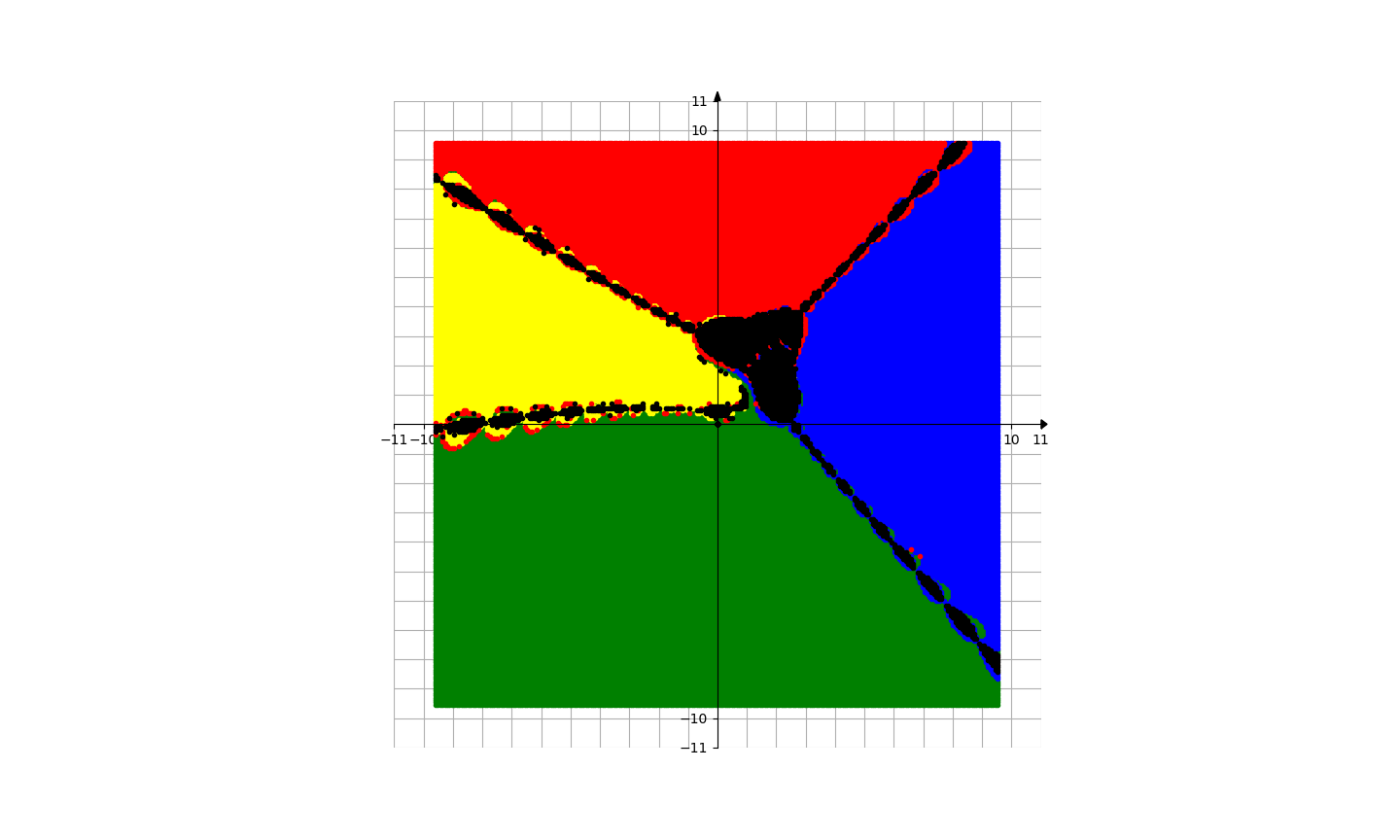}
    \includegraphics[width=3cm]{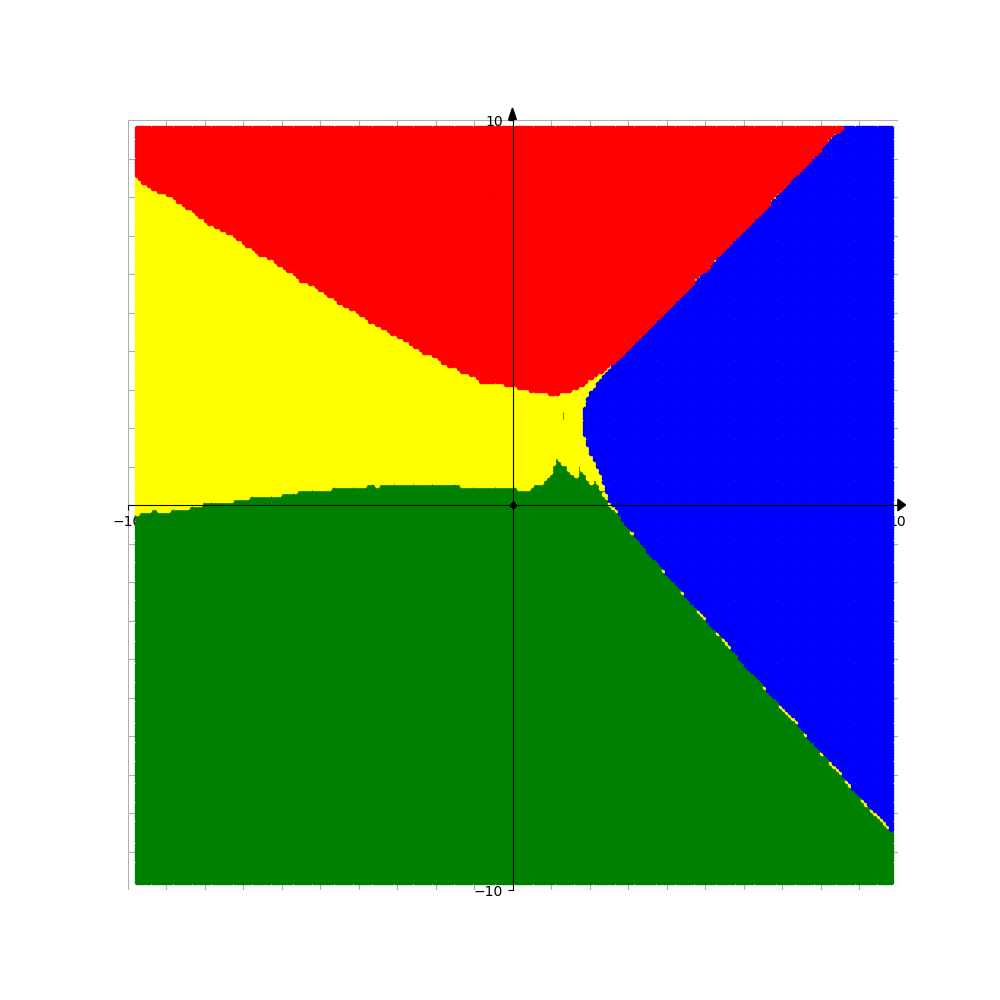}

    \bigskip
    \includegraphics[width=3cm]{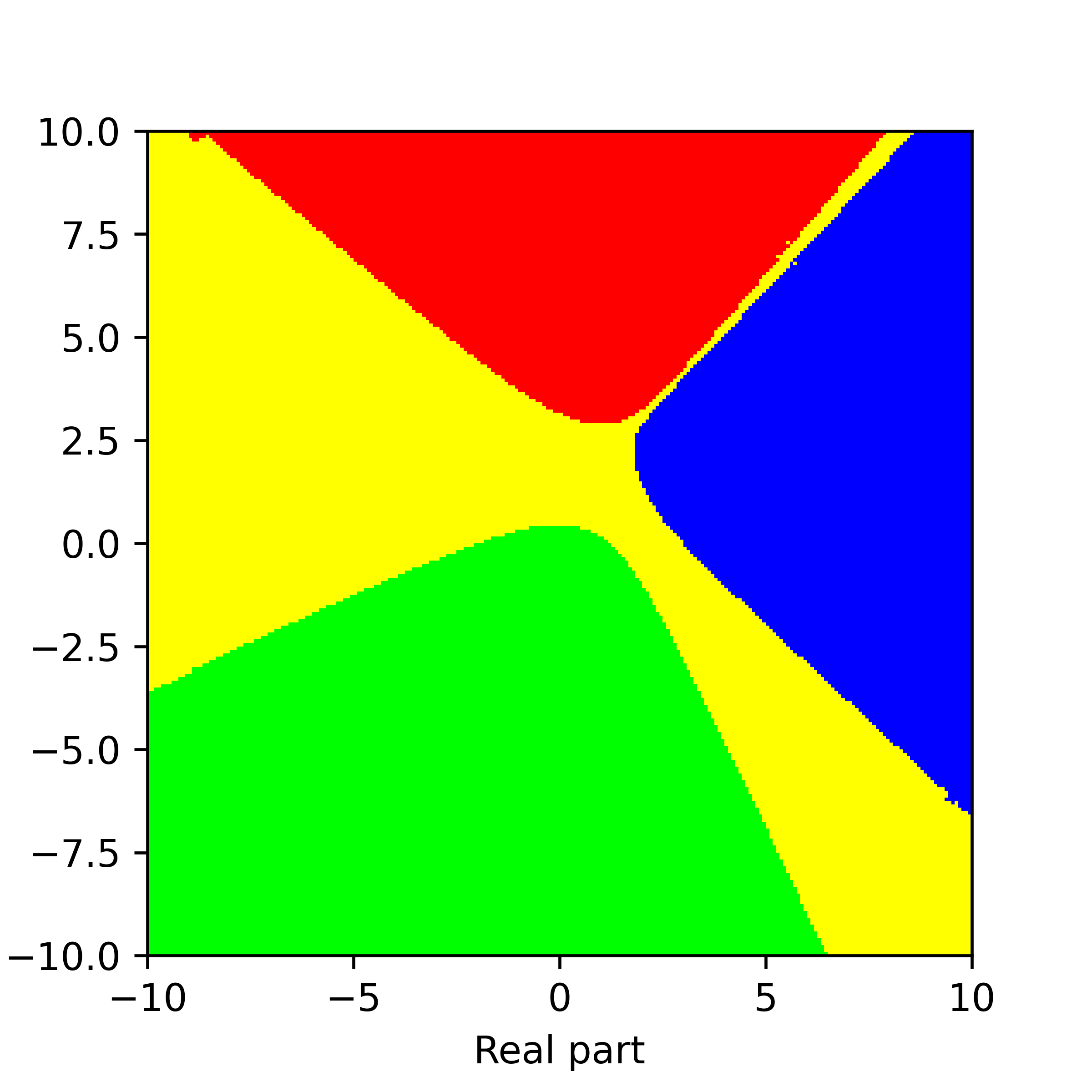}
    \includegraphics[width=3cm]{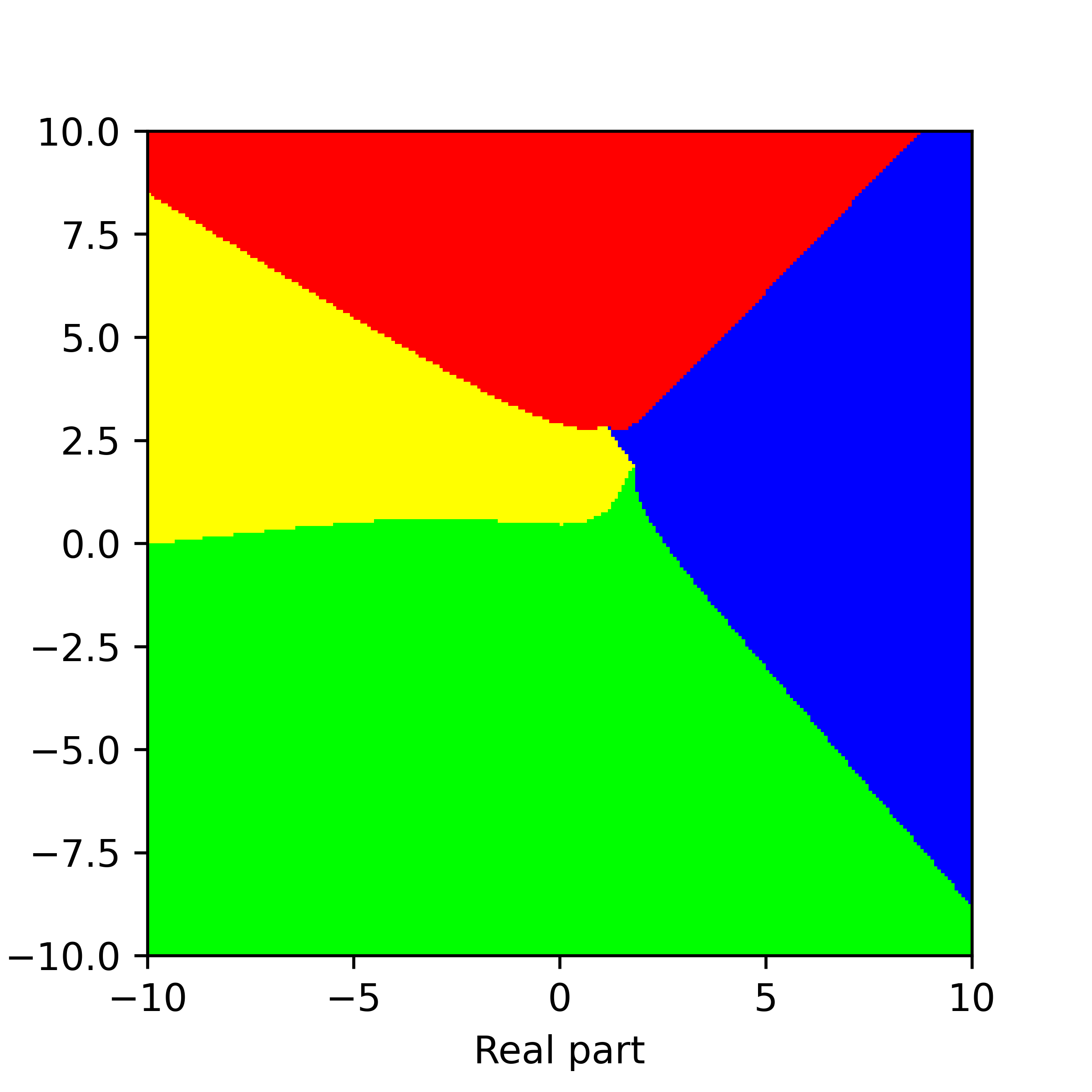}
    \includegraphics[width=3cm]{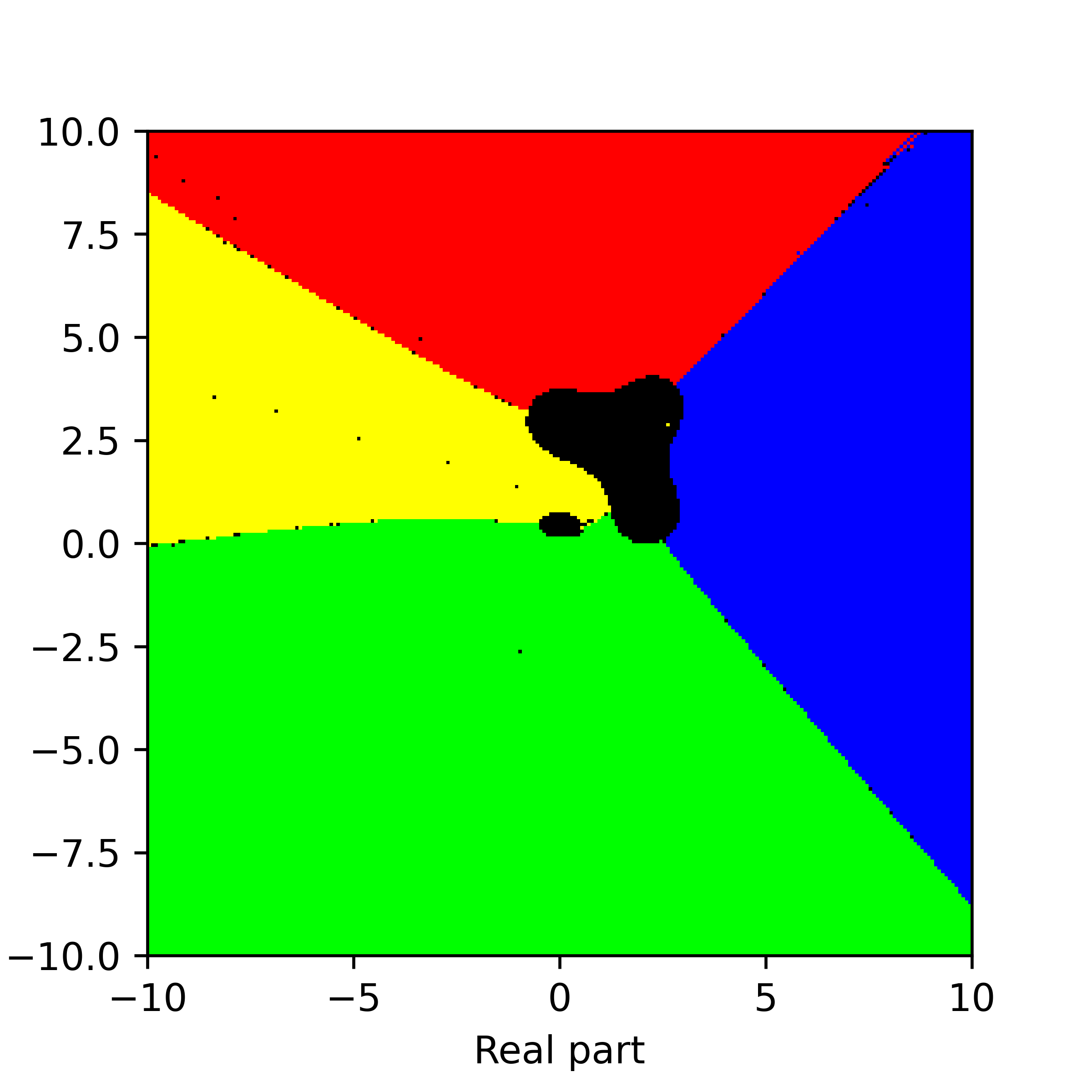}
    
    \caption{Basins of attraction for finding roots of the function $f_4$ by different methods. Pictures are referenced to from top to bottom, from left to right. Row 1: left picture is Voronoi's diagram, central picture is for Newton's method, right picture is for Random Relaxed Newton's method. Row 2: left picture is for Newton's method vOptimization, right picture is for BNQN. Row 3: left picture is for Newton's flow, central picture is for Newton's flow vFraction, right picture is for Newton's flow vOptimization. The black points in some of these pictures are those in the basin of attraction of critical points of $f_4$.}
    \label{fig:f4}
\end{figure}

\begin{figure}
    \centering
    \includegraphics[width=5cm]{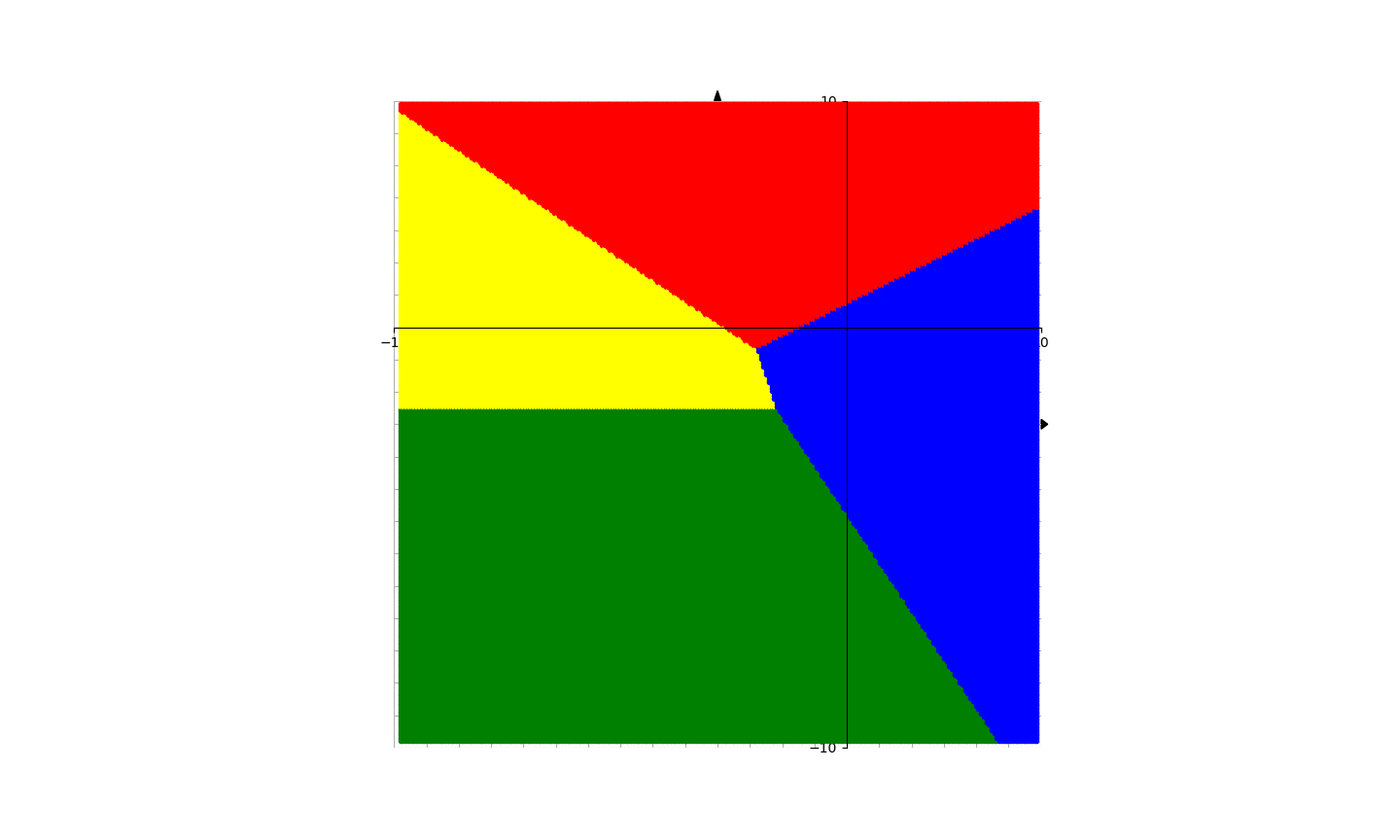}
    \includegraphics[width=3cm]{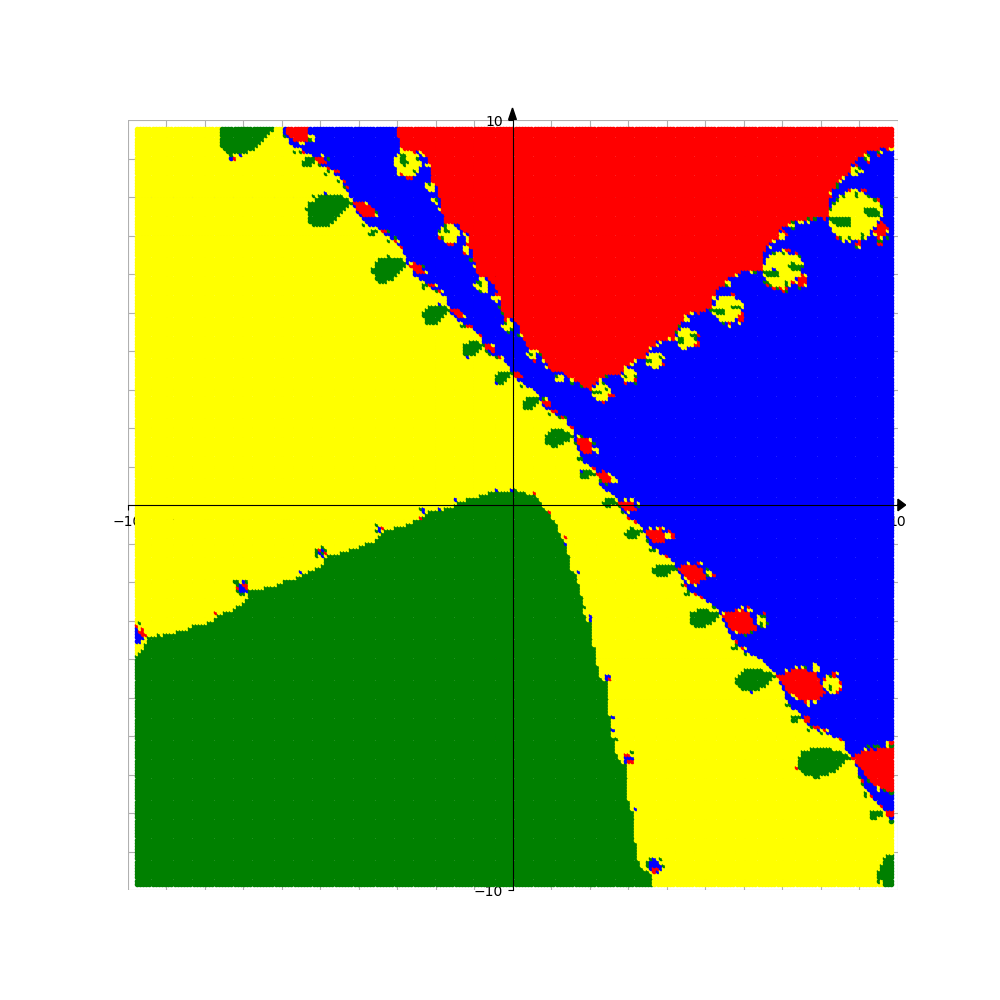}
    \includegraphics[width=3cm]{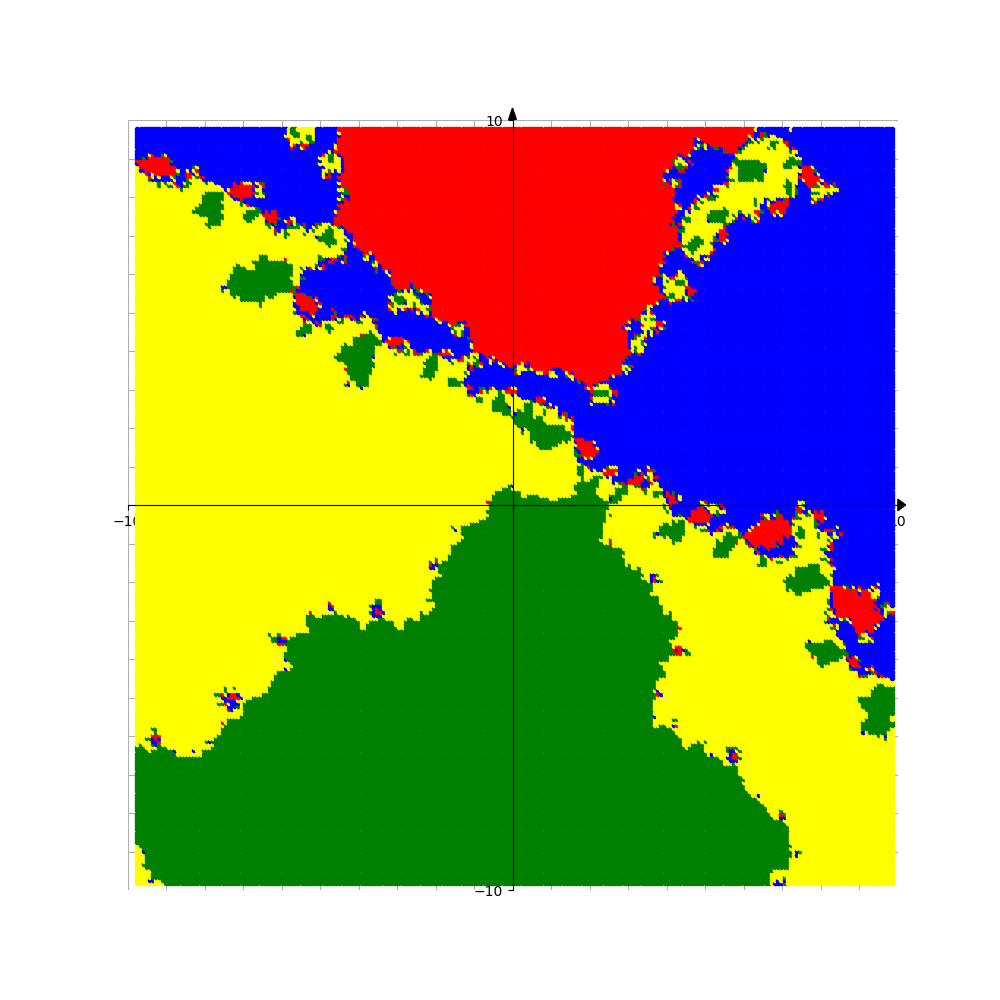}

    \bigskip
    \includegraphics[width=5.5cm]{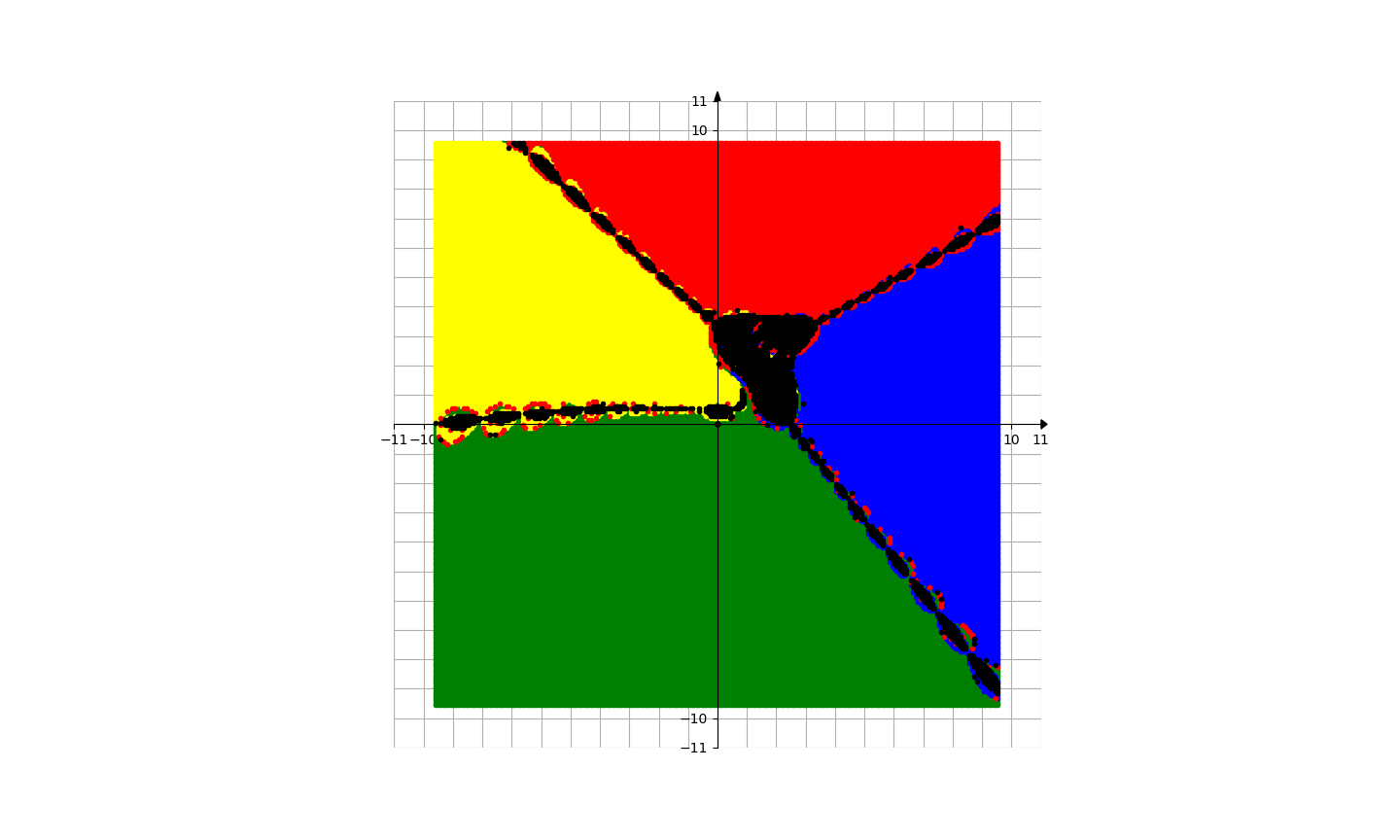}
    \includegraphics[width=3cm]{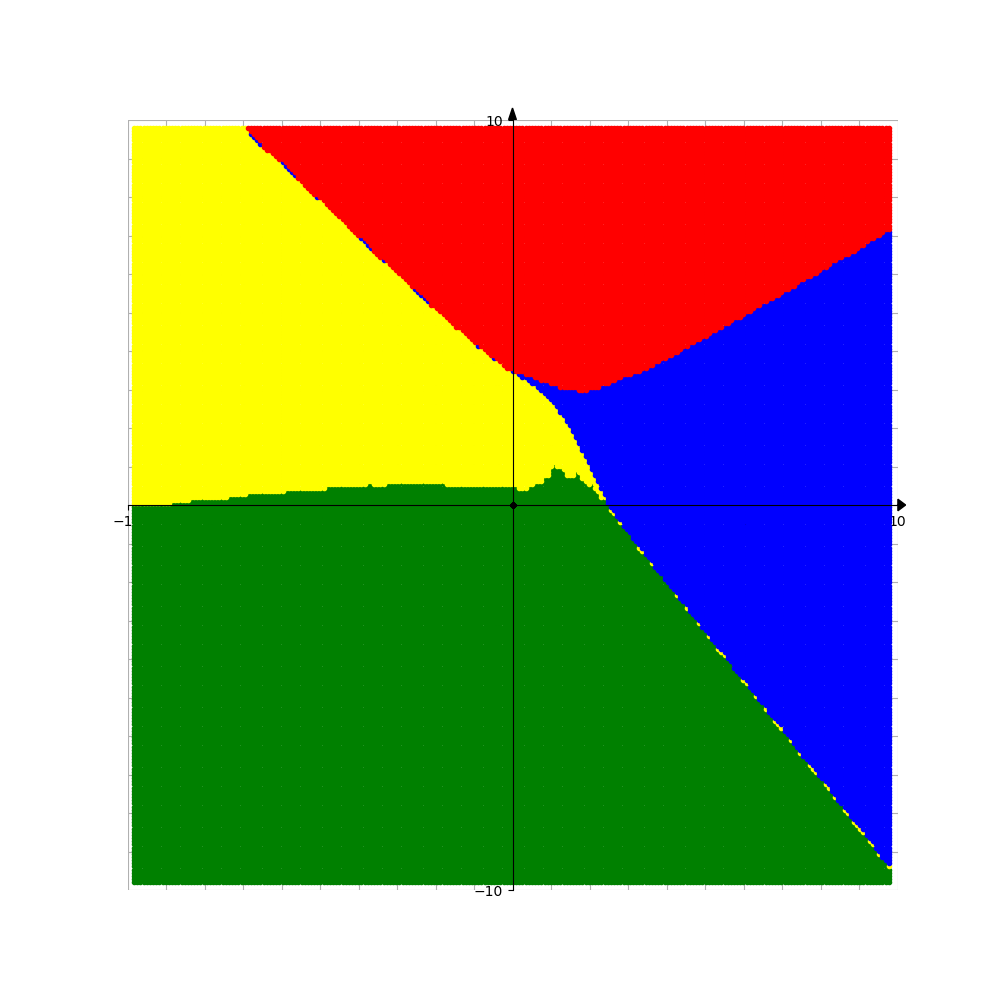}

    \bigskip
    \includegraphics[width=3cm]{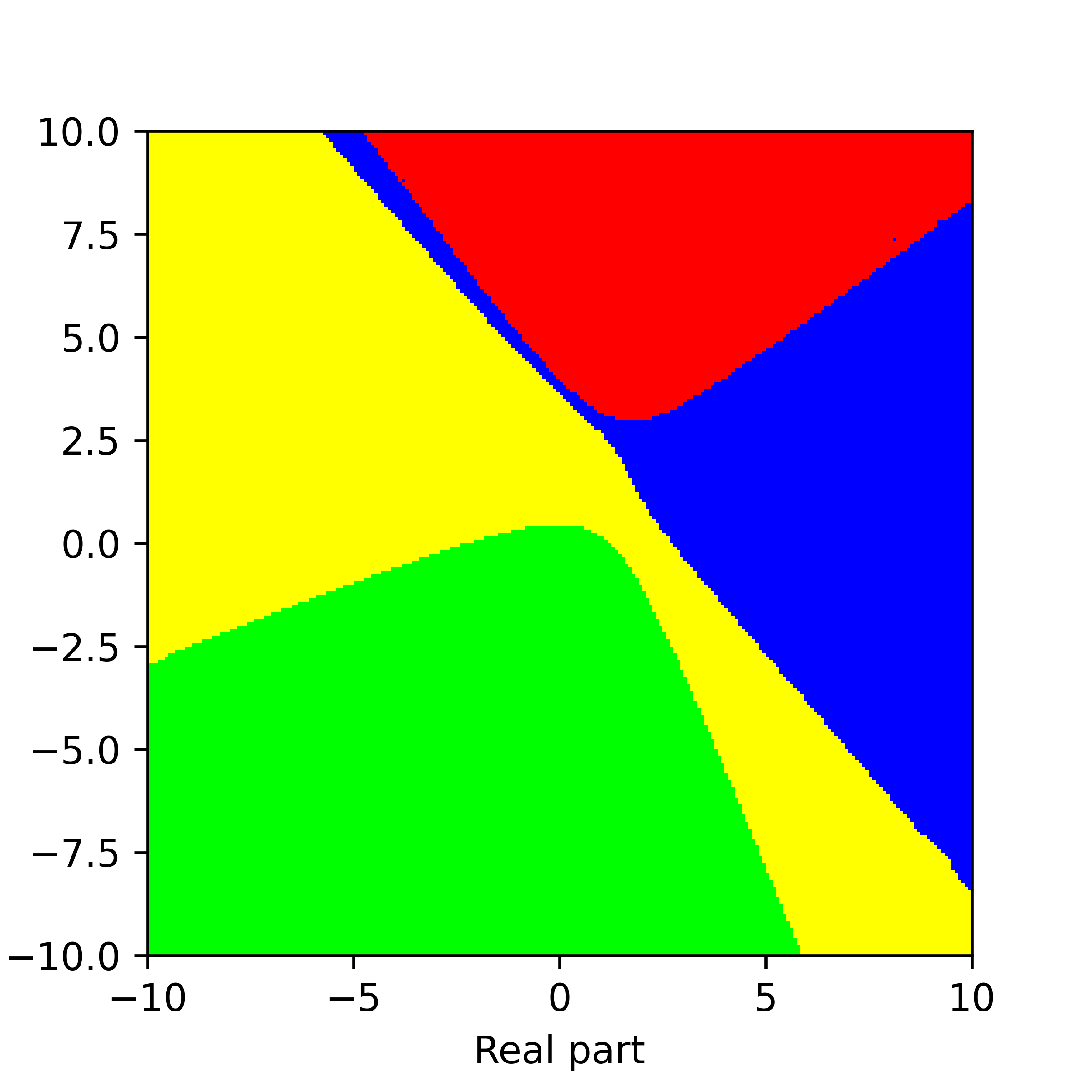}
    \includegraphics[width=3cm]{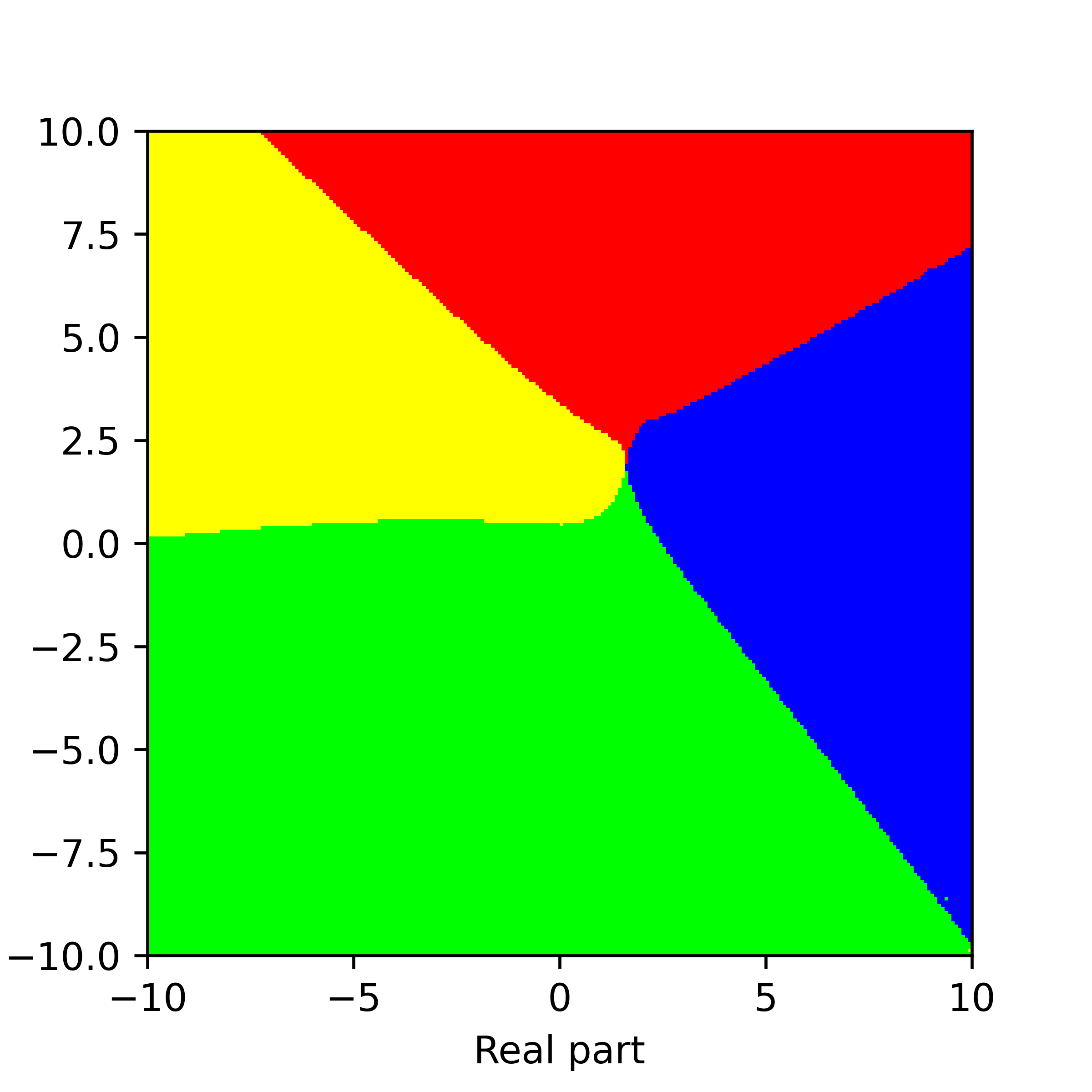}
    \includegraphics[width=3cm]{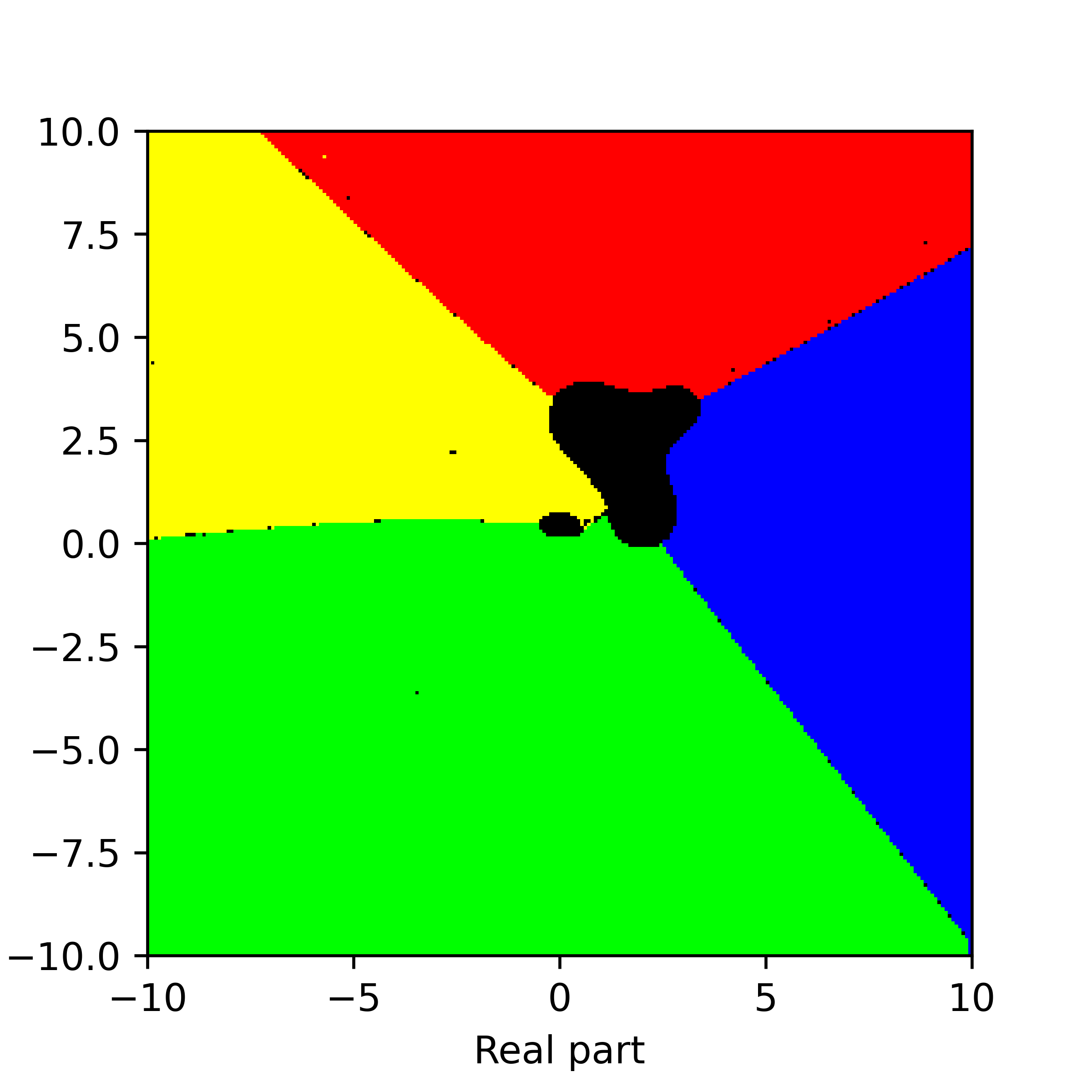}
    
    \caption{Basins of attraction for finding roots of the function $f_5$ by different methods. Pictures are referenced to from top to bottom, from left to right. Row 1: left picture is Voronoi's diagram, central picture is for Newton's method, right picture is for Random Relaxed Newton's method. Row 2: left picture is for Newton's method vOptimization, right picture is for BNQN. Row 3: left picture is for Newton's flow, central picture is for Newton's flow vFraction, right picture is for Newton's flow vOptimization. The black points in some of these pictures are those in the basin of attraction of critical points of $f_5$.}
    \label{fig:f5}
\end{figure}

\begin{figure}
    \centering
    \includegraphics[width=5cm]{Voronoif3Reduced.png}
    \includegraphics[width=3cm]{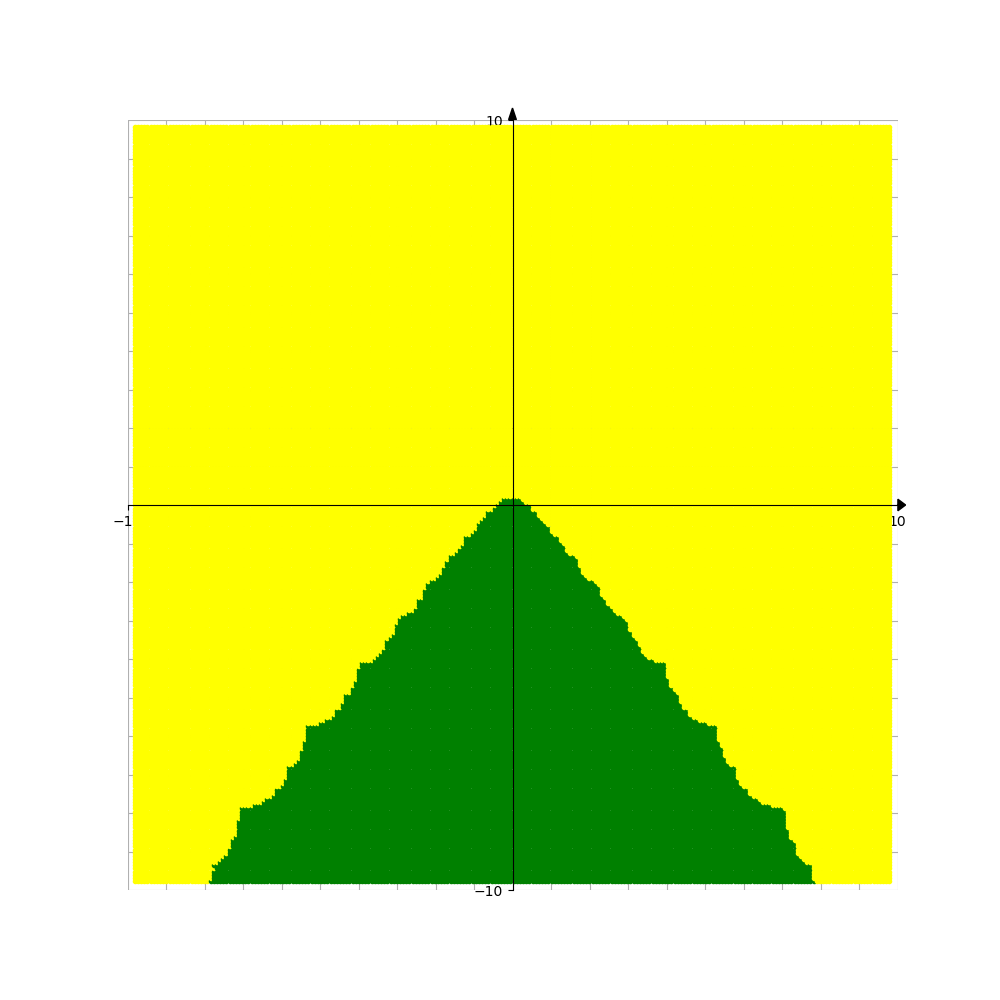}
    \includegraphics[width=3cm]{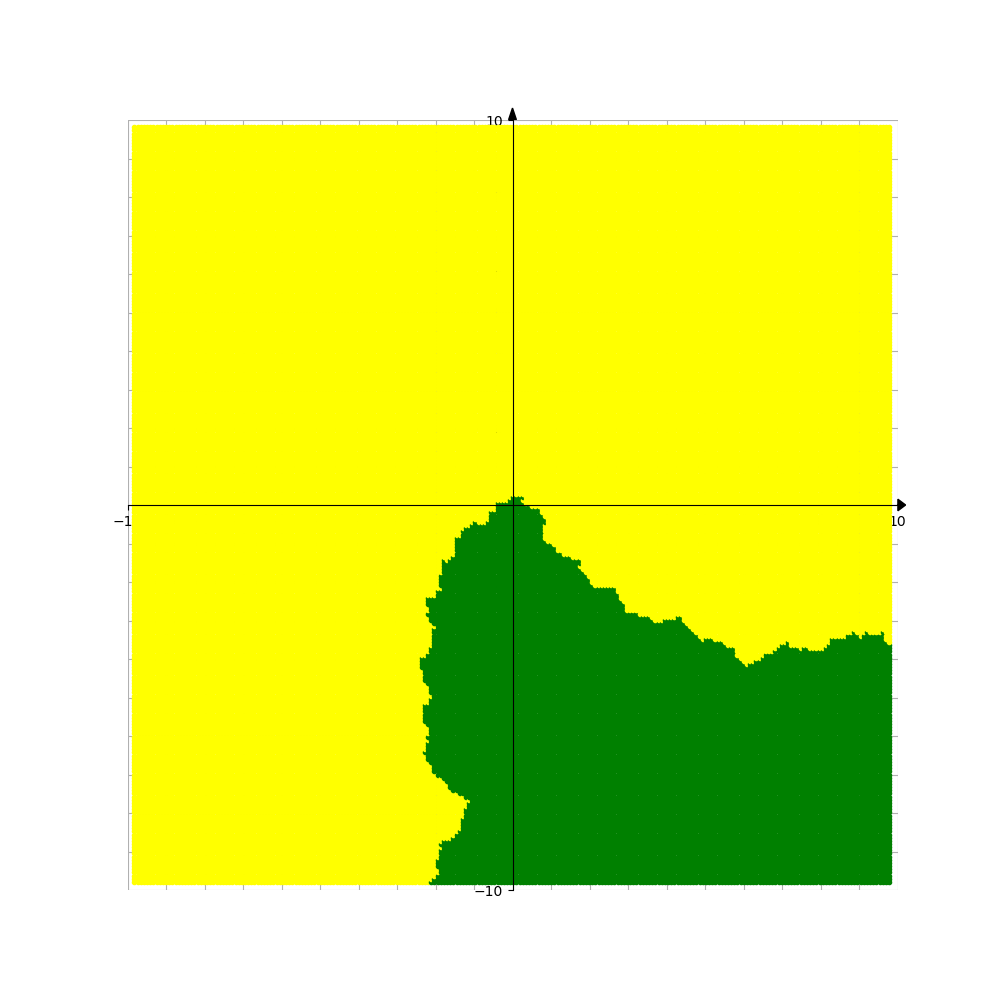}

    \bigskip
    \includegraphics[width=5.5cm]{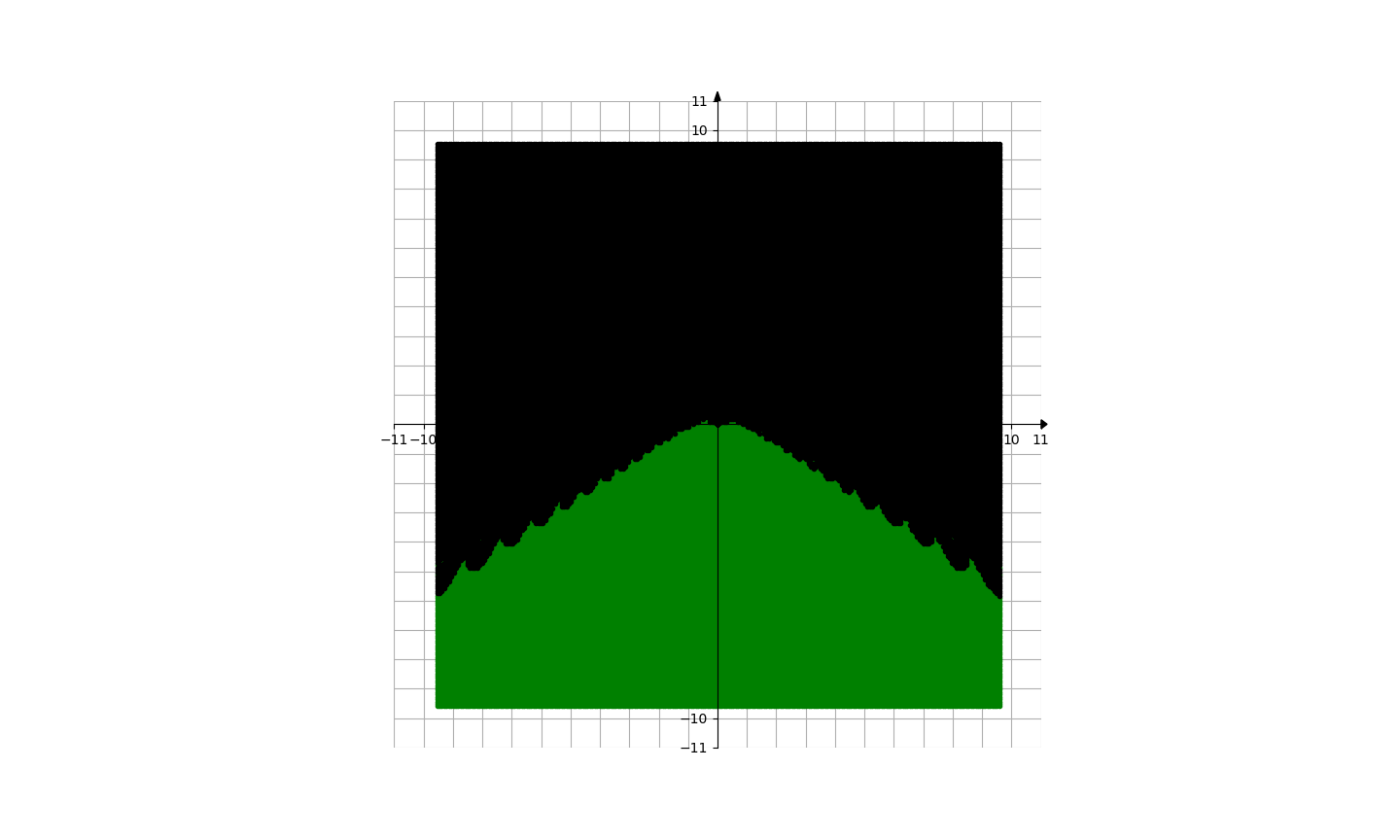}
    \includegraphics[width=5.5cm]{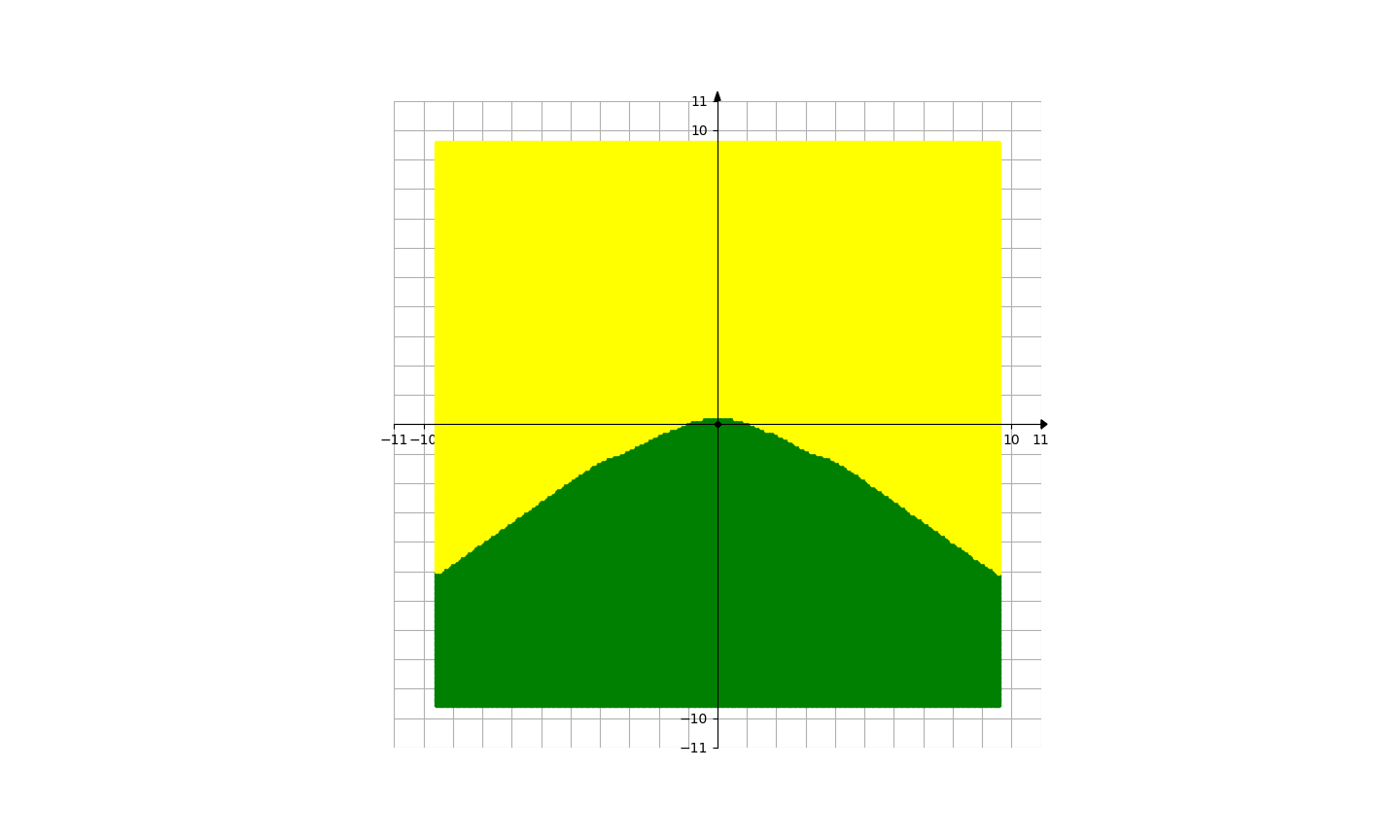}

    \bigskip
    \includegraphics[width=3cm]{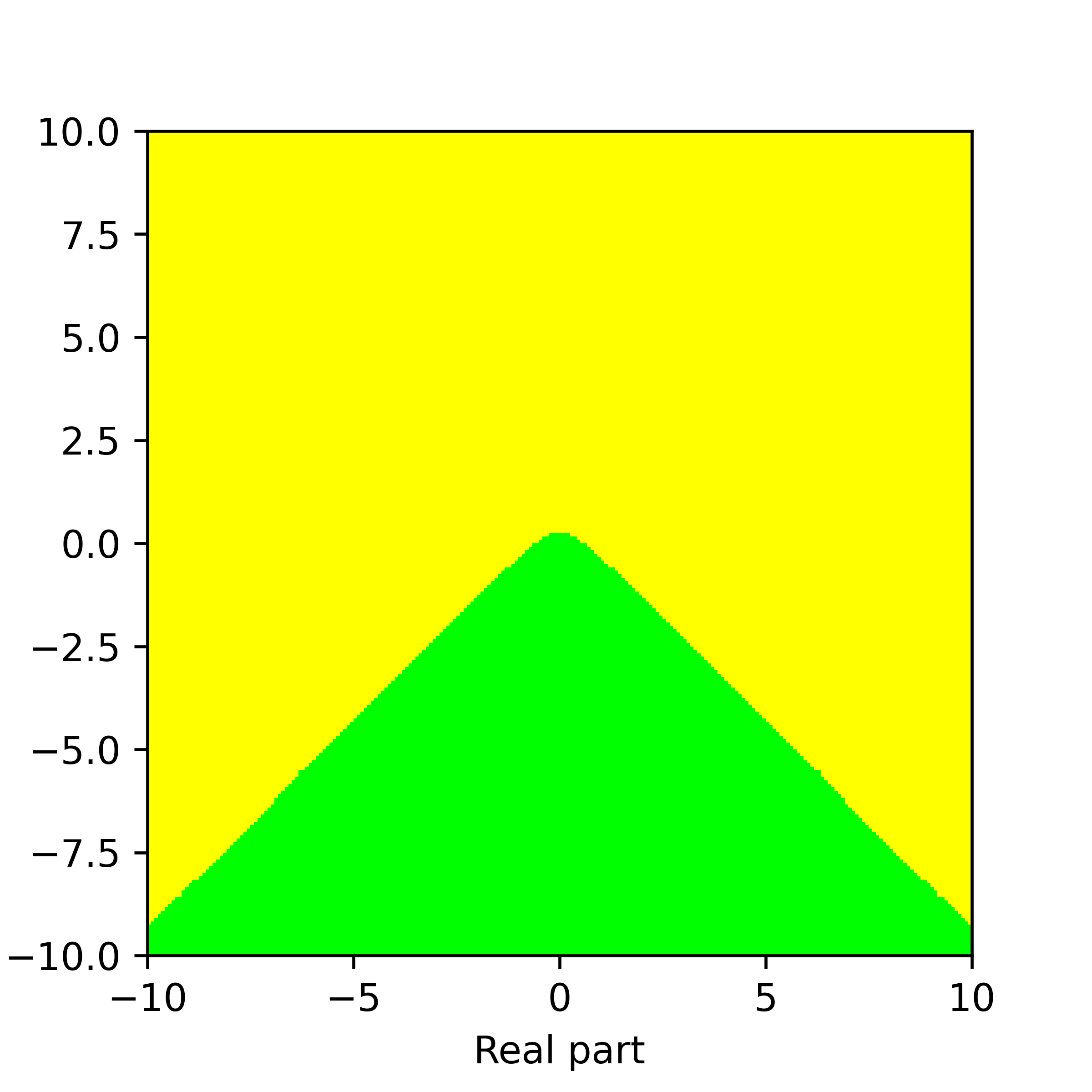}
    \includegraphics[width=3cm]{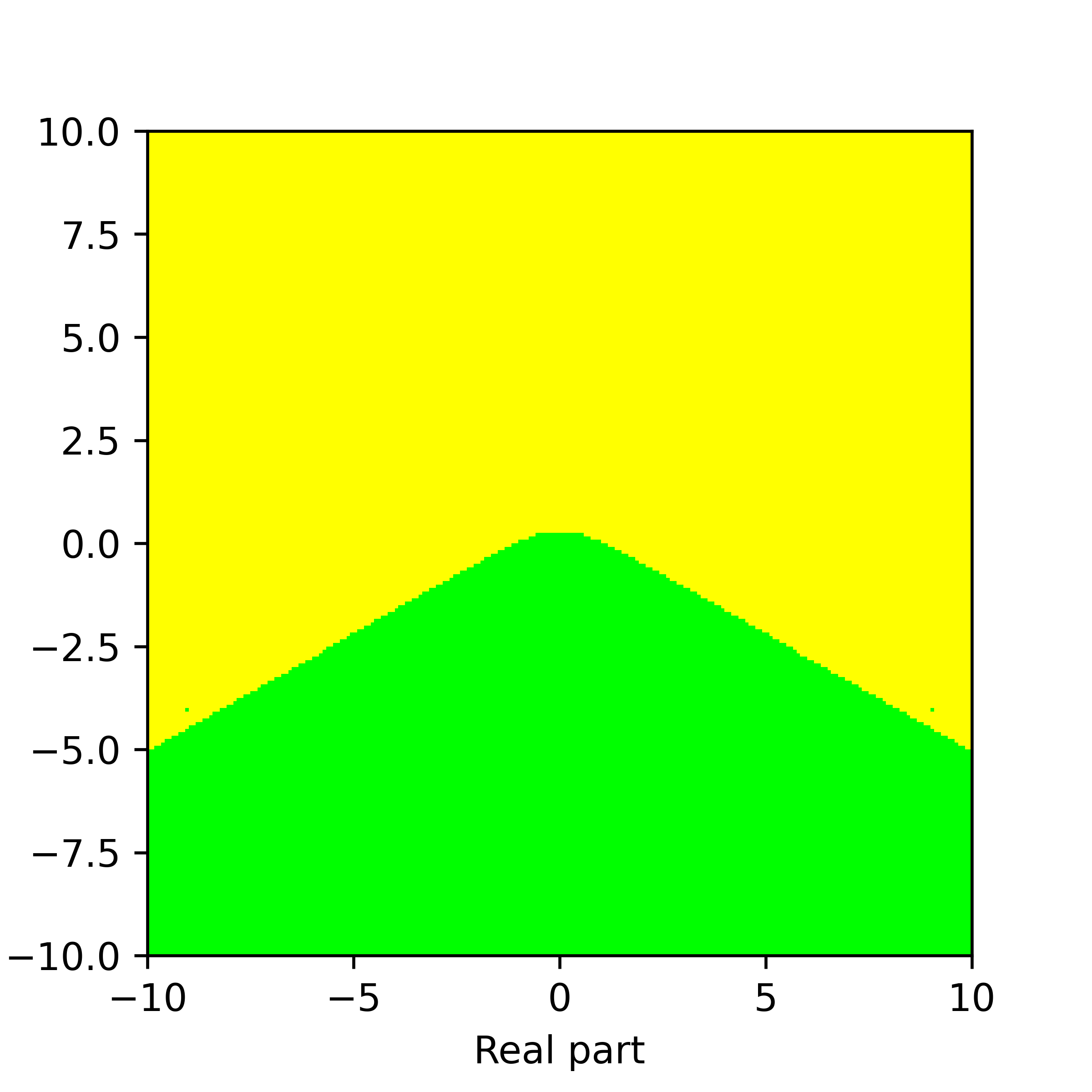}
    \includegraphics[width=3cm]{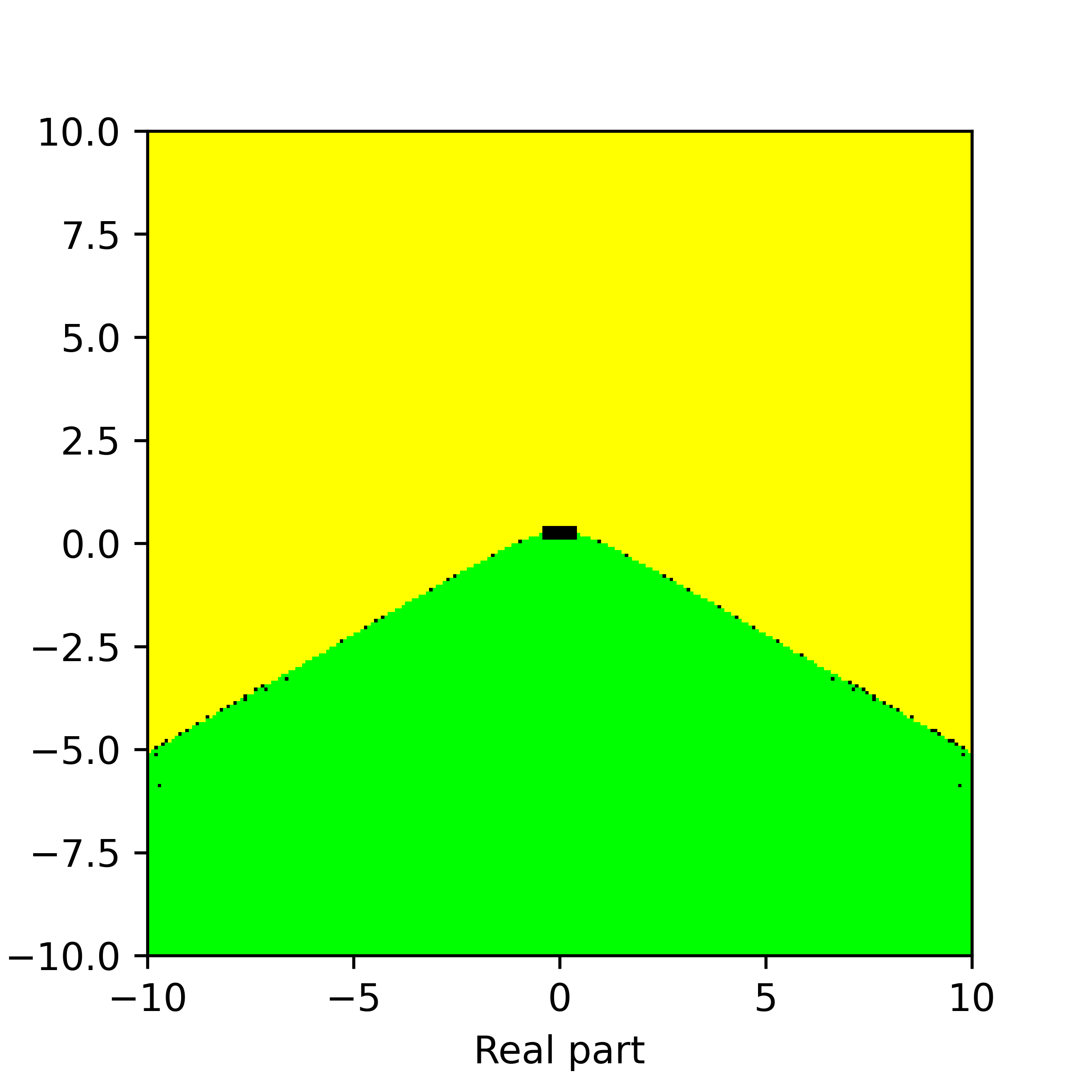}
    
    \caption{Basins of attraction for finding roots of the function $f_6$ by different methods. Pictures are referenced to from top to bottom, from left to right. Row 1: left picture is Voronoi's diagram, central picture is for Newton's method, right picture is for Random Relaxed Newton's method. Row 2: left picture is for Newton's method vOptimization, right picture is for BNQN. Row 3: left picture is for Newton's flow, central picture is for Newton's flow vFraction, right picture is for Newton's flow vOptimization. The black points for Newton's flow vOptimization are those in the basin of attraction of critical points of $f_6$ which are not roots, while the black points for Newton's method for Optimization seem to due to the fact that the cost function is highly degenerate.}
    \label{fig:f6}
\end{figure}

\begin{figure}
    \centering
    \includegraphics[width=5cm]{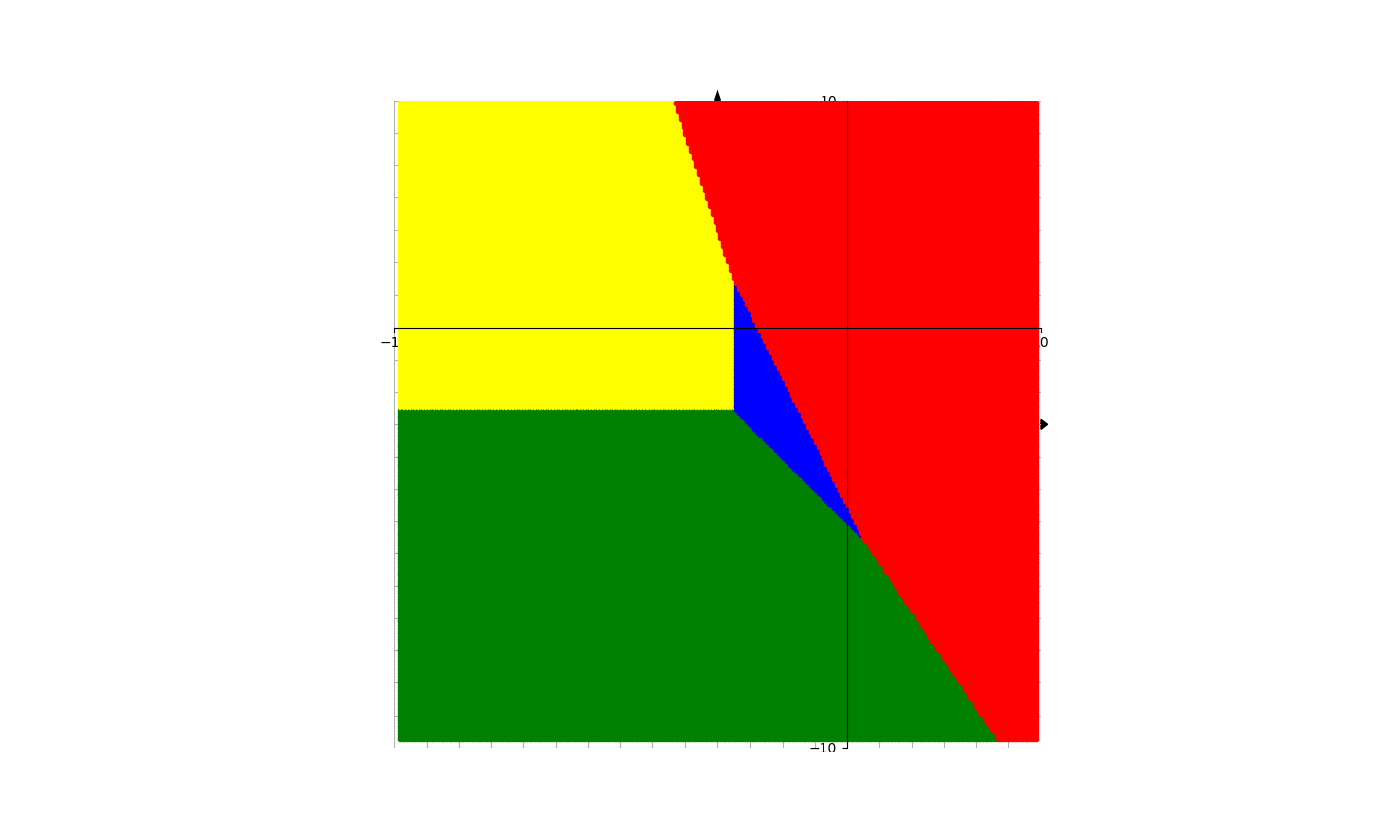}
    \includegraphics[width=3cm]{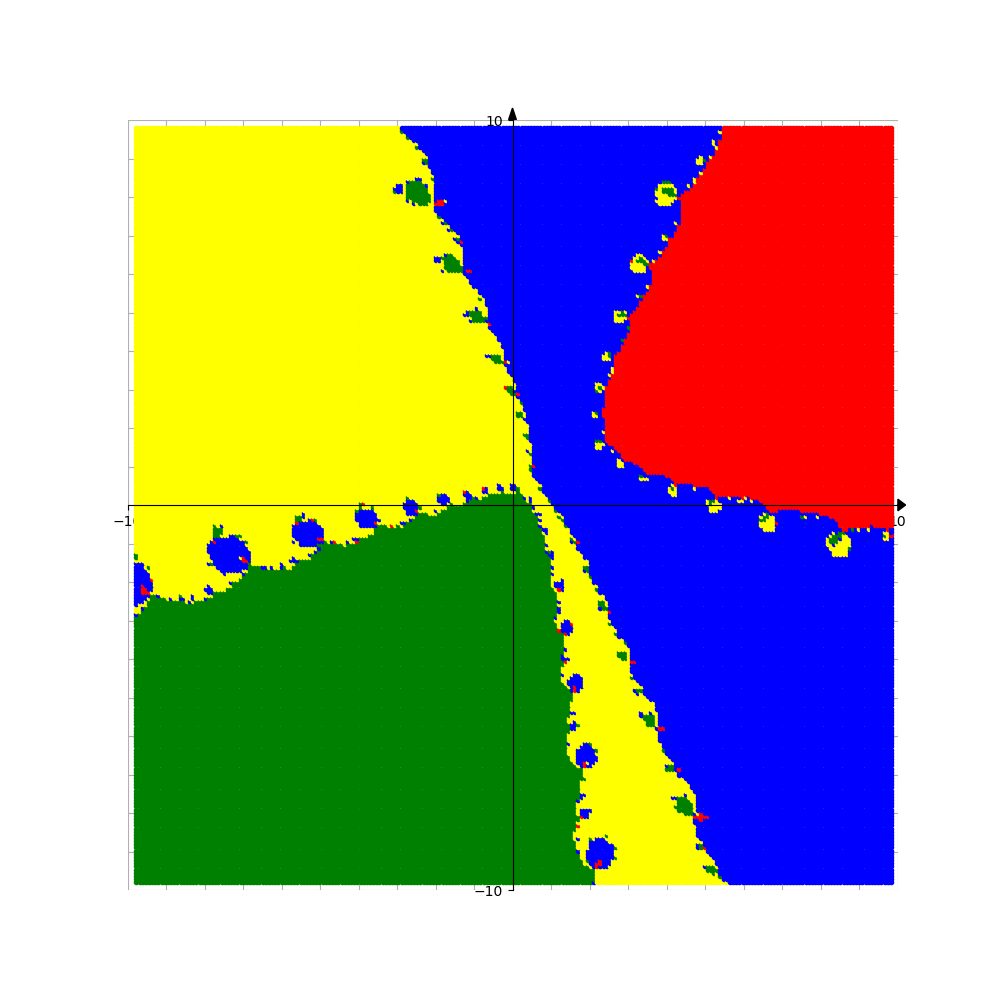}
    \includegraphics[width=3cm]{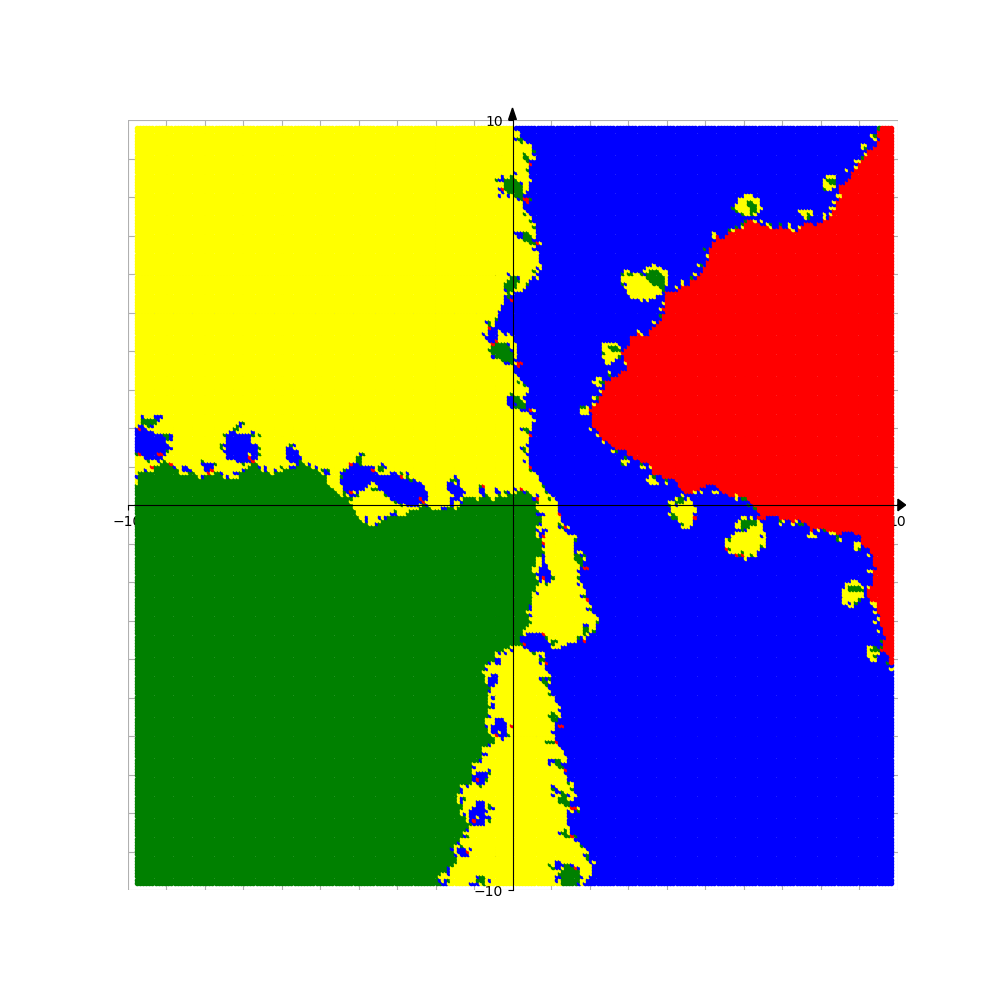}

    \bigskip
    \includegraphics[width=5.5cm]{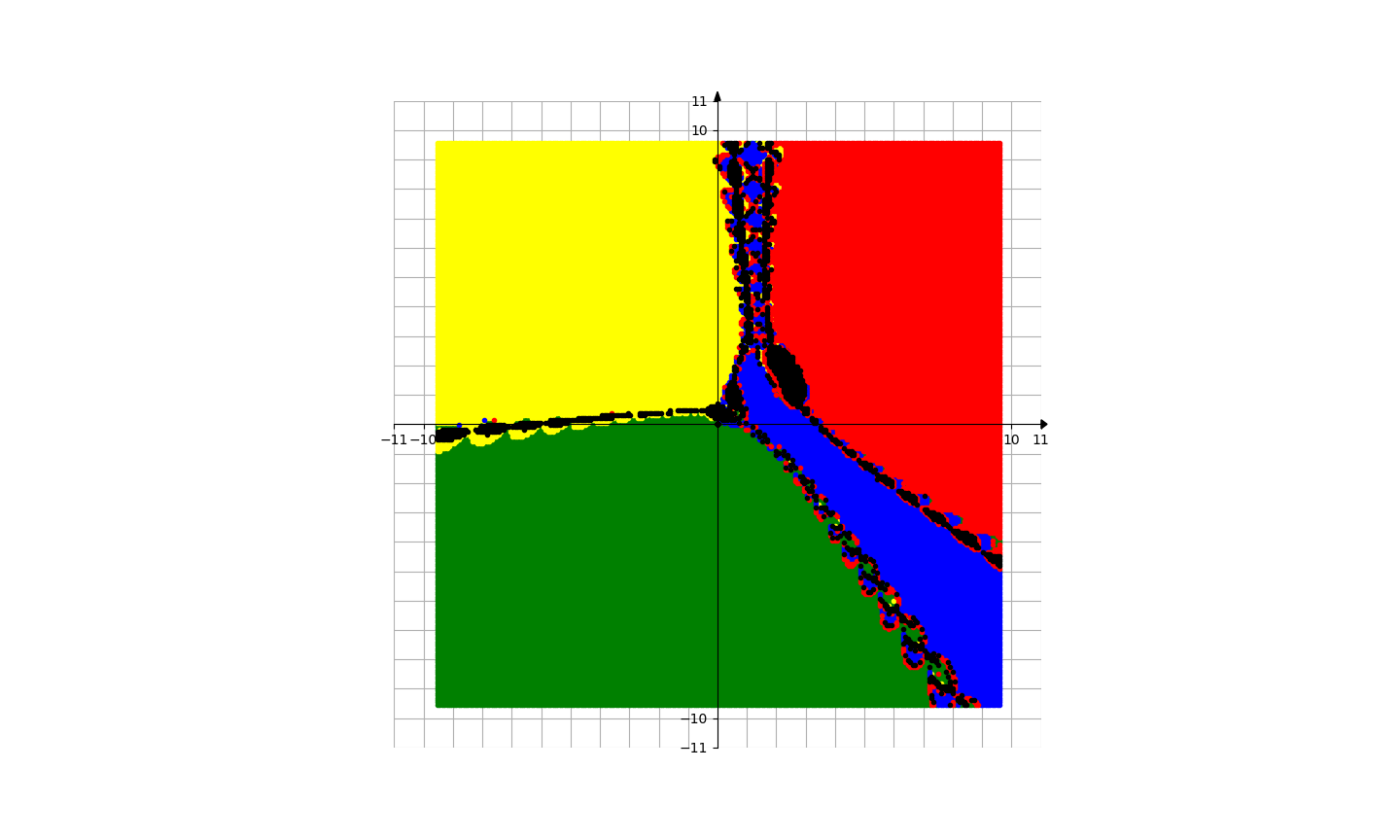}
    \includegraphics[width=3cm]{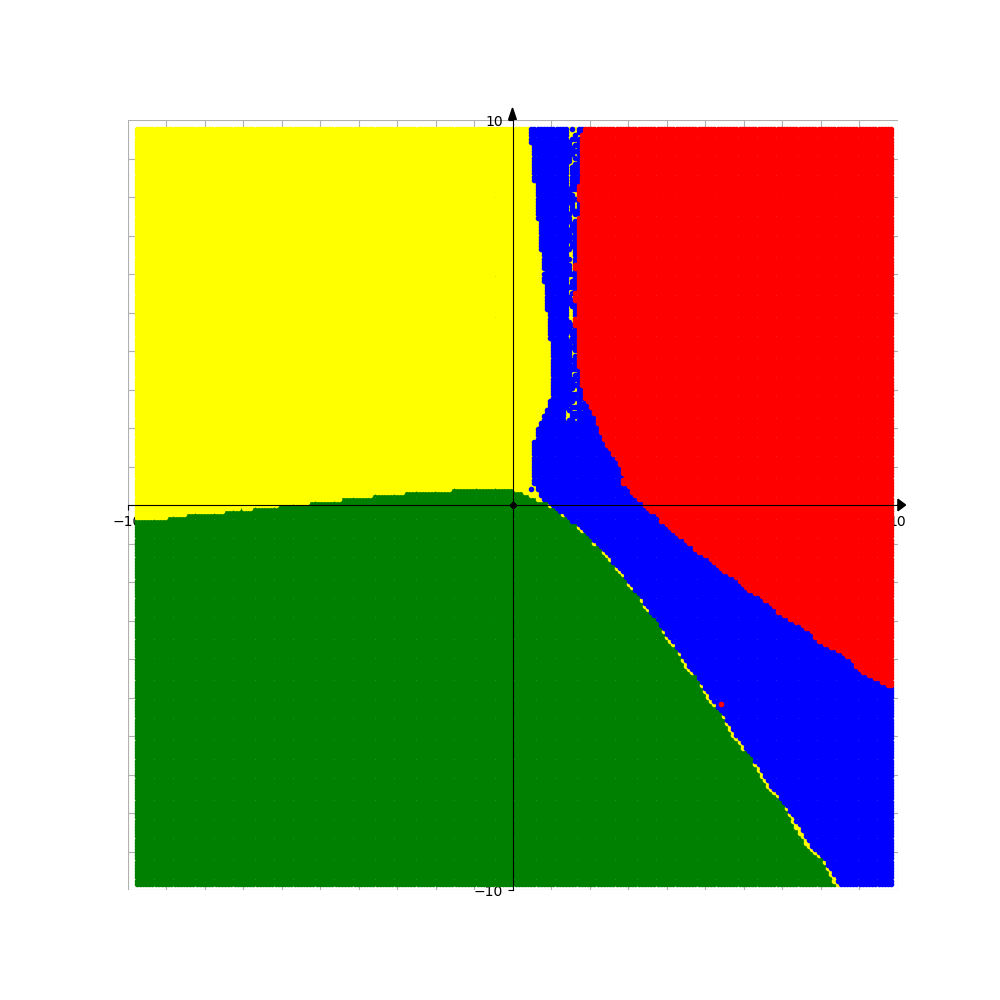}
    
    \bigskip
    \includegraphics[width=3cm]{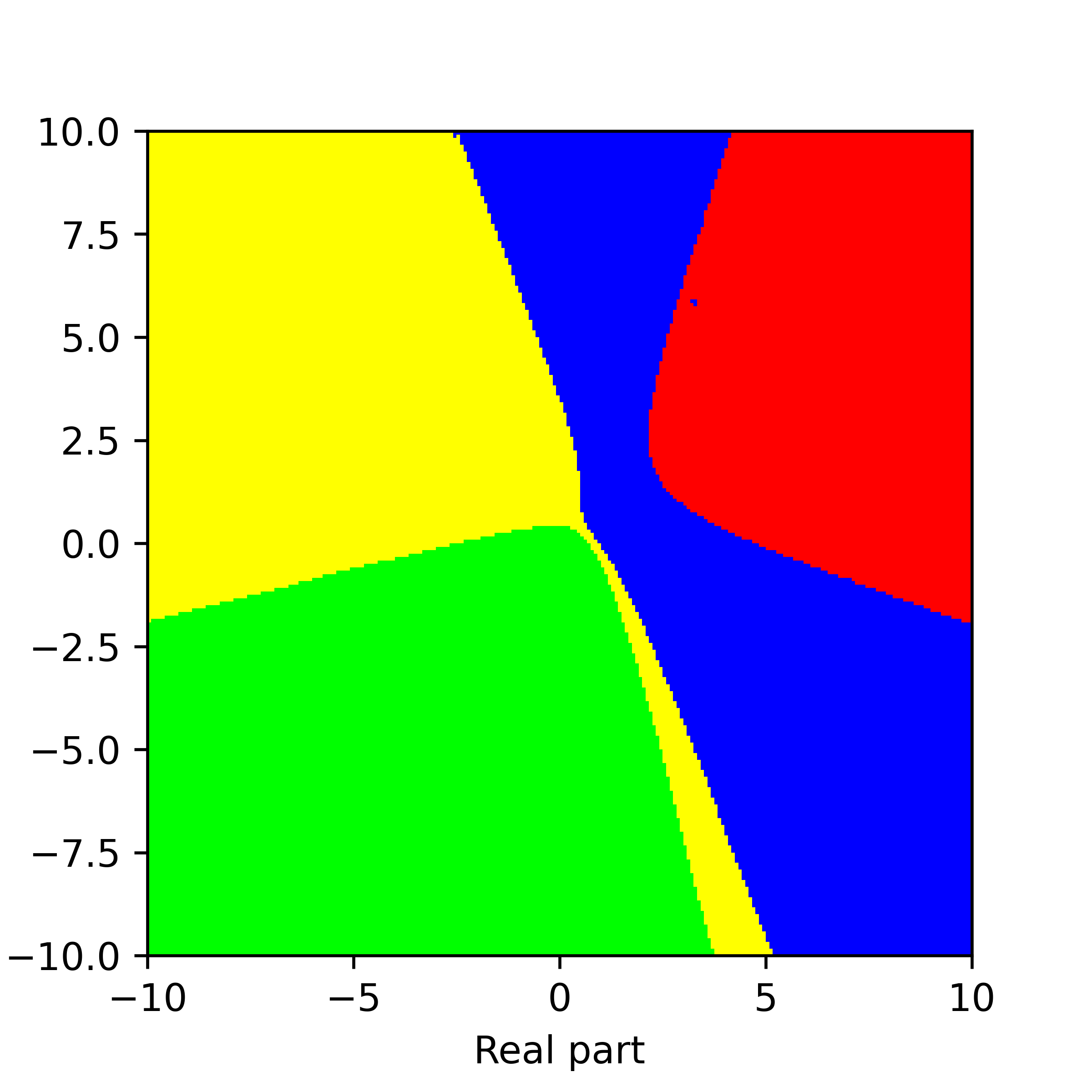}
    \includegraphics[width=3cm]{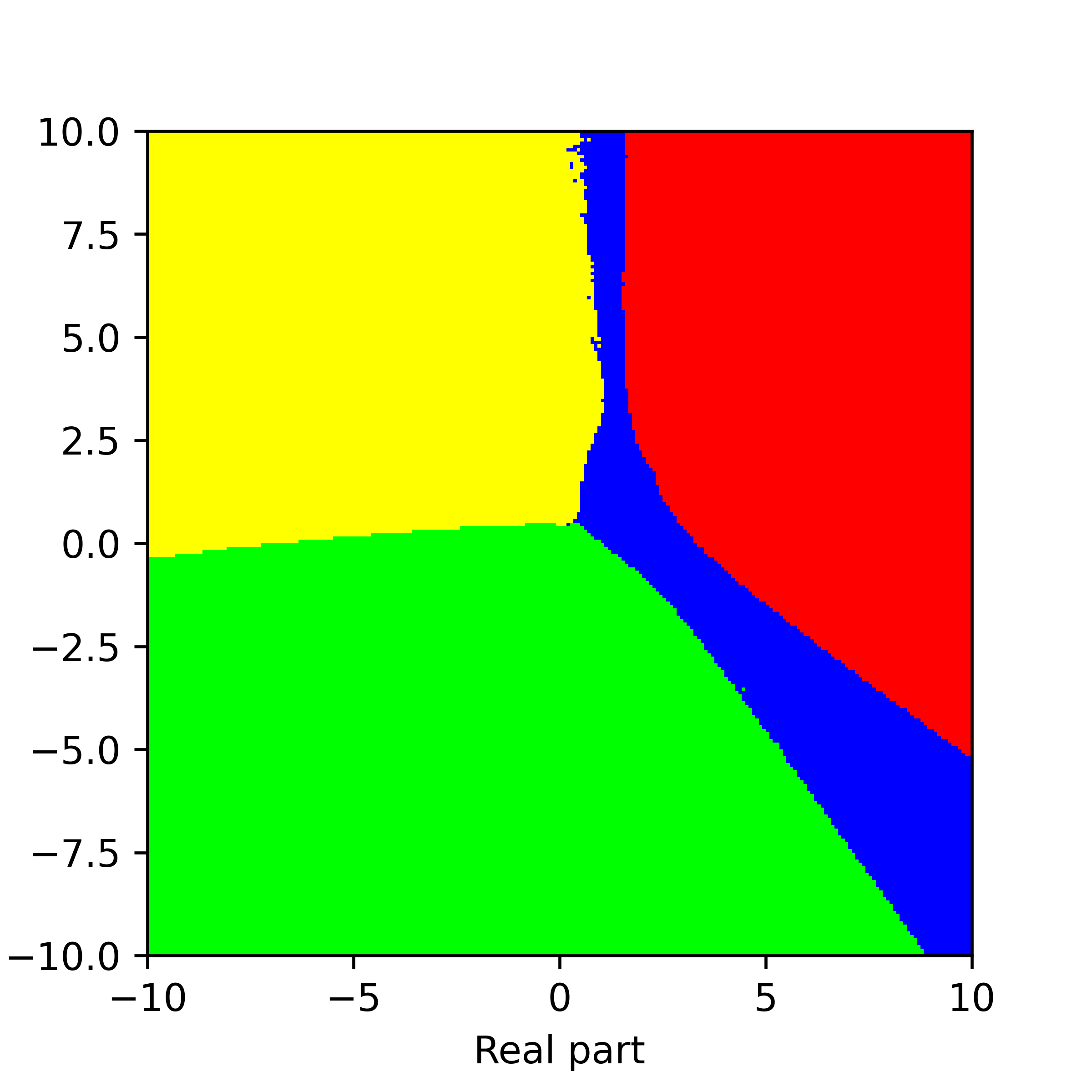}
    \includegraphics[width=3cm]{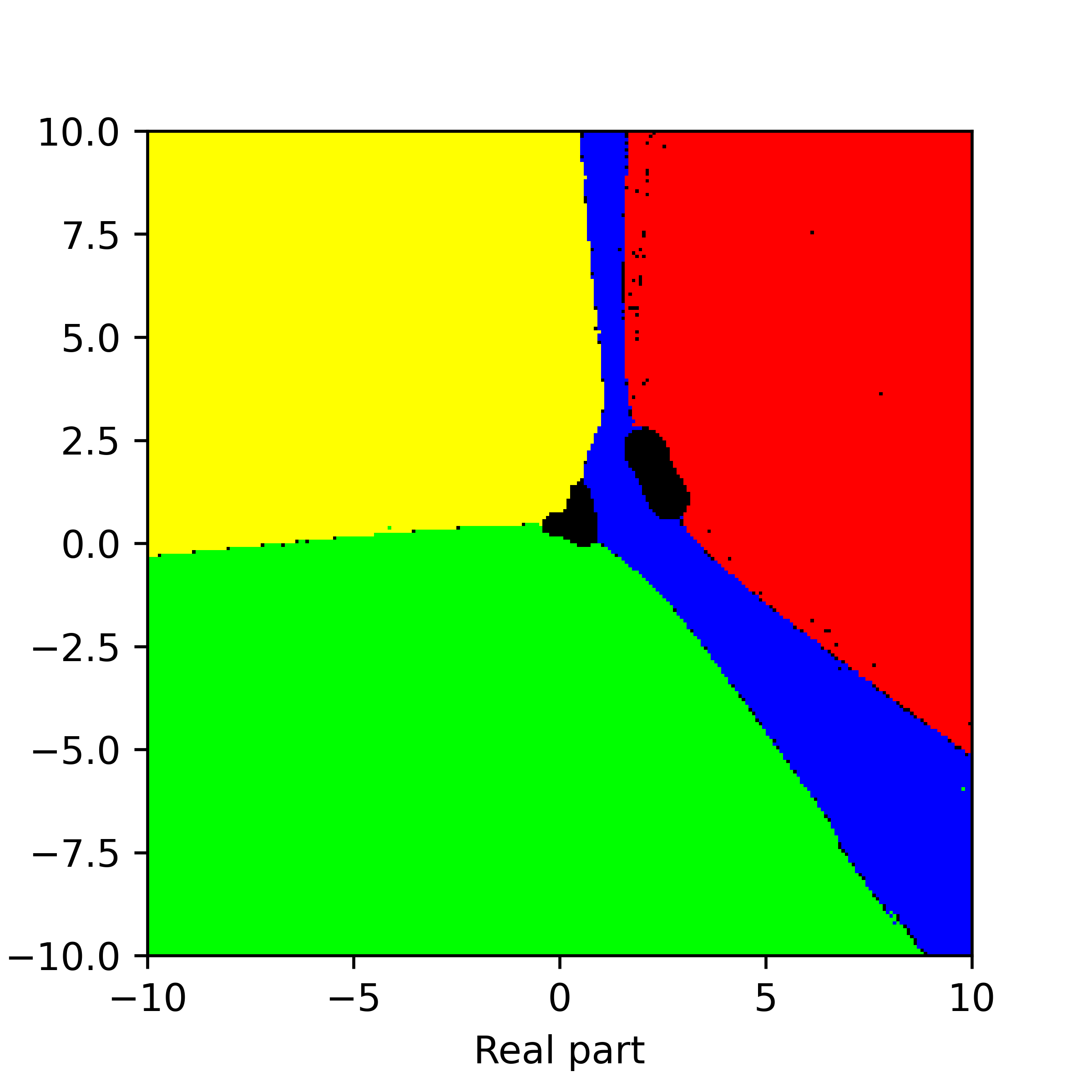}
    
    \caption{Basins of attraction for finding roots of the function $f_7$ by different methods. Pictures are referenced to from top to bottom, from left to right. Row 1: left picture is Voronoi's diagram, central picture is for Newton's method, right picture is for Random Relaxed Newton's method. Row 2: left picture is for Newton's method vOptimization, right picture is for BNQN. Row 3: left picture is for Newton's flow, central picture is for Newton's flow vFraction, right picture is for Newton's flow vOptimization. The black points in some of these pictures are those in the basin of attraction of critical points of $f_7$.}
    \label{fig:f7}
\end{figure}

\begin{figure}
    \centering
    \includegraphics[width=5cm]{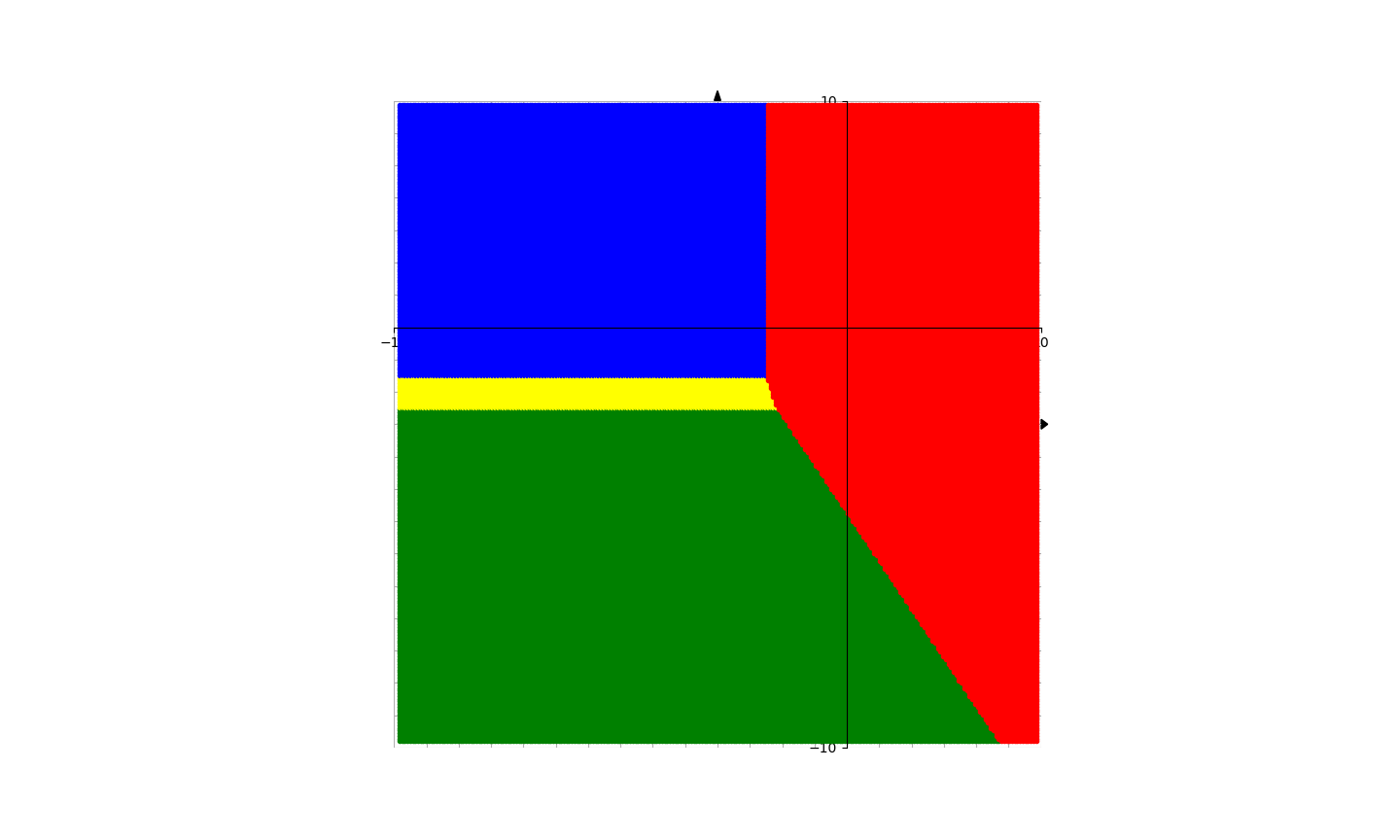}
    \includegraphics[width=3cm]{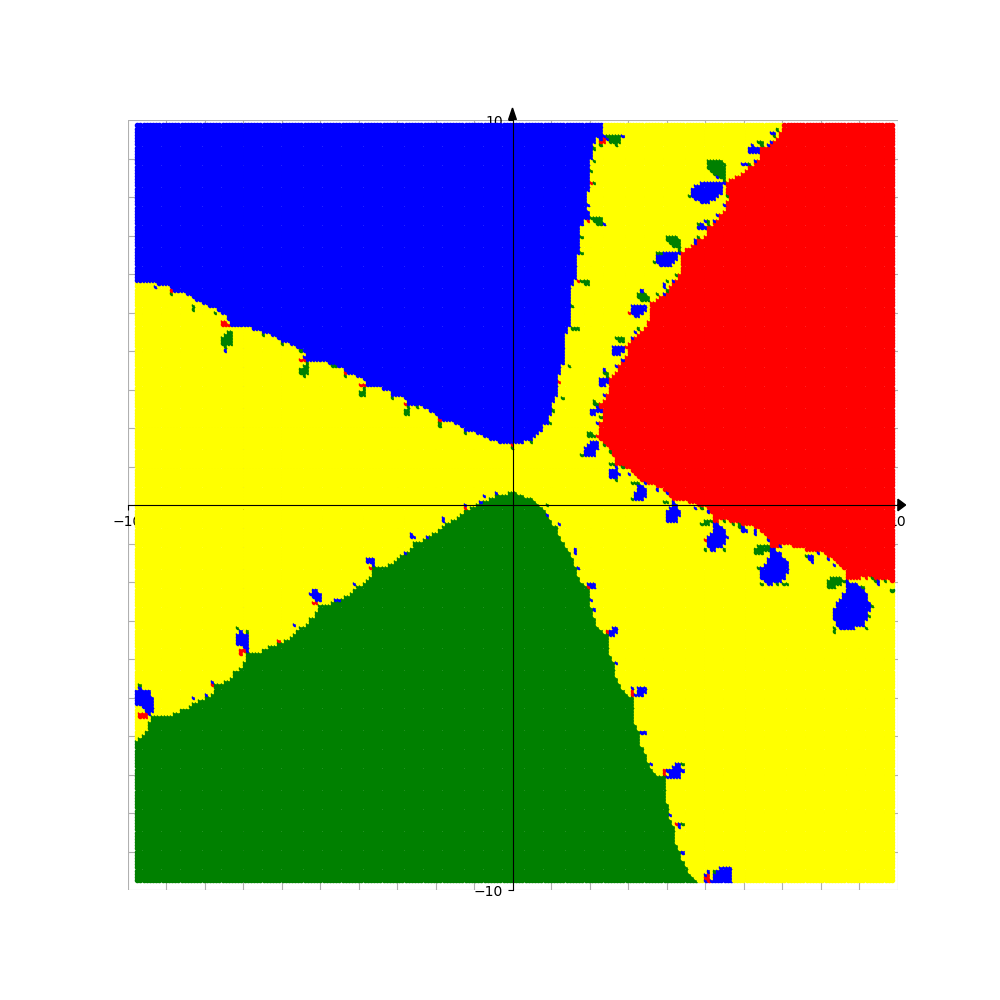}
    \includegraphics[width=3cm]{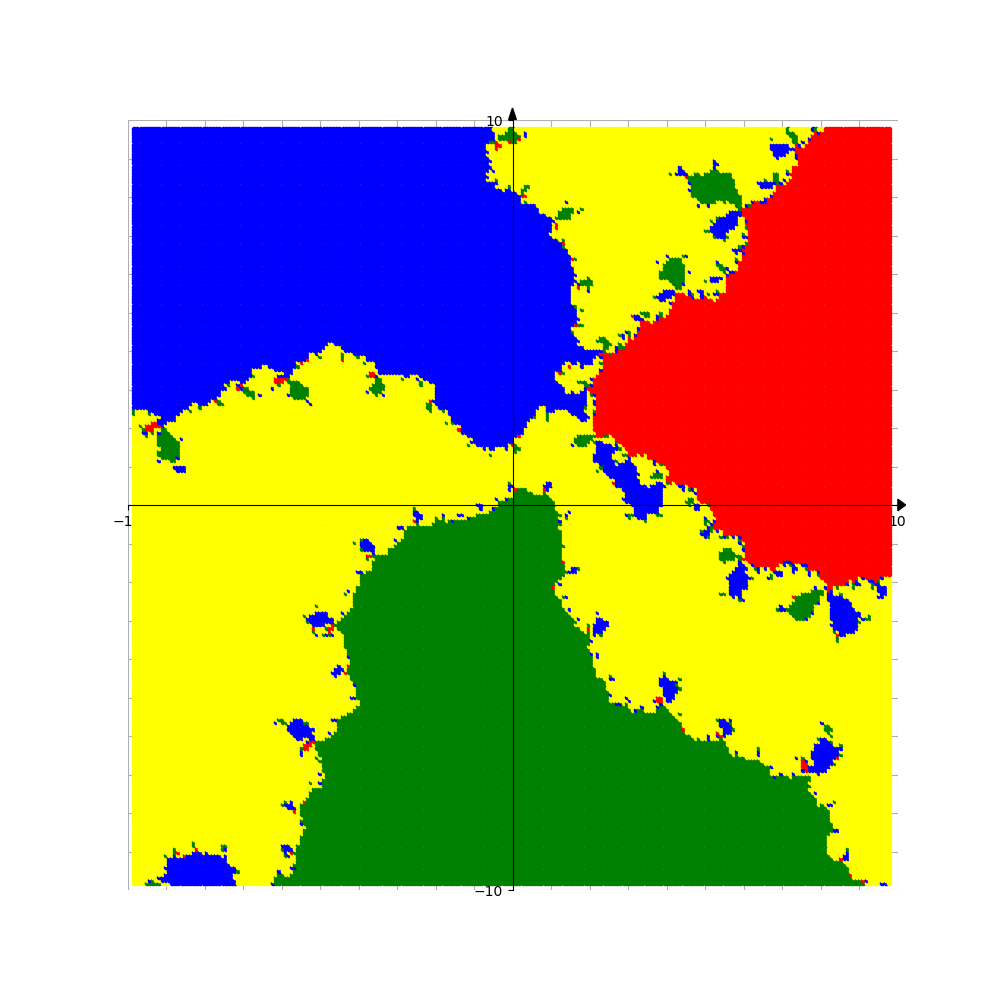}

    \bigskip
    \includegraphics[width=5.5cm]{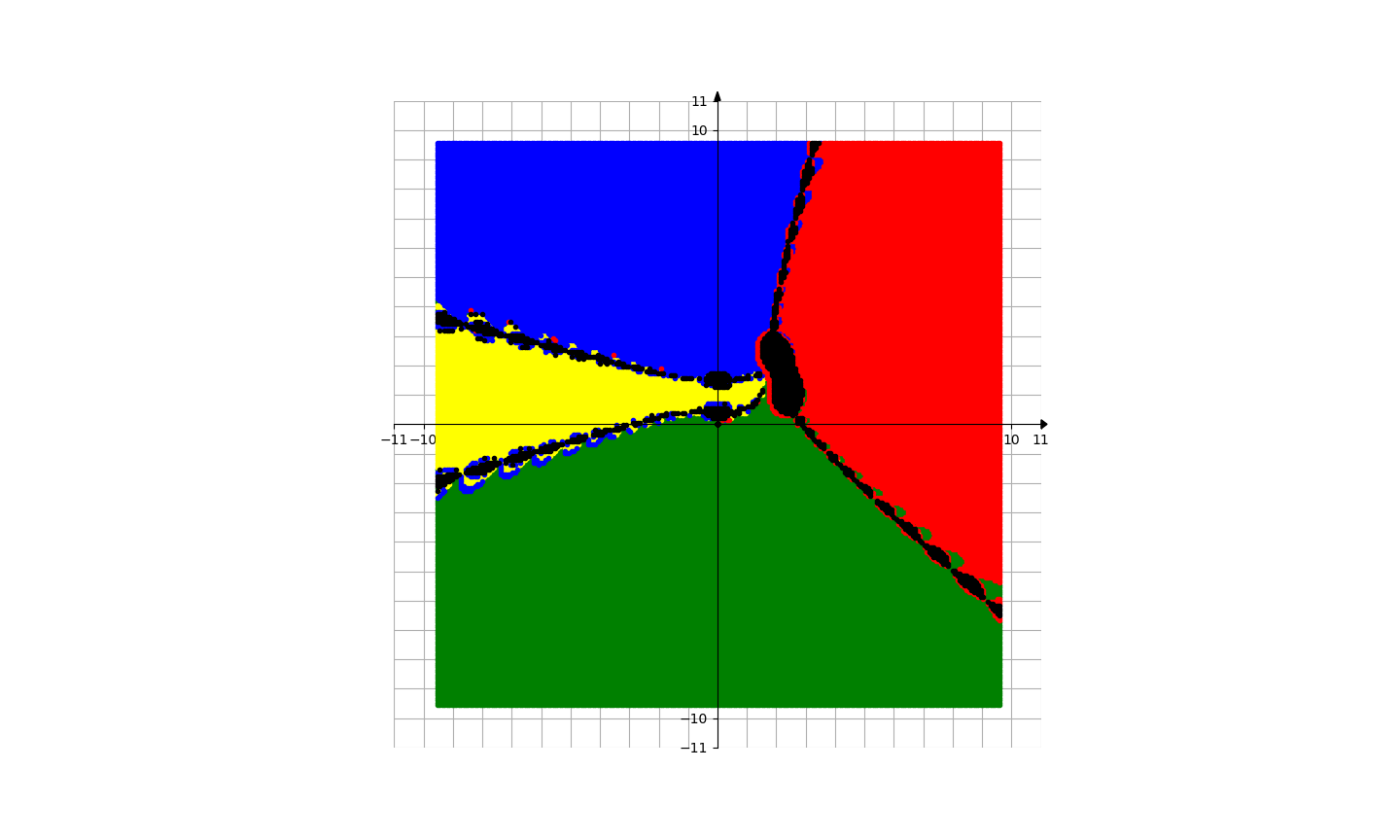}
    \includegraphics[width=3cm]{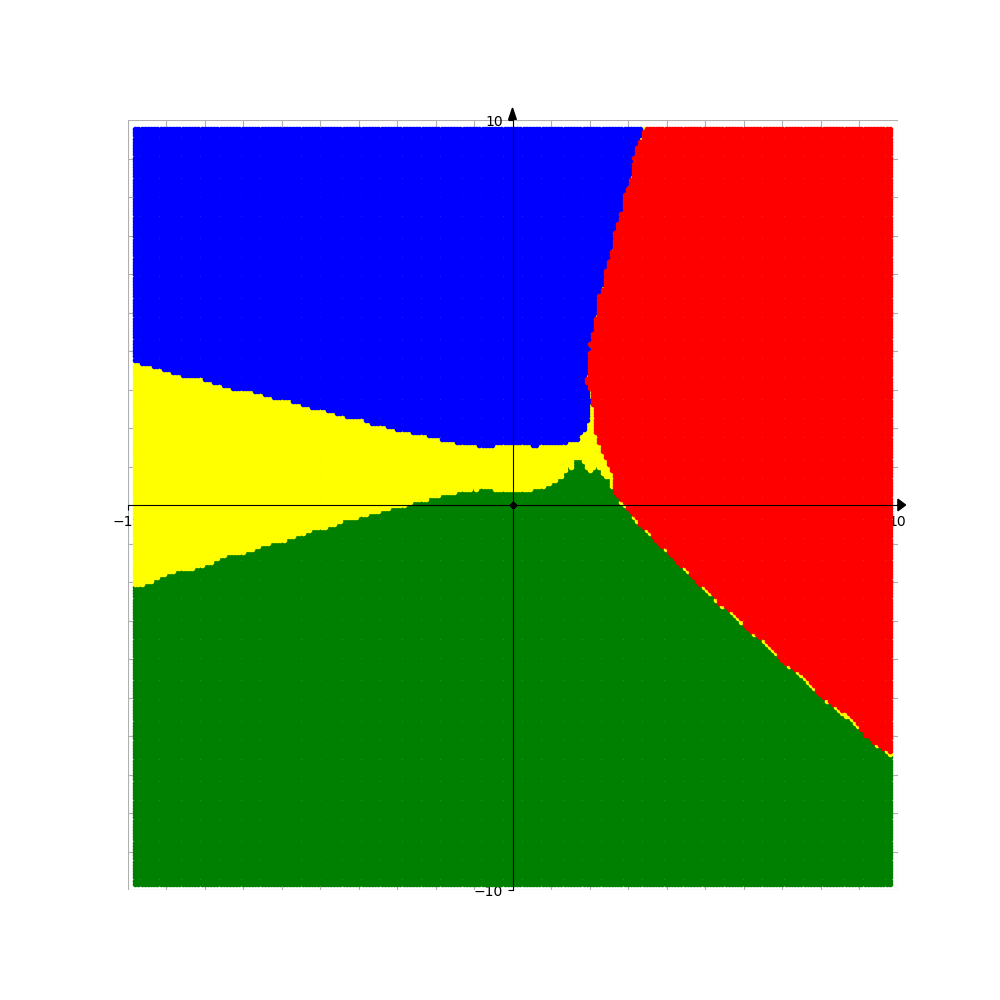}
    
    \bigskip
    \includegraphics[width=3cm]{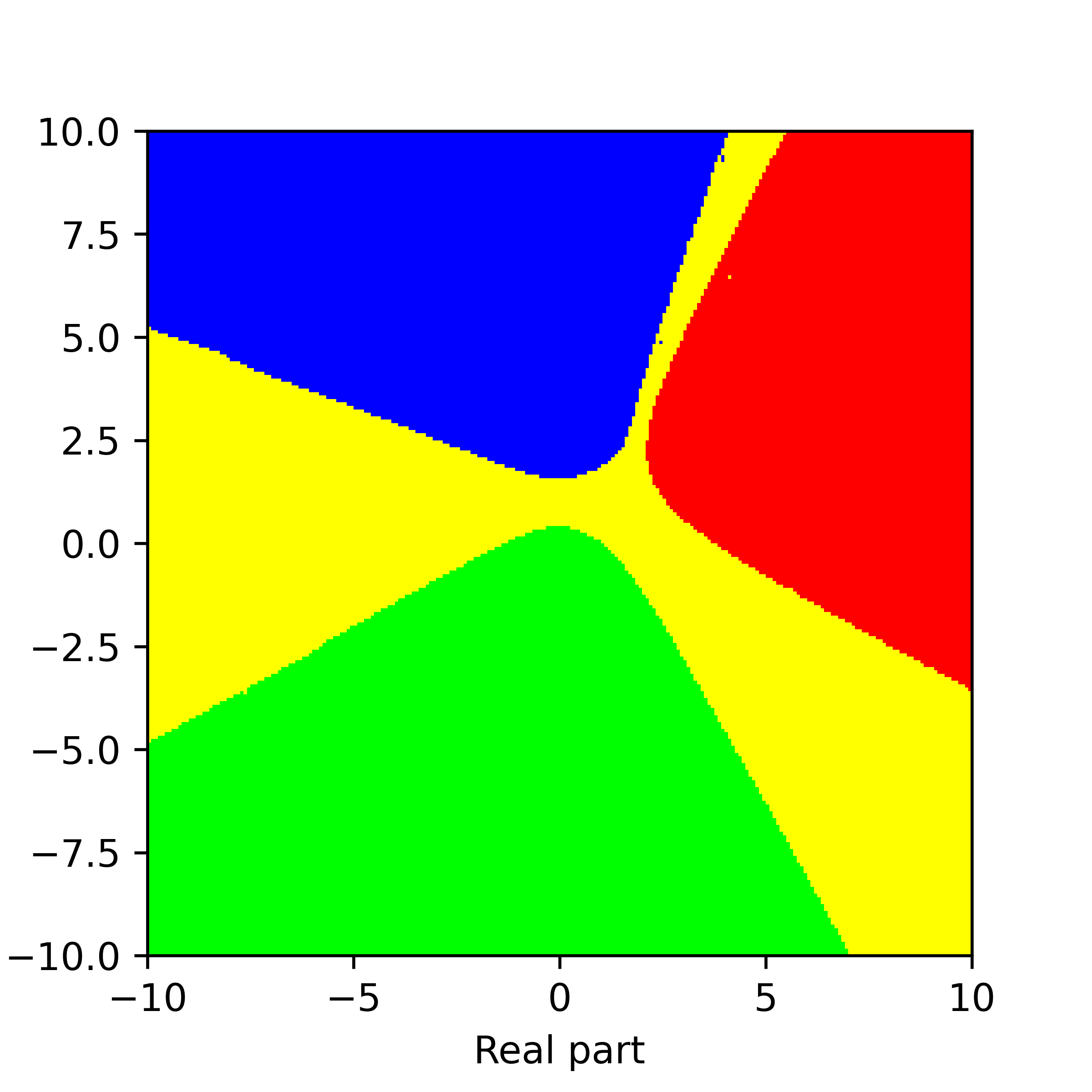}
    \includegraphics[width=3cm]{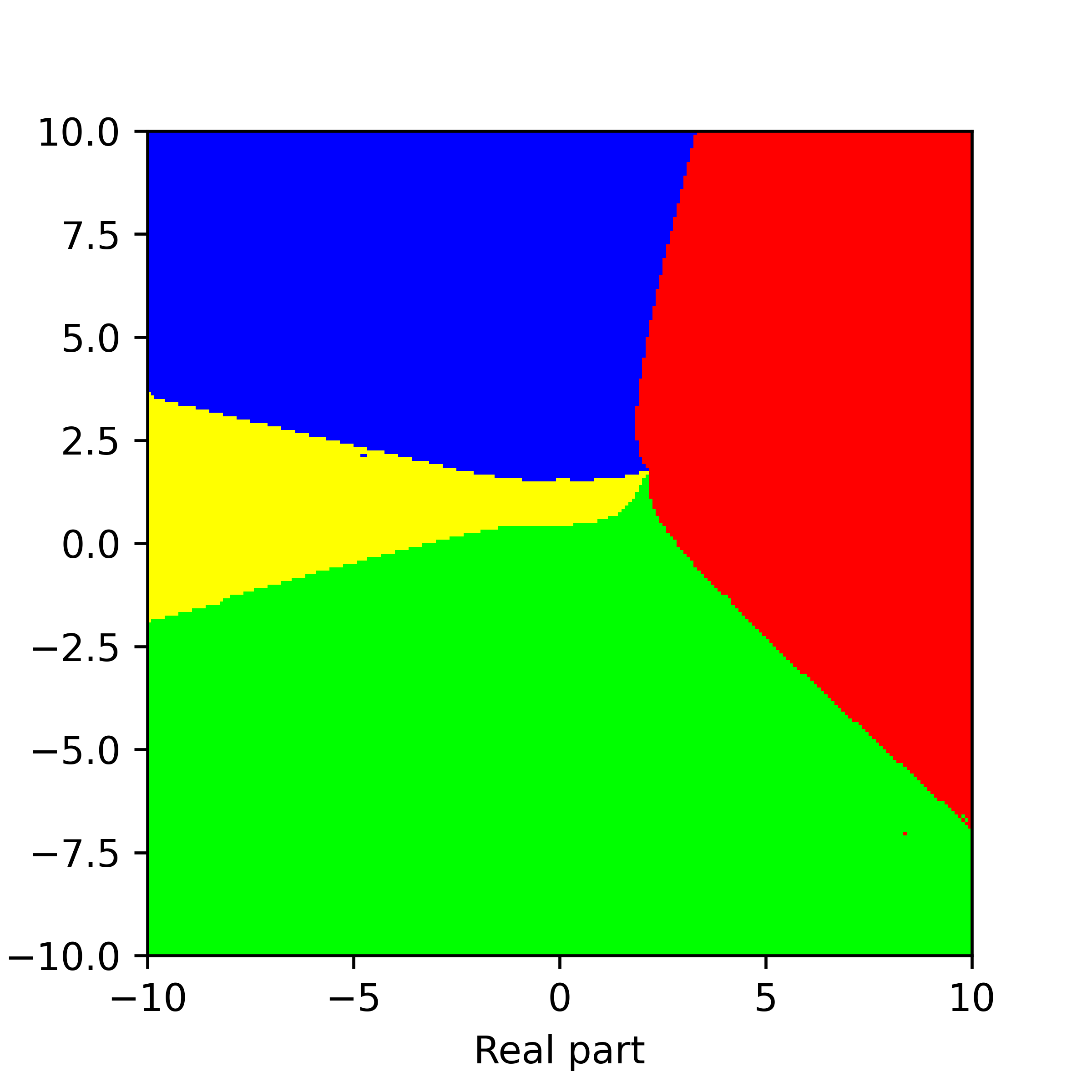}
    \includegraphics[width=3cm]{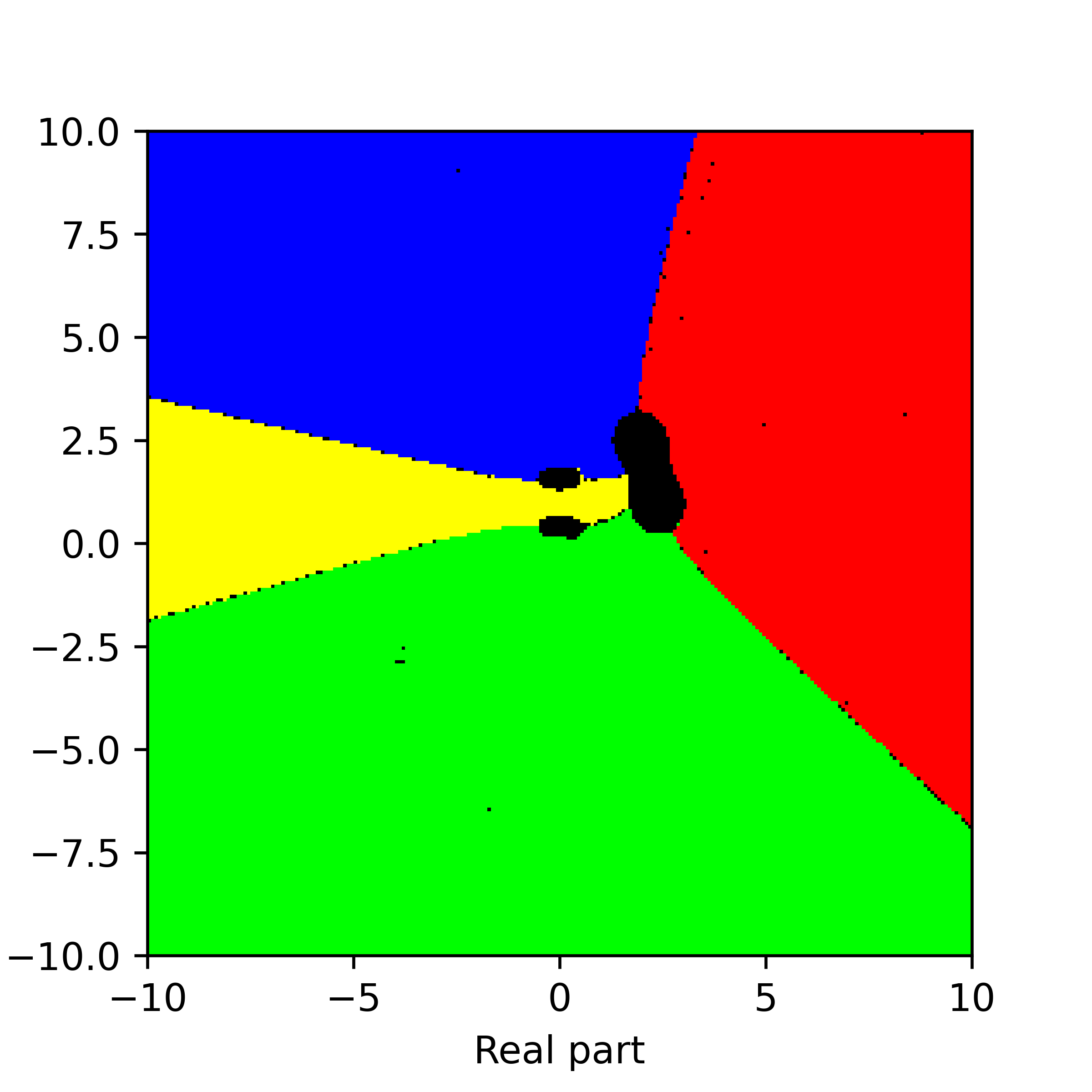}
    
    \caption{Basins of attraction for finding roots of the function $f_8$ by different methods. Pictures are referenced to from top to bottom, from left to right. Row 1: left picture is Voronoi's diagram, central picture is for Newton's method, right picture is for Random Relaxed Newton's method. Row 2: left picture is for Newton's method vOptimization, right picture is for BNQN. Row 3: left picture is for Newton's flow, central picture is for Newton's flow vFraction, right picture is for Newton's flow vOptimization. The black points in some of these pictures are those in the basin of attraction of critical points of $f_8$.}
    \label{fig:f8}
\end{figure}

\begin{figure}
    \centering
    \includegraphics[width=5cm]{Voronoif3Reduced.png}
    \includegraphics[width=3cm]{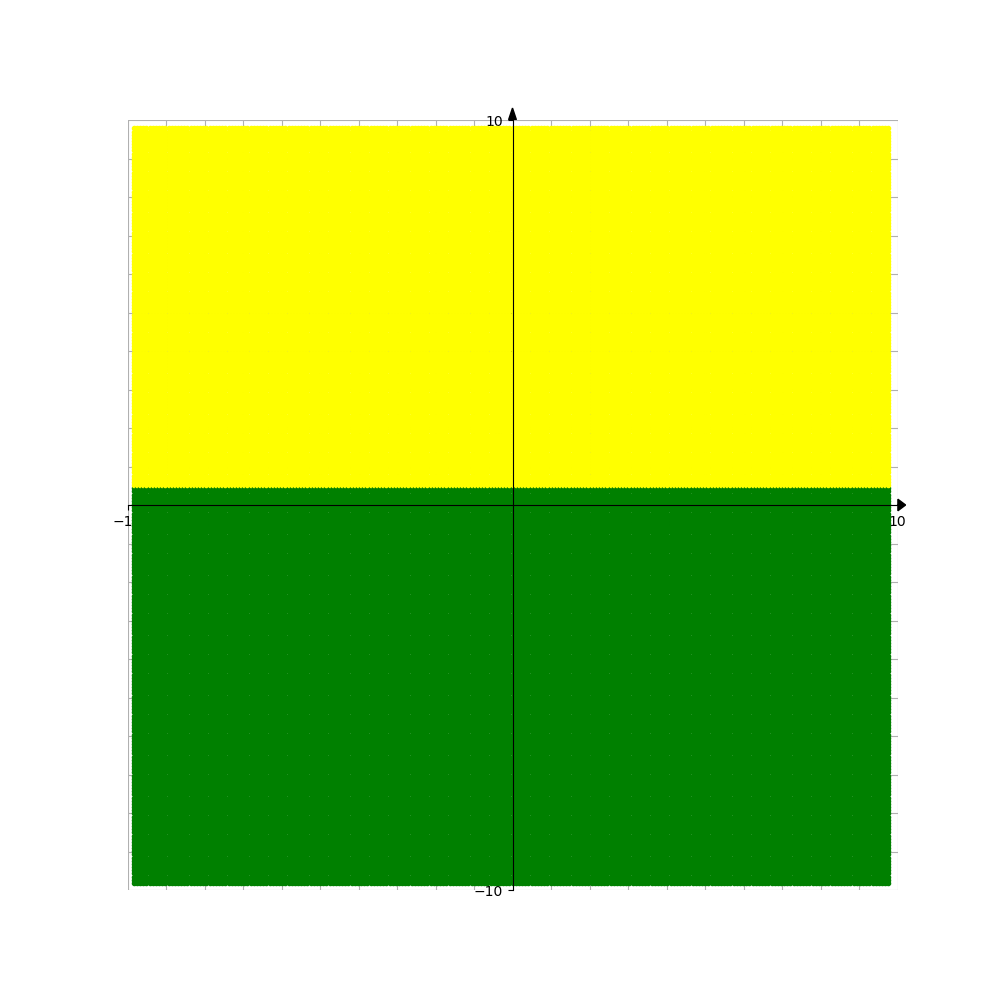}
    \includegraphics[width=3cm]{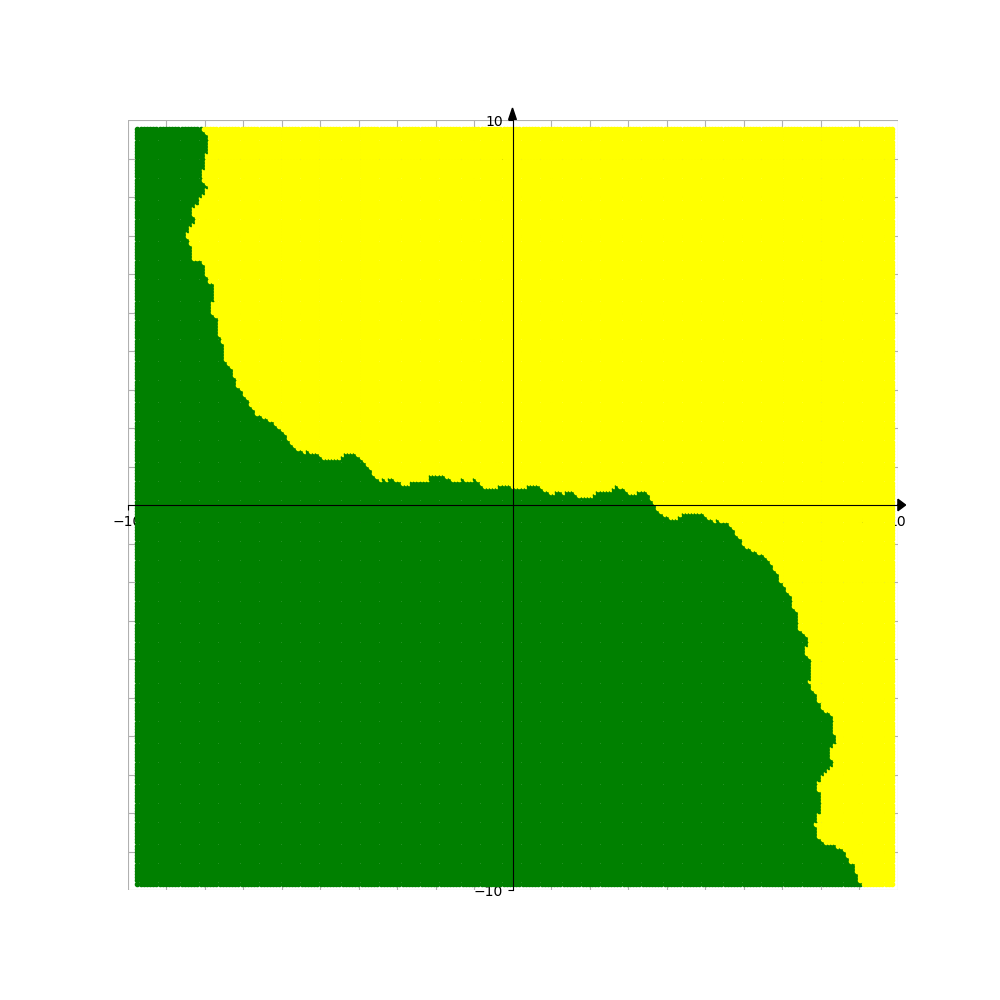}

    \bigskip
    \includegraphics[width=5.5cm]{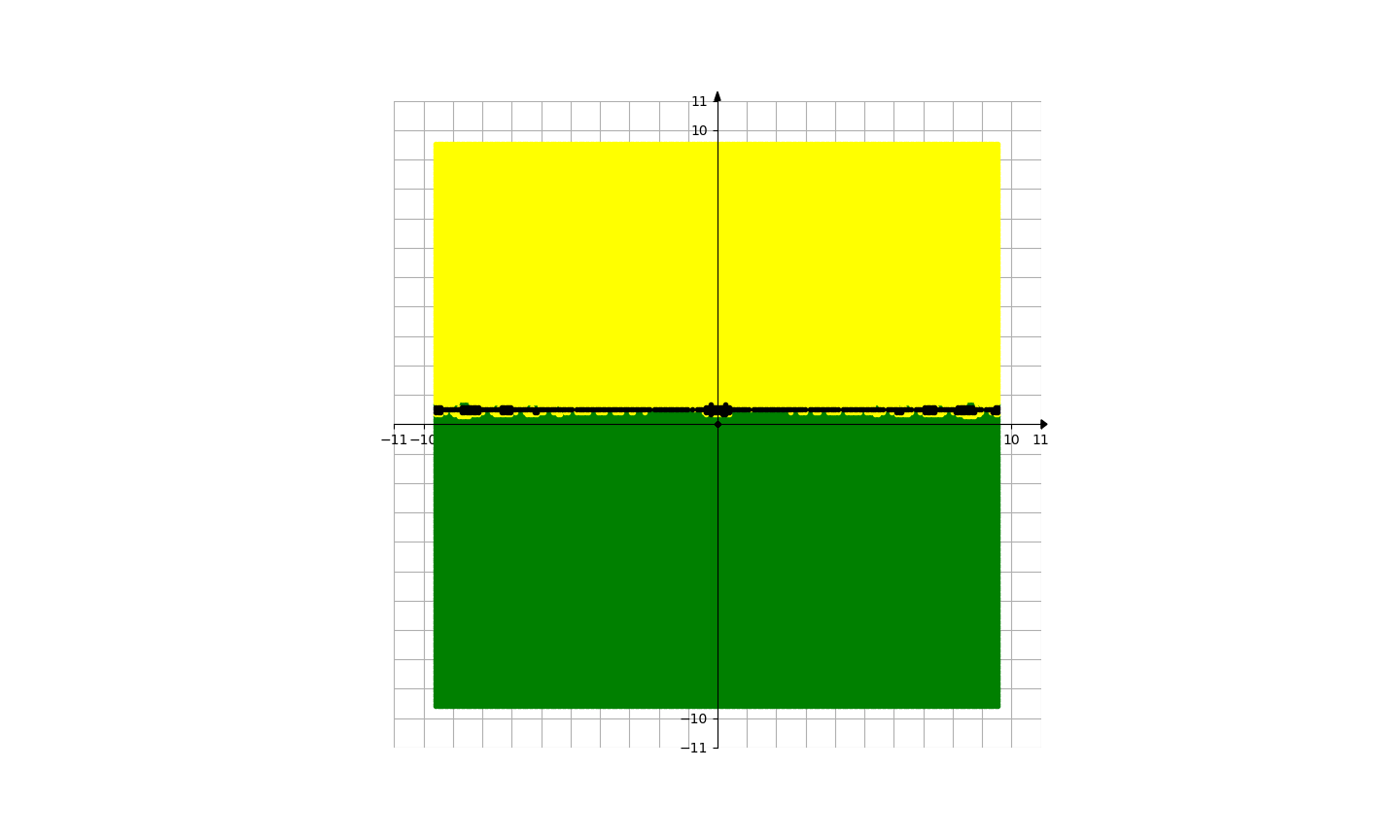}
    \includegraphics[width=3cm]{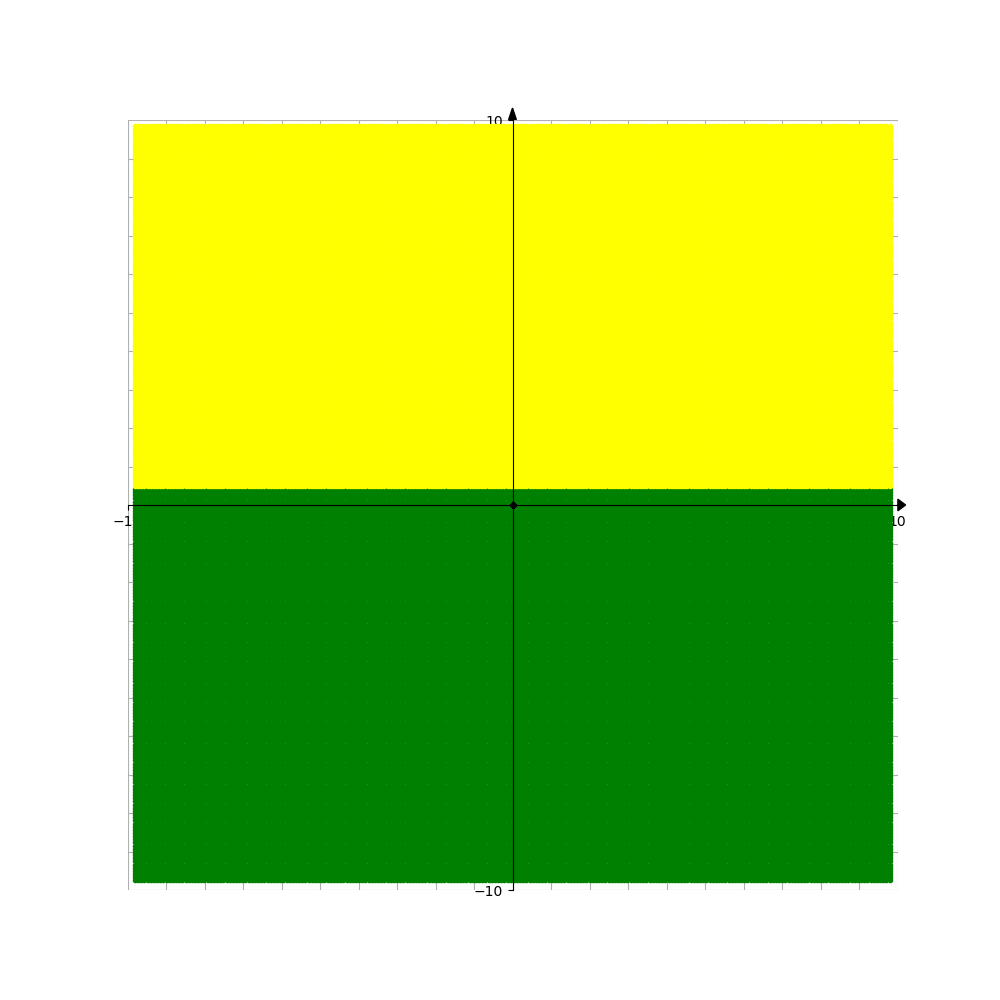}
    
    \bigskip
    \includegraphics[width=3cm]{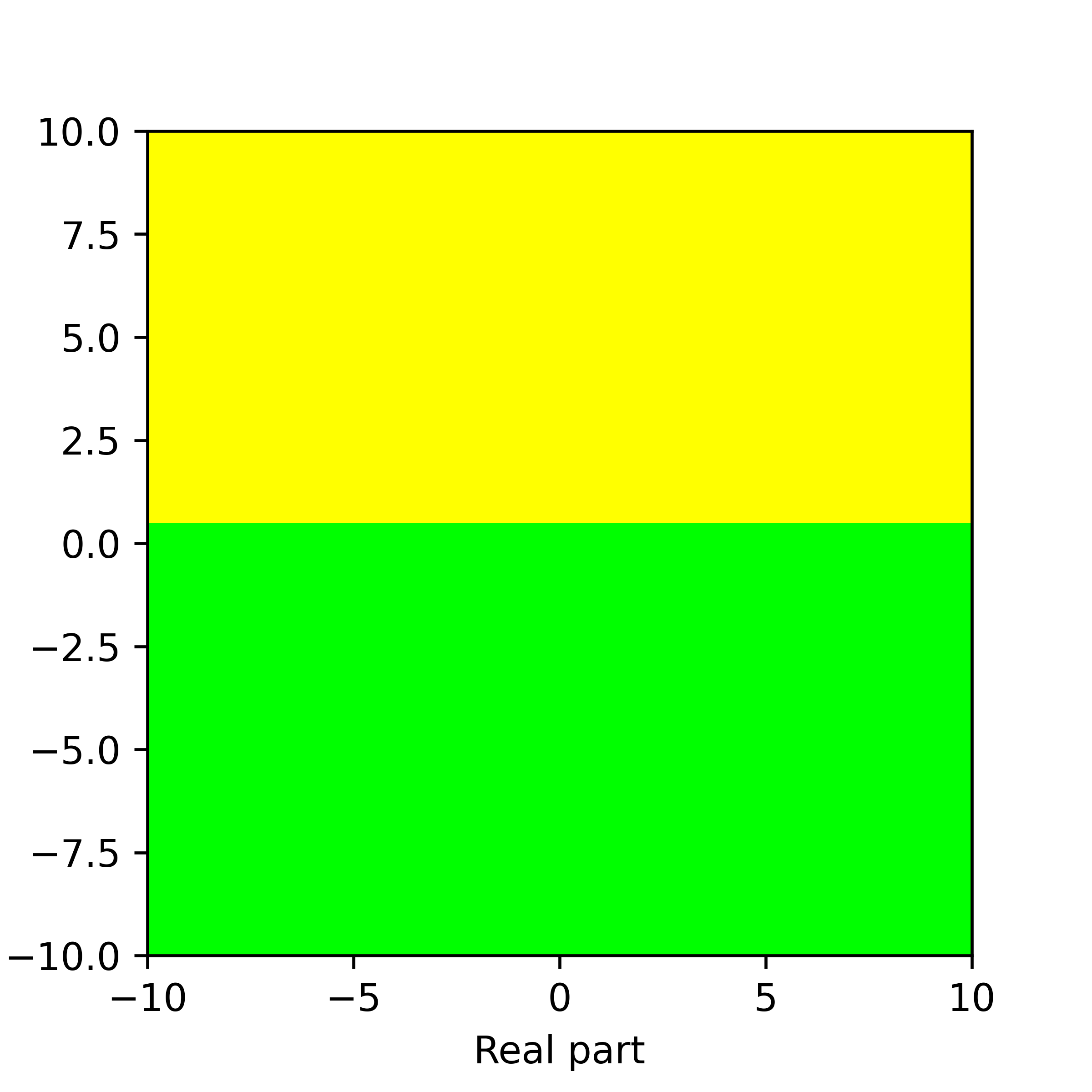}
    \includegraphics[width=3cm]{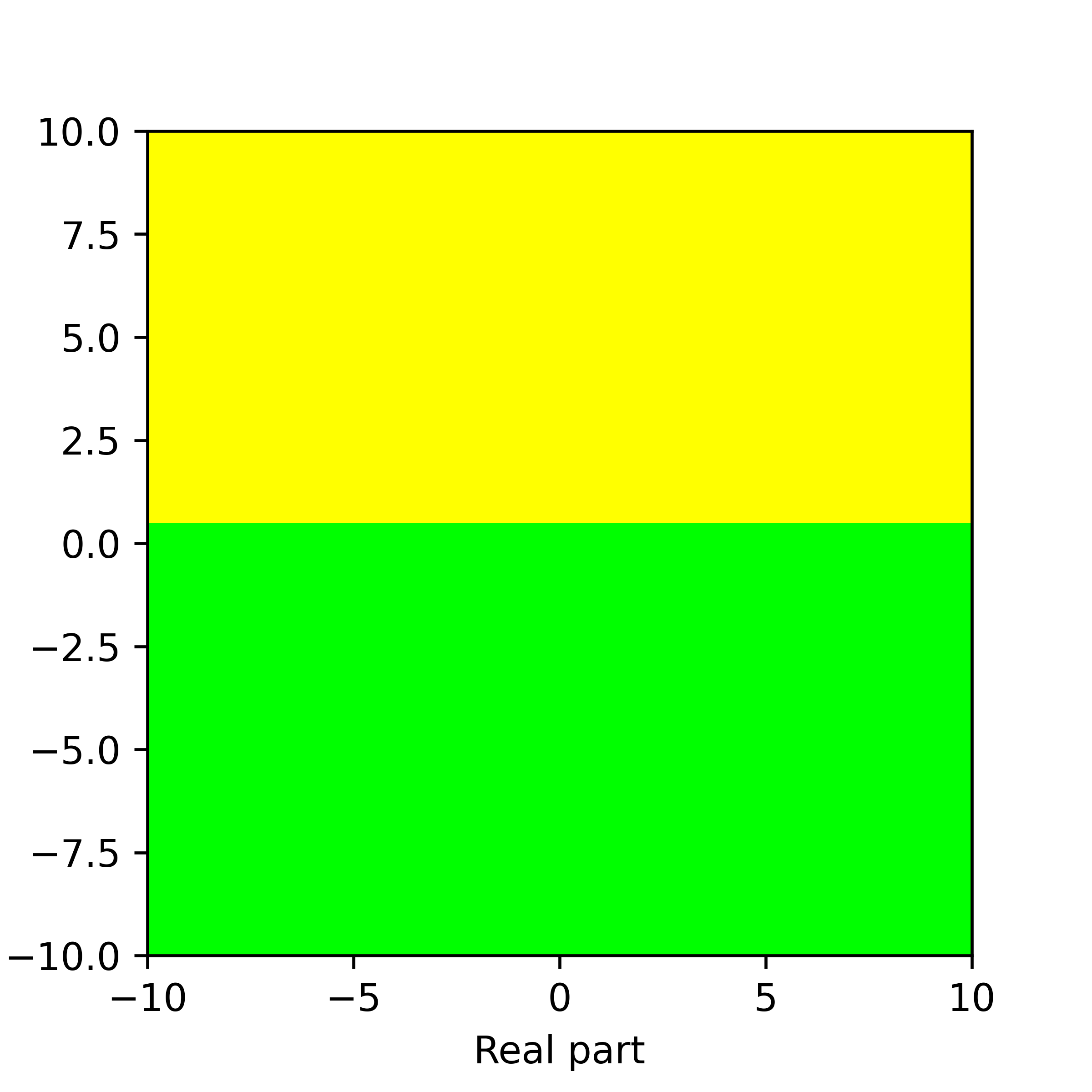}
    \includegraphics[width=3cm]{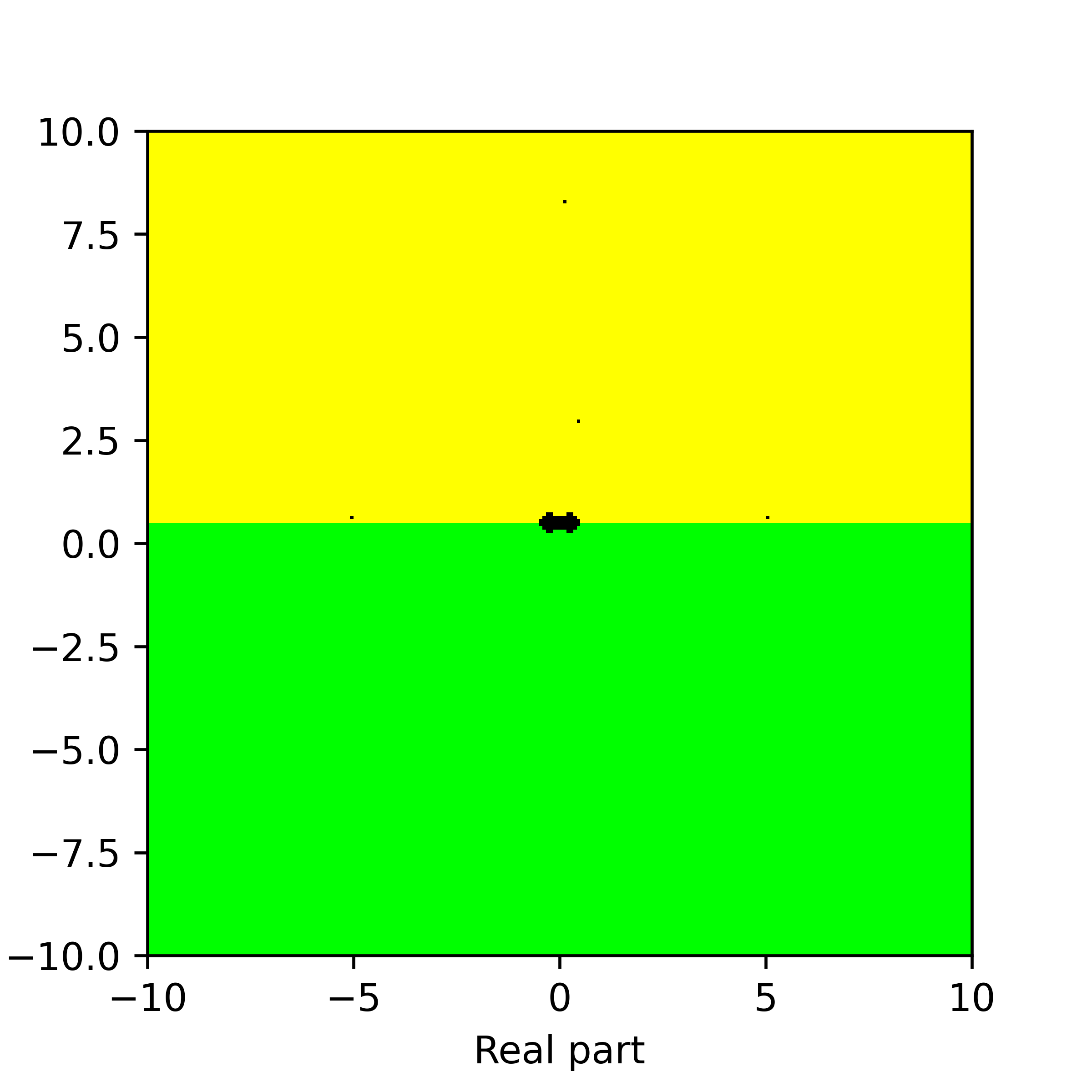}
    
    \caption{Basins of attraction for finding roots of the function $f_9$ by different methods. Pictures are referenced to from top to bottom, from left to right. Row 1: left picture is Voronoi's diagram, central picture is for Newton's method, right picture is for Random Relaxed Newton's method. Row 2: left picture is for Newton's method vOptimization, right picture is for BNQN. Row 3: left picture is for Newton's flow, central picture is for Newton's flow vFraction, right picture is for Newton's flow vOptimization. The black points in some of these pictures are those in the basin of attraction of critical points of $f_9$.}
    \label{fig:f9}
\end{figure}

\begin{figure}
    \centering
    \includegraphics[width=5cm]{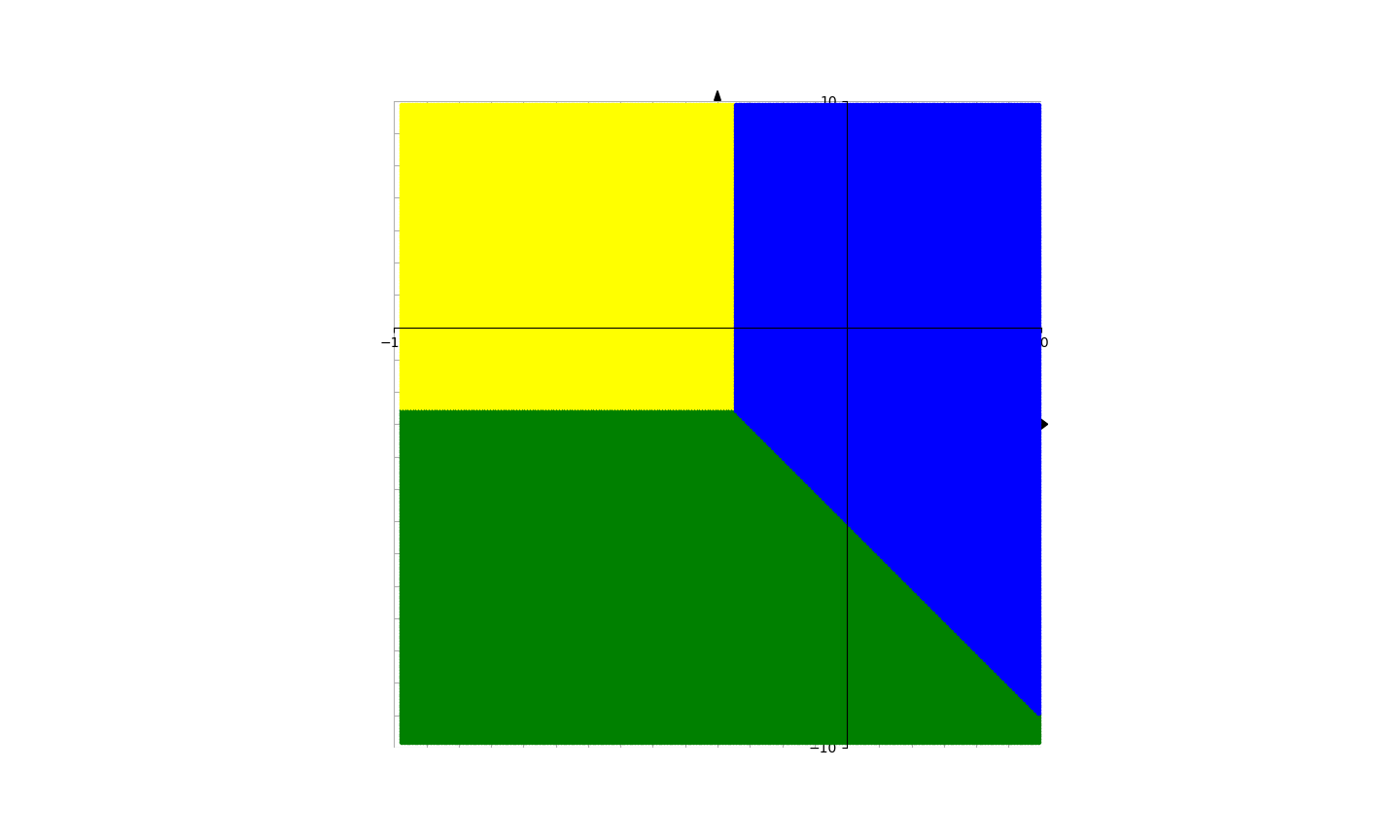}
    \includegraphics[width=3cm]{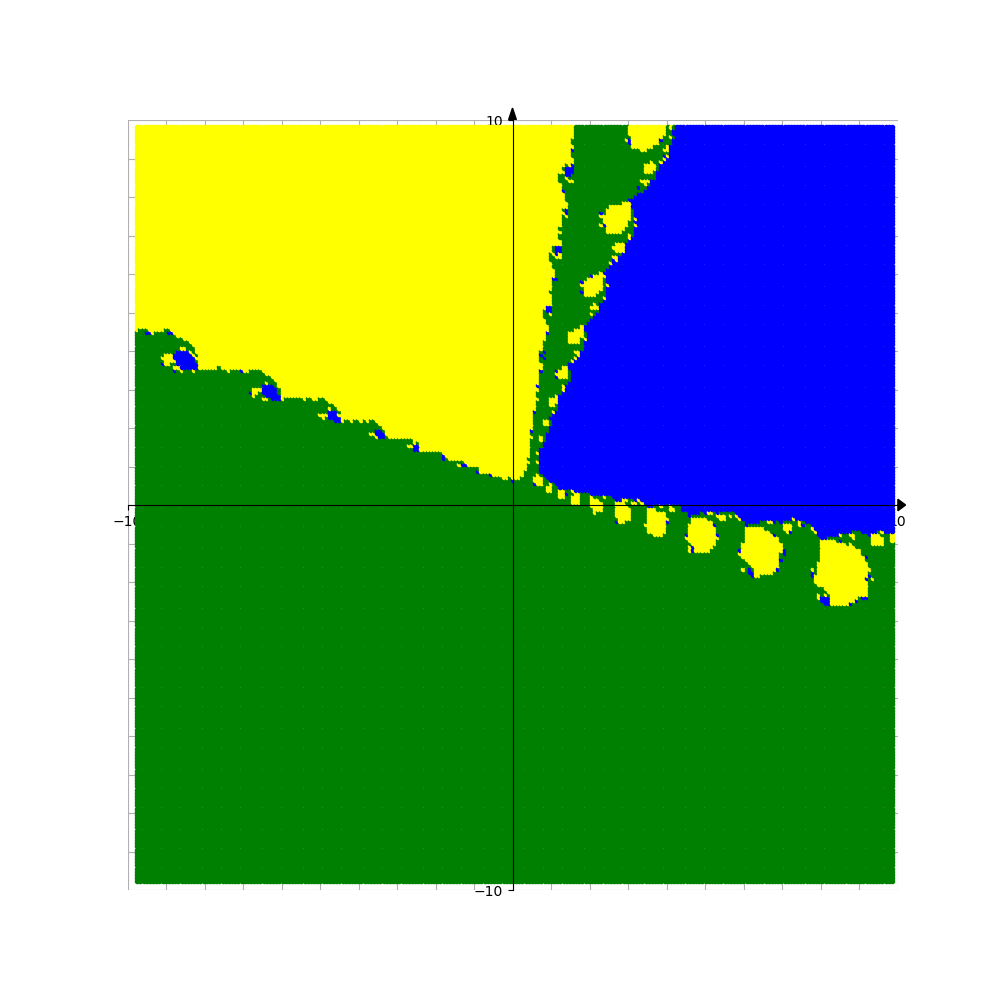}
    \includegraphics[width=3cm]{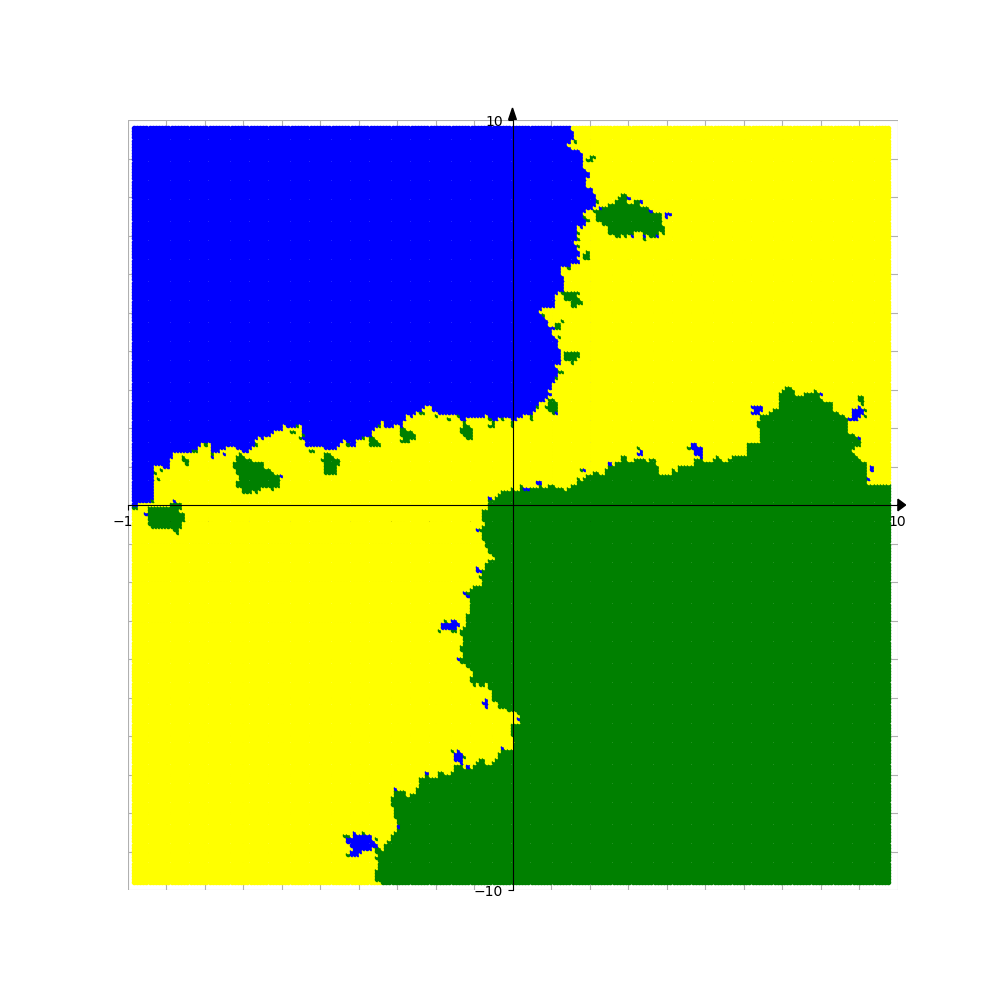}

    \bigskip
    \includegraphics[width=5.5cm]{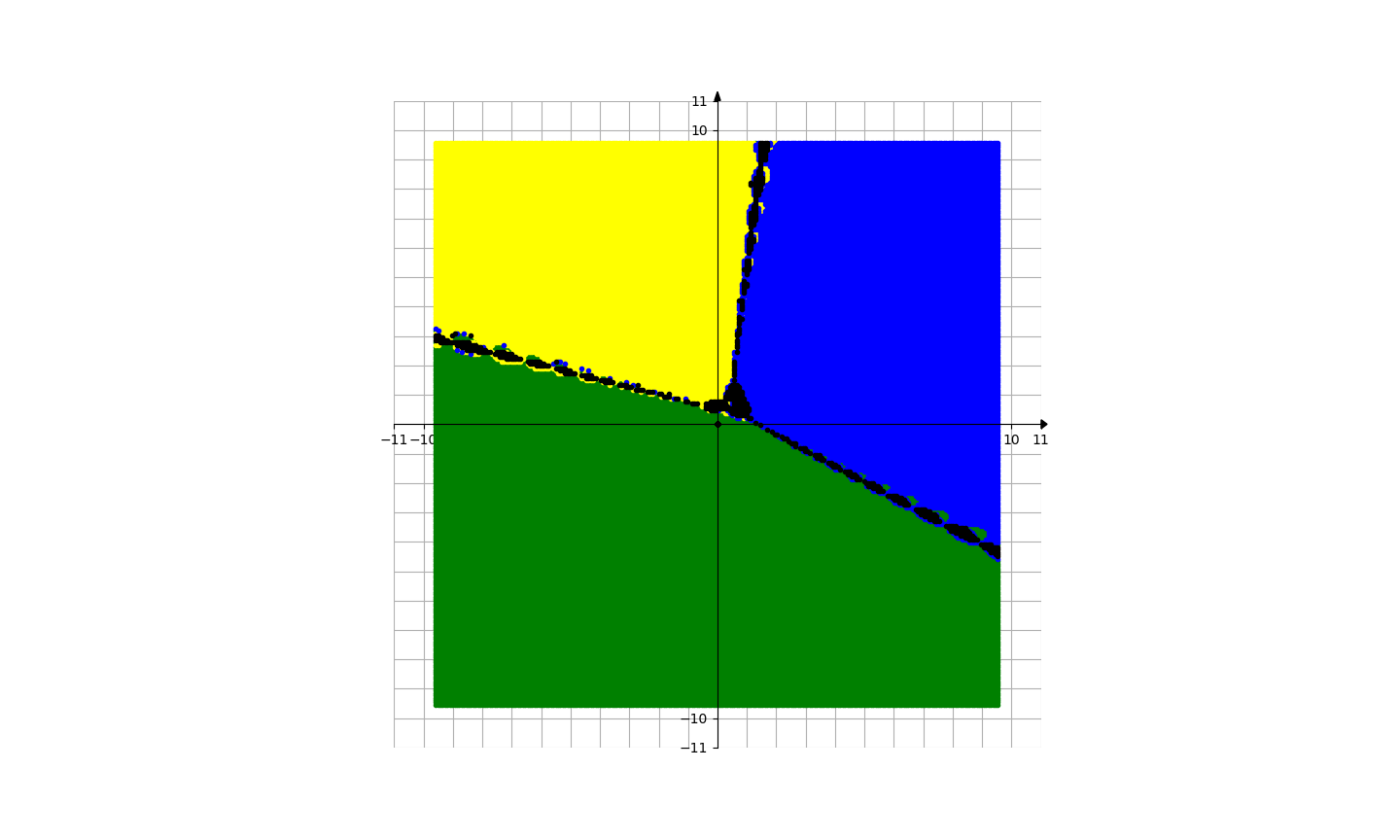}
    \includegraphics[width=3cm]{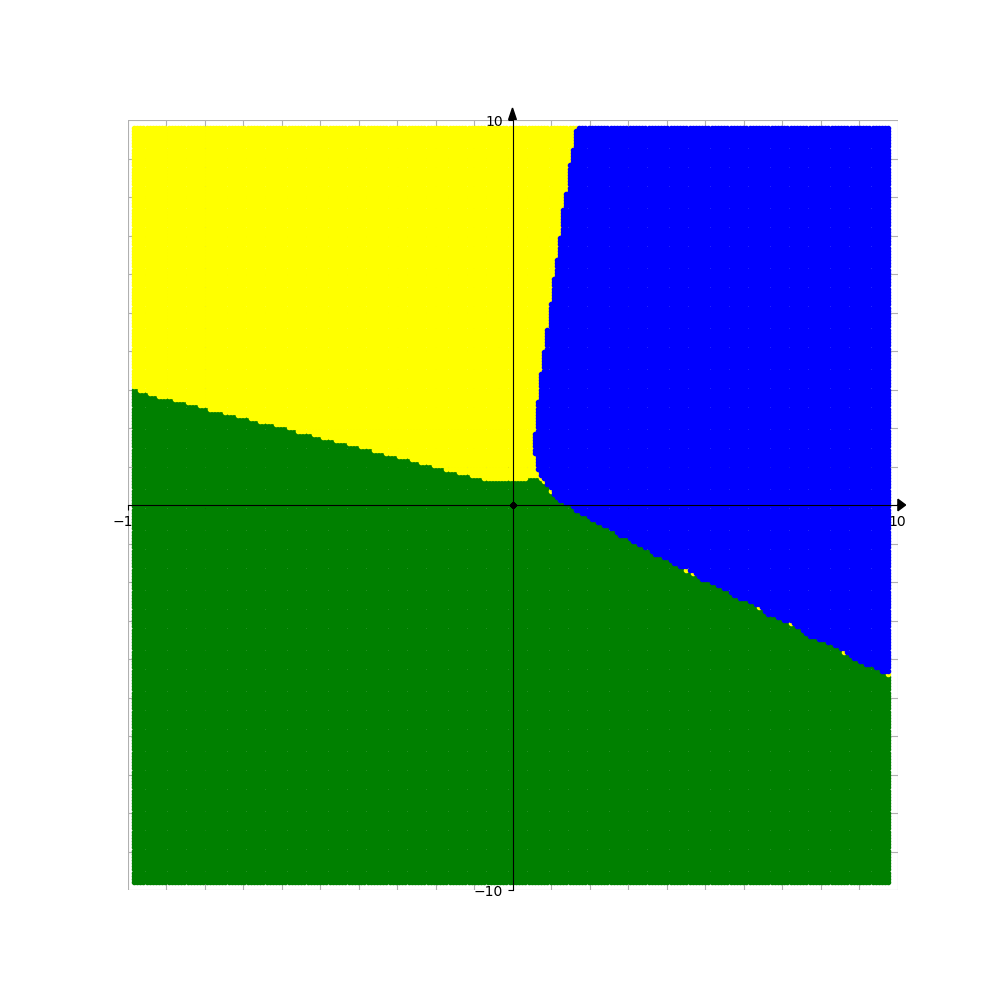}
    
    \bigskip
    \includegraphics[width=3cm]{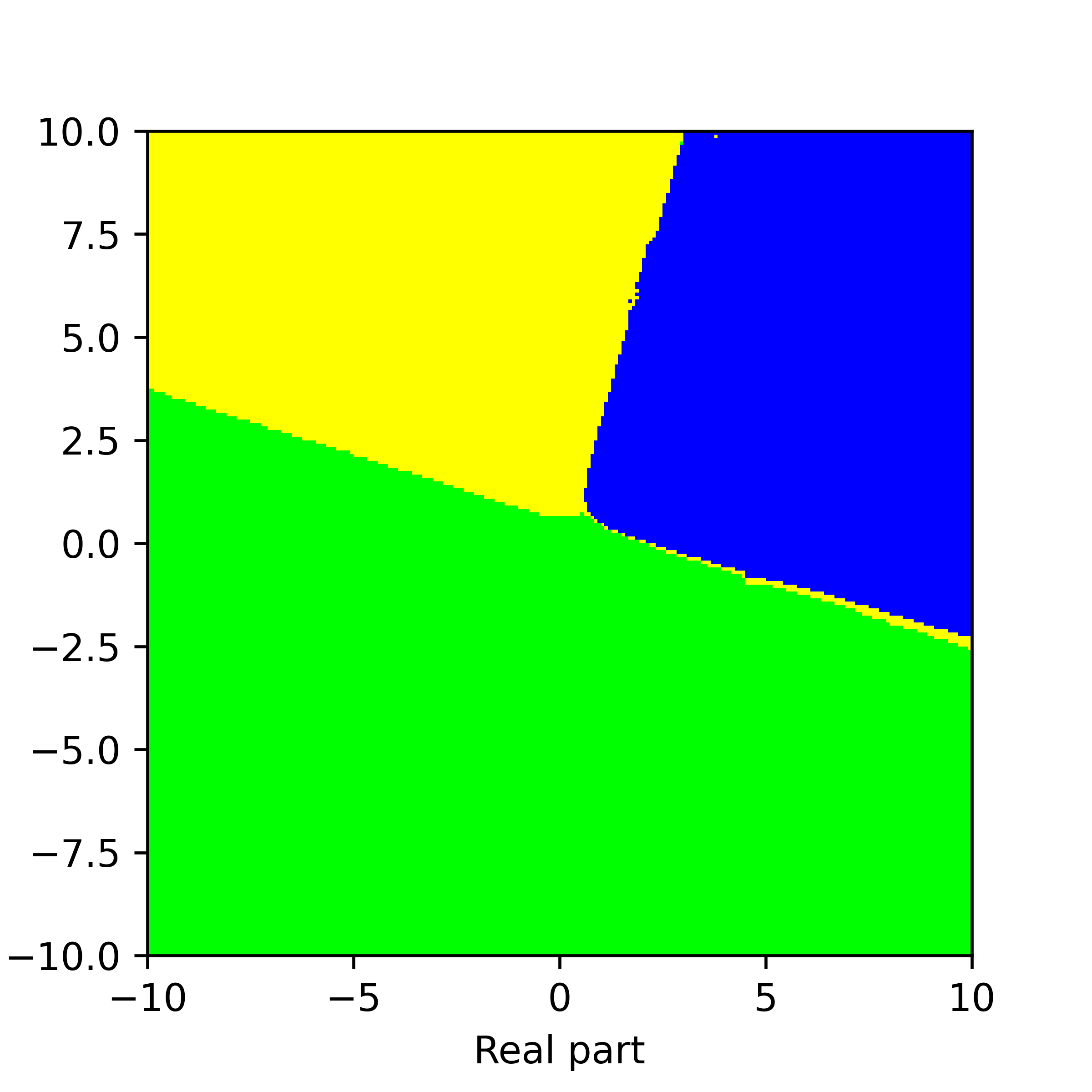}
    \includegraphics[width=3cm]{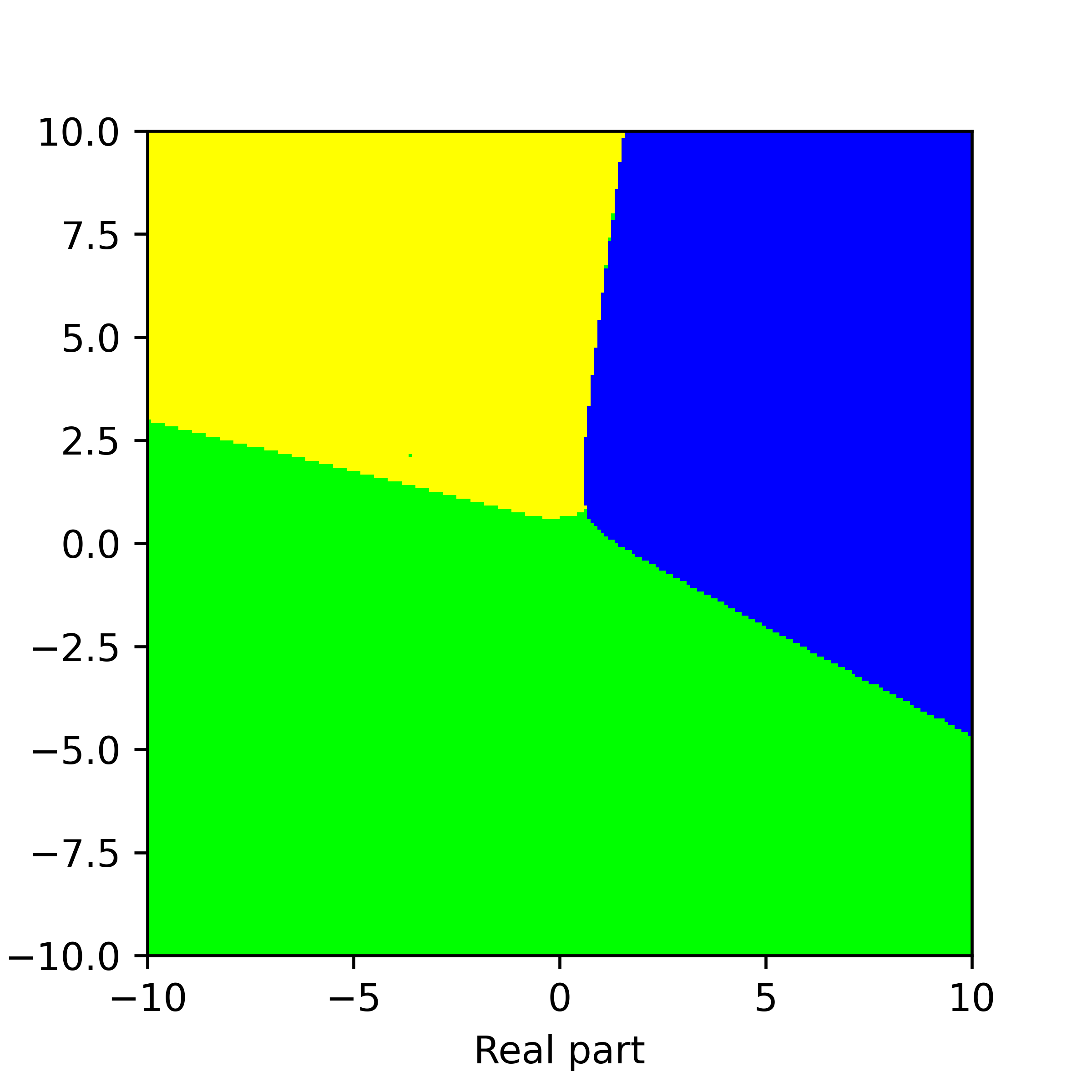}
    \includegraphics[width=3cm]{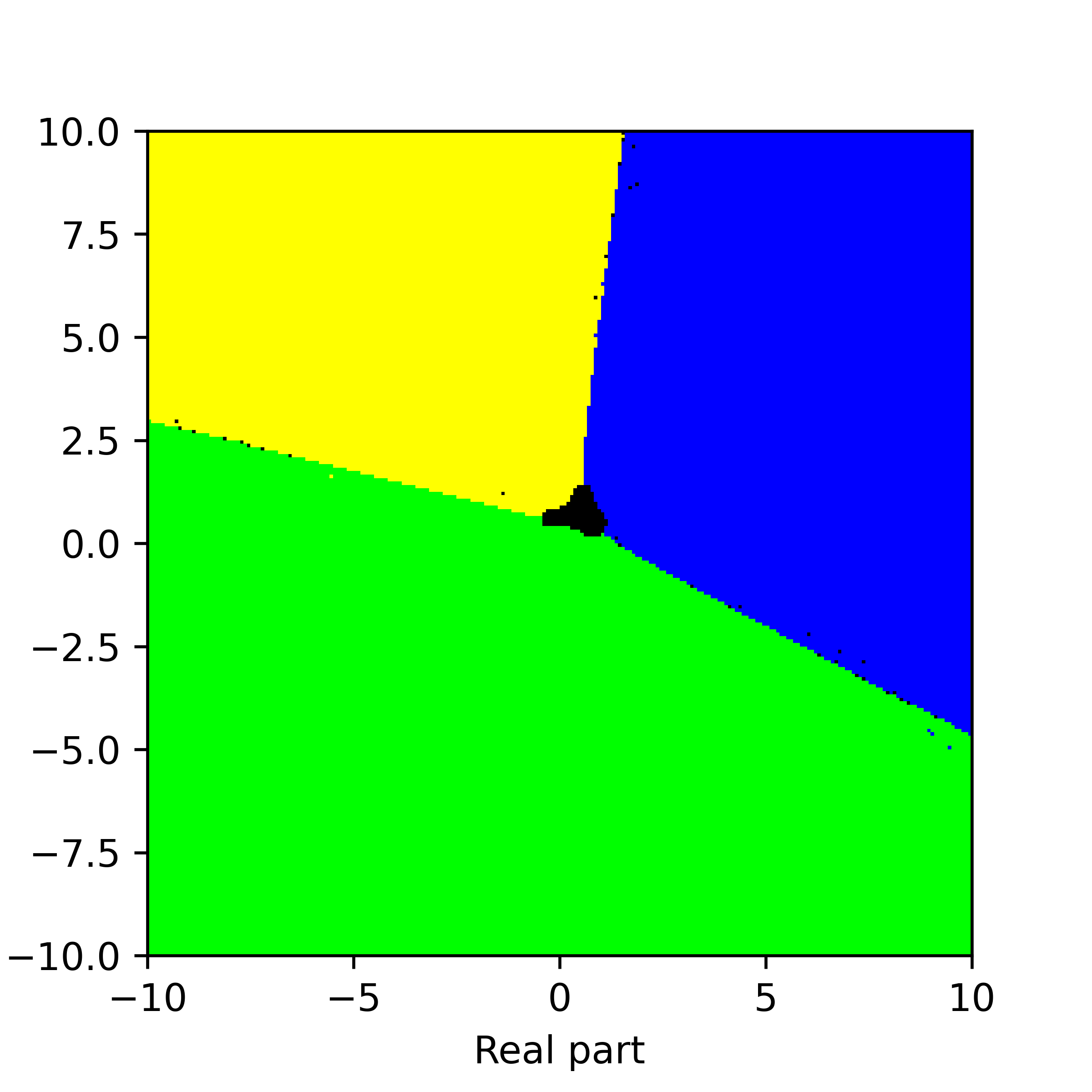}
    
    \caption{Basins of attraction for finding roots of the function $f_{10}$ by different methods. Pictures are referenced to from top to bottom, from left to right. Row 1: left picture is Voronoi's diagram, central picture is for Newton's method, right picture is for Random Relaxed Newton's method. Row 2: left picture is for Newton's method vOptimization, right picture is for BNQN. Row 3: left picture is for Newton's flow, central picture is for Newton's flow vFraction, right picture is for Newton's flow vOptimization. The black points in some of these pictures are those in the basin of attraction of critical points of $f_{10}$.}
    \label{fig:f10}
\end{figure}

\begin{figure}
    \centering
    \includegraphics[width=5cm]{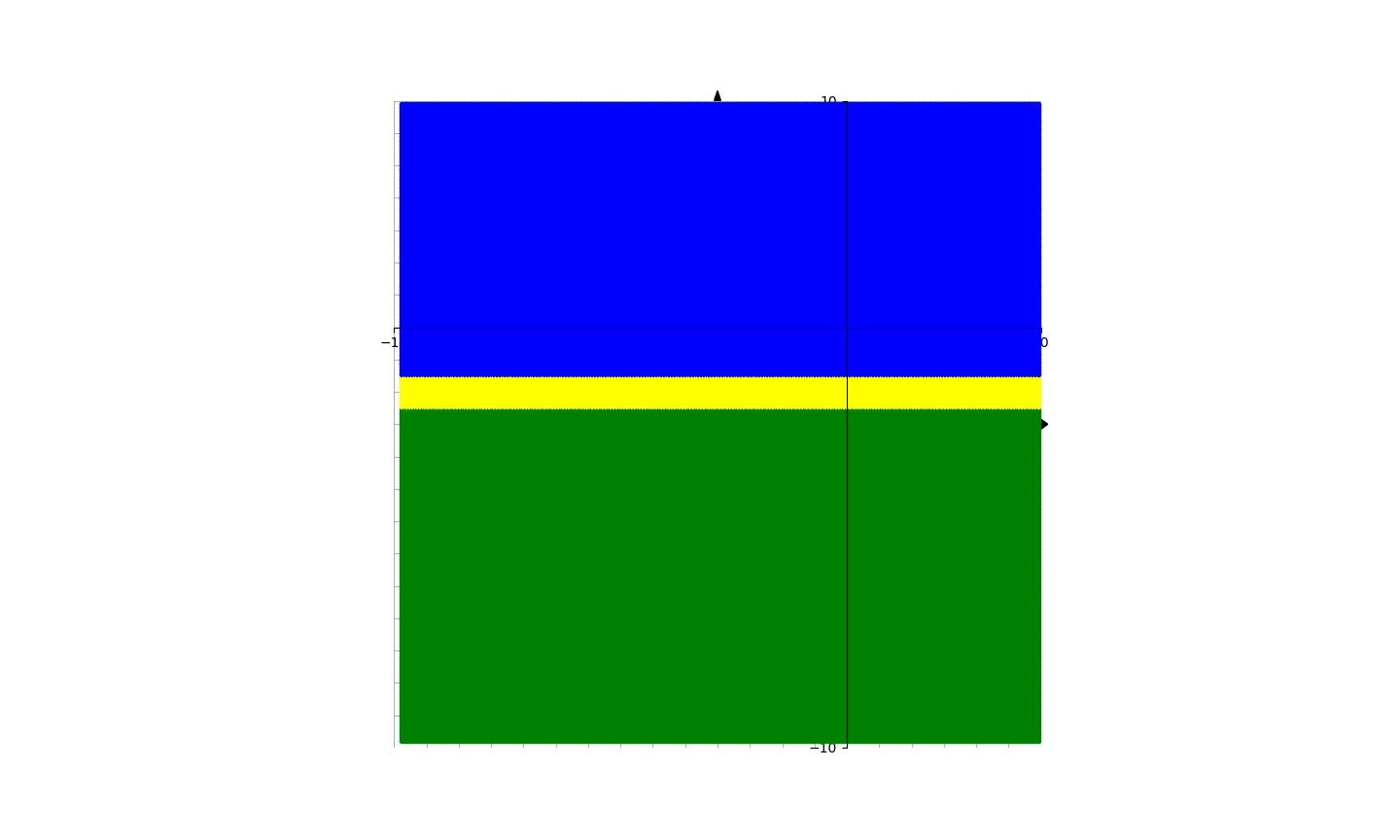}
    \includegraphics[width=3cm]{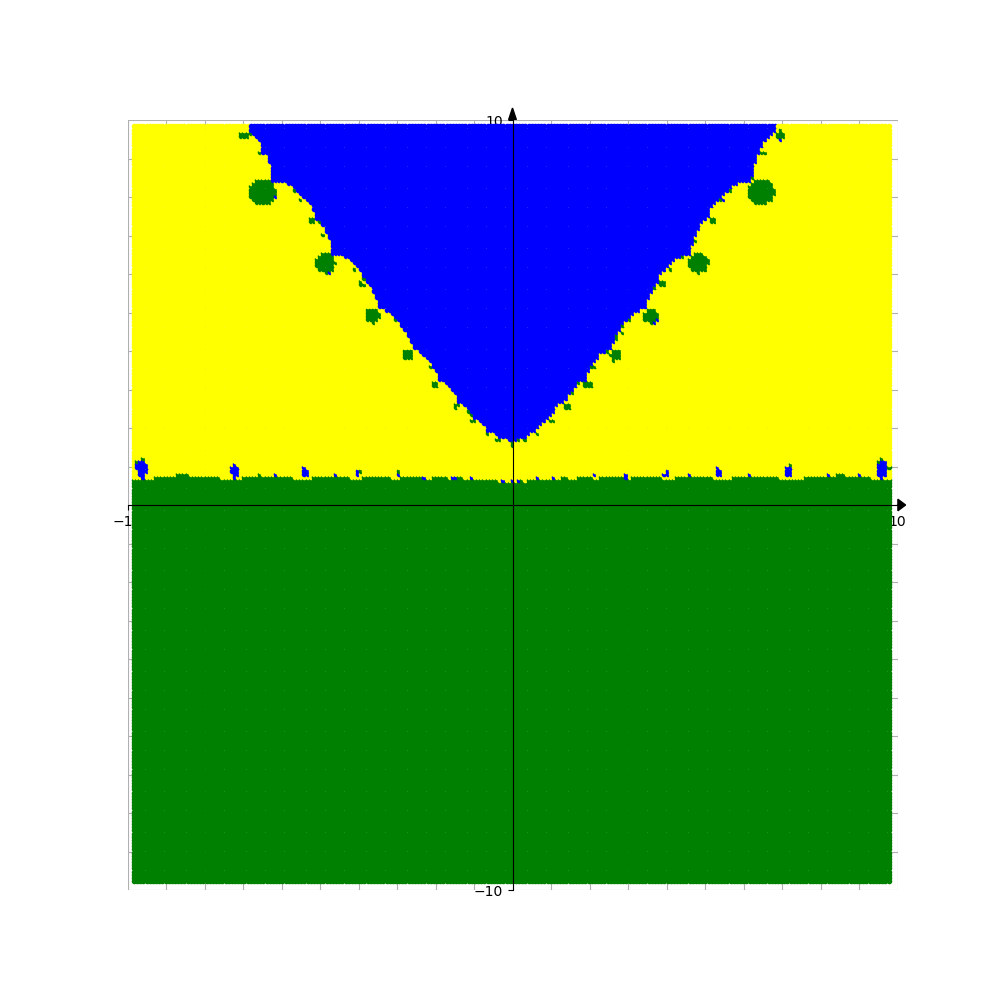}
    \includegraphics[width=3cm]{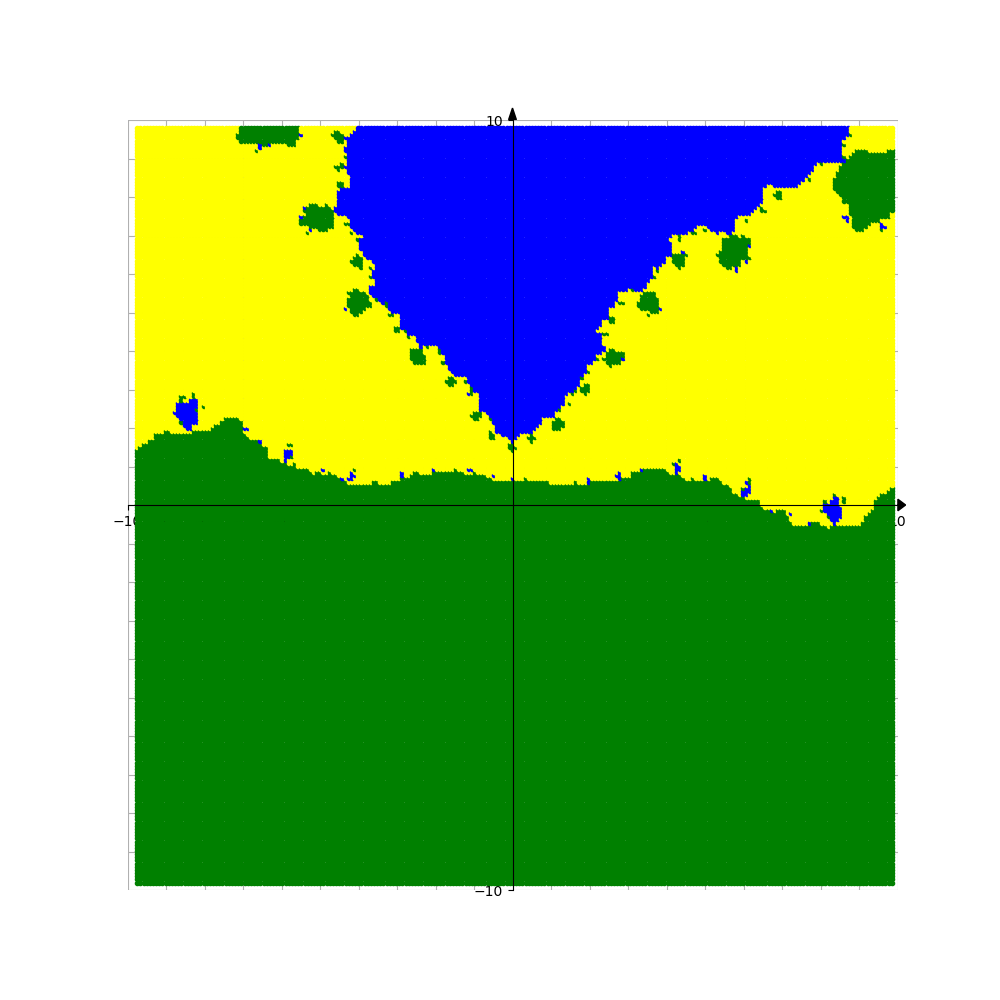}

    \bigskip
    \includegraphics[width=5.5cm]{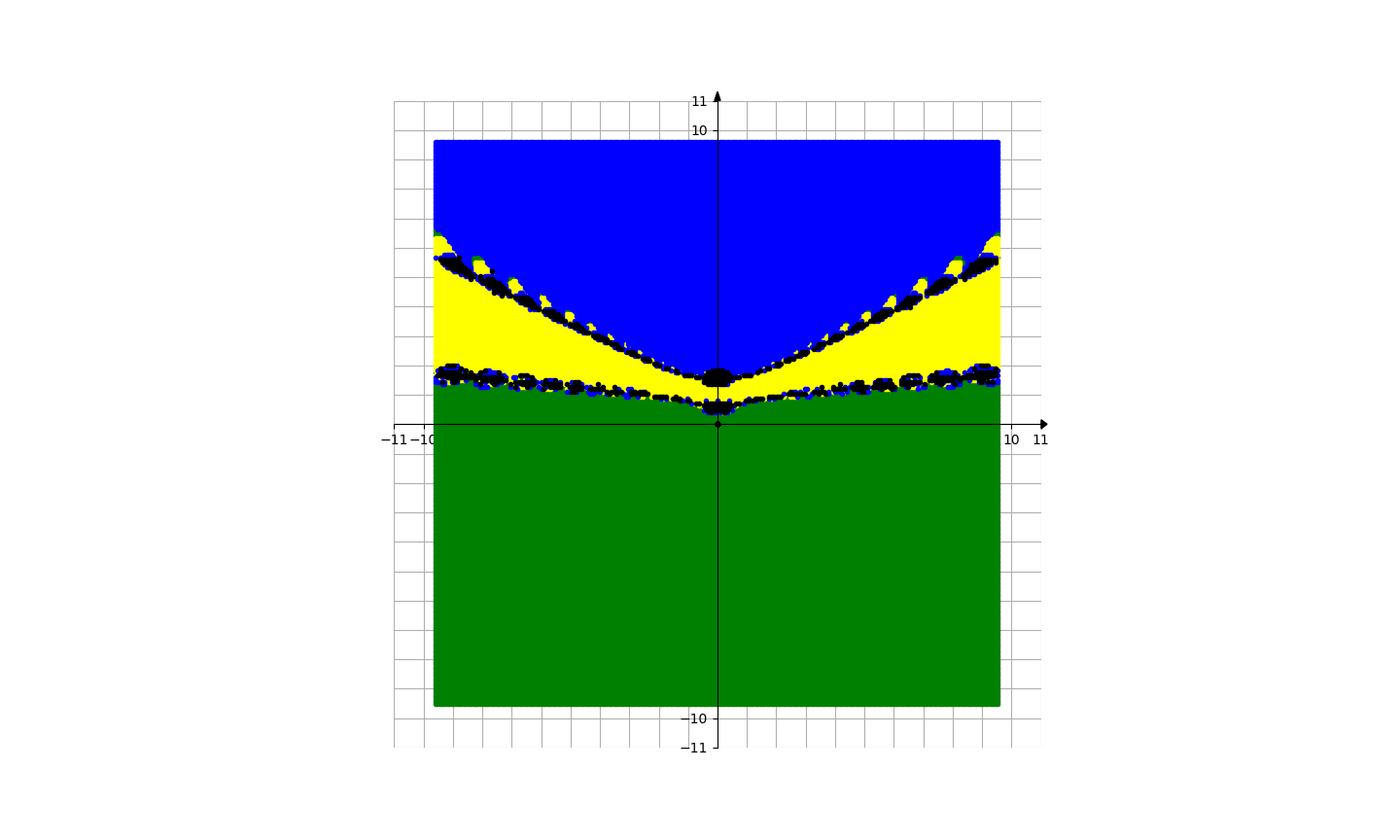}
    \includegraphics[width=3cm]{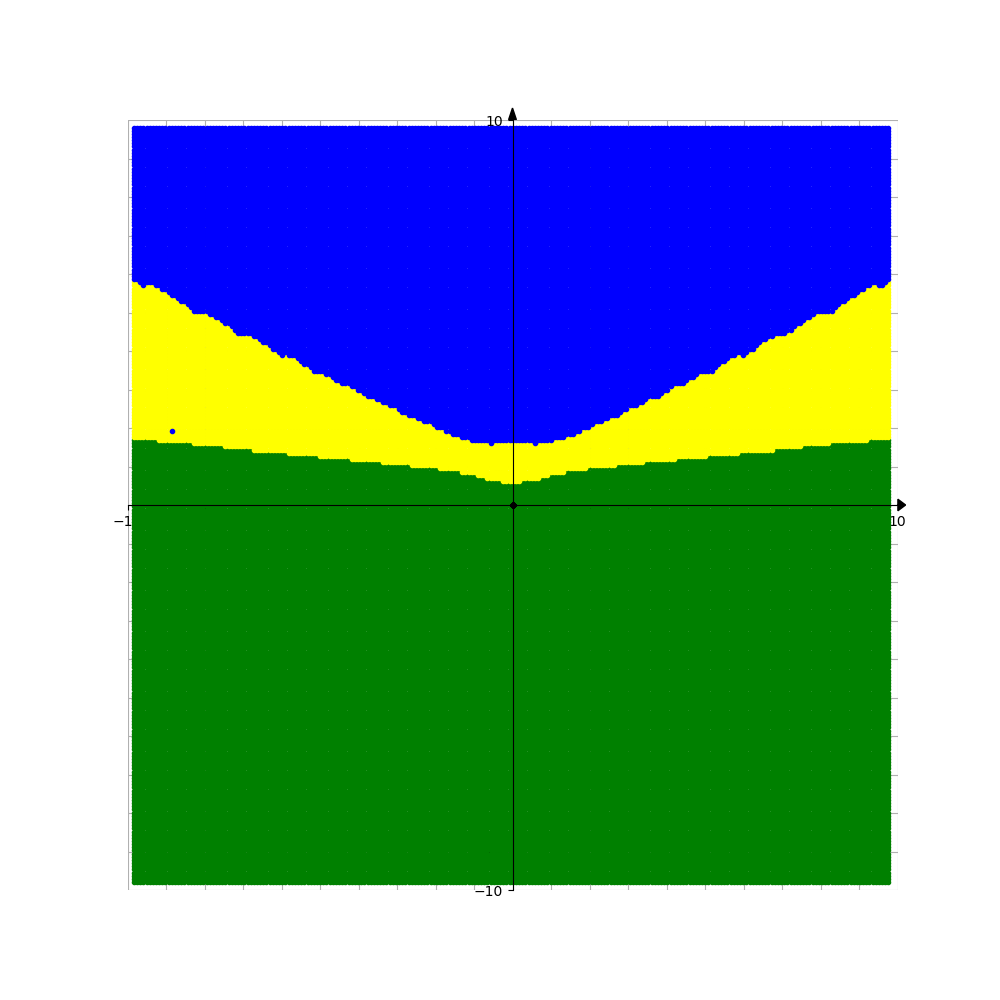}
    
    \bigskip
    \includegraphics[width=3cm]{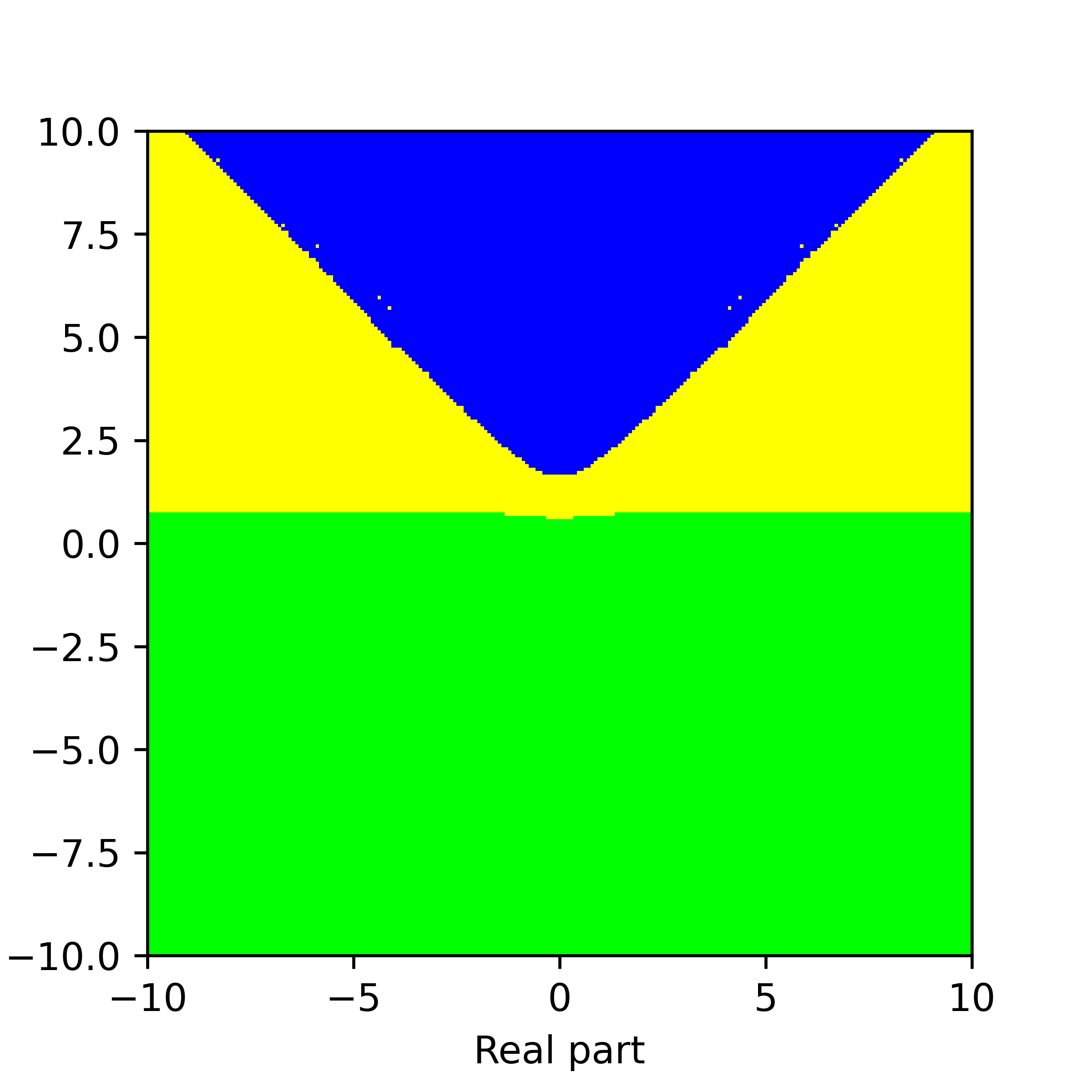}
    \includegraphics[width=3cm]{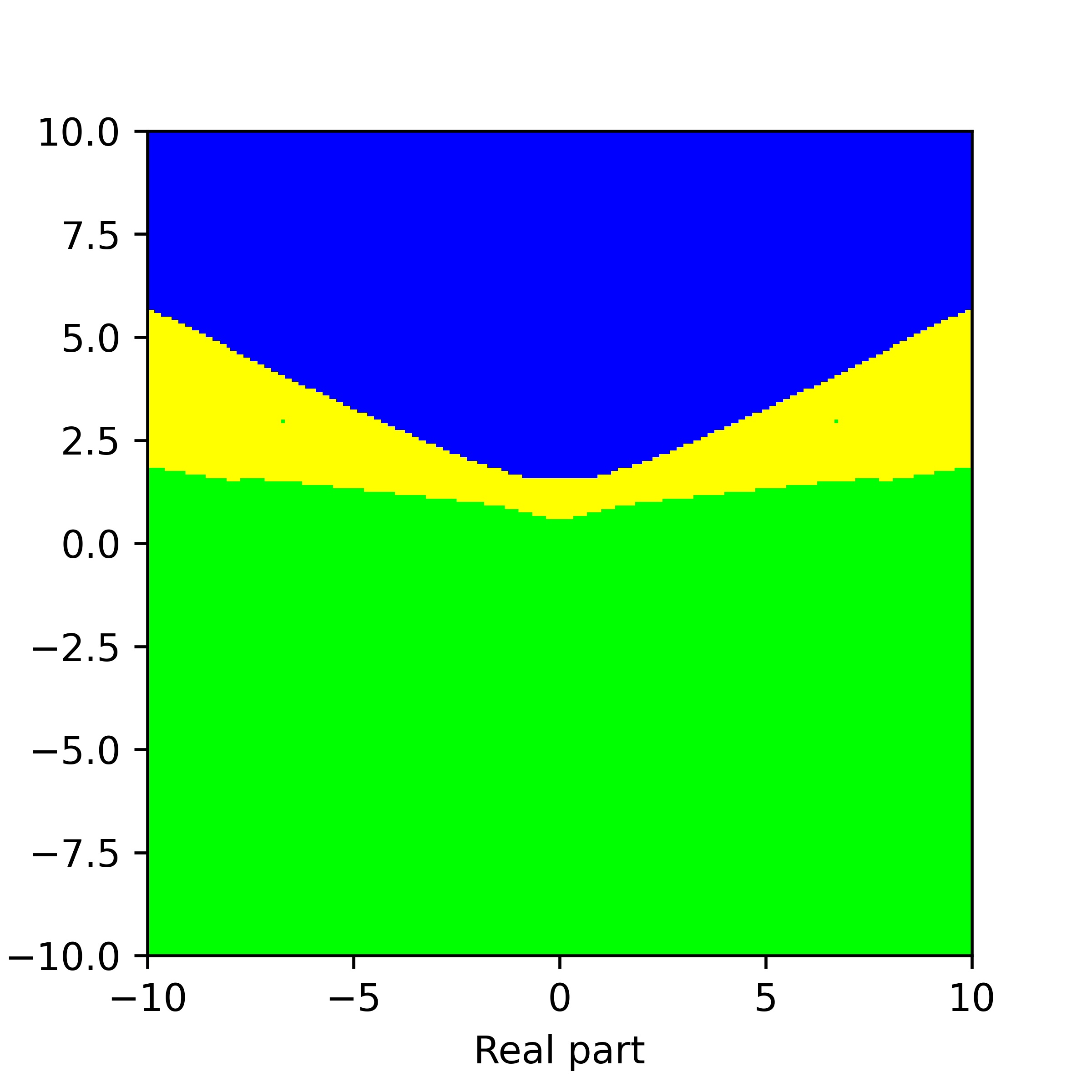}
    \includegraphics[width=3cm]{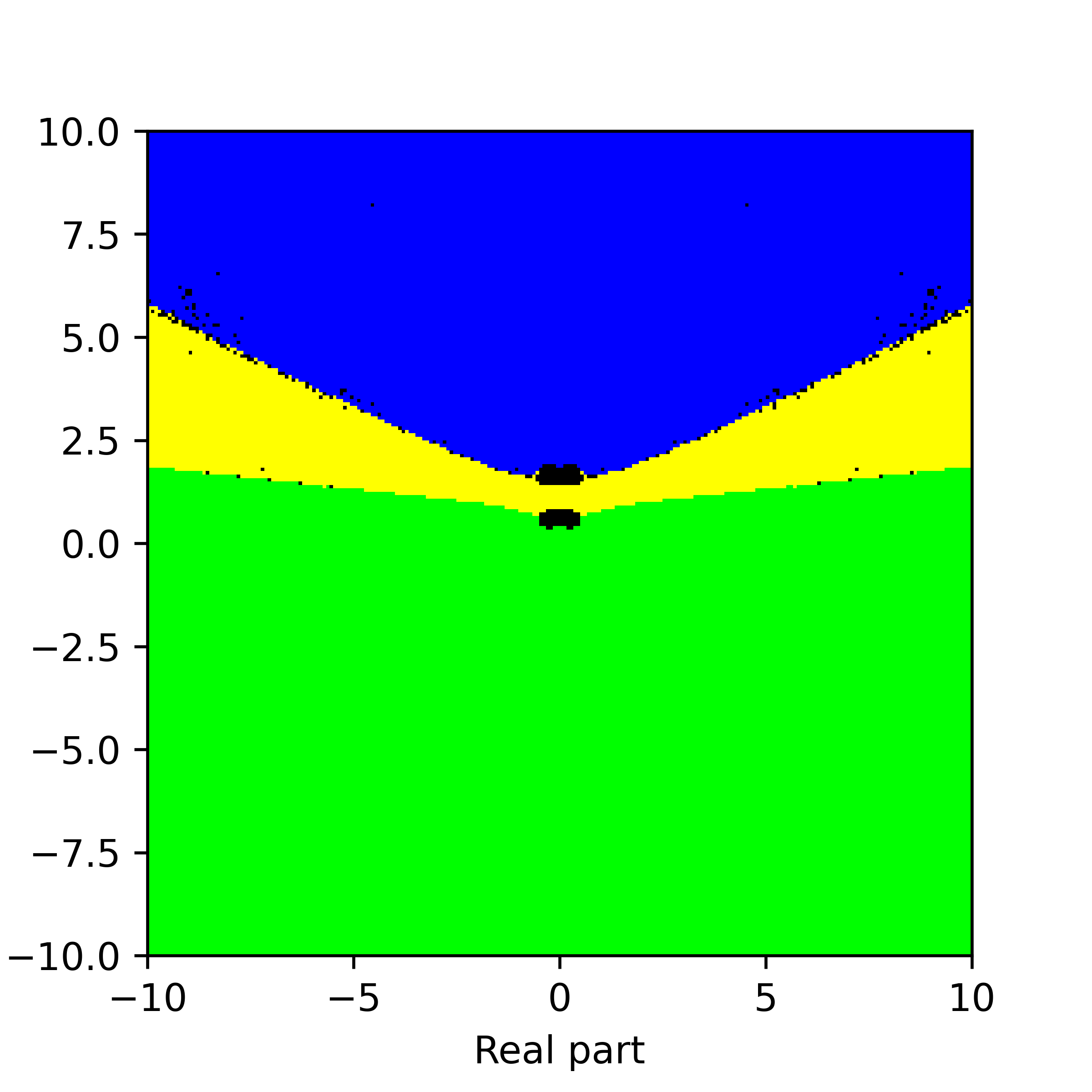}
    
    \caption{Basins of attraction for finding roots of the function $f_{11}$ by different methods. Pictures are referenced to from top to bottom, from left to right. Row 1: left picture is Voronoi's diagram, central picture is for Newton's method, right picture is for Random Relaxed Newton's method. Row 2: left picture is for Newton's method vOptimization, right picture is for BNQN. Row 3: left picture is for Newton's flow, central picture is for Newton's flow vFraction, right picture is for Newton's flow vOptimization. The black points in some of these pictures are those in the basin of attraction of critical points of $f_{11}$.}
    \label{fig:f11}
\end{figure}

\begin{figure}
    \centering
    \includegraphics[width=5cm]{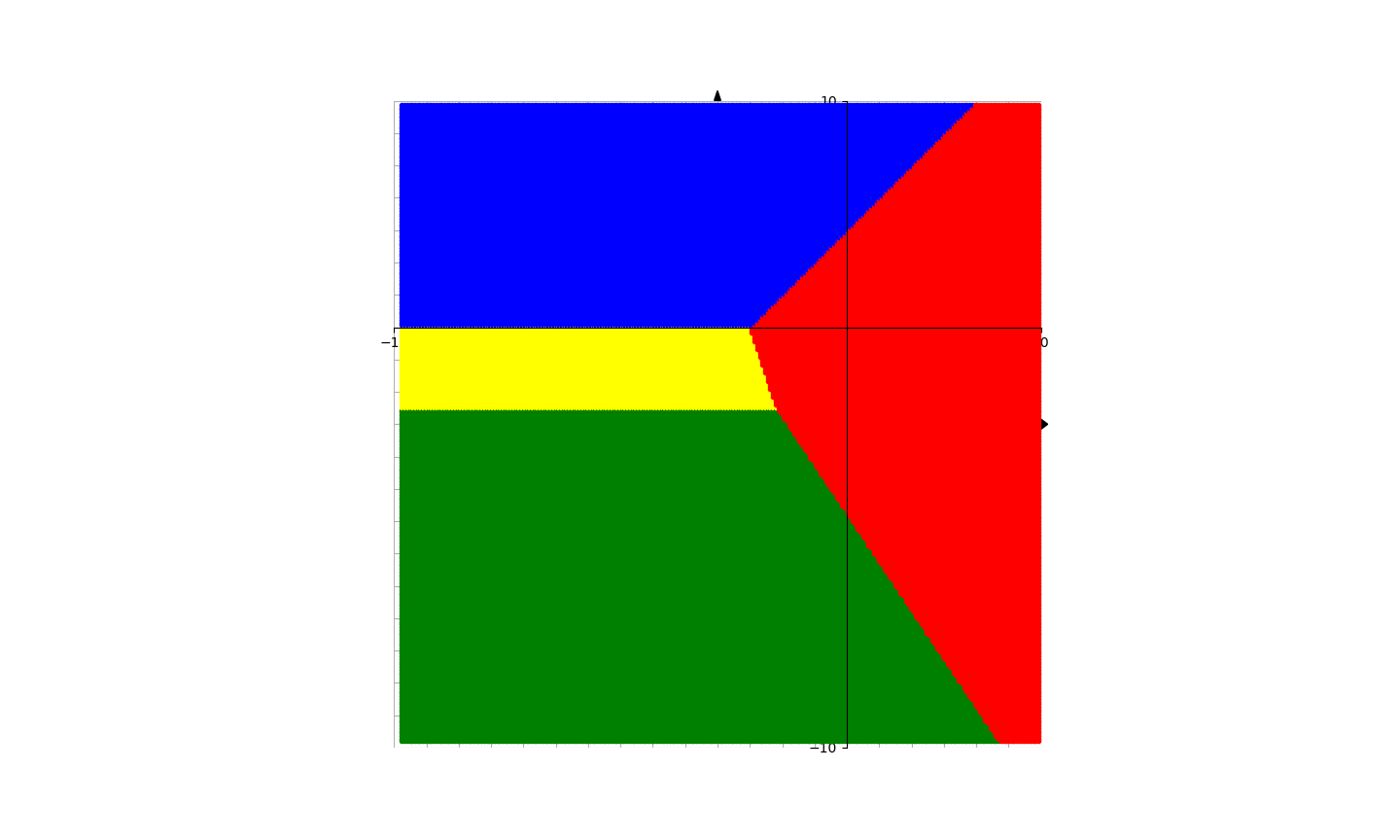}
    \includegraphics[width=3cm]{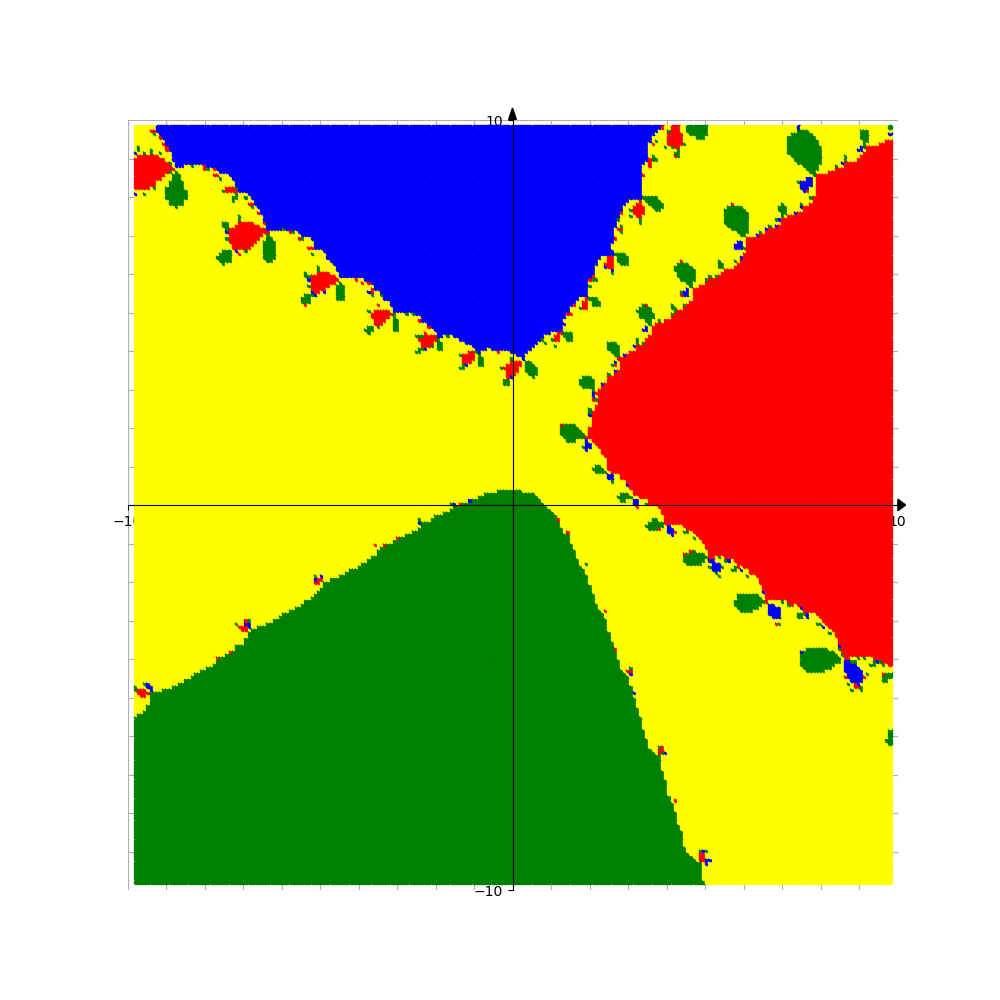}
    \includegraphics[width=3cm]{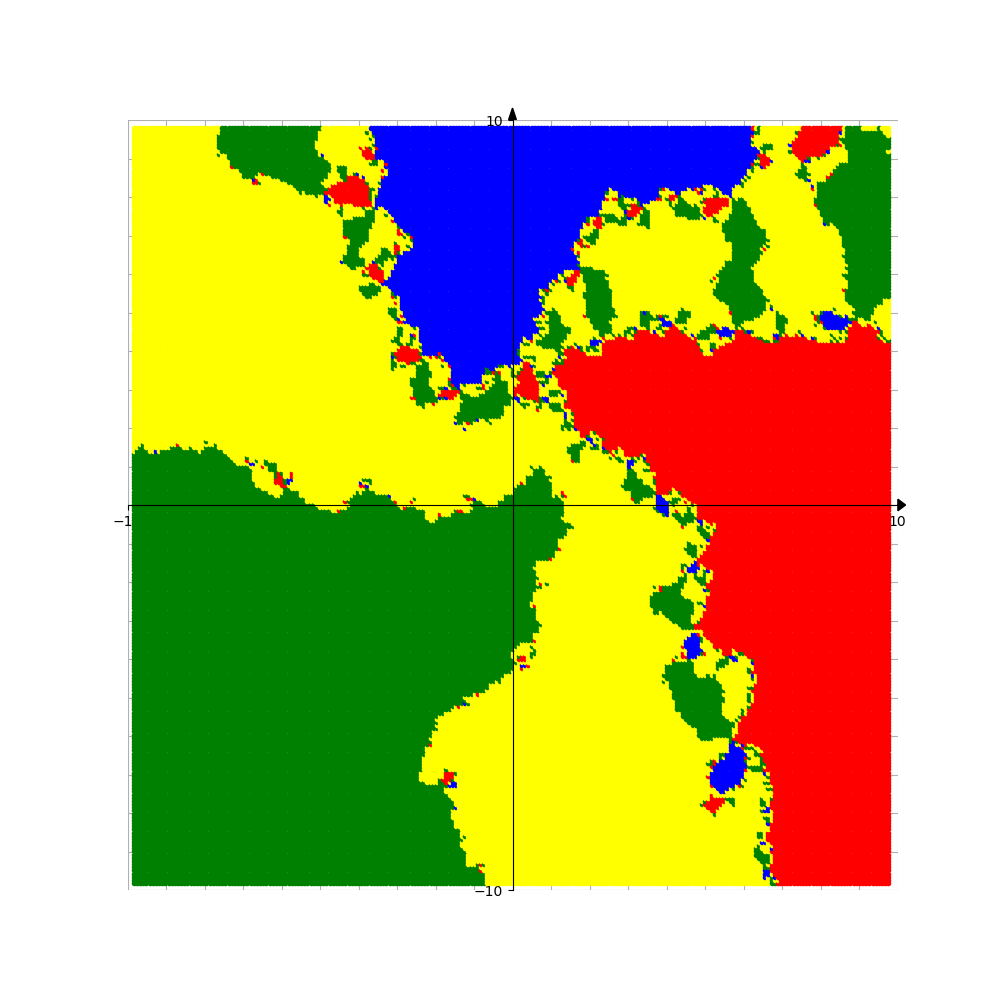}

    \bigskip
    \includegraphics[width=5.5cm]{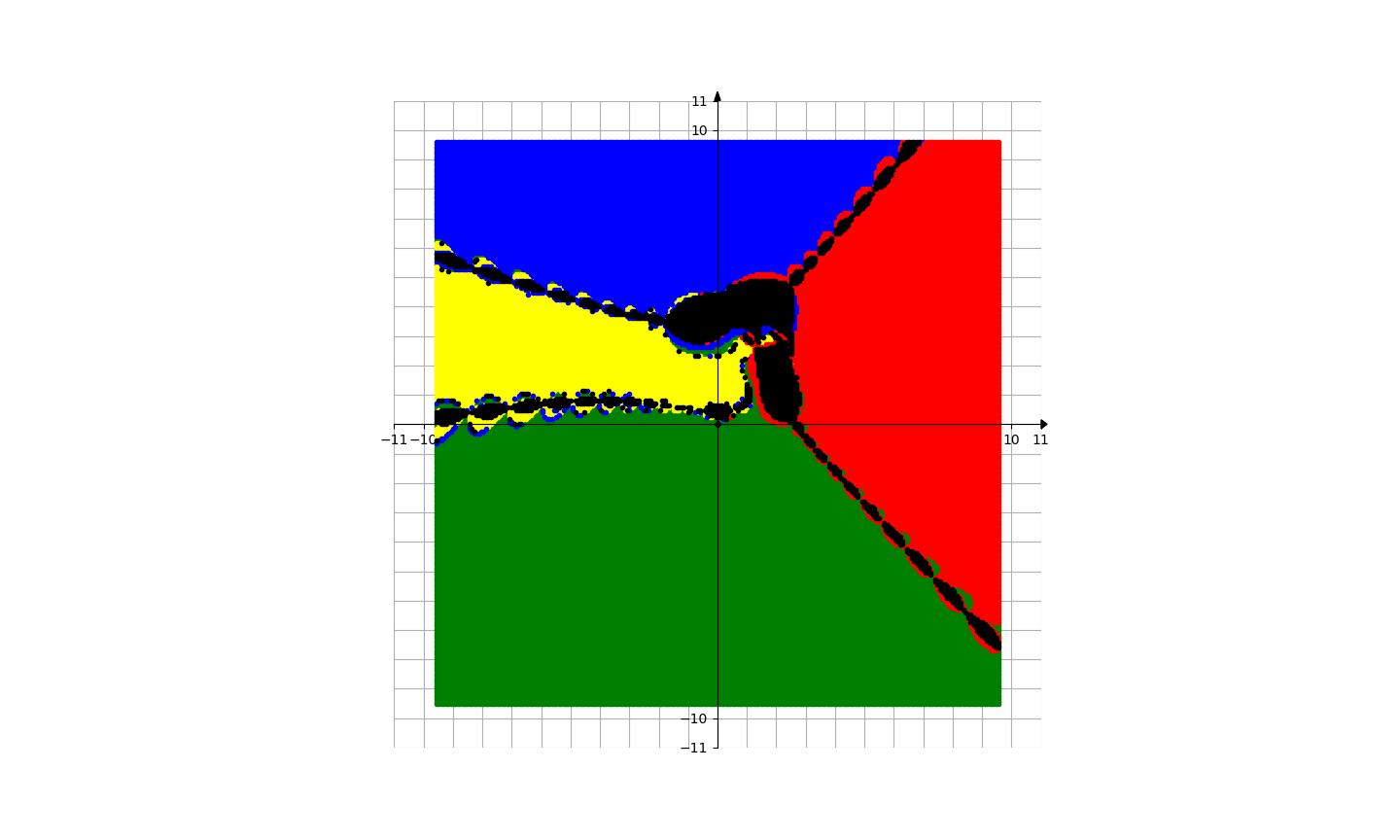}
    \includegraphics[width=3cm]{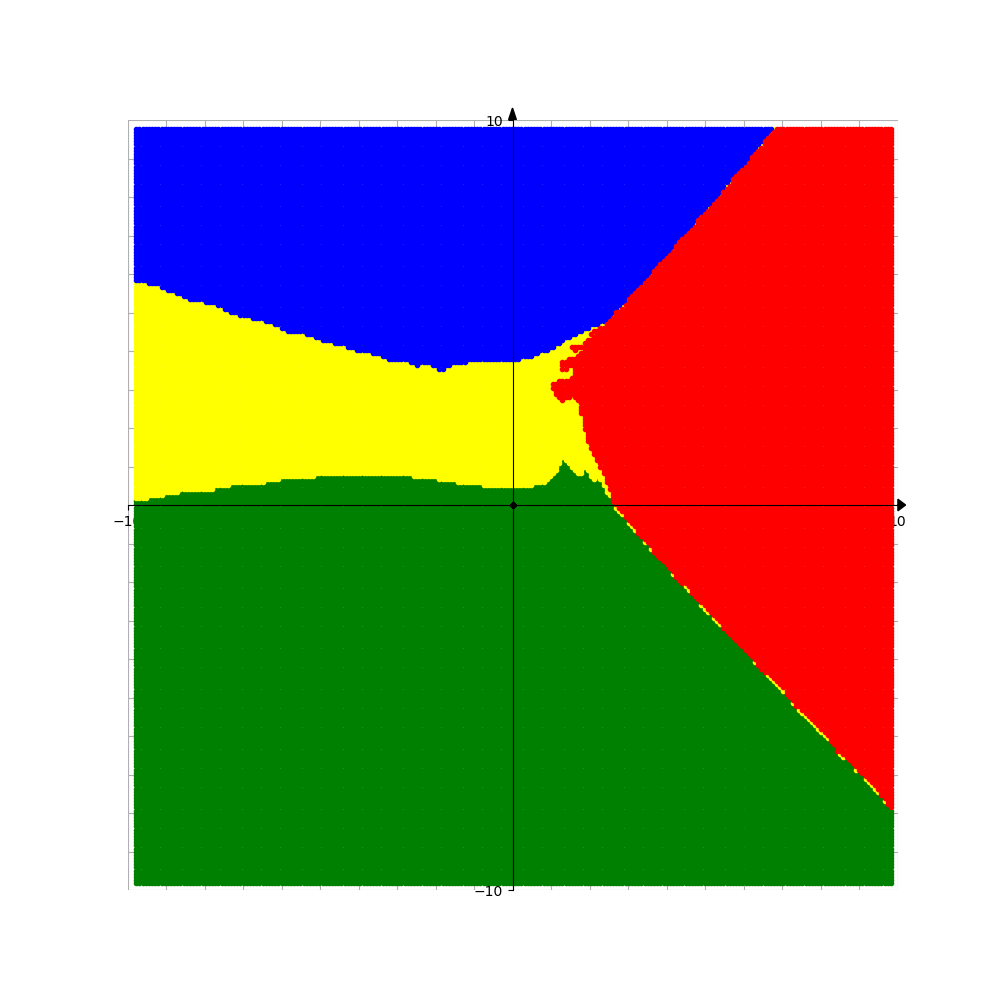}
    
    \bigskip
    \includegraphics[width=3cm]{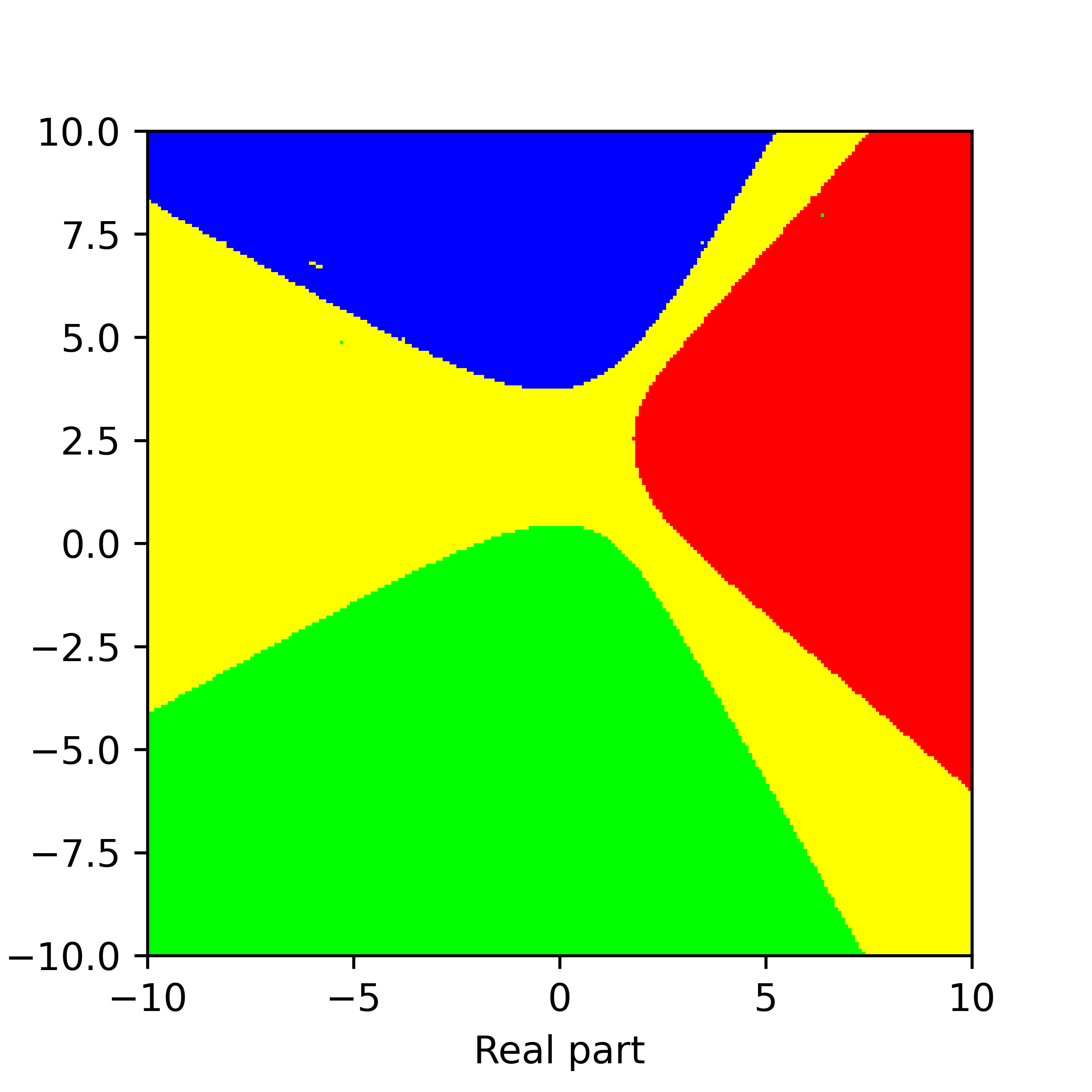}
    \includegraphics[width=3cm]{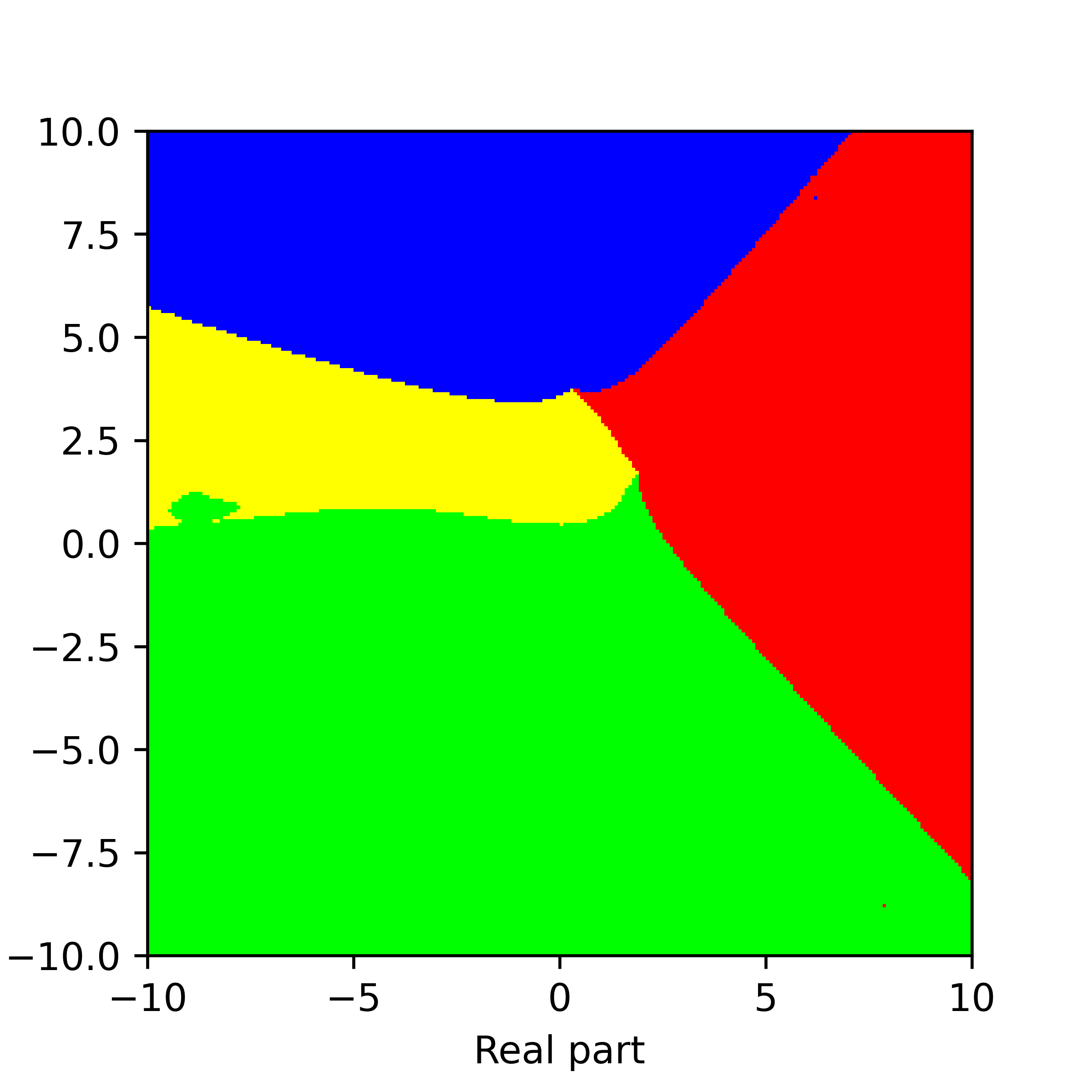}
    \includegraphics[width=3cm]{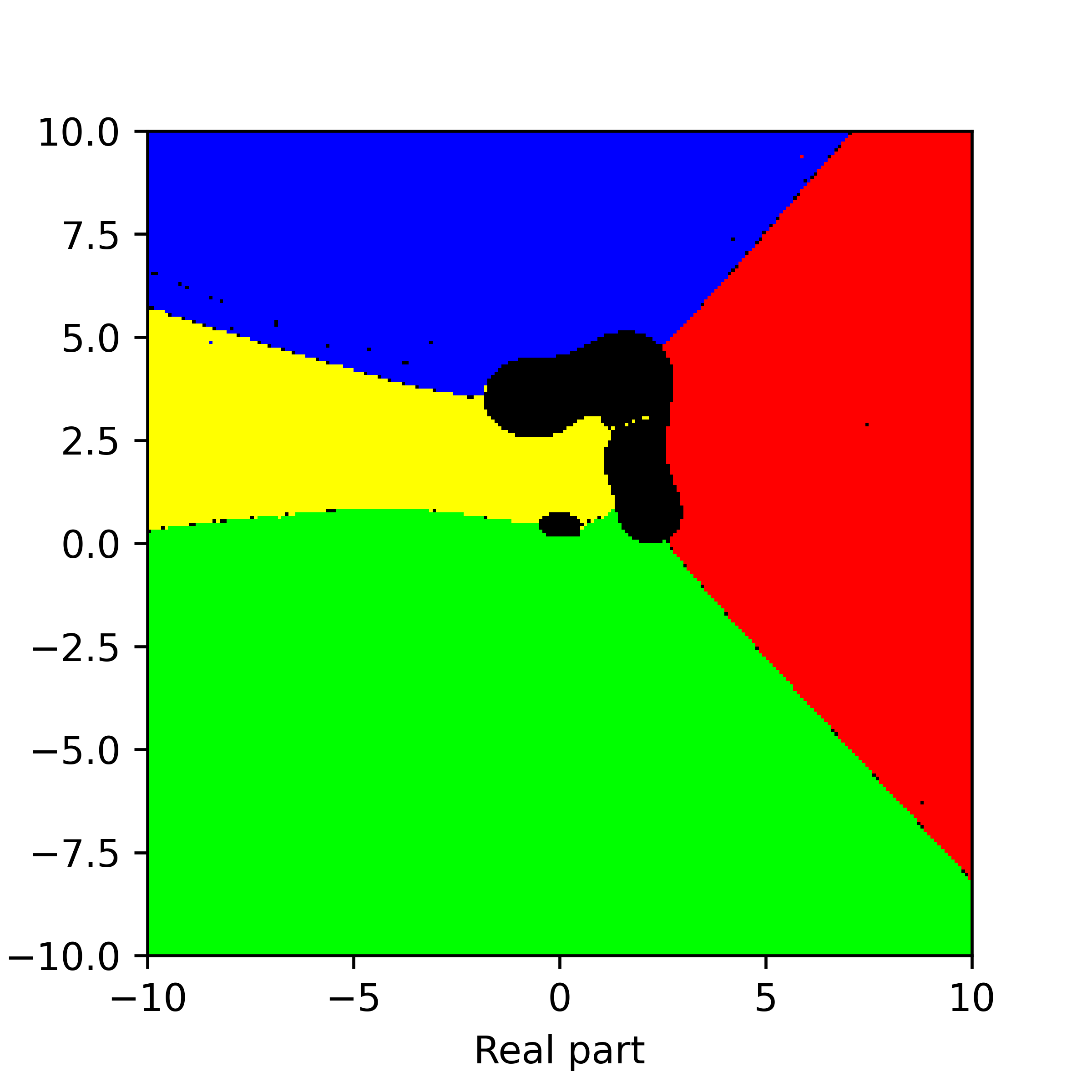}
    
    \caption{Basins of attraction for finding roots of the function $f_{12}$ by different methods. Pictures are referenced to from top to bottom, from left to right. Row 1: left picture is Voronoi's diagram, central picture is for Newton's method, right picture is for Random Relaxed Newton's method. Row 2: left picture is for Newton's method vOptimization, right picture is for BNQN. Row 3: left picture is for Newton's flow, central picture is for Newton's flow vFraction, right picture is for Newton's flow vOptimization. The black points in some of these pictures are those in the basin of attraction of critical points of $f_{12}$.}
    \label{fig:f12}
\end{figure}

\begin{figure}
    \centering
    \includegraphics[width=5cm]{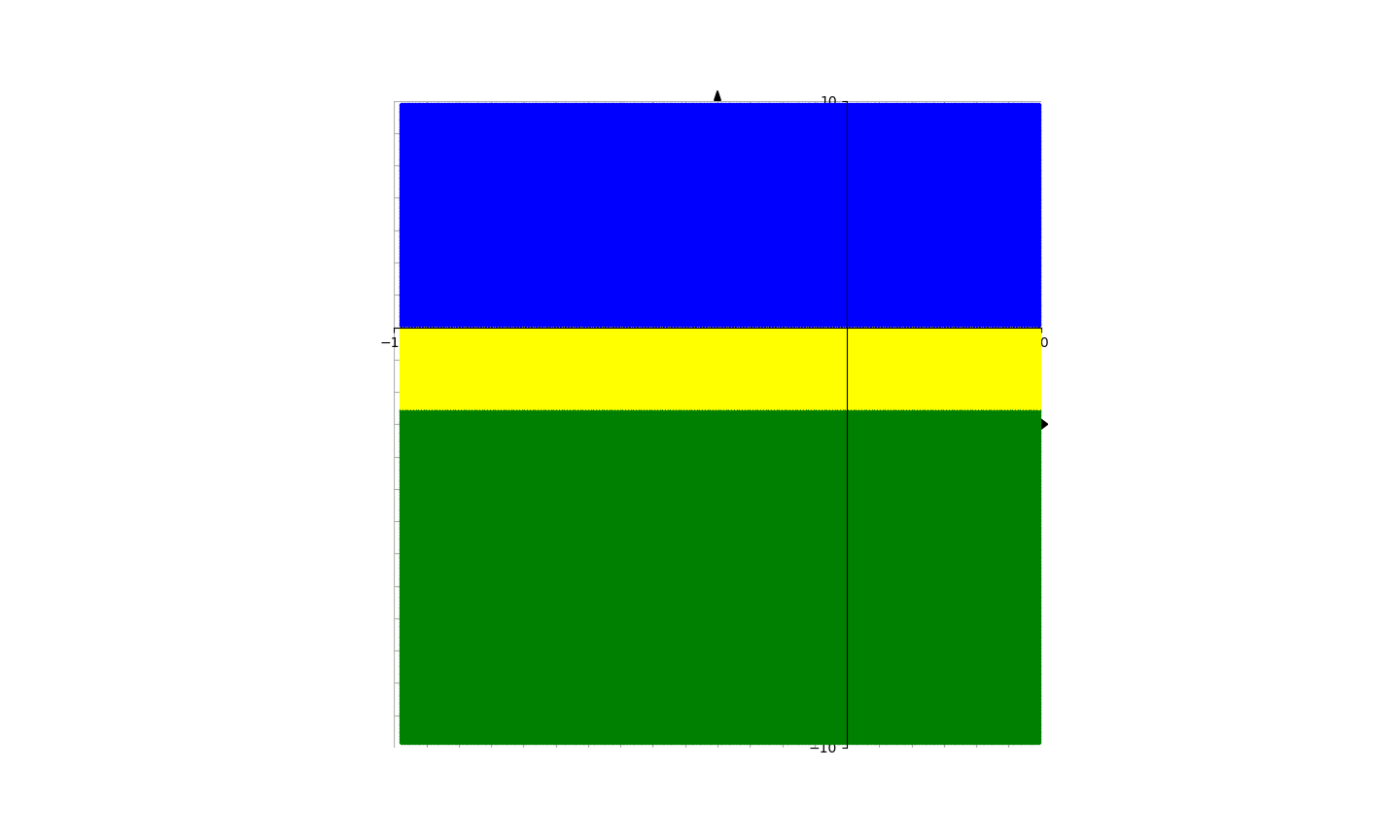}
    \includegraphics[width=3cm]{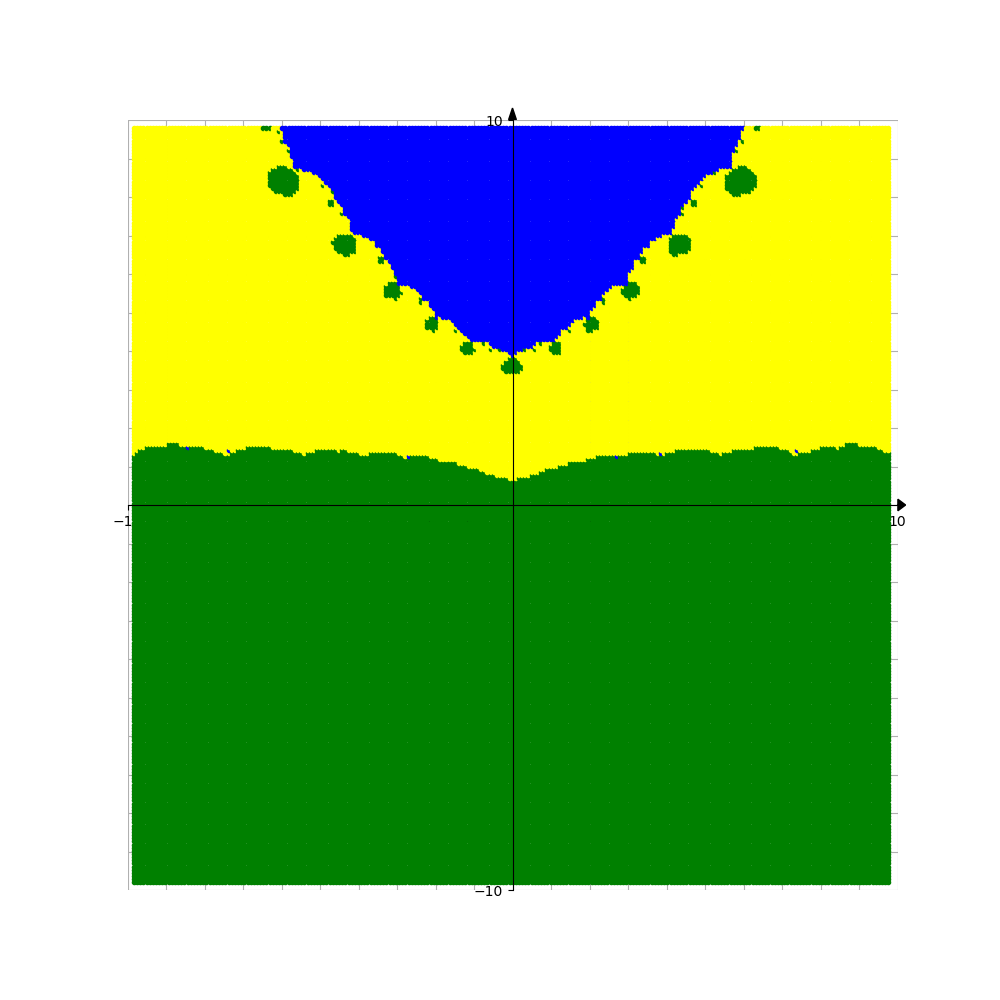}
    \includegraphics[width=3cm]{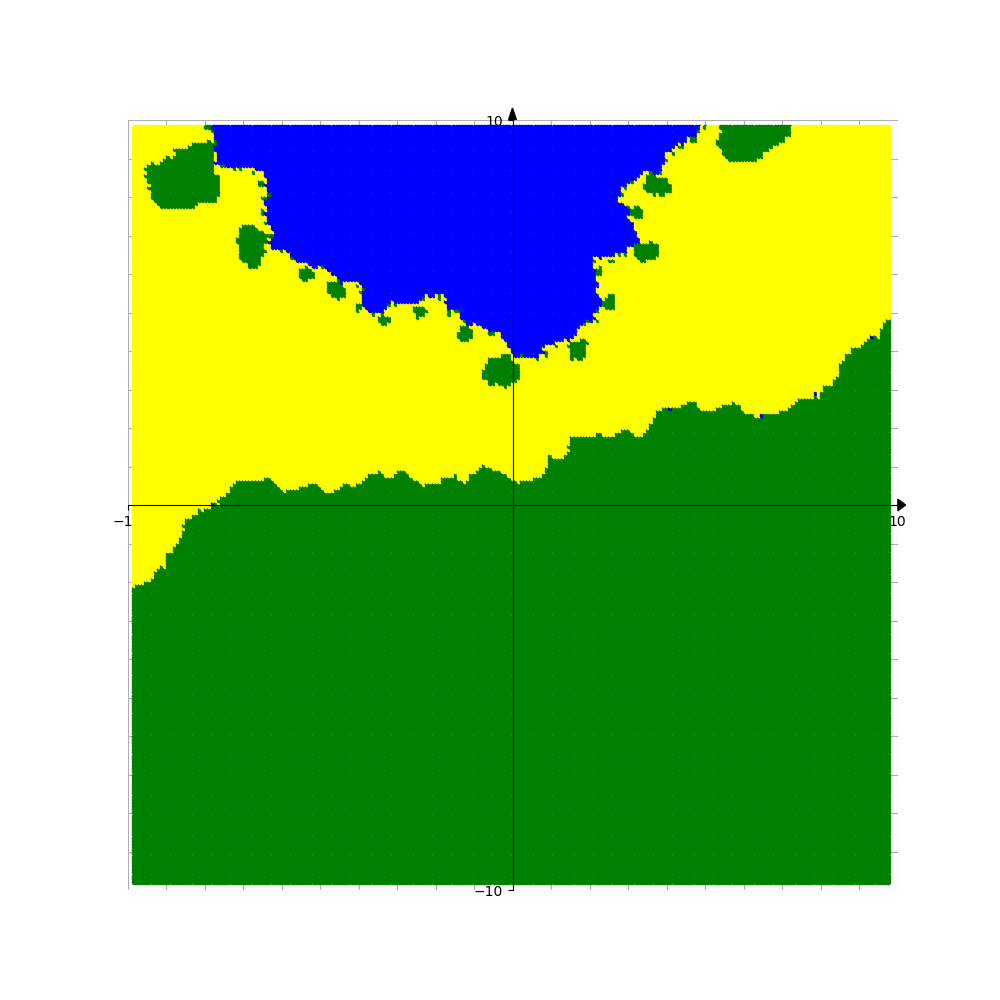}

    \bigskip
    \includegraphics[width=5.5cm]{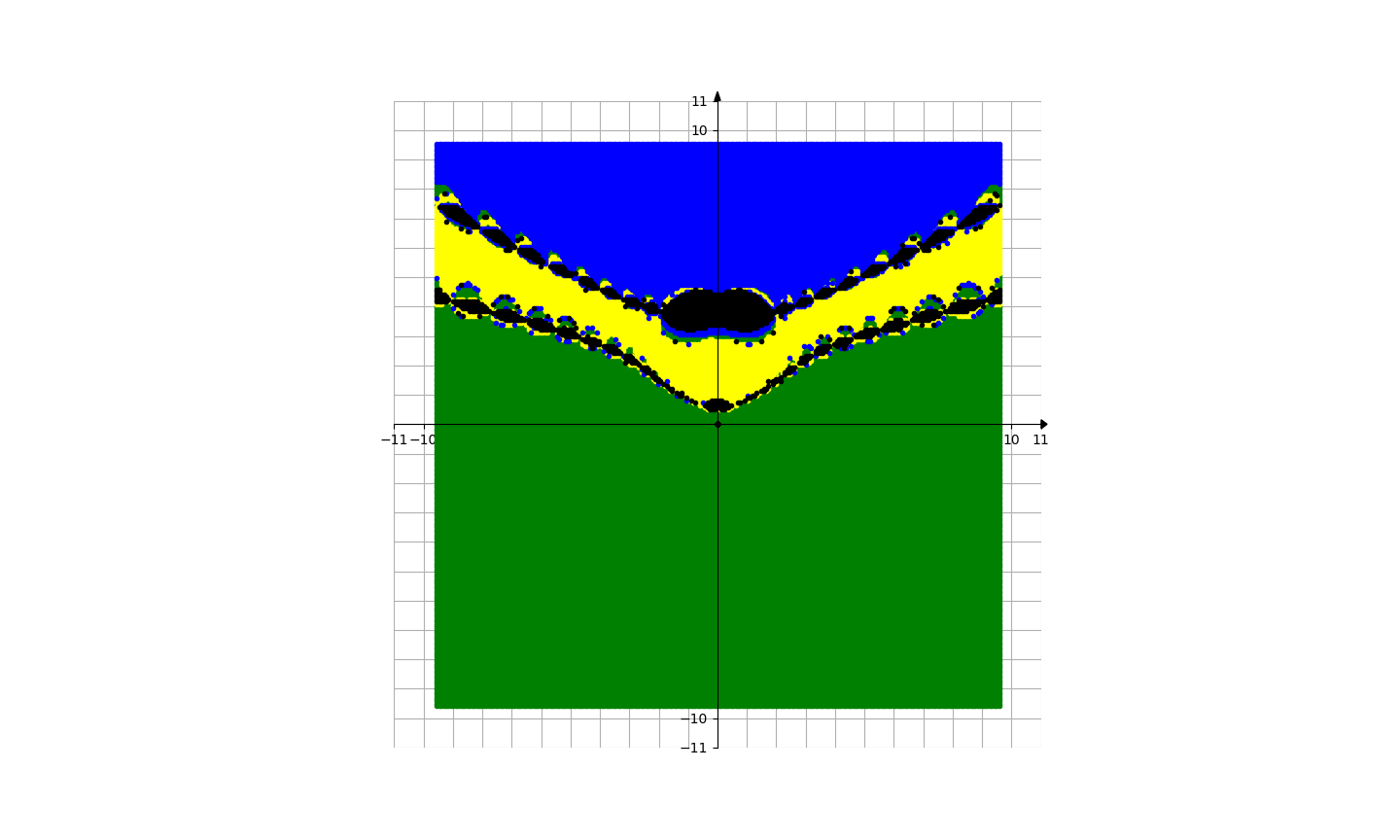}
    \includegraphics[width=3cm]{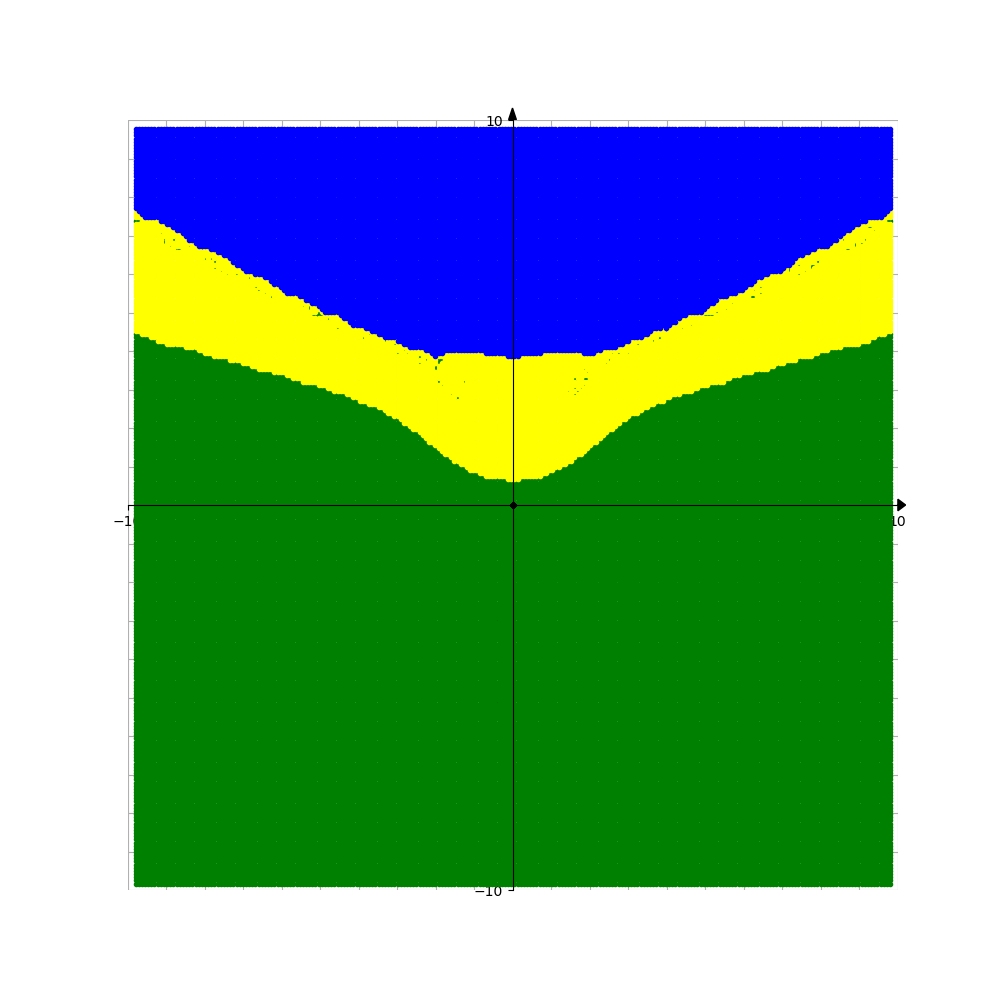}
    
    \bigskip
    \includegraphics[width=3cm]{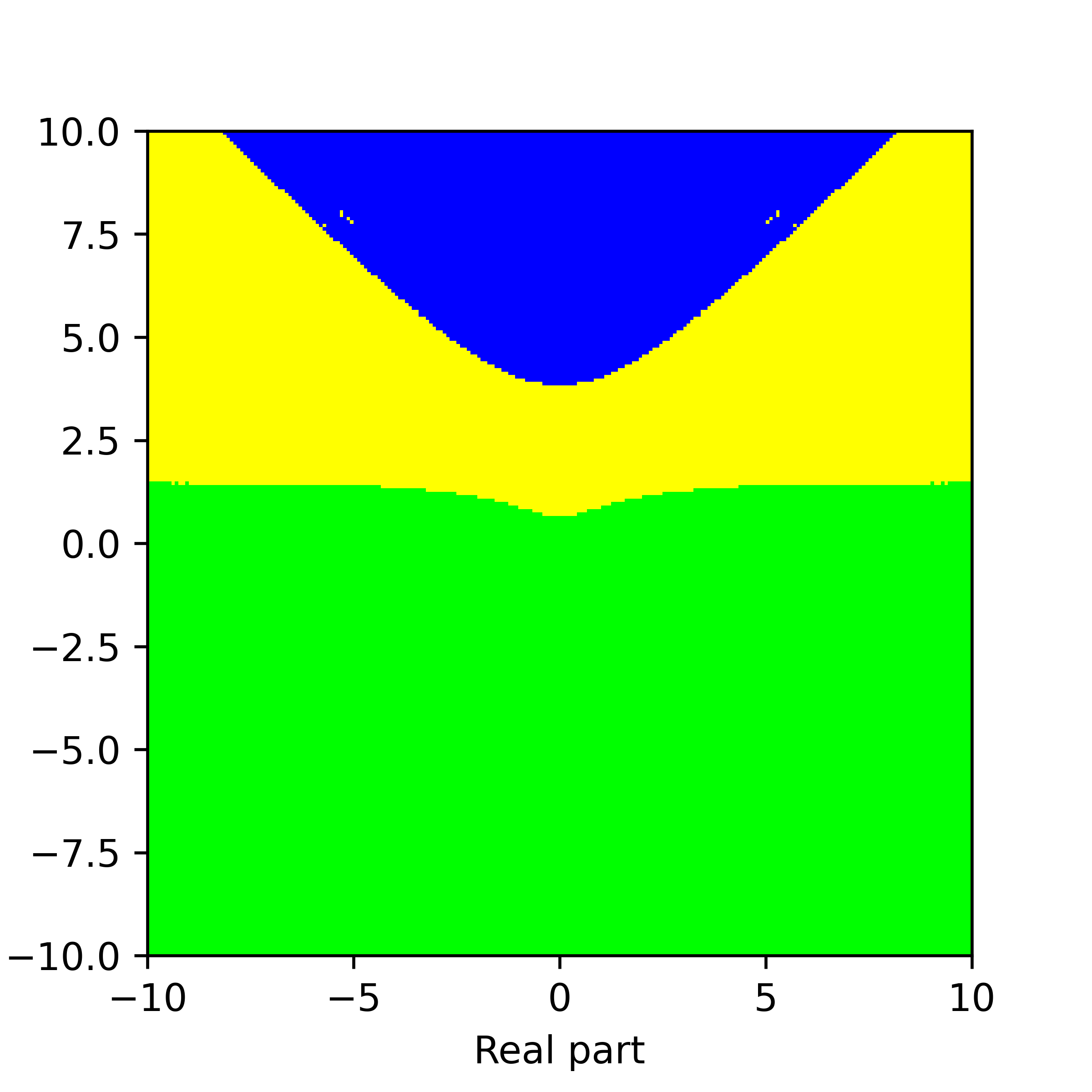}
    \includegraphics[width=3cm]{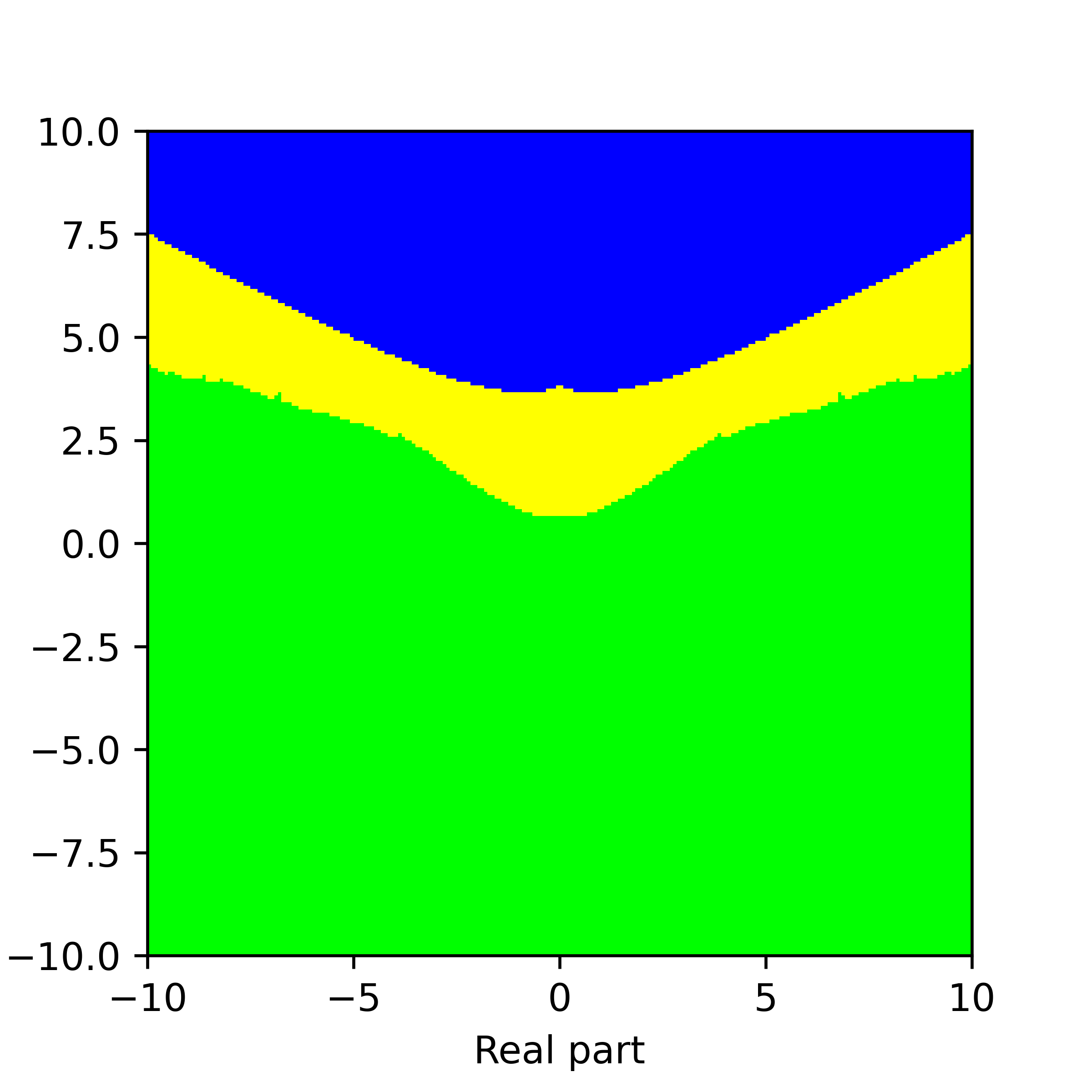}
    \includegraphics[width=3cm]{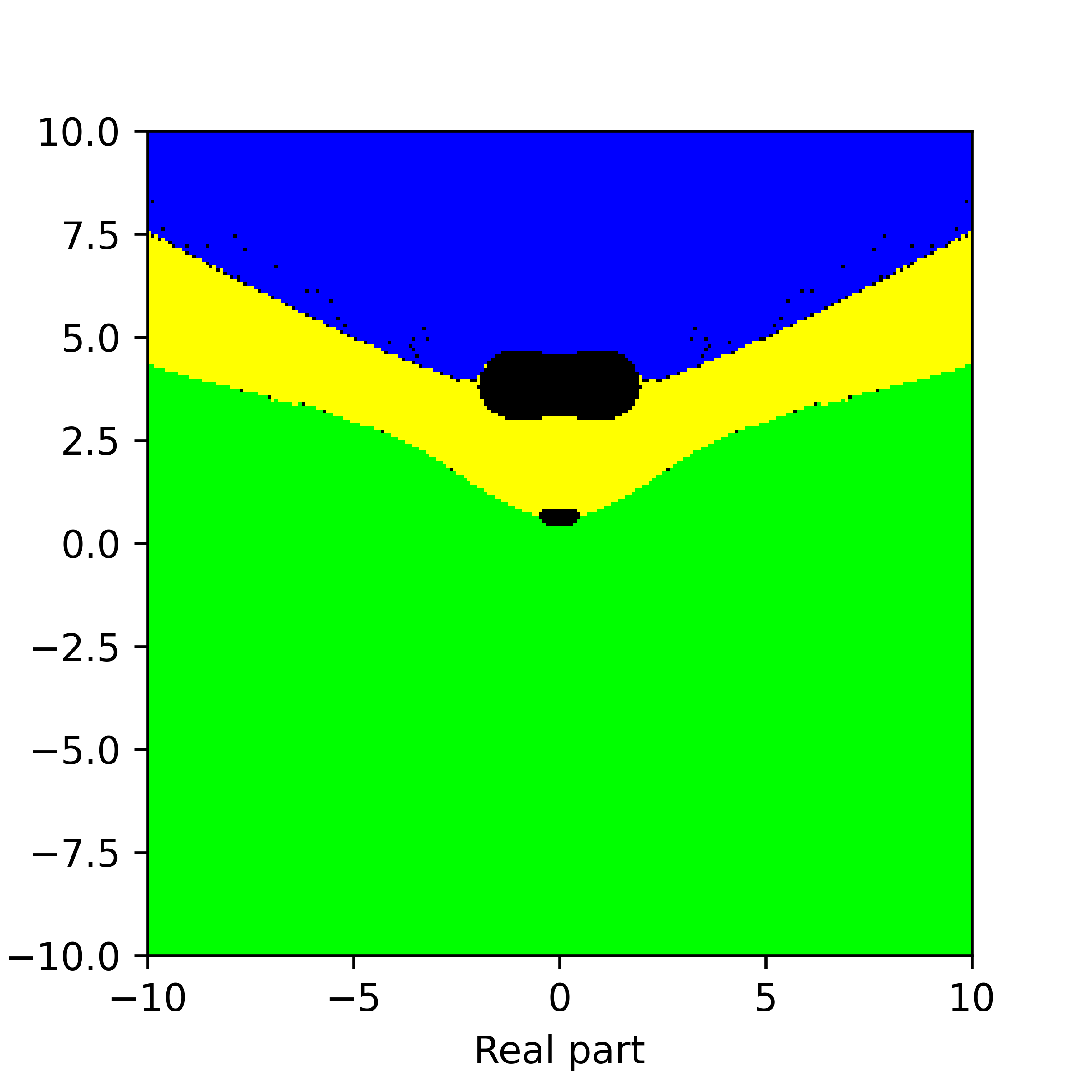}
    
    \caption{Basins of attraction for finding roots of the function $f_{13}$ by different methods. Pictures are referenced to from top to bottom, from left to right. Row 1: left picture is Voronoi's diagram, central picture is for Newton's method, right picture is for Random Relaxed Newton's method. Row 2: left picture is for Newton's method vOptimization, right picture is for BNQN. Row 3: left picture is for Newton's flow, central picture is for Newton's flow vFraction, right picture is for Newton's flow vOptimization. The black points in some of these pictures are those in the basin of attraction of critical points of $f_{13}$.}
    \label{fig:f13}
\end{figure}

\begin{figure}
    \centering
    \includegraphics[width=5cm]{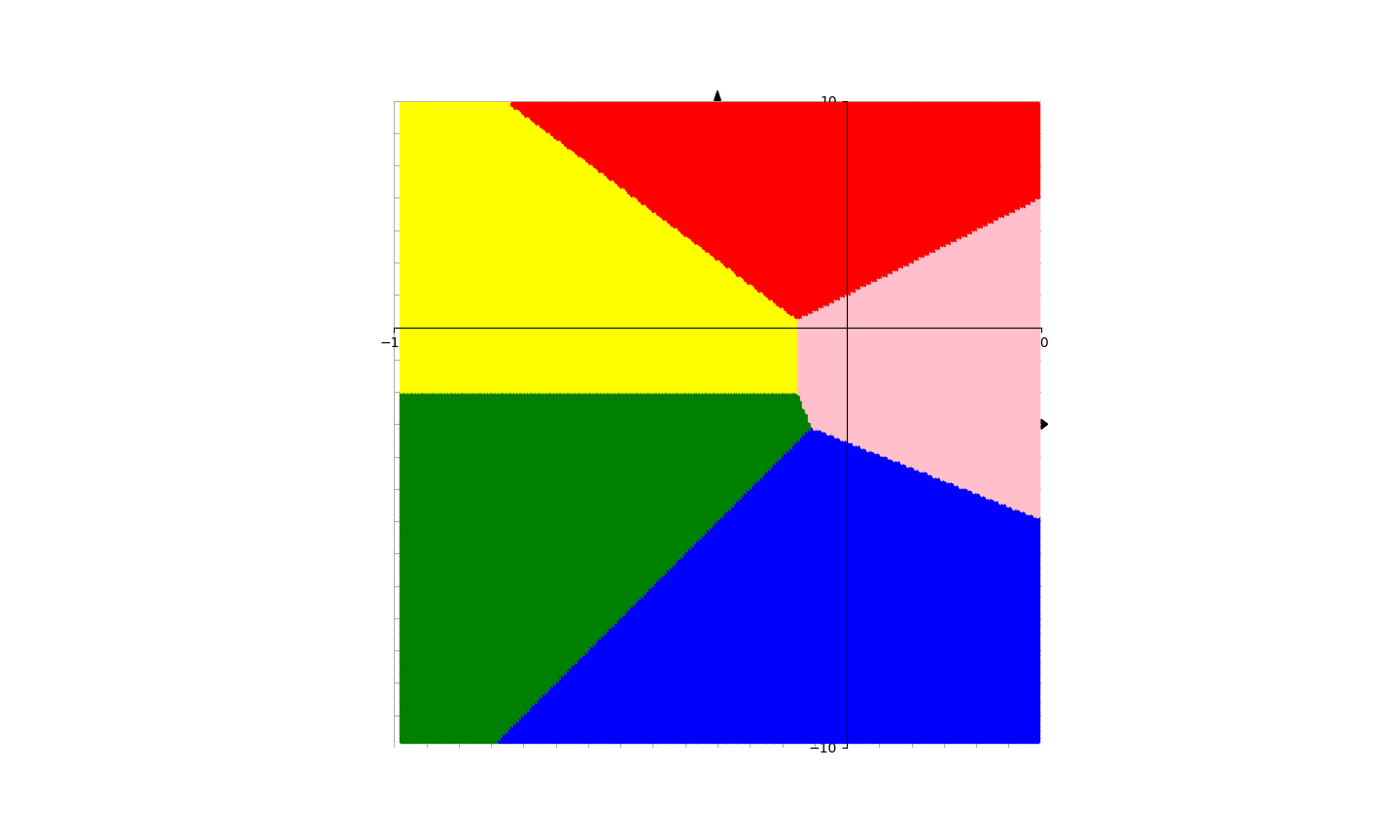}
    \includegraphics[width=3cm]{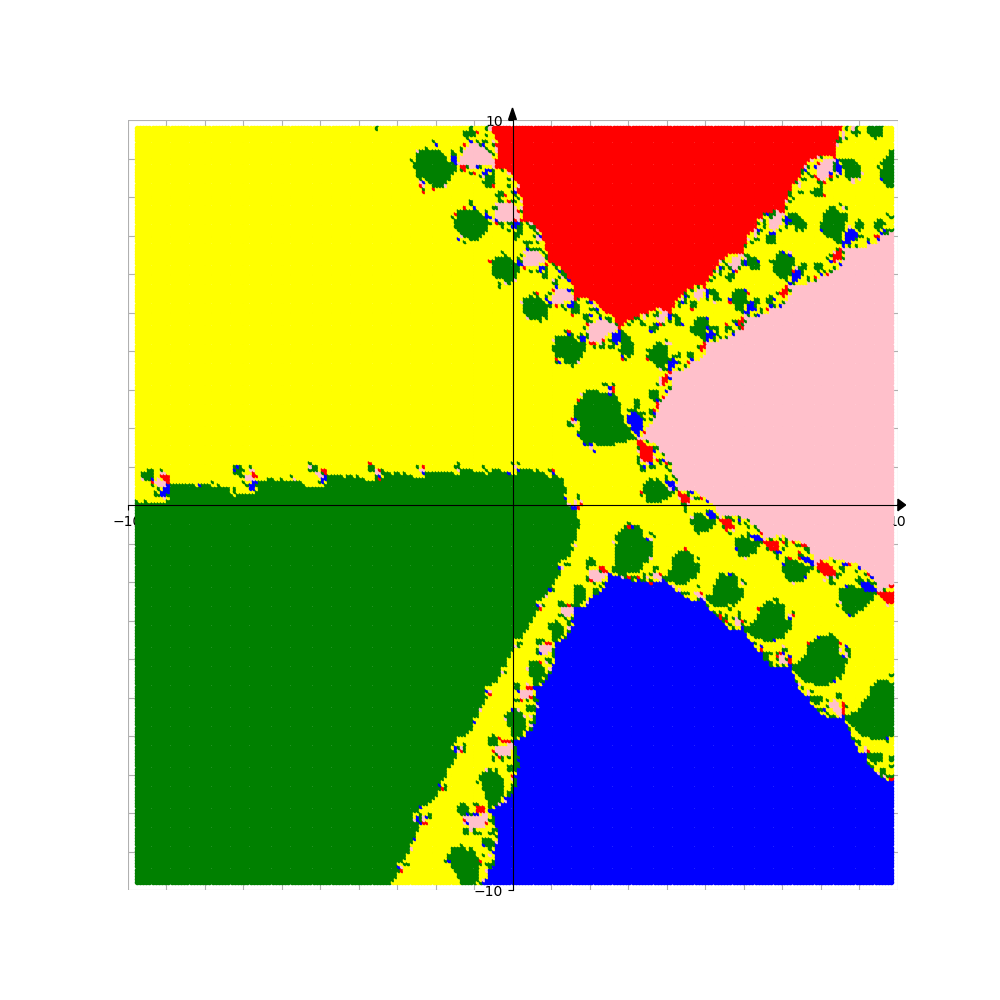}
    \includegraphics[width=3cm]{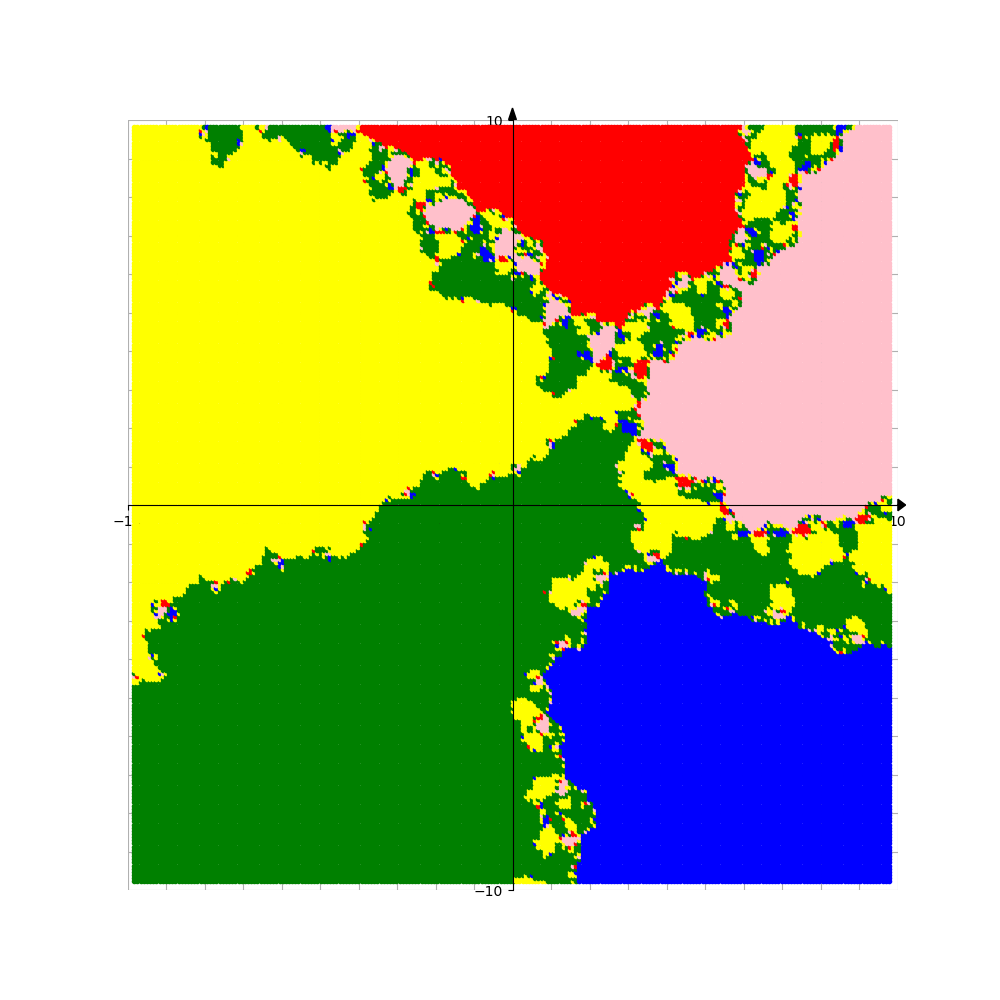}

    \bigskip
    \includegraphics[width=5.5cm]{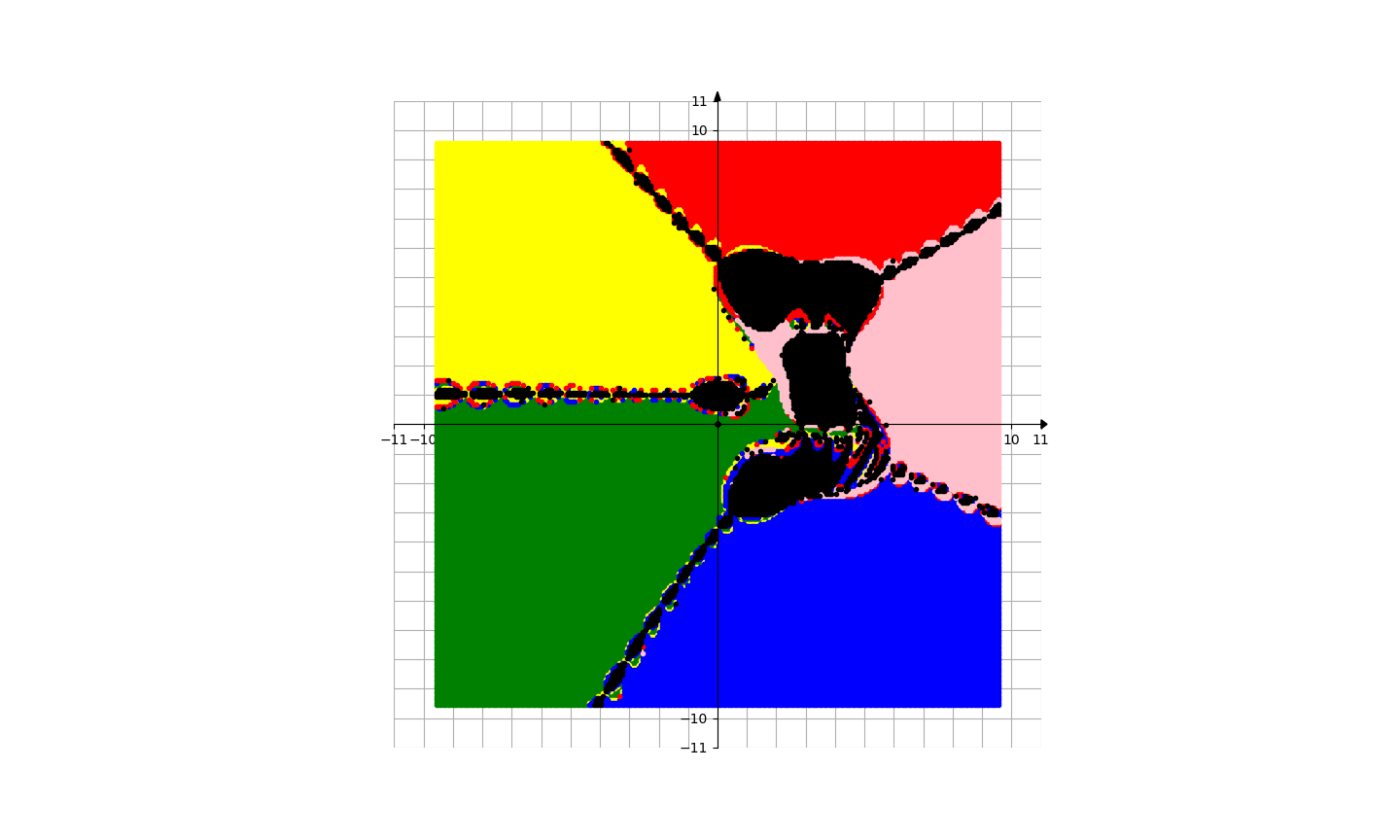}
    \includegraphics[width=3cm]{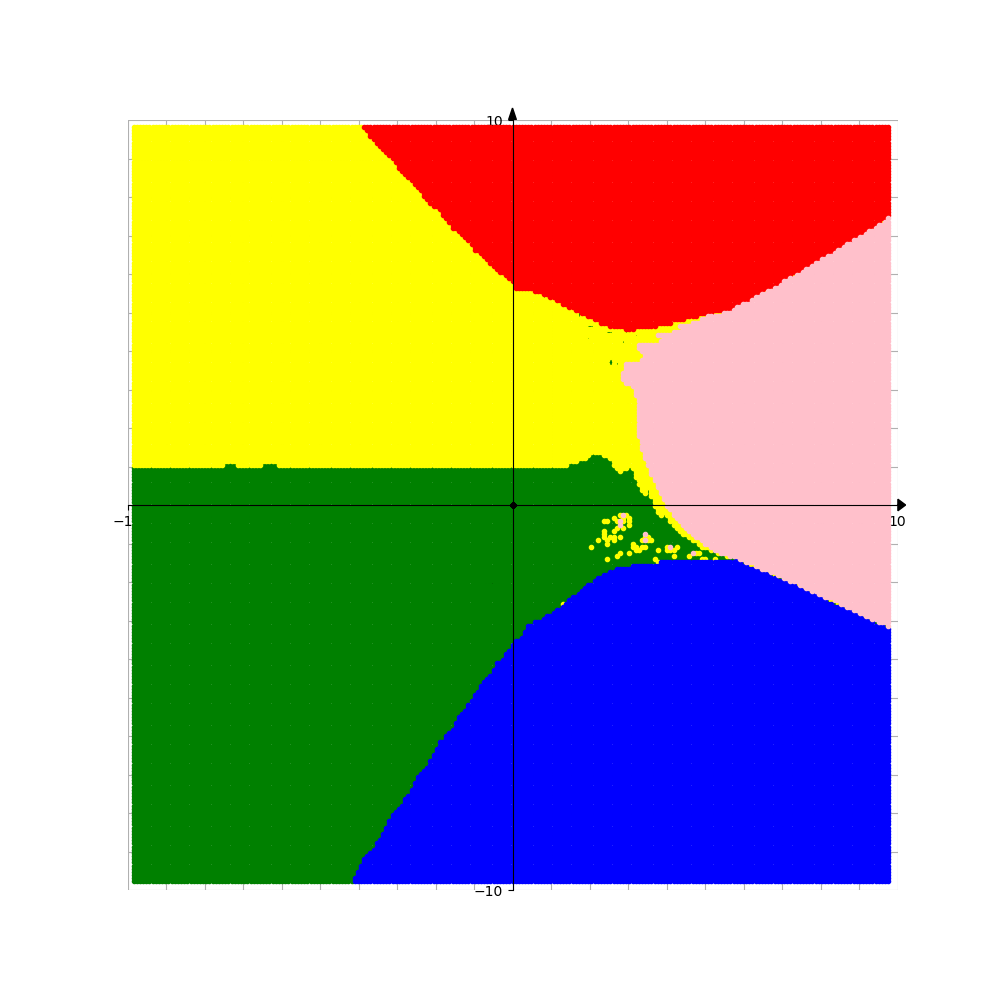}
    
    \bigskip
    \includegraphics[width=3cm]{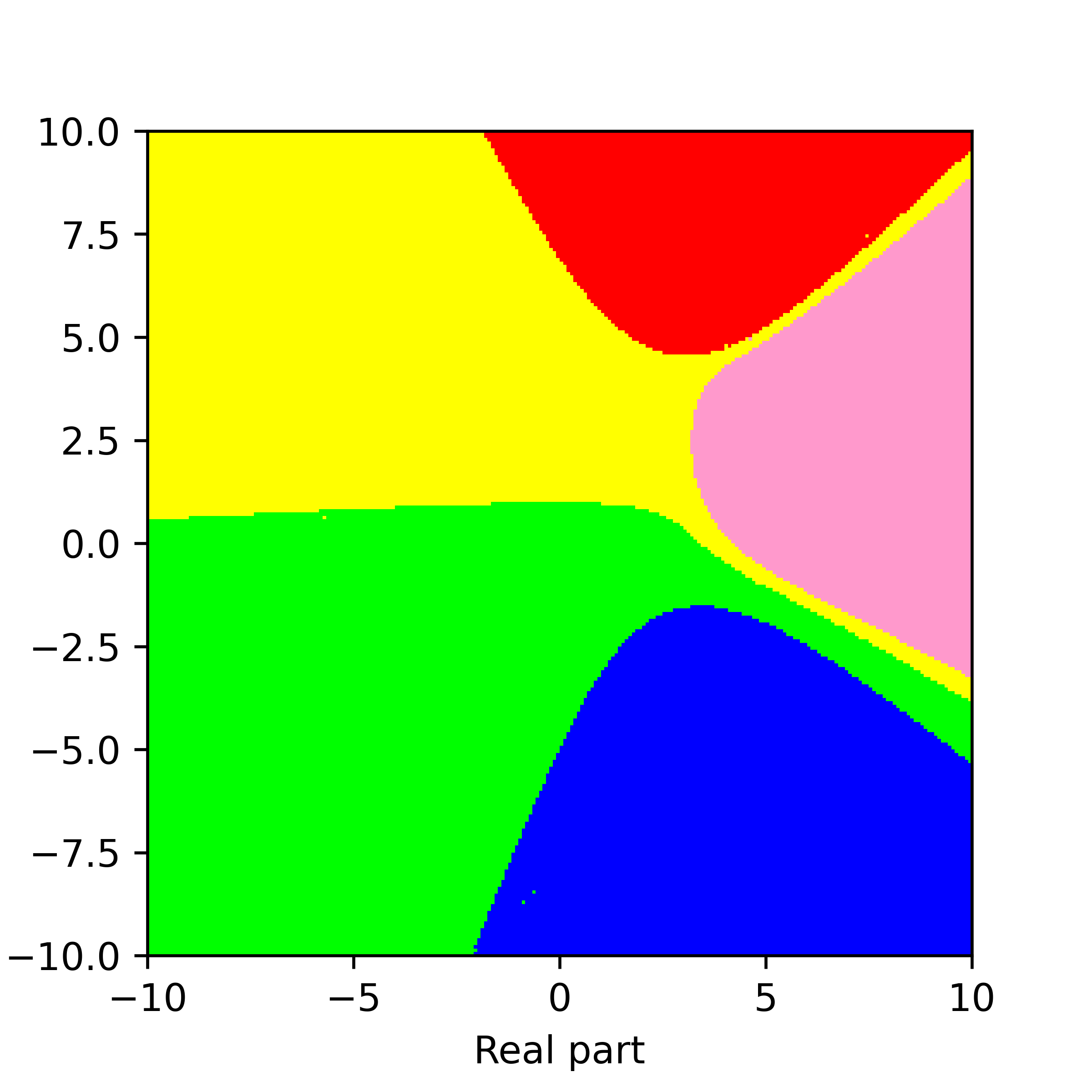}
    \includegraphics[width=3cm]{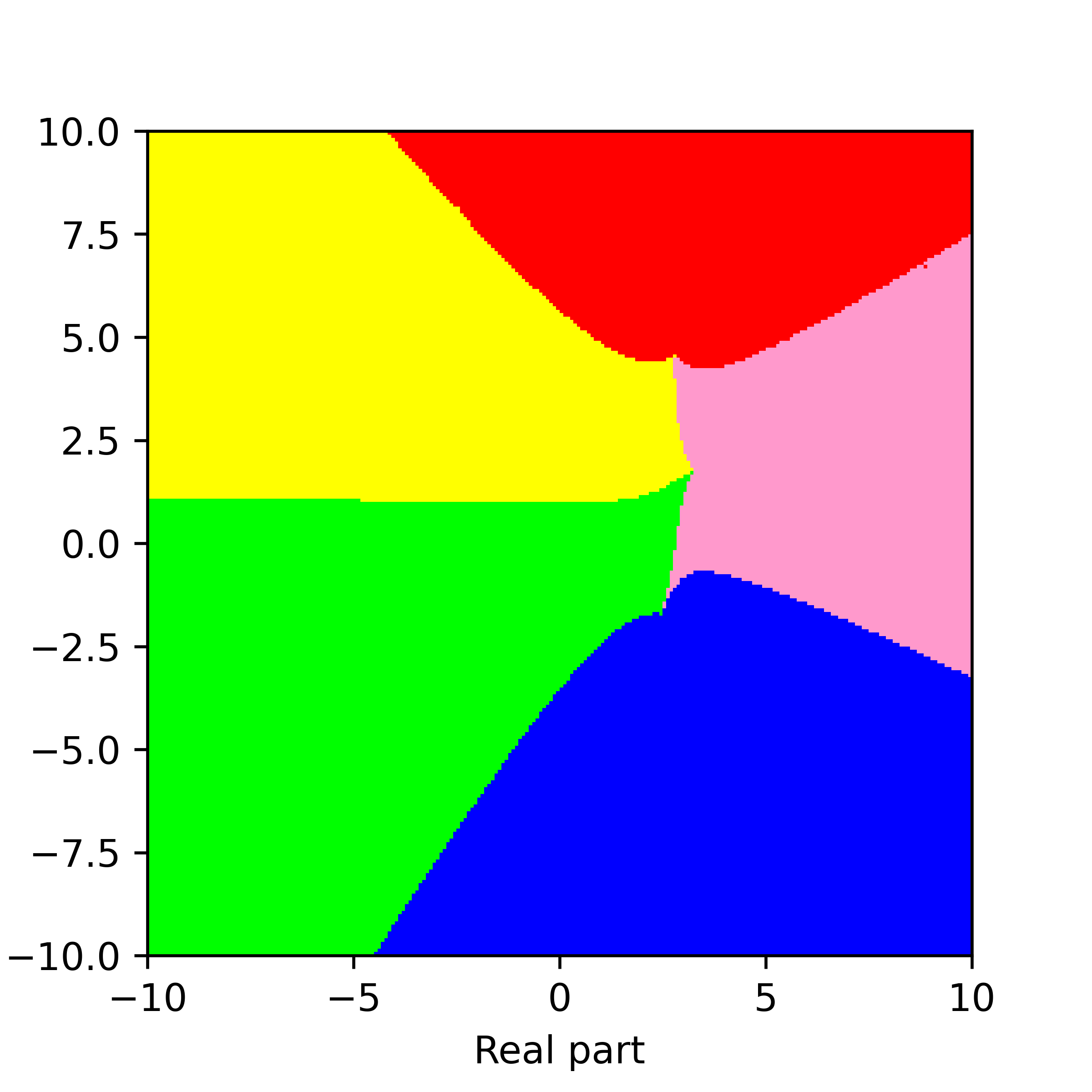}
    \includegraphics[width=3cm]{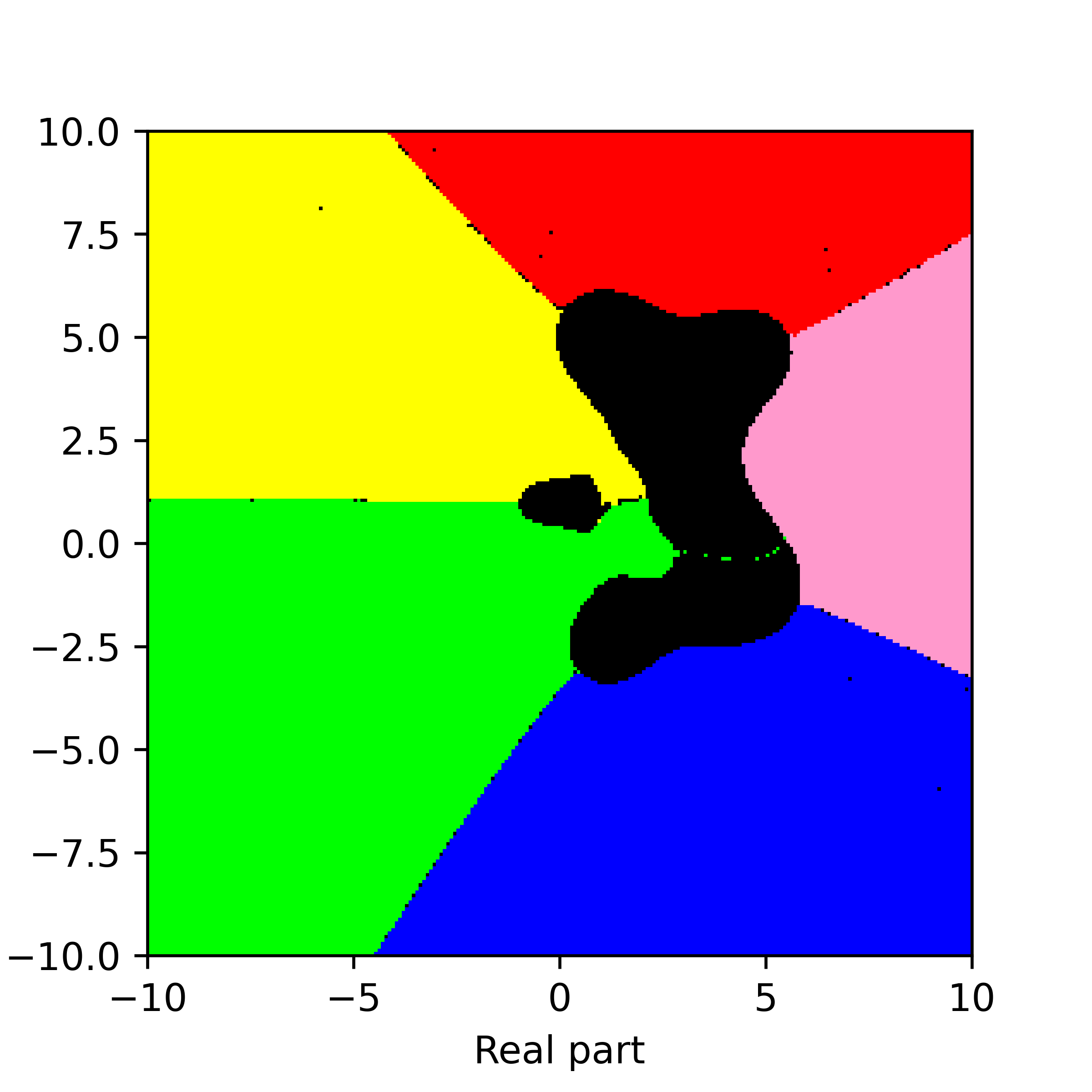}
    
    \caption{Basins of attraction for finding roots of the function $f_{14}$ by different methods. Pictures are referenced to from top to bottom, from left to right. Row 1: left picture is Voronoi's diagram, central picture is for Newton's method, right picture is for Random Relaxed Newton's method. Row 2: left picture is for Newton's method vOptimization, right picture is for BNQN. Row 3: left picture is for Newton's flow, central picture is for Newton's flow vFraction, right picture is for Newton's flow vOptimization. The black points in some of these pictures are those in the basin of attraction of critical points of $f_{14}$.}
    \label{fig:f14}
\end{figure}

\begin{figure}
    \centering
    \includegraphics[width=5cm]{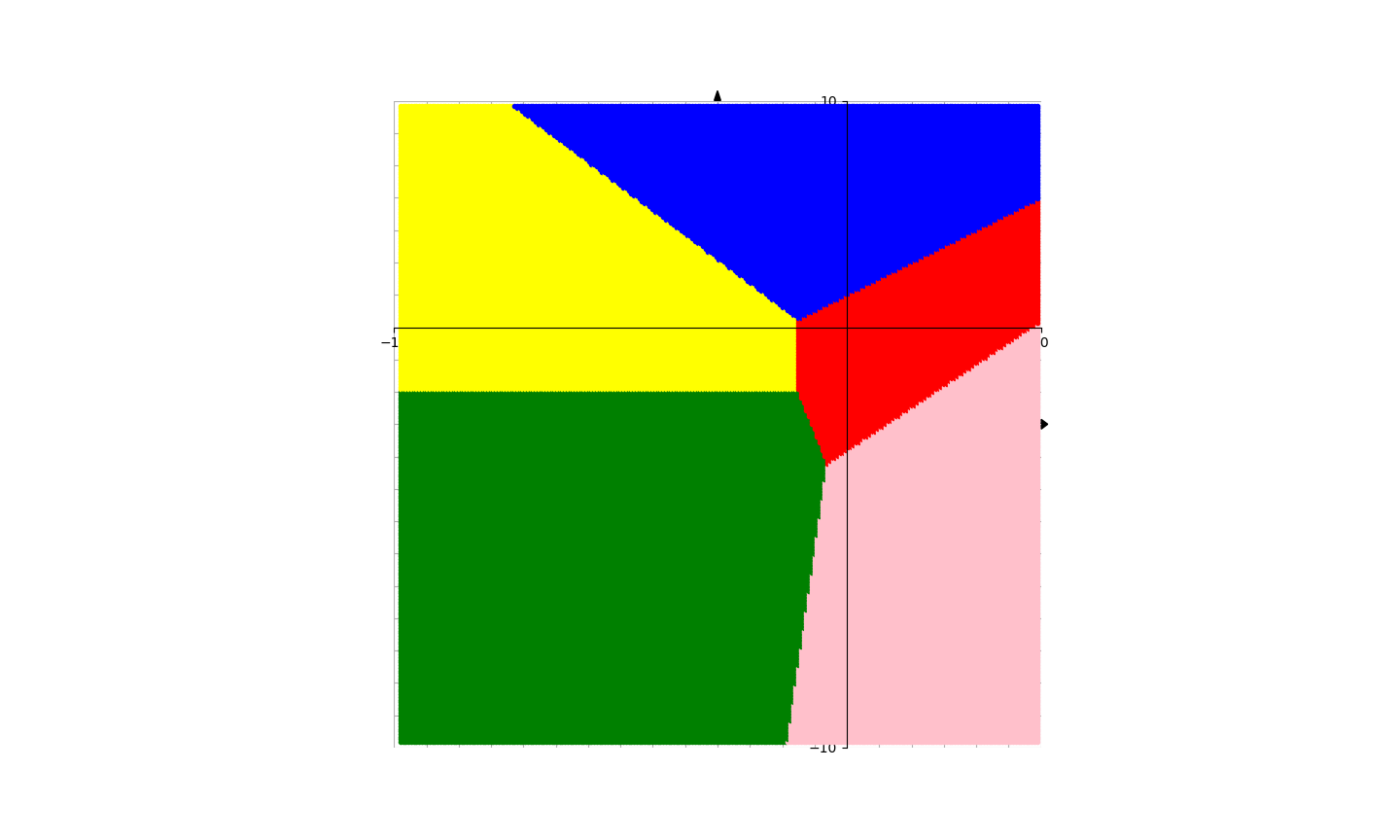}
    \includegraphics[width=3cm]{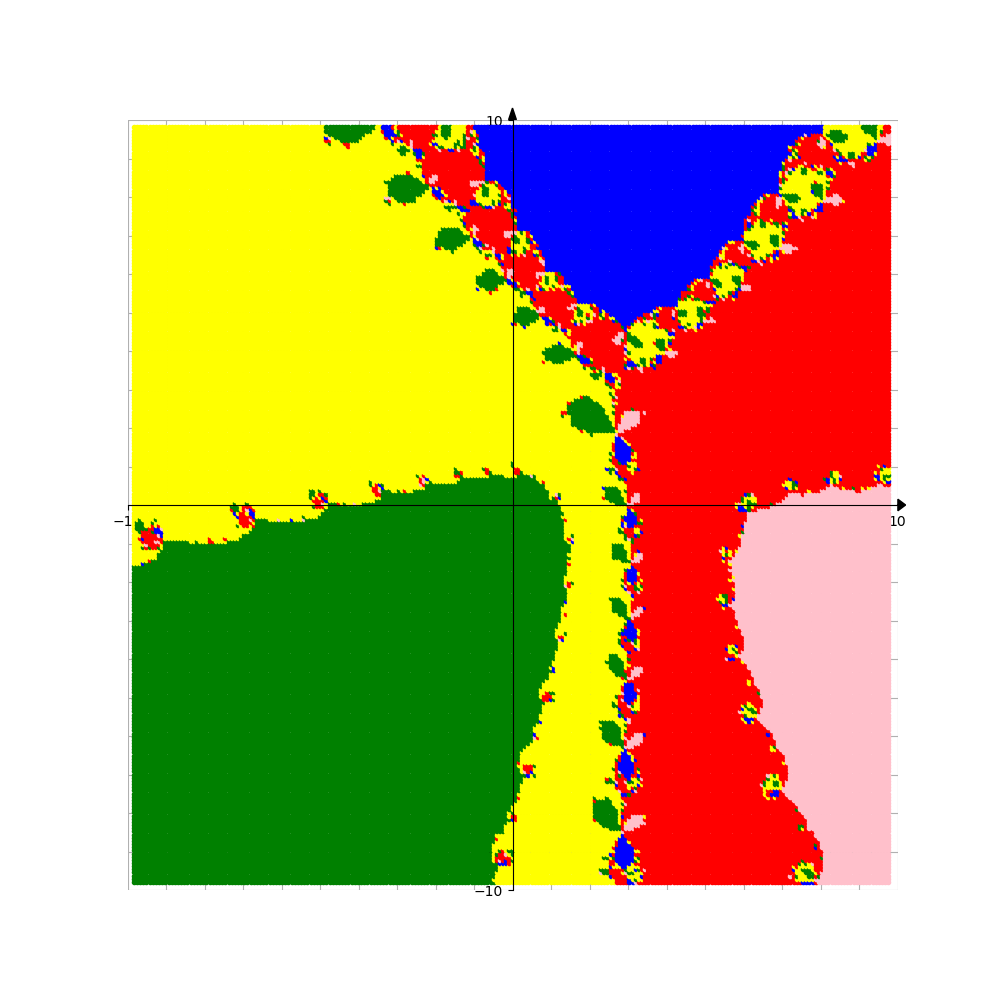}
    \includegraphics[width=3cm]{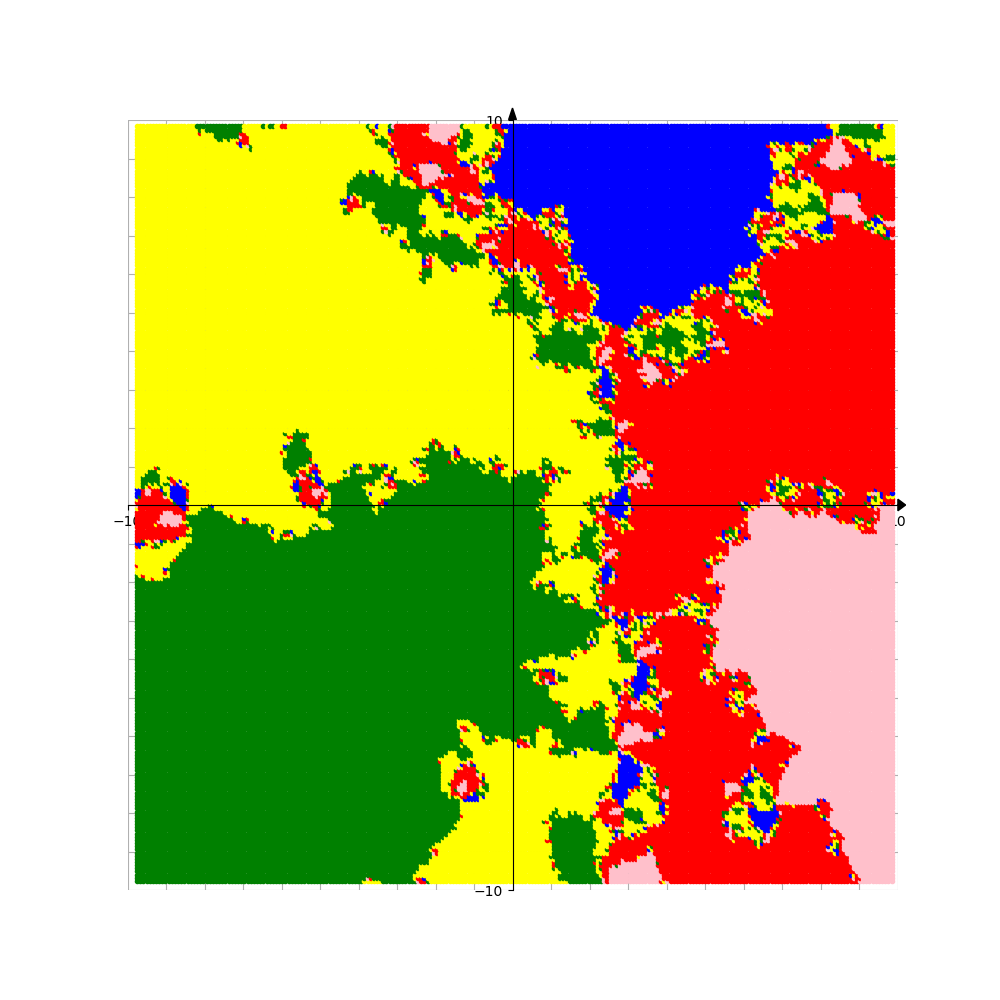}

    \bigskip
    \includegraphics[width=5.5cm]{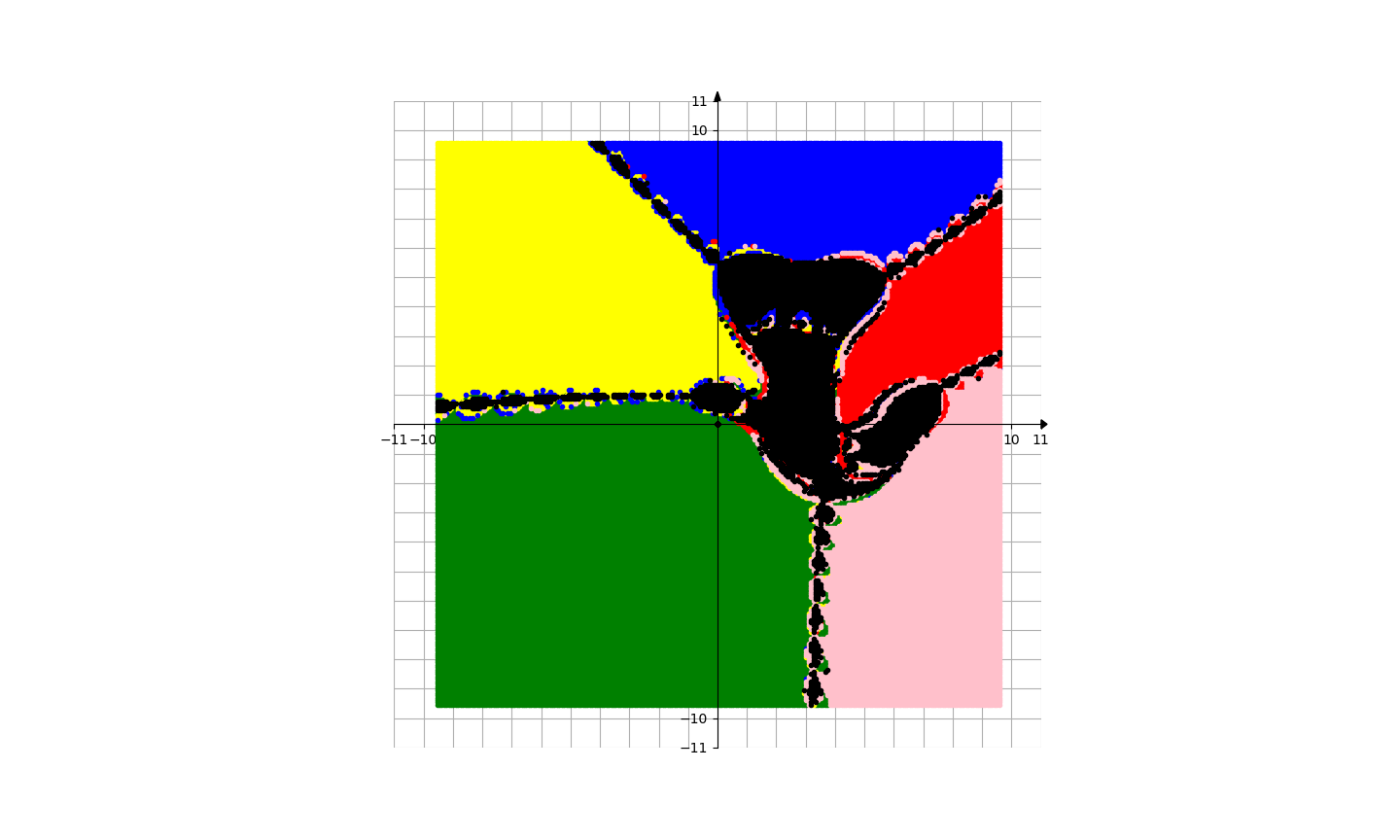}
    \includegraphics[width=3cm]{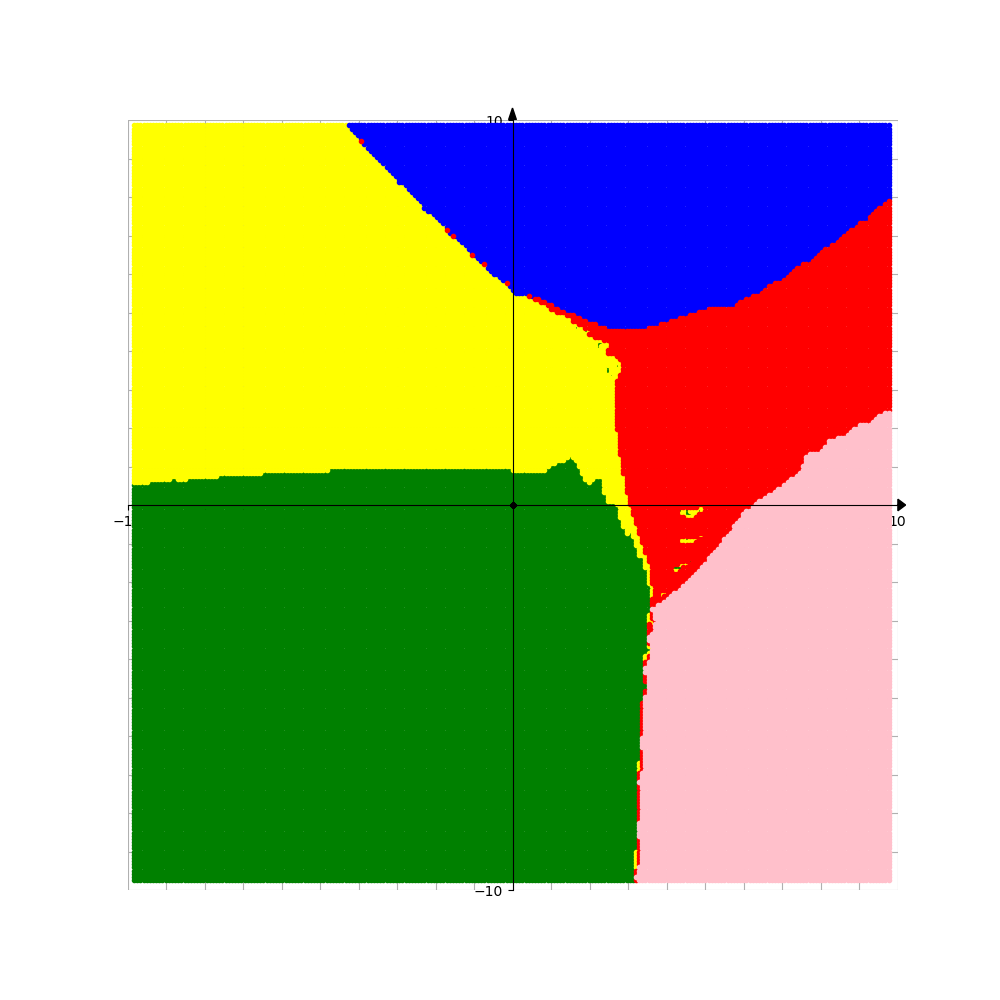}
    
    \bigskip
    \includegraphics[width=3cm]{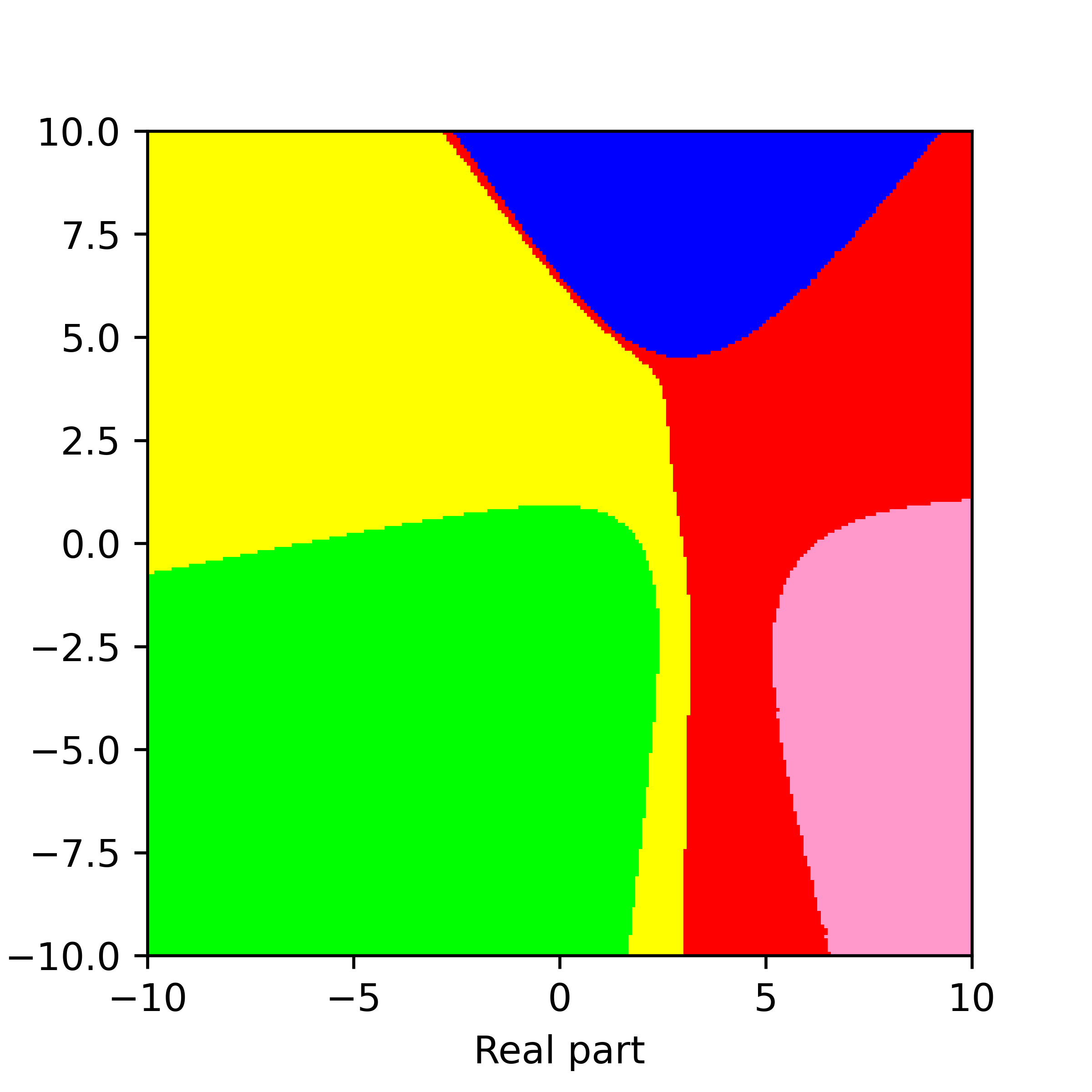}
    \includegraphics[width=3cm]{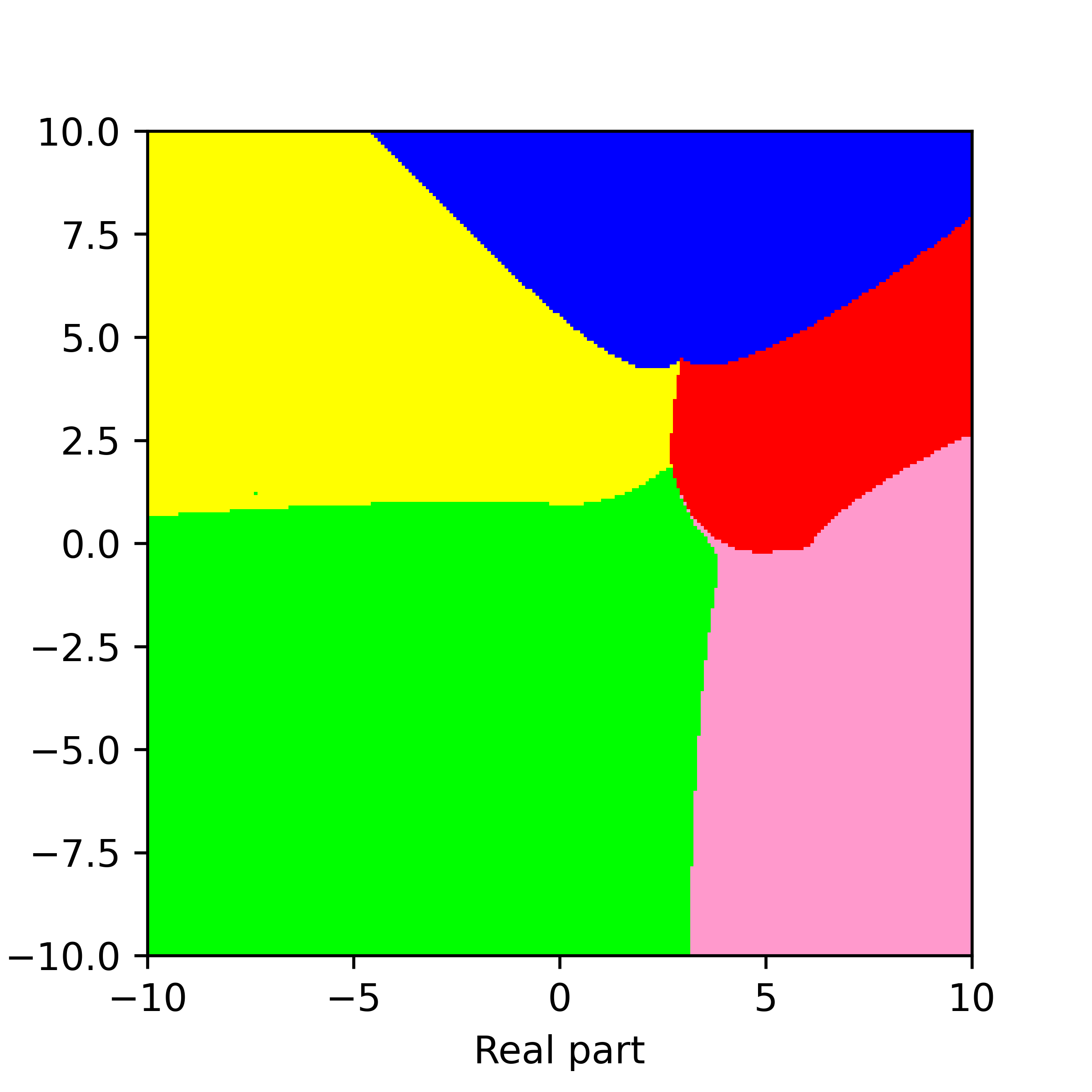}
    \includegraphics[width=3cm]{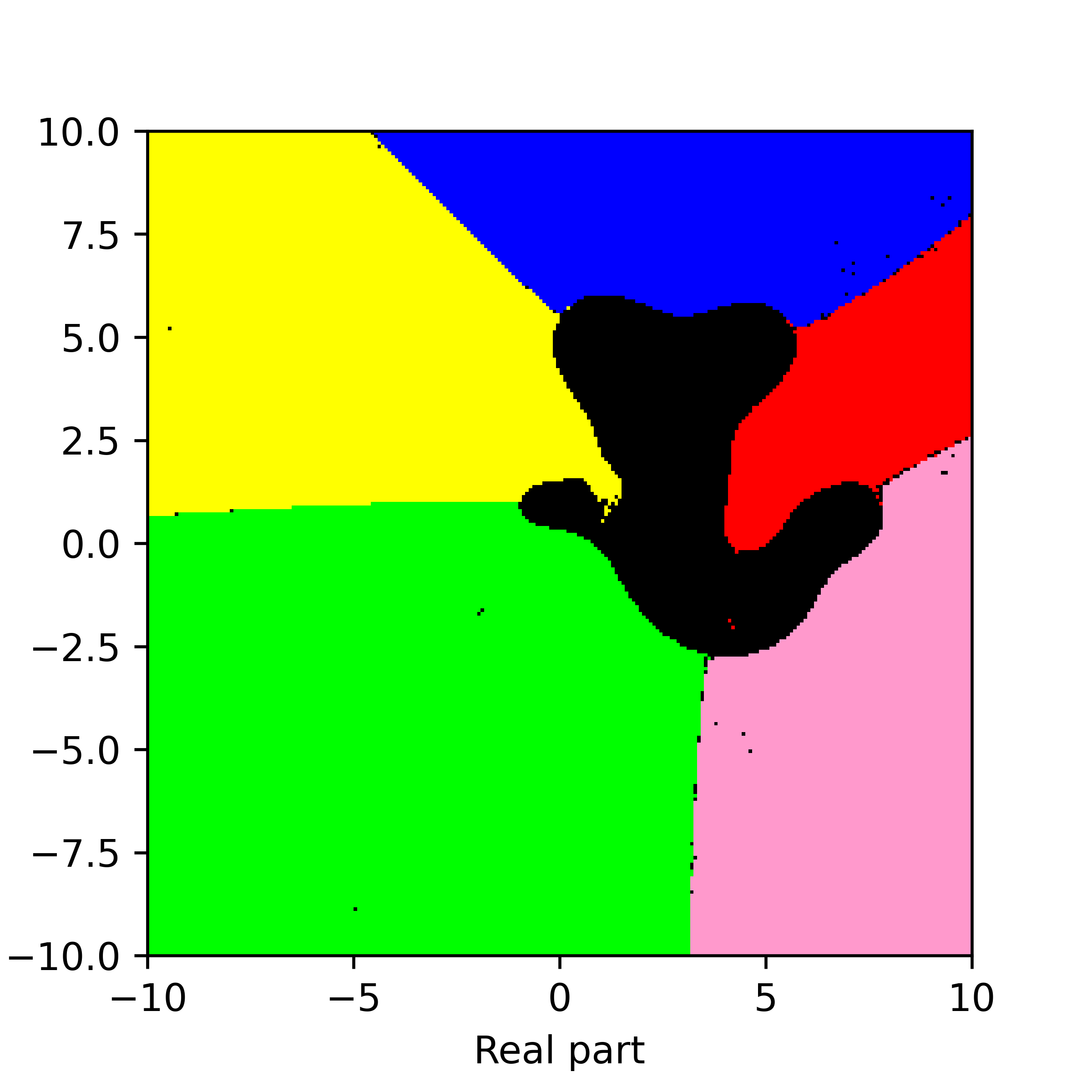}
    
    \caption{Basins of attraction for finding roots of the function $f_{15}$ by different methods. Pictures are referenced to from top to bottom, from left to right. Row 1: left picture is Voronoi's diagram, central picture is for Newton's method, right picture is for Random Relaxed Newton's method. Row 2: left picture is for Newton's method vOptimization, right picture is for BNQN. Row 3: left picture is for Newton's flow, central picture is for Newton's flow vFraction, right picture is for Newton's flow vOptimization. The black points in some of these pictures are those in the basin of attraction of critical points of $f_{15}$.}
    \label{fig:f15}
\end{figure}

\begin{figure}
    \centering
    \includegraphics[width=5cm]{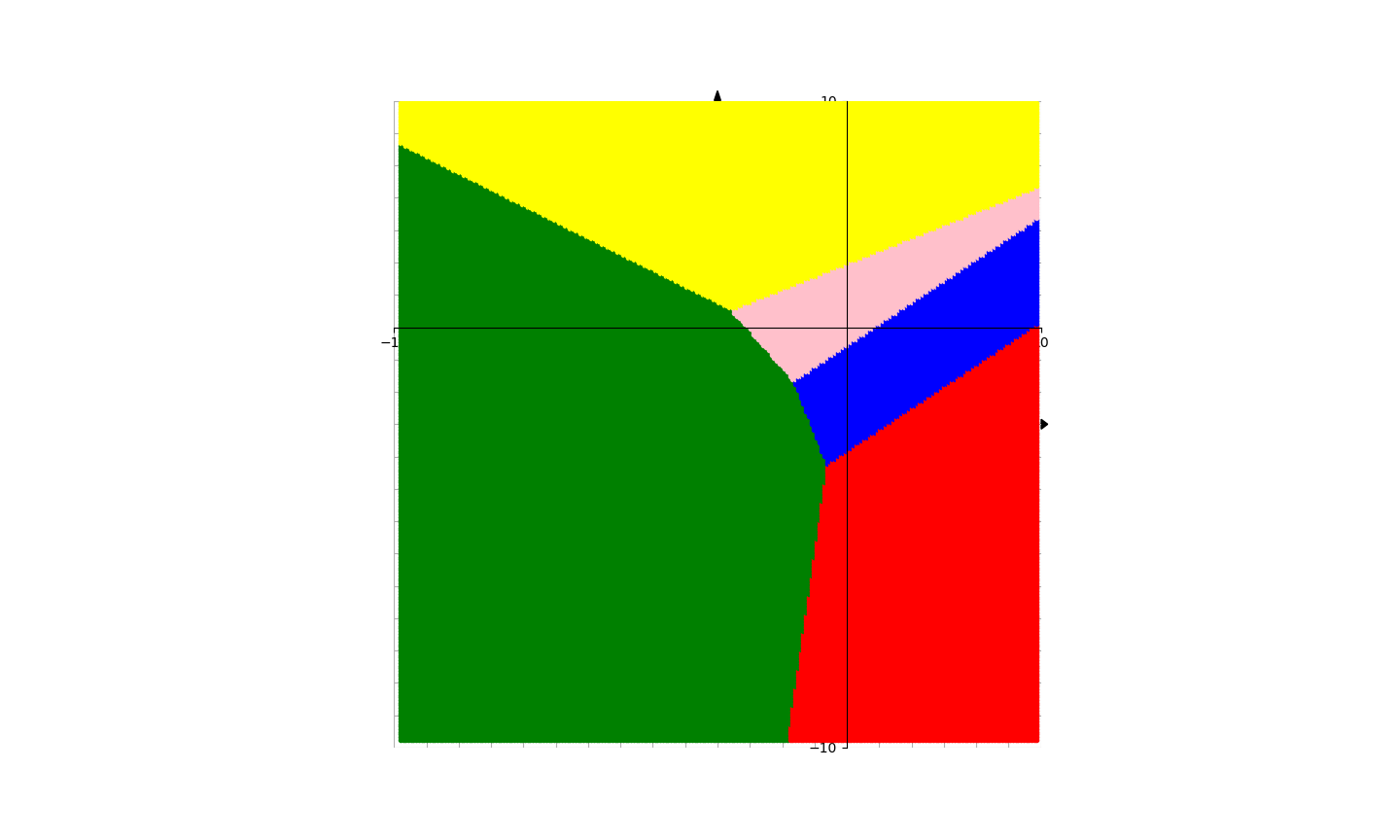}
    \includegraphics[width=3cm]{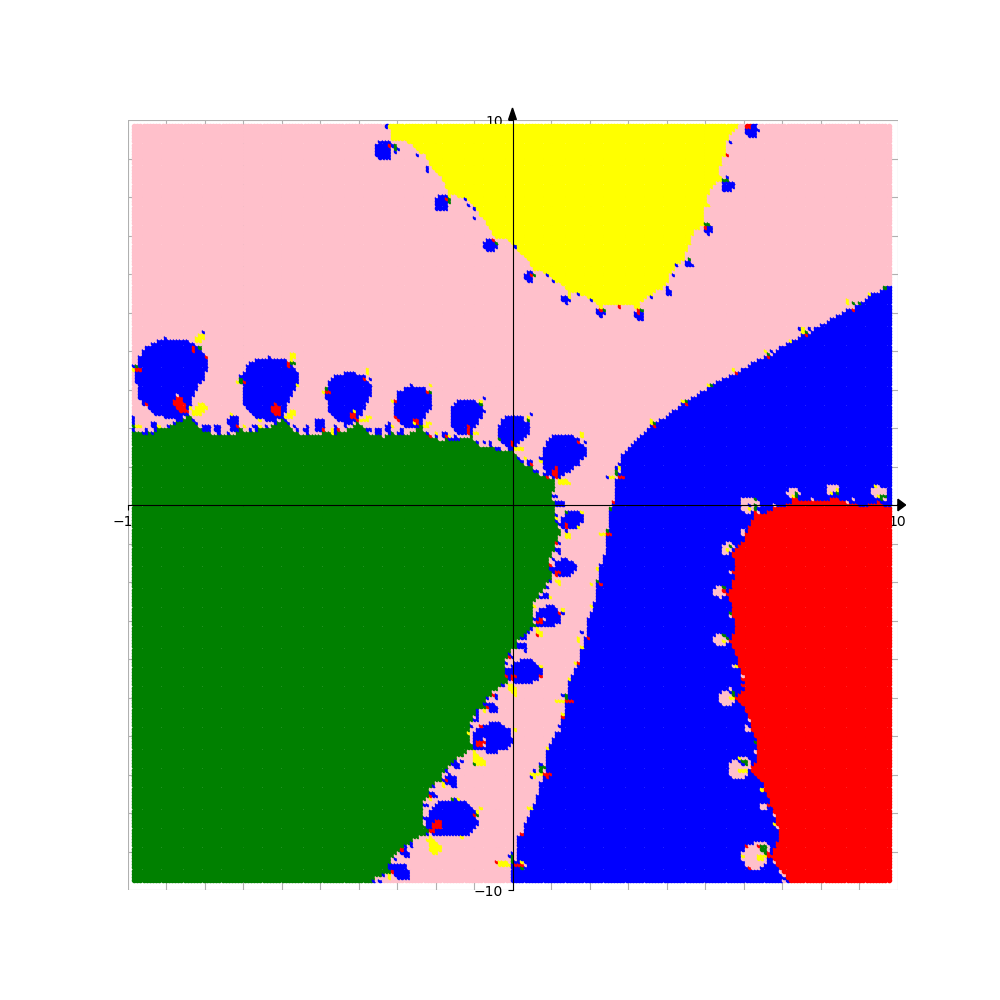}
    \includegraphics[width=3cm]{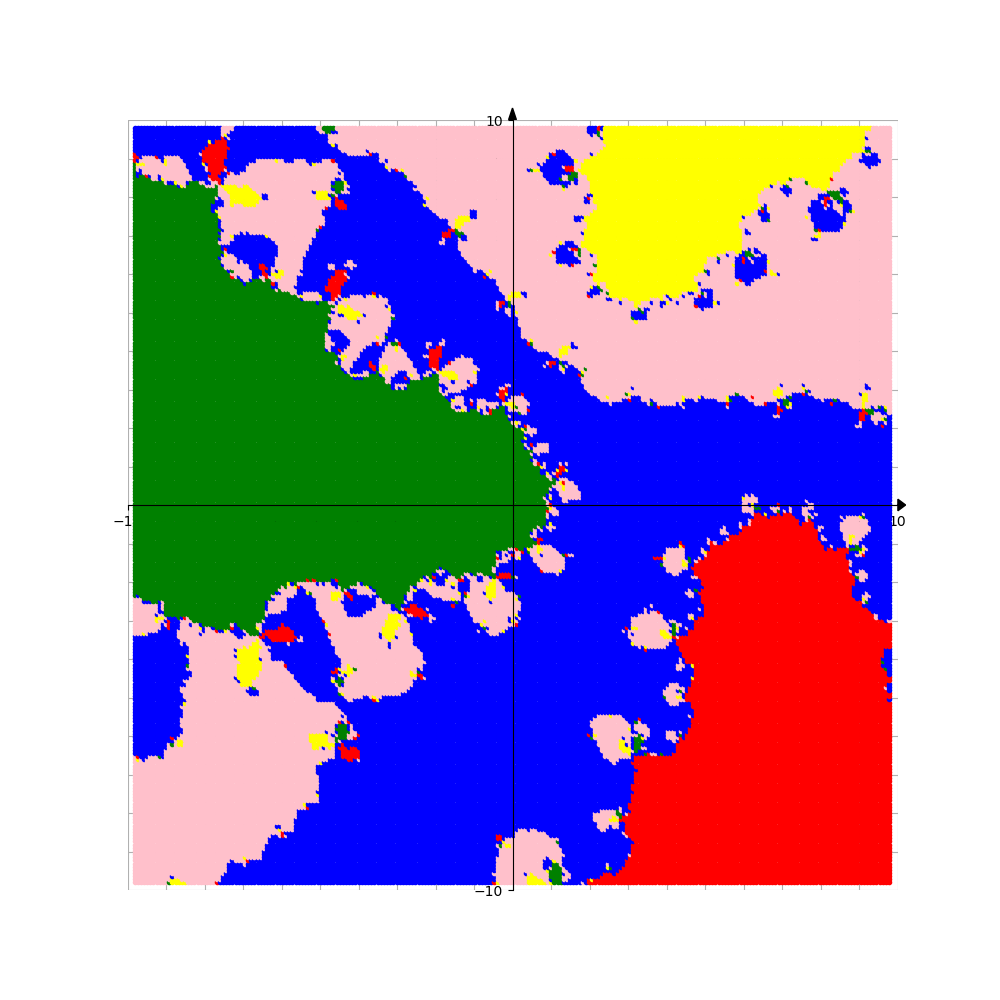}

    \bigskip
    \includegraphics[width=5.5cm]{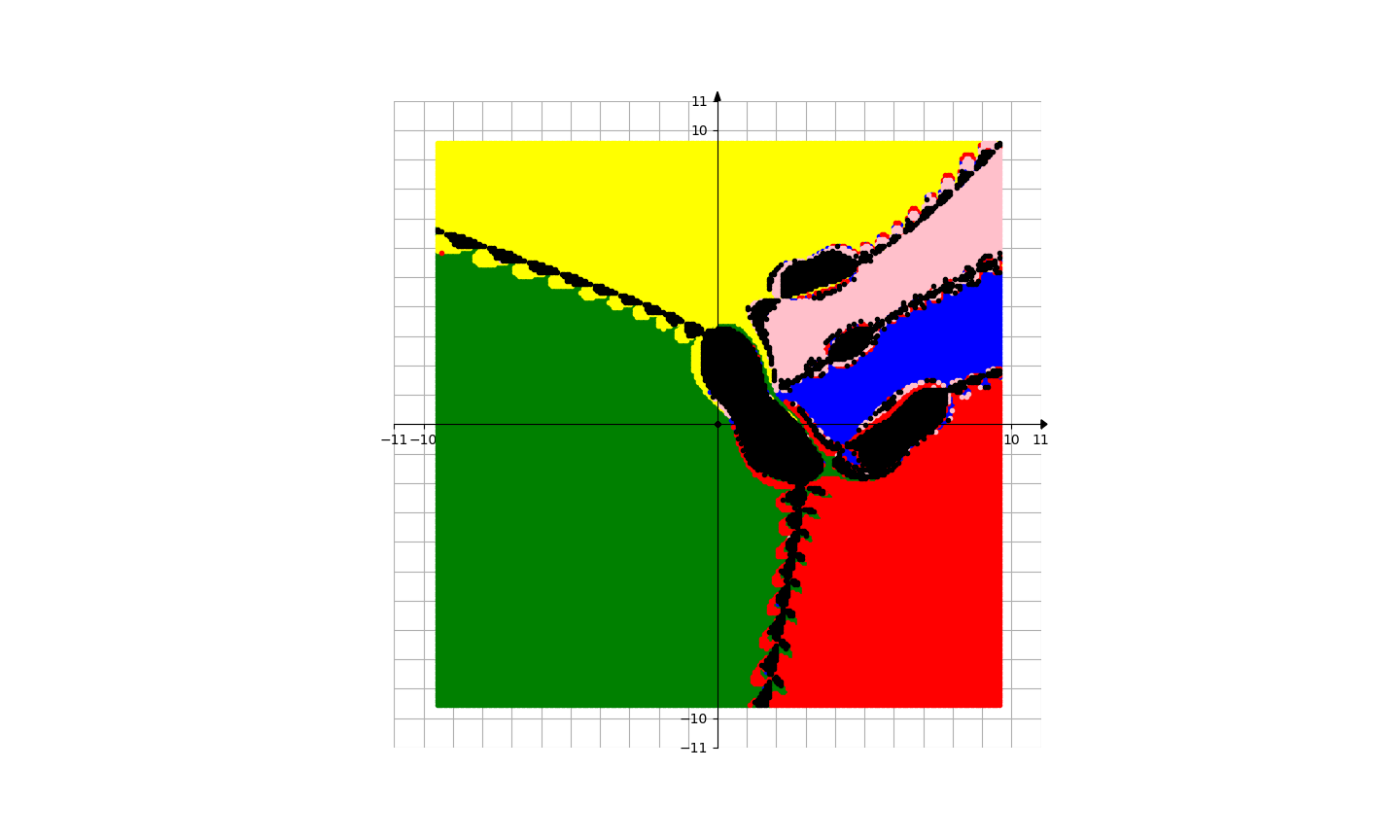}
    \includegraphics[width=3cm]{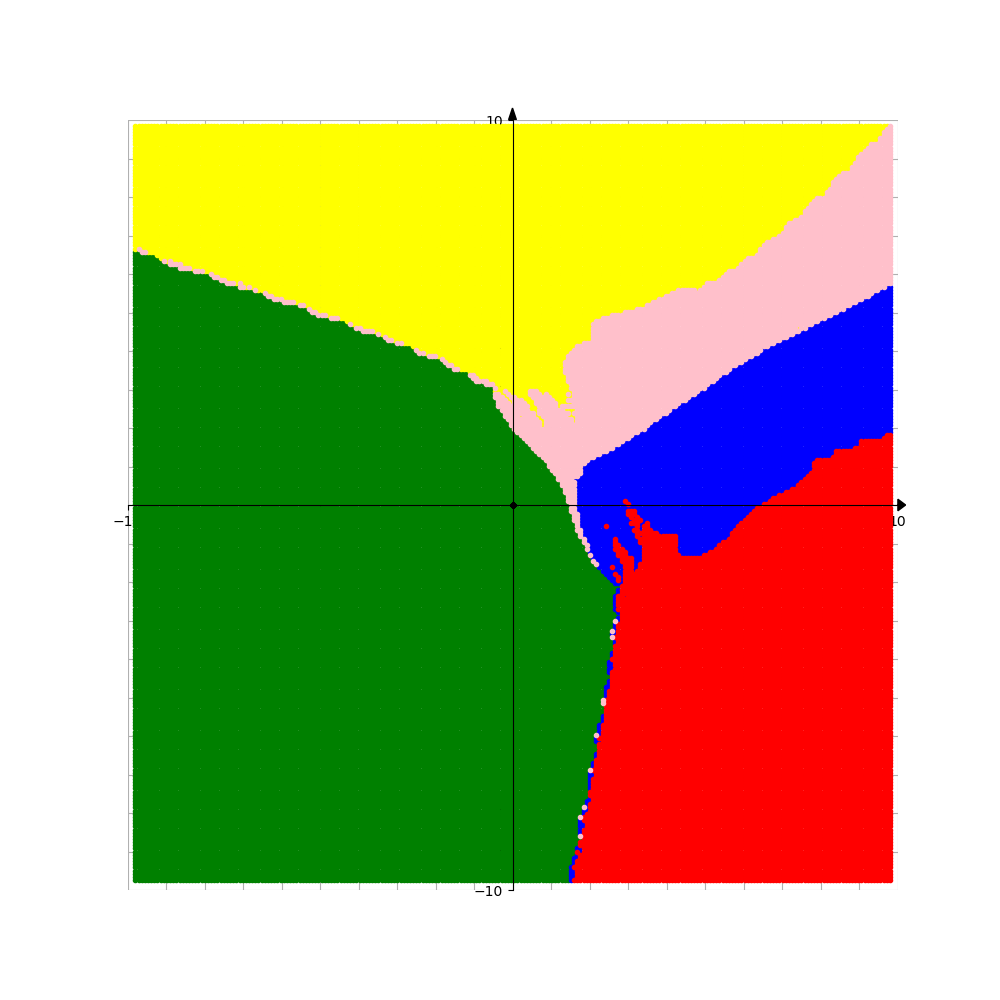}
    
    \bigskip
    \includegraphics[width=3cm]{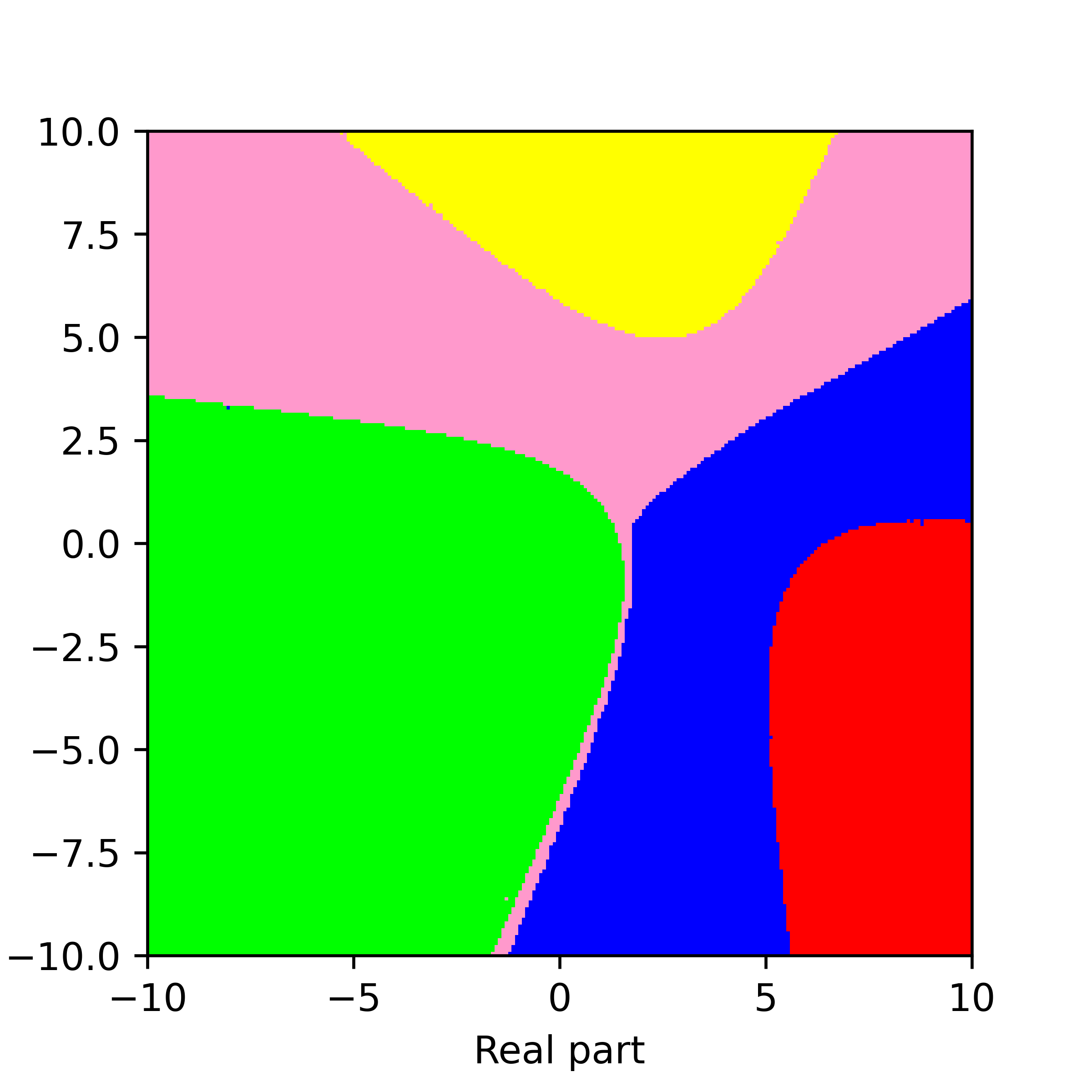}
    \includegraphics[width=3cm]{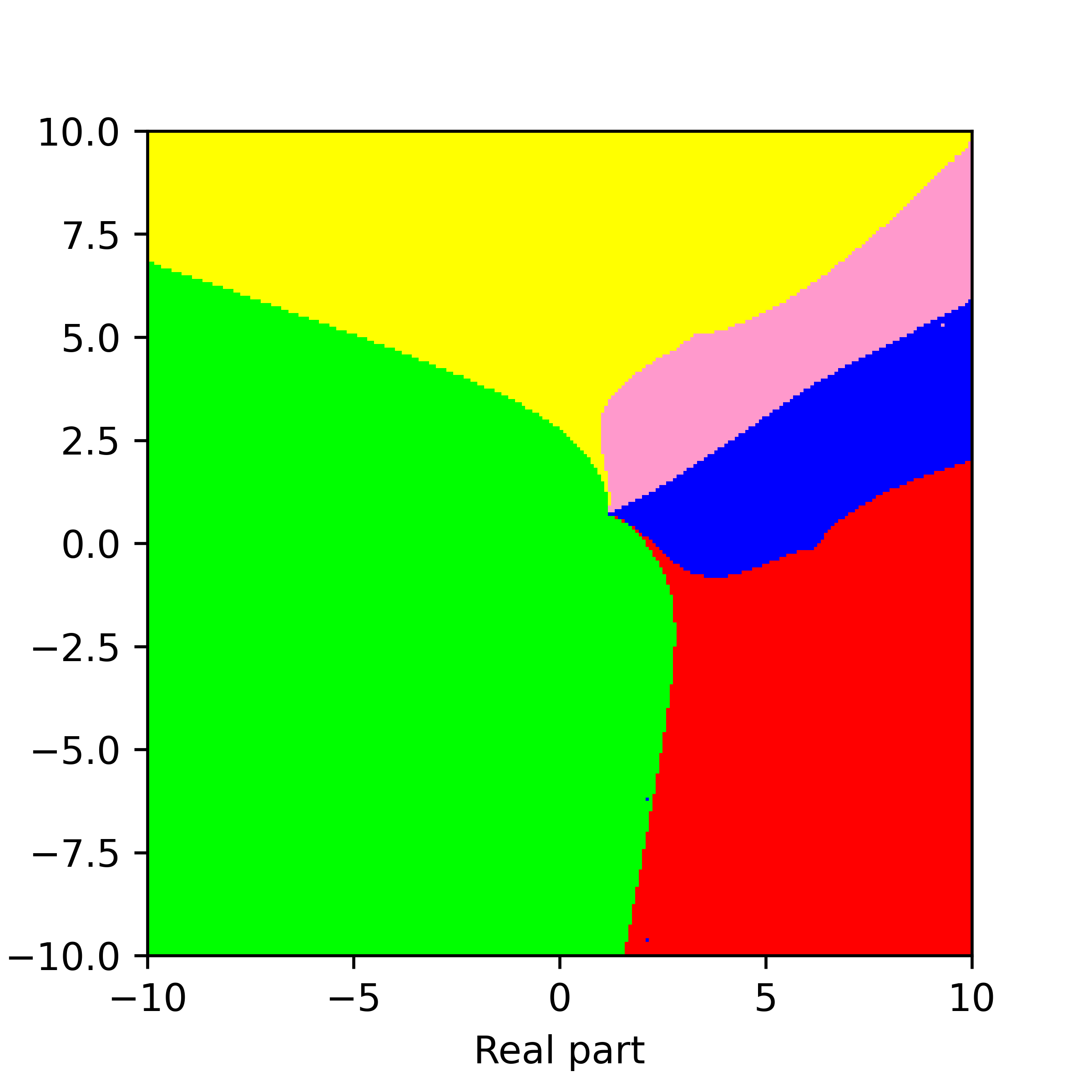}
    \includegraphics[width=3cm]{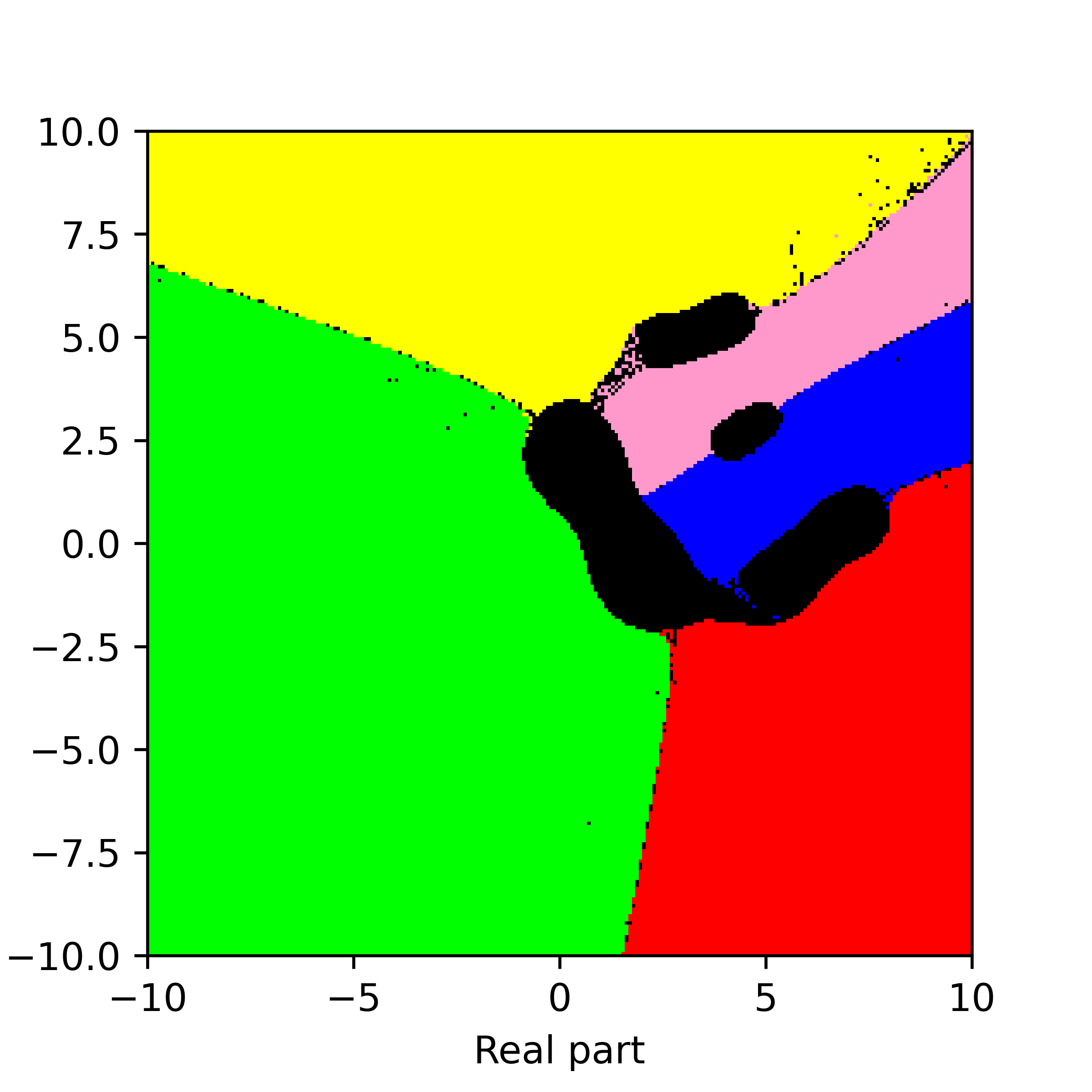}
    
    \caption{Basins of attraction for finding roots of the function $f_{16}$ by different methods. Pictures are referenced to from top to bottom, from left to right. Row 1: left picture is Voronoi's diagram, central picture is for Newton's method, right picture is for Random Relaxed Newton's method. Row 2: left picture is for Newton's method vOptimization, right picture is for BNQN. Row 3: left picture is for Newton's flow, central picture is for Newton's flow vFraction, right picture is for Newton's flow vOptimization. The black points in some of these pictures are those in the basin of attraction of critical points of $f_{16}$.}
    \label{fig:f16}
\end{figure}

\begin{figure}
    \centering
    \includegraphics[width=5cm]{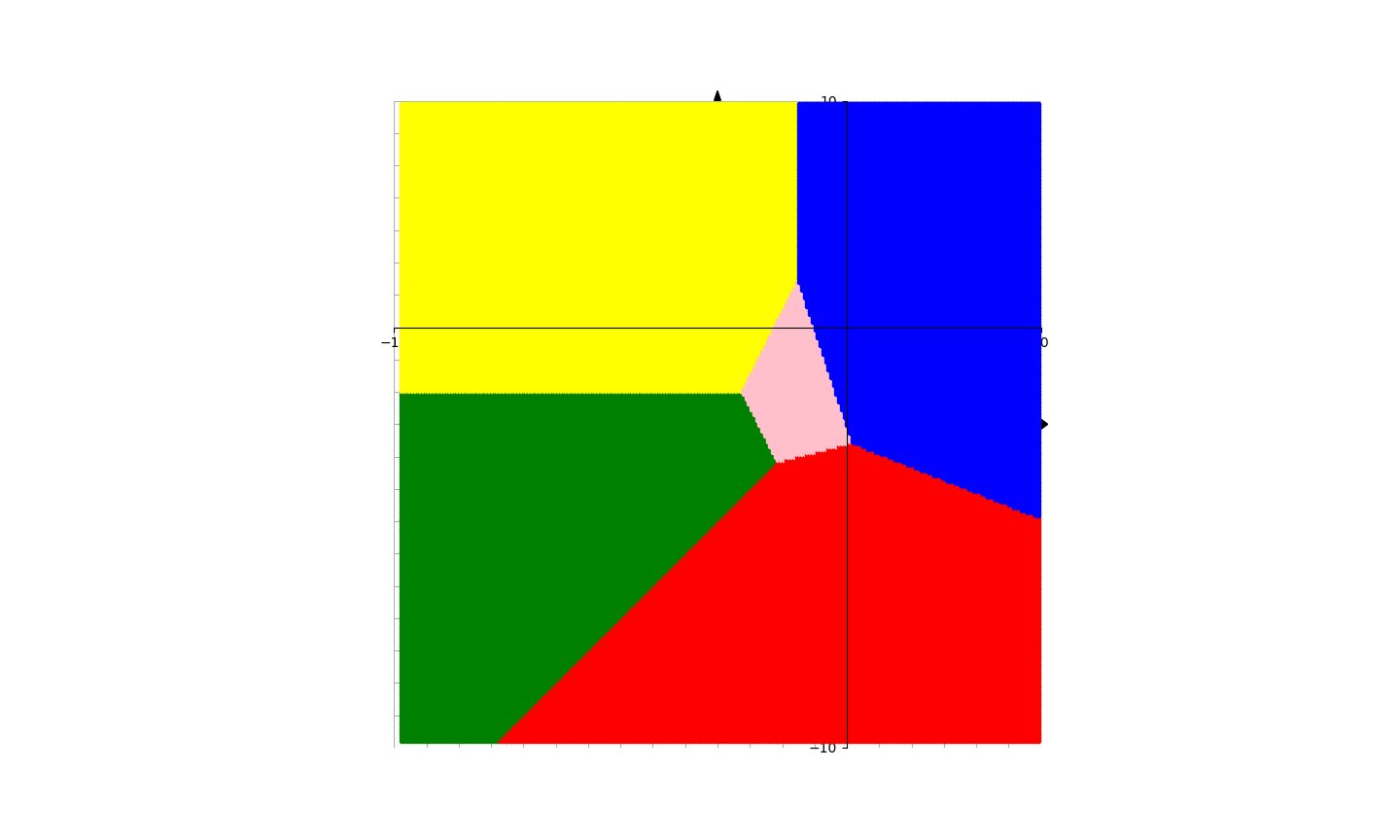}
    \includegraphics[width=3cm]{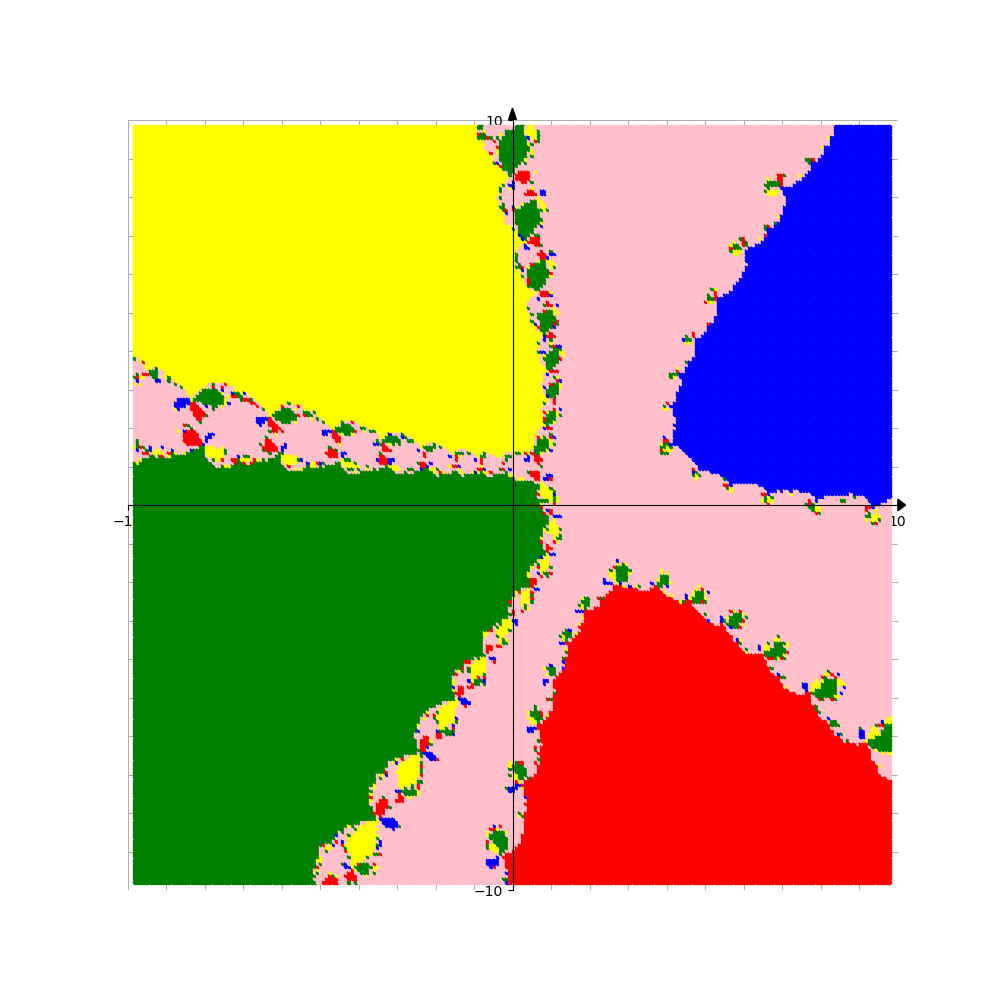}
    \includegraphics[width=3cm]{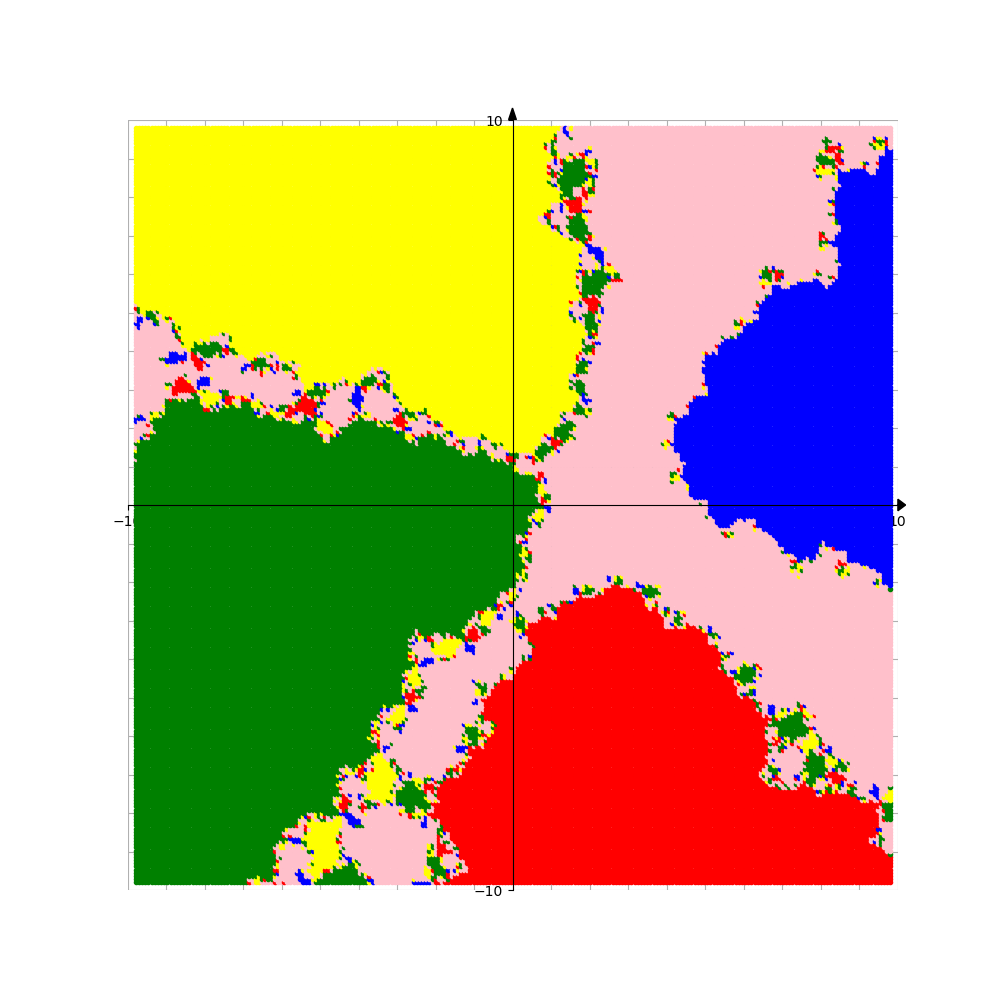}

    \bigskip
    \includegraphics[width=5.5cm]{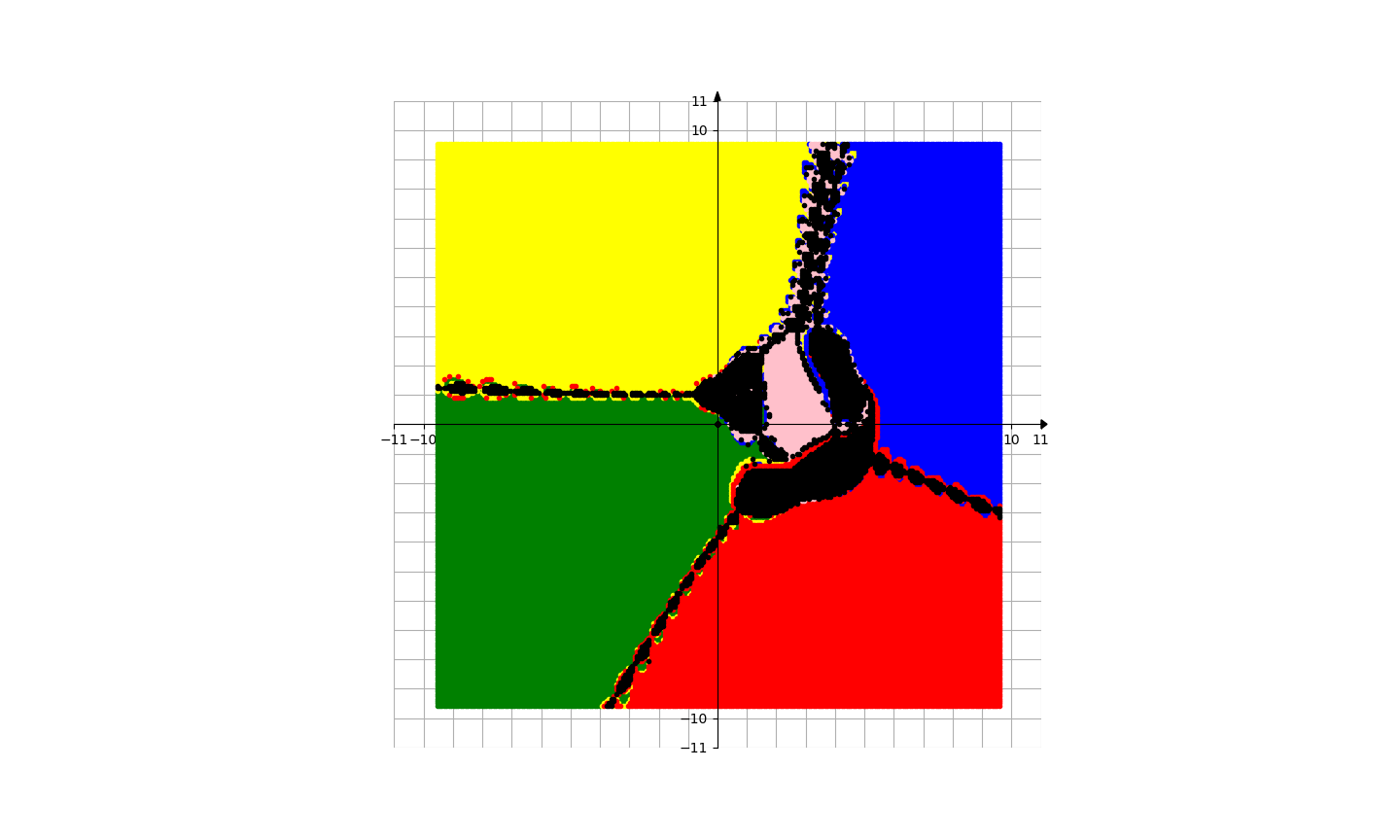}
    \includegraphics[width=3cm]{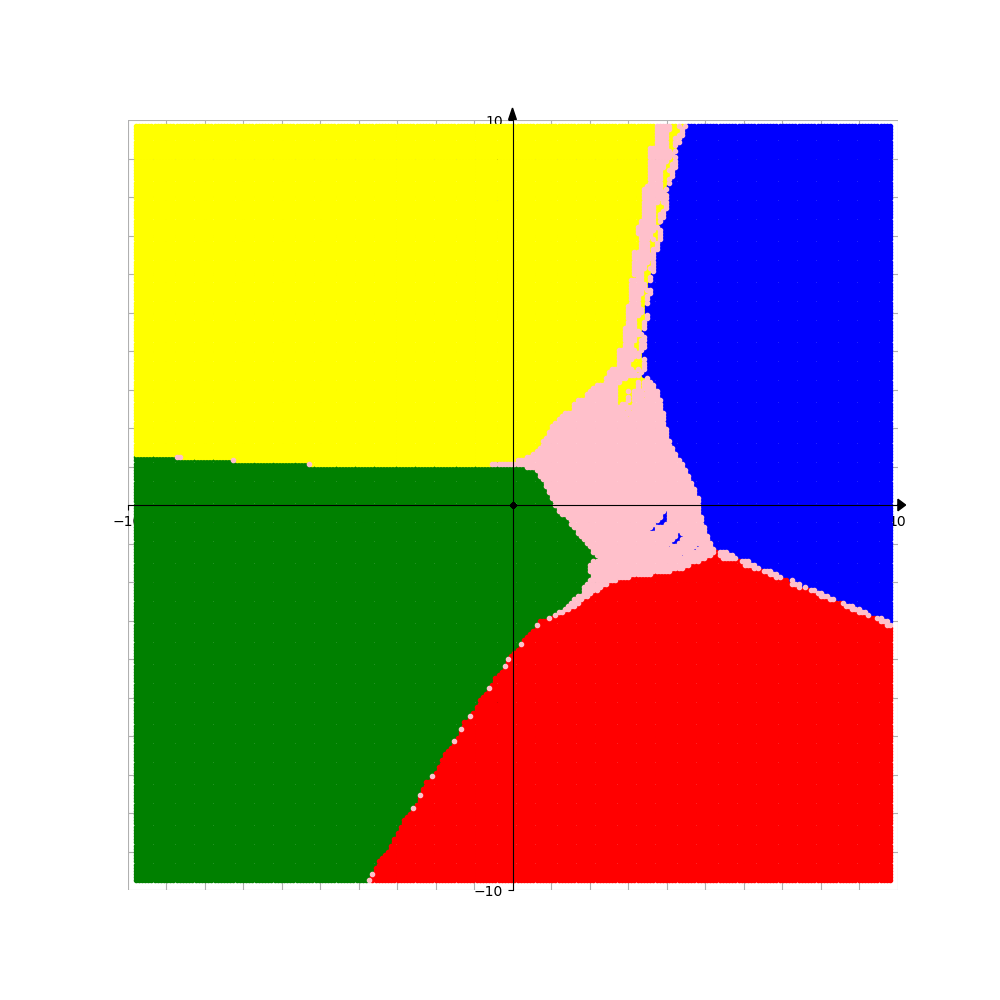}
    
    \bigskip
    \includegraphics[width=3cm]{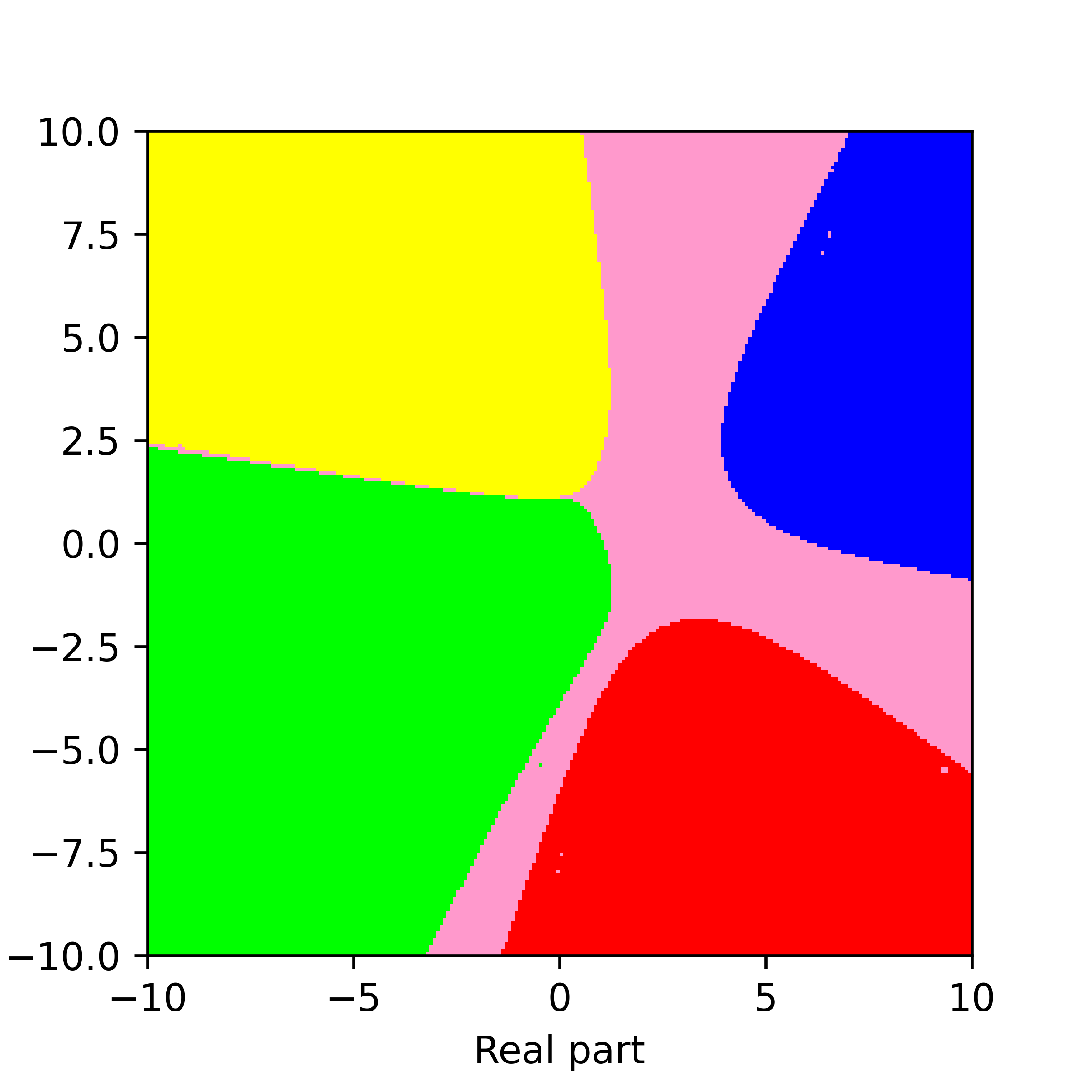}
    \includegraphics[width=3cm]{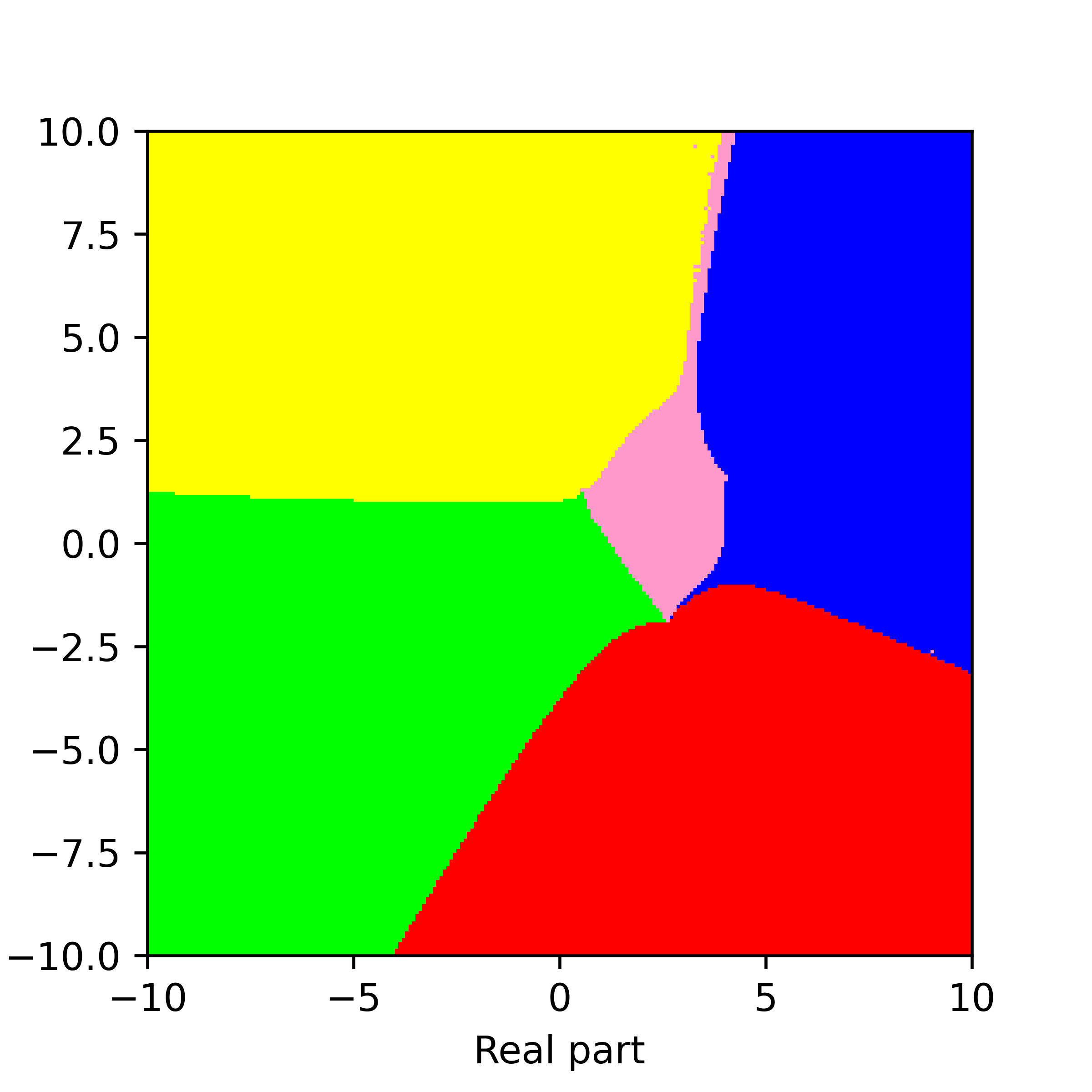}
    \includegraphics[width=3cm]{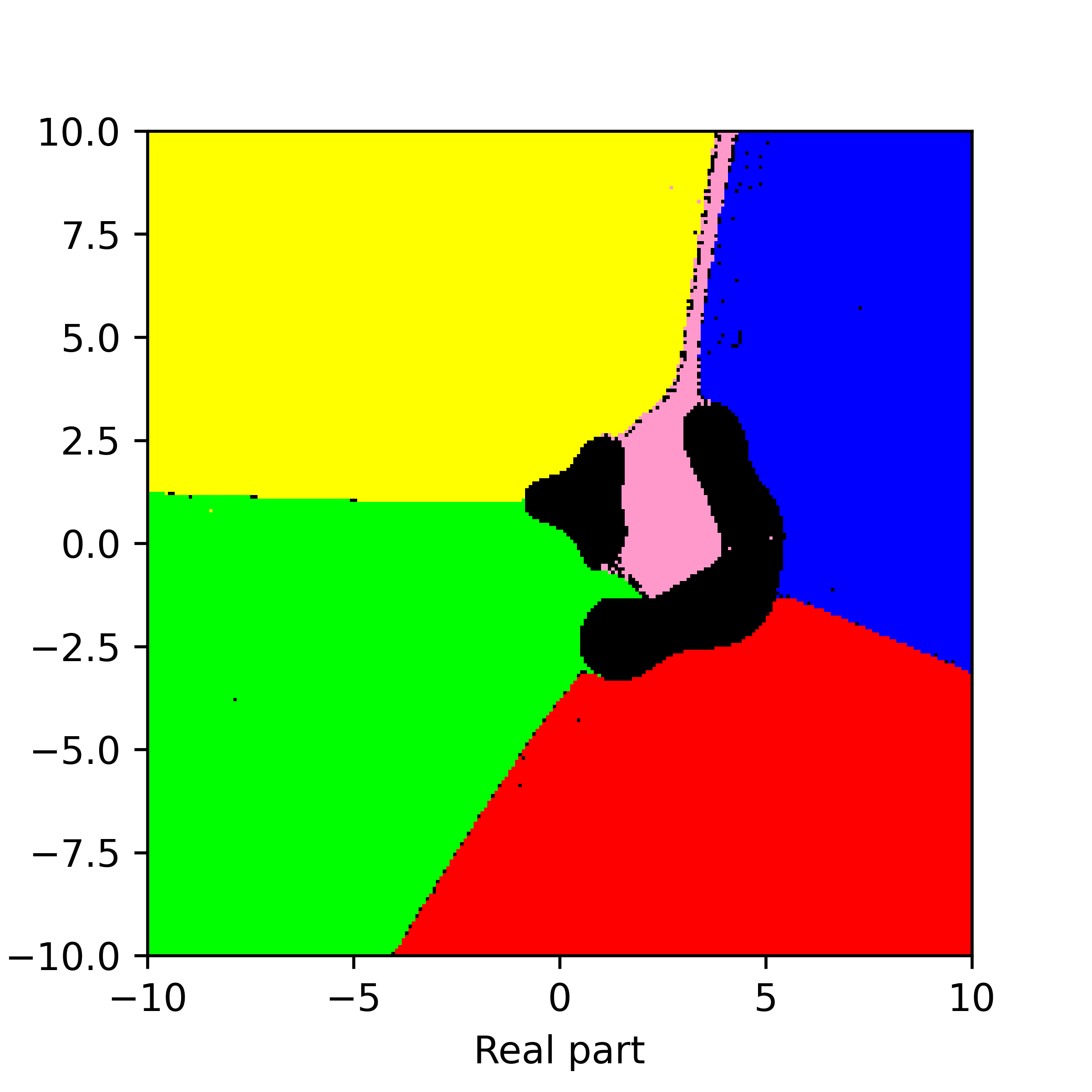}
    
    \caption{Basins of attraction for finding roots of the function $f_{17}$ by different methods. Pictures are referenced to from top to bottom, from left to right. Row 1: left picture is Voronoi's diagram, central picture is for Newton's method, right picture is for Random Relaxed Newton's method. Row 2: left picture is for Newton's method vOptimization, right picture is for BNQN. Row 3: left picture is for Newton's flow, central picture is for Newton's flow vFraction, right picture is for Newton's flow vOptimization. The black points in some of these pictures are those in the basin of attraction of critical points of $f_{17}$.}
    \label{fig:f17}
\end{figure}

\begin{figure}
    \centering
    \includegraphics[width=5cm]{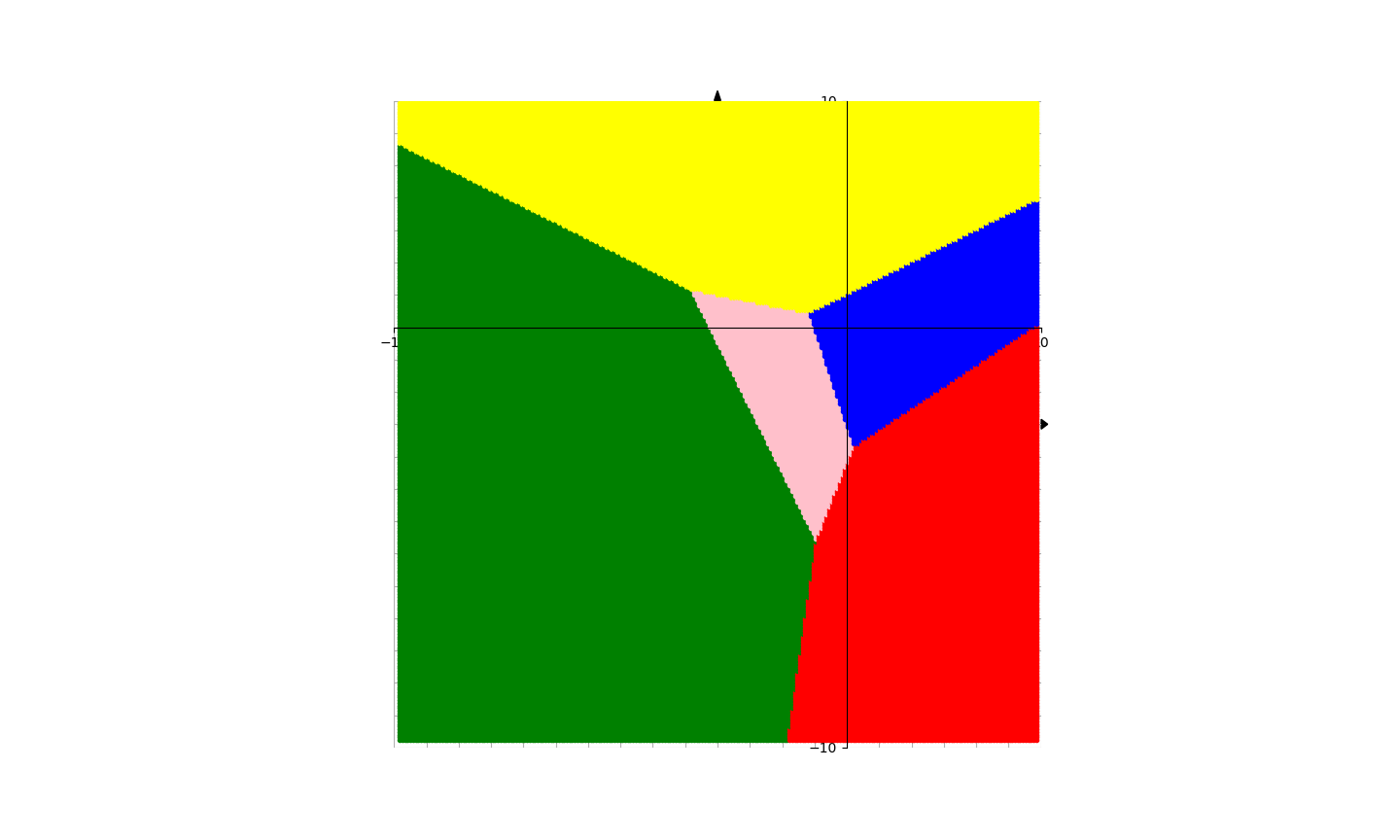}
    \includegraphics[width=3cm]{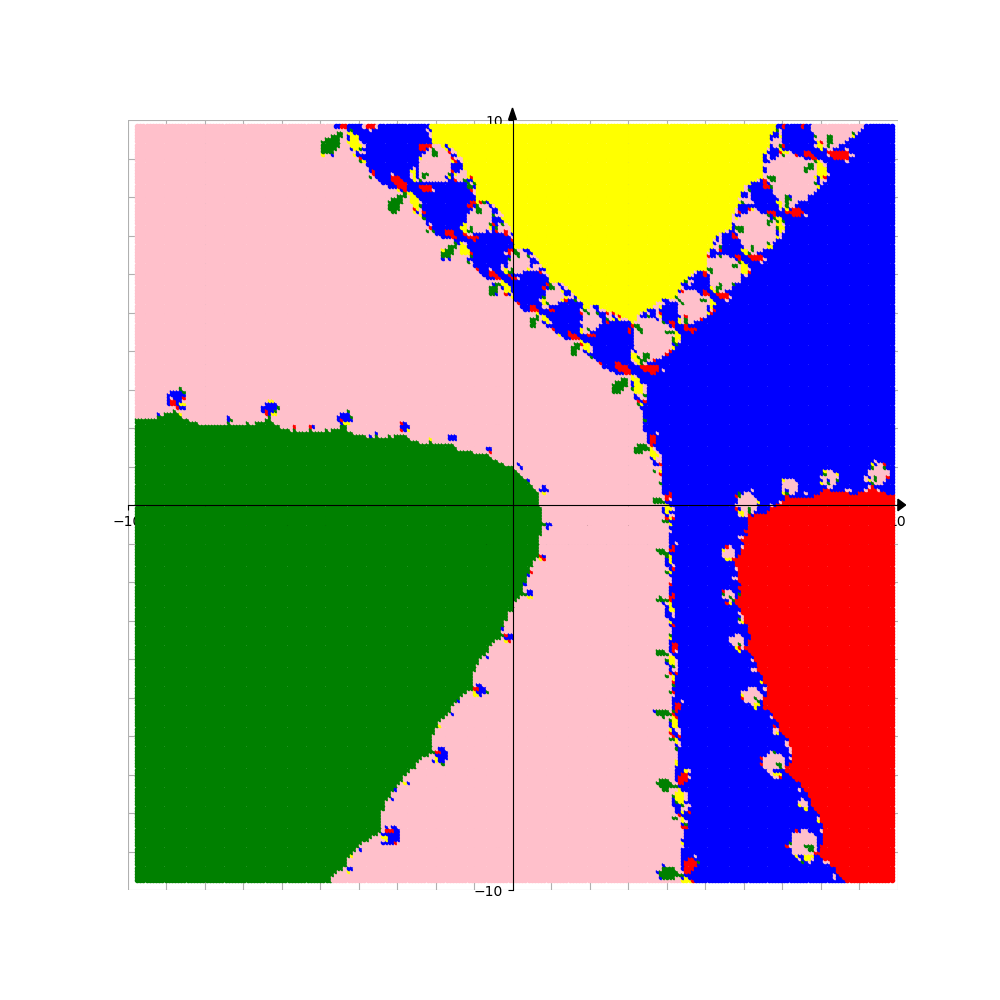}
    \includegraphics[width=3cm]{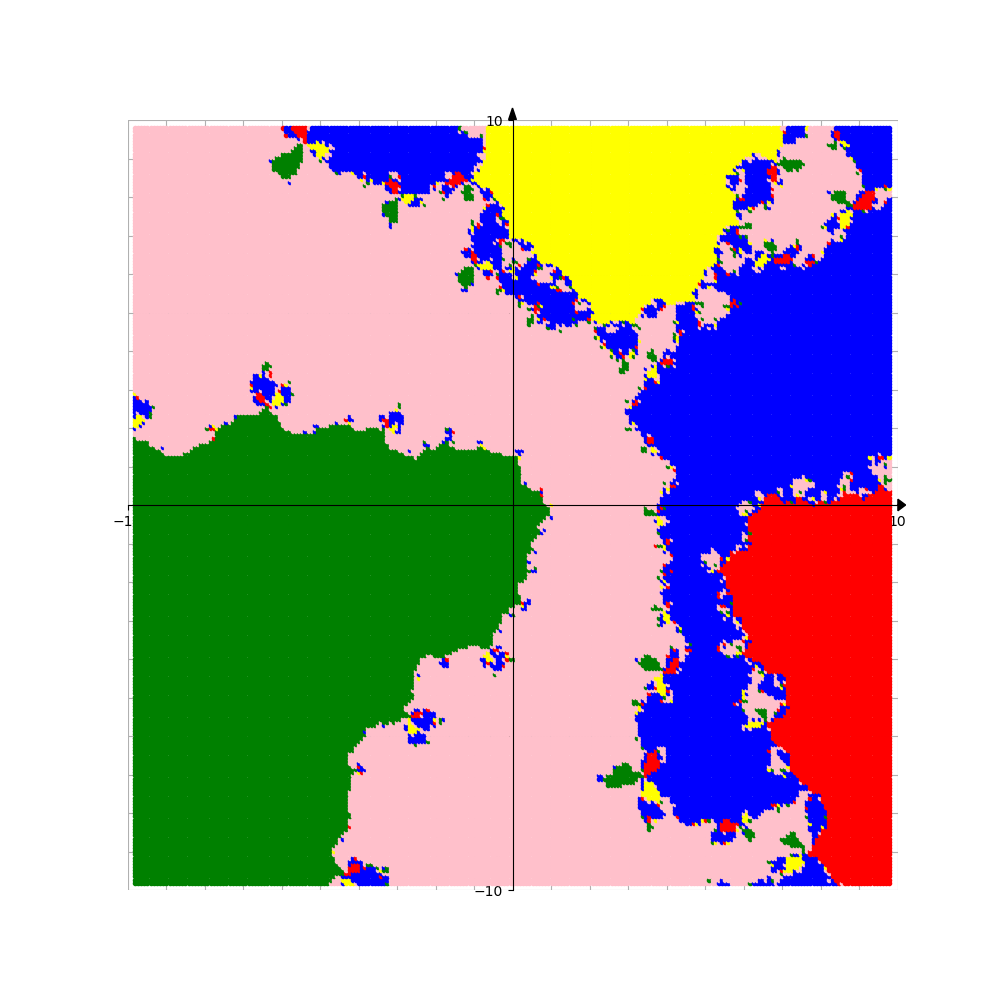}

    \bigskip
    \includegraphics[width=5.5cm]{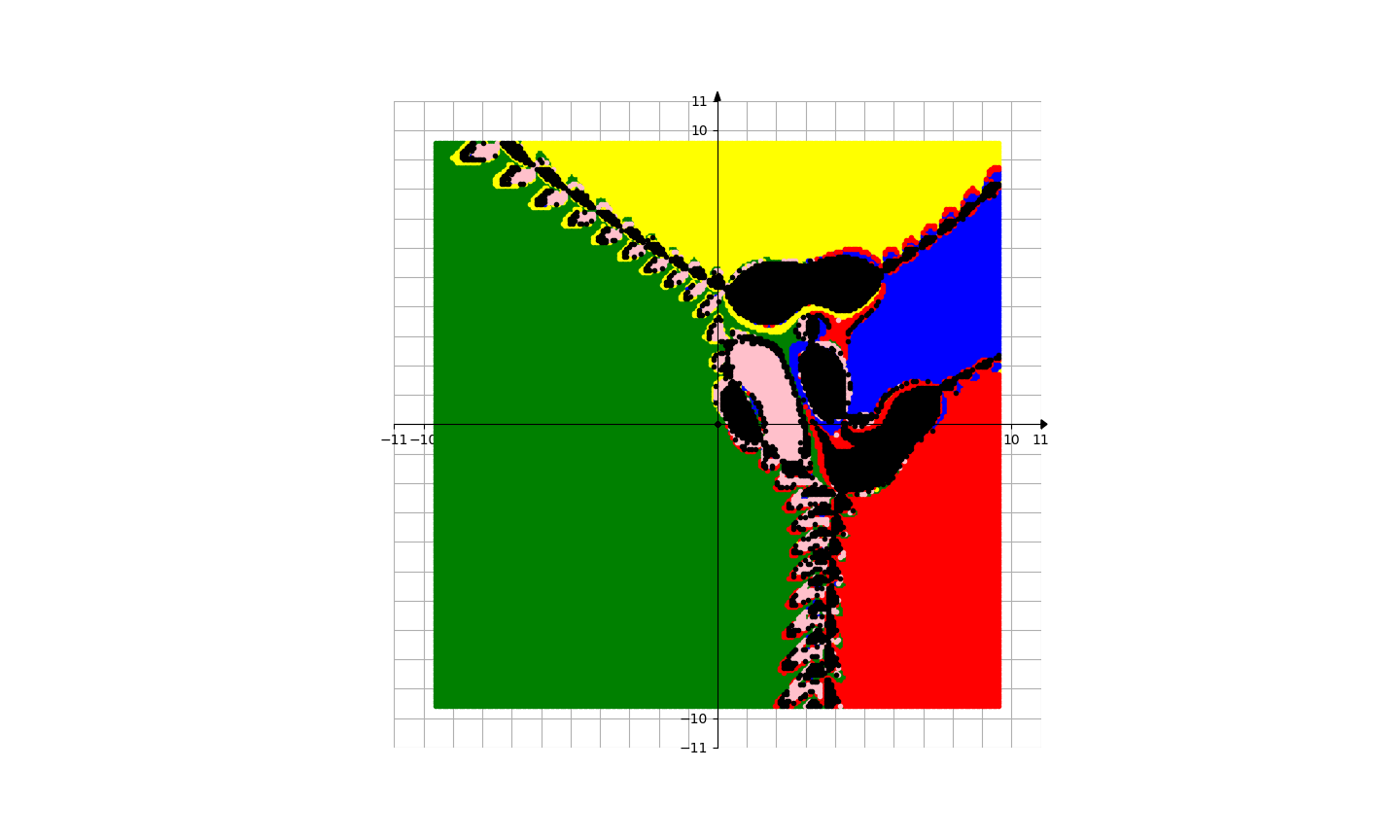}
    \includegraphics[width=3cm]{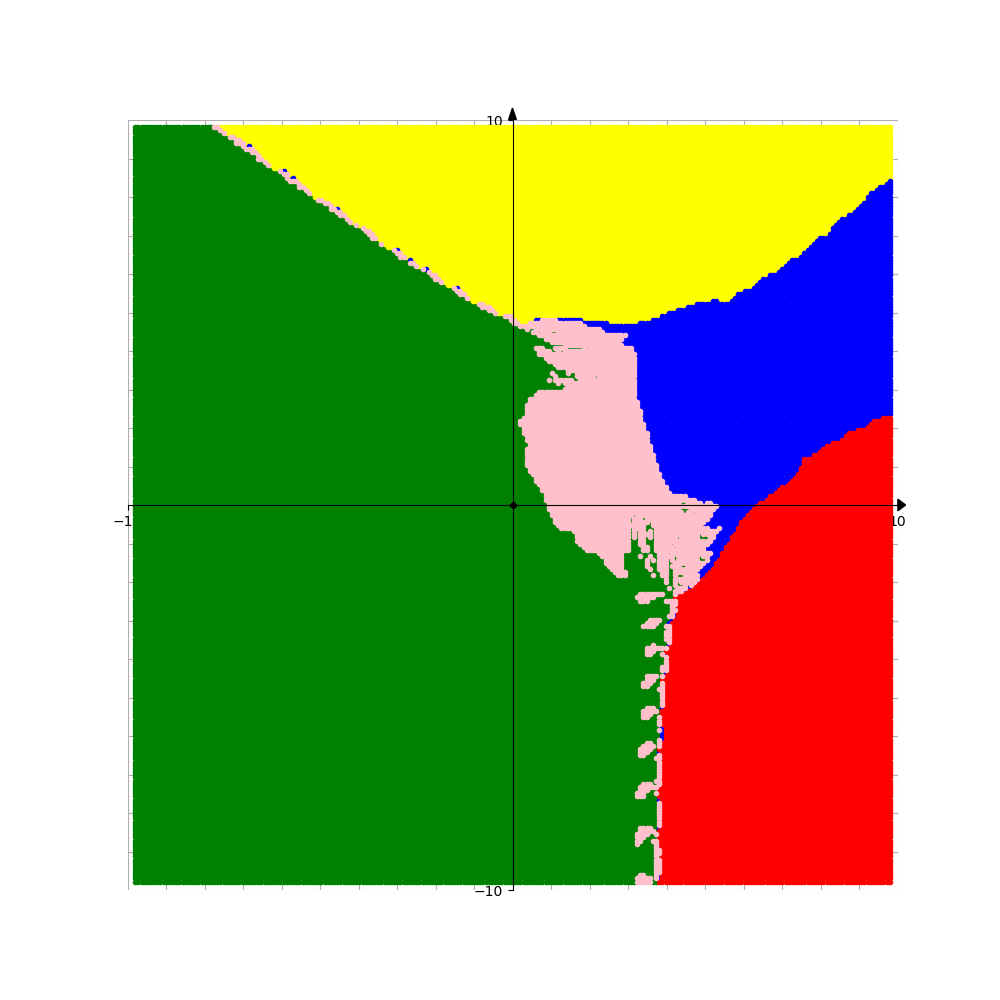}
    
    \bigskip
    \includegraphics[width=3cm]{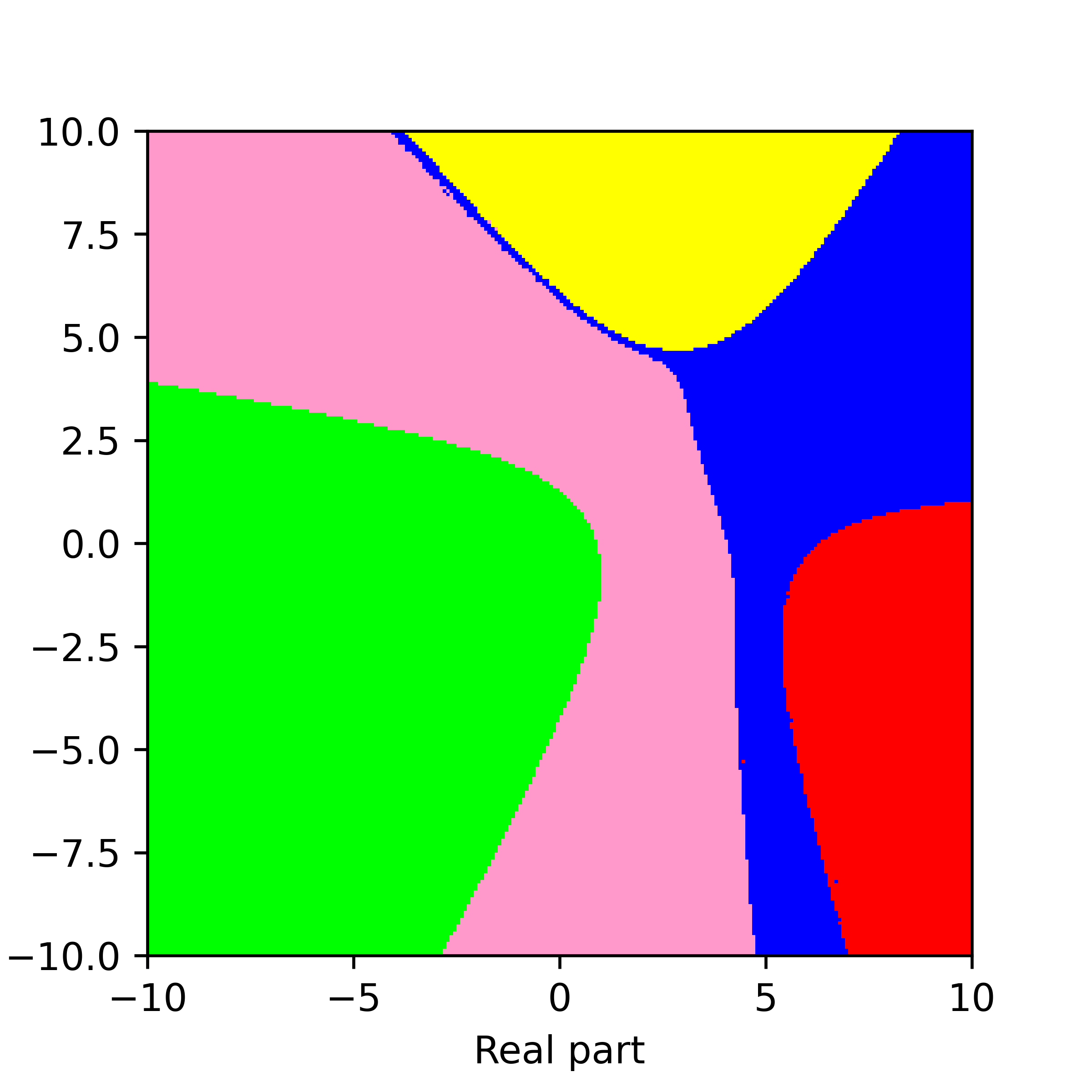}
    \includegraphics[width=3cm]{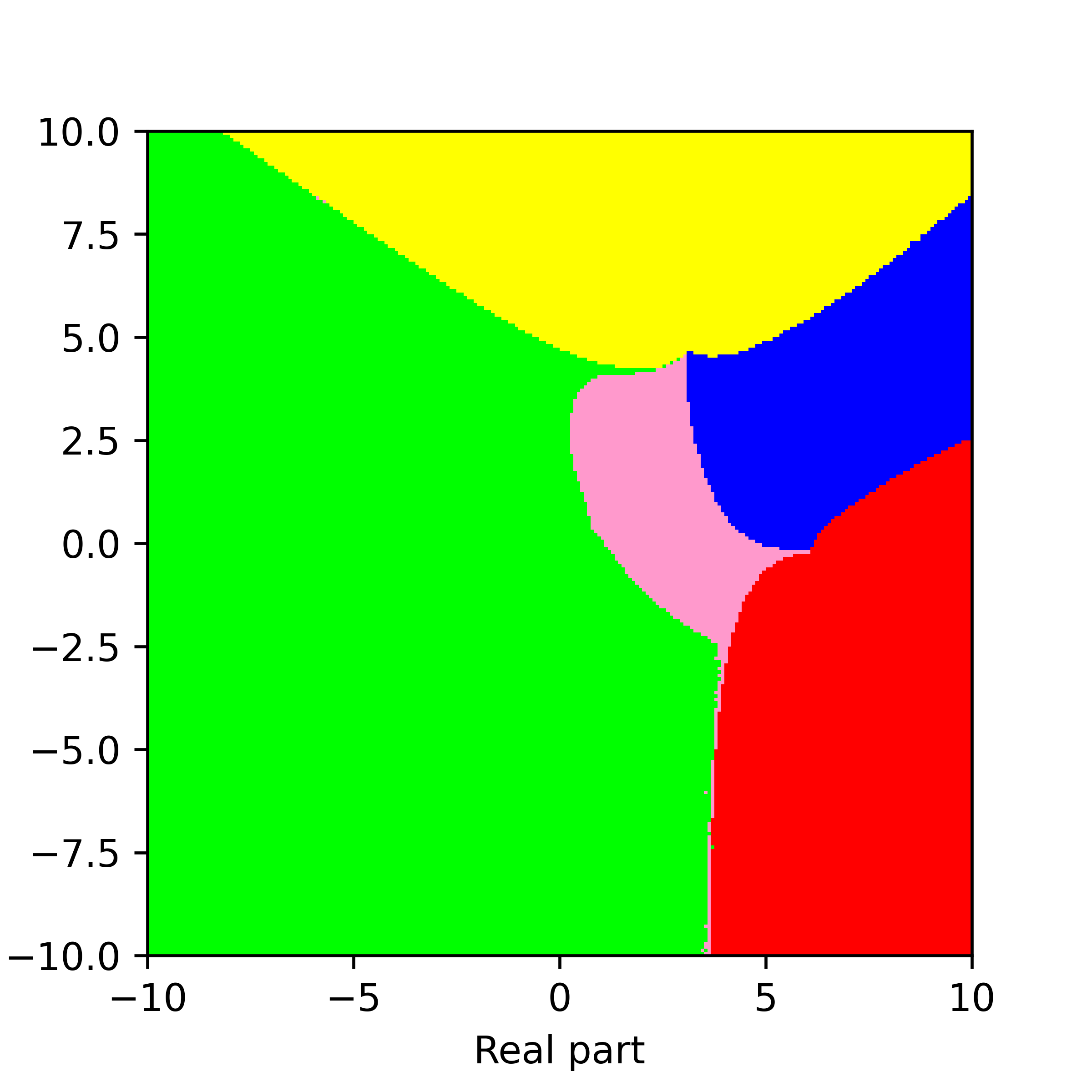}
    \includegraphics[width=3cm]{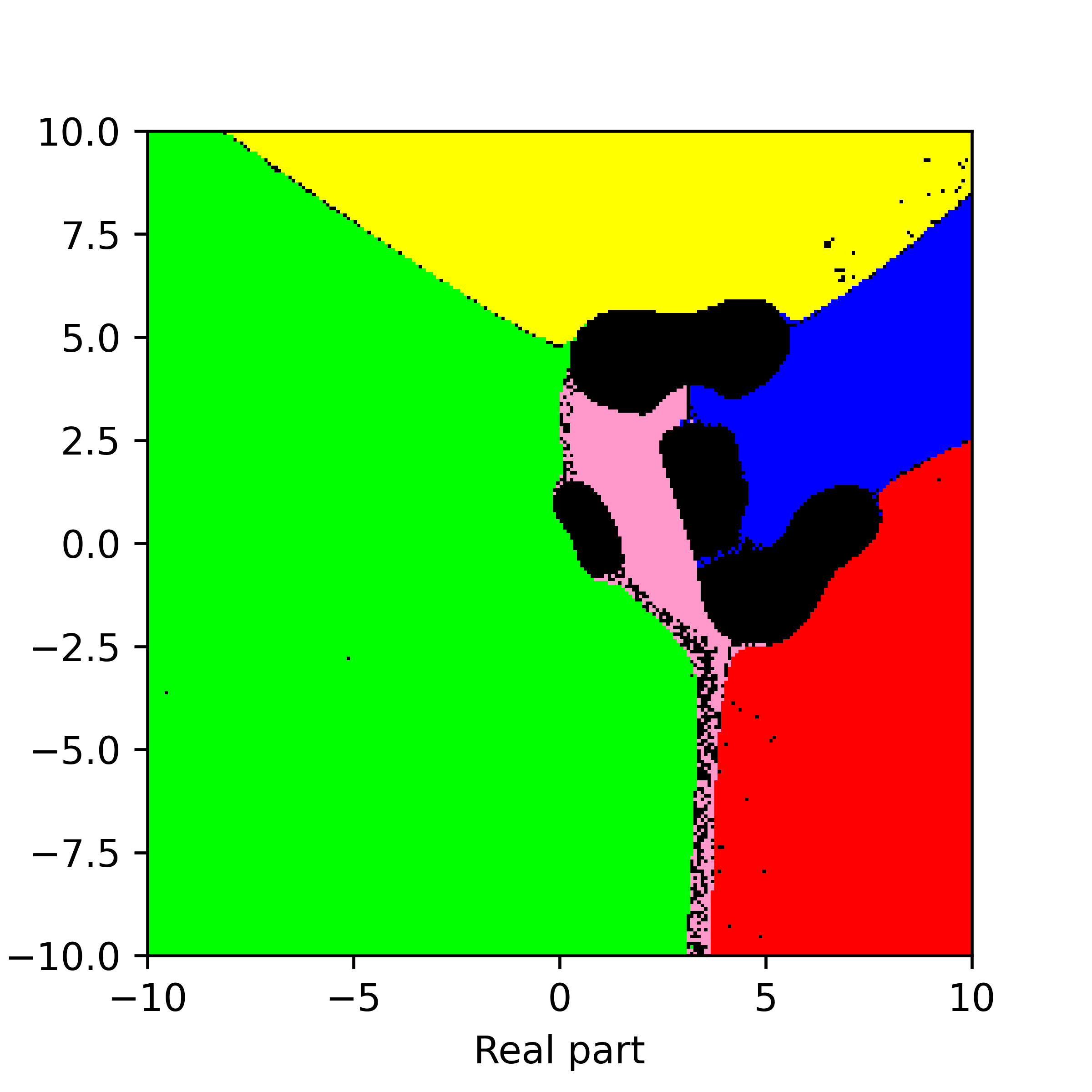}
    
    \caption{Basins of attraction for finding roots of the function $f_{18}$ by different methods. Pictures are referenced to from top to bottom, from left to right. Row 1: left picture is Voronoi's diagram, central picture is for Newton's method, right picture is for Random Relaxed Newton's method. Row 2: left picture is for Newton's method vOptimization, right picture is for BNQN. Row 3: left picture is for Newton's flow, central picture is for Newton's flow vFraction, right picture is for Newton's flow vOptimization. The black points in some of these pictures are those in the basin of attraction of critical points of $f_{18}$.}
    \label{fig:f18}
\end{figure}

\begin{figure}
    \centering
    \includegraphics[width=5cm]{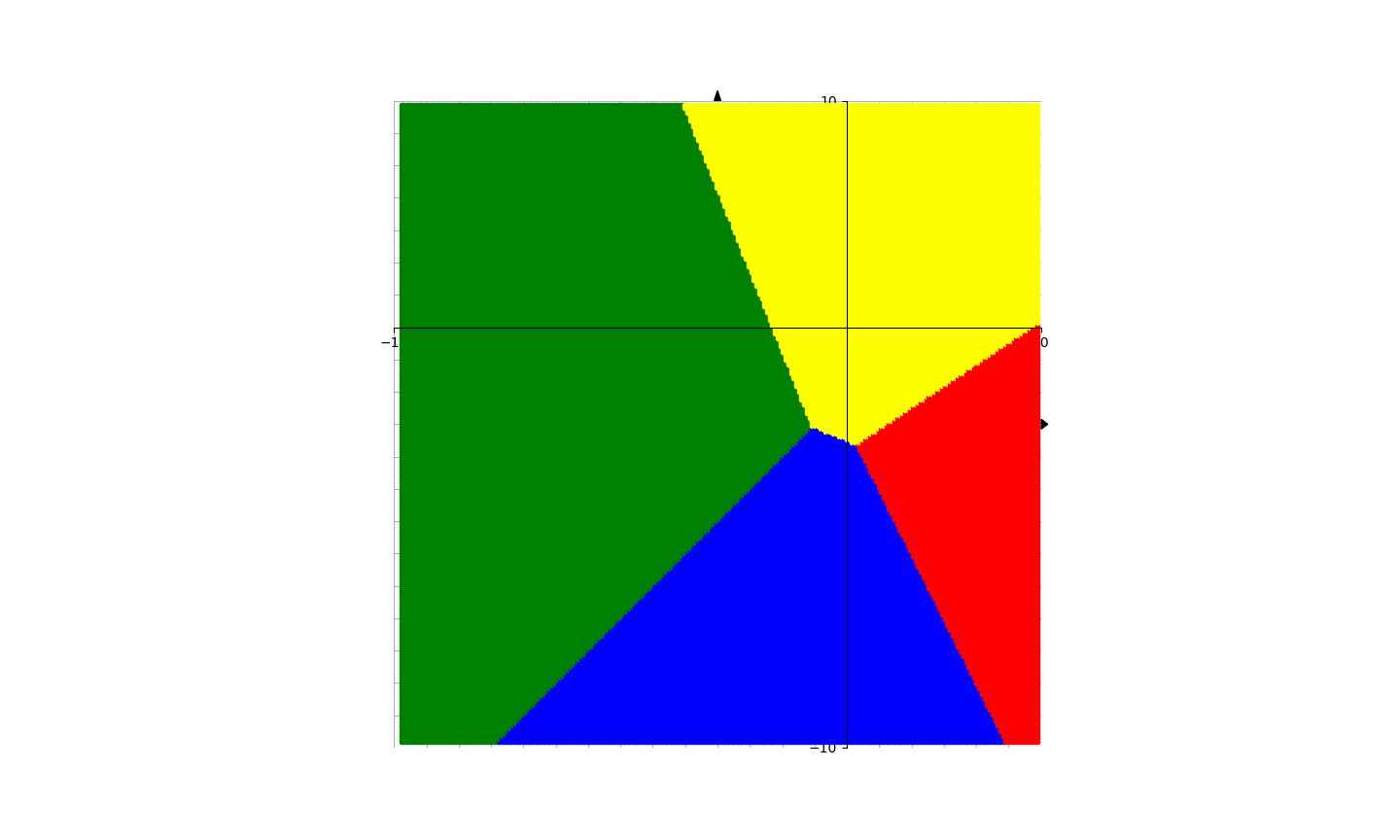}
    \includegraphics[width=3cm]{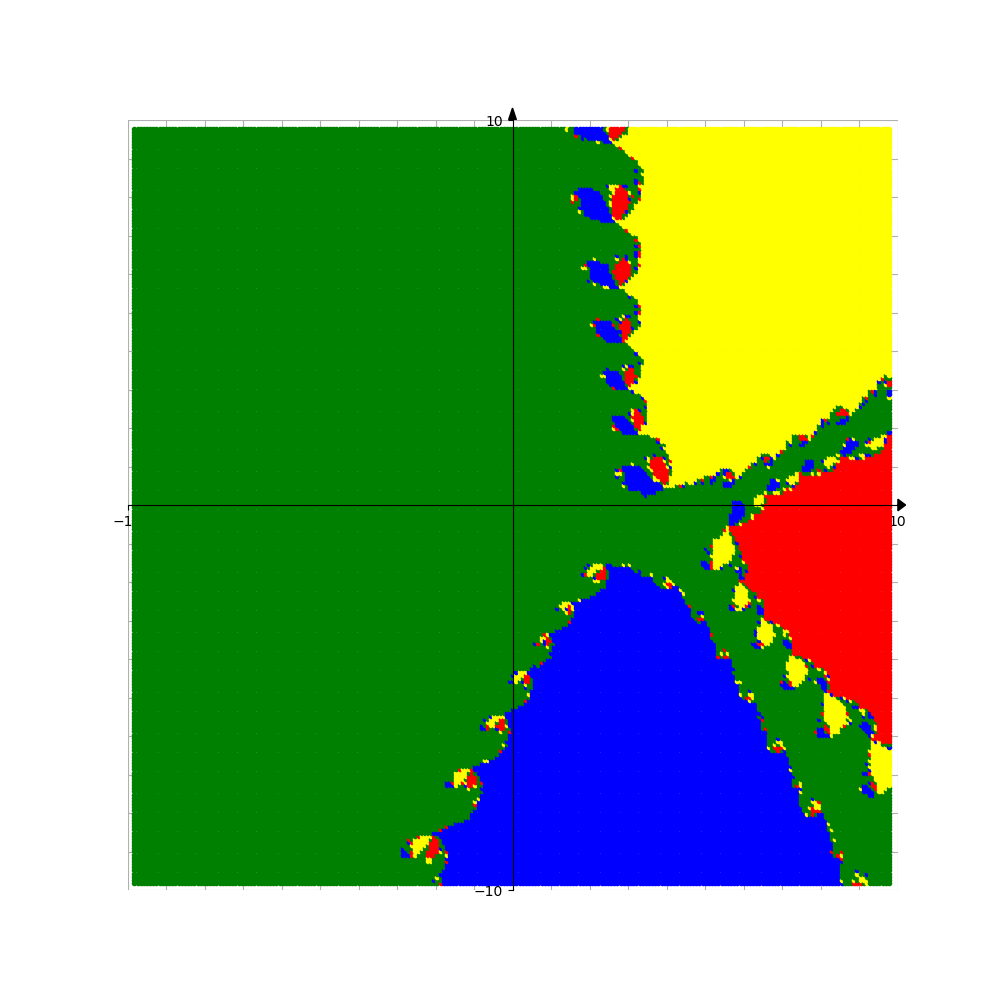}
    \includegraphics[width=3cm]{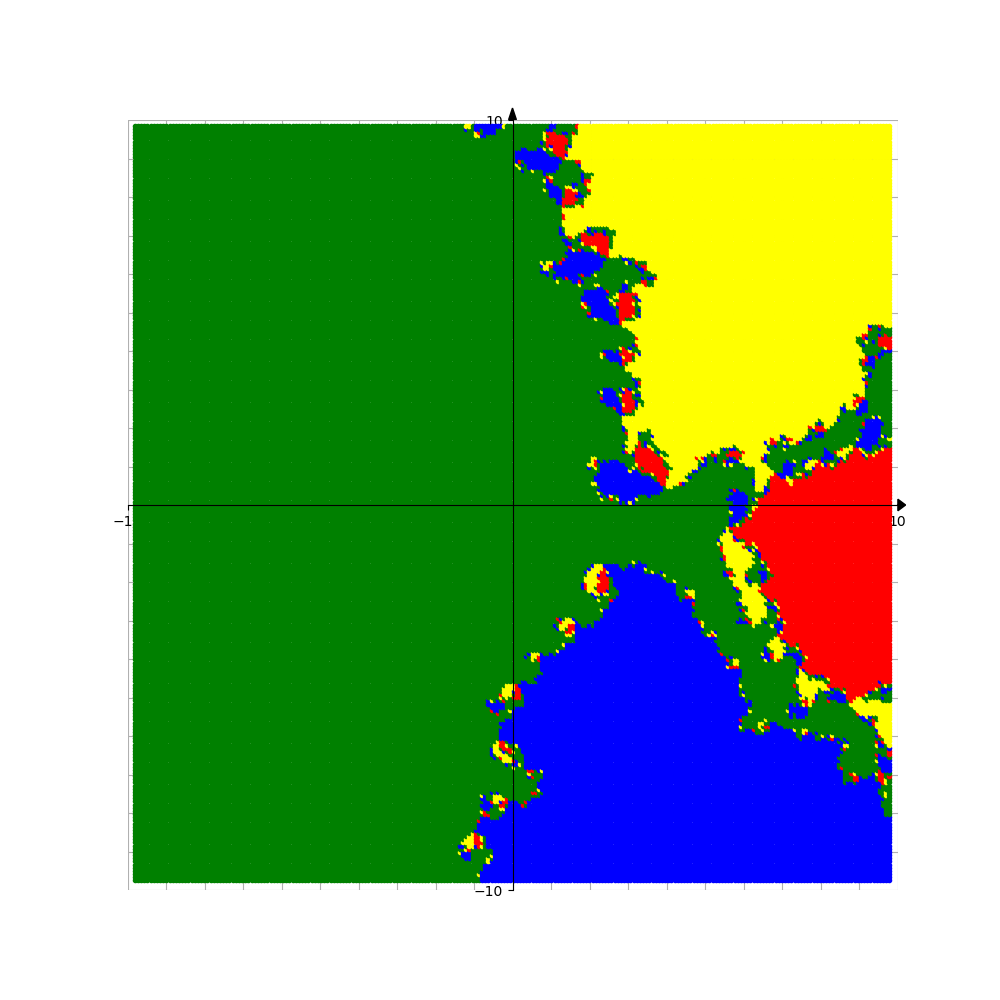}

    \bigskip
    \includegraphics[width=5.5cm]{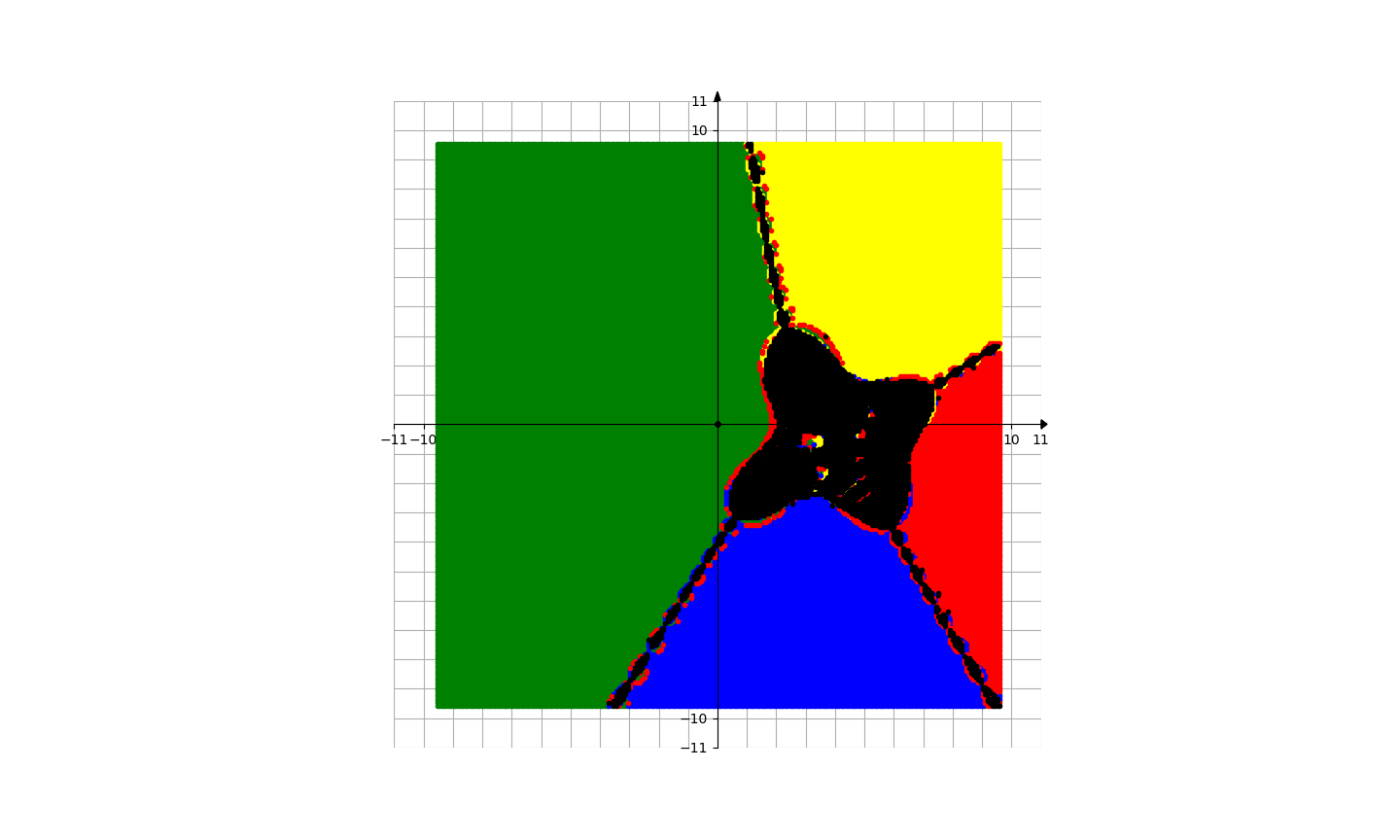}
    \includegraphics[width=3cm]{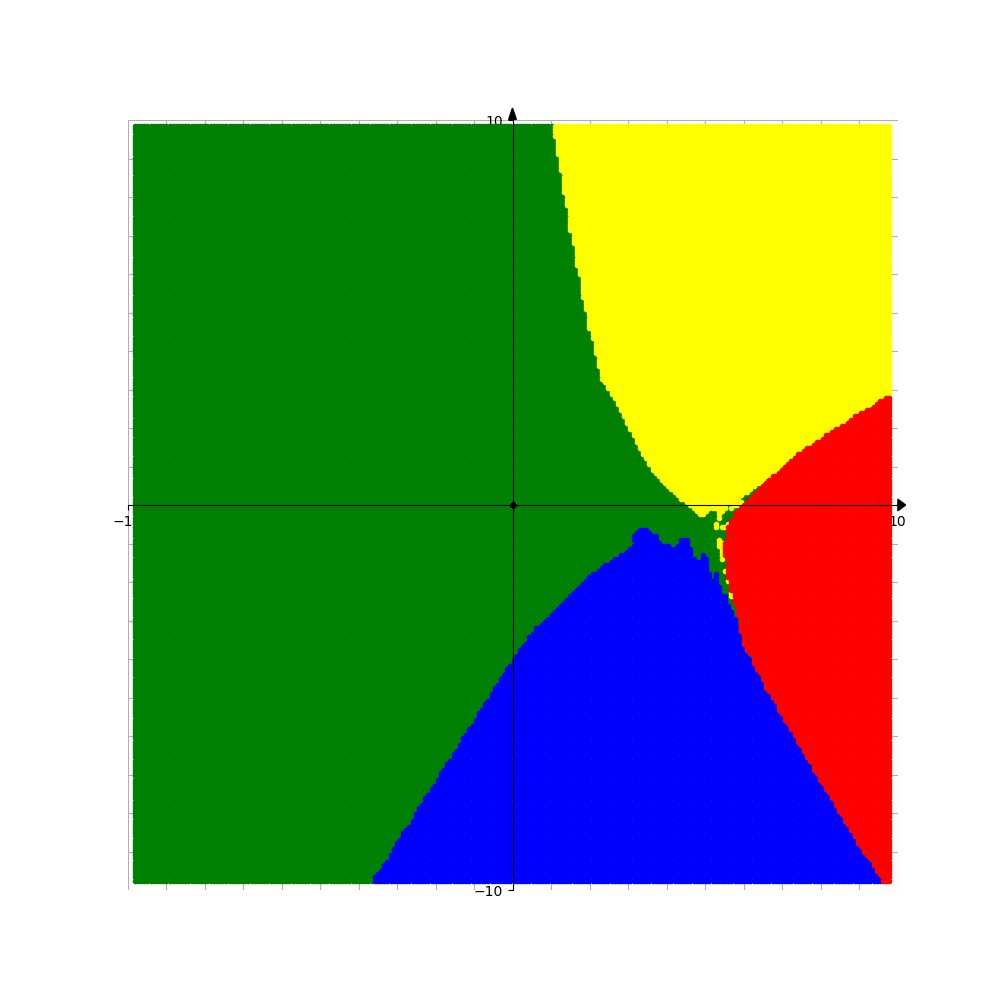}
    
    \bigskip
    \includegraphics[width=3cm]{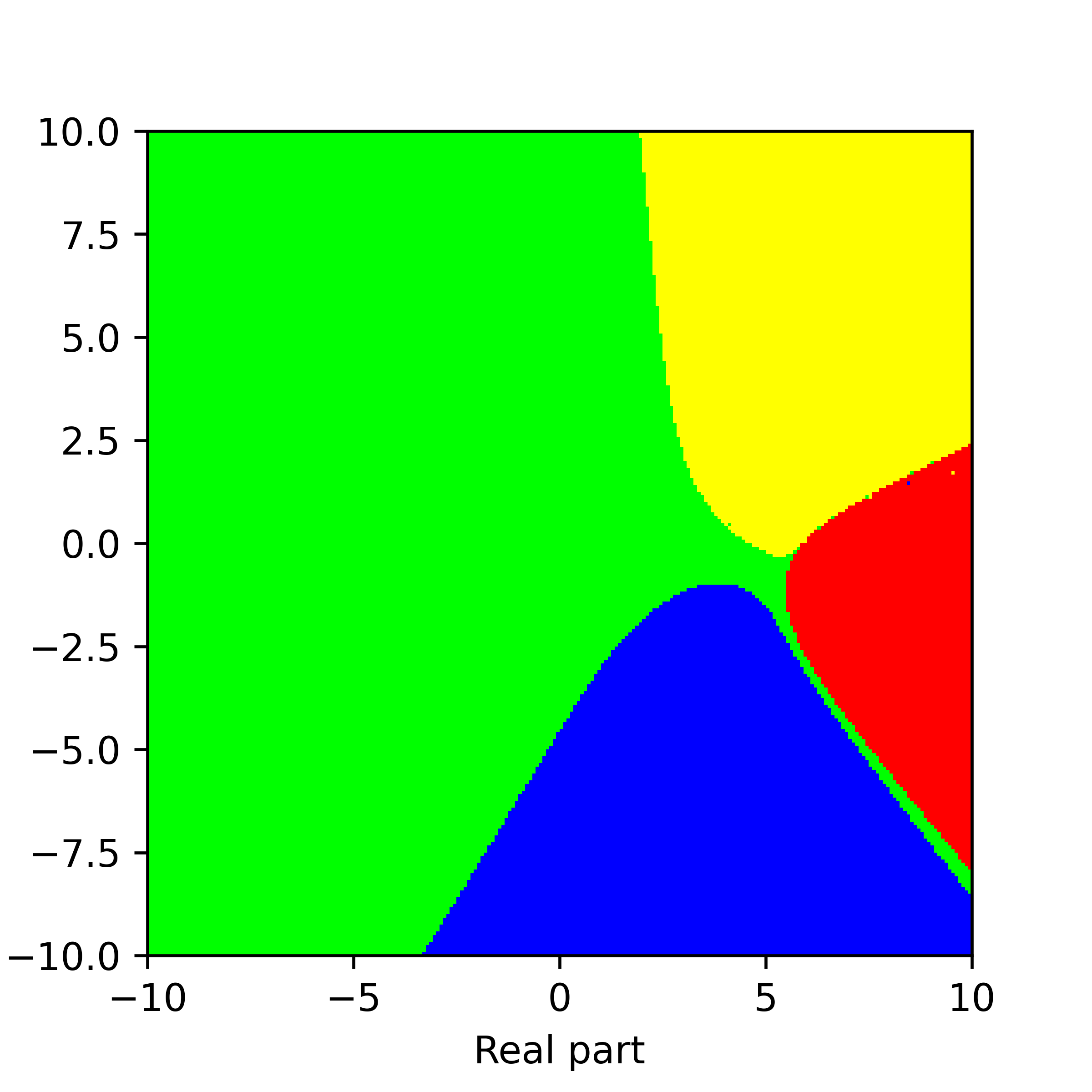}
    \includegraphics[width=3cm]{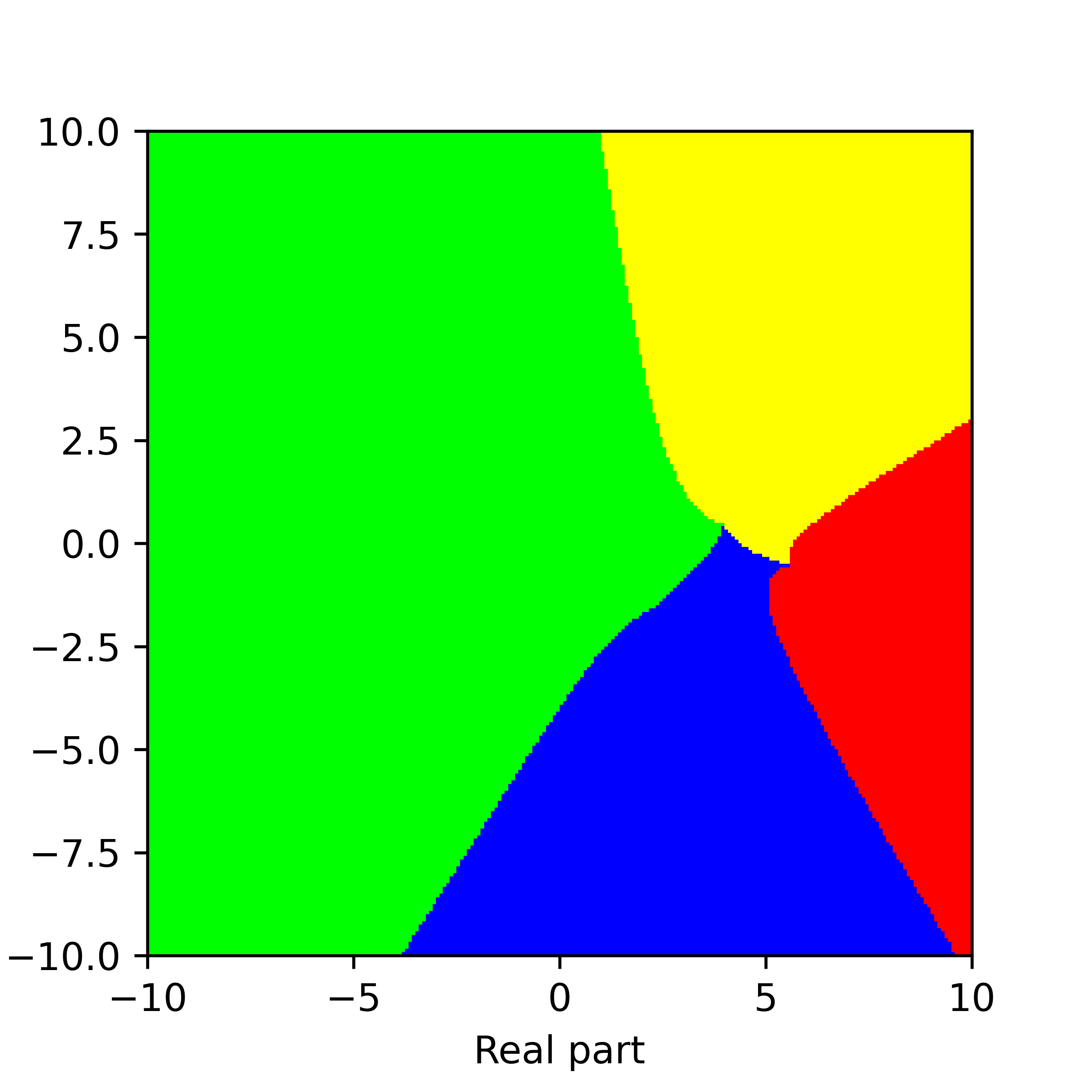}
    \includegraphics[width=3cm]{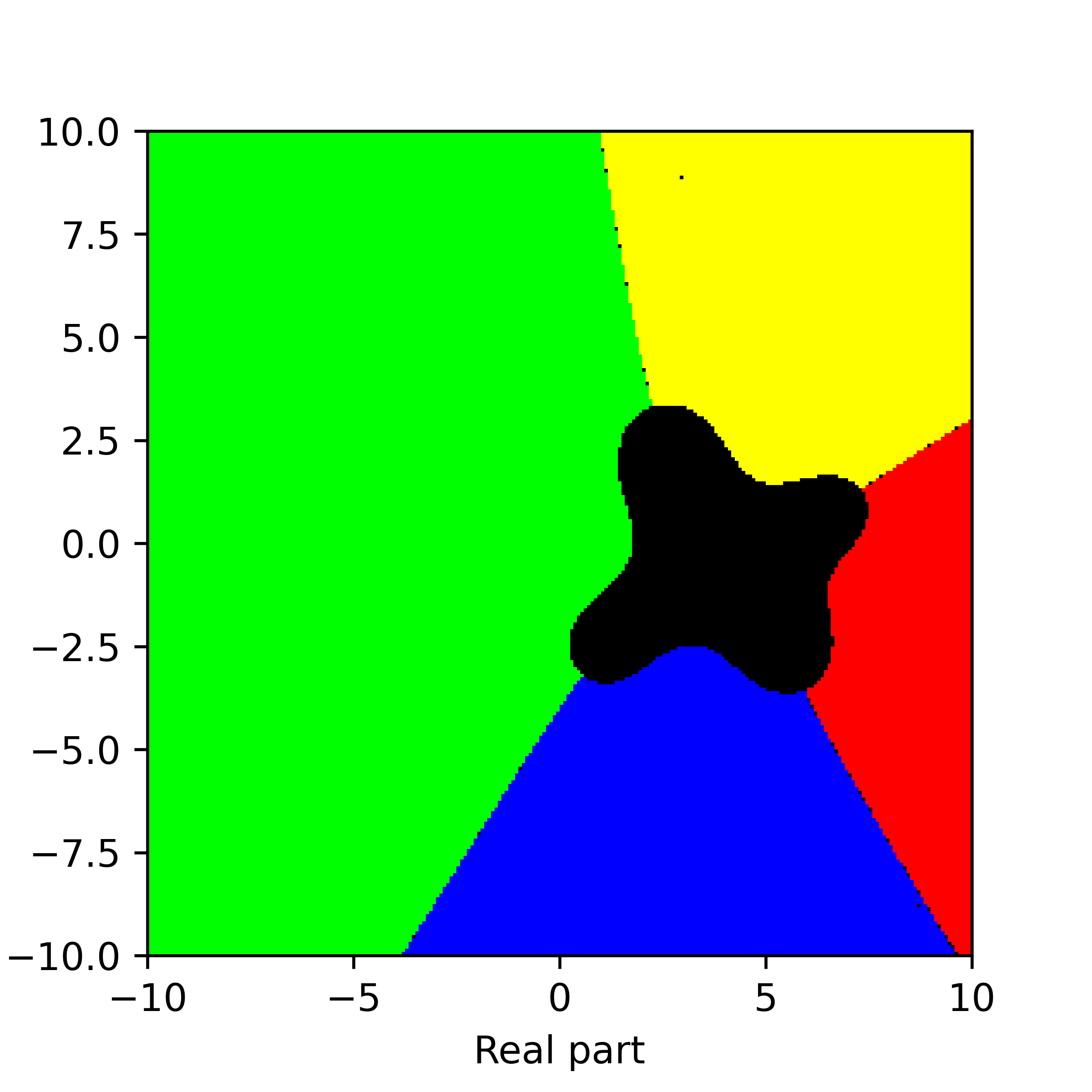}
    
    \caption{Basins of attraction for finding roots of the function $f_{19}$ by different methods. Pictures are referenced to from top to bottom, from left to right. Row 1: left picture is Voronoi's diagram, central picture is for Newton's method, right picture is for Random Relaxed Newton's method. Row 2: left picture is for Newton's method vOptimization, right picture is for BNQN. Row 3: left picture is for Newton's flow, central picture is for Newton's flow vFraction, right picture is for Newton's flow vOptimization. The black points in some of these pictures are those in the basin of attraction of critical points of $f_{19}$.}
    \label{fig:f19}
\end{figure}

\begin{figure}
    \centering
    \includegraphics[width=5cm]{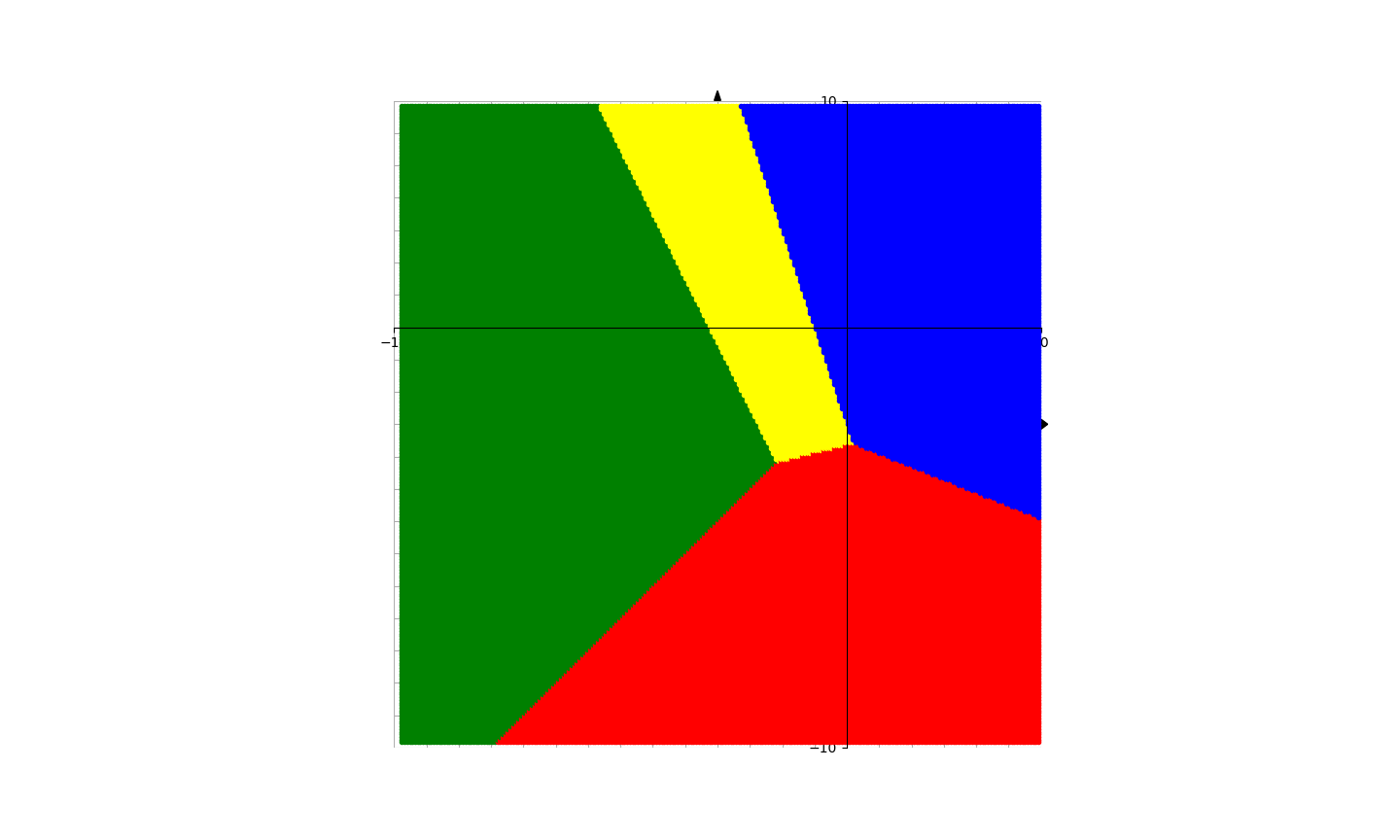}
    \includegraphics[width=3cm]{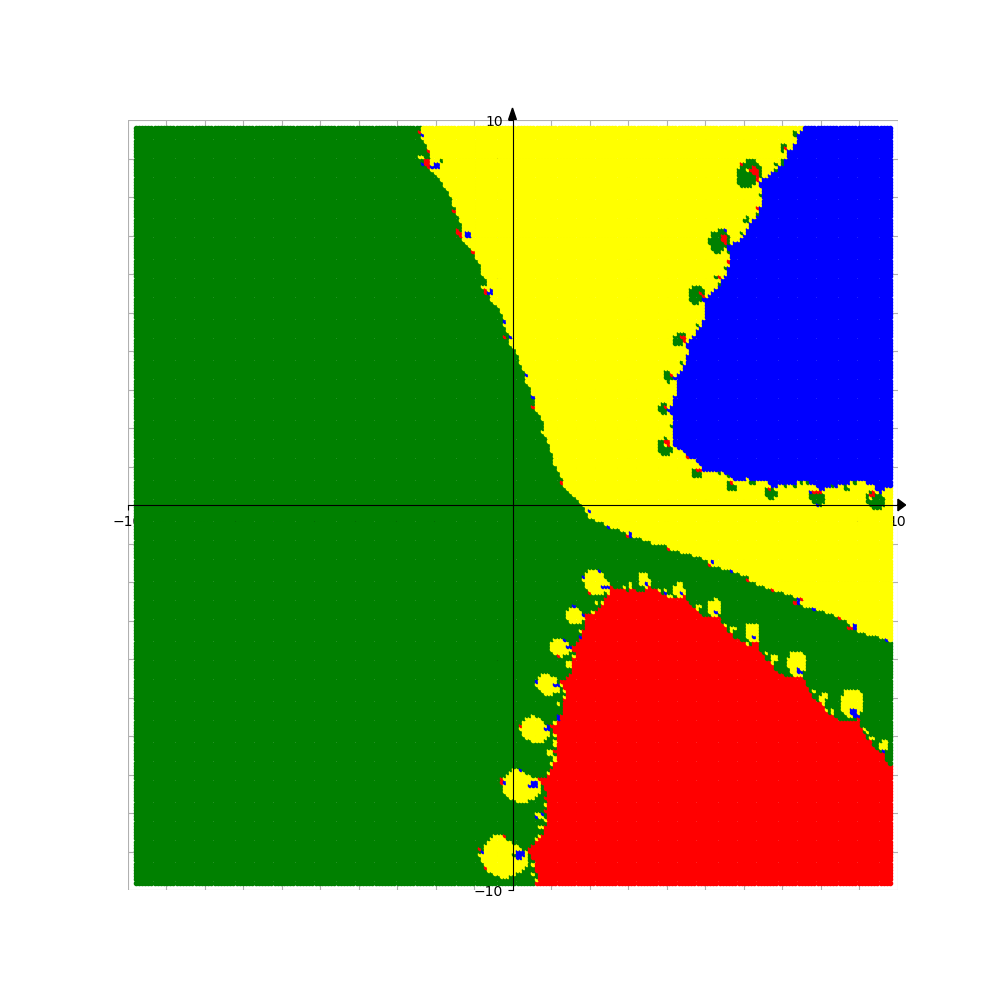}
    \includegraphics[width=3cm]{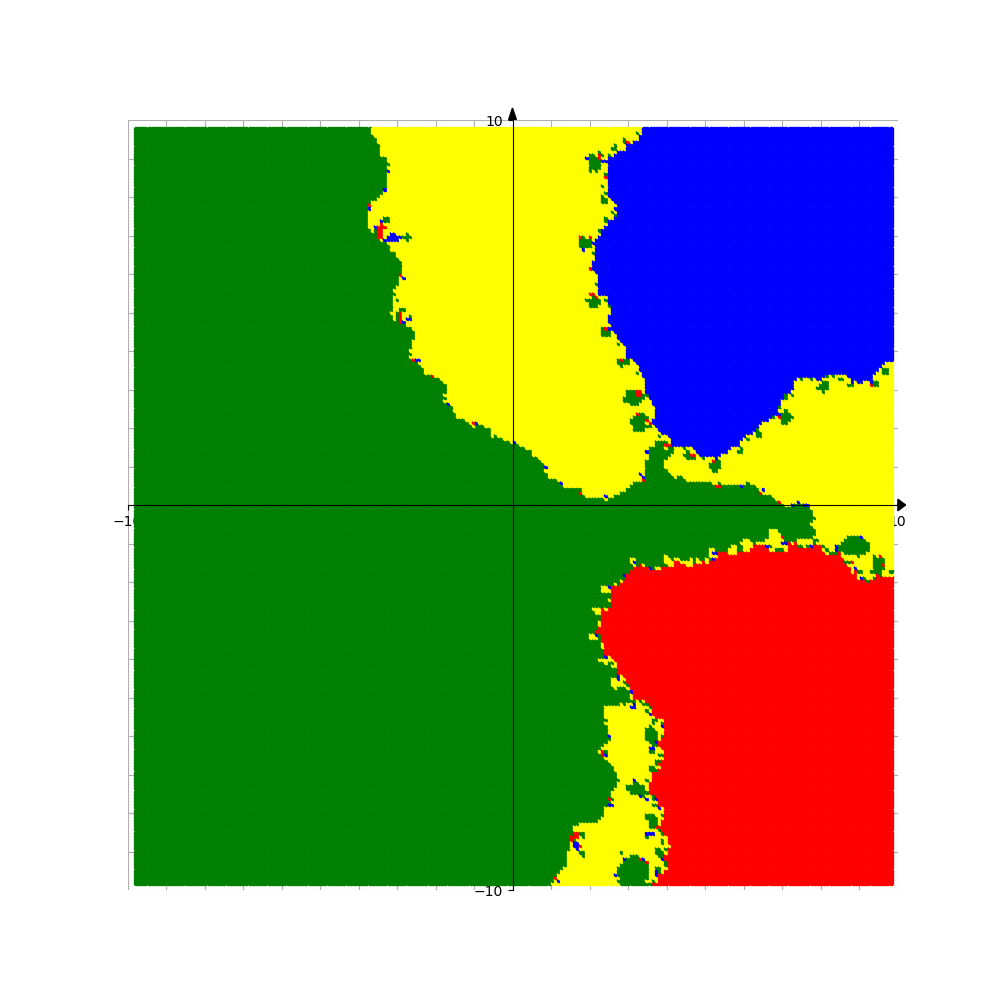}

    \bigskip
    \includegraphics[width=5.5cm]{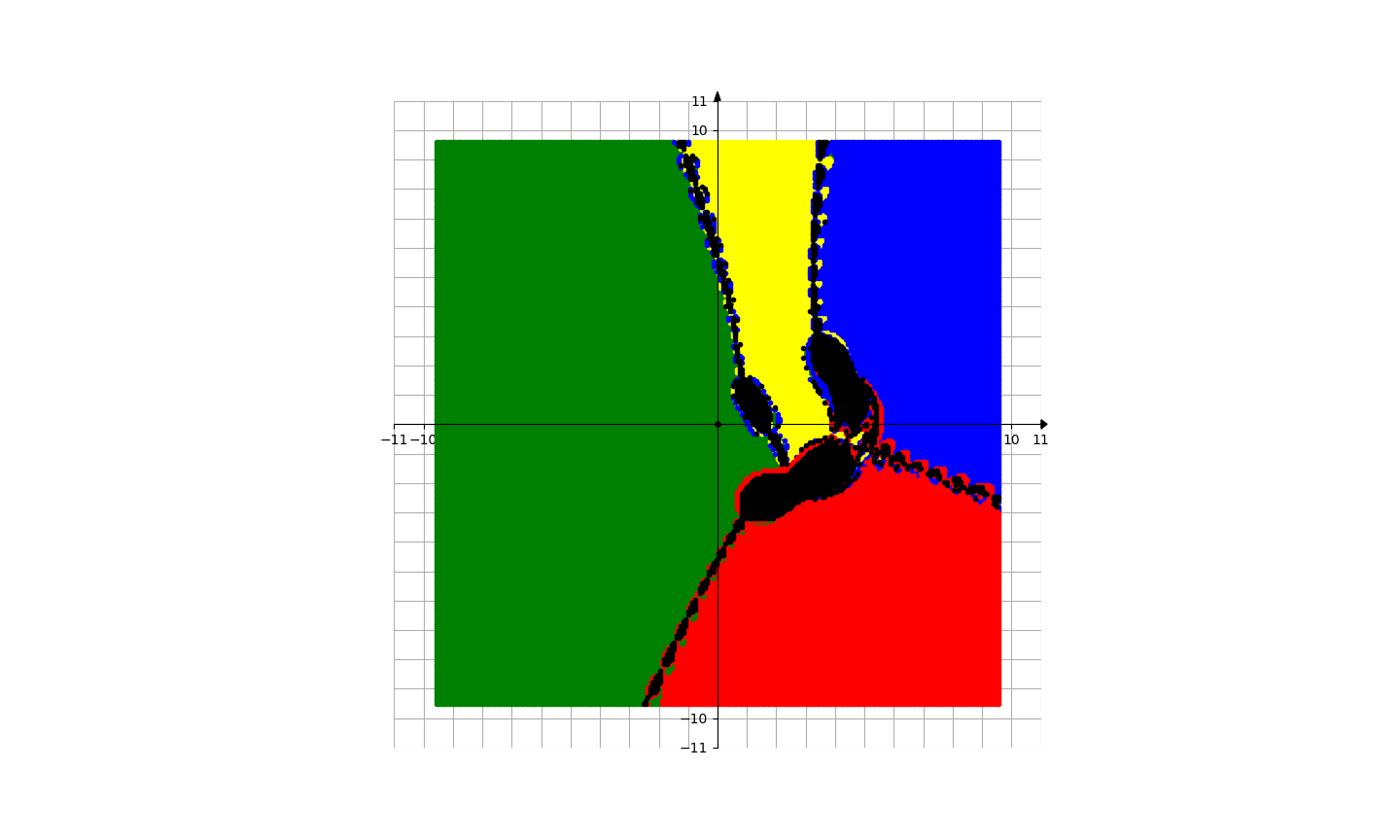}
    \includegraphics[width=3cm]{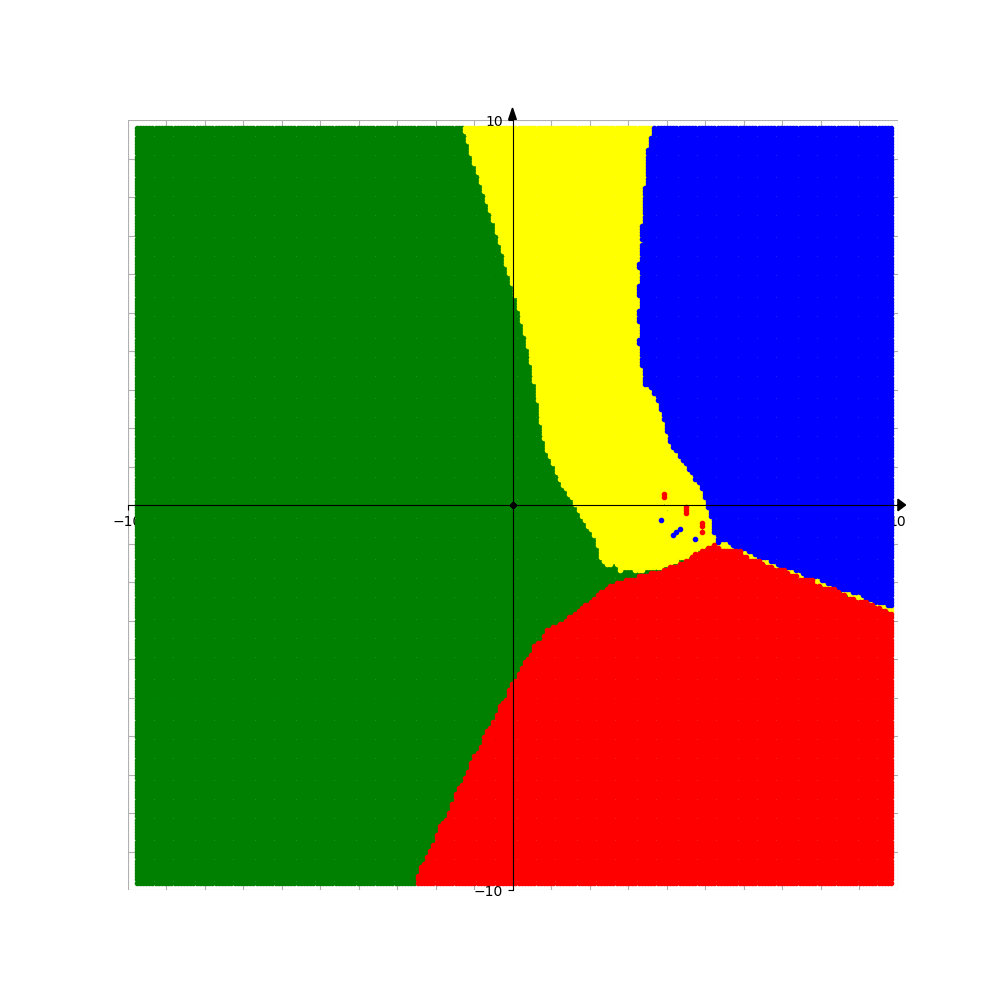}
    
    \bigskip
    \includegraphics[width=3cm]{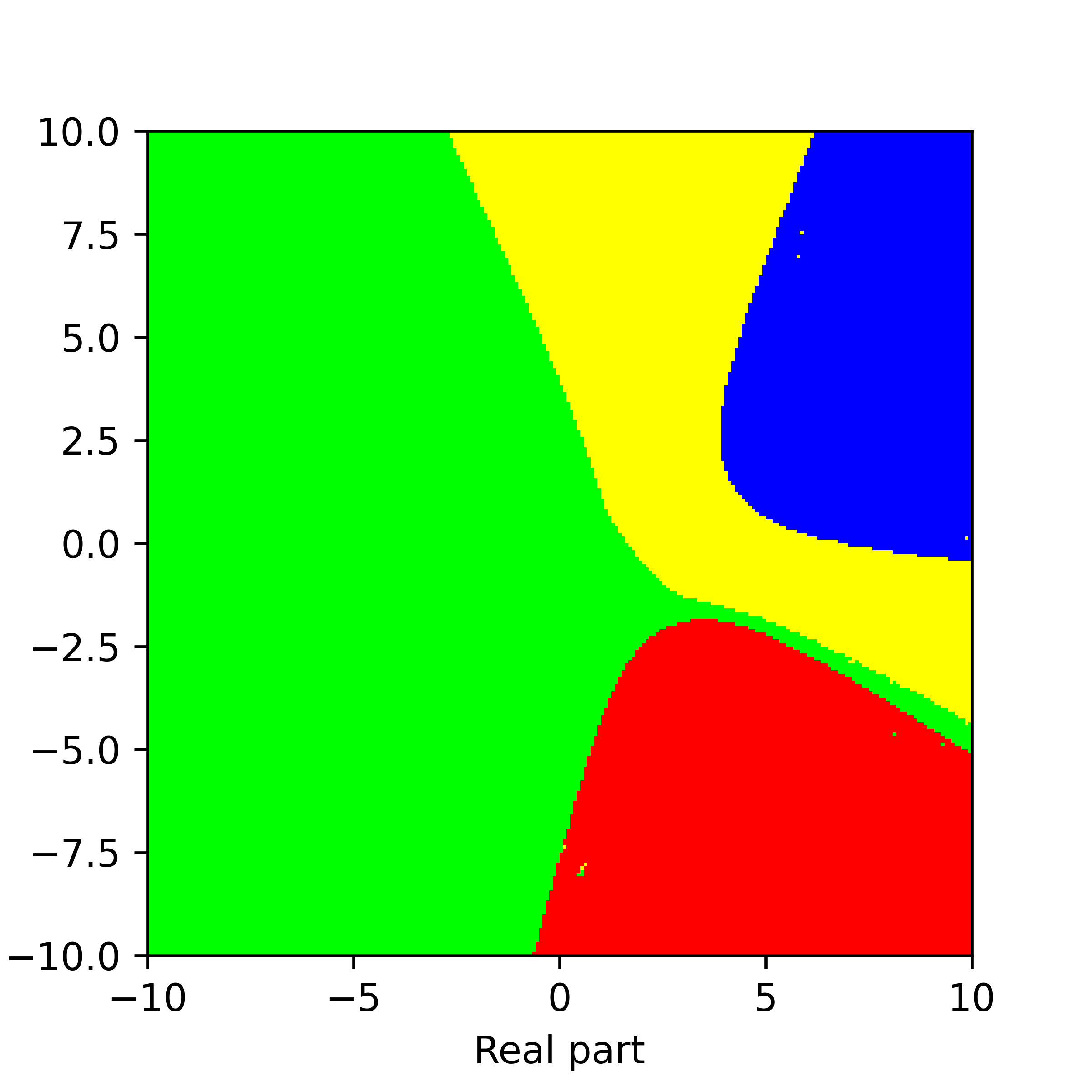}
    \includegraphics[width=3cm]{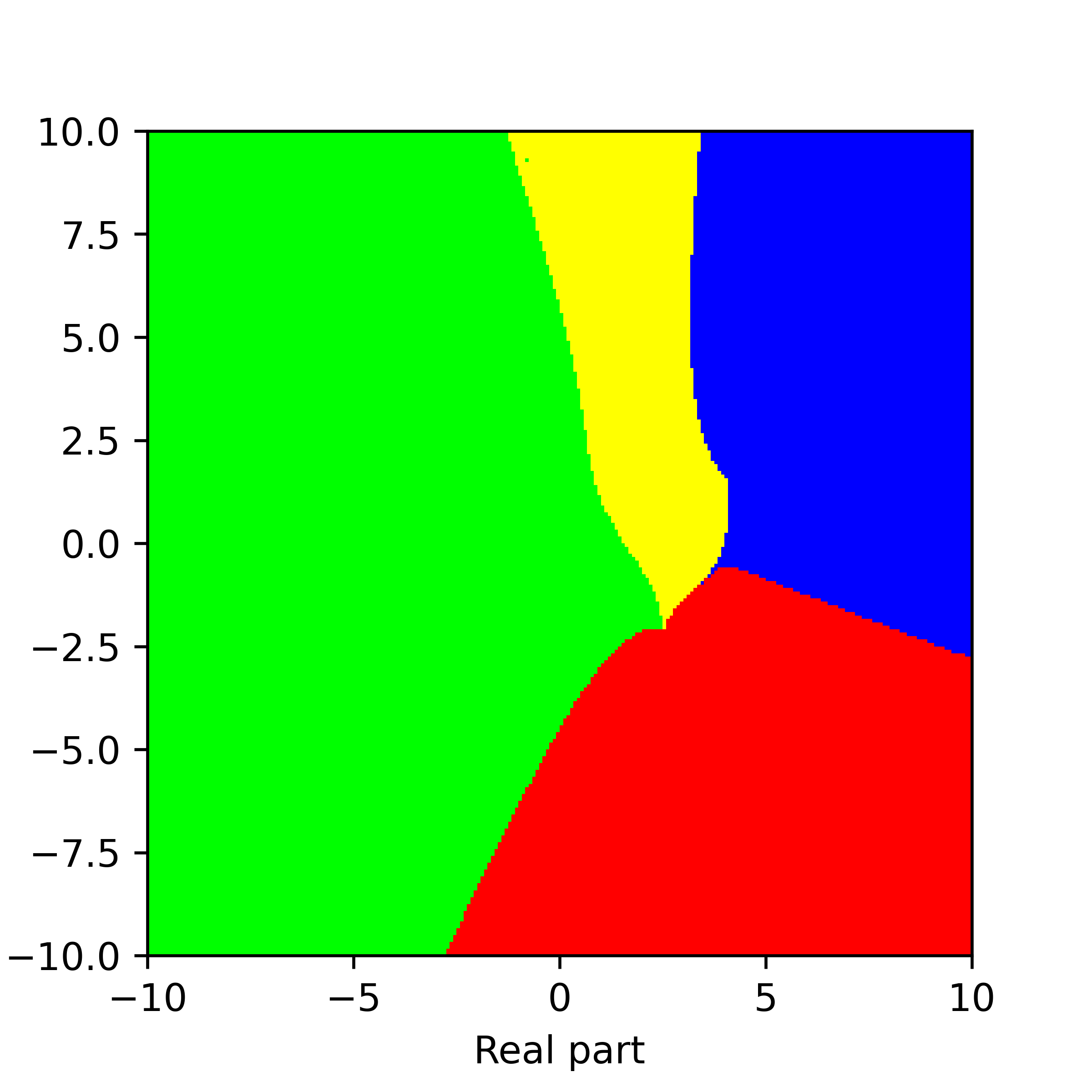}
    \includegraphics[width=3cm]{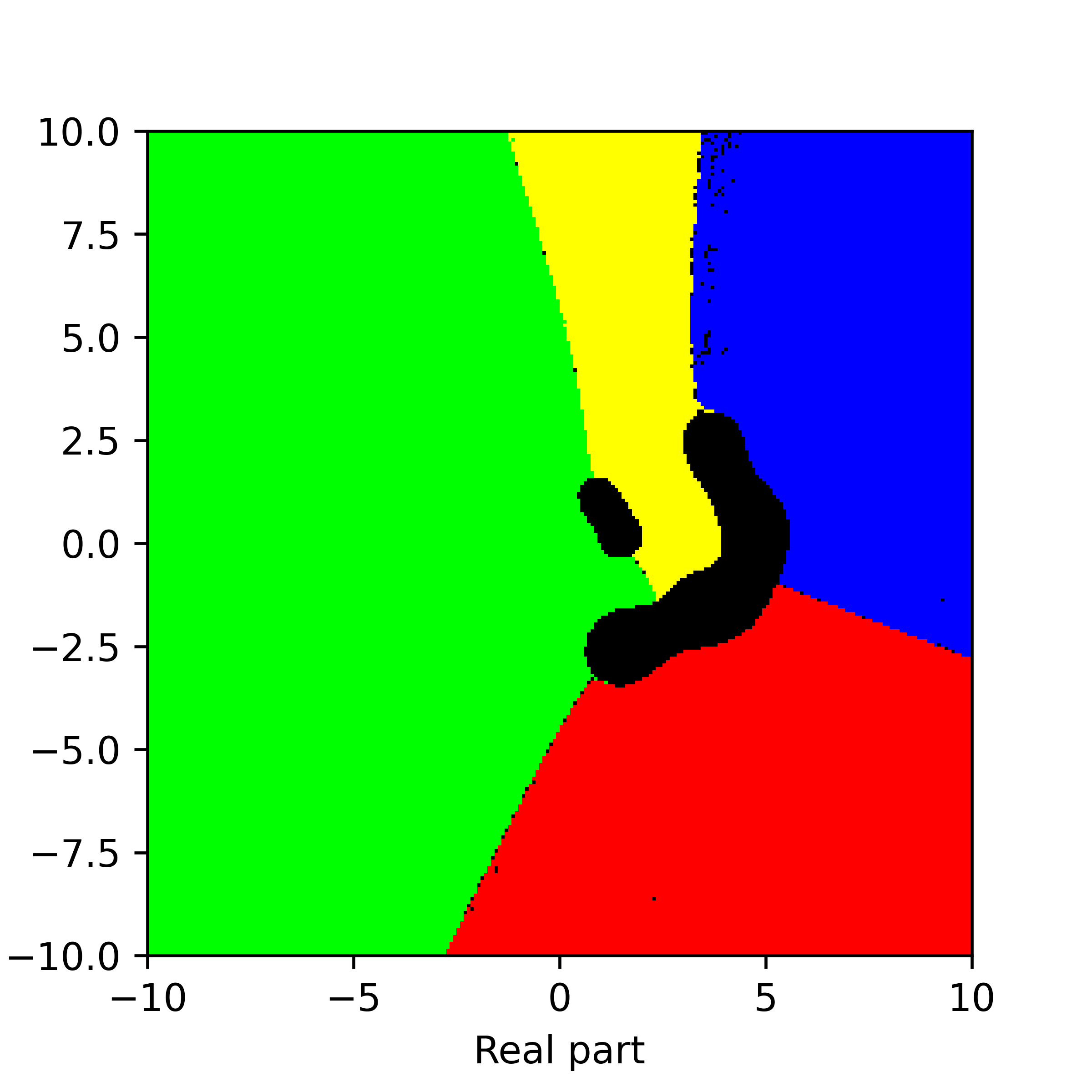}
    
    \caption{Basins of attraction for finding roots of the function $f_{20}$ by different methods. Pictures are referenced to from top to bottom, from left to right. Row 1: left picture is Voronoi's diagram, central picture is for Newton's method, right picture is for Random Relaxed Newton's method. Row 2: left picture is for Newton's method vOptimization, right picture is for BNQN. Row 3: left picture is for Newton's flow, central picture is for Newton's flow vFraction, right picture is for Newton's flow vOptimization. The black points in some of these pictures are those in the basin of attraction of critical points of $f_{20}$.}
    \label{fig:f20}
\end{figure}

\begin{figure}
    \centering
    \includegraphics[width=5cm]{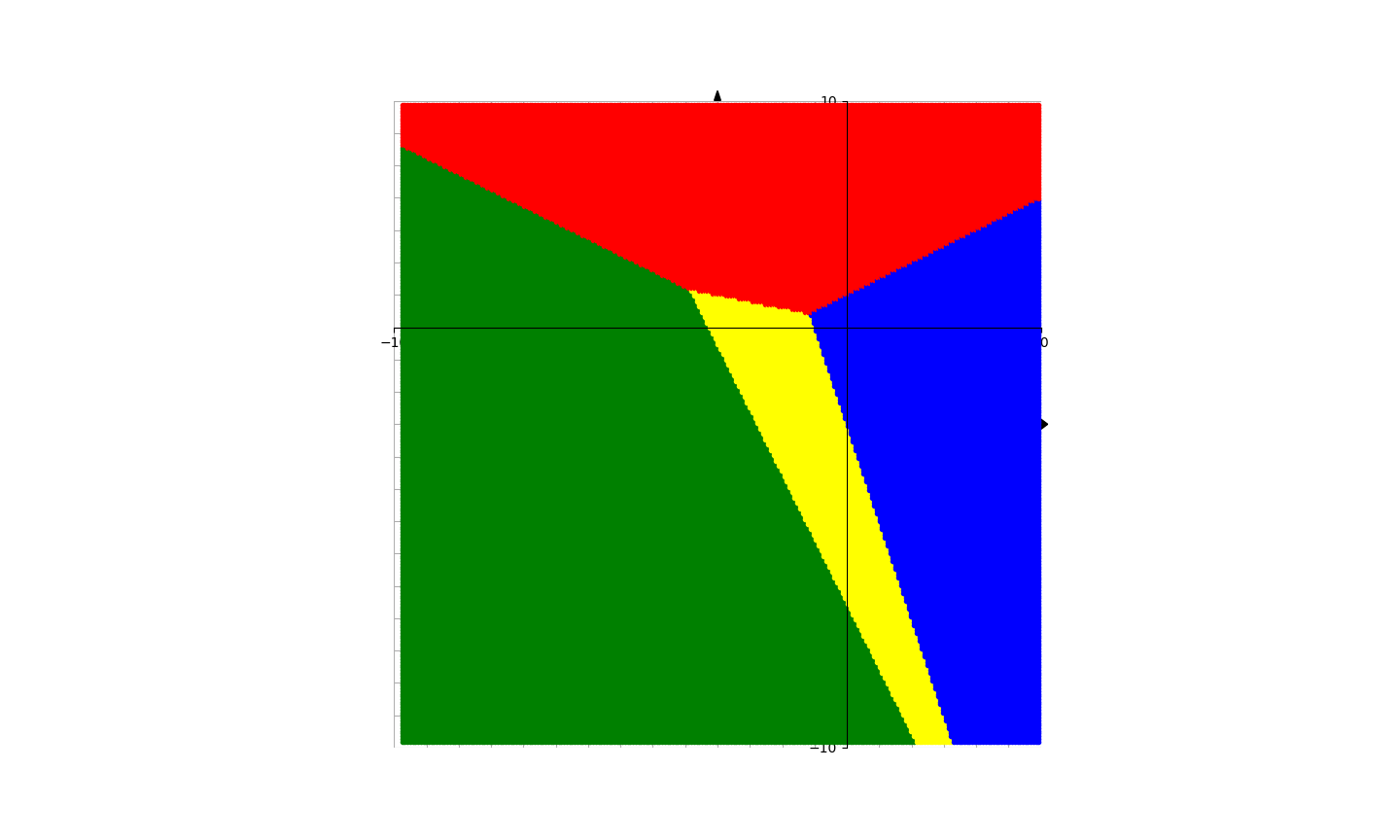}
    \includegraphics[width=3cm]{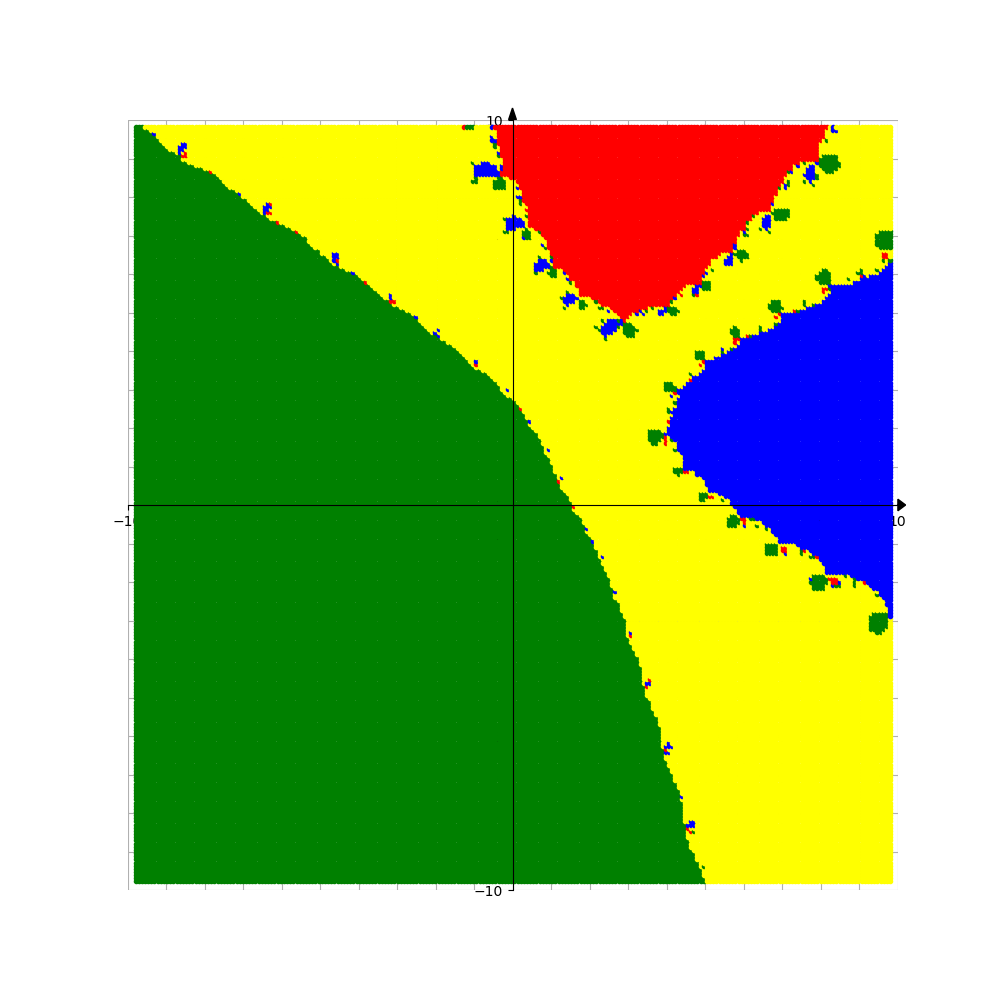}
    \includegraphics[width=3cm]{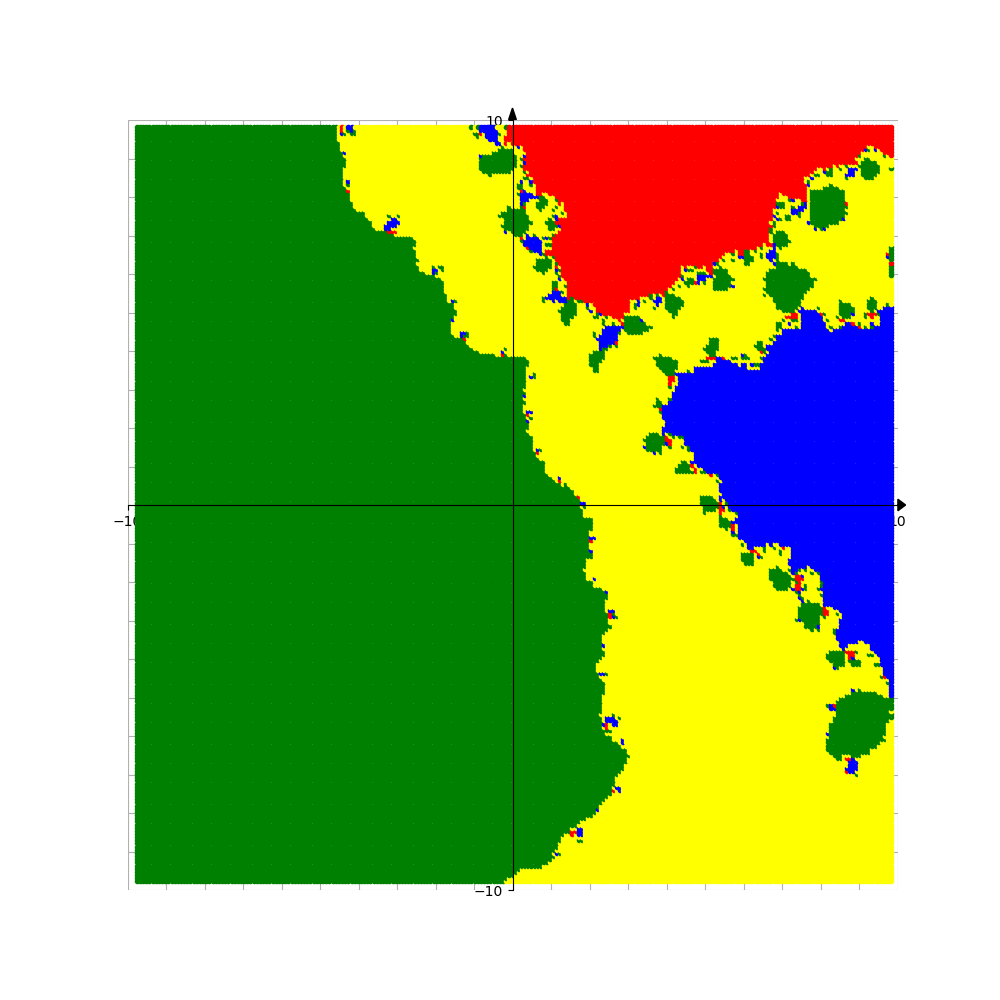}

    \bigskip
    \includegraphics[width=5.5cm]{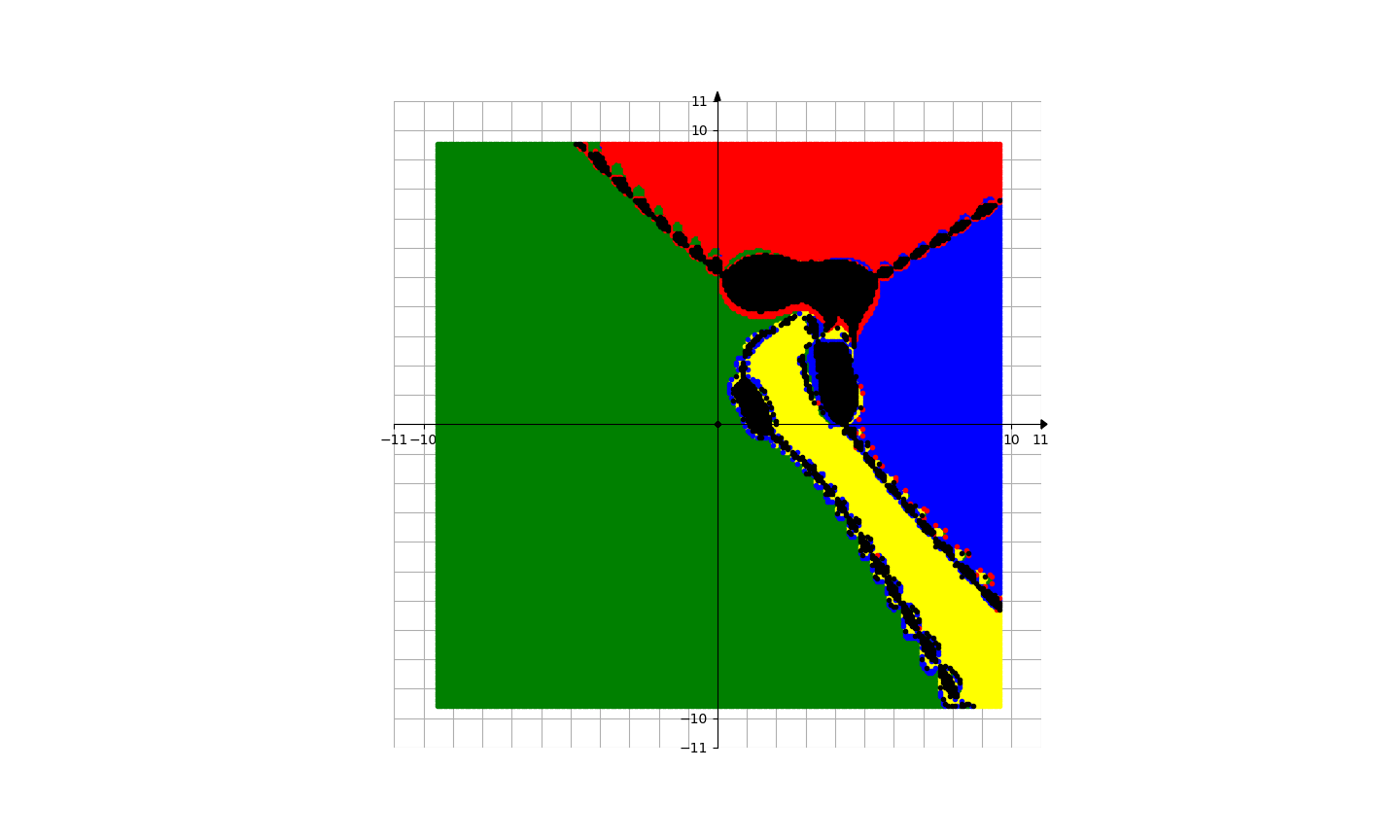}
    \includegraphics[width=3cm]{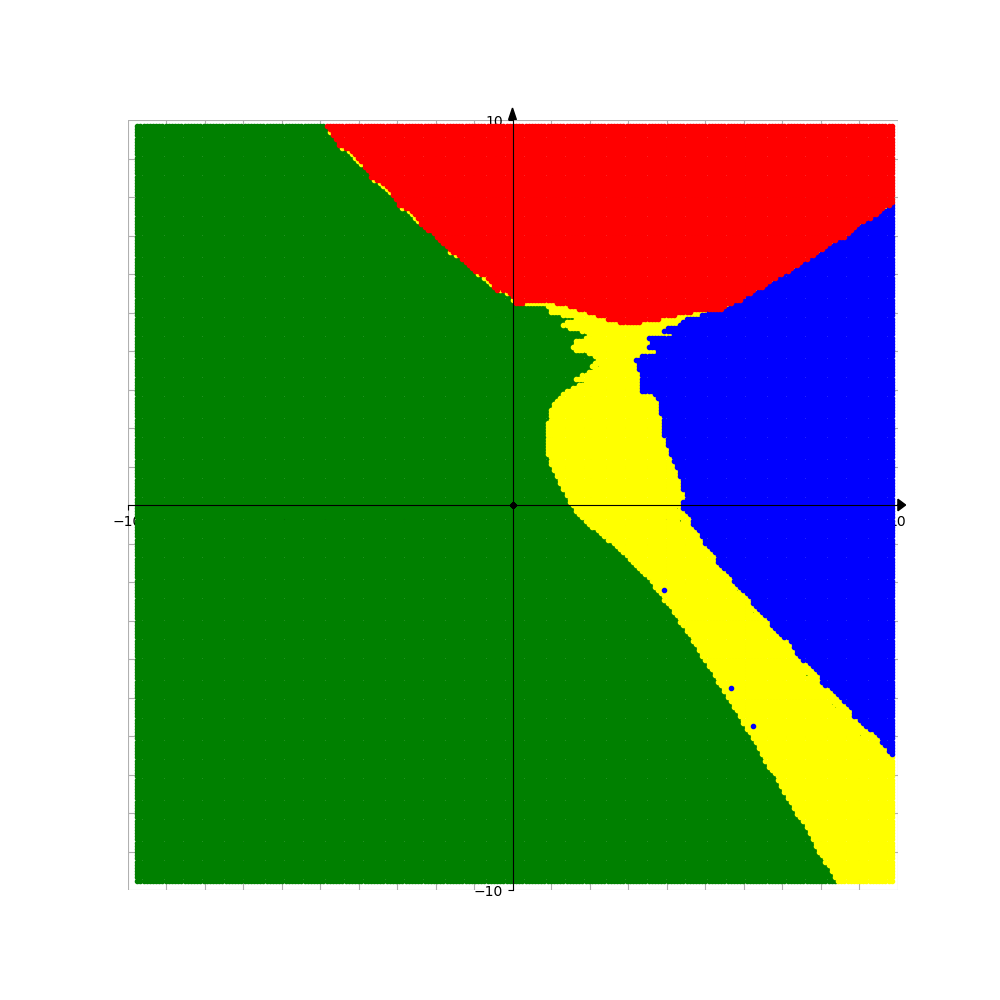}
    
    \bigskip
    \includegraphics[width=3cm]{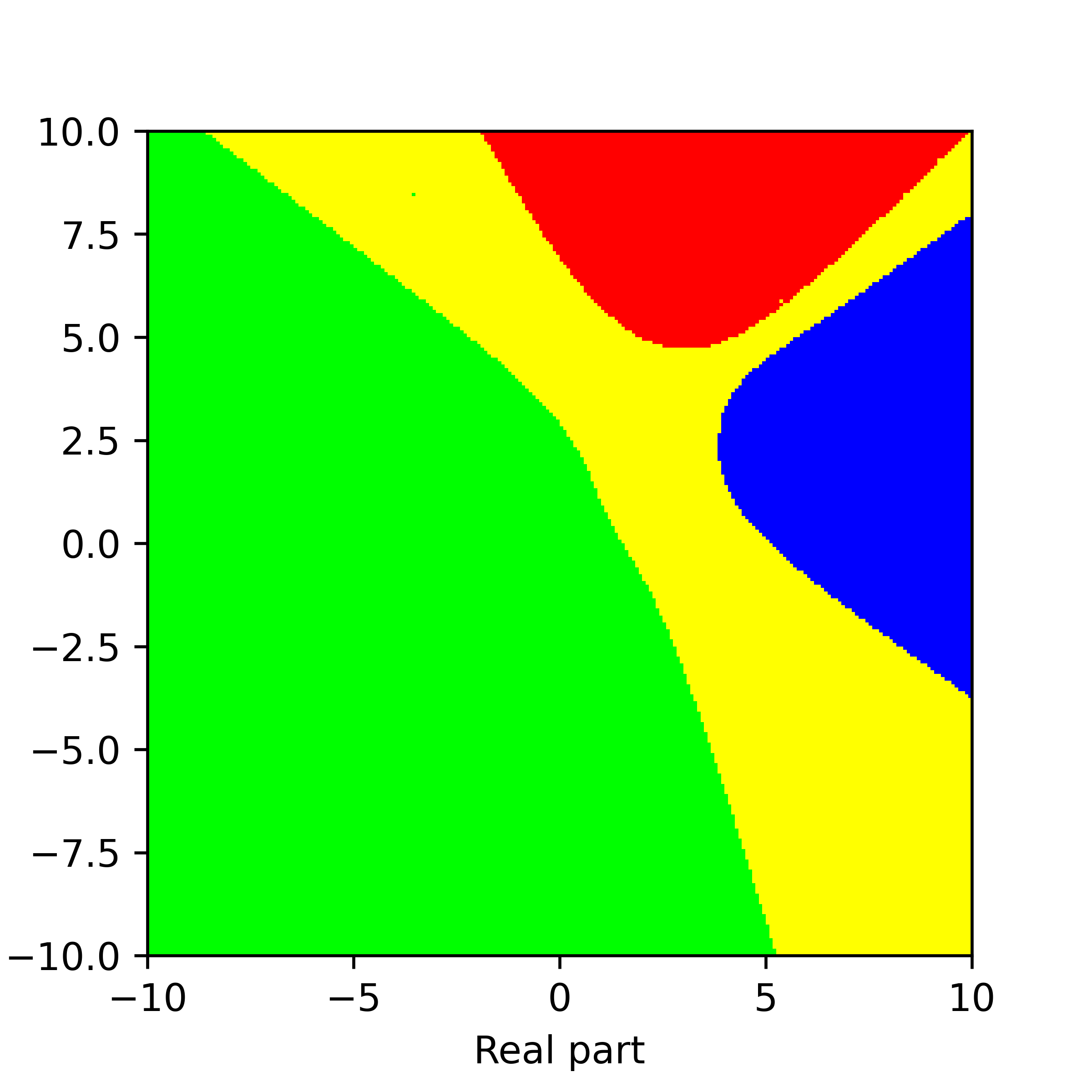}
    \includegraphics[width=3cm]{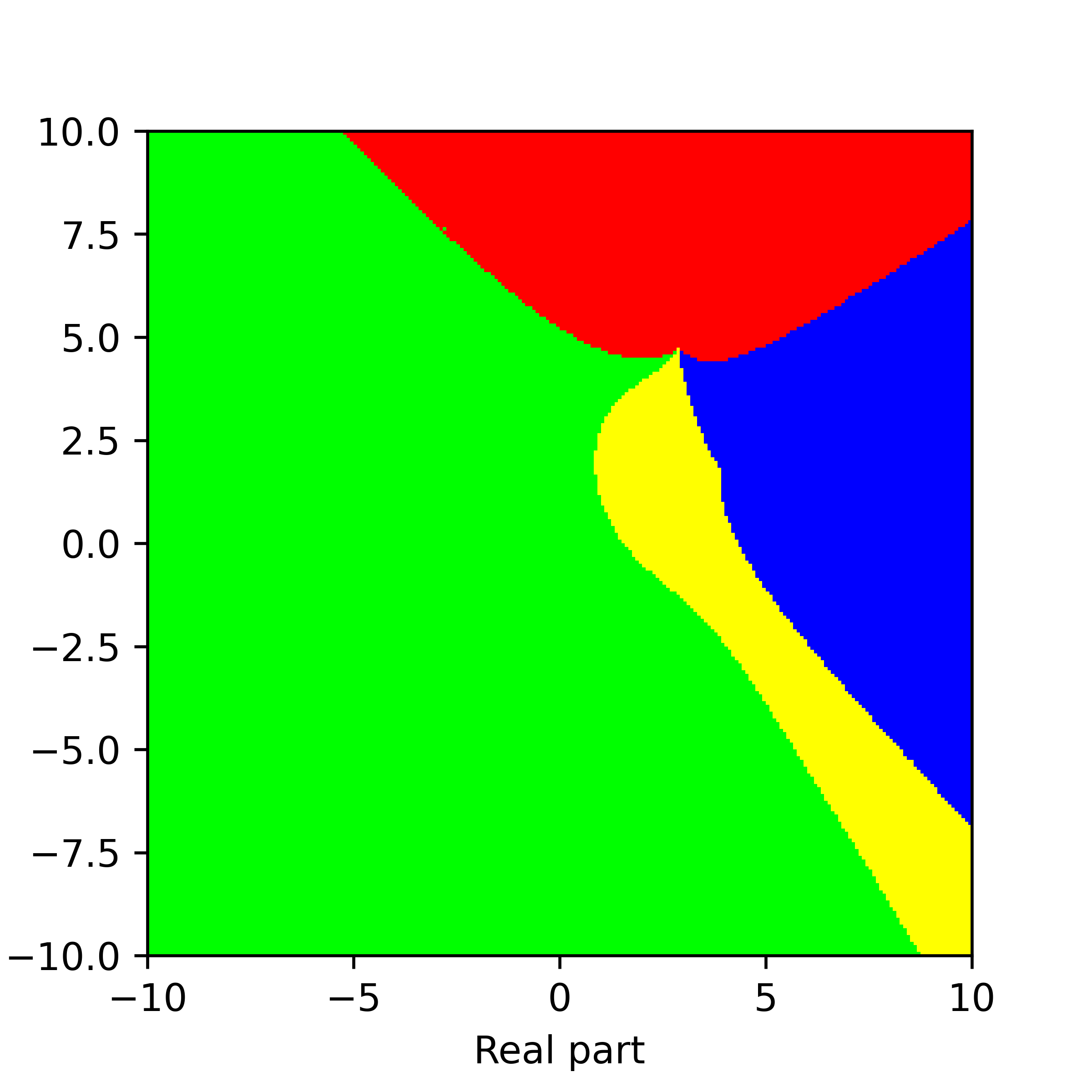}
    \includegraphics[width=3cm]{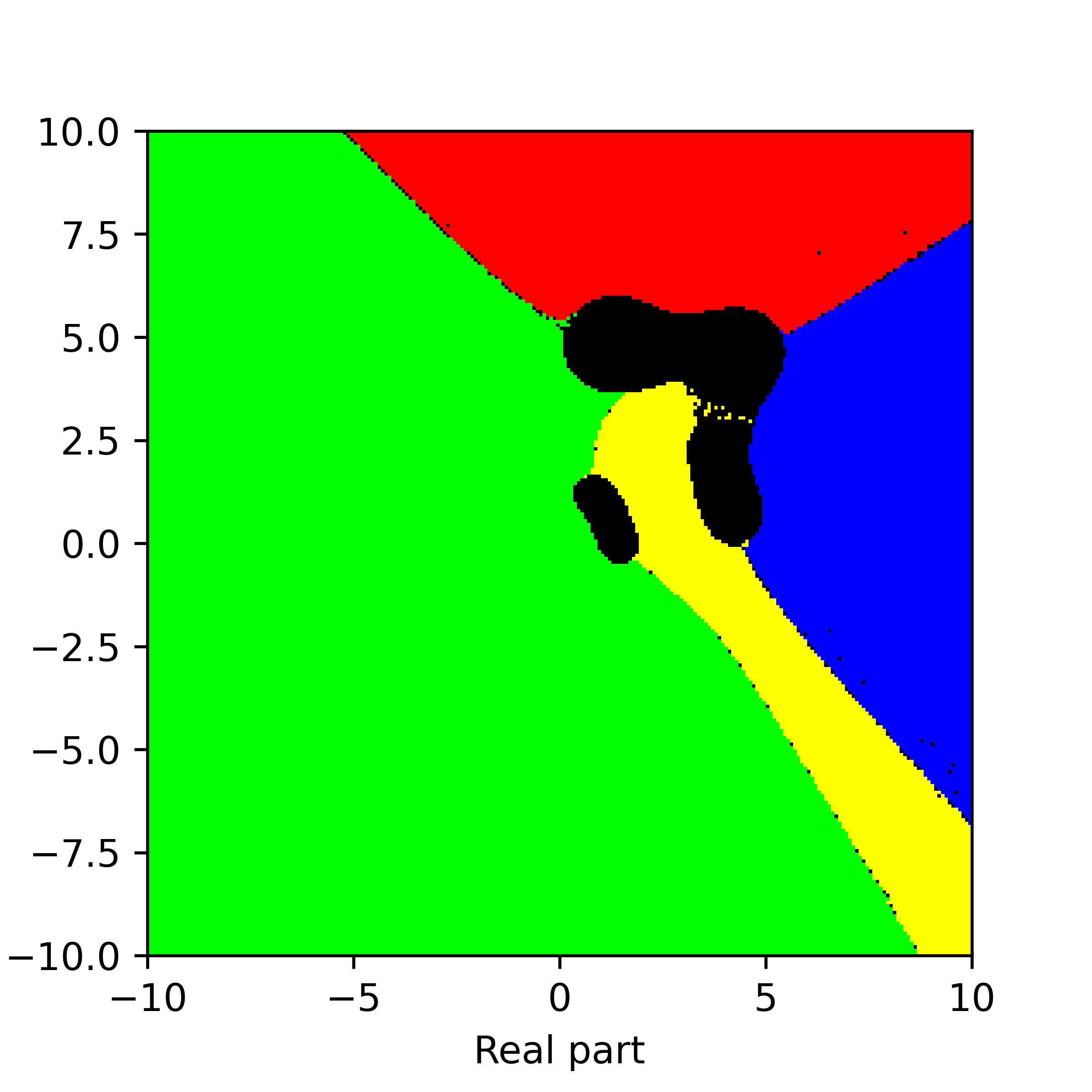}
    
    \caption{Basins of attraction for finding roots of the function $f_{21}$ by different methods. Pictures are referenced to from top to bottom, from left to right. Row 1: left picture is Voronoi's diagram, central picture is for Newton's method, right picture is for Random Relaxed Newton's method. Row 2: left picture is for Newton's method vOptimization, right picture is for BNQN. Row 3: left picture is for Newton's flow, central picture is for Newton's flow vFraction, right picture is for Newton's flow vOptimization. The black points in some of these pictures are those in the basin of attraction of critical points of $f_{21}$.}
    \label{fig:f21}
\end{figure}

\begin{figure}
    \centering
    \includegraphics[width=5cm]{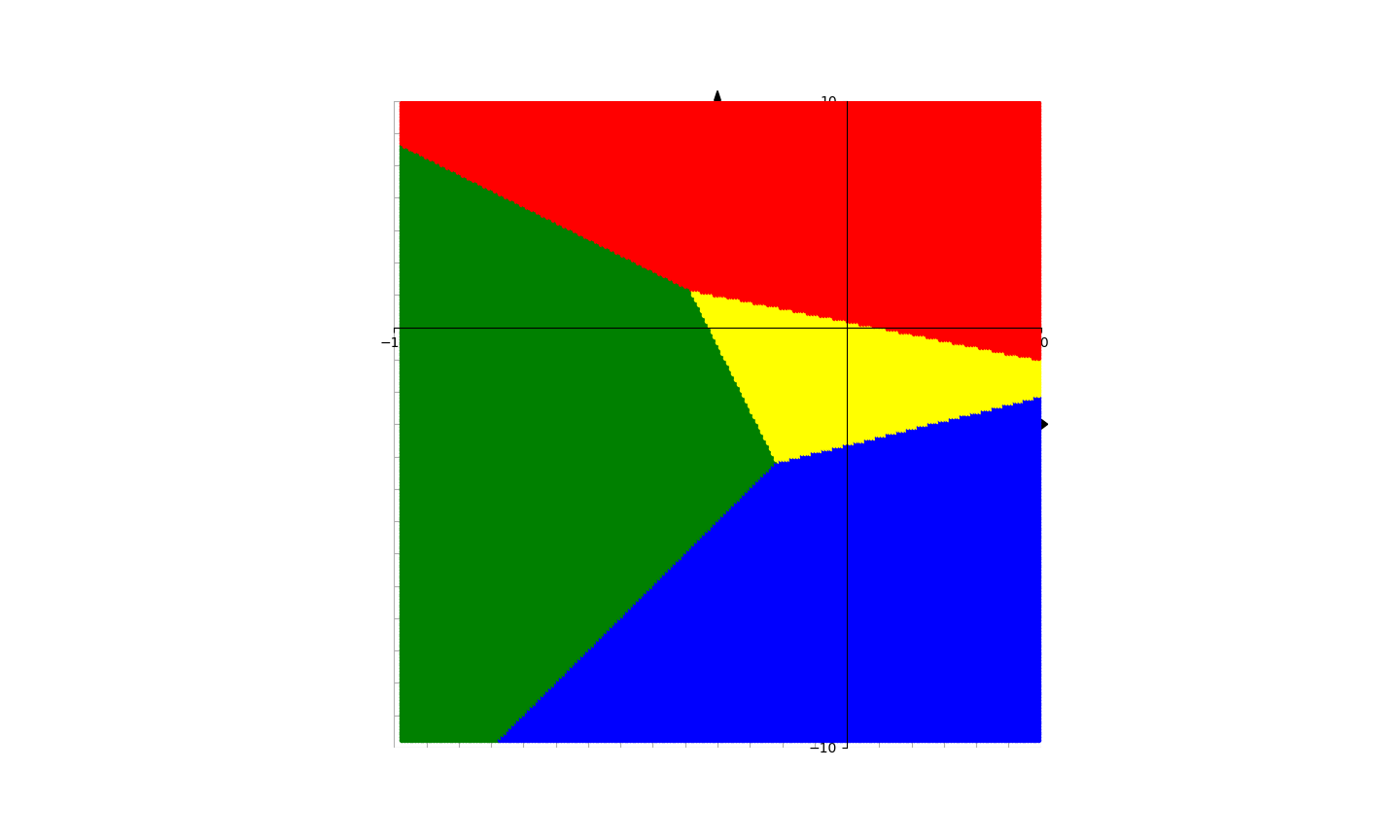}
    \includegraphics[width=3cm]{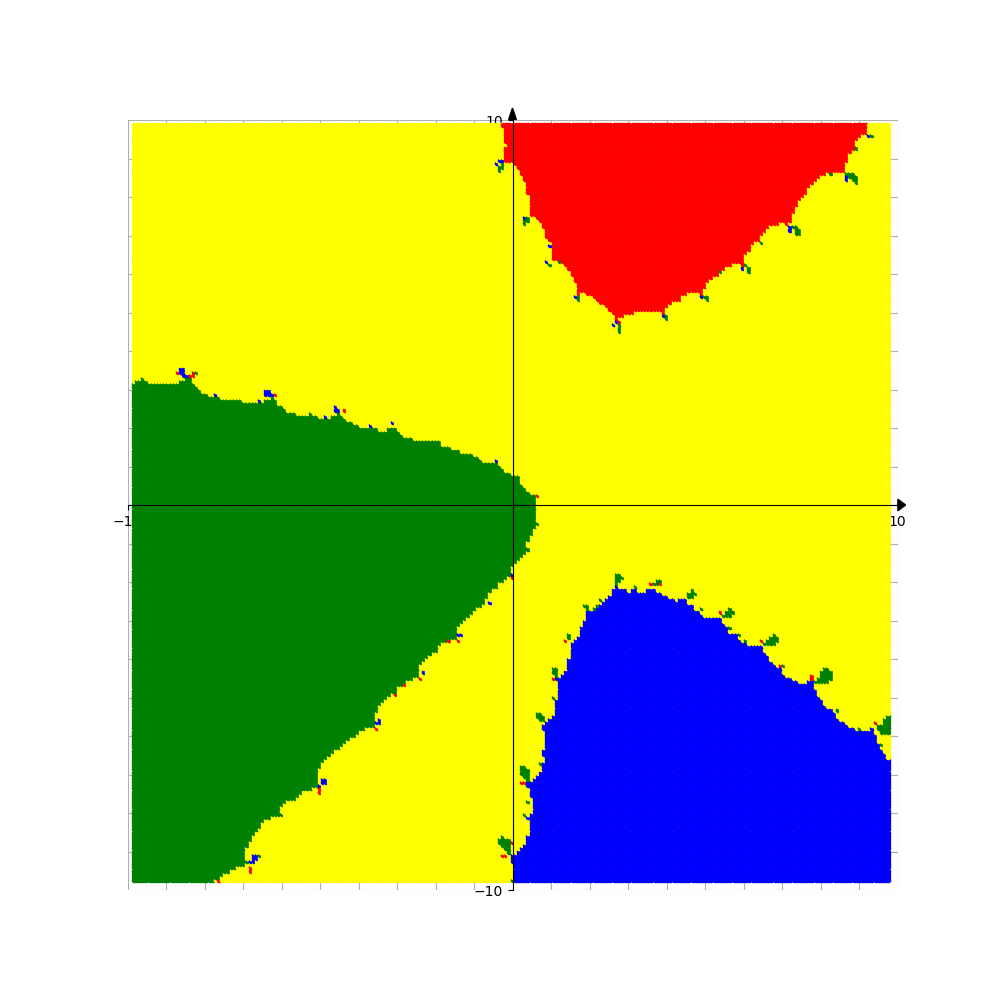}
    \includegraphics[width=3cm]{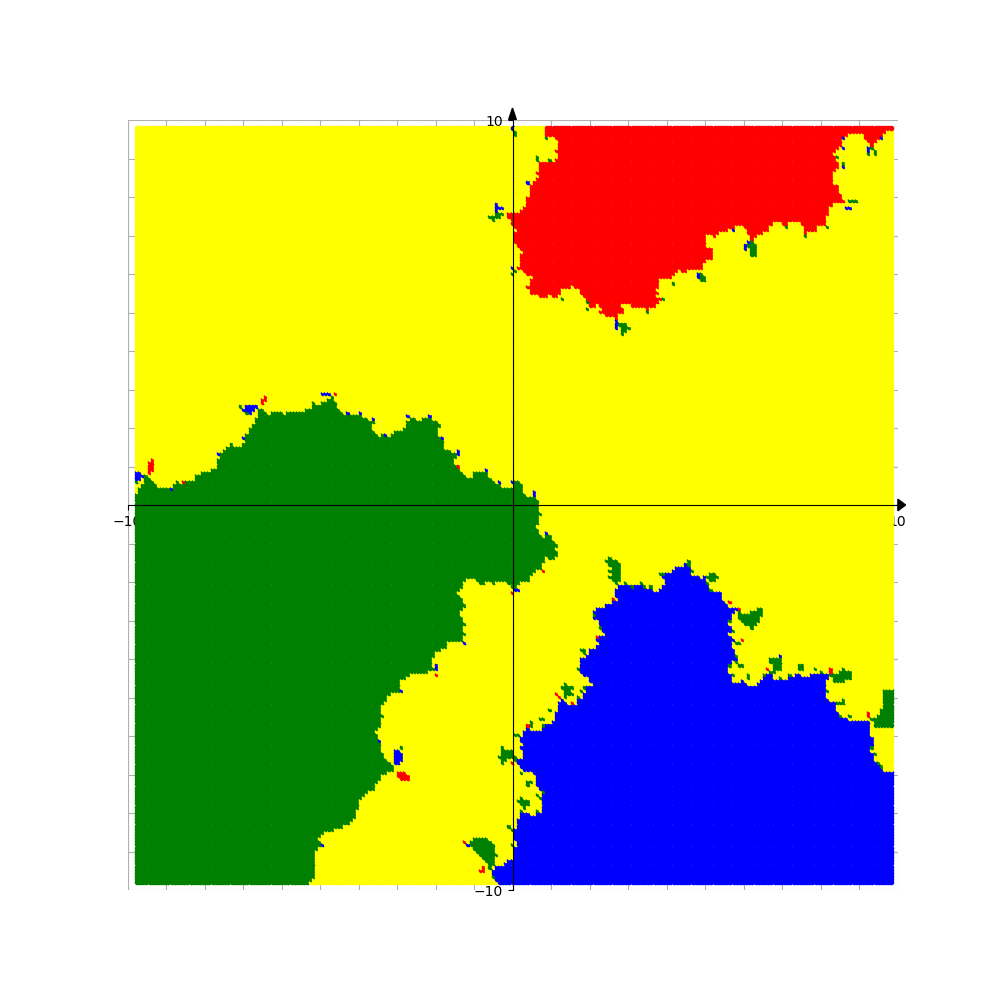}

    \bigskip
    \includegraphics[width=5.5cm]{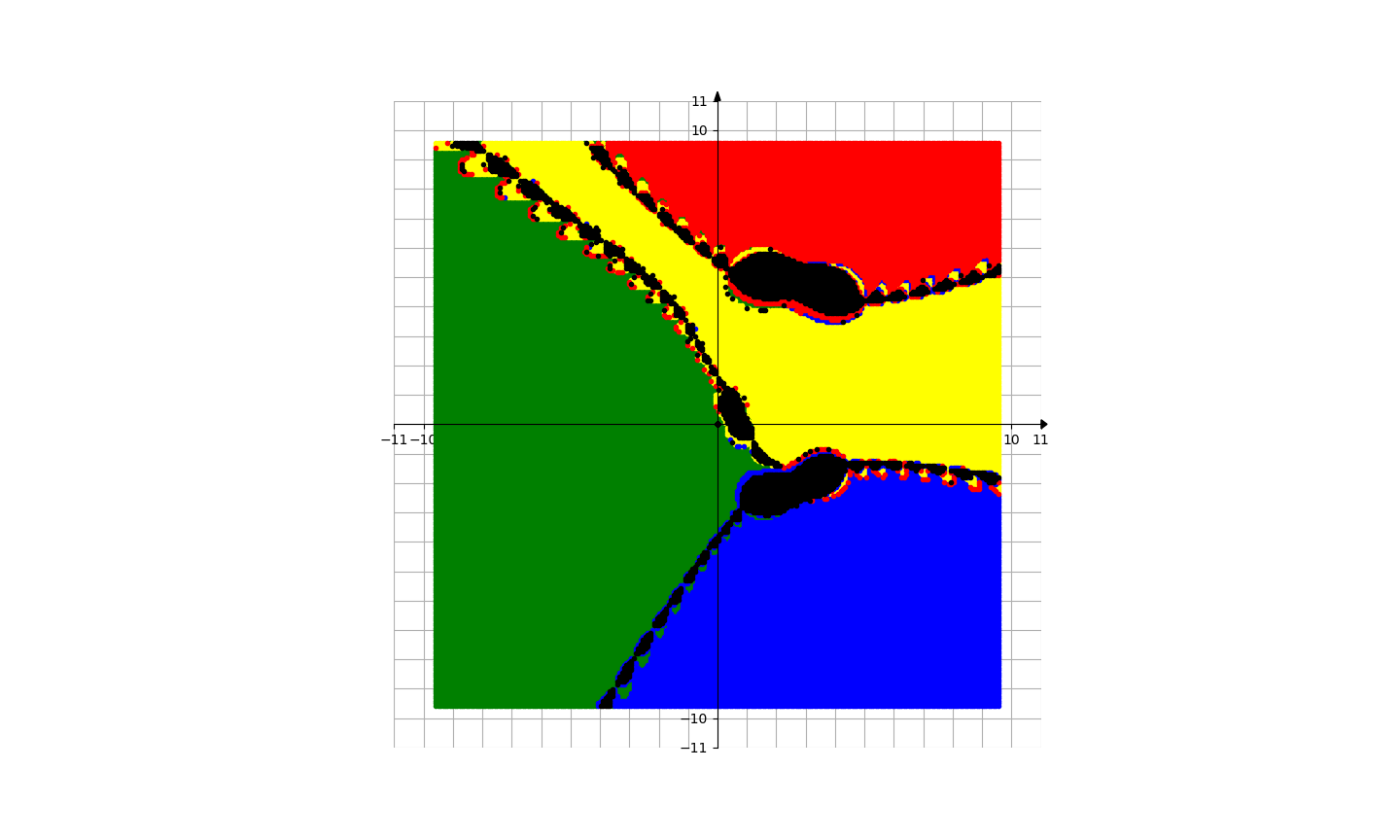}
    \includegraphics[width=3cm]{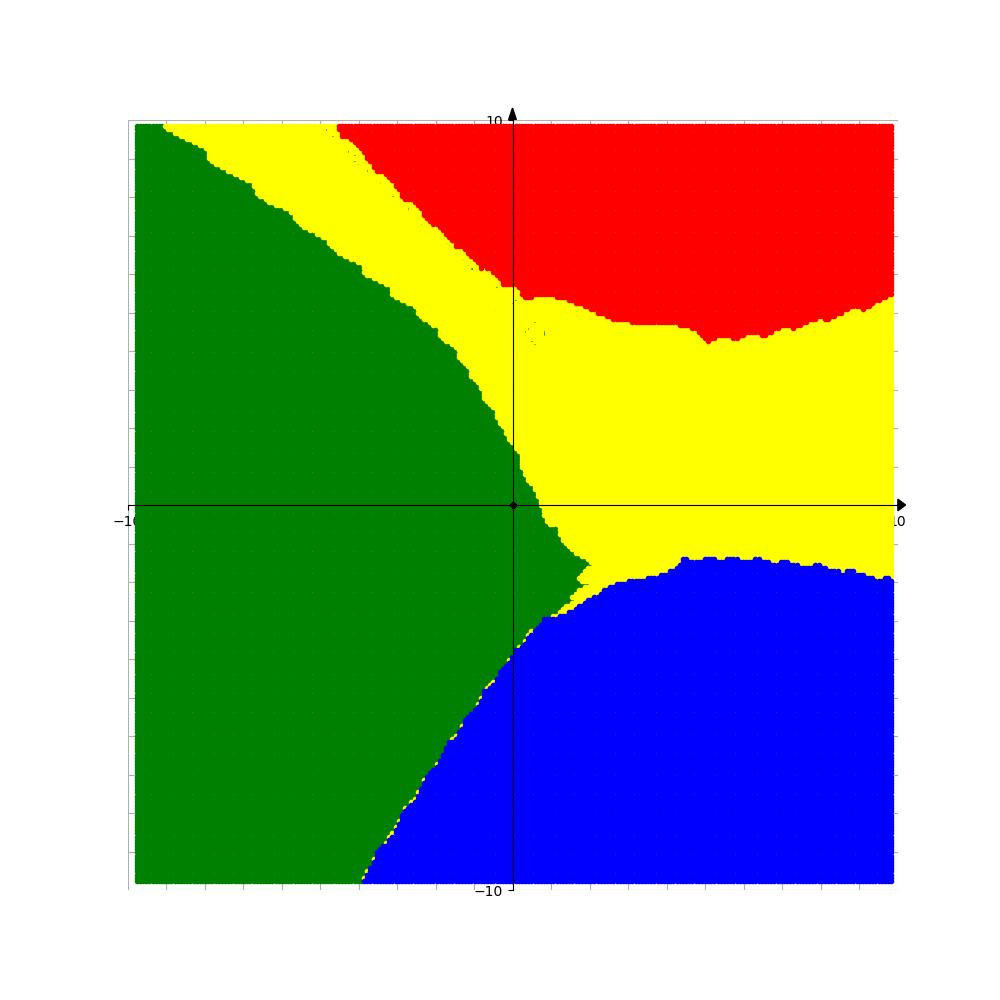}
    
    \bigskip
    \includegraphics[width=3cm]{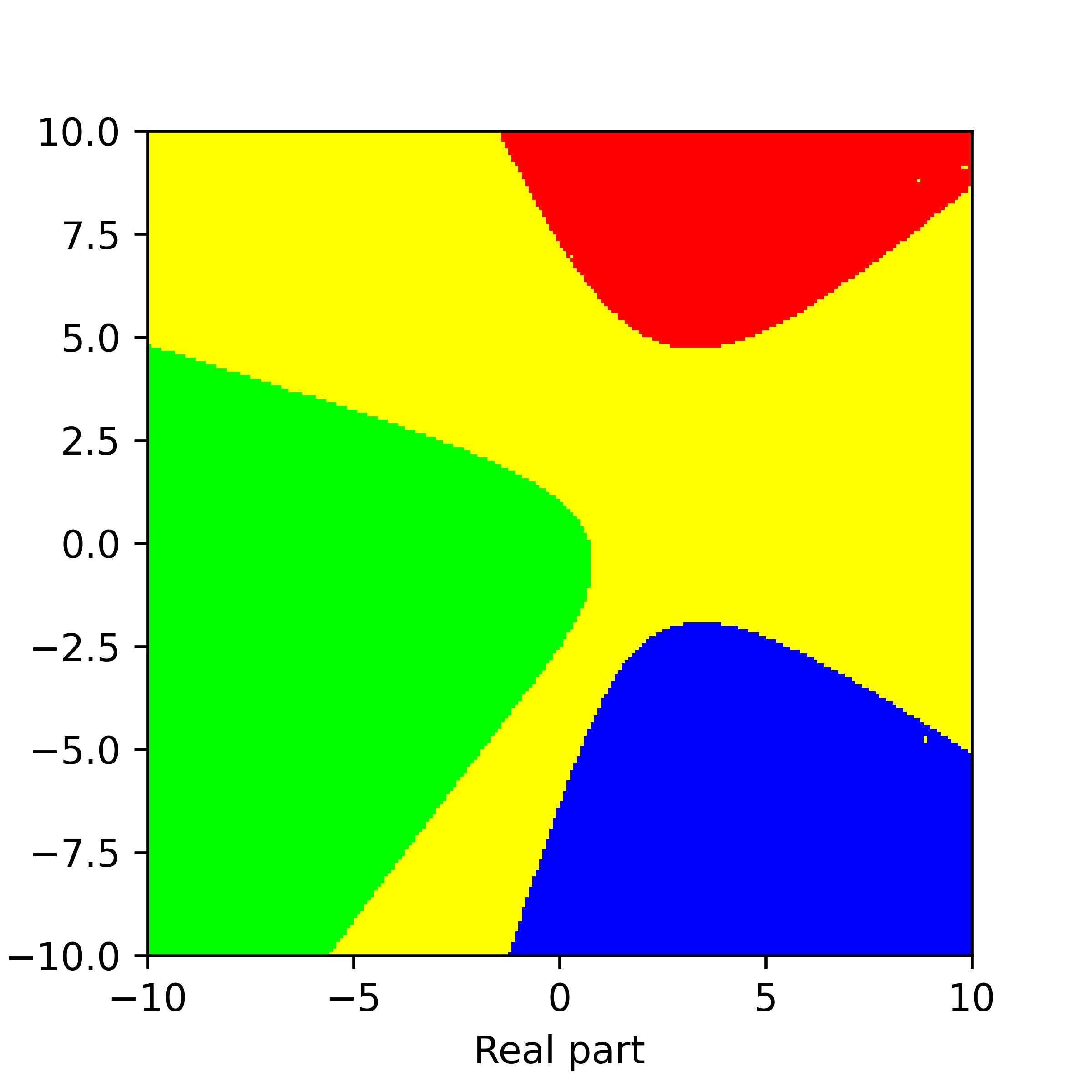}
    \includegraphics[width=3cm]{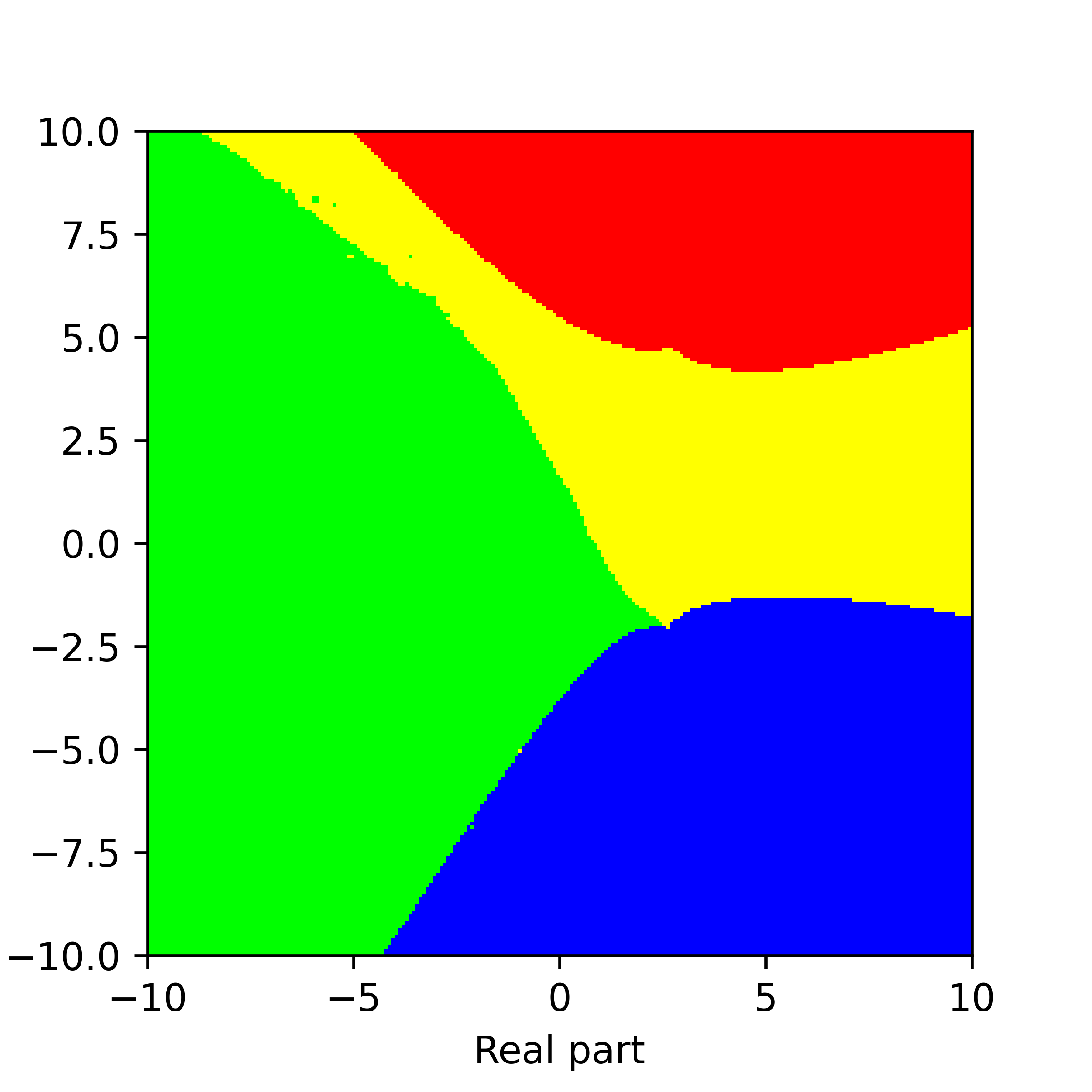}
    \includegraphics[width=3cm]{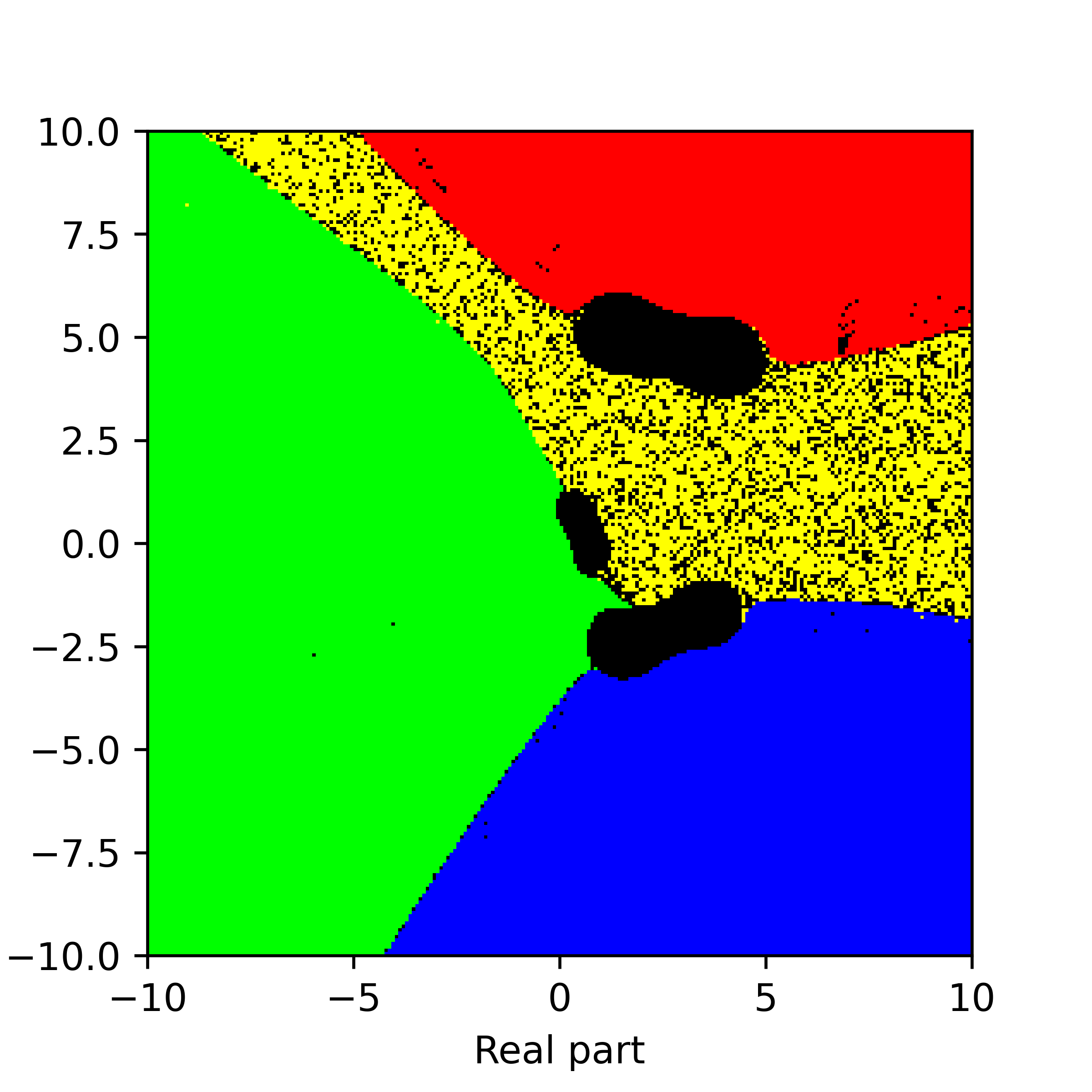}
    
    \caption{Basins of attraction for finding roots of the function $f_{22}$ by different methods. Pictures are referenced to from top to bottom, from left to right. Row 1: left picture is Voronoi's diagram, central picture is for Newton's method, right picture is for Random Relaxed Newton's method. Row 2: left picture is for Newton's method vOptimization, right picture is for BNQN. Row 3: left picture is for Newton's flow, central picture is for Newton's flow vFraction, right picture is for Newton's flow vOptimization. The black points in some of these pictures are those in the basin of attraction of critical points of $f_{22}$.}
    \label{fig:f22}
\end{figure}

\subsubsection{Transcendental functions}

\begin{figure}
    \centering
    \includegraphics[width=5cm]{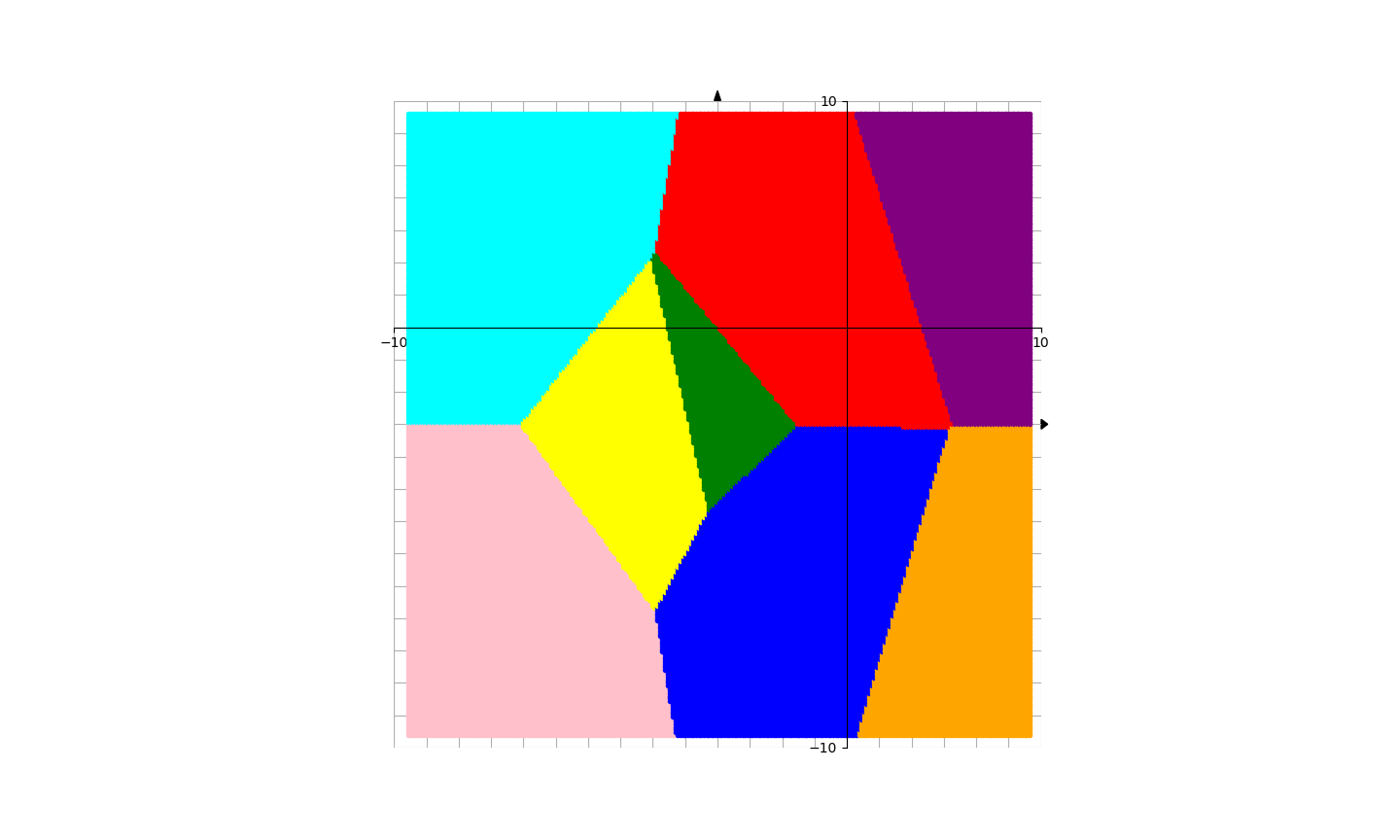}
   
    \bigskip
    
    \includegraphics[width=5cm]{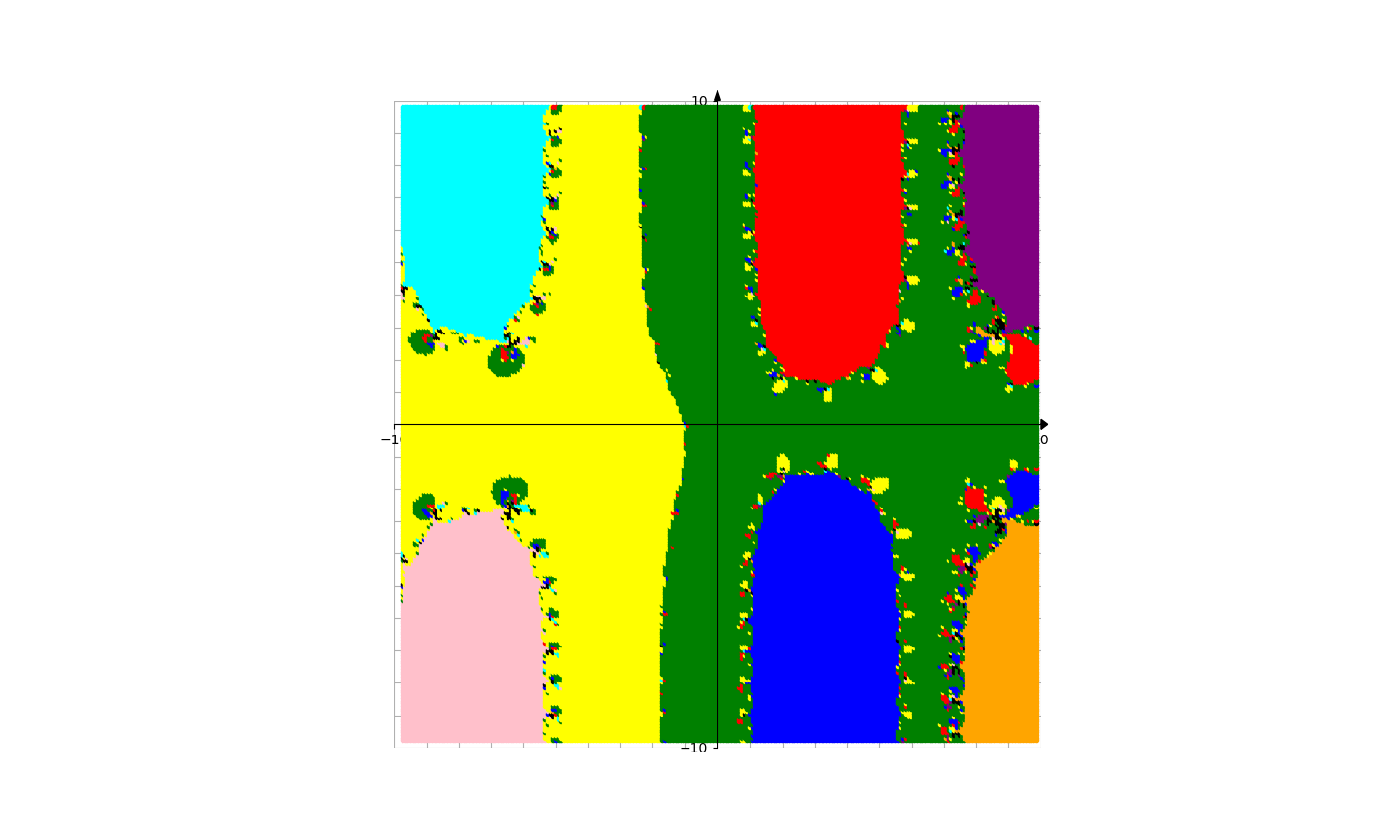}
    \includegraphics[width=5cm]{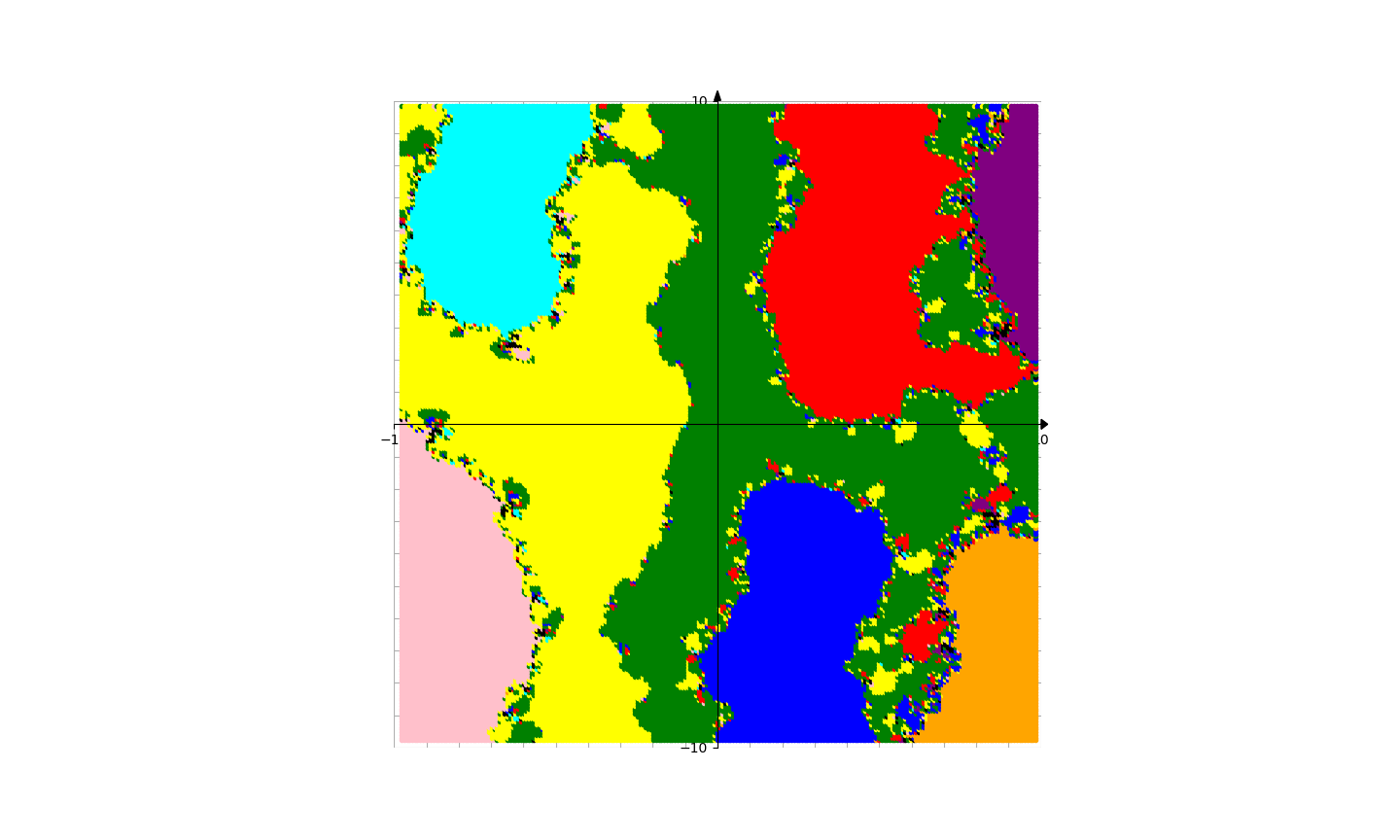}

    \bigskip
    \includegraphics[width=5cm]{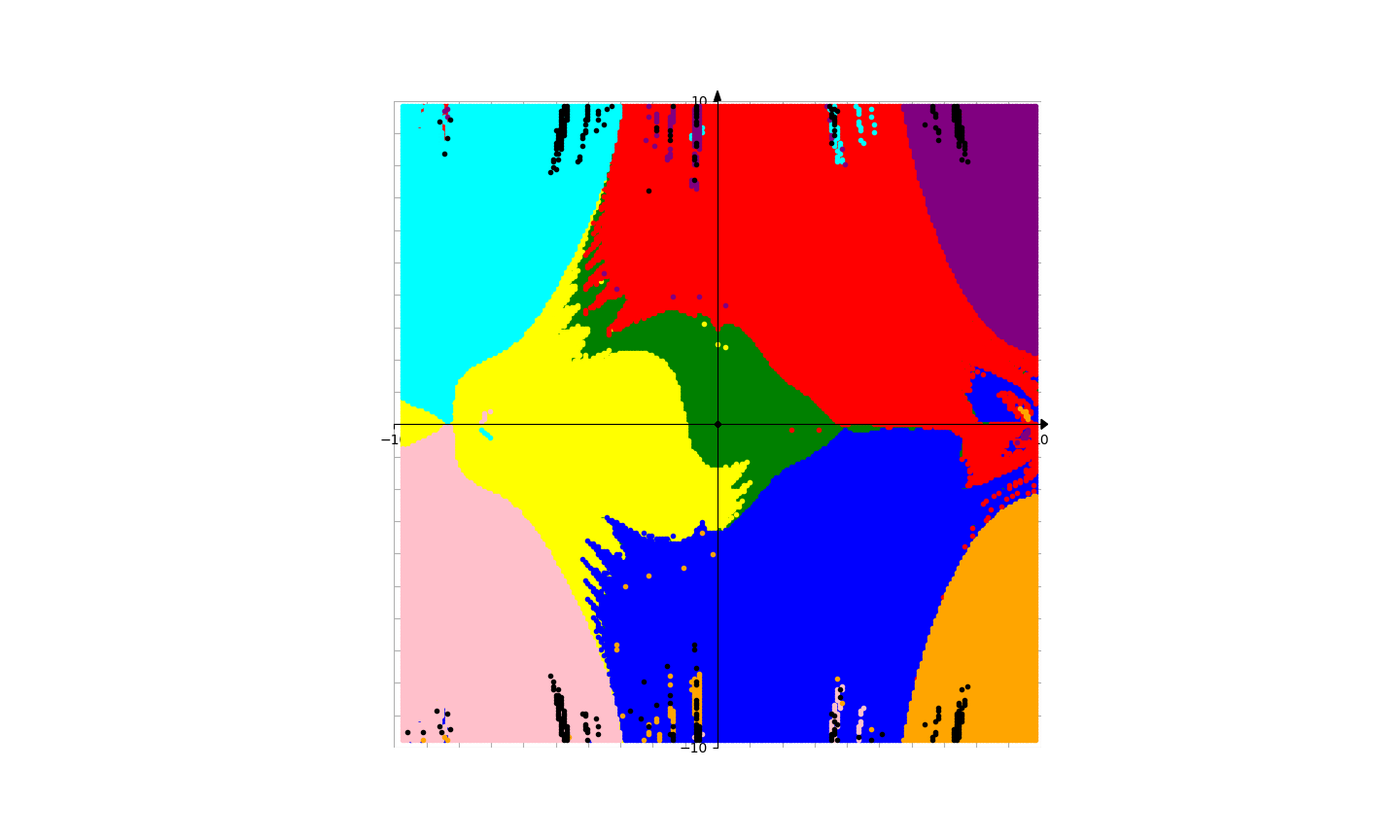}
    \includegraphics[width=5cm]{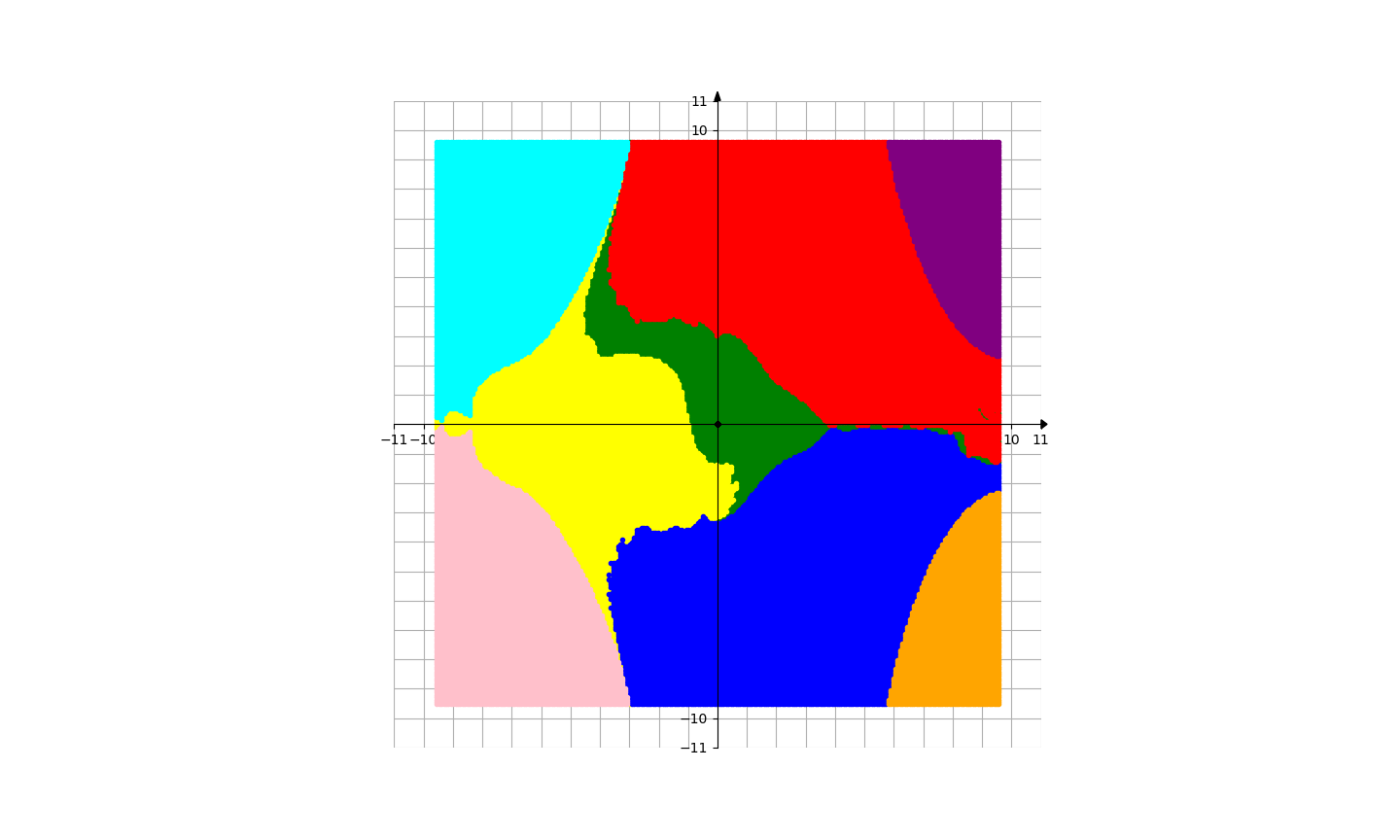}
    
    \bigskip
    \includegraphics[width=3cm]{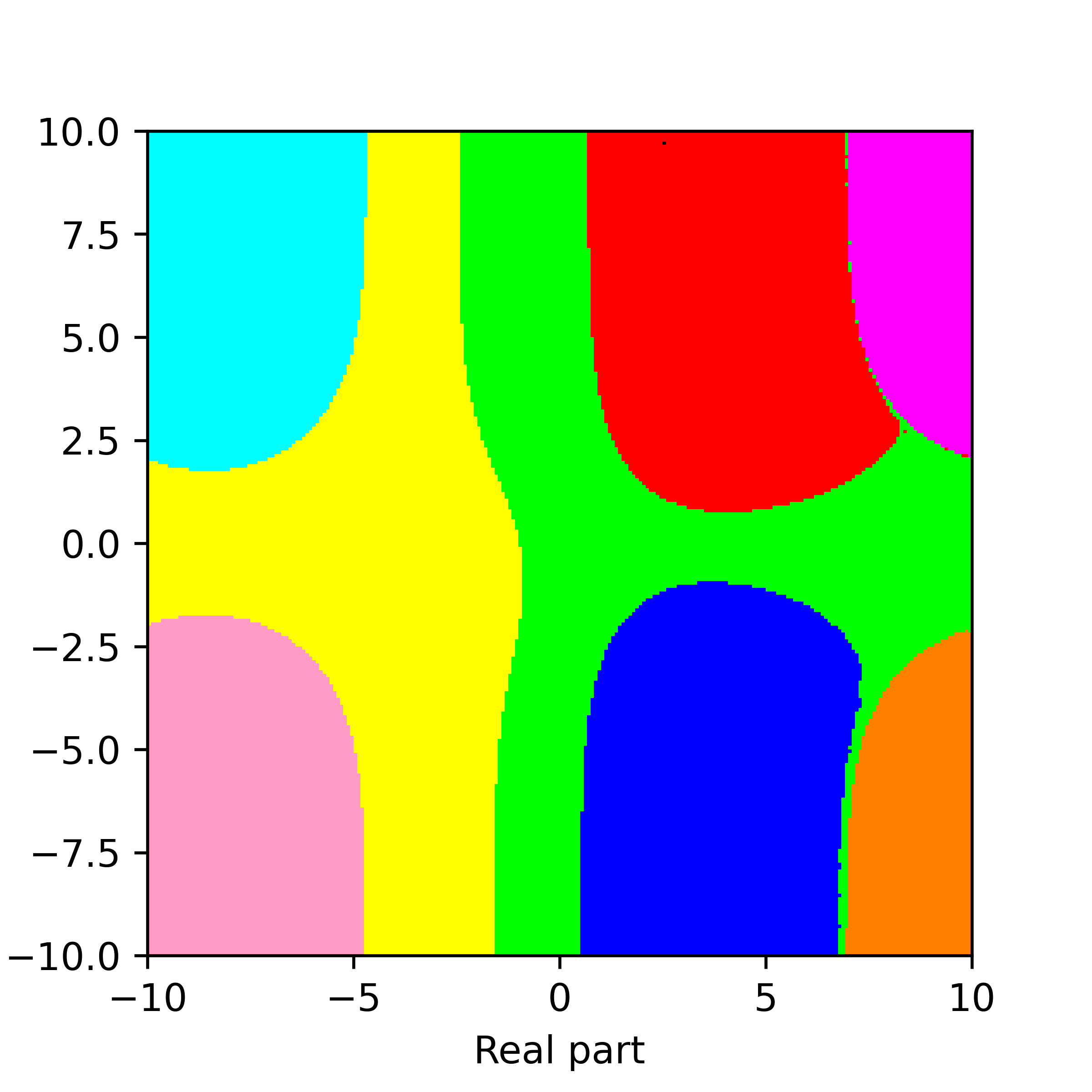}
    \includegraphics[width=3cm]{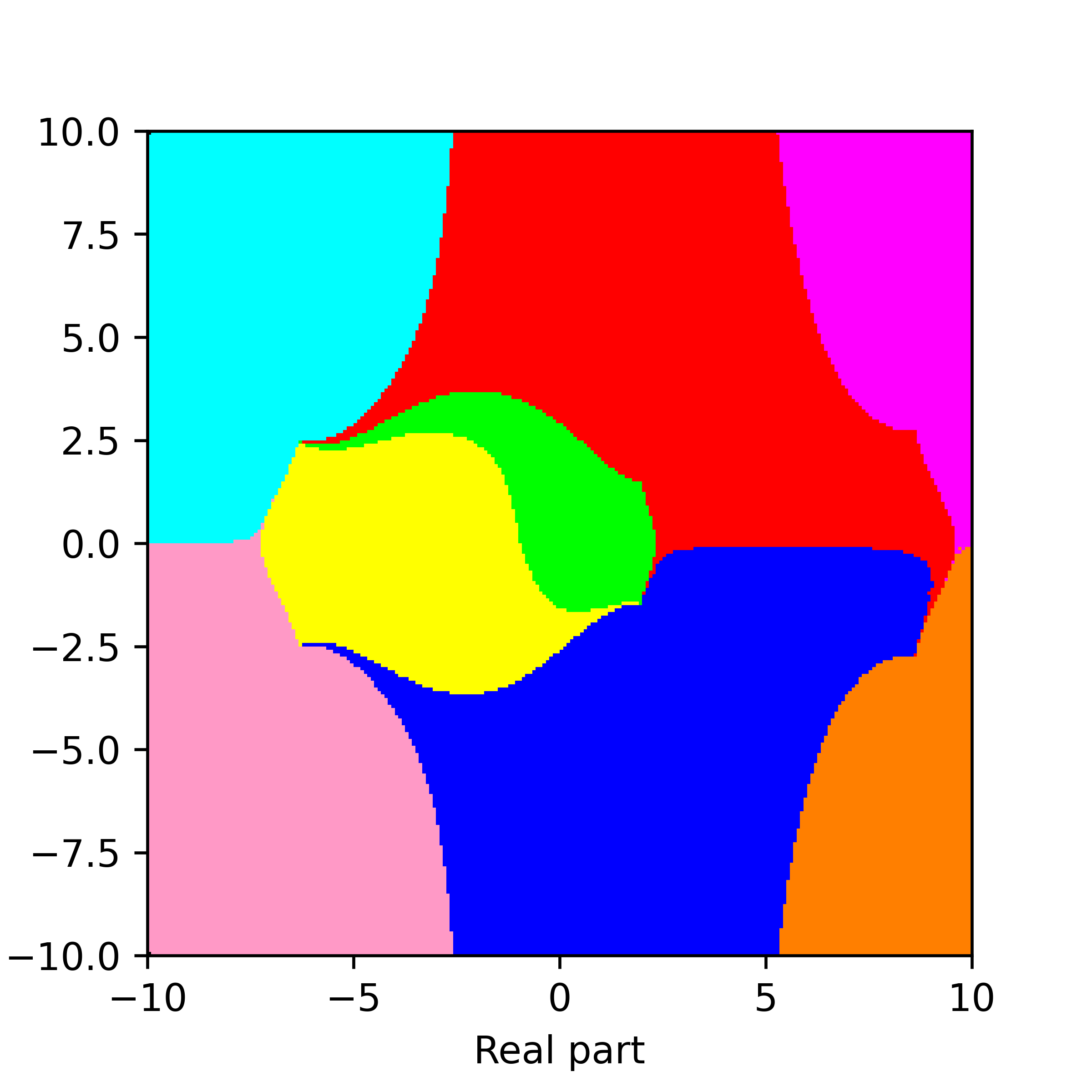}

    \caption{Basins of attraction for finding roots of the function $f_{23}$ by different methods. Pictures are referenced to from top to bottom, from left to right. Row 1: Voronoi's diagram (using only 8 roots inside the domain $(-10,10)\times (-10,10)$). Row 2: left picture is for Newton's method, right picture is for Random Relaxed Newton's method. Row 3: left picture is for BNQN, right picture is for BNQN v2. Row 4: left picture is for Newton's flow, right picture is for Newton's flow vFraction. The method Newton's method vOptimization encounters errors, while  Newton's flow vOptimization takes too long time to finish.}
    \label{fig:f23}
\end{figure}

\begin{figure}
    \centering
    \includegraphics[width=5cm]{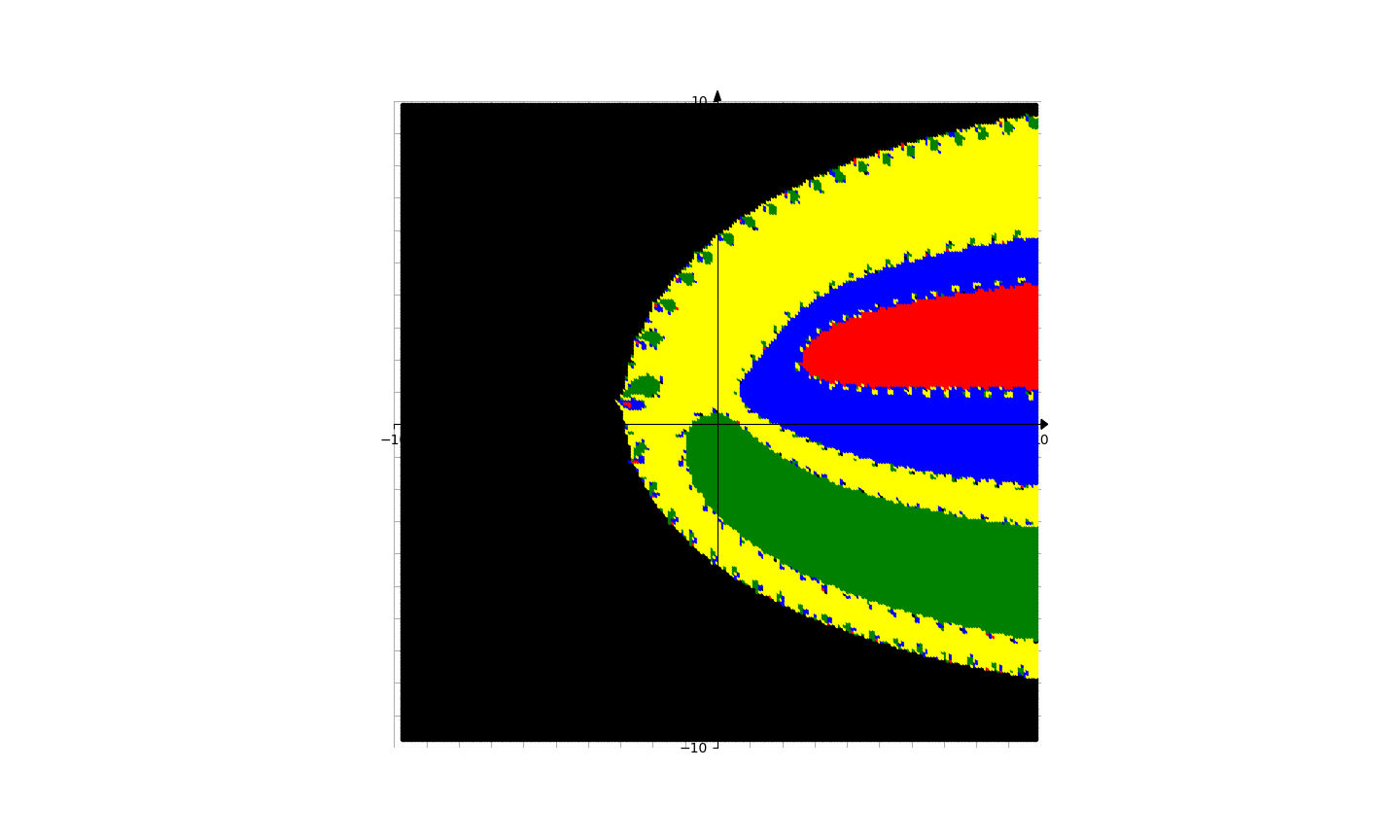}
     \includegraphics[width=5cm]{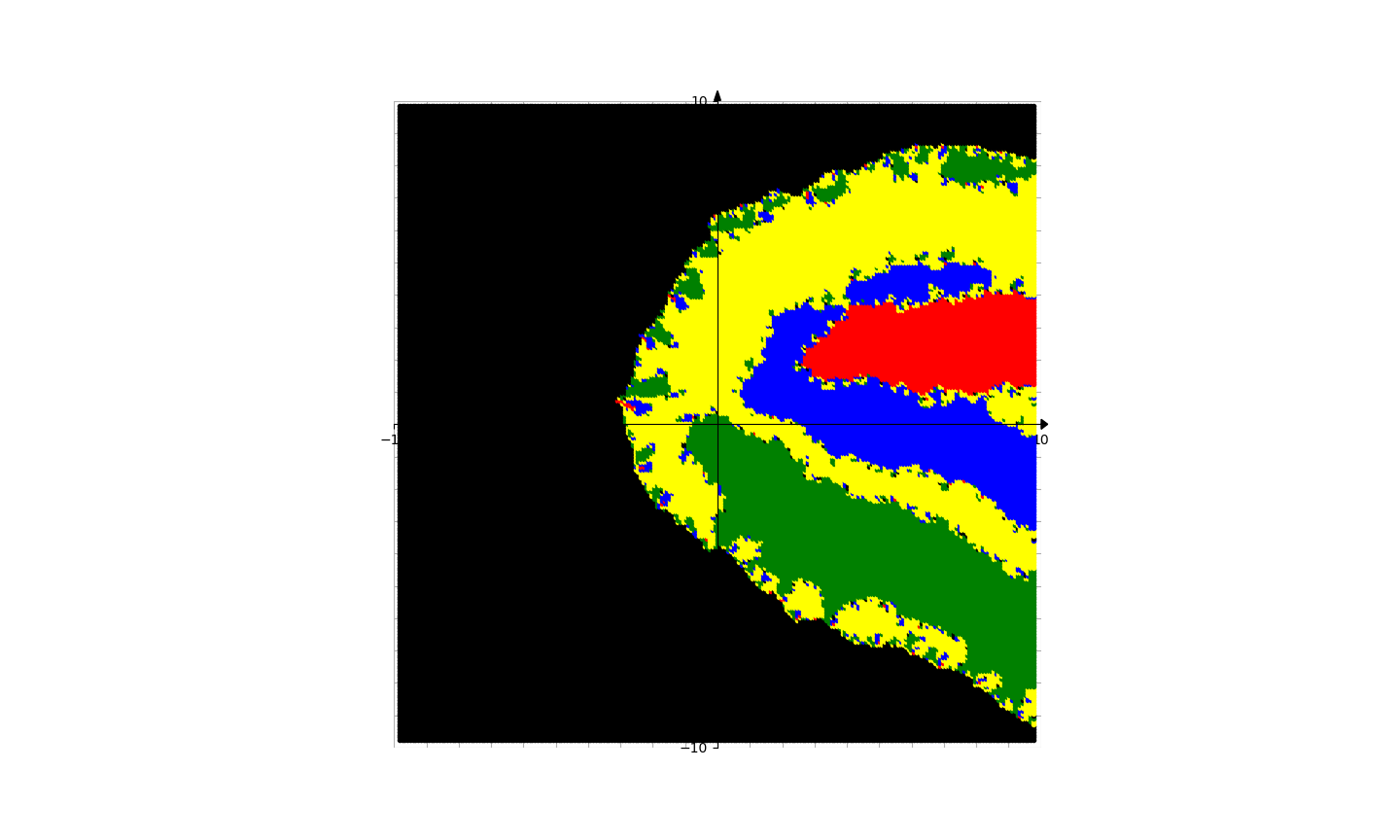}
   
    \bigskip
    
    \includegraphics[width=5cm]{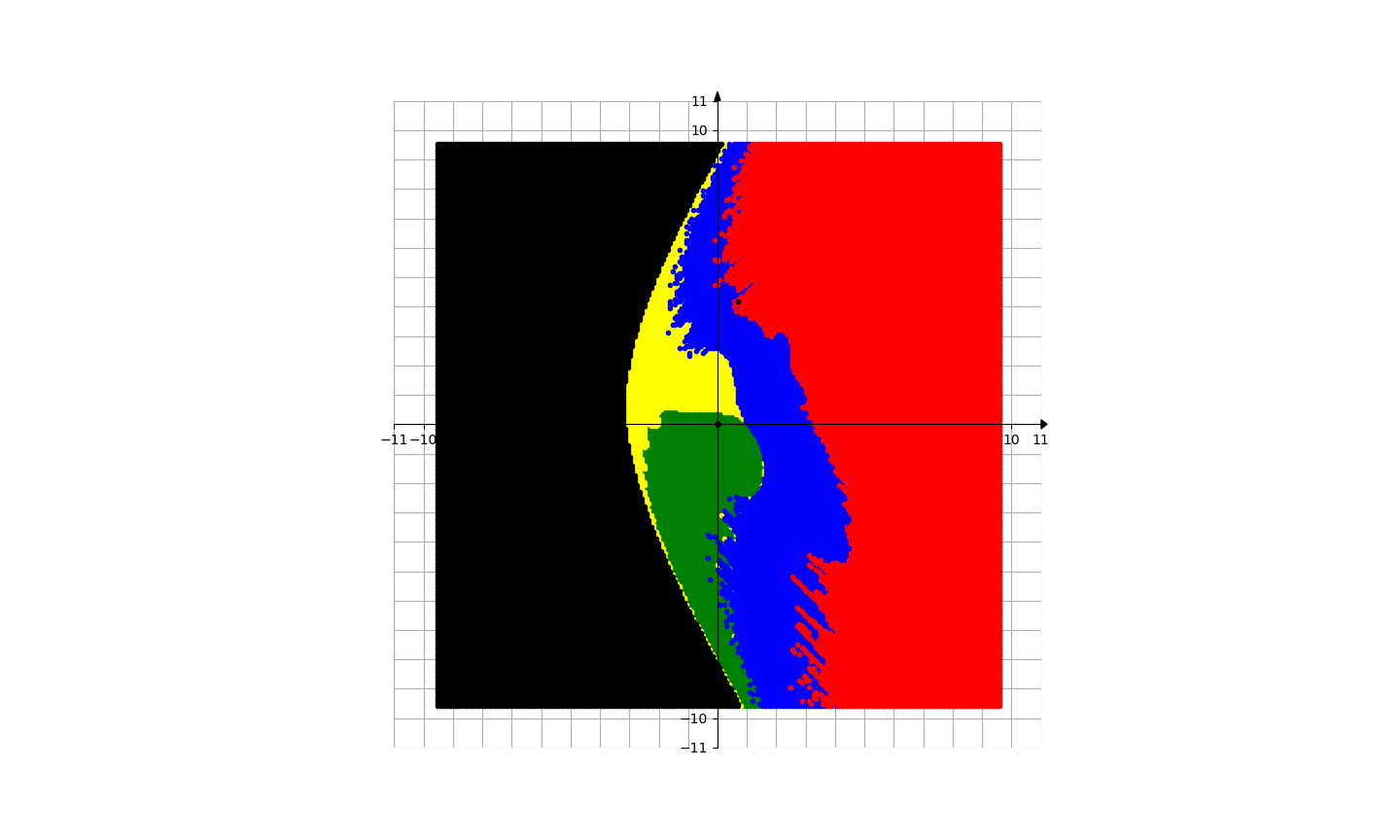}
    \includegraphics[width=5cm]{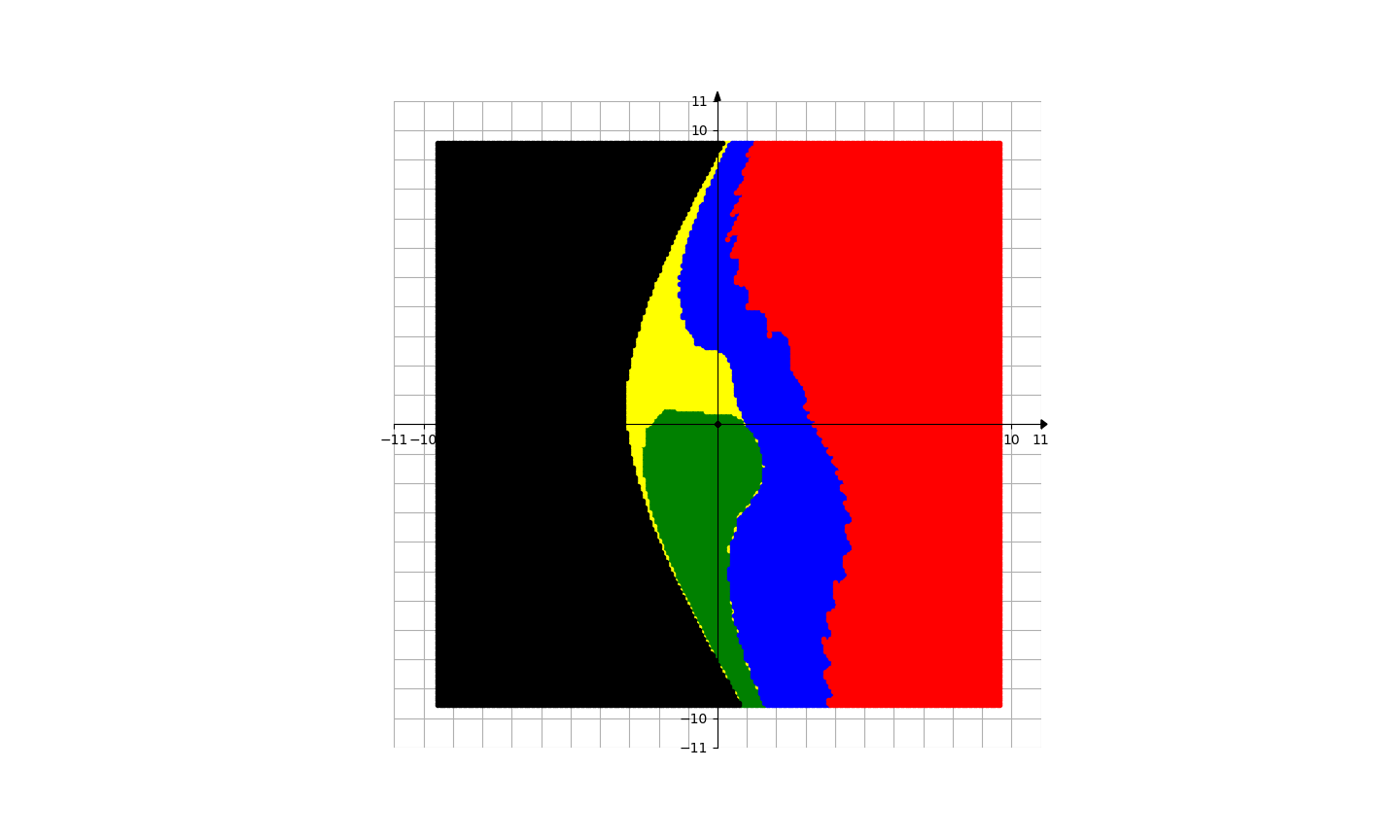}

    \bigskip
    \includegraphics[width=3cm]{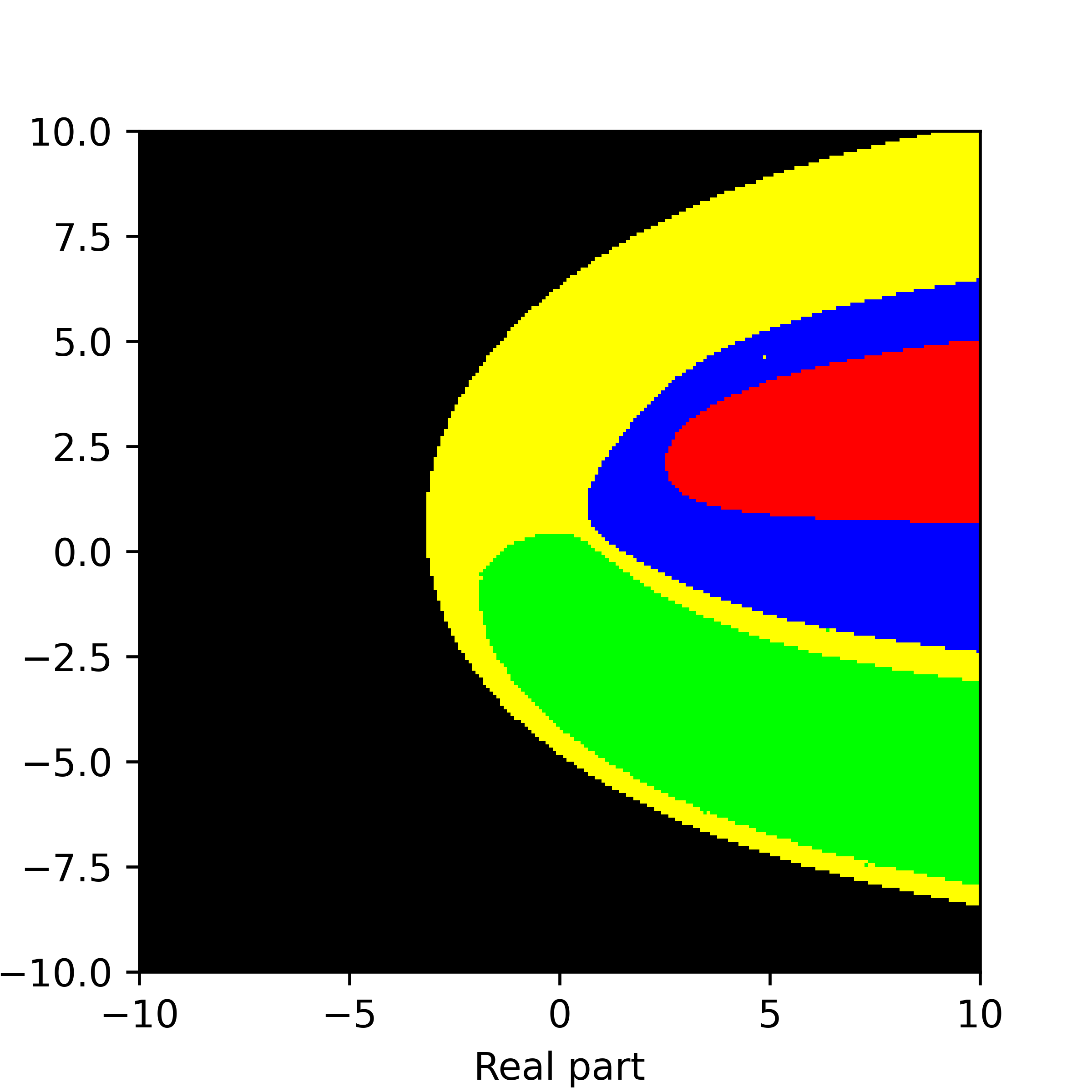}
    \includegraphics[width=3cm]{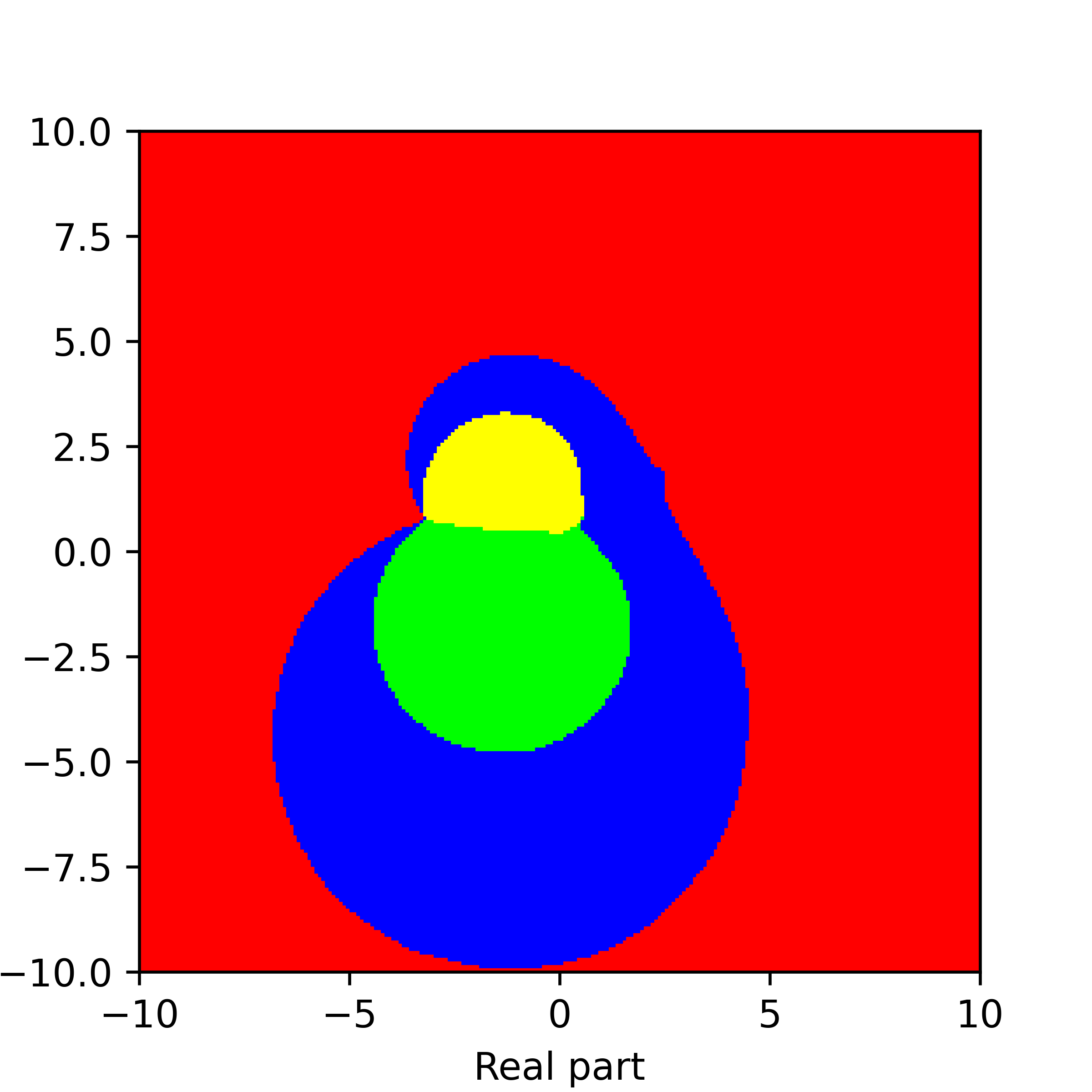}
    
    \bigskip
    \includegraphics[width=5cm]{VoronoiF7.png}

    \caption{Basins of attraction for finding roots of the function $f_{7}e^z$ by different methods. Pictures are referenced to from top to bottom, from left to right. Row 1: left picture is for Newton's method, right picture is for Random Relaxed Newton's method. Row 2: left picture is for BNQN, right picture is for BNQN v2. Row 3: left picture is for Newton's flow, right picture is for Newton's flow vFraction. Row 4: Voronoi's diagram. The method Newton's method vOptimization encounters errors, while  Newton's flow vOptimization takes too long time to finish. }
    \label{fig:f7ez}
\end{figure}

\begin{figure}
    \centering
    \includegraphics[width=5cm]{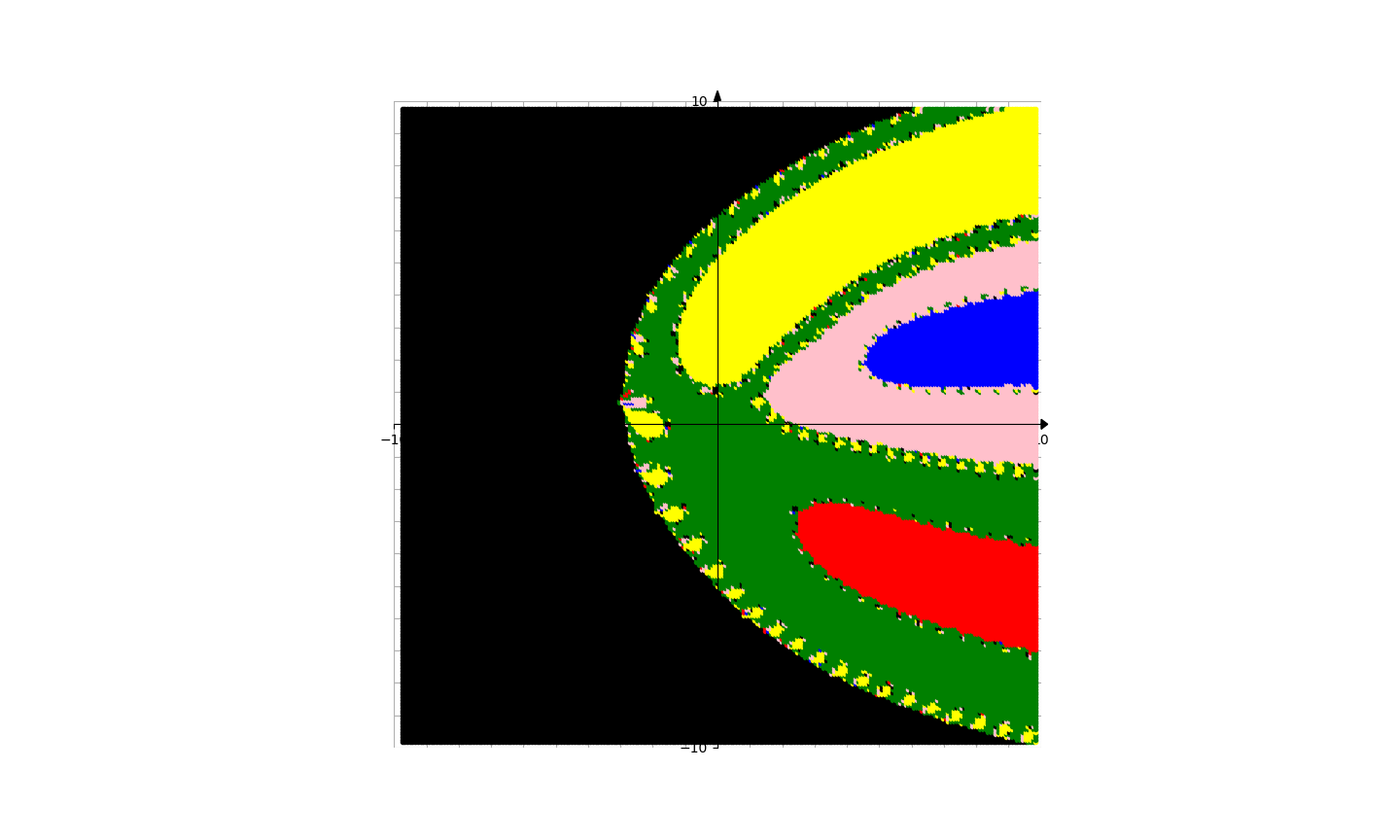}
     \includegraphics[width=5cm]{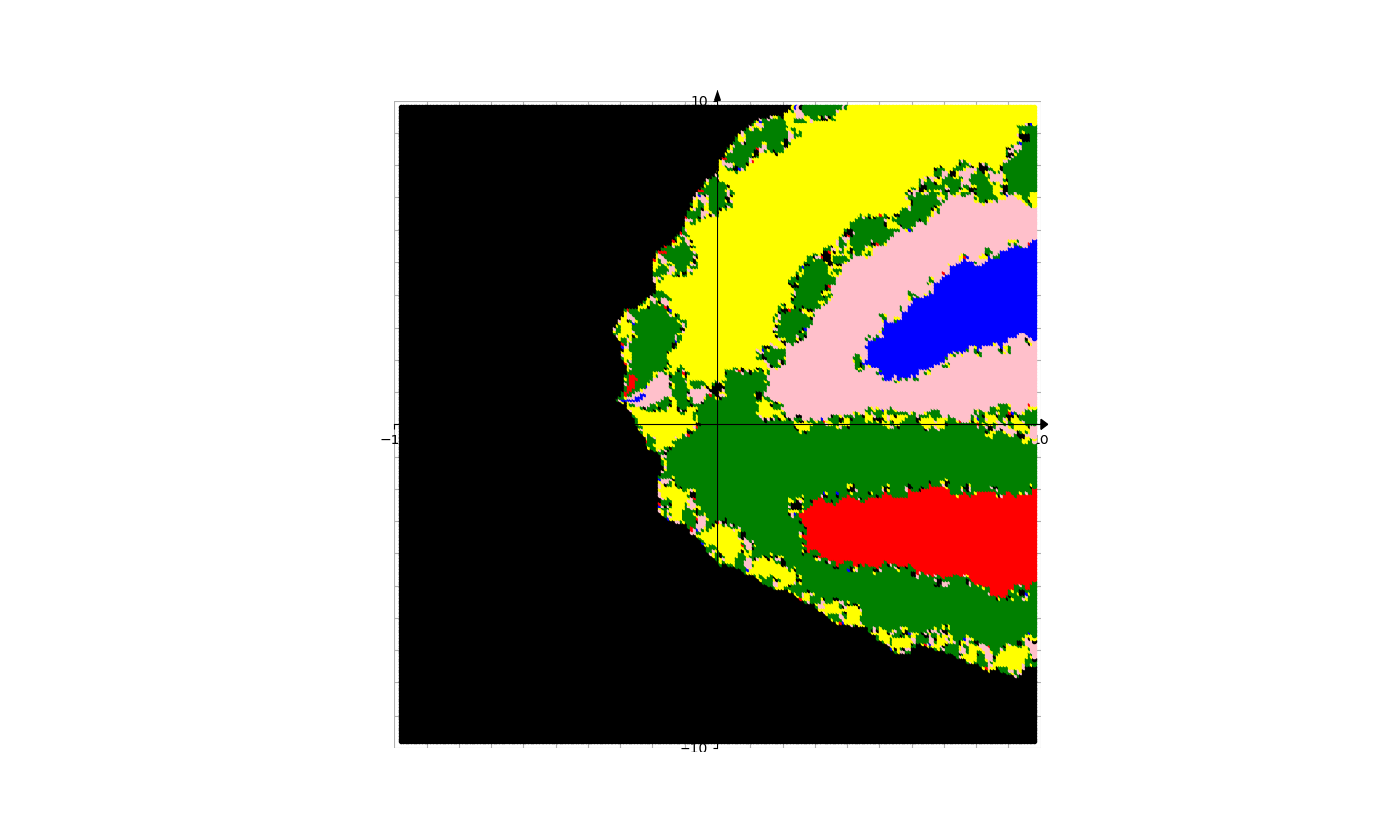}
   
    \bigskip
    
    \includegraphics[width=5cm]{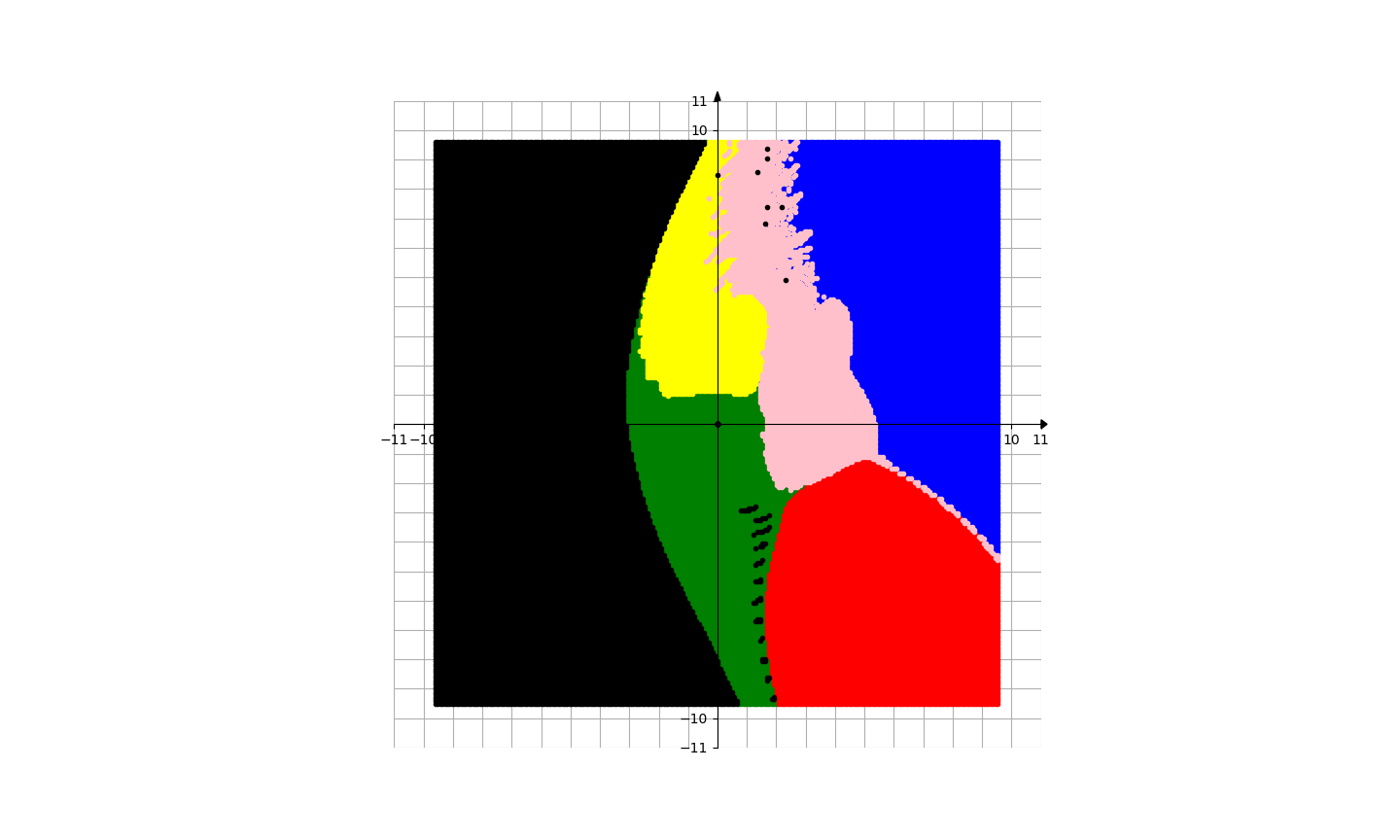}
    \includegraphics[width=5cm]{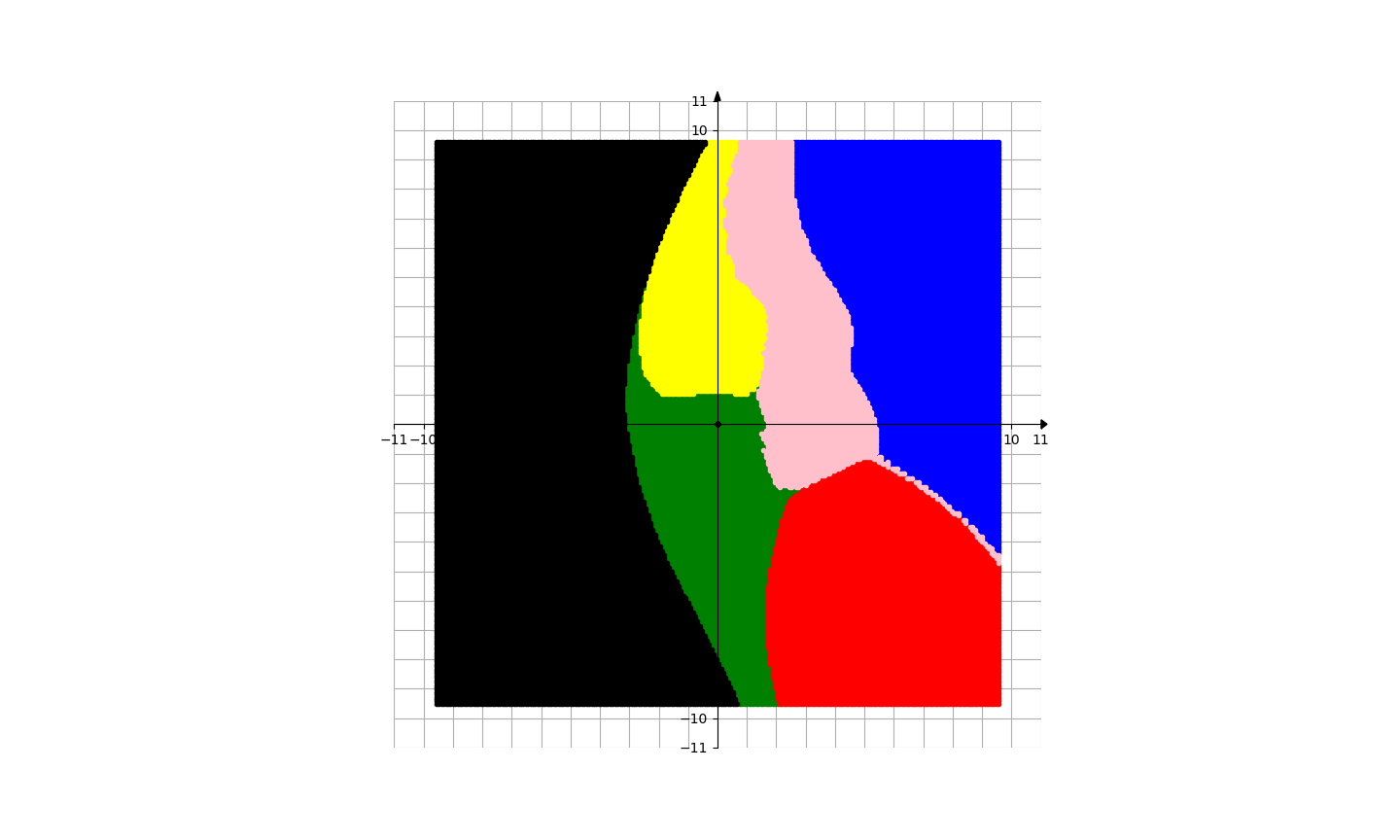}

    \bigskip
    \includegraphics[width=3cm]{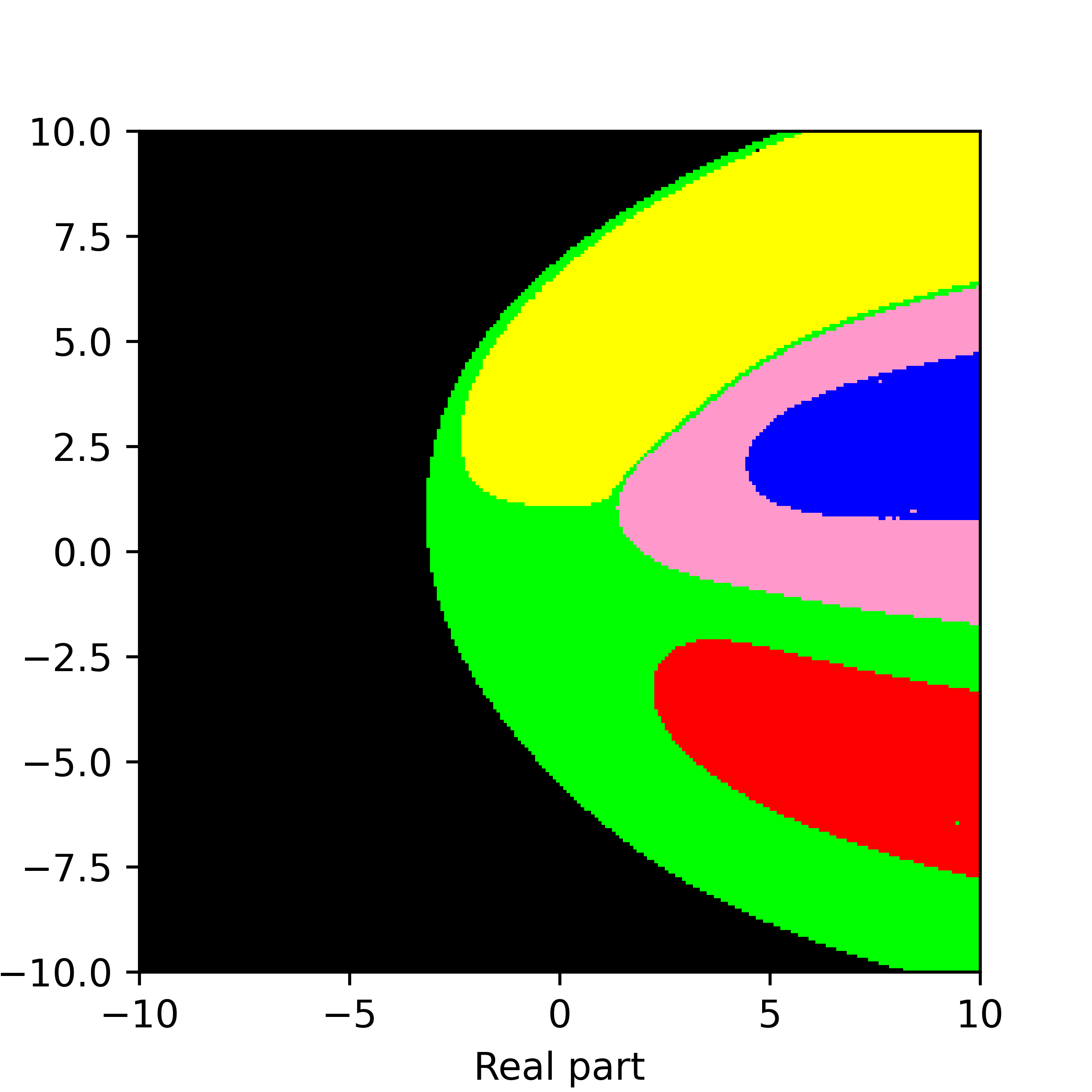}
    \includegraphics[width=3cm]{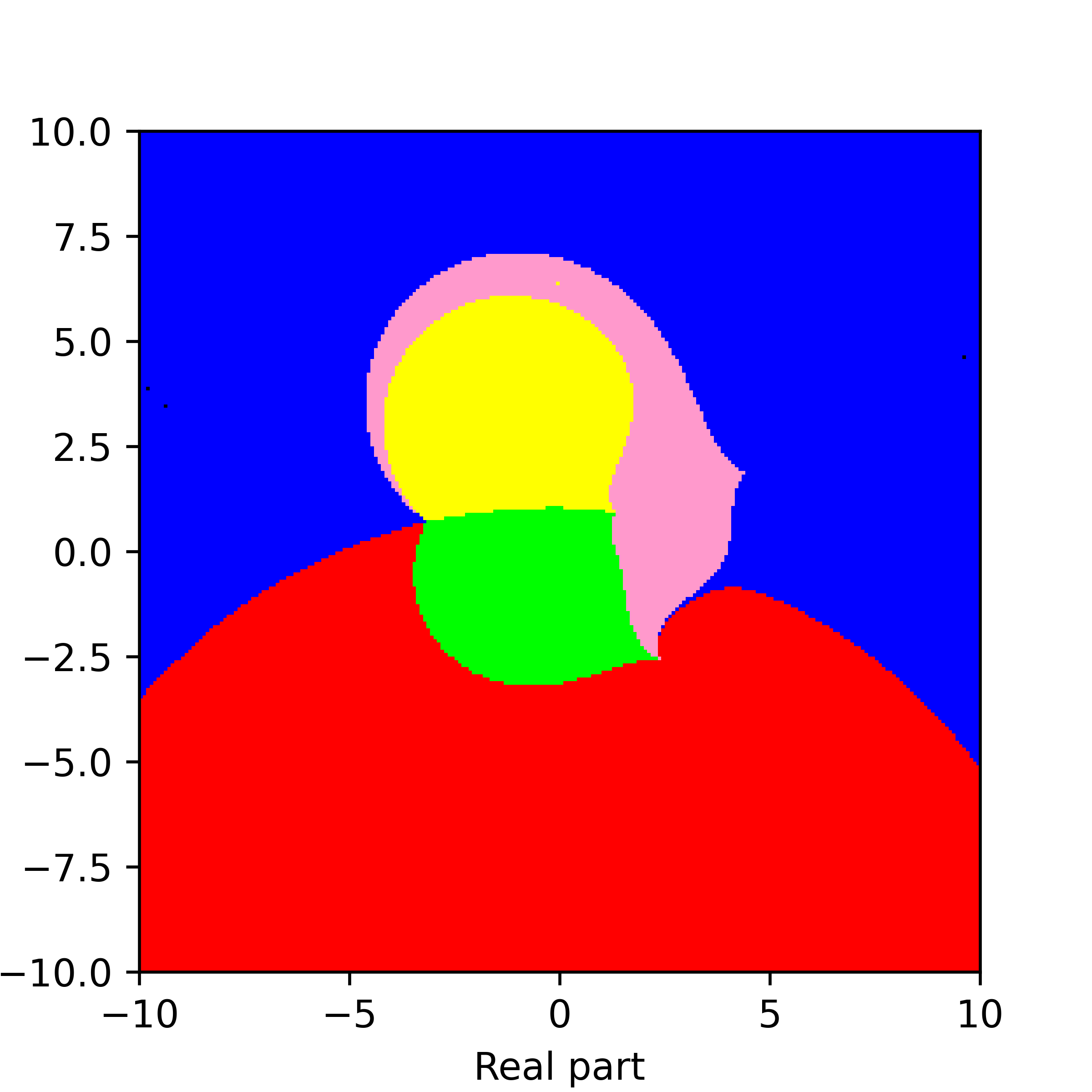}
    
    \bigskip
    \includegraphics[width=5cm]{VoronoiF17.png}

    \caption{Basins of attraction for finding roots of the function $f_{17}e^z$ by different methods. Pictures are referenced to from top to bottom, from left to right. Row 1: left picture is for Newton's method, right picture is for Random Relaxed Newton's method. Row 2: left picture is for BNQN, right picture is for BNQN v2. Row 3: left picture is for Newton's flow, right picture is for Newton's flow vFraction. The method Newton's method vOptimization encounters errors, while  Newton's flow vOptimization takes too long time to finish.}
    \label{fig:f17ez}
\end{figure}

\subsubsection{Quotient of a function and its derivative}

\begin{figure}
    \centering
    \includegraphics[width=3cm]{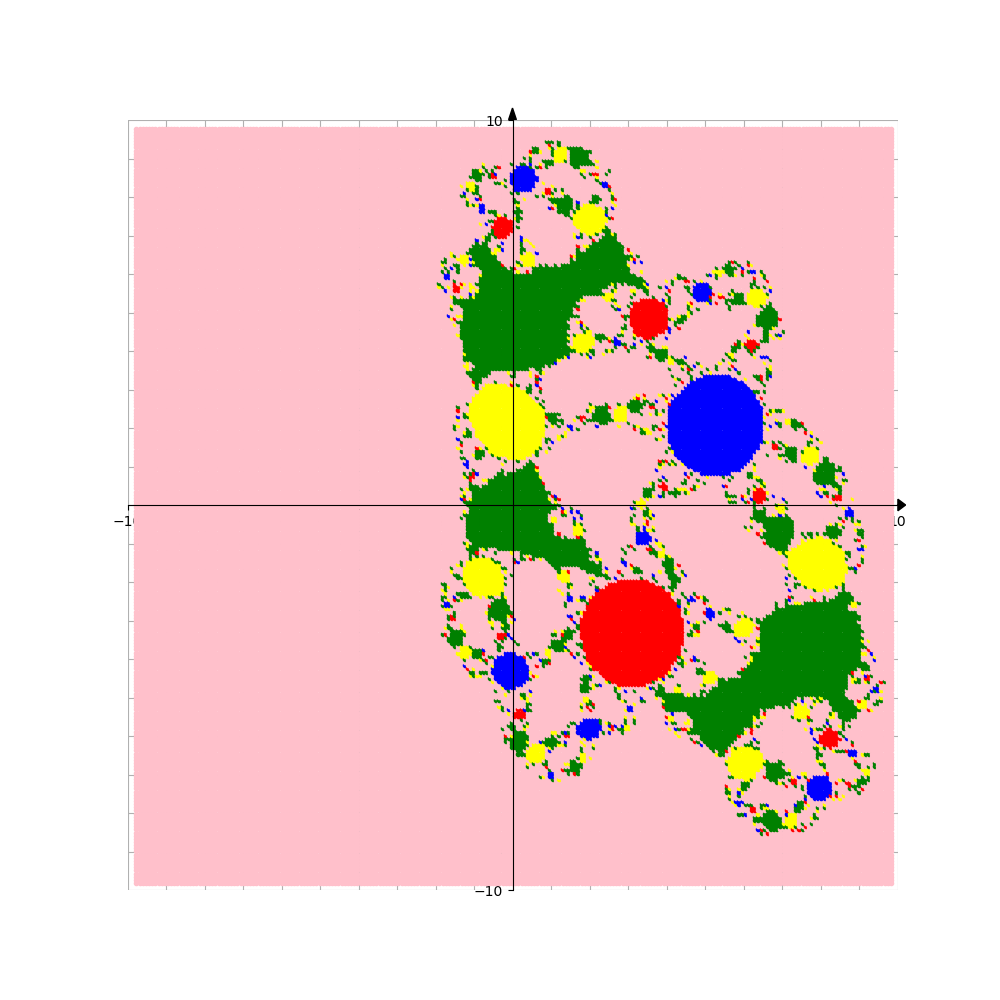}
     \includegraphics[width=3cm]{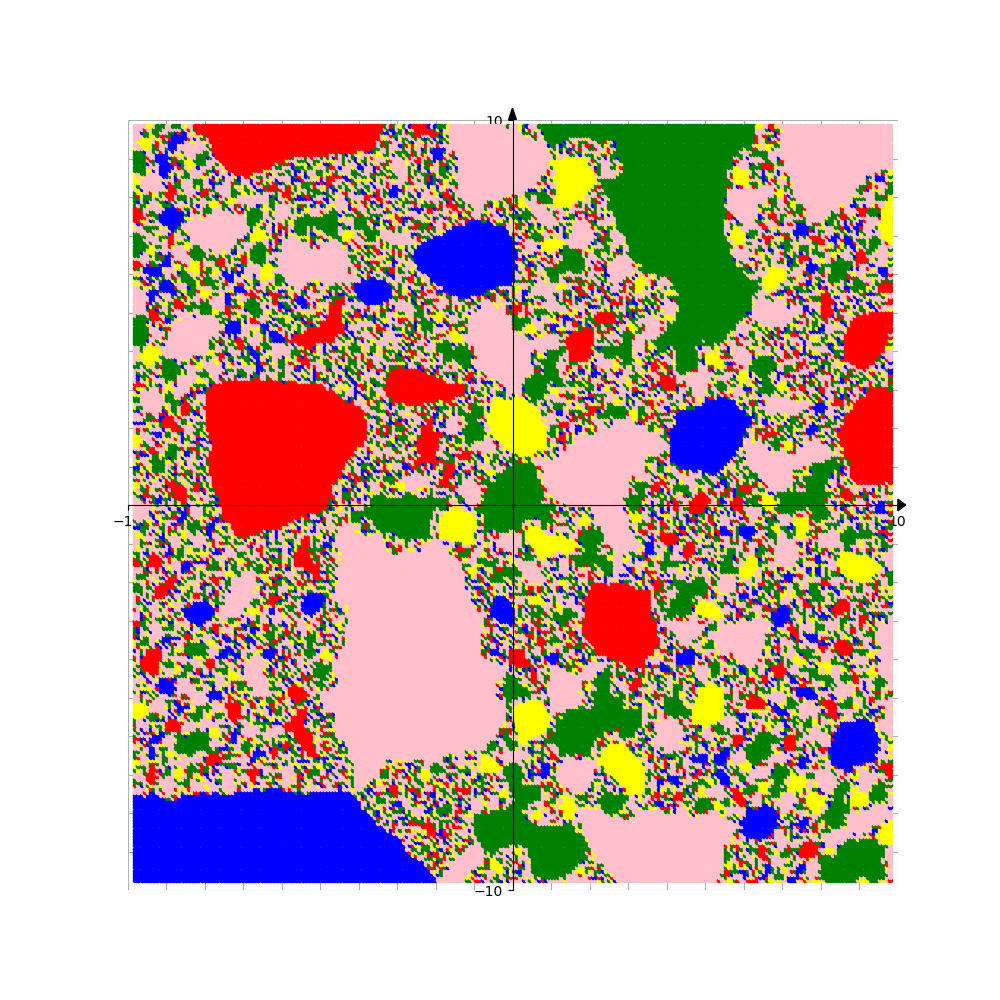}
        \includegraphics[width=5cm]{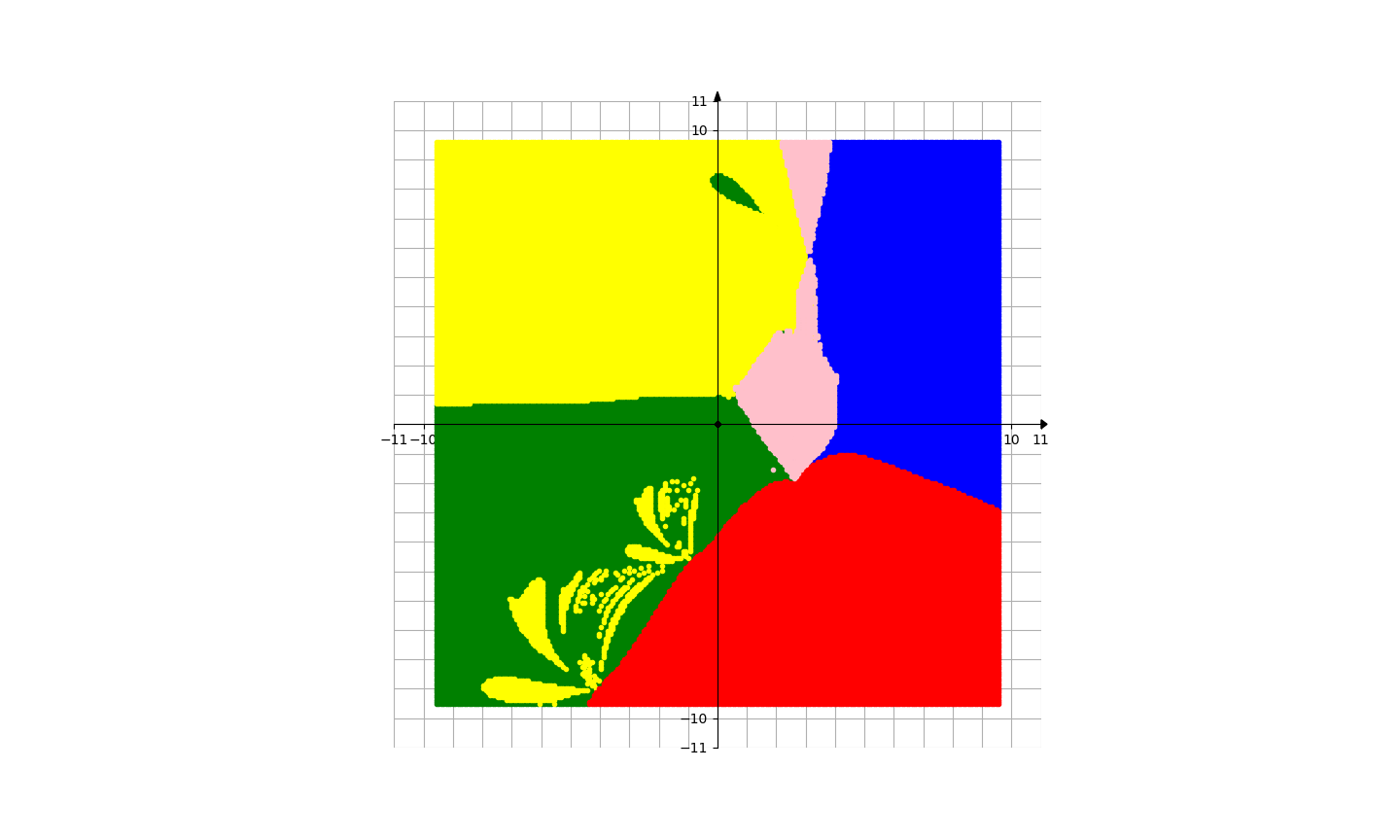}

    \caption{Basins of attraction for finding roots of the function $f_{17}/f_{17}'$ by different methods. Left picture is for Newton's method, central picture is for Random Relaxed Newton's method, right picture is for BNQN. }
    \label{fig:f17Fraction}
\end{figure}

\begin{figure}
    \centering
    \includegraphics[width=3cm]{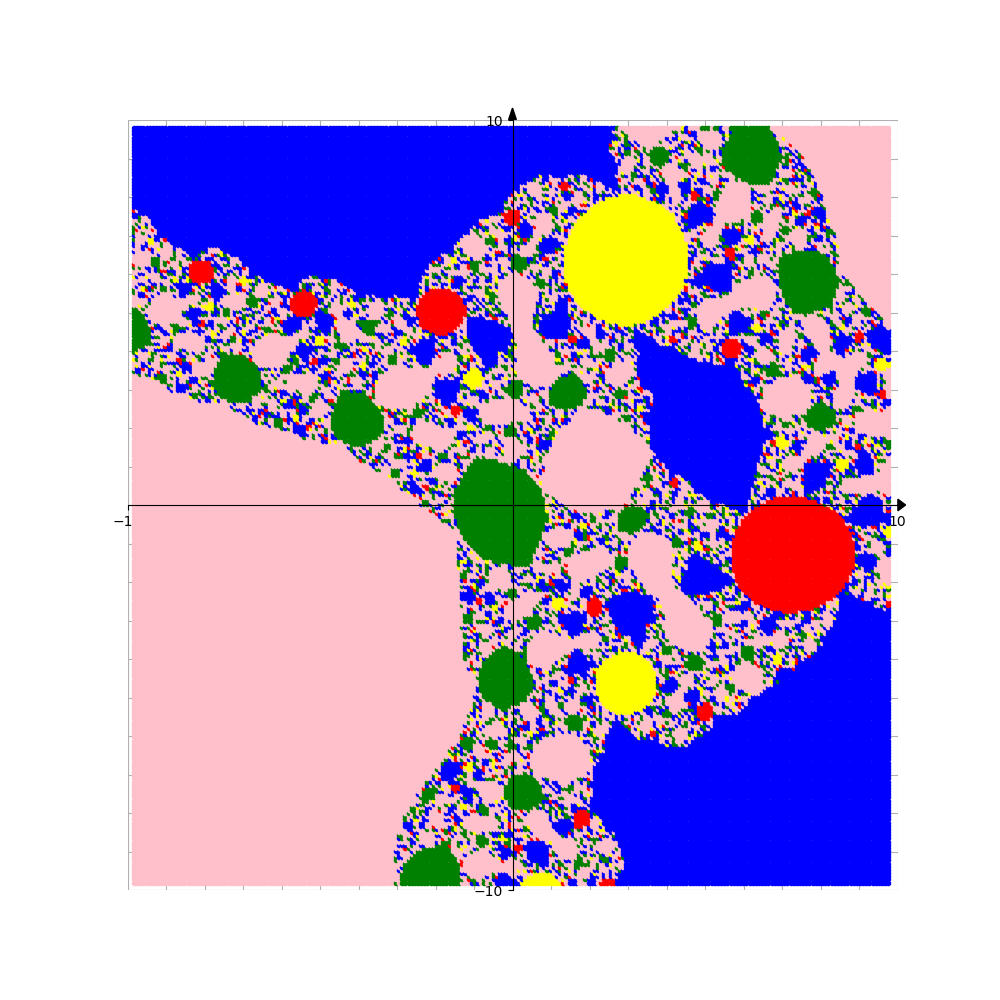}
     \includegraphics[width=3cm]{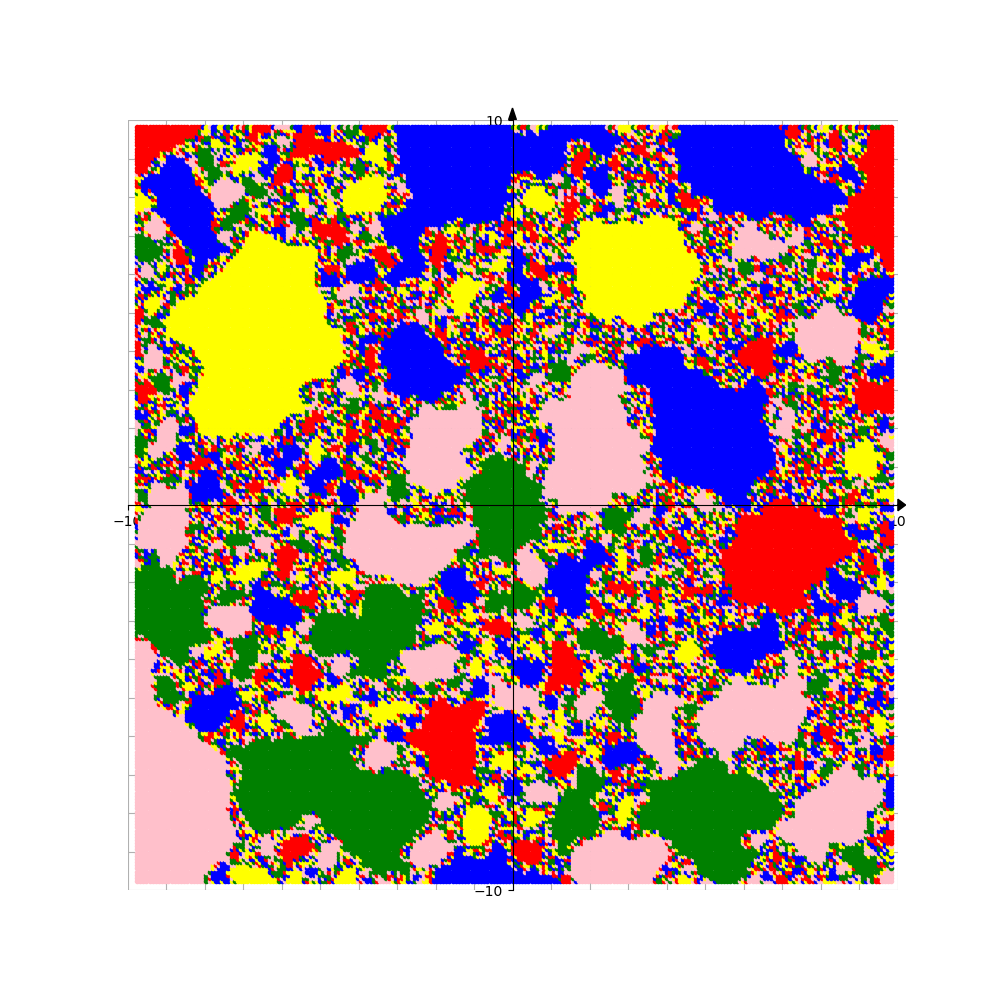}
        \includegraphics[width=5cm]{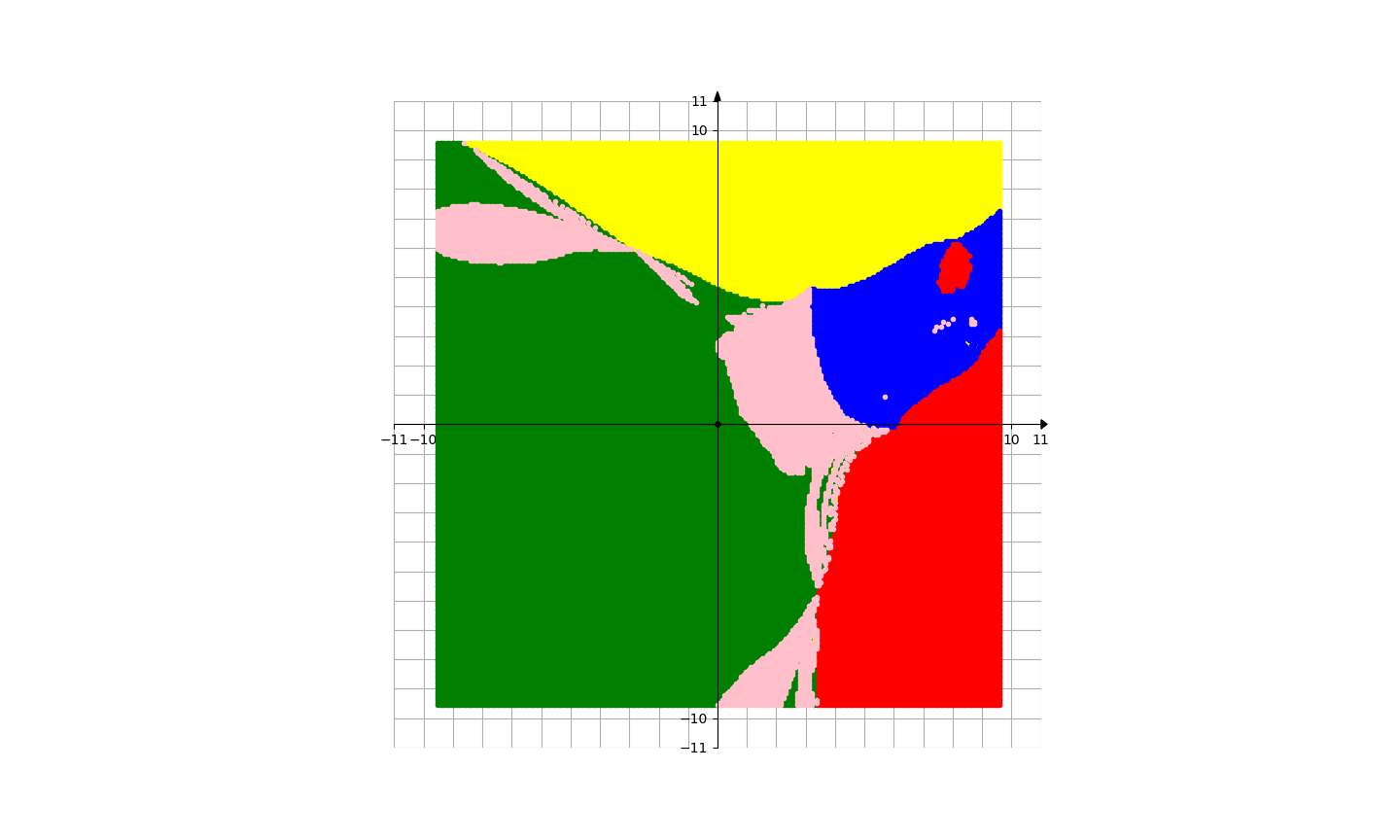}

    \caption{Basins of attraction for finding roots of the function $f_{18}/f_{18}'$ by different methods. Left picture is for Newton's method, central picture is for Random Relaxed Newton's method, right picture is for BNQN. }
    \label{fig:f18Fraction}
\end{figure}

\begin{figure}
    \centering
    \includegraphics[width=3cm]{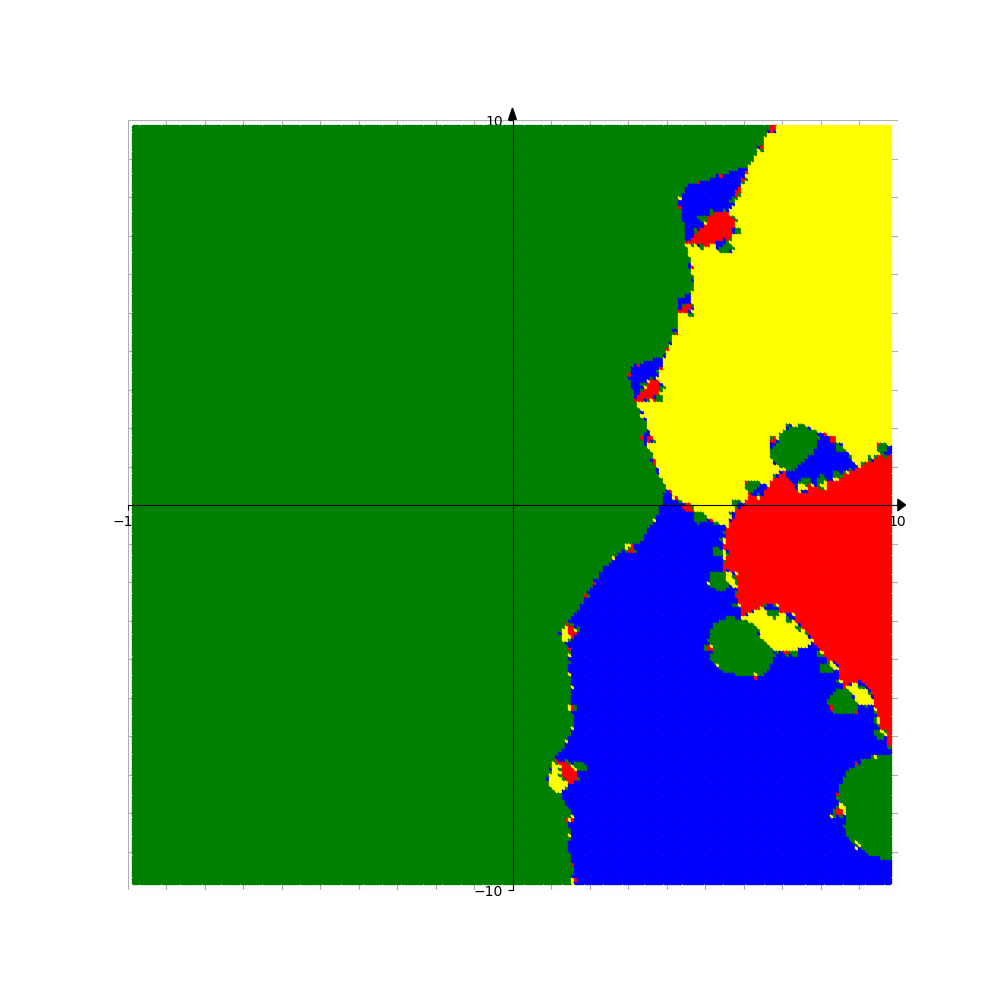}
     \includegraphics[width=3cm]{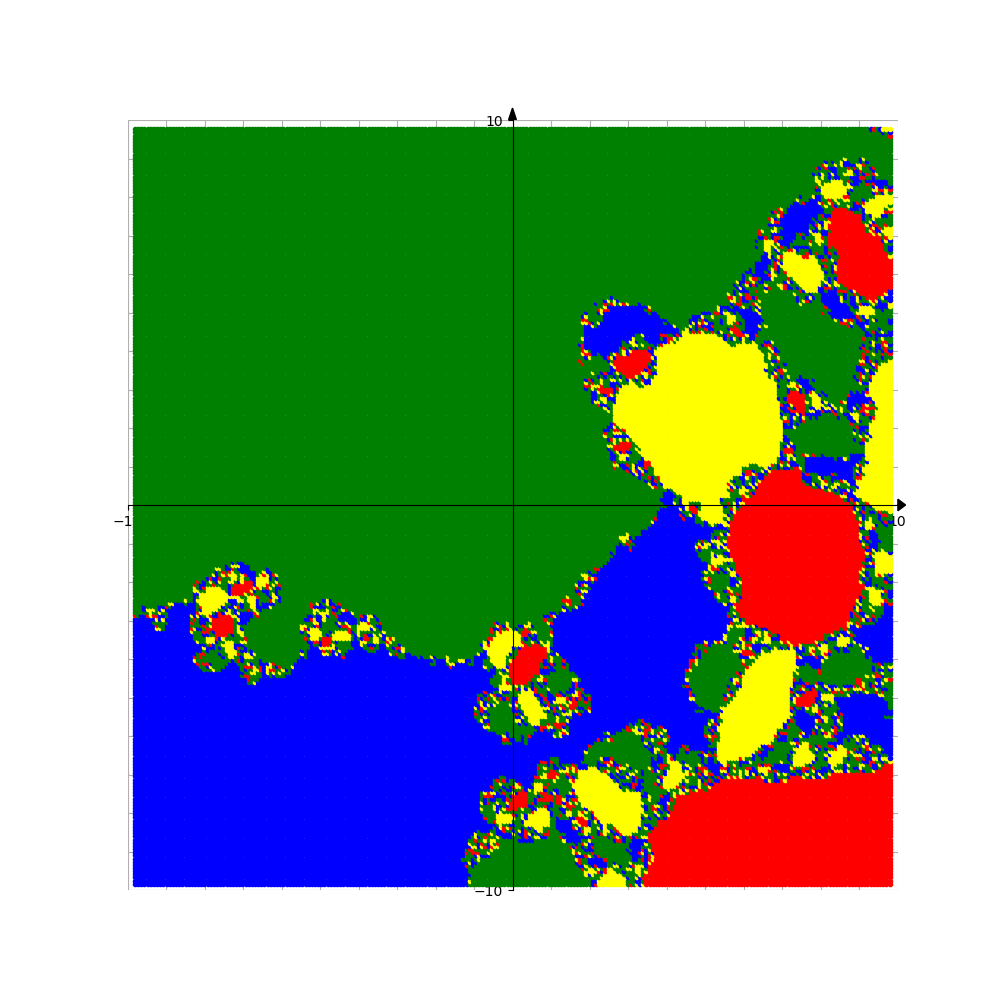}
        \includegraphics[width=5cm]{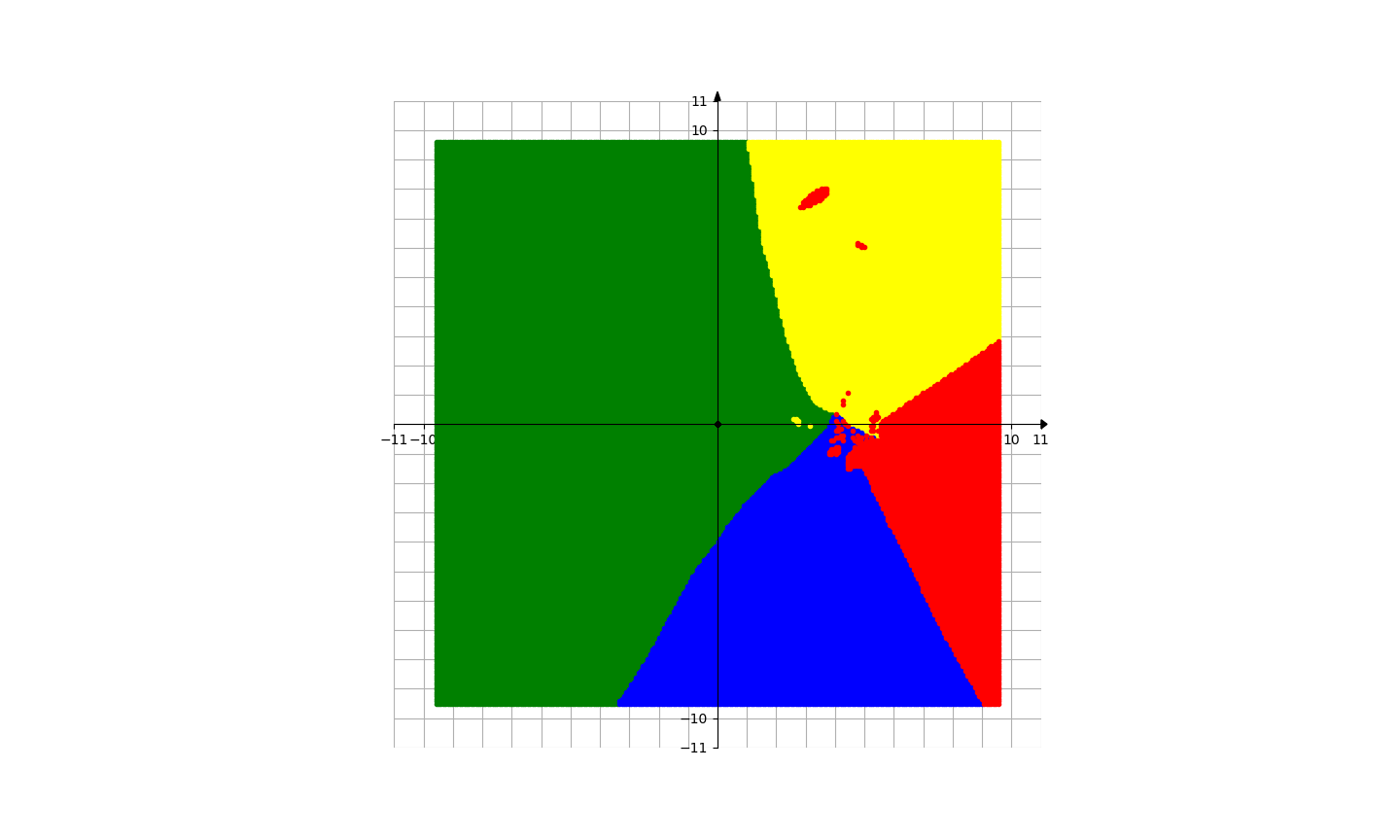}

    \caption{Basins of attraction for finding roots of the function $f_{19}/f_{19}'$ by different methods. Left picture is for Newton's method, central picture is for Random Relaxed Newton's method, right picture is for BNQN. }
    \label{fig:f19Fraction}
\end{figure}

\begin{figure}
    \centering
    \includegraphics[width=3cm]{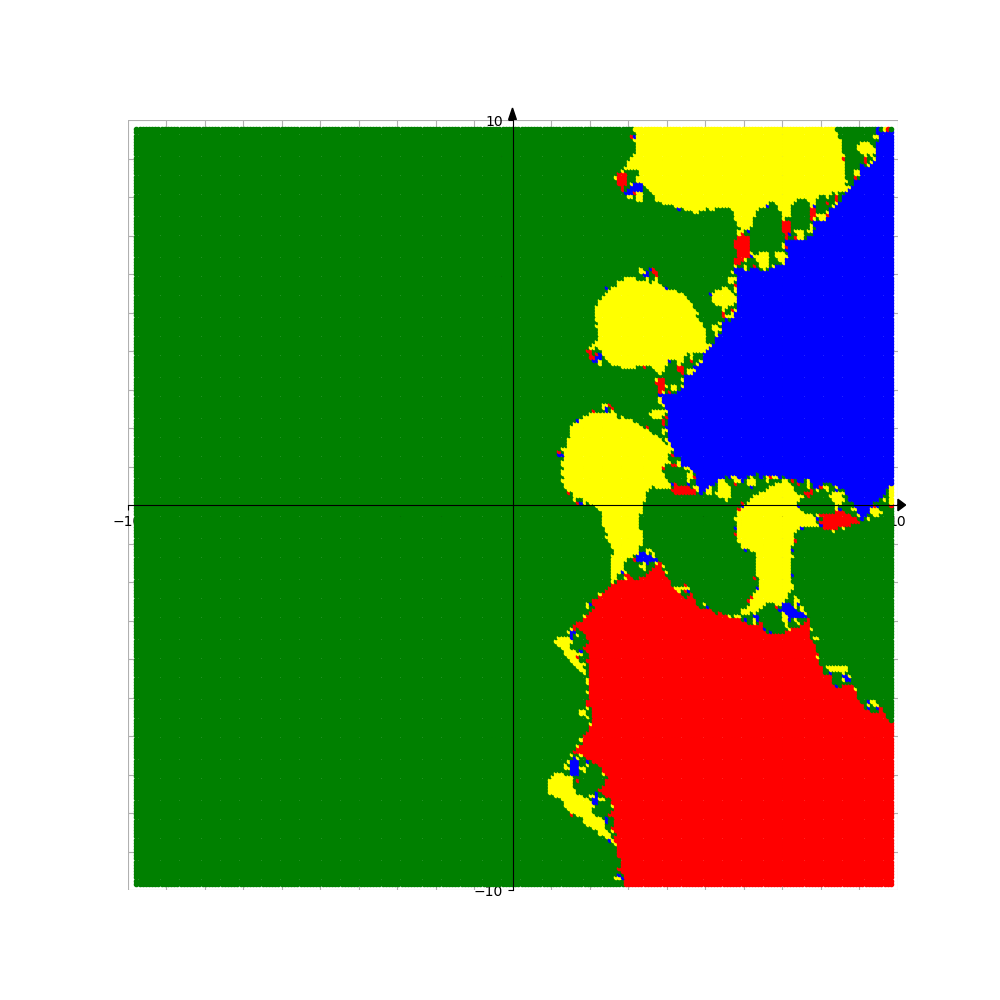}
     \includegraphics[width=3cm]{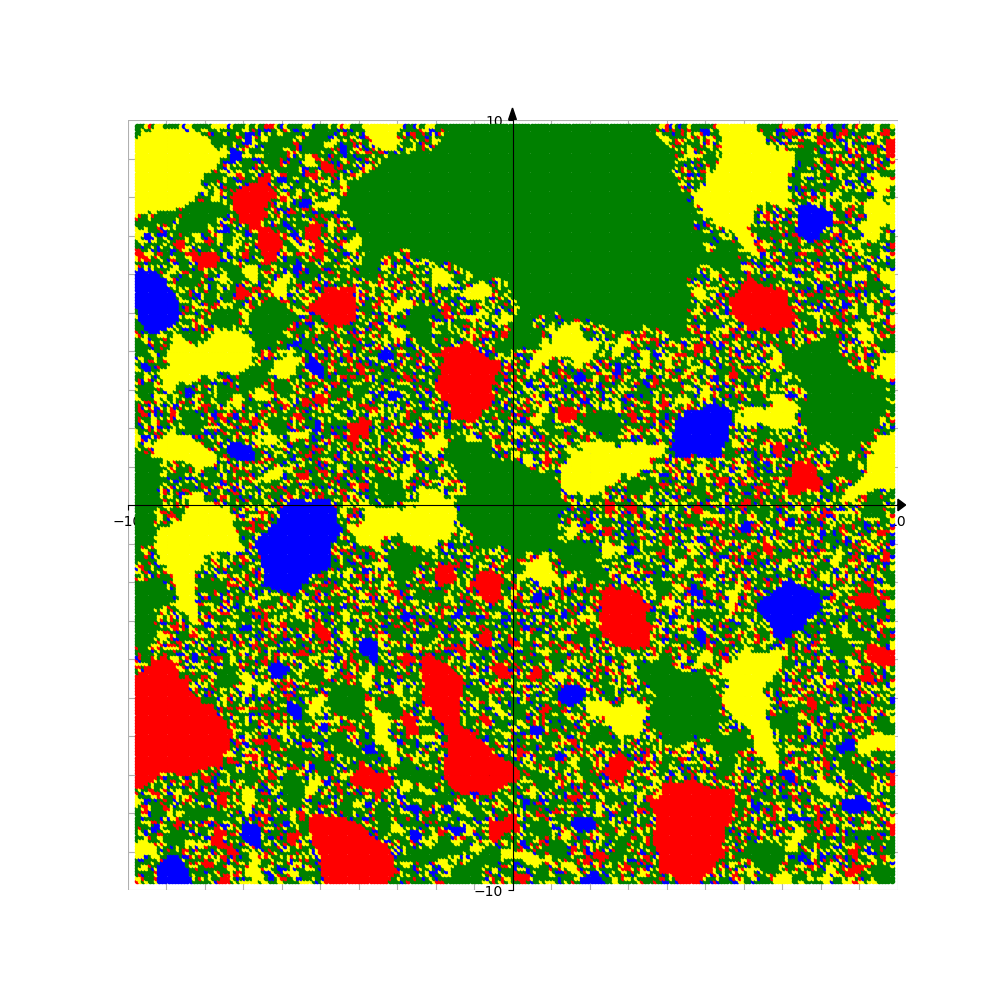}
        \includegraphics[width=5cm]{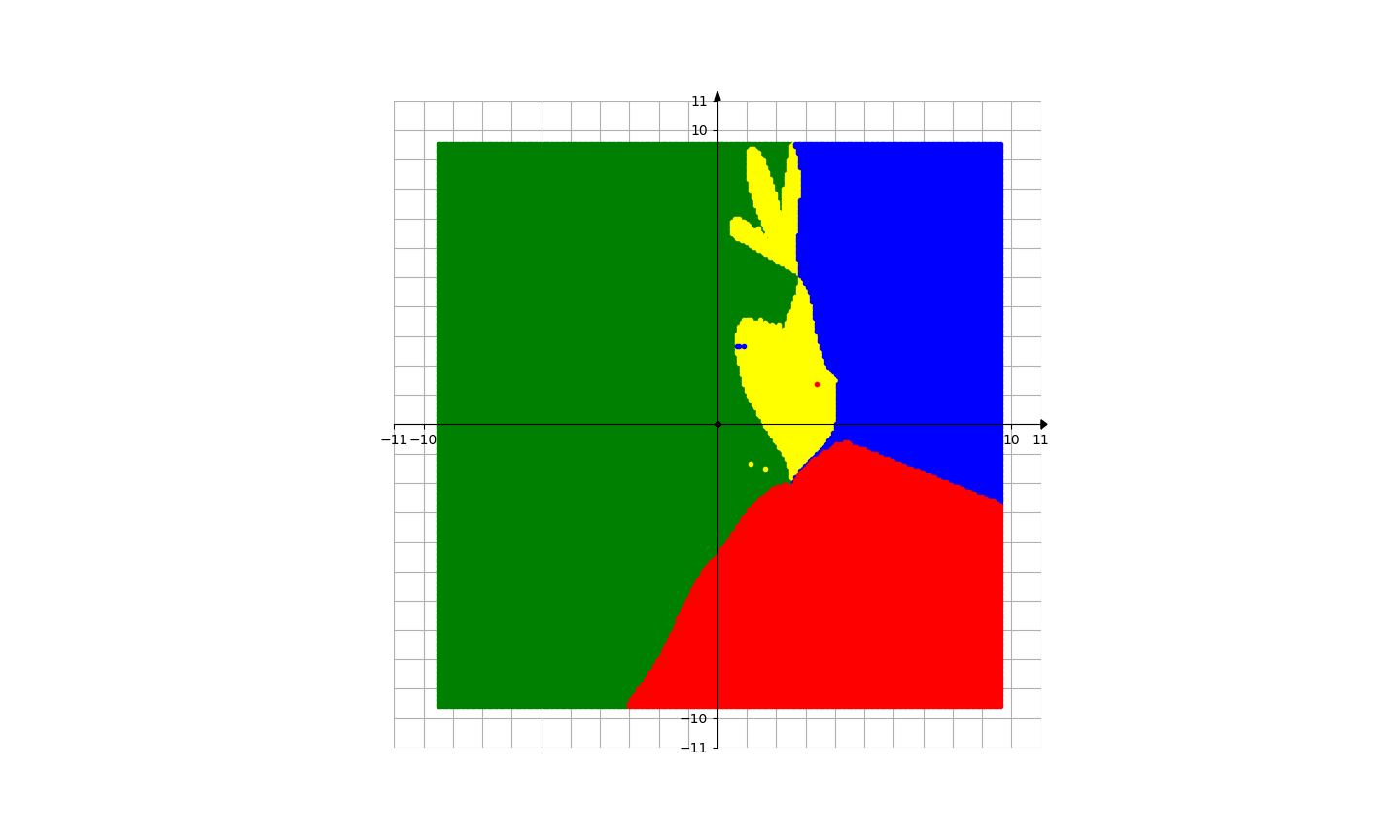}

    \caption{Basins of attraction for finding roots of the function $f_{20}/f_{20}'$ by different methods. Left picture is for Newton's method, central picture is for Random Relaxed Newton's method, right picture is for BNQN. }
    \label{fig:f20Fraction}
\end{figure}

\begin{figure}
    \centering
    \includegraphics[width=3cm]{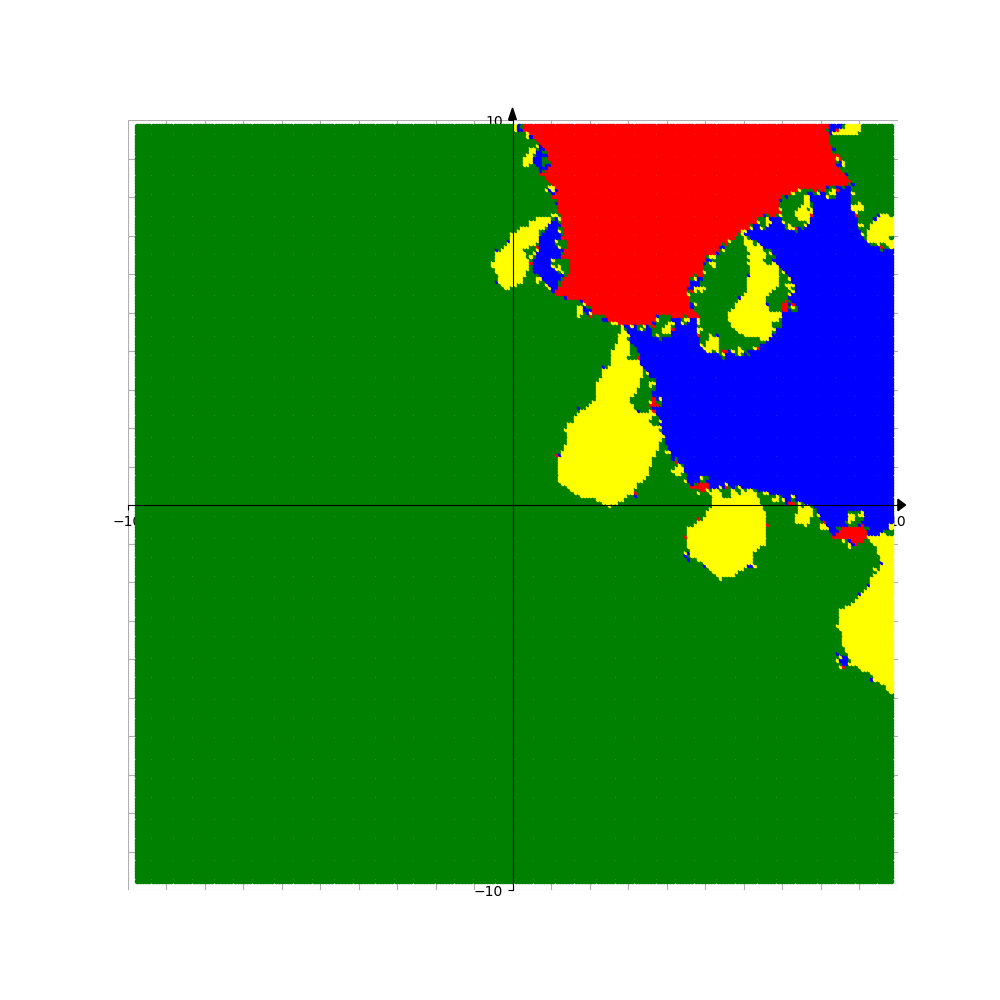}
     \includegraphics[width=3cm]{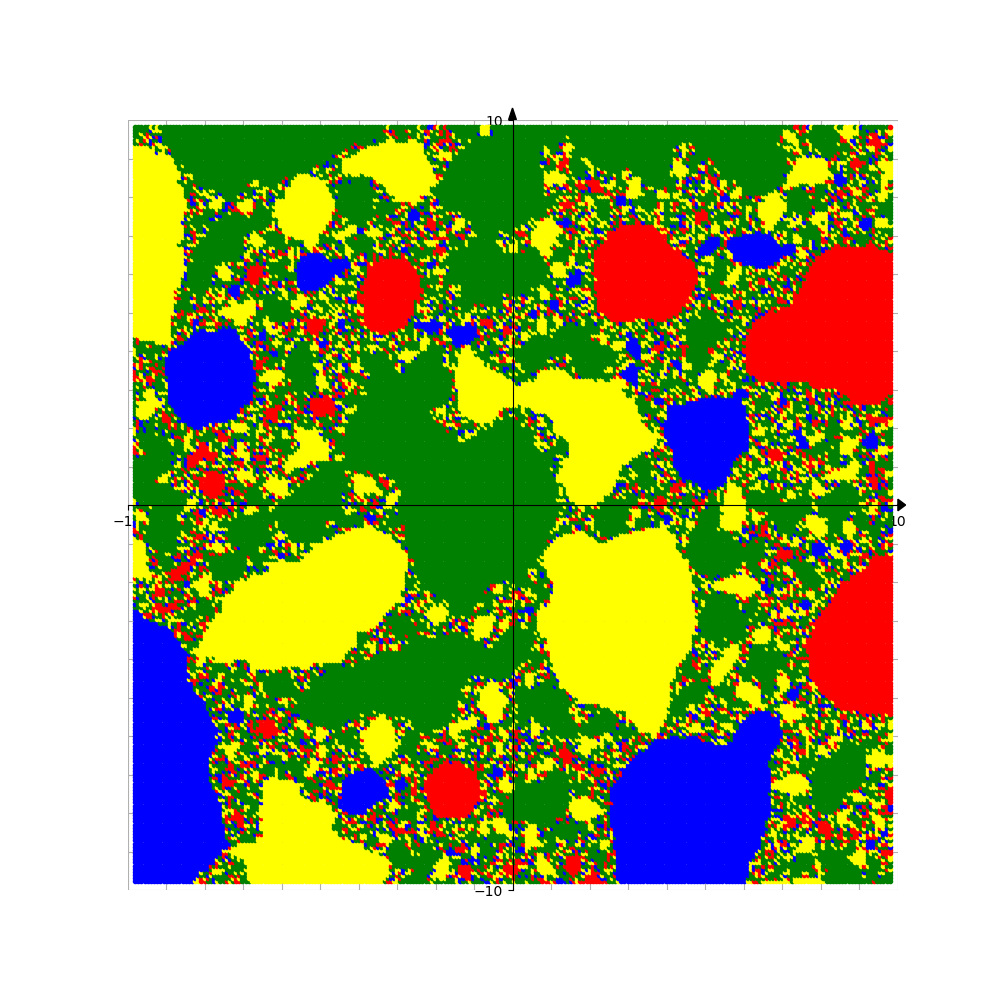}
        \includegraphics[width=5cm]{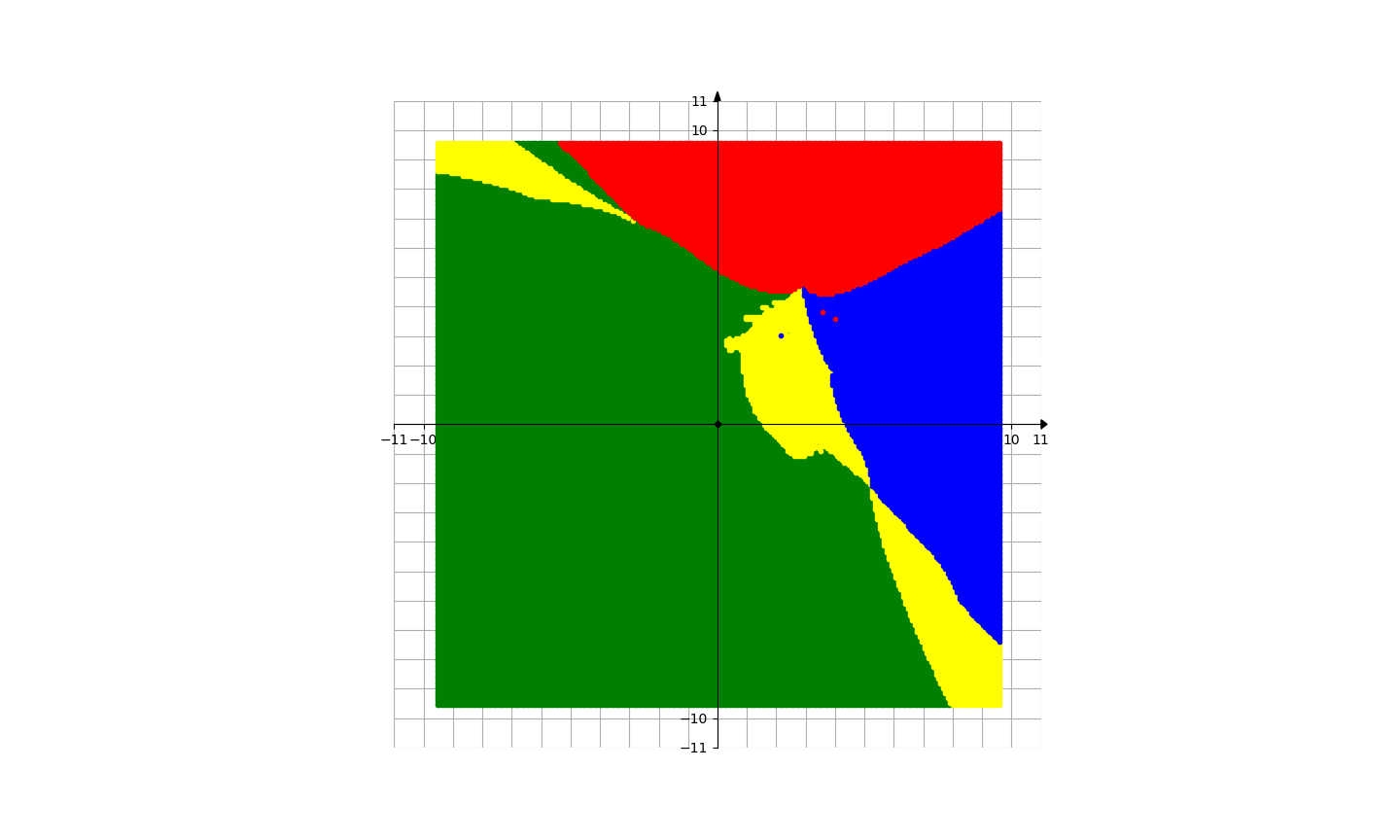}

    \caption{Basins of attraction for finding roots of the function $f_{21}/f_{21}'$ by different methods. Left picture is for Newton's method, central picture is for Random Relaxed Newton's method, right picture is for BNQN. }
    \label{fig:f21Fraction}
\end{figure}

\begin{figure}
    \centering
    \includegraphics[width=3cm]{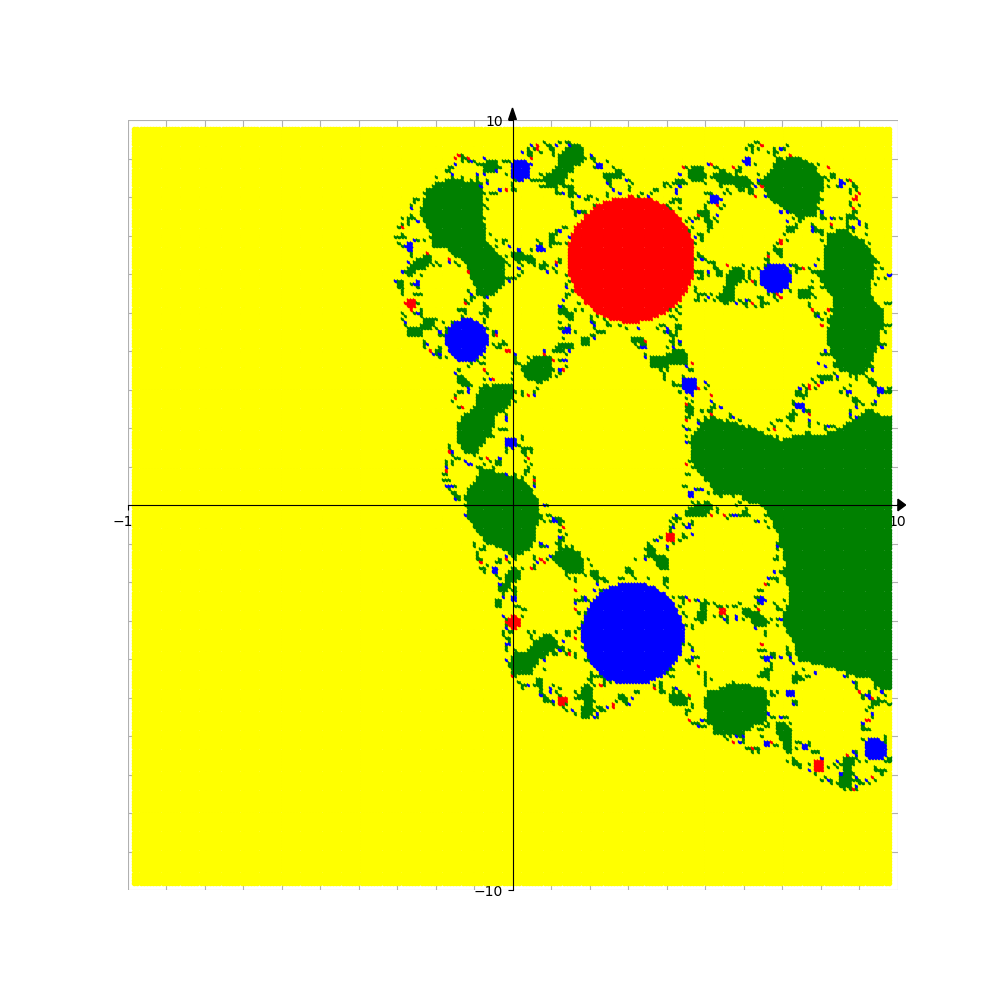}
     \includegraphics[width=3cm]{NewtonRandomf21overf21prime.png}
        \includegraphics[width=5cm]{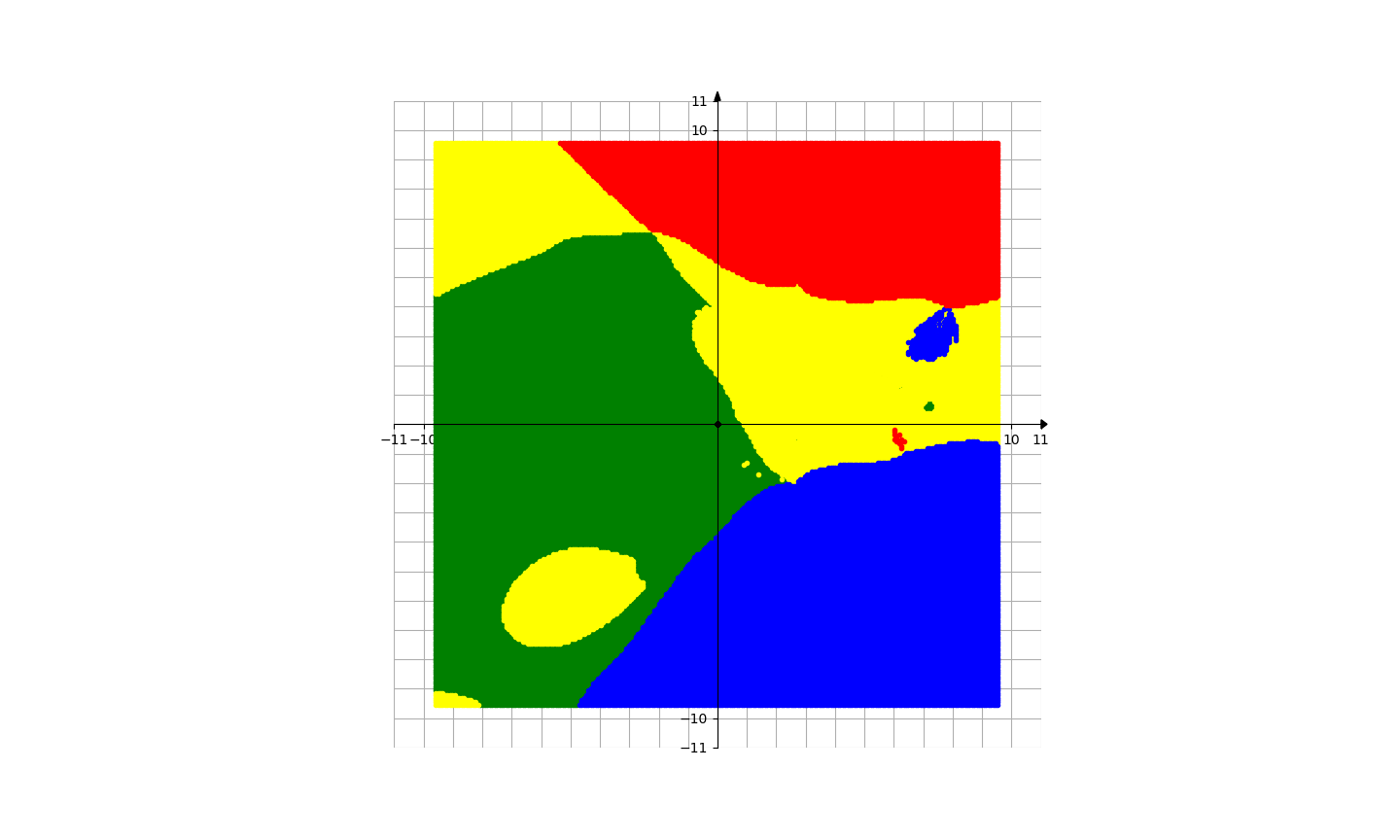}

    \caption{Basins of attraction for finding roots of the function $f_{22}/f_{22}'$ by different methods. Left picture is for Newton's method, central picture is for Random Relaxed Newton's method, right picture is for BNQN. }
    \label{fig:f22Fraction}
\end{figure}

\begin{figure}
    \centering

        \includegraphics[width=4cm]{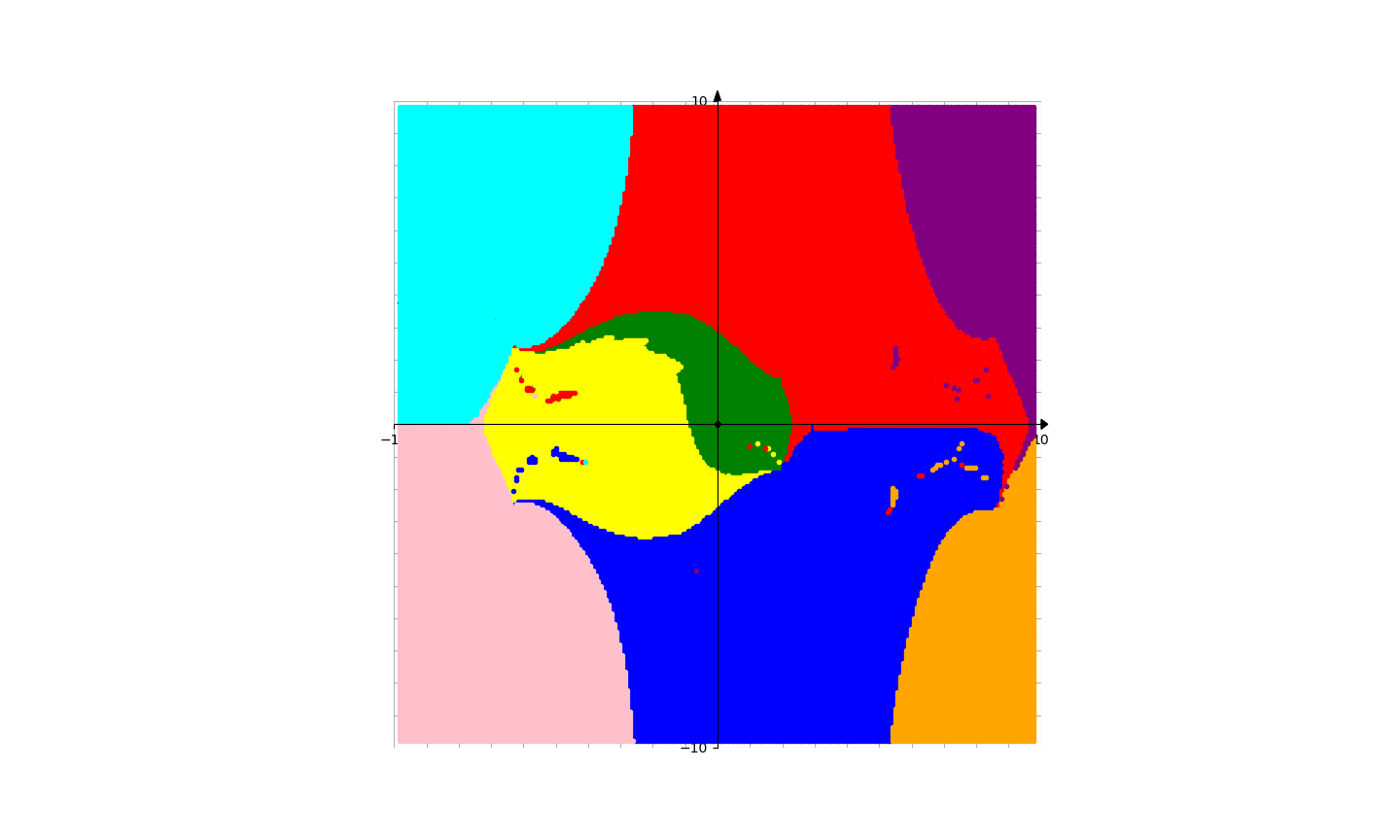}
        \includegraphics[width=4cm]{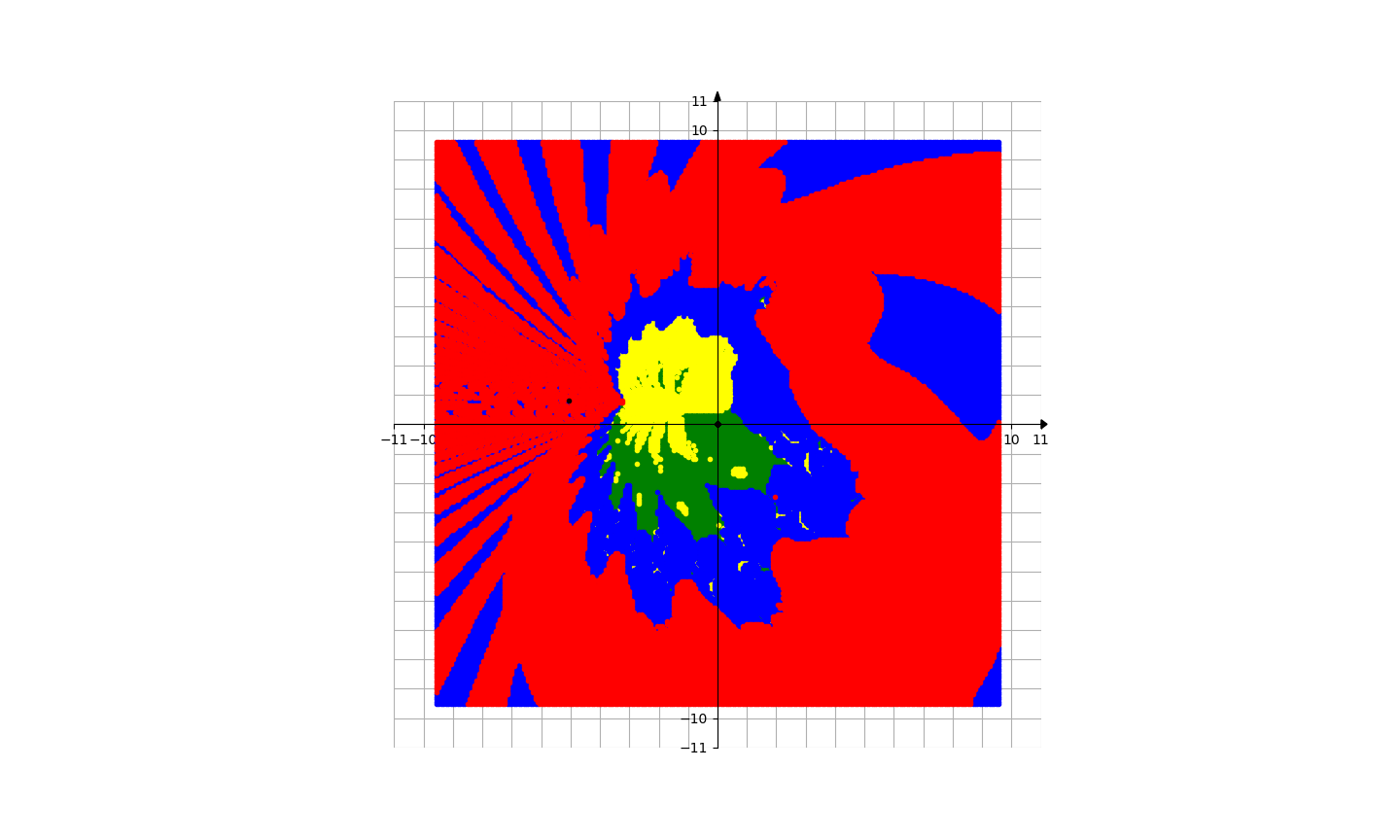}
        \includegraphics[width=4cm]{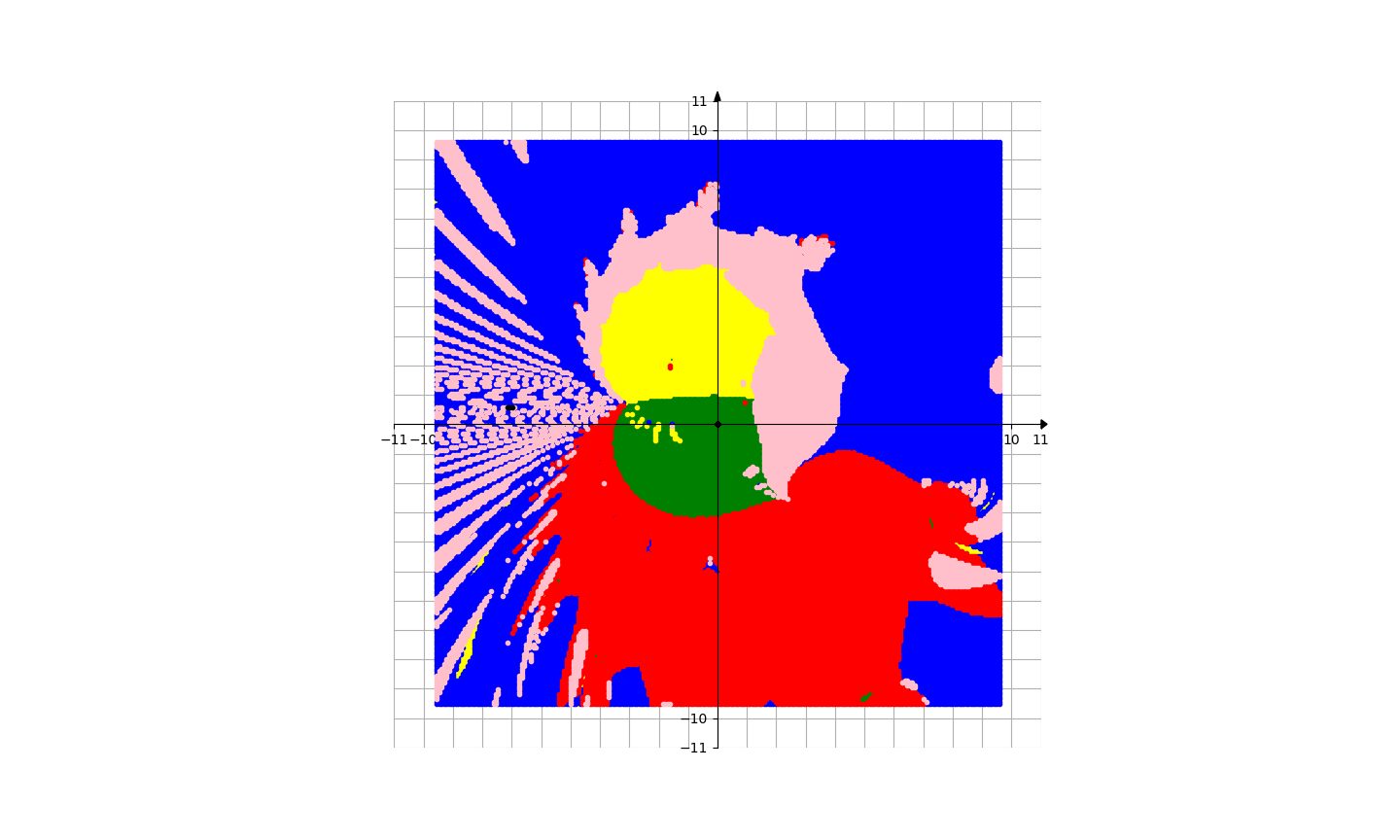}

    \caption{Basins of attraction for finding roots of the function $f_{23}/f_{23}'$, $f_{24}/f_{24}'$ and $f_{25}/f_{25}'$ by BNQN. For this case, Newton's method and  Random Relaxed Newton's method  either encounter errors or take very long time to finish. }
    \label{fig:f23f24f25Fraction}
\end{figure}

\begin{figure}
    \centering

        \includegraphics[width=4cm]{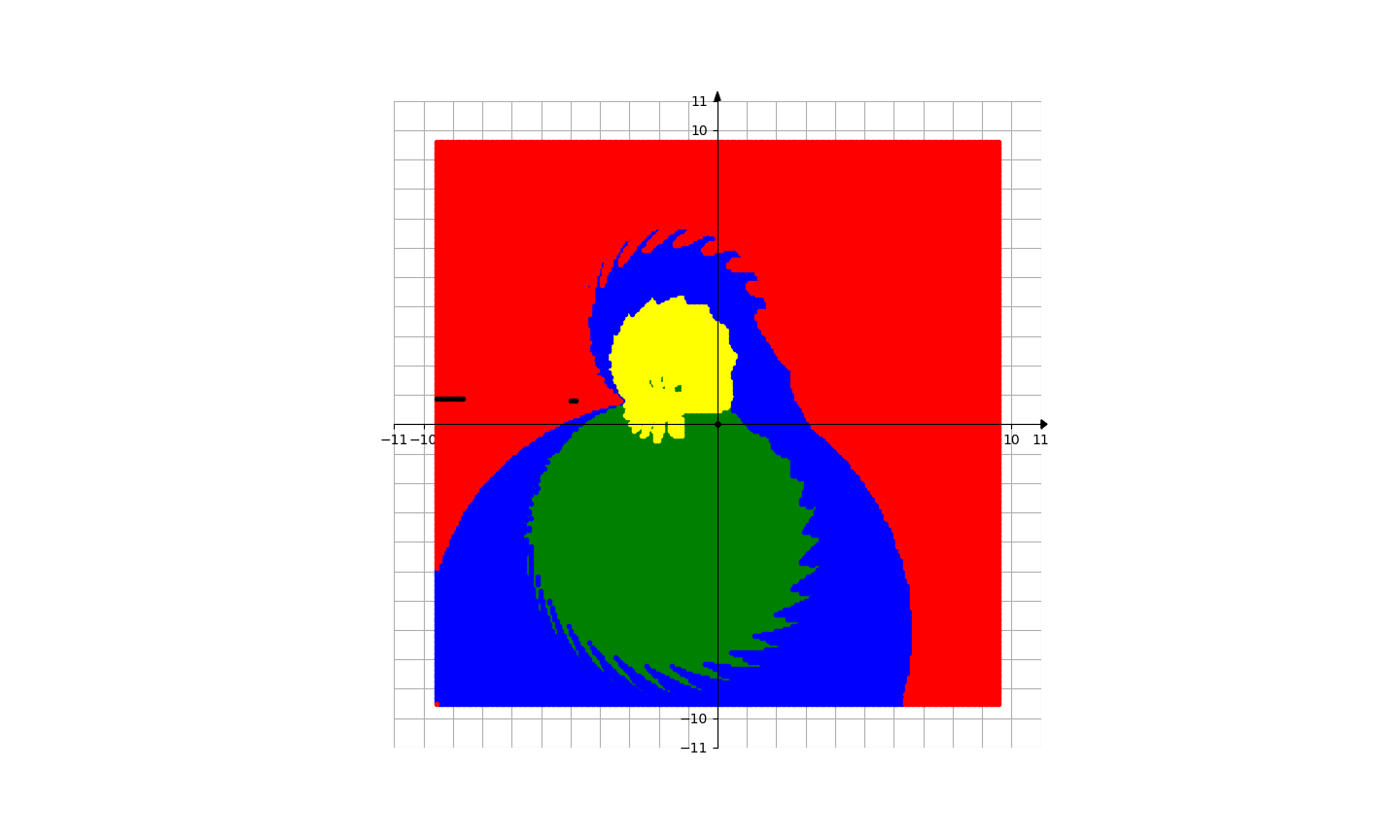}
        \includegraphics[width=4cm]{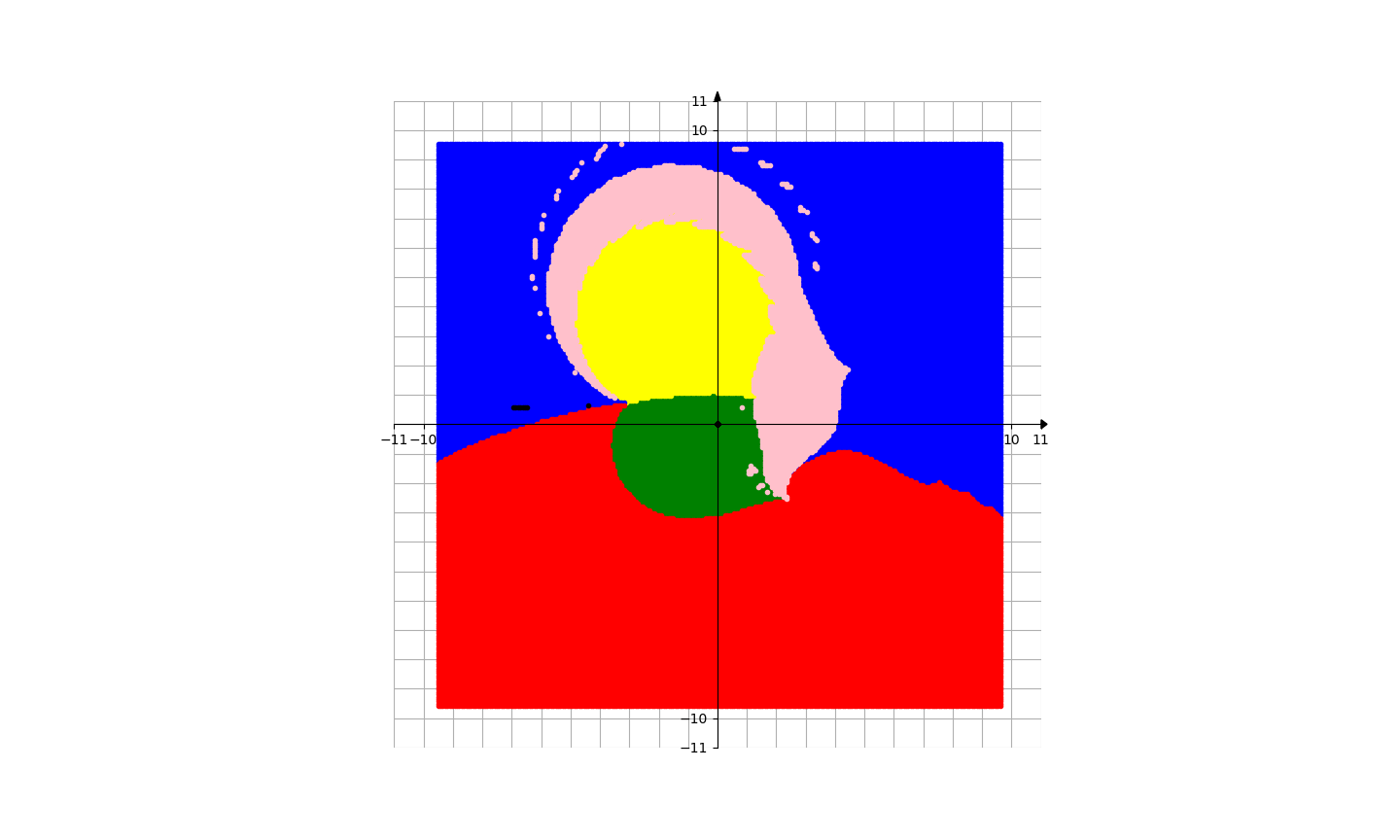}

    \caption{Basins of attraction for finding roots of the function $f_{24}/f_{24}'$ and $f_{25}/f_{25}'$ by BNQN v2, which are more smooth than that by BNQN in Figure \ref{fig:f23f24f25Fraction}. }
    \label{fig:f23f24f25FractionV2}
\end{figure}

\subsubsection{Stochastic root finding}

\begin{figure}
    \centering

        \includegraphics[width=4cm]{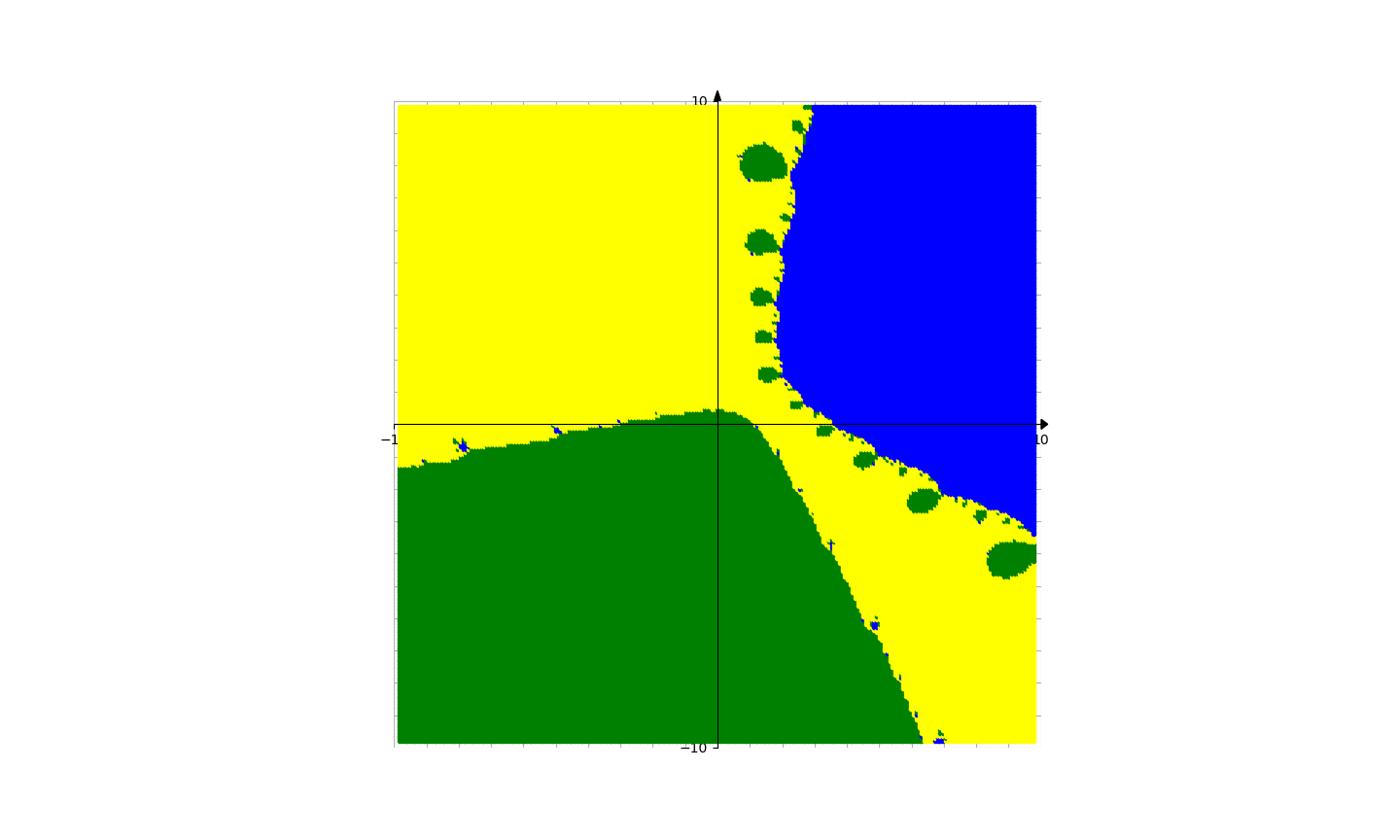}
        \includegraphics[width=4cm]{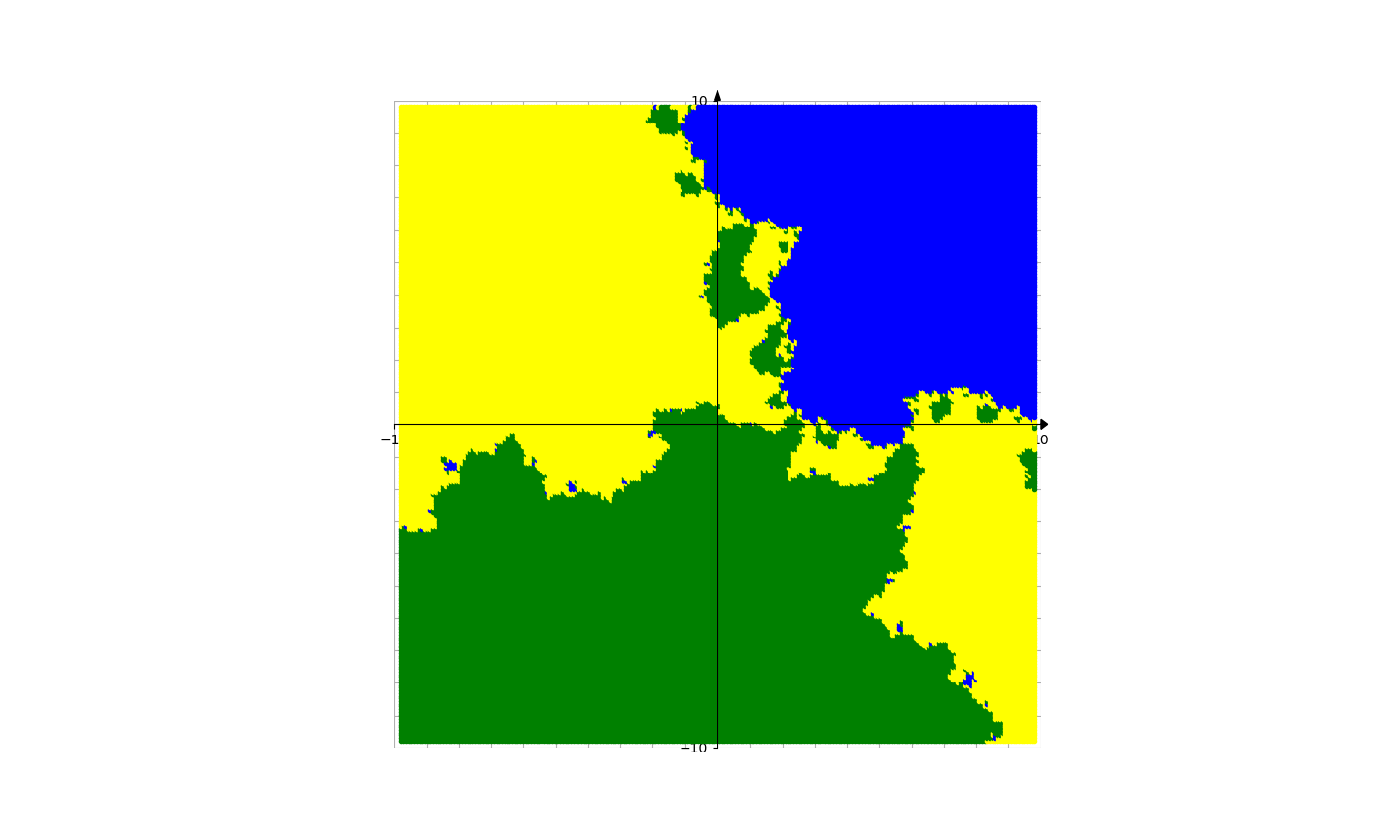}
   
    \bigskip
    
     \includegraphics[width=4cm]{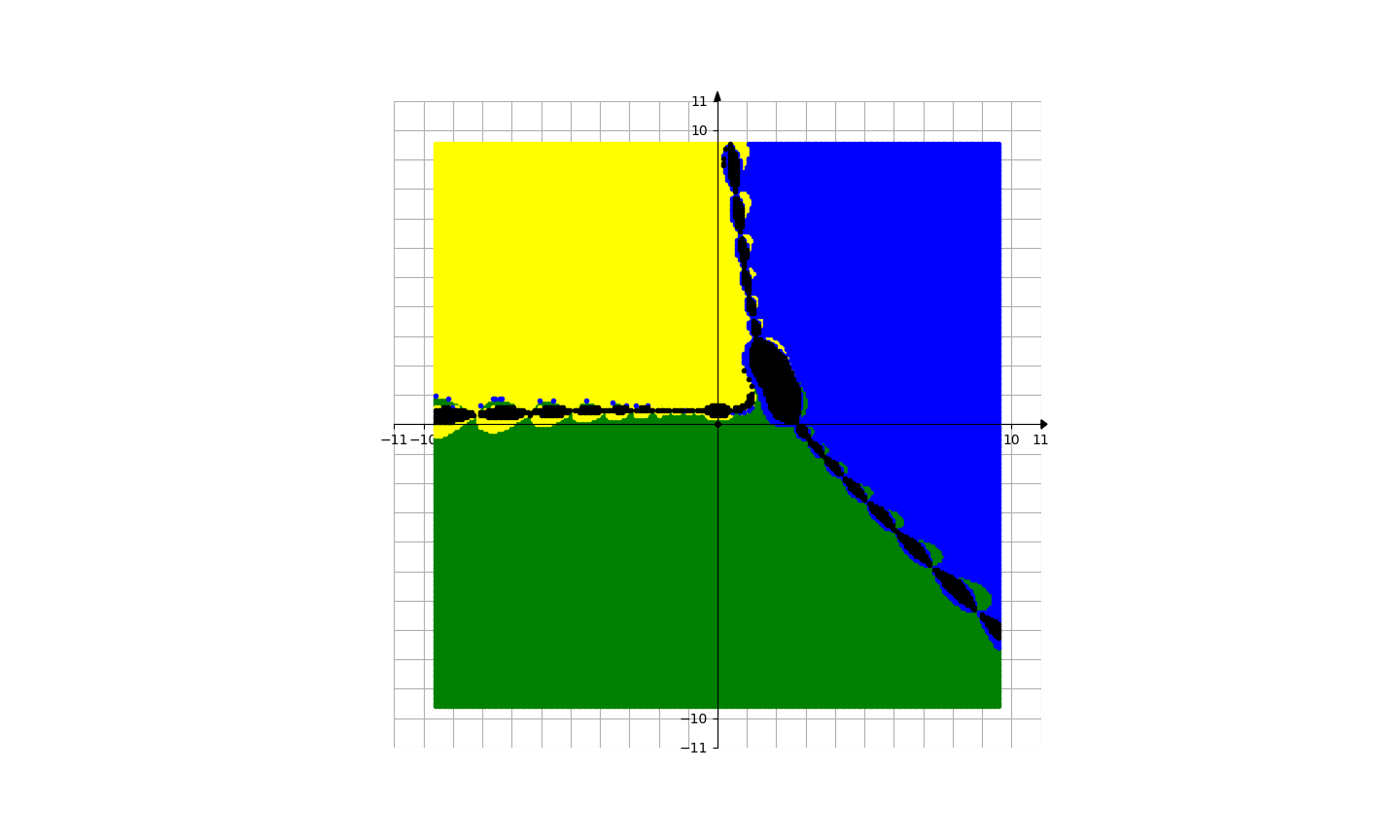}
      \includegraphics[width=4cm]{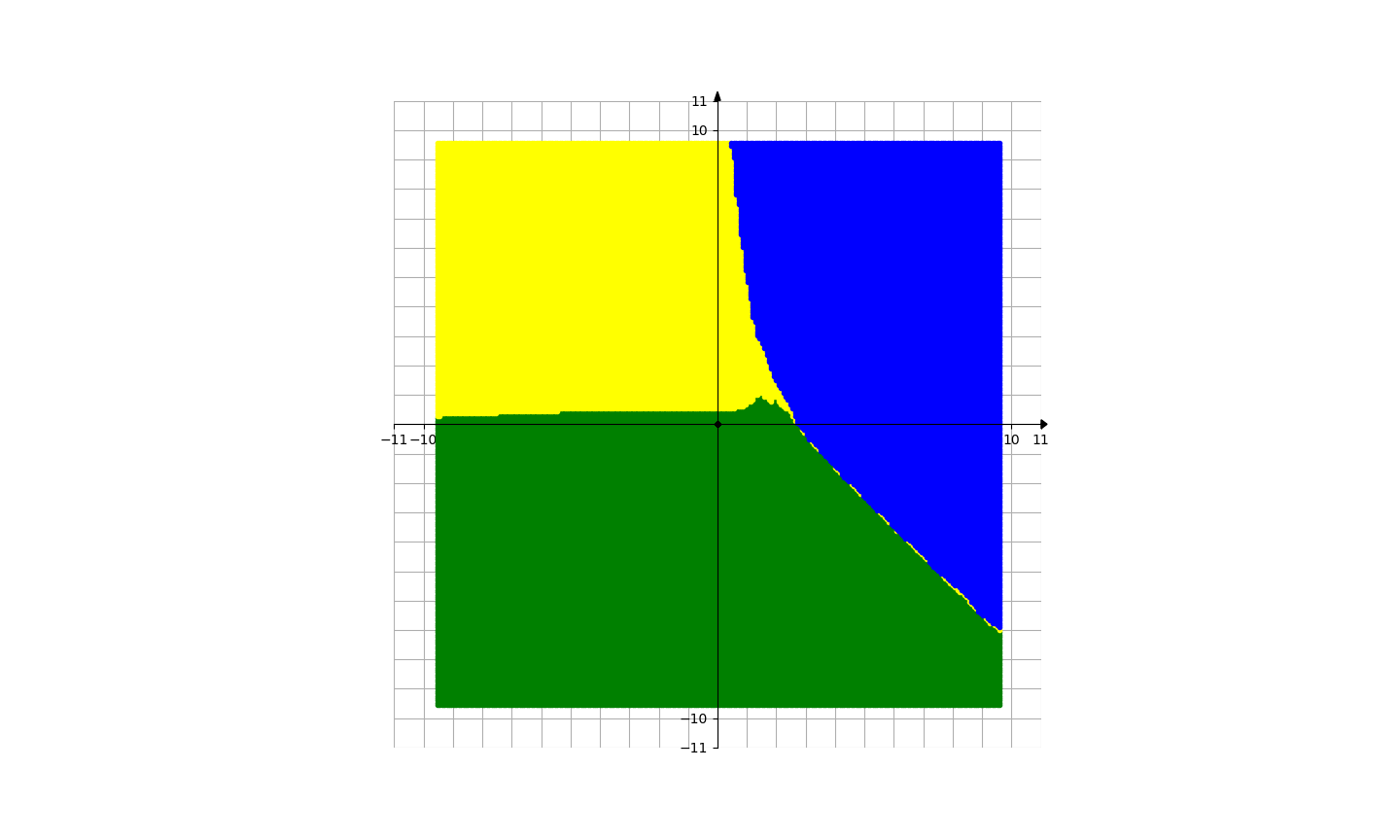}
     \includegraphics[width=4cm]{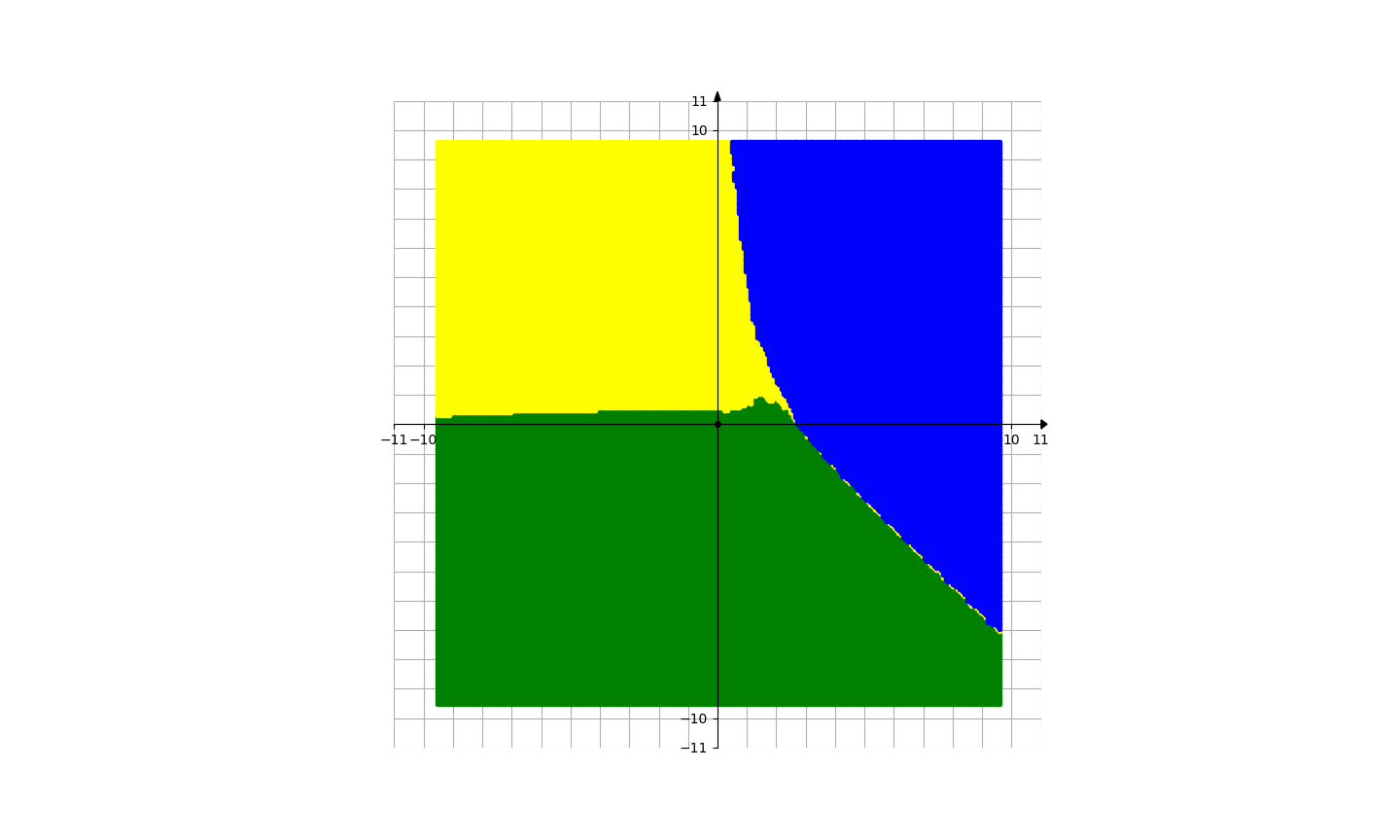}

    \caption{Basins of attraction for finding roots of the stochastic function $f_{1}+\epsilon \xi (z^3+2z-5)$ by different methods. Pictures are referenced to from top to bottom, from left to right. Row 1: left picture is for Newton's method, right picture is for Random Relaxed Newton's method. Row 2: left picture is for Newton's method vOptimization, central picture is for BNQN, right picture is for BNQN v2.}
    \label{fig:f1Stochastic}
\end{figure}

\begin{figure}
    \centering

        \includegraphics[width=4cm]{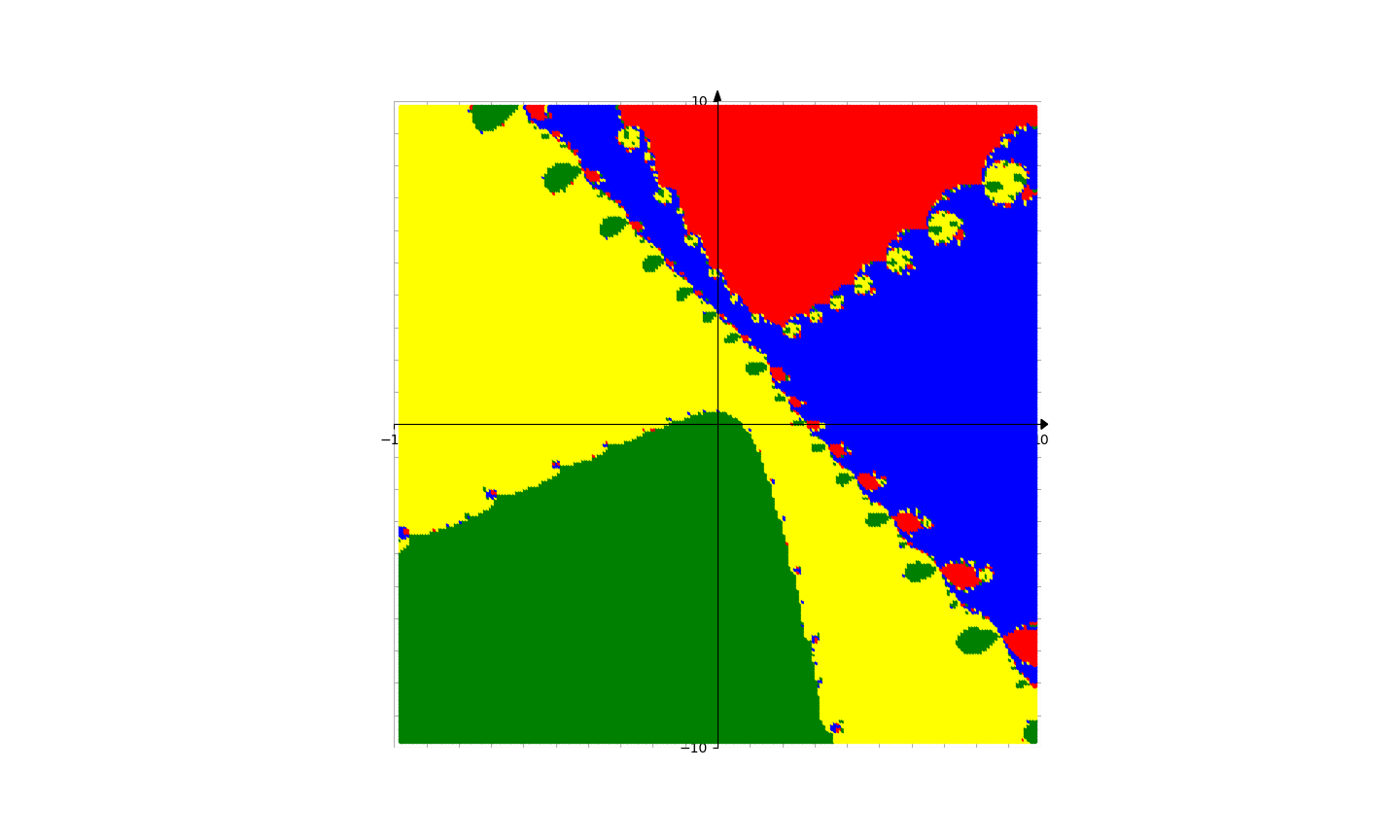}
        \includegraphics[width=4cm]{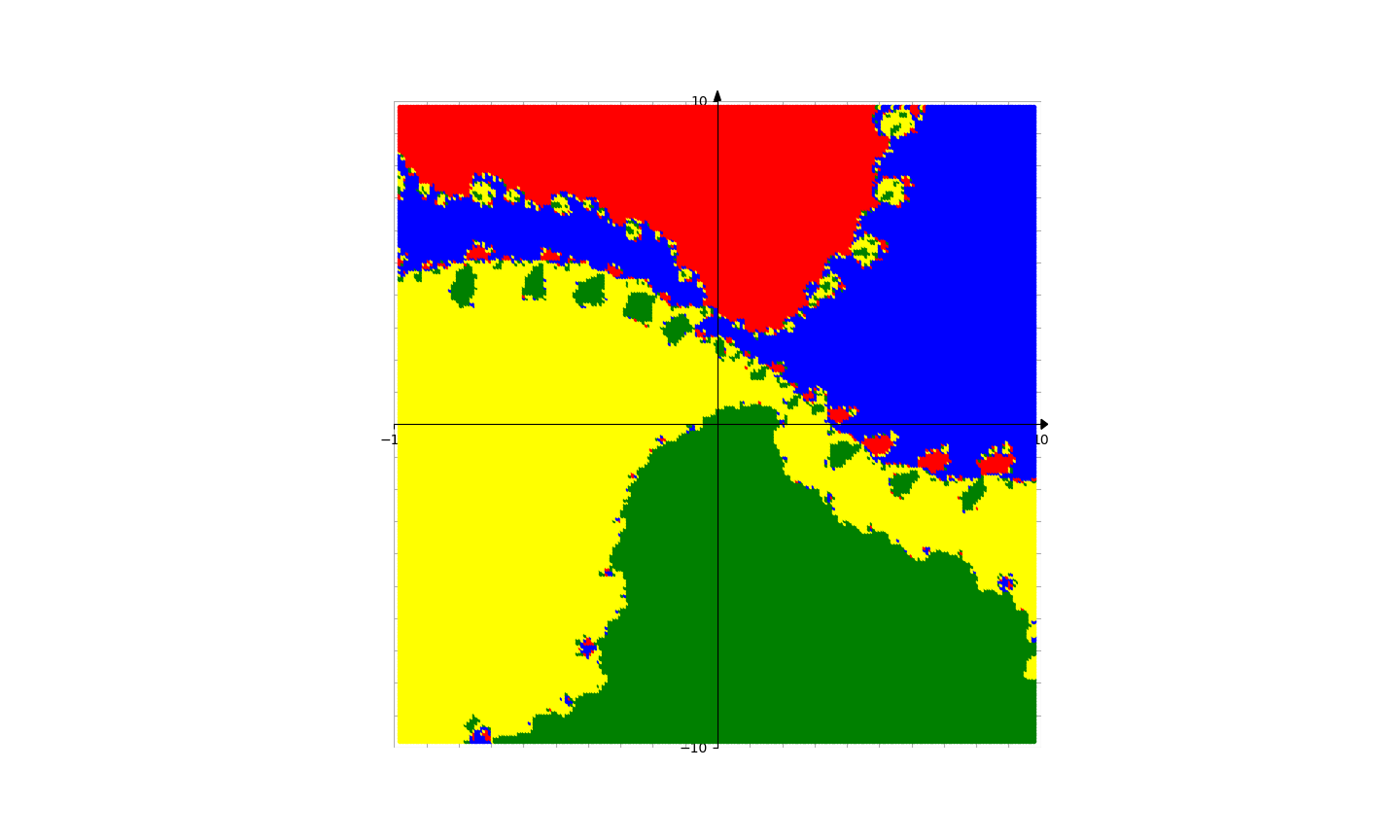}
   
    \bigskip
    
     \includegraphics[width=4cm]{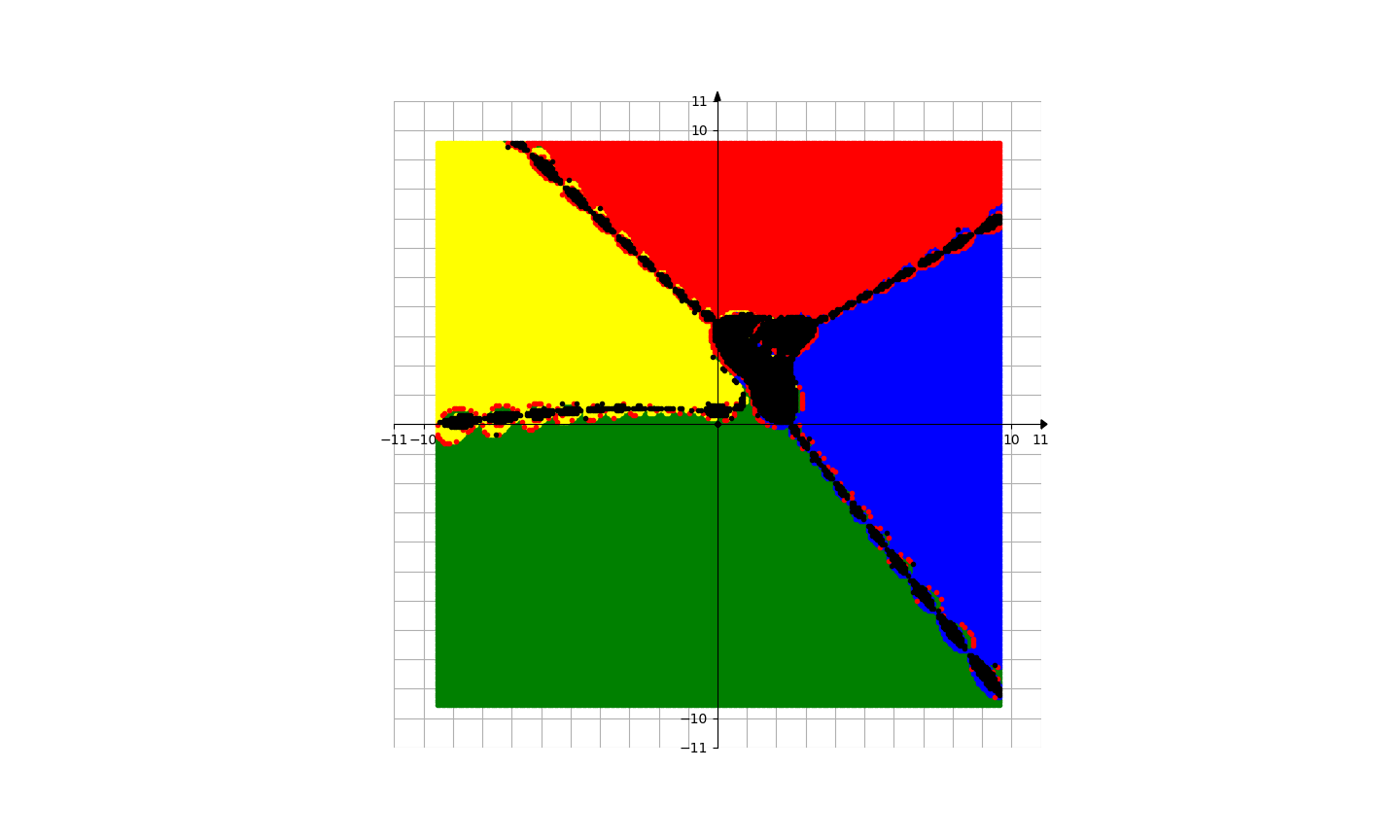}
      \includegraphics[width=4cm]{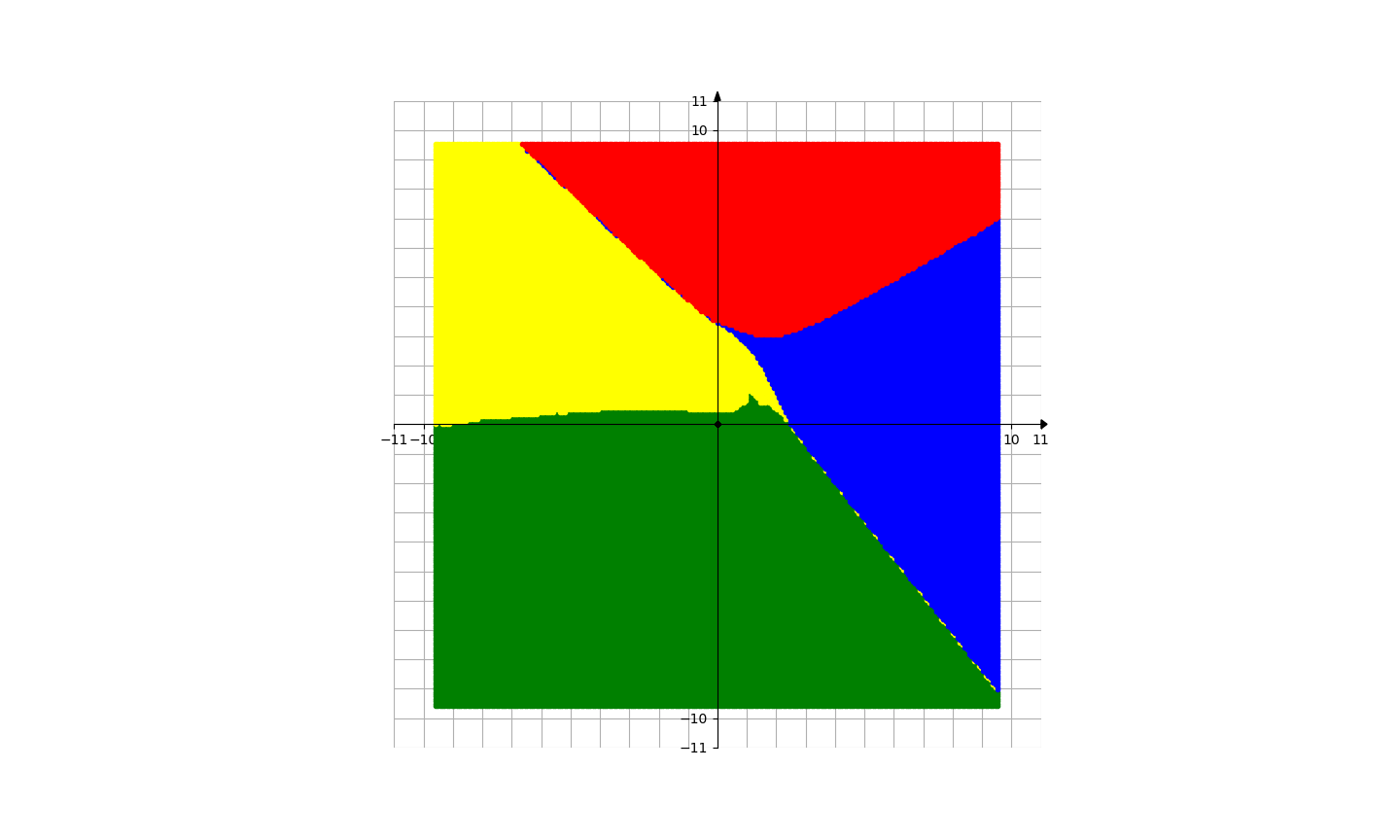}
     \includegraphics[width=4cm]{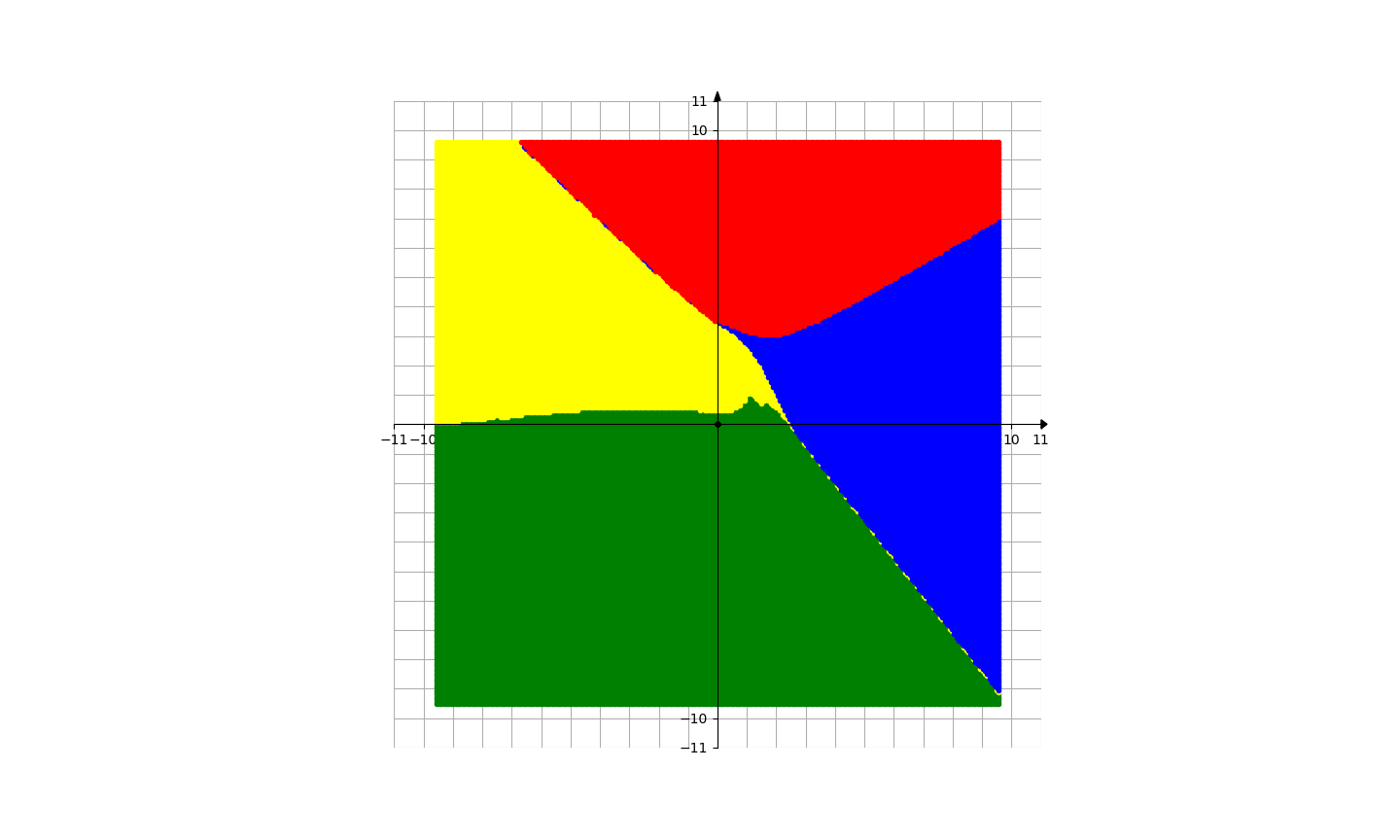}

    \caption{Basins of attraction for finding roots of the stochastic function $f_{5}+\epsilon \xi (z^3+2z-5)$ by different methods. Pictures are referenced to from top to bottom, from left to right. Row 1: left picture is for Newton's method, right picture is for Random Relaxed Newton's method. Row 2: left picture is for Newton's method vOptimization, central picture is for BNQN, right picture is for BNQN v2.}
    \label{fig:f5Stochastic}
\end{figure}

\begin{figure}
    \centering

        \includegraphics[width=4cm]{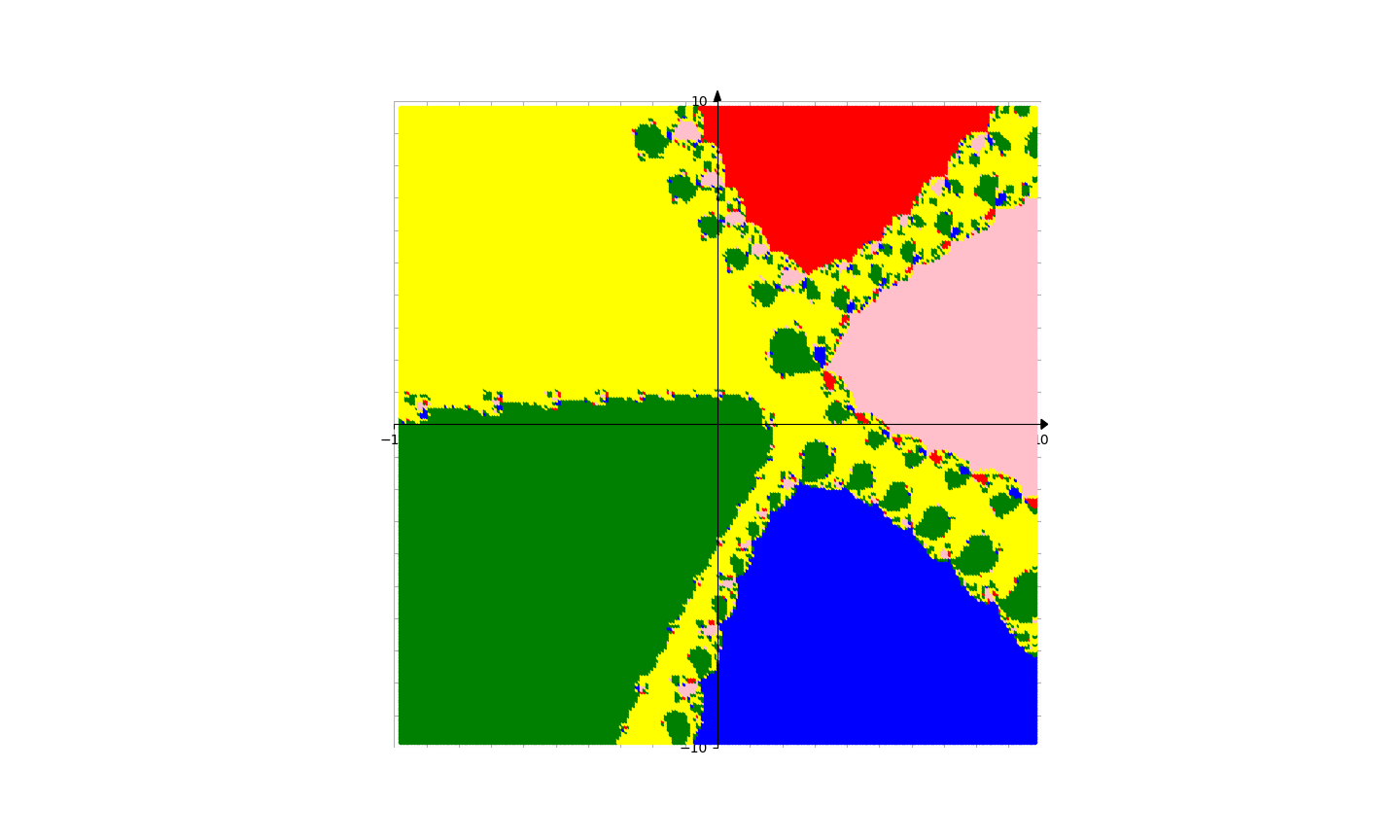}
        \includegraphics[width=4cm]{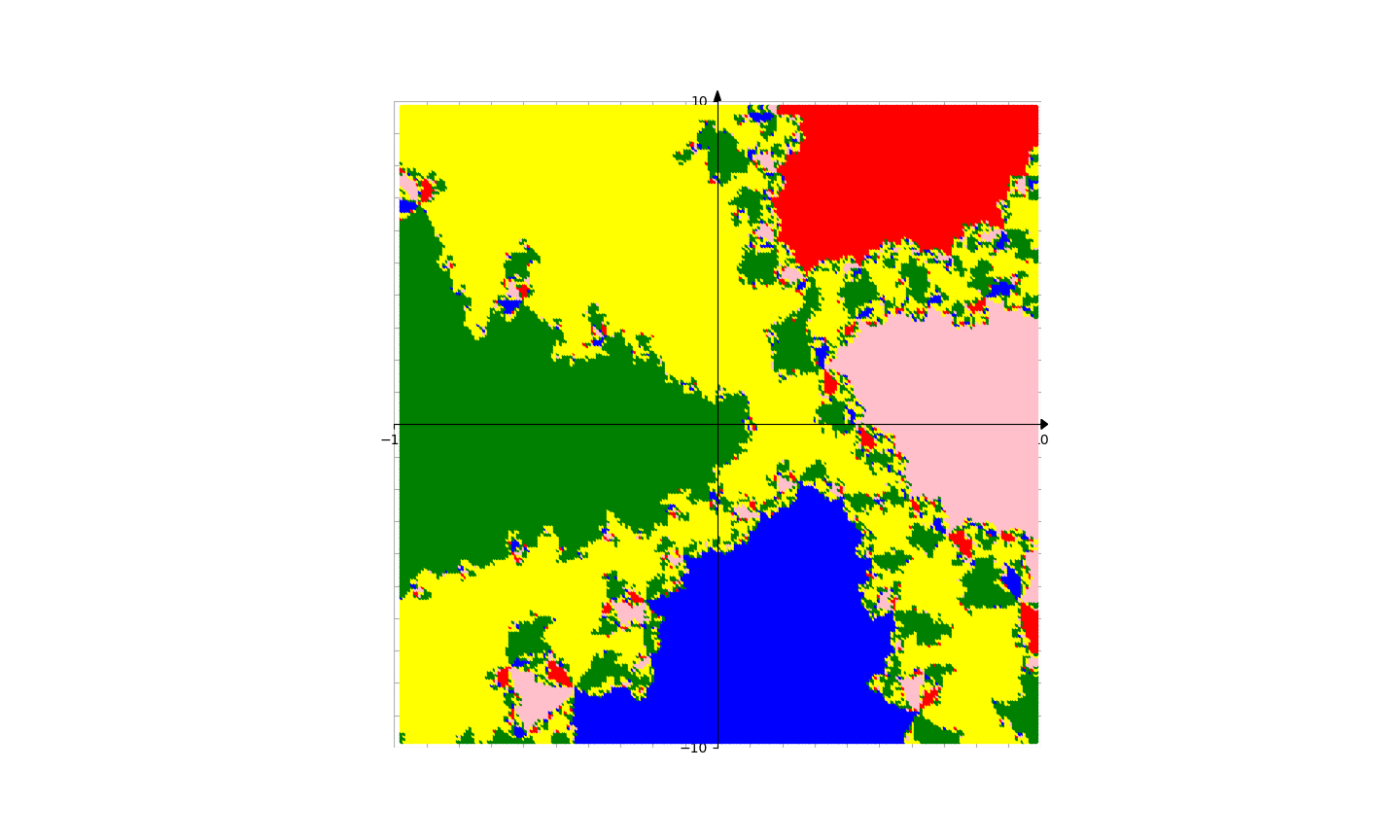}
   
    \bigskip
    
     \includegraphics[width=4cm]{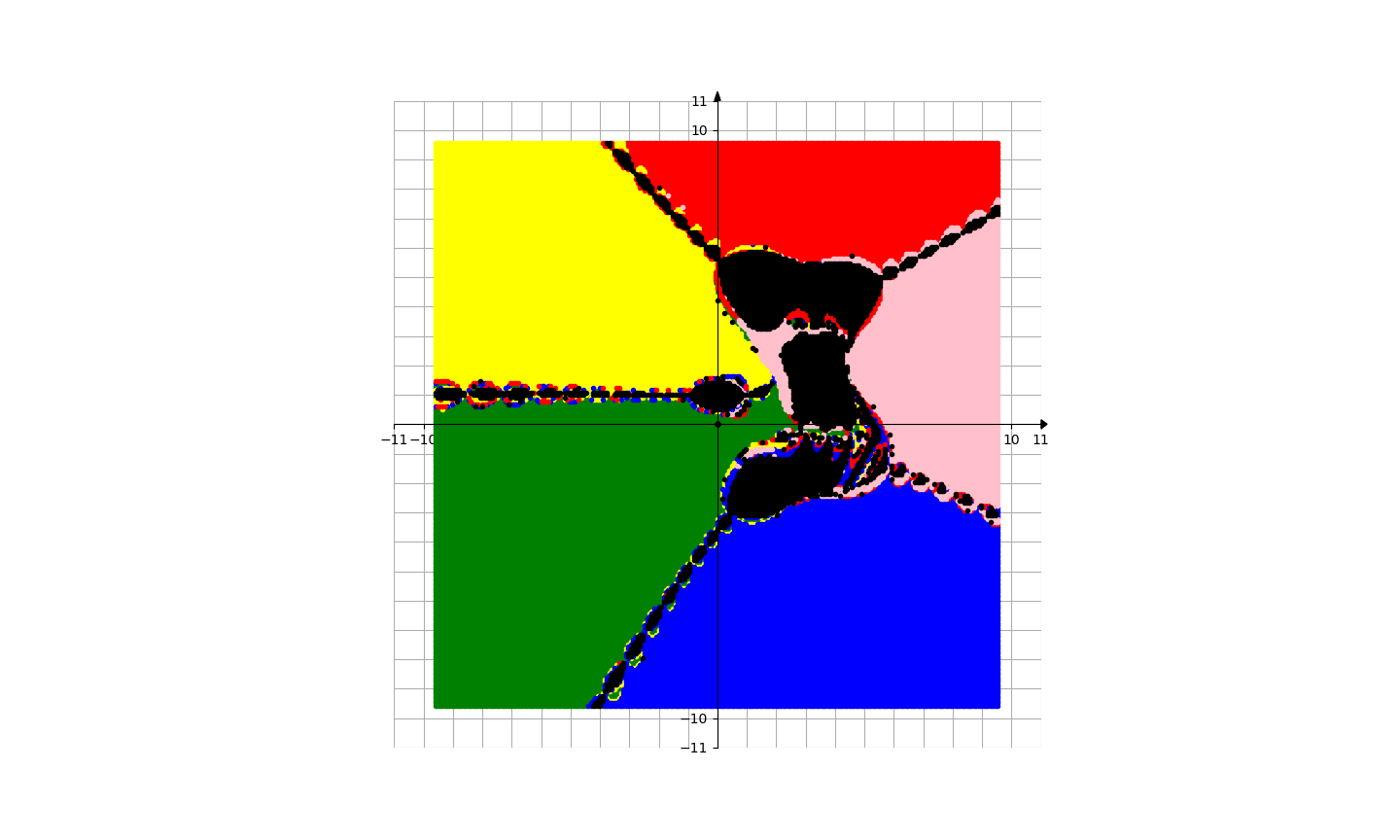}
      \includegraphics[width=4cm]{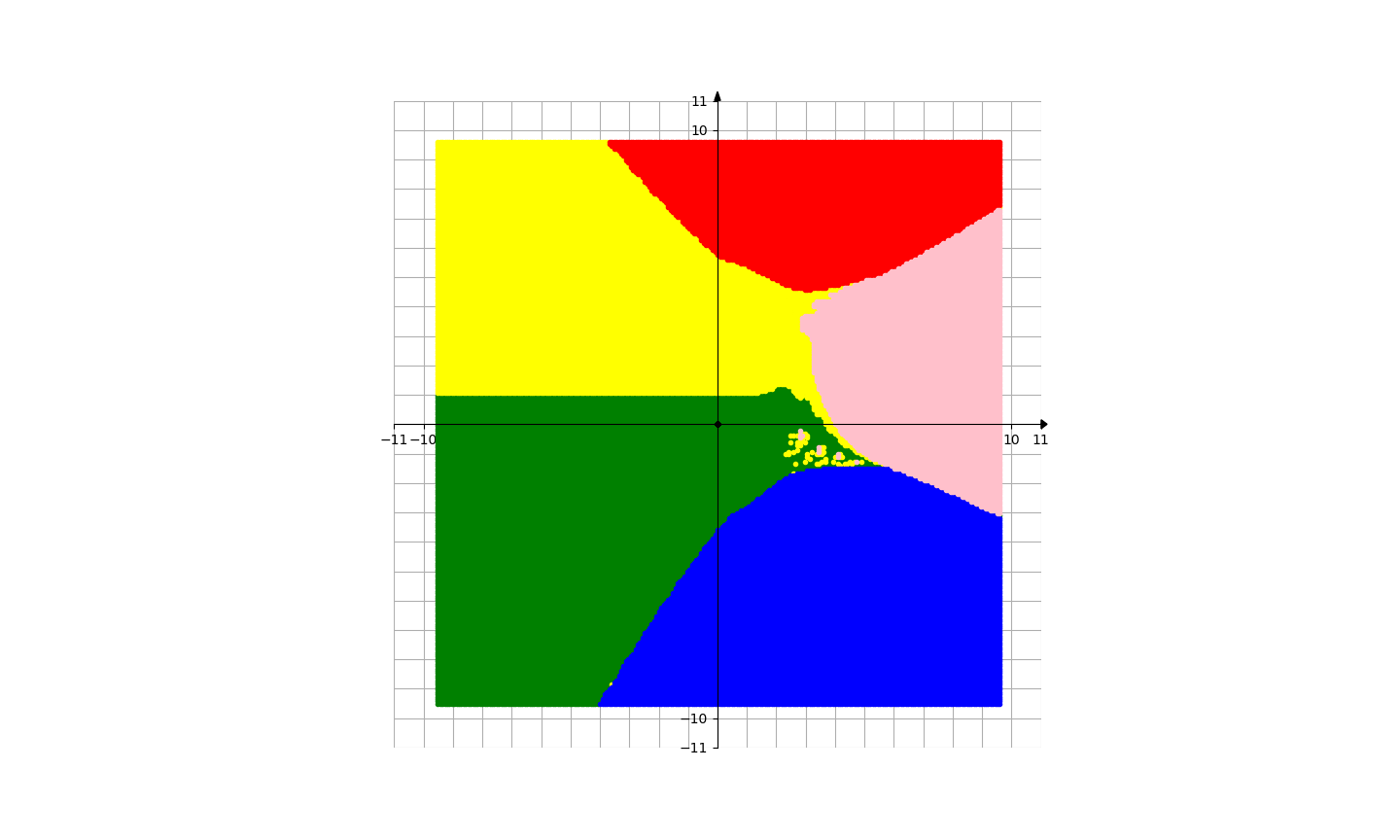}
     \includegraphics[width=4cm]{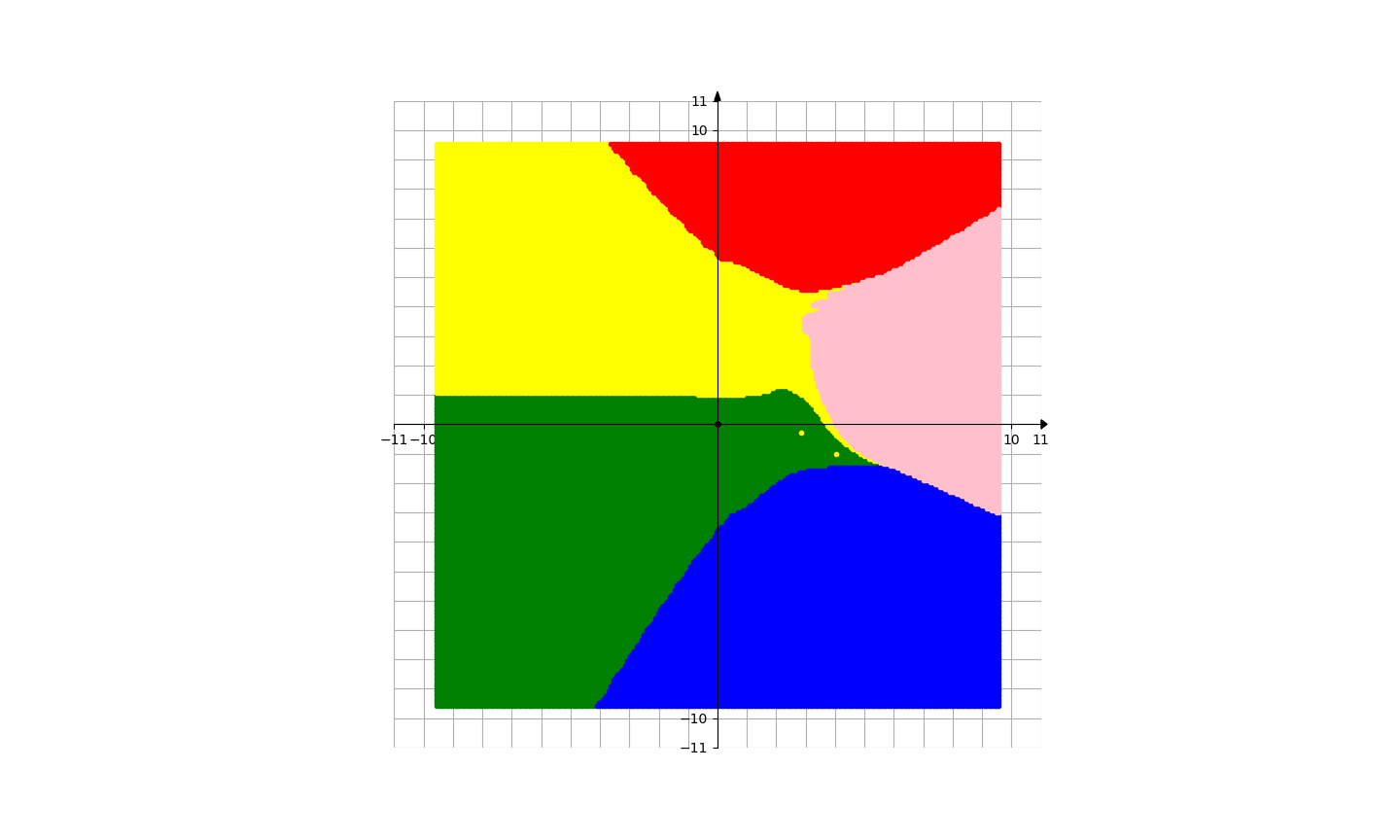}

    \caption{Basins of attraction for finding roots of the stochastic function $f_{14}+\epsilon \xi (z^3+2z-5)$ by different methods. Pictures are referenced to from top to bottom, from left to right. Row 1: left picture is for Newton's method, right picture is for Random Relaxed Newton's method. Row 2: left picture is for Newton's method vOptimization, central picture is for BNQN, right picture is for BNQN v2.}
    \label{fig:f14Stochastic}
\end{figure}

\begin{figure}
    \centering

        \includegraphics[width=4cm]{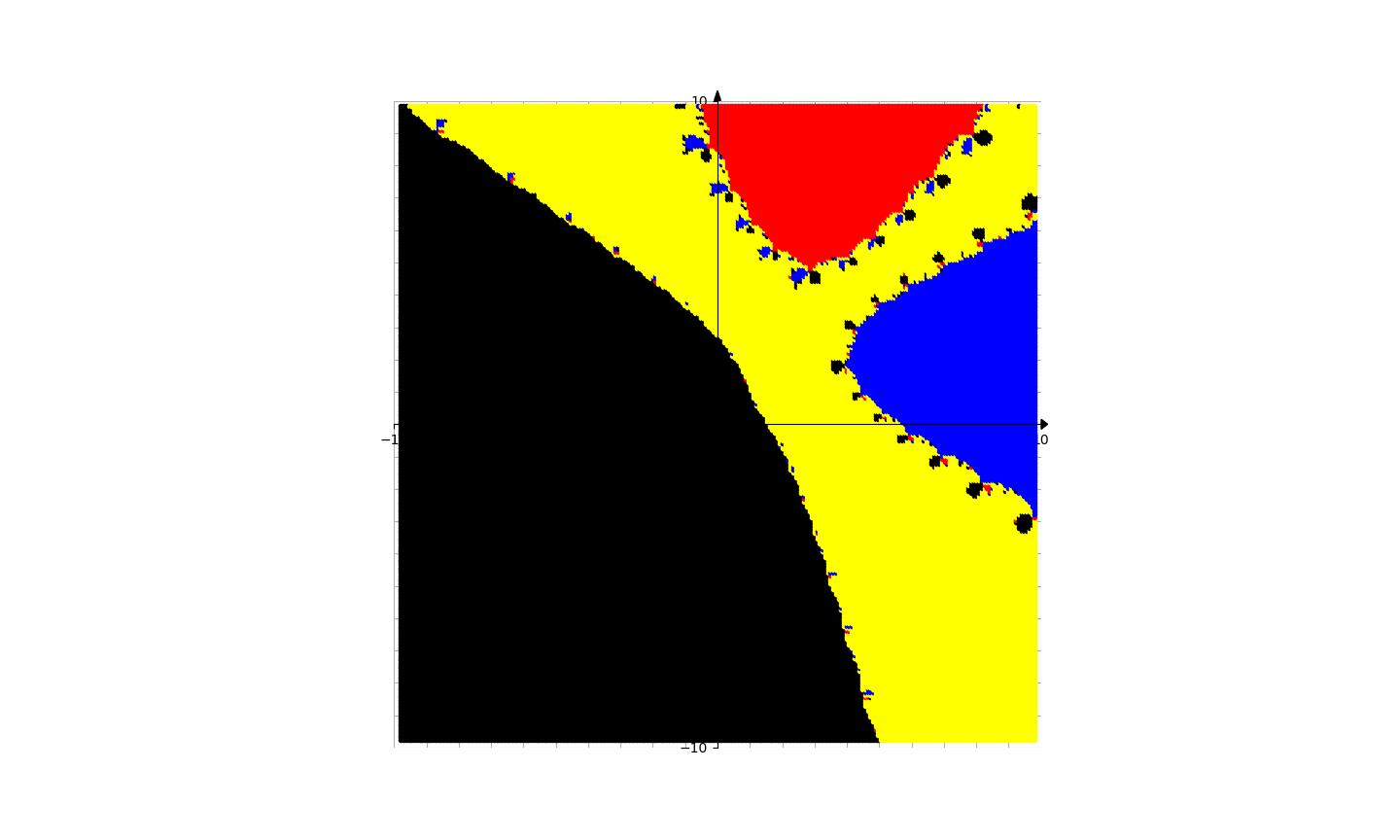}
        \includegraphics[width=4cm]{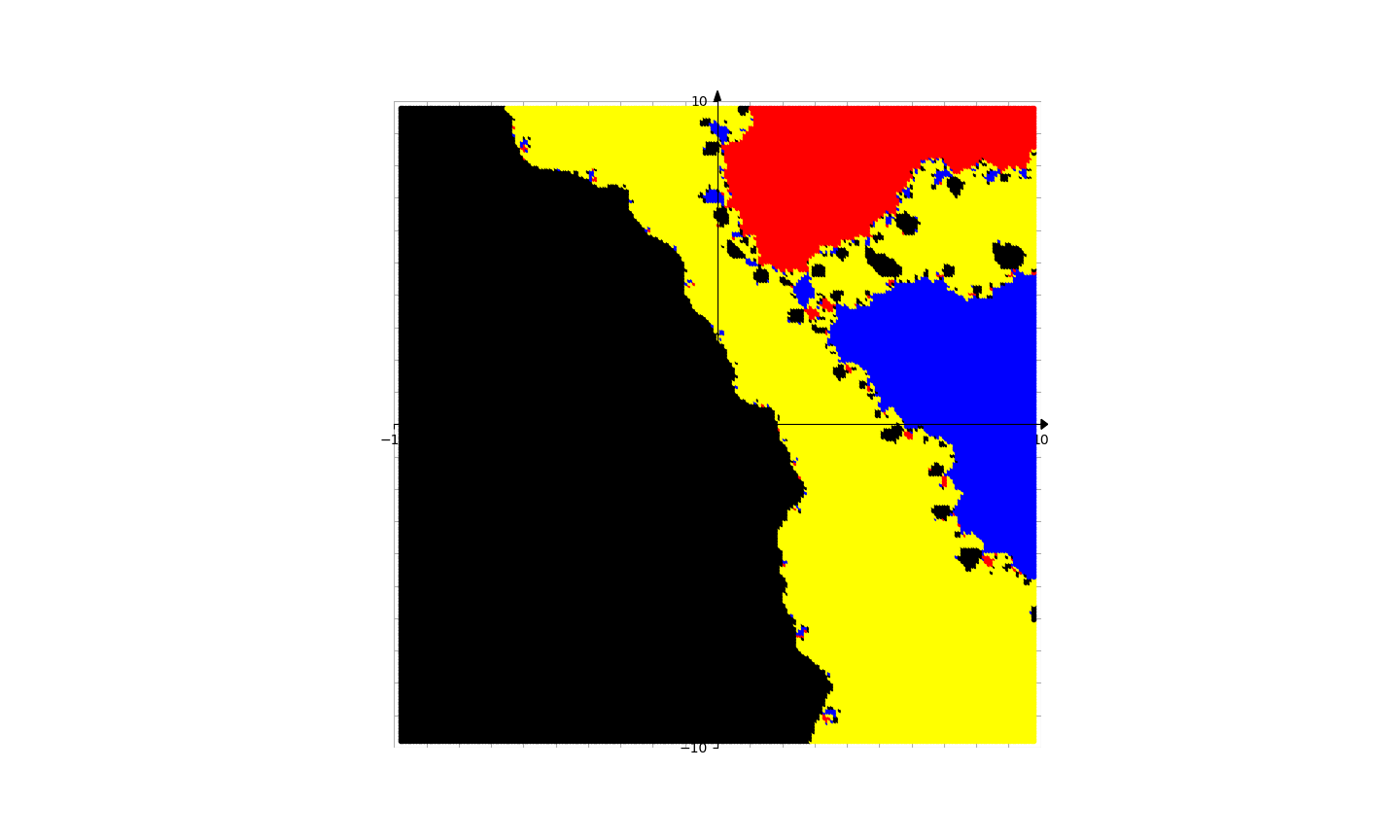}
   
    \bigskip
    
     \includegraphics[width=4cm]{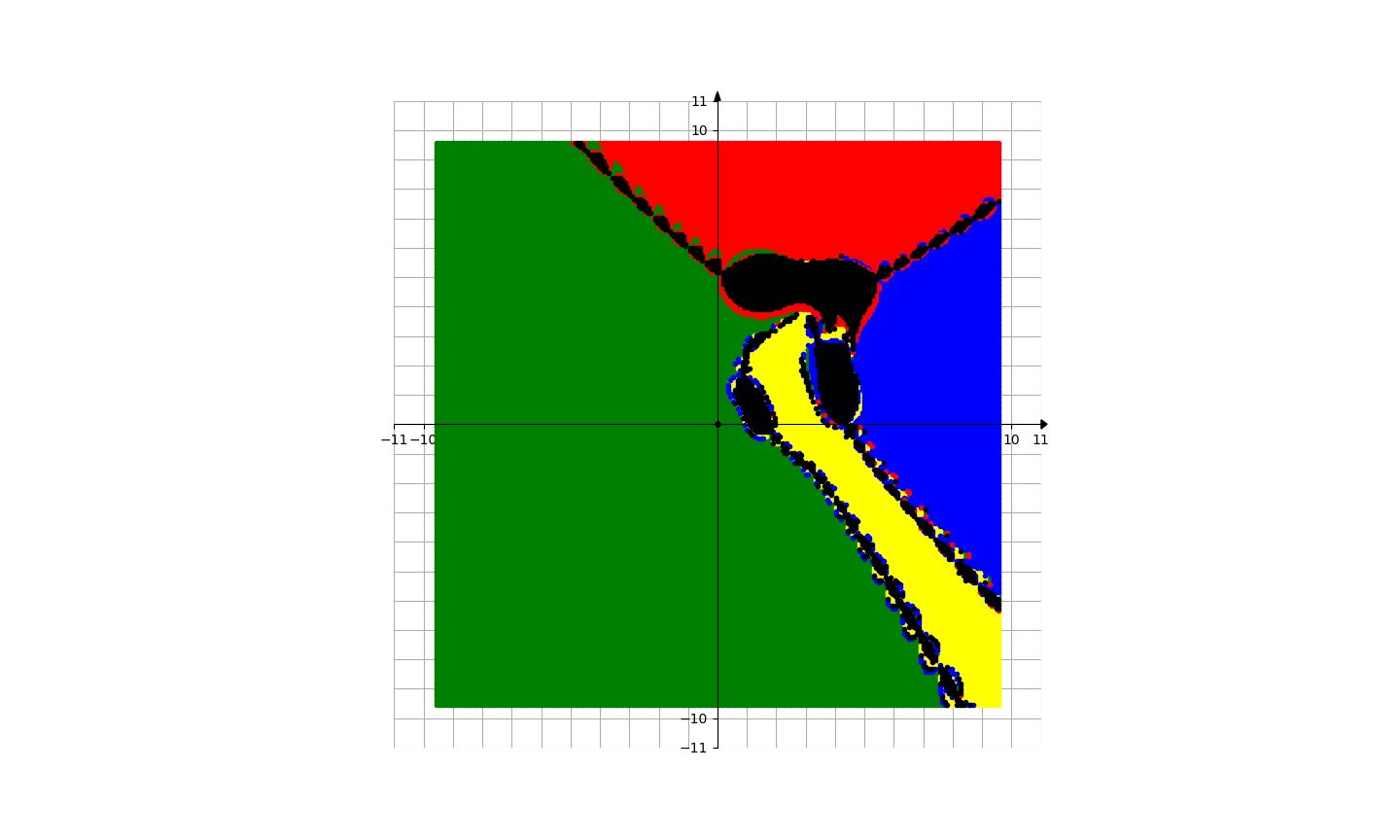}
      \includegraphics[width=4cm]{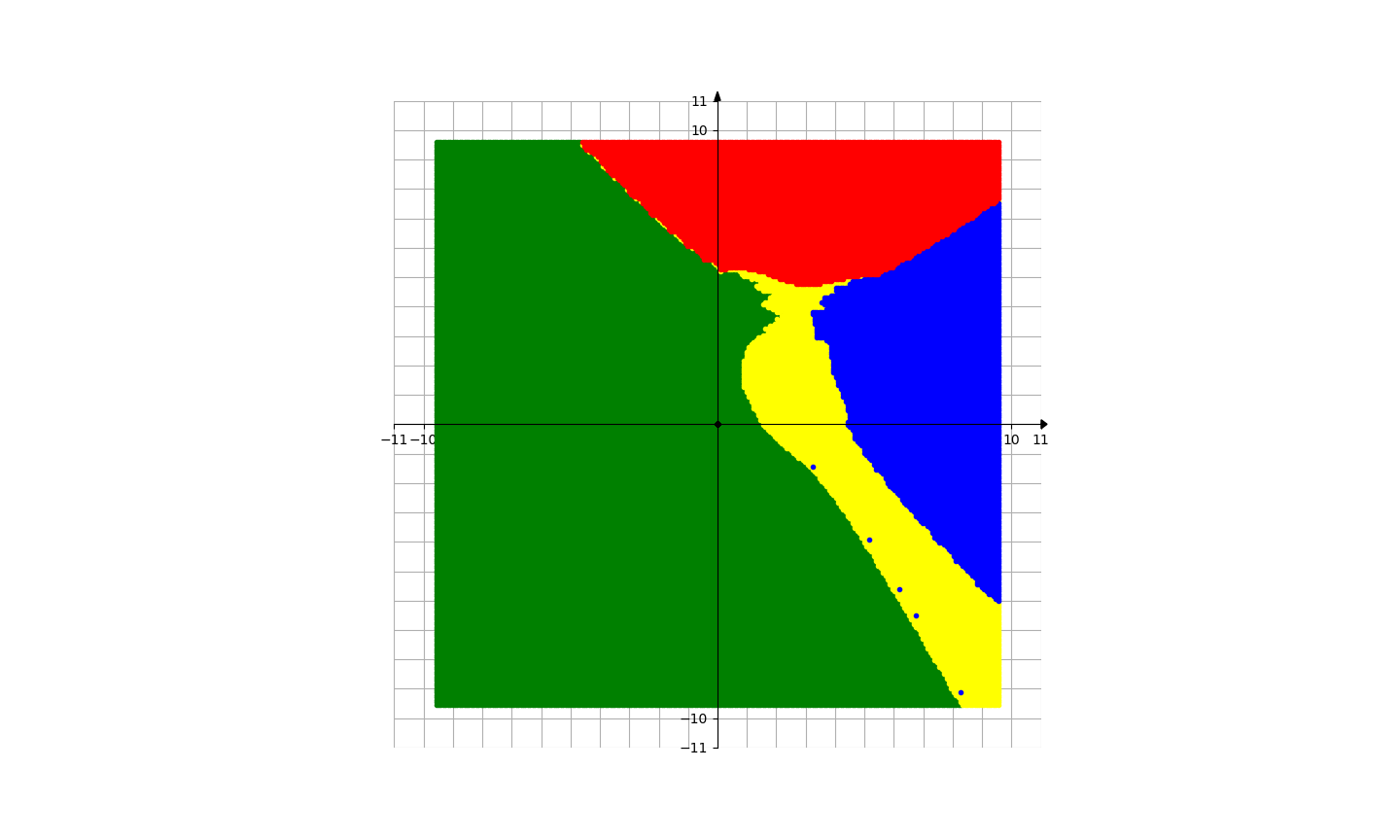}
     \includegraphics[width=4cm]{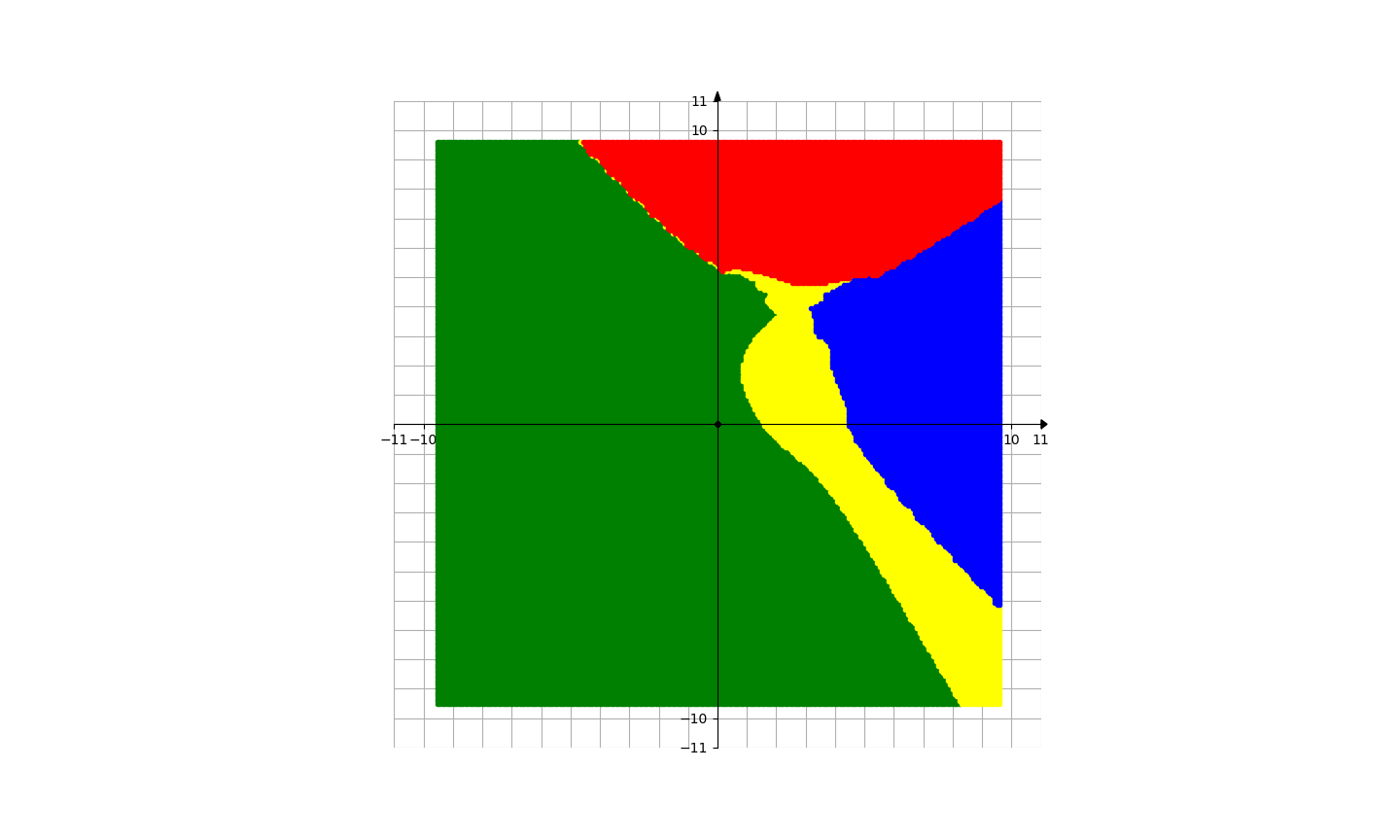}

    \caption{Basins of attraction for finding roots of the stochastic function $f_{21}+\epsilon \xi (z^3+2z-5)$ by different methods. Pictures are referenced to from top to bottom, from left to right. Row 1: left picture is for Newton's method, right picture is for Random Relaxed Newton's method. Row 2: left picture is for Newton's method vOptimization, central picture is for BNQN, right picture is for BNQN v2.}
    \label{fig:f21Stochastic}
\end{figure}

\section{Conclusions and further directions}

From experiments, we draw the following conclusions, when finding the roots of $f$: 

- Newton's method does not reflect well the geometric configuration of the roots of $f$ (i.e. Voronoi's diagrams). Newton's flow method, while more smooth than Newton's method, also does not reflect well the geometric configuration of the roots. The basins of attraction for Random Relaxed Newton's method are distorted versions of Newton's method. 

Both Newton's method and Random Relaxed Newton's method seem not work well (take too long time to run) with functions of the form $f/f'$, where $f$ is transcendental. A reason could be that in these cases, it can take many iterates for these methods to converge, or they are more likely to diverge, in both cases it takes more time to verify the STOP conditions. 

- Newton's method for Optimization performs poorly.

 Newton's flow for Optimization, except for some small open sets which are basins of attraction of the critical points of $f$, produces pictures which match well with Voronoi's diagrams. This is in big contrast to its discrete version. However, the method is quite heavy, and takes too long time for the transcendental examples.   

- Newton's flow for $f/f'$, on the other hand, reflects well Voronoi's diagrams, in contrast to Newton's flow for $f$. The reason could be as follows: The RHS in Newton's flow for $f/f'$ is 
\begin{eqnarray*}
(f/f')/(f/f')'=\frac{ff'}{(f')^2-ff"},
\end{eqnarray*}
which involves the second derivative and hence reflects better the curvature of the function landscape. 

- BNQN reflects well Voronoi's diagram. In particular, BNQN v2 (where we choose $\theta =1$) is more similar to Newton's flow for $f/f'$. 

- If some of the roots are contained inside the convex hull of other roots, then the similarity between these iterative methods and Voronoi's diagrams becomes less apparent. A reason could be that indeed one needs to use a more appropriate metric than the usual Euclidean metric when constructing Voronoi's diagrams. 

- BNQN is more robust against stochastic root finding than Newton's method, Random Relaxed Newton's method and Newton's method for Optimization. 

The above conclusions suggest a couple of natural follow ups, which are beyond the current paper and are left for future explorations: 

- While Newton's flow for $f/f'$ and BNQN both involve the second derivatives of $f$, how they depend on the second derivatives is different. In BNQN, we need to change the signs of negative eigenvalues, and we also need to run an Armijo's Backtracking line search. What is the reason that they produce similar pictures?

- Is it true that in BNQN New Variant, if we choose bigger values for the parameter $\theta$, then we will obtain more smooth pictures? 

- What should be the metric to be chosen in constructing Voronoi's diagrams, so that we will obtain more similar pictures to basins of attraction? 

- For the case of functions with multiple roots, are there  corresponding Voronoi's diagrams to the basins of attraction?

\end{document}